\documentclass[reqno]{amsart}
\usepackage{amsthm,amsfonts,amssymb,euscript}
\usepackage{latexsym, multicol, fancybox}
\usepackage{graphicx}
\usepackage{color}
\usepackage{amsmath, amsthm, amssymb, bm}
\usepackage{epstopdf}
\usepackage{caption}
\usepackage{psfrag}

\usepackage{mathrsfs}
 \usepackage{xcolor}
 \usepackage[citebordercolor={green}]{hyperref}
 \usepackage{tikz}
 \usepackage{bbm}
 \usepackage{dsfont}
 \usepackage{bbold}
 
\numberwithin{equation}{section}

 \usepackage{tikz}

\newtheorem{theorem}{Theorem}[section]
\newtheorem{lemma}[theorem]{Lemma}
\newtheorem{proposition}[theorem]{Proposition}
\newtheorem{corollary}[theorem]{Corollary}
\newtheorem{definition}[theorem]{Definition}
\newtheorem{remark}[theorem]{Remark}
\newtheorem{conjecture}[theorem]{Conjecture}

\newcommand{\bea}{\begin{eqnarray}}
\newcommand{\eea}{\end{eqnarray}}
\def\beaa{\begin{eqnarray*}}
\def\eeaa{\end{eqnarray*}}
\def\ba{\begin{array}}
\def\ea{\end{array}}
\def\be#1{\begin{equation} \label{#1}}
\def \eeq{\end{equation}}

\newcommand{\bsub}{\begin{subequations}}
\newcommand{\esub}{\end{subequations}}

\newcommand{\nn}{\nonumber}

\def\a{{\alpha}}

\def\b{{\beta}}
\def\be{{\beta}}
\def\ga{\gamma}
\def\Ga{\Gamma}
\def\de{\delta}
\def\De{\Delta}
\def\ep{\epsilon}

\def\la{\lambda}
\def\La{\Lambda}

\def\Si{\Sigma}
\def\om{\omega}
\def\Om{\Omega}

\def\vphi{\varphi}

\def\th{\theta}

\def\ze{\zeta}

\def\nab{\nabla}

\def\pr{{\partial}}
\def\les{\lesssim}
\def\c{\cdot}

\def\AA{{math\cal A}}

\def\MM{{\mathcal M}}
\def\NN{{\mathcal N}}

\def\LL{{\mathcal L}}
\def\II{{\mathcal I}}

\def\HH{{\mathcal H}}

\def\JJ{{\mathcal J}}
\def\KK{{\mathcal K}}
\def\Lie{{\mathcal L}}

\def\DD{{\mathcal D}}
\def\PP{{\mathcal P}}
\def\RR{{\mathcal R}}
\def\QQ{{\mathcal Q}}
\def\AA{{\mathcal A}}

\def\HH{{\mathcal H}}

\def\Lie{{\mathcal L}}

\def\lap{{\triangle}}

\def\A{{\bf A}}
\def\B{{\bf B}}

\def\D{{\bf D}}
\def\E{{\bf E}}
\def\F{{\bf F}}

\def\H{{\bf H}}

\def\M{{\bf M}}
\def\L{{\bf L}}
\def\O{{\bf O}}

\def\R{{\bf R}}

\def\K{{\bf K}}

\def\T{T}
\def\Z{Z}

\def\g{{\bf g}}

\def\CCC{{\Bbb C}}
\def\f12{{\frac 1 2}}

\def\dual{{\,\,^*}}
\def\div{{\mbox div\,}}
\def\curl{{\mbox curl\,}}

\def\Hb{\,\underline{H}}

\def\Xh{\,^{(h)}X}

\def\trch{{\mbox tr}\chi}
\def\chih{{\widehat \chi}}
\def\chib{{\underline \chi}}
\def\chibh{{\underline{\chih}}}

\def\etab{{\underline \eta}}
\def\omb{{\underline{\om}}}
\def\bb{{\underline{\b}}}
\def\aa{\protect\underline{\a}}
\def\xib{{\underline \xi}}

\def\Xib{\underline{\Xi}}

\def\Ab{\protect\underline{A}}
\def\Bb{\protect\underline{B}}
\def\Xb{\protect\underline{X}}

\def\Xh{\widehat{X}}
\def\Xbh{\widehat{\Xb}}

\def\tr{\mbox{tr}}
\def\atr{\,^{(a)}\mbox{tr}}

\def\trchb{{\tr\chib}}

\def\Div{\mbox{Div}}

\def\Hh{{\widehat H}}

\def\atrch{\atr\chi}
\def\atrchb{\atr\chib}

\def\hot{\widehat{\otimes}}
\def\rhod{\,\dual\hspace{-2pt}\rho}

\def\piX{\, ^{(X)}\pi}

\def\err{{\mbox{Err}}}
\def\ov{\overline}

\def\f12{\frac 1 2}
\def\lab{\label}
\def\nabc{\,^{(c)}\nab}

\def\bsplit{\begin{split}}

\newcommand{\Mext}{{\,{}^{(ext)}\mathcal{M}}}

\newcommand{\Mint}{{\protect \,{}^{(int)}\mathcal{M}}}

\def\Mint{{\protect \, ^{(int)}\MM}}

\def\qf{\mathfrak{q}}
\def\qfb{\protect\underline{\qf}}

\def\Jk{\mathfrak{J}}

\def\sk{\mathfrak{s}}
\def\dkb{ \, \mathfrak{d}     \mkern-9mu /}
\def\dk{\mathfrak{d}}

\DeclareFontFamily{U}{mathx}{\hyphenchar\font45}
\DeclareFontShape{U}{mathx}{m}{n}{
      <5> <6> <7> <8> <9> <10>
      <10.95> <12> <14.4> <17.28> <20.74> <24.88>
      mathx10
      }{}
\DeclareSymbolFont{mathx}{U}{mathx}{m}{n}
\DeclareFontSubstitution{U}{mathx}{m}{n}
\DeclareMathAccent{\widecheck}{0}{mathx}{"71}

\def\Zc{\widecheck{Z}}
\def\Hc{\widecheck{H}}
\def\Hbc{\widecheck{\Hb}}
\def\trXc{\widecheck{\tr X}}
\def\trXbc{\widecheck{\tr\Xb}}
\def\Pc{\widecheck{P}}

\def\Gac{\widecheck{\Ga}}

\def\omc{\widecheck \omega}

\def\DDc{\,^{(c)} \DD}

\newcommand{\deh}{\delta_{\mathcal{H}}}
\newcommand{\dec}{\delta_{dec}}

\newcommand{\Lieb}{\Lie \mkern-10mu /\,}

\def\Rdot{\dot{\R}}

\def\DDc{\,^{(c)} \DD}

\def\DDb{\ov{\DD}}
\def\DDbc{\ov{\DDc}}

\def\Lied{\dot{\Lie}}

\def\Db{\dot{\D}}
\def\Ddot{\dot{\D}}
\def\squared{\dot{\square}}

\def\Ddot{\dot{\D}}

\def\Ddot{\dot{\D}}

\def\Rhat{{\widehat{R}}}

\def\nabc{\,^{(c)}\nab}

\DeclareFontFamily{U}{mathx}{\hyphenchar\font45}
\DeclareFontShape{U}{mathx}{m}{n}{
      <5> <6> <7> <8> <9> <10>
      <10.95> <12> <14.4> <17.28> <20.74> <24.88>
      mathx10
      }{}
\DeclareSymbolFont{mathx}{U}{mathx}{m}{n}
\DeclareFontSubstitution{U}{mathx}{m}{n}
\DeclareMathAccent{\widecheck}{0}{mathx}{"71}

\def\trXc{\widecheck{\tr X}}
\def\trXbc{\widecheck{\tr\Xb}}
\def\Hc{\widecheck{H}}
\def\Zc{\widecheck{Z}}
\def\Pc{\widecheck{P}}

      \def\ntrap{trap\mkern-18 mu\big/\,}
          \def\Mtrap{\,\MM_{trap}}
\def\Mntrap{{\MM_{\ntrap}}}

    \def\DDc{\,^{(c)} \DD}

      \def\ntrap{trap\mkern-18 mu\big/\,}
\def\Mntrap{{\MM_{\ntrap}}}

 \def\NNt{\widetilde{\NN}}

    \def\DDc{\,^{(c)} \DD}

\def\und{\underline}

\def\N{\mathbf{N}}

\newcommand{\psiplus}[1]{\pmb\psi_{+2}^{(#1)}}
\newcommand{\psiminus}[1]{\pmb\psi_{-2}^{(#1)}}
\newcommand{\psis}[1]{\pmb\psi_{s}^{(#1)}}
\newcommand{\psiss}[2]{\psi_{s,#1}^{(#2)}}

\newcommand{\phiplus}[1]{\pmb\phi_{+2}^{(#1)}}
\newcommand{\phiminus}[1]{\pmb\phi_{-2}^{(#1)}}
\newcommand{\phis}[1]{\pmb\phi_{s}^{(#1)}}

\newcommand{\phipluss}[2]{\phi_{+2,#1}^{(#2)}}
\newcommand{\phiminuss}[2]{\phi_{-2,#1}^{(#2)}}
\newcommand{\phiss}[2]{\phi_{s,#1}^{(#2)}}

\def\gam{\g_{a,m}}
\def\dhor{\delta_{\HH}}
\def\dbl{\delta_{\textbf{BL}}}
\def\tmod{t_{\text{mod}}}
\def\phimod{\phi_{\text{mod}}}
\def\tt{\tau}

\def\tphi{\varphi}

\def\Xcal{\mathcal{X}}

\def\dred{\delta_{\text{red}}}

\def\EMF{{\bf EMF}}
\def\EM{{\bf EM}}
\def\EF{{\bf EF}}

\def\EF{{\bf EF}}

\def\NNh{\mathring{\NN}}
\def\MMh{\mathring{\MM}}
\def\Mtraph{\,\mathring{\MM}_{trap}}
\def\Mntraph{{\mathring{\MM}_{\ntrap}}}
\def\gh{\mathring{\g}}
\def\Sih{\mathring{\Sigma}}
\def\Eh{{\bf \mathring{E}}}
\def\Fh{{\bf \mathring{F}}}
\def\Mh{{\bf \mathring{M}}}
\def\EMFh{{\bf \mathring{EMF}}}
\def\EMh{{\bf \mathring{EM}}}
\def\EFh{{\bf \mathring{EF}}}
\def\MFh{{\bf \mathring{EMF}}}

\def\NNt{\widetilde{\NN}}
\def\NNtlede{\NNt_{\text{le},\de}}

\newcommand{\Rmic}{R_0}
\newcommand{\Nmic}{N_0}
\newcommand{\tmic}{\tau_{\Nmic}}
\newcommand{\Iti}{I_{\Nmic}}

\def\gcheck{\widecheck{\g}}

\def\qs{|q|^2}

\def\tauu{\underline{\tau}}
\def\tauut{\widetilde{\tauu}}

\def\reg{\mathbf{k}}

\def\Reals{\mathbb{R}}

\def\Errdefect{\E_{\mathrm{defect}}}
\def\Errdefects{\F_{\mathrm{defect,*}}}

\def\NNtlocal{\NNt_{\text{local}}}
\def\Err{{\bf{Err}}}

\newcommand{\Bulkxw}[1]{\mathbf{Bulk}_{#1,(X,w)}[\pmb\psi]}

\def\NNtaux{\widetilde{\NN}_{\text{aux}}}
\def\NNtmora{\NNt_{\text{Mora}}}
\def\NNtdemora{\NNt_{\text{Mora}, \de}}
\def\NNtener{\NNt_{\text{Ener}}}

\newcommand{\IE}[1]{{\bf{IE}}[#1]}

\newcommand{\EMFtotalh}[2]{\widetilde{\EMF}_{#1,\de, \text{total}}[#2\pmb\phi_{#1}]}
\newcommand{\EMFtotalhps}[2]{\EMF_{#1,\de, \text{total}}[#2\pmb\phi_{#1}]}

\newcommand{\NNttotalph}[2]{\NNt_{#1, \de, \text{total}}[#2\pmb\phi_{#1}]}

\newcommand{\EMFtotalp}[1]{\widetilde{\EMF}_{#1, \de, \text{total}}[\pmb\phi_{#1}]}

\newcommand{\EMFtotalps}[1]{{\EMF}_{#1, \de, \text{total}}[\pmb\phi_{#1}]}

\newcommand{\NNttotalp}[1]{\NNt_{#1,\de, \text{total}}[\pmb\phi_{#1}]}

\newcommand{\psish}[2]{\pmb\psi_{s}^{(#1), #2}}

\newcommand{\IEde}[1]{{\bf{IE}_{\de}}[#1]}

\def\prtphihat{\widehat{\pr}_{\tphi}}

\def\Ao{\protect\overline{\A}}

\begin{document}

\title{Teukolsky equations in perturbations of Kerr}
\author{J\'{e}r\'{e}mie Szeftel}

\begin{abstract}
The Kerr stability conjecture has been proved in the slowly rotating case, i.e., $|a|\ll m$, in the sequence of works \cite{KS-GCM1} \cite{KS-GCM2} \cite{KS:Kerr} by Sergiu Klainerman and the author, \cite{GKS22} by Elena  Giorgi, Sergiu Klainerman and the author, and \cite{Shen} by Dawei Shen, and extended in \cite{Sze}, by the author, to the full subextremal range $|a|<m$, thereby completing the proof of the Kerr stability conjecture. The proof in \cite{Sze} crucially relies on the two companion papers \cite{MaSz24} \cite{MaSz26} by Siyuan Ma and the author, in which we prove energy-Morawetz estimates respectively for the scalar wave equation and for Teukolsky equations on perturbations of Kerr with $|a|<m$. In addition, two results in \cite{Sze}, concerning the derivation of the Teukolsky wave-transport system in perturbations of Kerr and energy-Morawetz estimates for Teukolsky in the setting of \cite{Sze}, are stated without proofs. The goal of the present companion paper to \cite{Sze} is to provide the proof of these two results. 
\end{abstract}

\maketitle

\tableofcontents

%%%%%%%%%%%%%%%%%%%%

\section{Introduction}

%%%%%%%%%%%%%%%%%%%%%

%%%%%%%%%%%%%%%%%%%%%%%%%%%%%%%%%%%%%%%%%%%%%%

\subsection{Kerr stability conjecture and rough statement of the main result in \cite{Sze}}
  
%%%%%%%%%%%%%%%%%%%%%%%%%%%%%%%%%%%%%%%%%%%%%% 

%%%%%%%%%%%%%%%%%%%%%%%%%%%

\subsubsection{Einstein vacuum equations}

%%%%%%%%%%%%%%%%%%%%%%%%%%%

The Einstein vacuum equations (EVE) in a {Lorentzian manifold} $(\MM, \g)$ take the form
\bea\lab{eq:EVE:intro}
\mathbf{R}_{\a\b}=0,
\eea
where $\mathbf{R}_{\a\b}$ {denotes} the Ricci curvature tensor of the metric $\g$. Foundational contributions by Choquet-Bruhat \cite{CB52}, and by  
Choquet-Bruhat and Geroch \cite{CBG69}, formulate EVE as an evolution problem of hyperbolic type and associate to any suitable sufficiently regular initial data set a unique (up to diffeomorphisms) maximal Cauchy development.

%%%%%%%%%%%%%%%%%%%%%%%%%%%

\subsubsection{Kerr solution}

%%%%%%%%%%%%%%%%%%%%%%%%%%%

The EVE admit a family of explicit solutions, found by Kerr \cite{Kerr63} in 1963, which describe asymptotically flat, stationary, axially symmetric black hole spacetimes. The metrics of Kerr spacetimes are parameterized by an angular momentum per unit mass $a$ and a mass $m$,  satisfying $|a|\leq m$, and take the following form in the Boyer--Lindquist coordinates $(t,r,\th, \phi)$
\bea\lab{eq:expressionofKerrmetricinBLcoordinates:intro}
\gam=-\frac{\Delta \qs}{\Sigma^2} dt^2 + \frac{\sin^2\th\Sigma^2 }{\qs}\bigg(d\phi - \frac{2amr}{\Sigma^2} dt\bigg)^2 +\frac{\qs}{\Delta} dr^2 + \qs d\th^2,
\eea
where
\bea
\Delta = r^2 - 2mr +a^2, \quad \qs=r^2+a^2\cos^2\th, \quad \Sigma^2=(r^2+a^2)^2 - a^2\sin^2\th \Delta.
\eea
Note that the particular case $a=0$ with $m>0$ corresponds to the family of Schwarzschild spacetimes, introduced by Schwarzschild \cite{Sch16} in 1916.

We consider in this work the family of \textit{subextremal} Kerr spacetimes, in which the two parameters $(a,m)$ satisfy the strict inequality $|a|<m$. Such a subextremal Kerr spacetime contains a black hole $\{r<r_+\}$ with a nondegenerate event horizon located at $\{r=r_+\}$ where $r_+:=m+\sqrt{m^2-a^2}$ is the larger root of $\Delta=\De(r)$, see Figure \ref{fig:penrosediagramofKerr} for the corresponding Penrose diagram.

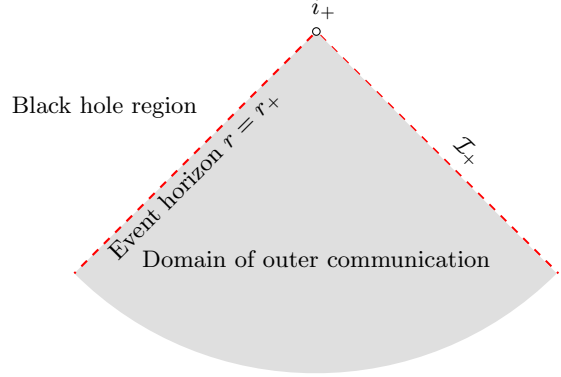
\begin{figure}[htbp]
  \begin{center}
\begin{tikzpicture}[scale=1]
\tikzstyle{every node}=[font=\small]
      \draw[dashed, color=red, thick] (0.05,3.95) -- (3.2,0.8);
  \fill[lightgray!50] (-0.05,3.95)--(-3.2,0.8) arc(225:315: 4.51 and 4.51) -- (0.05,3.95);
   \draw[dashed, color=red, thick] (-0.05,3.95)--(-3.2,0.8) ;
     \node at (-2.8,3) {Black hole region};
     \node at (0.1,4.3) {$i_+$};
       \node[rotate=315]  at (2.0,2.4) {$\II_+$};
       \node[rotate=45] at (-1.6,2.1) {Event horizon $r=r_+$};
         \draw[] (0,4) circle (0.05);
         \node at (0,1) {Domain of outer communication};  
\end{tikzpicture}
\end{center}
\caption{\footnotesize{Penrose diagram of subextremal Kerr spacetimes.}}
\lab{fig:penrosediagramofKerr}
\end{figure}

%%%%%%%%%%%%%%%%%%%%%%%%%%%

\subsubsection{Kerr stability conjecture}

%%%%%%%%%%%%%%%%%%%%%%%%%%%

The \textit{black hole stability conjecture} is one of the central open problems in general relativity. We provide a rough statement below.

\begin{conjecture}[Kerr stability conjecture]
The maximal Cauchy development of any initial data set for EVE, that is sufficiently close to a subextremal Kerr initial data set in a suitable sense, has a complete future null infinity and a domain of outer communication\footnote{The domain of outer communication is the complement of the black hole region.} which is asymptotic to a nearby member of the subextremal Kerr family. 
\end{conjecture}

%%%%%%%%%%%%%%%%%%%%%%%%%%%%%%%%
  
\subsubsection{Statement of the main result in \cite{Sze}}

%%%%%%%%%%%%%%%%%%%%%%%%%%%%%%%%
 
The following theorem proves the Kerr stability conjecture.
   \begin{theorem}[Main theorem in \cite{Sze}, rough version]
\lab{MainThm-firstversion}
Let $a_0$ and $m_0$ be real constants such that $|a_0|<m_0$. The future globally hyperbolic   development  of  a general,   asymptotically  flat, initial data set, sufficiently close to a   $Kerr(a_0, m_0) $   initial data set in the sense that 
\bea\lab{eq:rmk:remarkonroughtopologyinitialdatamainTh}
\g=\g_{a_0,m_0}+O(\ep_0r^{-\frac{3}{2}-\de})\quad\textrm{as}\quad r\to +\infty\quad\textrm{along the spacelike initial hypersurface}\quad\Si_0,
\eea
for $\de>0$ and for $\ep_0>0$ small enough, has a complete    future null infinity  $\II^+$ and converges in  its causal past  $\JJ^{-1}(\II^{+})$  to another  nearby Kerr spacetime $Kerr(a_f, m_f)$ with parameters    $(a_f, m_f)$ close to the initial ones $(a_0, m_0)$.
 \end{theorem}

 \begin{figure}[ht!]
 \lab{fig0-introd}
\centering
\includegraphics[scale=0.35]{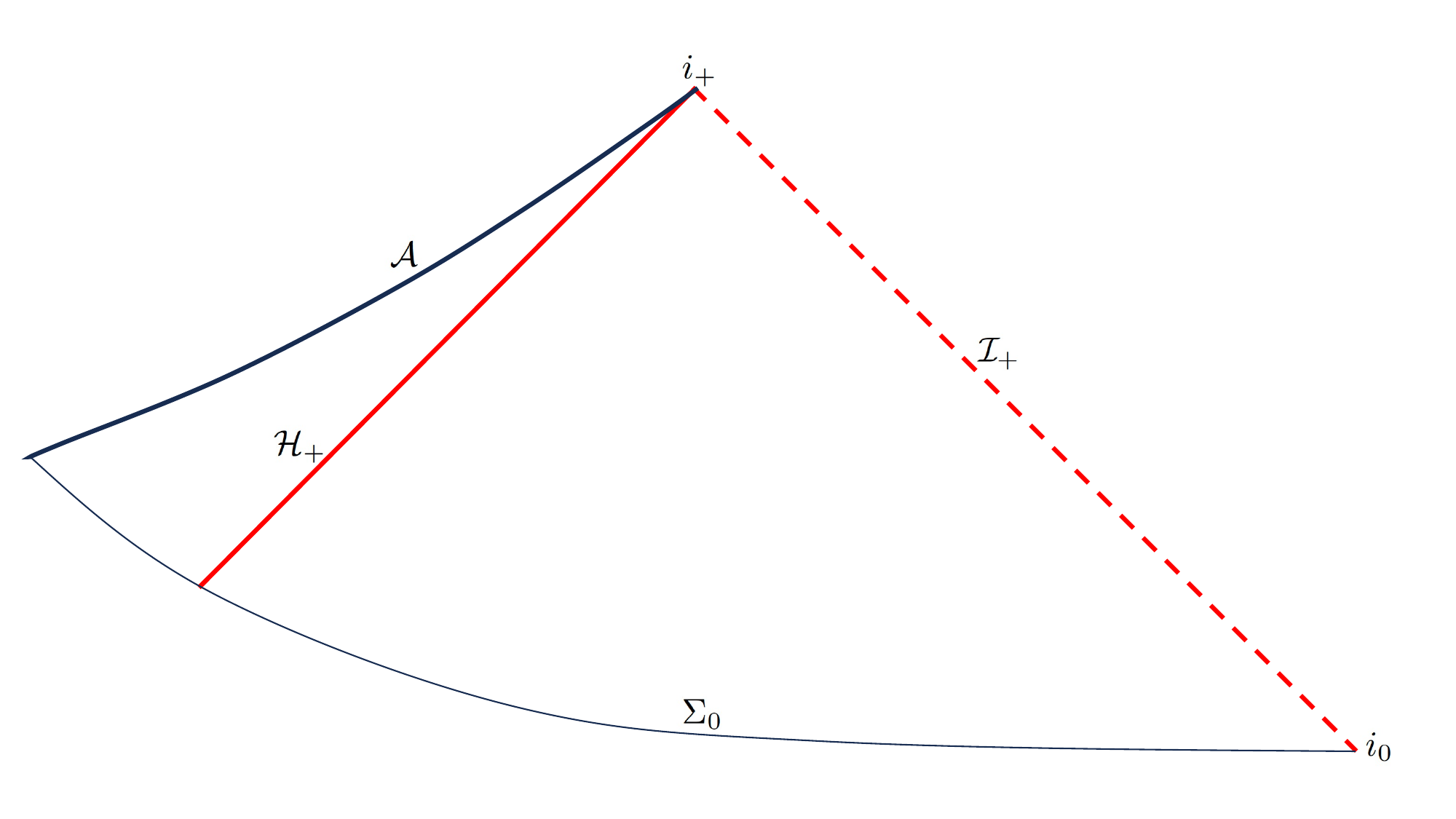}
\caption{\footnotesize{The Penrose diagram of the  final   space-time in  Theorem \ref{MainThm-firstversion} with initial spacelike hypersurface $\Si_0$,  future space-like boundary $\AA$, and $\II^+$ the complete future null infinity.  The hypersurface $\HH_+$} denotes the future event horizon.}
\end{figure} 

Prior to this work, the Kerr stability conjecture has been addressed in the case of Schwarzschild (i.e., for $a=0$) in polarized symmetry by Klainerman-Szeftel \cite{KS}, in the case of Schwarzschild for a codimension 3 subset of initial data by Dafermos-Holzegel-Rodnianski-Taylor \cite{DHRT}, and in the case of Kerr for $|a|\ll m$ in the sequence of works \cite{KS-GCM1} \cite{KS-GCM2} \cite{KS:Kerr} by Klainerman and the author, \cite{GKS22} by Giorgi, Klainerman and the author, and \cite{Shen} by Shen, where the initial data in \cite{KS} and \cite{GKS22} \cite{KS-GCM1} \cite{KS-GCM2} \cite{KS:Kerr} \cite{Shen} are as in \eqref{eq:rmk:remarkonroughtopologyinitialdatamainTh}, while the ones in \cite{DHRT} fall off like $O(\ep_0r^{-\frac{5}{2}})$ as $r\to +\infty$ along $\Si_0$. Recently, it has also been proved for all $|a|<m$ by Hintz \cite{Hi2} \cite{Hi3} \cite{Hi1} for a class of well-prepared initial data satisfying
\bea\lab{eq:intro:wellpreparedinitialdatainresultbyPeterHintz:0}
\g-\g_{a_0,m_0} = \textrm{finite polyhomogeneous expansion} + O(\ep_0r^{-3-\de})\quad\textrm{as}\quad r\to +\infty,
\eea
for $\de>0$ and for $\ep_0>0$ small enough, which allows for a reduction to finite-dimensional gauge modifications\footnote{More precisely, the initial data in  \eqref{eq:intro:wellpreparedinitialdatainresultbyPeterHintz:0} are tailored to avoid the difficulties associated with the infinite dimensional kernel of the linearized operator generated by general covariance, see Section 1.3.4 in \cite{Sze} and Remark 1.3 in \cite{Hi1} for more details.}. We refer to the introduction in our companion paper \cite{Sze} for a detailed review of the literature concerning the Kerr stability conjecture, and in particular to Section 1.3.4 in \cite{Sze} for a comparison of the result in \cite{Hi1} with Theorem \ref{MainThm-firstversion}. The rest of the introduction concerns Teukolsky equations which is the  focus of the present paper.

%%%%%%%%%%%%%%%%%%%%%%%%%%%%

\subsection{Teukolsky wave/transport system}
\lab{subsect:Teukolskywavetransport:intro}

%%%%%%%%%%%%%%%%%%%%%%%%%%%%

Kerr spacetimes possess a distinguished pair of null vectorfields known as the principal null pair, see \eqref{def:e3e4inKerr}, that diagonalizes the curvature tensor. In perturbations of Kerr, we consider a pair of null vectorfields $(e_3, e_4)$, normalized by $\g(e_3, e_4)=-2$, which is a suitable perturbation of the principal null pair of Kerr. We then consider an orthonormal pair of spacelike vectorfields $(e_1, e_2)$ spanning the horizontal bundle $\{e_3, e_4\}^\perp$, see \eqref{def:e1e2inKerr} in Kerr, so that $(e_3, e_4, e_a)$, $a=1,2$, forms a null frame of the spacetime $(\MM, \g)$. As in \cite{GKS20} \cite{GKS22}, we associate to the null pair $(e_3, e_4)$ horizontal tensors and denote in particular by $\sk_2(\Reals)$ the set of symmetric traceless horizontal real $2$-tensors, see Section \ref{subsection:review-horiz.structures}.

Next, we denote the curvature components $\a, \aa\in\sk_2(\Reals)$ by
\beaa
\a_{ab}={\bf R}_{a4b4},\qquad \aa_{ab}={\bf R}_{a3b3},\qquad a,b=1,2,
\eeaa
where ${\bf R}_{\a\b\mu\nu}$ denotes the curvature tensor of the spacetime $(\MM, \g)$. Also, we define the complexified curvature components $A, \Ab\in \sk_2(\mathbb{C})$ as
\beaa
A=\a+i\dual\a, \qquad \Ab=\aa+i\dual\aa,
\eeaa
where $\sk_2(\mathbb{C})$ is introduced in Definition \ref{def:skC:horizontaltensors}  as the set of symmetric, traceless, anti-self-dual horizontal complex $2$-tensors. The \textit{Teukolsky equations} \cite{Teuk}, the governing equations for these curvature components, are, in the tensorial formalism introduced in \cite{GKS20} \cite{GKS22}, given by
\bea
\lab{eq:Teu:intro}
\mathcal{T}_{+2, \g} A = \N_A, \qquad \mathcal{T}_{-2,\g}\Ab=\N_{\Ab},
\eea
where $\mathcal{T}_{\pm 2,\g}$ are tensorial Teukolsky wave operators in $(\MM, \g)$ and where $\N_{A}$ and $\N_{\Ab}$ are source terms\footnote{For the explicit formulas of $\N_{A}$ and $\N_{\Ab}$ in terms of the Ricci coefficients and curvature components of a perturbation of Kerr $(\MM, \g)$, see Sections 5.1.1 and 5.3.1 in \cite{GKS22}.}. 

The heart of the analysis in this paper relies on a \textit{wave/transport hierarchy} constructed from the Teukolsky equations \eqref{eq:Teu:intro}. 
Following \cite{Ma}, {with the construction} adapted to the tensorial formalism of \cite{GKS22},  we consider  tensors $\pmb\phi_s^{(p)}\in\sk_2(\mathbb{C})$, $s=\pm 2$\footnote{In this paper, $s$ refers to the spin weight of the tensors.}, $p=0,1,2$, with $\pmb\phi_s^{(0)}$ given by 
\bsub
\lab{def:TensorialTeuScalars:wavesystem:Kerrperturbation:intro}
\bea
\pmb\phi_{+2}^{(0)}=\frac{\ov{q}}{q}A, \qquad \pmb\phi_{-2}^{(0)}=\frac{q}{\ov{q}}\left(\frac{\De}{|q|^2}\right)^2\Ab
\eea
and with $\phis{p}$, $s=\pm 2$, $p=0,1,2$, satisfying the Teukolsky transport equations
 \bea
\lab{def:TensorialTeuScalars:wavesystem:Kerrperturbation:+2:intro}
\nab_3 \left(\frac{r\bar{q}}{q}\left(\frac{r^2}{|q|^2}\right)^{p-2}\pmb\phi_{+2}^{(p)}\right)&=&\frac{\bar{q}}{rq}\left(\frac{r^2}{|q|^2}\right)^{p-1}\pmb\phi_{+2}^{(p+1)}+\N_{T,+2}^{(p)}, \quad p=0,1,\\
\lab{def:TensorialTeuScalars:wavesystem:Kerrperturbation:-2:intro}
\nab_4\left(\frac{rq}{\bar{q}}\left(\frac{{r^2}}{|q|^2}\right)^{p-2}\pmb\phi_{-2}^{(p)}\right)&=&\frac{q}{r\bar{q}}\left(\frac{r^2}{|q|^2}\right)^{p-1}\frac{\De}{\qs}\pmb\phi_{-2}^{(p+1)}+\N_{T,-2}^{(p)}, \,\,\,\, p=0,1,
\eea
\esub
where $q:=r+ia\cos\th$ and where $\N_{T, s}^{(p)}$, $s=\pm 2$, $p=0,1$, are source terms in the transport equations. These tensors $\phis{p}$ satisfy the coupled Teukolsky wave equations
\bsub
\lab{eq:TensorialTeuSysandlinearterms:rescaleRHScontaine2:general:Kerrperturbation:intro}
\bea
\lab{eq:TensorialTeuSys:Kerrpert:intro}
\bigg(\squared_2 -\frac{4ia\cos\th}{|q|^2}\nab_{\pr_t}- \frac{4-2\de_{p0}}{\qs}\bigg){\phis{p}} = \L_{s}^{(p)}[\pmb\phi_{s}]+\N_{W,s}^{(p)}, \quad s=\pm2, \quad p=0,1,2,
\eea
where $\squared_2$ is the tensorial wave operator for tensors in $\sk_2(\mathbb{C})$, see \eqref{eq:def=squared-2}, and where the linear coupling terms $\L_{s}^{(p)}[{\pmb\phi_s}]$ have the following schematic form
\bea
\lab{eq:tensor:Lsn:onlye_2present:general:Kerrperturbation:intro}
\bsplit
{\L_{s}^{(0)}[\pmb\phi_{s}]}={}& (2sr^{-3} +O(mr^{-4}))\phis{1}+ O(mr^{-3}) \nab_{\Xcal_s}^{\leq 1}\phis{0},\\
{\L_{s}^{(1)}[\pmb\phi_{s}]}={}& (sr^{-3} +O(mr^{-4}))\phis{2}+ O(mr^{-3}) \nab_{\Xcal_s}^{\leq 1}  \phis{1}+O(mr^{-2})\nab_{\pr_{\phi}+a\pr_t}^{\leq 1}\phis{0},\\
{\L_{s}^{(2)}[\pmb\phi_{s}]}={}&O(mr^{-3})\phis{2}+O(mr^{-2})\nab_{\pr_{\phi}+a\pr_t}^{\leq 1}\phis{1}+O(m^2 r^{-2})\phis{0},
\end{split}
\eea
\esub
with $\Xcal_s$, $s=\pm 2$, being regular vectorfields that are horizontal in the case of Kerr\footnote{See \eqref{eq:formofregularhorizontalvectorfieldmathcalXs:Kerrperturbation} for the form of $\Xcal_s$, $s=\pm 2$.}. The equations \eqref{def:TensorialTeuScalars:wavesystem:Kerrperturbation:intro} \eqref{eq:TensorialTeuSysandlinearterms:rescaleRHScontaine2:general:Kerrperturbation:intro} correspond to the tensorial \textit{Teukolsky wave/transport system} in perturbations of Kerr considered throughout this paper.

%%%%%%%%%%%%%%%%%%%%%%%%%%%%%%%%%%%%%%%%%%%%%%%%%

\subsection{State of the art on energy-Morawetz estimates for Teukolsky equations}
\lab{subsect:literature:Teukolsky:intro}

%%%%%%%%%%%%%%%%%%%%%%%%%%%%%%%%%%%%%%%%%%%%%%%%%

The analysis of the Teukolsky equations is central to understanding the dynamical evolution of Kerr spacetimes and fundamentally builds upon the framework developed for the scalar wave equation. For an in-depth review of the literature concerning scalar waves, we direct the reader to the introduction in our previous  paper \cite{MaSz24}. In this section, we review the literature pertaining to energy-Morawetz estimates for solutions to Teukolsky equations.

%%%%%%%%%%%%%%%%%%%%%%%%%%

\subsubsection{Teukolsky equations in Kerr}
\lab{subsubsect:literature:TeukolskyinKerr:intro}

%%%%%%%%%%%%%%%%%%%%%%%%%%

In order to derive energy-Morawetz estimates, one must first address the question of mode stability for solutions to Teukolsky equations in Kerr spacetimes. The absence of exponentially growing mode solutions was proved in the seminal work of Whiting \cite{Whit}. It was later extended in \cite{AMPW17}, see also \cite{TdC20}, to show the absence of non-trivial mode solutions with real frequencies. Recently, an unconditional bound for the horizon flux using a generalized physical space version of Whiting's transform has been derived in \cite{HeKl2}.

In Schwarzschild,  energy-Morawetz estimates for the Teukolsky equations were first obtained by Dafermos-Holzegel-Rodnianski \cite{D-H-R}. The proof relies on a physical-space analog of the Chandrasekhar's transformation \cite{Chand2} that converts the Teukolsky equations into a Regge-Wheeler type wave equation \cite{RW57}, to which the techniques developed for the scalar wave equation can be directly applied. Generalizations to Kerr spacetimes were achieved in the slowly rotating case by Ma \cite{Ma} and Dafermos-Holzegel-Rodnianski \cite{D-H-R-Kerr}, and for the full subextremal range by Millet \cite{Millet}\footnote{While \cite{Millet} derives sharp decay estimates for solutions to Teukolsky equations (for any half-integer spin) that do not rely on energy-Morawetz estimates, one can easily adapt the methodology in that paper to derive a weak Morawetz estimate, though with a loss of several derivatives, see Section 12 in \cite{MaSz26}.} and Shlapentokh-Rothman-Teixeira da Costa \cite{SRTdC20, SRTdC23}.

%%%%%%%%%%%%%%%%%%%%%%%%%%%%%%%%%%%%%%%%%

\subsubsection{Teukolsky equations in perturbations of Kerr with $|a|\ll m$}
\lab{subsubsect:literature:TeukolskyinpertKerr:intro}

%%%%%%%%%%%%%%%%%%%%%%%%%%%%%%%%%%%%%%%%%

To address the nonlinear stability of Kerr, it is important to extend the energy-Morawetz estimates for Teukolsky equations in Kerr reviewed in Section \ref{subsubsect:literature:TeukolskyinKerr:intro} to small perturbations of Kerr. This has been achieved in the context of the recent proofs of the nonlinear stability of Schwarzschild and of Kerr spacetimes for $|a|\ll m$: see Chapter 10 of \cite{KS} in the context of the nonlinear stability of Schwarzschild under polarized axisymmetry, Chapters 12 and 13 of \cite{DHRT} in the context of the nonlinear stability of Schwarzschild spacetimes for a codimension-3 set of initial data, and Chapter 9 of \cite{GKS22} in the context of the nonlinear stability of slowly rotating Kerr, i.e., with $|a|\ll m$. The proof of energy-Morawetz estimates for Teukolsky equations on perturbations of Kerr for $|a|<m$ has been recently obtained by Ma and the author \cite{MaSz26}, building in particular on \cite{MaSz24}, and plays a crucial role in  the present paper.

%%%%%%%%%%%%%%%%%%%%%%%

\subsection{Main results}

%%%%%%%%%%%%%%%%%%%%%%%

This paper proves two results used in the companion paper \cite{Sze}. The first concerns the derivation of a Teukolsky wave-transport system in perturbations of Kerr, and the second one concerns an extension of the energy-Morawetz estimates in perturbations of Kerr for Teukolsky equations in \cite{MaSz26} to the setting of \cite{Sze}.

%%%%%%%%%%%%%%%%%%%%%%%%%%%%%%%%%%%%%%%%%%%%%%

\subsubsection{Derivation of the Teukolsky wave-transport system in perturbations of Kerr}

%%%%%%%%%%%%%%%%%%%%%%%%%%%%%%%%%%%%%%%%%%%%%%%

In Section \ref{sec:sectionontheteukolskywavetransportsystem}, we derive a Teukosky wave-transport system in perturbations of Kerr, see Section \ref{sec:definitionofpmbphisspinpertubrtionsofKerr} for the definition of the Teukolsky quantities $\pmb\phi_s^{(p)}$, $s=\pm 2$, $p=0,1,2$, Theorem \ref{thm:derivationoftheTeukolskytensorialwavesystemfors=plusminus2:kerrpert:alternateformnullframeinsteadcoordvectorfield} for the wave part of the Teukolsky wave-transport system, and Proposition \ref{prop:transportequationphis=plusminus2p=0and1} for the transport part of the wave-transport system. 

The proof of Theorem \ref{thm:derivationoftheTeukolskytensorialwavesystemfors=plusminus2:kerrpert:alternateformnullframeinsteadcoordvectorfield} and Proposition \ref{prop:transportequationphis=plusminus2p=0and1} relies on the properties of the formalism for non-integrable structures introduced in \cite{GKS20} \cite{GKS22} and reviewed in Section \ref{sec:nonintergrableformalism}, and in particular on the tensorial wave equations for\footnote{With the notations of \cite{GKS22}, we have $\pmb\phi_{+2}^{(2)}=\qf$ and $\pmb\phi_{-2}^{(2)}=\qfb$, see Remarks \ref{rmk:compasisionphiplus2p=0withMaSz26andphiplus2p=2withqfGKS22} and \ref{rmk:compasisionphiminus2p=0withMaSz26andphiminus2p=2withqfGKS22}.} $\pmb\phi_{s}^{(2)}$, $s=\pm 2$, derived in Section 5 of \cite{GKS22}. 

\begin{remark}
The wave part  of the Teukolsky wave-transport system in Theorem \ref{thm:derivationoftheTeukolskytensorialwavesystemfors=plusminus2:kerrpert:alternateformnullframeinsteadcoordvectorfield} is the analog of the one in \eqref{eq:TensorialTeuSysandlinearterms:rescaleRHScontaine2:general:Kerrperturbation:intro}, where coordinates derivatives are replaced by their analog in the null frame $(e_3, e_4, e_1, e_2)$. The form of the wave part  of the Teukolsky wave-transport system using coordinates derivatives in then easily derived from Theorem \ref{thm:derivationoftheTeukolskytensorialwavesystemfors=plusminus2:kerrpert:alternateformnullframeinsteadcoordvectorfield} in Corollary 3.20 in \cite{Sze}. 
\end{remark}

%%%%%%%%%%%%%%%%%%%%%%%%%%%%%%%%%%%%%%%%%%%%%%%%%

\subsubsection{Energy-Morawetz estimates for Teukolsky equations in the setting of \cite{Sze}}

%%%%%%%%%%%%%%%%%%%%%%%%%%%%%%%%%%%%%%%%%%%%%%%%%

Energy-Morawetz estimates for Teukolsky equations in perturbations of Kerr have been established in \cite{MaSz26} for spacetimes extending all the way to $\II_+$. For the companion paper \cite{Sze}, we need to extend these results to the case of a spacetime\footnote{$\MM$ is the bootstrap spacetime used in the proof of Theorem \ref{MainThm-firstversion} in \cite{Sze}.} $\MM$ which is a perturbation of Kerr ending at a spacelike hypersurface $\Si_*$, see section \ref{section:SpacetimeMM-chap6} for the properties of $\Si_*$ and Theorem \ref{thm:main:MaSz26} for the main energy-Morawetz estimates for Teukolsky equations in perturbations of Kerr proved in this paper. This requires to revisit the proofs in \cite{MaSz26} by working in an expanding space-time $\MMh$ extending all the way to $\II_+$, and then to deal with new terms coming from a cut-off procedure in a causal spacetime region $\RR_*$ close to $\Si_*$, see Section \ref{sec:derfintionofRRstarandMMh} for the definitions of $\MMh$ and $\Si_*$.

%%%%%%%%%%%%%%%%%%%%%%%%%%%%%%%%%%%%

\subsection{Organisation of the paper}

%%%%%%%%%%%%%%%%%%%%%%%%%%%%%%%%%%%%

In Section \ref{sec:nonintergrableformalism}, we review the non-integrable formalism introduced in \cite{GKS20} \cite{GKS22}. Next, Section \ref{sec:sectionontheteukolskywavetransportsystem} is dedicated to the derivation of the Teukolsky wave system in perturbations of Kerr. Then, Section \ref{sec:energyMorawetzesitmatesforTeukoslkyonMM:upto15derivatives} adapts the energy-Morawetz estimates for Teukolsky equations in \cite{MaSz26} to our setting.

%%%%%%%%%%%%%%%%%%%%%%%%%%%%%%%%%%%%

\subsection{Acknowledgements}

%%%%%%%%%%%%%%%%%%%%%%%%%%%%%%%%%%%%

The author would like to emphasize that this work particularly owes to the contributions of his collaborators Elena Giorgi, Sergiu Klainerman and Siyuan Ma, as well as the ones of his former student Dawei Shen. 
In particular, he wishes to express his deepest gratitude to Sergiu Klainerman without whom this program would have been unthinkable, and he would also like to acknowledge the essential contribution of Siyuan Ma in  managing the last (but not least) hurdles. The author is supported by the ERC grant ERC-2023 AdG 101141855 BlaHSt.

%%%%%%%%%%%%%%%%%%%%

\section{Non-integrable formalism}
\lab{sec:nonintergrableformalism}

%%%%%%%%%%%%%%%%%%%%%

In this section, we review the formalism for non-integrable structures introduced in \cite{GKS20} \cite{GKS22}. This will be used throughout the paper.

%%%%%%%%%%%%%%%%%%%%%%%%%%%%

\subsection{Definition of the non-integrable formalism}

%%%%%%%%%%%%%%%%%%%%%%%%%%%%

%%%%%%%%%%%%%%%%%%%%%%%%%%

\subsubsection{Null pairs and horizontal structures}
\lab{subsection:review-horiz.structures}

%%%%%%%%%%%%%%%%%%%%%%%%%%

Consider a fixed null pair $e_3, e_4$, i.e., $\g(e_3, e_3)=\g(e_4, e_4)=0$,   $\g(e_3, e_4)=-2,$ and
 denote  by  $\O(\MM)$ the vectorspace  of horizontal vectorfields $X$  on $\MM$, i.e., $\g(e_3, X)= \g(e_4, X)=0$.
  Given a fixed   orientation  on $\MM$,  with corresponding  volume form  $\in$,  we define  the induced 
 volume form on   $\O(\MM)$ by,  
 \beaa
  \in(X, Y):=\frac 1 2\in(X, Y, e_3, e_4). 
  \eeaa
 A null  frame on $\MM$ consists of a choice of horizontal vectorfields  $e_1, e_2$, such that\footnote{We use greek 
 indices $\a, \b, \ga$ for $1,2,3,4$ and latin indices $a,b$ for $1,2$.}
 \beaa
 \g(e_a, e_b)=\de_{ab}\qquad  a, b=1,2.
 \eeaa  
 The commutator $[X,Y]$ of two horizontal vectorfields
may fail however to be horizontal. We say that the pair $(e_3, e_4 )$ is integrable if   $\O(\MM)$  forms an integrable distribution, i.e., $X, Y\in\O(\MM) $ implies that $[X,Y]\in\O(\MM)$. As is well-known,  the  principal null pair in Kerr fails to be integrable.
Given an arbitrary vectorfield $X$, we denote by $^{(h)}X$
its  horizontal projection, 
\beaa
{}^{(h)}X := X+ \frac 1 2 \g(X,e_3)e_4+ \frac 1 2   \g(X,e_4) e_3. 
\eeaa
A  $k$-covariant tensor-field $U$ is said to be horizontal,  $U\in \O_k(\MM)$,
if  for any $X_1,\ldots, X_k$ we have 
\beaa
U(X_1,\ldots, X_k)=U( ^{(h)} X_1,\ldots, {}^{(h)}X_k).
\eeaa

\begin{definition}\label{definition-SS-real}
We denote by $\sk_0=\sk_0(\MM, \mathbb{R})$ the set of pairs of real scalar functions on $\MM$,  by $\sk_1=\sk_1(\MM, \mathbb{R})$ the  set of real horizontal $1$-forms  on $\MM$   and by $\sk_2=\sk_2(\MM, \mathbb{R})$ the set of symmetric traceless   horizontal real $2$-tensors on $\MM$.
\end{definition}

\begin{definition}\label{definition-hodge-duals}
We define the dual of $\xi\in\sk_1$ and $U\in\sk_2$ by
\beaa
\dual \xi_{a}:=\in_{ab}\xi_b,\qquad \dual U_{ab}:=\in_{ac} U_{cb}.
\eeaa
\end{definition}

Note that given $\xi, \eta\in\sk_1$ and $U\in\sk_2$, we have 
\beaa
\dual(\dual \xi)=-\xi, \qquad \dual (\dual U)=-U,\qquad \dual\xi \c  \eta=-\xi\c\dual\eta.
\eeaa
 Also, given  $\xi, \eta\in\sk_1 $,  $U, V\in \sk_2$  we denote
\beaa
\xi\c \eta&:=&\de^{ab} \xi_a\eta_b,\qquad 
\xi\wedge\eta:=\in^{ab} \xi_a\eta_b=\xi\c\dual \eta,\qquad 
(\xi\hot \eta)_{ab}:=\xi_a \eta_b +\xi_b \eta_a-\de_{ab} \xi\c \eta,\\
(\xi\c U)_a&:=&\de^{bc} \xi_b U_{ac}, \qquad  (U\wedge V)_{ab} := \in^{ab}U_{ac}V_{cb}.
\eeaa

For any $ X, Y\in \O(\MM)$ we define  the induced metric $g(X, Y):=\g(X, Y)$ and the null second fundamental forms
\bea
\chi(X,Y):=\g(\D_Xe_4 ,Y), \qquad \chib(X,Y):=\g(\D_Xe_3,Y).
\eea
Observe that  $\chi$ and $\chib$  are  symmetric if and only if   the horizontal structure is integrable. Indeed this  follows easily from the following formulas
 \beaa
 \chi(X,Y)-\chi(Y,X)&=&\g(\D_X e_4, Y)-\g(\D_Ye_4,X)=-\g(e_4, [X,Y]),\\
 \chib(X,Y)-\chib(Y,X)&=&\g(\D_X e_3, Y)-\g(\D_Ye_3,X)=-\g(e_3, [X,Y]).
\eeaa
  Note  that  we  can view  $\chi$ and $\chib$ as horizontal 2-covariant tensor-fields
 by extending their definition to arbitrary vectorfields  $X, Y$  by setting  $\chi(X, Y)= \chi( ^{(h)}X, ^{(h)}Y)$,  $\chib(X, Y)= \chib( ^{(h)}X, ^{(h)}Y)$.
 Given an horizontal 2-tensor $U$  we define its trace $\tr U$  and anti-trace $\atr U$
\beaa
\tr (U):=\de^{ab}U_{ab}, \qquad \atr U:=\in^{ab} U_{ab}.
\eeaa
Accordingly we  decompose $\chi, \chib$ as follows,
\beaa
\chi_{ab}=\chih_{ab} +\frac 1 2 \de_{ab} \trch+\frac 1 2 \in_{ab}\atrch,\qquad \chib_{ab}=\chibh_{ab} +\frac 1 2 \de_{ab} \trchb+\frac 1 2 \in_{ab}\atrchb,
\eeaa
where $\chih$ and $\chibh$ denote respectively the symmetric traceless part of $\chi$ and $\chib$.

We define the horizontal covariant operator $\nab$ as follows. Given $X, Y\in \O(\MM)$
 \bea
 \nab_X Y&:=&^{(h)}(\D_XY)=\D_XY- \frac 1 2 \chib(X,Y)e_4 -  \frac 1 2 \chi(X,Y) e_3.
 \eea
 In particular, for  all  $X,Y, Z\in \O(\MM)$,
 \beaa
 Z g (X,Y)=g(\nab_Z X, Y)+ g(X, \nab_ZY).
 \eeaa

In the integrable case, $\nab$ coincides with the Levi-Civita connection
 of the metric induced on the integral surfaces of   $\O(\MM)$.  
 Given $X$ horizontal, $\D_4X$ and $\D_3 X$ are in general not horizontal.
 We define $\nab_4 X$ and $\nab_3 X$  to be the horizontal projections
 of the former.  More precisely,
 \beaa
 \nab_4 X&:=&^{(h)}(\D_4 X)=\D_4 X- \frac 1 2 \g(X, \D_4 e_3 ) e_4- \frac 1 2  \g(X, \D_4 e_4)  e_3 ,\\
 \nab_3 X&:=&^{(h)}(\D_3 X)=\D_3 X-   \frac 1 2 \g(X, \D_3e_3) e_3 - \frac 1 2   \g(X, \D_3 e_4 ) e_3. 
 \eeaa
The definition can be easily extended to arbitrary  $  \O_k(\MM) $ tensor-fields  $U$ 
\beaa
 \nab_4U(X_1,\ldots, X_k)&:=&e_4 (U(X_1,\ldots, X_k))- \sum_i U( X_1,\ldots, \nab_4 X_i, \ldots, X_k),\\
  \nab_3 U(X_1,\ldots, X_k)&:=&e_3 (U(X_1,\ldots, X_k)) -\sum_i U( X_1,\ldots, \nab_3 X_i, \ldots, X_k).
 \eeaa

%%%%%%%%%%%%%%%%%%%%%%%%%

\subsubsection{Ricci and curvature  coefficients}

%%%%%%%%%%%%%%%%%%%%%%%%%

Given a null frame $(e_1, e_2, e_3, e_4)$ we define the following connection coefficients,
 \bea
 \begin{split}
\chib_{ab}&:=\g(\D_ae_3, e_b),\qquad \,\,\,\,\,\,\,\chi_{ab}:=\g(\D_ae_4, e_b),\\
\xib_a&:=\frac 1 2 \g(\D_3 e_3 , e_a),\qquad\,\,\,\,\, \xi_a:=\frac 1 2 \g(\D_4 e_4, e_a),\\
\omb&:=\frac 1 4 \g(\D_3e_3 , e_4),\qquad\,\,\,\,\,\,\, \om:=\frac 1 4 \g(\D_4 e_4, e_3),\qquad \\
\etab_a&:=\frac 1 2\g(\D_4 e_3, e_a),\qquad \quad \eta_a:=\frac 1 2 \g(\D_3 e_4, e_a),\qquad\\
 \ze_a&:=\frac 1 2 \g(\D_{e_a}e_4,  e_3),
 \end{split}
\eea
which account for all the  connection coefficients except $\g(\D_{e_\mu} e_b, e_a)$, $\mu=1,2,3,4$, $a, b=1,2$.

We have the Ricci formulas 
\bea
\lab{eq:Ricciformula}
\D_a e_b&=&\nab_a e_b+\frac 1 2 \chi_{ab} e_3+\frac 1 2  \chib_{ab}e_4,\nn\\
\D_a e_4&=&\chi_{ab}e_b -\ze_a e_4,\nn\\
\D_a e_3&=&\chib_{ab} e_b +\ze_ae_3,\nn\\
\D_3 e_a&=&\nab_3 e_a +\eta_a e_3+\xib_a e_4,\nn\\
\D_3 e_3&=& -2\omb e_3+ 2 \xib_b e_b,\label{ricci}\\
\D_3 e_4&=&2\omb e_4+2\eta_b e_b,\nn\\
\D_4 e_a&=&\nab_4 e_a +\etab_a e_4 +\xi_a e_3,\nn\\
\D_4 e_4&=&-2 \om e_4 +2\xi_b e_b,\nn\\
\D_4 e_3&=&2 \om e_3+2\etab_b e_b.\nn
\eea 

For a given horizontal   1-form $\xi$, we  define the frame independent   operators
\bea\lab{eq:defintiondivcurlandnabhot}
\div\xi:=\de^{ab}\nab_b\xi_a,\qquad 
\curl\xi:=\in^{ab}\nab_a\xi_b,\qquad 
(\nab\hot \xi)_{ba}:=\nab_b\xi_a+\nab_a  \xi_b-\de_{ab}( \div \xi).
\eea
We also define the curvature components 
\bea
\a_{ab}:={\bf R}_{a4b4},\quad \b_a:=\frac 12 {\bf R}_{a434}, \quad \rho:=\frac 1 4 {\bf R}_{3434}, \quad\rhod:=\frac 1 4 \dual{\bf R}_{3434},\quad \bb_a:=\frac 1 2 {\bf R}_{a334}, \quad \aa_{ab}:={\bf R}_{a3b3},
\eea
where $\dual{\bf R}$ denotes the Hodge dual of the curvature tensor ${\bf R}$.

%%%%%%%%%%%%%%%%%%%%%%%%%%%%%%%%%%%%

\subsubsection{The tensorial wave operator}

%%%%%%%%%%%%%%%%%%%%%%%%%%%%%%%%%%%%

In order to define the tensorial wave operator in a covariant way, we first introduce the covariant derivative $\Ddot$ acting on mixed tensors of the type $\T_k (\MM)\otimes   \O_l (\MM)$, i.e., tensors  of the form  $U_{\nu_1\ldots \nu_k,  a_1\ldots a_l}$, 
for which we define
\beaa
\Ddot_\mu U_{\nu_1\ldots \nu_k,  a_1\ldots a_l}&:=& e_\mu(U_{\nu_1\ldots \nu_k,  a_1\ldots a_l}) -U_{\D_\mu e_{\nu_1}\ldots \nu_k,  a_1\ldots a_l}-\ldots- U_{\nu_1\ldots  \D_\mu e_{\nu_k},  a_1\ldots a_l}\\
&-& U_{\nu_1\ldots \nu_k,   ^{(h)}(\D_\mu e_{a_1})\ldots a_l}-  U_{\nu_1\ldots \nu_k,   a_1 \ldots ^{(h)}(\D_\mu e_{a_l})}.
\eeaa

\begin{proposition}
\lab{Proposition:commutehorizderivatives}
For a tensor $\Psi\in \O_1 (\MM)$, we   have  the following formula
 \bea
( \Ddot _\mu\Ddot_\nu -\Ddot_\nu\Ddot _\mu)\Psi_a=\Rdot_{a b  \mu\nu}\Psi^b
 \eea
with an immediate generalization to tensors $\Psi\in \O_l (\MM)$, where, with $(\La_\a)_{\b\ga}= \g(\D_\a e_\ga, e_\b)$,
 \bea
 \lab{eq:DefineRdot}
 \bsplit
 \Rdot_{ab   \mu\nu}&:= {\bf R}_{ab    \mu\nu}+ \frac 1 2  \B_{ab   \mu\nu},\\
  \B_{ab   \mu\nu} &:=  (\La_\mu)_{3a} (\La_\nu)_{b4}+  (\La_\mu)_{4a} (\La_\nu)_{b3}- (\La_\nu)_{3a} (\La_\mu)_{b4}-  (\La_\nu)_{4a} (\La_\mu)_{b3}.
  \end{split}
 \eea 
 \end{proposition}

\begin{proof}
See Proposition 2.1.27 in \cite{GKS22}.
\end{proof}

\begin{proposition}
\lab{proposition:componentsofB}
The components of $\B$   are given   by the following formulas:
\bea
\begin{split}
\B_{ a   b  c 3}&=     -  \trchb  \big( \de_{ca}\eta_b-  \de_{cb} \eta_a\big)  -  \atrchb \big( \in_{ca}  \eta_b -  \in_{cb}  \eta_a\big) \\
&+ 2 \big(- \chibh_{ca}  \eta_b + \chibh_{cb} \eta_a-  \chi_{ca} \xib_b+  \chi_{cb} \xib_a\big),\\
\B_{ a   b  c 4}&=     -  \trch  \big( \de_{ca}\etab_b-  \de_{cb} \etab_a\big)  -  \atrch \big( \in_{ca}  \etab_b -  \in_{cb}  \etab_a\big) \\
&+ 2 \big(- \chih_{ca}  \etab_b + \chih_{cb} \etab_a-  \chib_{ca} \xi_b+  \chib_{cb} \xi_a\big),\\
\B_{ a   b  3 4} &=4\big(-\eta_a \etab_b+\etab_a\eta_b -\xib_a \xi_b+\xi_a \xib_b\big),\\
\B_{abcd} &= \left(- \frac 12  \trch \trchb-\frac 1 2 \atrch \atrchb+\chih \c \chibh\right)\in_{ab}\in_{cd}.
\end{split}
\eea
\end{proposition}

\begin{proof}
See Proposition 2.2.4 in \cite{GKS22}.
\end{proof}

Then, we  define  the wave operator for $\psi \in \sk_k(\mathbb{C})$, $k=0,1,2$, to be, see Definition 2.3.1 in \cite{GKS22},  
 \bea\label{eq:def=squared-2}
 \squared_k\psi:= \g^{\mu\nu} \Ddot_\mu\Ddot_ \nu \psi.
\eea

%%%%%%%%%%%%%%%%%%%%%%%%%

\subsubsection{Commutation formulas}
\lab{sec:generalcommutationformulasrealcase}

%%%%%%%%%%%%%%%%%%%%%%%%%

We start with the following general commutation formulas.
\begin{lemma}
   \lab{LEMMA:COMM-GEN-B}
Let $U_{A}= U_{a_1\ldots a_k} $ be a general $k$-horizontal  tensorfield, and let $\B_{ab\mu\nu}$ be given by \eqref{eq:DefineRdot}. Then:
\begin{enumerate}
\item  We have
\bea
\,[\nab_3, \nab_b] U_A =- \chib_{bc} \nab_c U_A+( \eta_b-\ze_b) \nab_3 U_A +\xib_b \nab_4 U_A +\sum_{i=1}^k\Big(-\in_{a_i c} \dual\bb_b +  \frac 1 2 \B_{a_i c 3b} \Big) U_{a_1\ldots }\,^ c \,_{\ldots a_k}. 
\eea

\item We have
\bea
\,[\nab_4, \nab_b] U_A =- \chi_{bc} \nab_c U_a+( \etab_b+\ze_b) \nab_4 U_a +\xi_b \nab_3 U_a +\sum_{i=1}^k\Big(\in_{a_i c} \dual\b_b +  \frac 1 2 \B_{a_i c 4b} \Big) U_{a_1\ldots }\,^ c \,_{\ldots a_k}.
\eea

\item We have
\bea
\,[\nab_4, \nab_3] U_A = 2(\etab_b-\eta_b ) \nab_b U_A + 2 \om \nab_3 U_A -2\omb \nab_4 U_A+ \sum_{i=1}^k\Big( - \in_{a_i b}\dual \rho+  \frac 1 2 \B_{a_i b 43} \Big) U_{a_1\ldots}\,^b\,_{\ldots a_k}. 
\eea
\end{enumerate}
\end{lemma}

\begin{proof}
See the proof of Lemma 2.2.7 in \cite{GKS22}.
\end{proof}

Using the values of $\B_{ab\mu\nu}$ given by Proposition \ref{proposition:componentsofB}, we specialize  to the case of $\sk_0$, $\sk_1$ and $\sk_2$.

 \begin{lemma}
   \lab{lemma:comm}
   The following commutation formulas hold true:
   \begin{enumerate}
\item Given   $f \in \sk_0$, we have
       \bea\label{eq:comm-nab3-nab4-naba-f-general}
       \begin{split}
        \,[\nab_3, \nab_a] f &=-\frac 1 2 \left(\trchb \nab_a f+\atrchb \dual \nab_a f\right)+(\eta_a-\ze_a) \nab_3 f-\chibh_{ab}\nab_b f  +\xib_a \nab_4 f,\\
         \,[\nab_4, \nab_a] f &=-\frac 1 2 \left(\trch \nab_a f+\atrch \dual \nab_a f\right)+(\etab_a+\ze_a) \nab_4 f-\chih_{ab}\nab_b f  +\xi_a \nab_3 f, \\
         \, [\nab_4, \nab_3] f&= 2(\etab-\eta ) \c \nab f + 2 \om \nab_3 f -2\omb \nab_4 f. 
         \end{split}
       \eea

  \item   Given  $u\in \sk_1$, we have
    \bea\label{commutator-3-a-u-b}\label{commutator-u-in-SS1}
         \bsplit            
\,  [\nab_3,\nab_a] u_b    &=-\frac 1 2 \trchb \big( \nab_a u_b+\eta_b u_a-\de_{ab} \eta \c u \big) -\frac 1 2 \atrchb \big( \dual \nab_a u_b+\eta_b \dual u_a-\in_{ab} \eta\c u\big) \\
&+(\eta-\ze)_a \nab_3 u_b+\err_{3ab}[u],\\
  \err_{3ab}[u] &=-\dual \bb_a\dual u_b+\xib_a\nab_4 u_b-\xib_b \chi_{ac} u_c+\chi_{ab} \,\xib\c u-\chibh_{ac}\nab_c u_b-\eta_b\chibh_{ac}u_c+\chibh_{ab}\eta\c u,
   \end{split}
   \eea
   \bea\label{commutator-4-a-u-b}
   \bsplit
\,  [\nab_4,\nab_a] u_b    &=-\frac 1 2 \trch \big( \nab_a u_b+\etab_b u_a-\de_{ab} \etab \c u \big) -\frac 1 2 \atrch \big( \dual \nab_a u_b+\etab_b \dual u_a-\in_{ab} \etab\c u\big)\\
& +(\etab +\ze)_a \nab_4 u_b +\err_{4ab}[u],\\
   \err_{4ab}[u]&=\dual \b_a\dual u_b+\xi_a\nab_3 u_b-\xi_b \chib_{ac} u_c+\chib_{ab} \,\xi\c u-\chih_{ac}\nab_c u_b-\etab_b\chih_{ac}u_c+\chih_{ab}\etab\c u, 
      \end{split}
   \eea
   \bea
   \bsplit
 \, [\nab_4, \nab_3] u_a&=2 \om \nab_3 u_a -2\omb \nab_4 u_a+ 2(\etab_b-\eta_b ) \nab_b u_a +2(\etab \c u ) \eta_{a} -2 (\eta \c u )\etab_{a}\\
 &  -2 \dual \rho \dual u_a +\err_{43a}[u],\\
 \err_{43a}[u]&= 2 \big( \xib_{a}  \xi_b- \xi_{a}  \xib_b )u^b.
\end{split}
\eea

\item  Given  $u\in \sk_2$, we have 
    \bea\label{commutator-u-in-SS2}\label{commutator-3-a-u-bc}
         \bsplit            
\,  [\nab_3,\nab_a] u_{bc}    &=-\frac 1 2 \trchb\, (\nab_a u_{bc}+\eta_bu_{ac}+\eta_c u_{ab}-\de_{a b}(\eta \c u)_c-\de_{a c}(\eta \c u)_b )\\
&-\frac 1 2 \atrchb\, (\dual \nab_a u_{bc} +\eta_b\dual u_{ac}+\eta_c\dual u_{ab}- \in_{a b}(\eta \c u)_c- \in_{a c}(\eta \c u)_b )\\
&+(\eta_a-\ze_a)\nab_3 u_{bc}+\err_{3abc}[u],\\
\err_{3abc}[u]&= -2\dual \bb_a \dual u_{bc}+\xib_a \nab_4 u_{bc} -\xib_b\chi_{ad}u_{dc} -\xib_c\chi_{ad}u_{bd}+\chi_{ab}\xib_d u_{dc} \\
&+\chi_{ac}\xib_d u_{bd}-\chibh_{ad} \nab_d u_{bc} -\eta_b\chibh_{ad}u_{dc} - \eta_c\chibh_{ad}u_{bd}+\chibh_{ab}\eta_du_{dc} +\chibh_{ac}\eta_du_{bd},
   \end{split}
   \eea
   \bea\label{commutator-4-a-u-bc}
   \bsplit
\,  [\nab_4,\nab_a] u_{bc}    &=-\frac 1 2 \trch\, (\nab_a u_{bc}+\etab_bu_{ac}+\etab_c u_{ab}-\de_{a b}(\etab \c u)_c-\de_{a c}(\etab \c u)_b )\\
&-\frac 1 2 \atrch\, (\dual \nab_a u_{bc} +\etab_b\dual u_{ac}+\etab_c\dual u_{ab}- \in_{a b}(\etab \c u)_c- \in_{a c}(\etab \c u)_b )\\
&+(\etab_a+\ze_a)\nab_4 u_{bc}+\err_{4abc}[u],\\
\err_{4abc}[u]&= 2\dual \b_a \dual u_{bc}+\xi_a \nab_3 u_{bc} -\xi_b\chib_{ad}u_{dc} -\xi_c\chib_{ad}u_{bd}+\chib_{ab}\xi_d u_{dc} +\chib_{ac}\xi_d u_{bd}\\
& -\chih_{ad} \nab_d u_{bc} -\etab_b\chih_{ad}u_{dc} - \etab_c\chih_{ad}u_{bd}+\chih_{ab}\etab_du_{dc} +\chih_{ac}\etab_du_{bd}, 
     \end{split}
   \eea
   \bea\label{commutator-4-3-u-bc}
   \bsplit
   \, [\nab_4, \nab_3] u_{ab} &=2 \om \nab_3 u_{ab} -2\omb \nab_4 u_{ab} + 2(\etab_c-\eta_c ) \nab_c u_{ab} + 4 \eta \hot (\etab \c u)  \\
   &-4 \etab \hot (\eta \c u)-4 \dual \rho \dual u_{ab}+\err_{43ab}[u],\\
\err_{43ab}[u]&= 2 \big( \xib_{a}  \xi_c- \xi_{a}  \xib_c )u^c\,_{b}+2 \big( \xib_{b}  \xi_c- \xi_{b}  \xib_c )u_{a} \,^c.
   \end{split}
\eea
       \end{enumerate}
 \end{lemma}

\begin{proof}
See the proof of Lemma 2.2.8 in \cite{GKS22}.
\end{proof}

%%%%%%%%%%%%%%%%%%%%%%%%%%%%%

\subsubsection{Horizontal Lie derivatives}
\label{subsect:horizontalliederivatives}

%%%%%%%%%%%%%%%%%%%%%%%%%%%%%

Recall that the Lie derivative of a $k$-covariant tensor $U$ relative to a vectorfield  $X$ is given by
\beaa
\Lie_X{U}\big(e_{\a_1}, \ldots , e_{\a_k}\big) = X\big(U_{\a_1\ldots\a_k}\big) -  U\big(\Lie_Xe_{\a_1}, \ldots,e_{\a_k}\big) - U\big(e_{\a_1}, \ldots, \Lie_Xe_{\a_k}\big),
\eeaa
where $\Lie_X Y=[X, Y]$. We define horizontal Lie derivatives as follows, see Definition 2.2.12 in \cite{GKS22}.

\begin{definition}[Horizontal Lie derivatives]\label{definition:hor-Lie-derivative}
Given   vectorfields  $X$, $Y$,  the horizontal Lie  derivative  $\Lieb_XY$ is given by
 \beaa
 \Lieb_X Y :=\Lie_X Y+ \frac 1 2 \g(\Lie_XY, e_3) e_4+  \frac 1 2 \g(\Lie_XY, e_4) e_3.
 \eeaa
 Given  a horizontal covariant k-tensor $U$,  the horizontal  Lie derivative $\Lieb_X U $ is defined   to be the projection of $\Lie_X U$ to the  horizontal space, i.e.,
  \beaa
 \Lieb_X U\big(e_{a_1}, \ldots,e_{a_k}\big) := X\big(U_{a_1\ldots a_k}\big)- U\big(\Lieb_Xe_{a_1},\ldots, e_{a_k} \big)-\ldots -U\big(e_{a_1},\ldots,  \Lieb_Xe_{a_k} \big).
  \eeaa 

Also, given a mixed tensor $U$ of the type $\T_k (\MM)\otimes   \O_l (\MM)$, we define the general horizontal derivative $\Lied_XU$ as follows 
\beaa
&& \Lied_XU\big(e_{\a_1},\ldots, e_{\a_k},  e_{a_1},\ldots, e_{a_l}\big) \\
&:=& X\big(U_{\a_1\ldots \a_k,  a_1\ldots a_l}\big) -U\big(\Lie_Xe_{\a_1},\ldots, e_{\a_k},  e_{a_1}, \ldots, e_{a_l}\big) -\ldots  - U\big(e_{\a_1},\ldots.  \Lie_Xe_{\a_k},  e_{a_1},\ldots, e_{a_l}\big)\\
&& -U\big(e_{\a_1}, \ldots, e_{\a_k},  \Lieb_Xe_{a_1},\ldots, e_{a_l}\big) - \ldots -  U\big(e_{\a_1},\ldots, e_{\a_k},   e_{a_1}, \ldots, \Lieb_X e_{a_l}\big).
\eeaa
\end{definition}

%%%%%%%%%%%%%%%%%%%%%%%%%%%%%%%%%%%%%%%%%%%%    
    
\subsection{Main equations using complex notations}\label{section:complex-notations}

%%%%%%%%%%%%%%%%%%%%%%%%%%%%%%%%%%%%%%%%%%%%

%%%%%%%%%%%%%%%%%%%%%%%%

\subsubsection{Complex notations}
\lab{sec:complexnotationsRicciandcurvature}

%%%%%%%%%%%%%%%%%%%%%%%%

Recall Definition \ref{definition-SS-real} of the set of real horizontal tensors $\sk_k=\sk_k(\MM, \mathbb{R})$ on $\MM$ for $k=0,1,2$. We now define the corresponding complexified versions.

\begin{definition} 
\lab{def:skC:horizontaltensors}
We denote by $\sk_k(\mathbb{C})$, $k=0,1,2$, the following set of  horizontal tensors on $\MM$: 
\beaa
& a+ i b \in \sk_0(\mathbb{C})\,\,\,\textrm{if}\,\,\, (a, b) \in \sk_0, \qquad F= f+ i \dual f  \in \sk_1(\mathbb{C})\,\,\,\textrm{if}\,\,\, f \in \sk_1, \qquad U=u + i \dual u \in \sk_2(\mathbb{C})\,\,\,\textrm{if}\,\,\, u \in \sk_2,
\eeaa
where $F\in\sk_1(\mathbb{C})$ and $U\in\sk_2(\mathbb{C})$ are anti-self dual, i.e., $\dual F=-iF$ and $\dual U=-iU$. 
\end{definition}

\begin{definition}\lab{def:complexRicciandcurvaturecoefficients}
We define the following complexified curvature components 
\beaa
A:=\a+i\dual\a, \quad B:=\b+i\dual\b, \quad P:=\rho+i\dual\rho,\quad \Bb:=\bb+i\dual\bb, \quad \Ab:=\aa+i\dual\aa,
\eeaa      
with $A, \Ab\in\sk_2(\mathbb{C})$, $B, \Bb\in\sk_1(\mathbb{C})$, $P\in\sk_0(\mathbb{C})$, and the following complexified Ricci coefficients     
\beaa
&& X:=\chi+i\dual\chi, \quad \Xb:=\chib+i\dual\chib, \quad H:=\eta+i\dual \eta, \quad \Hb:=\etab+i\dual \etab,  \\ 
&& Z:=\ze+i\dual\ze, \quad \Xi:=\xi+i\dual\xi, \quad \Xib:=\xib+i\dual\xib,
\eeaa    
with $\widehat{X}, \Xbh\in\sk_2(\mathbb{C})$, $H, \Hb, Z, \Xi, \Xib\in\sk_1(\mathbb{C})$, where $\widehat{X}, \Xbh$, as well as $\tr X, \tr\Xb$ are given by
\beaa
\tr X := \trch-i\atrch, \quad \widehat{X}:=\chih+i\dual\chih, \quad \tr\Xb:=\trchb -i\atrchb, \quad \Xbh:=\chibh+i\dual\chibh.
\eeaa
\end{definition}

\begin{definition}
We define derivatives of complex quantities as follows
\begin{itemize}
\item For two scalar functions $a$ and $b$, we define
\beaa
\DD(a+ib) &:=& (\nabla+i\dual\nabla)(a+ib).
\eeaa

\item For a 1-form $f$, we define
\beaa
\ov{\DD}\c(f+i\dual f) &:=& (\nabla-i\dual\nabla)\c(f+i\dual f)
\eeaa
and  
\beaa
\DD\hot(f+i\dual f) &:=& (\nabla+i\dual\nabla)\hot(f+i\dual f).
\eeaa

\item For a symmetric traceless 2-tensor $u$, we define
\beaa
\ov{\DD}\c(u+i\dual u) &:=& (\nabla-i\dual\nabla)\c(u+i\dual u).
\eeaa
\end{itemize}
\end{definition}

We will use the following two lemmas involving complex notations.
\begin{lemma}\label{dot-hot-complex} 
Let  $E, F \in \sk_1(\mathbb{C})$ and $U\in \sk_2(\mathbb{C}) $.  Then
\bea\label{simil-Leibniz}
E \hot ( \ov{F} \c U)+F \hot ( \ov{E} \c U)&=&2 ( E \c \ov{F}+ \ov{E} \c F) \ U.
\eea
Also, for $E=e+i \dual e$, $F=f+i \dual f$
\bea\label{simil-Leib-half}
E \hot ( \ov{F} \c U) &=& 4 \big(e \c f - i e \wedge f \big) U.
\eea
\end{lemma}

\begin{proof}
See Lemma 2.4.5 in \cite{GKS22}.
\end{proof}

\begin{lemma}\label{SIMPLIFICATION-ANGULAR}
Let $h$ be a scalar function, $F \in \sk_1(\mathbb{C}) $,  $U\in \sk_2(\mathbb{C})$.  Then
\bea
\bsplit
 \ov{\DD} \c (h F)&= h \ov{\DD} \c F+  \ov{\DD}(h) \c F \lab{ov-HH-hF},\\
 \DD\hot (h F)&= h  \DD\hot  F+ \DD( h) \hot F  \lab{DD-hot-hF},\\
\ov{\DD}\c( h U)&=  \ov{\DD}(h) \c U  + h(\ov{\DD} \c U) \lab{ov-DD-hu},\\
 \DD\hot(\ov{F}\c U)&=  2   (\DD\c \ov{F} ) U + 2   (\ov{F}\c \DD) U\label{Leibniz-hot},\\
U \c \ov{\DD} F&= U (\ov{\DD} \c F).\label{rule-1}
\end{split}
\eea
Also,
\bea
\bsplit
F \hot (\DDb \c U)&= 2(F\c \DDb) U= 4F \c \nab U , \lab{Leib-eq-DDb} \\
(F \c \DDb) U + (\ov{F} \c \DD) U&= 4f \c \nab U= 2(F +\ov{F}) \c \nab U. \lab{Leib-eq-DDb-DD-nab}\lab{relation0angular=der}
\end{split}
\eea
\end{lemma}

\begin{proof}
See Lemma 2.4.6 in \cite{GKS22}.
\end{proof}

%%%%%%%%%%%%%%%%%%%%%%%%%%%%%

\subsubsection{Main equations  in complex form}

%%%%%%%%%%%%%%%%%%%%%%%%%%%%%

The complex notations allow us to rewrite the Ricci equations in  a more compact  form. 
\begin{proposition}
\label{prop-nullstr:complex} {We have}
\beaa
\nab_3\tr\Xb +\frac{1}{2}(\tr\Xb)^2+2\omb\,\tr\Xb &=& \DD\c\ov{\Xib}+\Xib\c\ov{\Hb}+\ov{\Xib}\c(H-2Z)-\frac{1}{2}\Xbh\c\ov{\Xbh},\\
\nab_3\Xbh+\Re(\tr\Xb) \Xbh+ 2\omb\,\Xbh&=&\frac 1 2  \DD\hot \Xib+\frac 1 2   \Xib\hot(H+\Hb-2Z)-\Ab,
\eeaa
\beaa
\nab_3\tr X +\frac{1}{2}\tr\Xb\tr X-2\omb\tr X &=& \DD\c\ov{H}+H\c\ov{H}+2P+\Xib\c\ov{\Xi}-\frac{1}{2}\Xbh\c\ov{\Xh},\\
\nab_3\widehat{X} +\frac{1}{2}\tr\Xb\, \widehat{X} -2\omb\widehat{X} &=&\frac 1 2  \DD\hot H  +\frac 1 2 H\hot H -\frac{1}{2}\ov{\tr X} \widehat{\Xb}+\frac{1}{4}\Xib\hot\Xi,
\eeaa
\beaa
\nab_4\tr\Xb +\frac{1}{2}\tr X\tr\Xb -2\om\tr\Xb &=& \DD\c\ov{\Hb}+\Hb\c\ov{\Hb}+2\ov{P}+\Xi\c\ov{\Xib}-\frac{1}{2}\Xh\c\ov{\Xbh},\\
\nab_4\widehat{\Xb} +\frac{1}{2}\tr X\, \widehat{\Xb} -2\om\widehat{\Xb} &=&\frac  12  \DD\hot\Hb  +\frac 1 2 \Hb\hot\Hb -\frac{1}{2}\ov{\tr\Xb} \widehat{X}+\frac{1}{4}\Xi\hot\Xib,
\eeaa
\beaa
\nab_4\tr X +\frac{1}{2}(\tr X)^2+2\om\tr X &=& \DD\c\ov{\Xi}+\Xi\c\ov{H}+\ov{\Xi}\c({\Hb}+2Z)-\frac{1}{2}\Xh\c\ov{\Xh},\\
\nab_4\Xh+\Re(\tr X)\Xh+ 2\om\Xh&=&\frac 1 2  \DD\hot \Xi+\frac 12   \Xi\hot(\Hb+H+2Z)-A.
\eeaa
Also,
\beaa
\nab_3Z +\frac{1}{2}\tr\Xb(Z+H)-2\omb(Z-H) &=& -2\DD\omb -\frac{1}{2}\widehat{\Xb}\c(\ov{Z}+\ov{H})+\frac{1}{2}\tr X\Xib+2\om\Xib -\Bb+\frac{1}{2}\ov{\Xib}\c\Xh,\\
\nab_4Z +\frac{1}{2}\tr X(Z-\Hb)-2\om(Z+\Hb) &=& 2\DD\om +\frac{1}{2}\widehat{X}\c(-\ov{Z}+\ov{\Hb})-\frac{1}{2}\tr\Xb\Xi-2\omb\Xi -B-\frac{1}{2}\ov{\Xi}\c\Xbh,\\
\nab_3\Hb -\nab_4\Xib &=&  -\frac{1}{2}\ov{\tr\Xb}(\Hb-H) -\frac{1}{2}\Xbh\c(\ov{\Hb}-\ov{H}) -4\om\Xib+\Bb,\\
\nab_4H -\nab_3\Xi &=&  -\frac{1}{2}\ov{\tr X}(H-\Hb) -\frac{1}{2}\Xh\c(\ov{H}-\ov{\Hb}) -4\omb\Xi-B,
\eeaa
and
\beaa
\nab_3\om+\nab_4\omb -4\om\omb -\xi\c \xib -(\eta-\etab)\c\ze +\eta\c\etab&=&   \rho.
\eeaa
Also,
\beaa
\frac{1}{2}\ov{\DD}\c\Xh +\frac{1}{2}\Xh\c\ov{Z} &=& \frac{1}{2}\DD\ov{\tr X}+\frac{1}{2}\ov{\tr X}Z-i\Im(\tr X)H-i\Im(\tr \Xb)\Xi-B,\\
\frac{1}{2}\ov{\DD}\c\Xbh -\frac{1}{2}\Xbh\c\ov{Z} &=& \frac{1}{2}\DD\ov{\tr\Xb}-\frac{1}{2}\ov{\tr\Xb}Z-i\Im(\tr\Xb)\Hb-i\Im(\tr X)\Xib+\Bb,
\eeaa
and,
\beaa
\curl\ze&=&-\frac 1 2 \chih\wedge\chibh   +\frac 1 4 \big(  \trch\atrchb-\trchb\atrch   \big)+\om \atrchb -\omb\atrch+\dual \rho.
\eeaa
\end{proposition}

The complex notations allow us to rewrite the Bianchi identities as follows.  
  \begin{proposition}\label{prop:bianchi:complex} 
    We have,
 \beaa
 \nab_3A -\frac 1 2 \DD\hot B &=& -\frac{1}{2}\tr\Xb A+4\omb A +\frac 1 2 (Z+4H)\hot B -3\ov{P}\Xh,\\
\nab_4B -\frac{1}{2}\ov{\DD}\c A &=& -2\ov{\tr X} B -2\om B +\frac{1}{2}A\c  (\ov{2Z +\Hb})+3\ov{P} \,\Xi,\\
\nab_3B-\DD\ov{P} &=& -\tr\Xb B+2\omb B+\ov{\Bb}\c \Xh+3\ov{P}H +\frac{1}{2}A\c\ov{\Xib},\\
\nab_4P -\frac{1}{2}\DD\c \ov{B} &=& -\frac{3}{2}\tr X P +\frac{1}{2}(2\Hb+Z)\c\ov{B} -\ov{\Xi}\c\Bb -\frac{1}{4}\Xbh\c \ov{A}, \\
\nab_3P +\frac{1}{2}\ov{\DD}\c\Bb &=& -\frac{3}{2}\ov{\tr\Xb} P -\frac{1}{2}(\ov{2H-Z})\c\Bb +\Xib\c \ov{B} -\frac{1}{4}\ov{\Xh}\c\Ab, \\
\nab_4\Bb+\DD P &=& -\tr X\Bb+2\om\Bb+\ov{B}\c \Xbh-3P\Hb -\frac{1}{2}\Ab\c\ov{\Xi},\\
\nab_3\Bb +\frac{1}{2}\ov{\DD}\c\Ab &=& -2\ov{\tr\Xb}\,\Bb -2\omb\,\Bb -\frac{1}{2}\Ab\c (\ov{-2Z +H})-3P \,\Xib,\\
\nab_4\Ab +\frac{1}{2}\DD\hot\Bb &=& -\frac{1}{2}\tr X \Ab+4\om\Ab +\frac{1}{2}(Z-4\Hb)\hot \Bb -3P\Xbh.
\eeaa
    \end{proposition} 
    
    \begin{proof} 
    See Proposition 2.4.12 in \cite{GKS22}.
    \end{proof}

%%%%%%%%%%%%%%%%%%%%%%%%%%%%%%%%%%%%%

\subsubsection{Main complex equations using conformal derivatives}

%%%%%%%%%%%%%%%%%%%%%%%%%%%%%%%%%%%%%

Following sections 2.2.9 and 2.4.3 in \cite{GKS22}, we introduce conformal derivatives. Consider  frame transformations of the form
\beaa
e_3'=\la^{-1} e_3, \qquad  e'_4 = \la e_4 , \qquad e_a'= e_a.
\eeaa
Note that  under   the   above  mentioned  frame transformation we have
\beaa
 \Xb'&=&\la^{-1} \Xb, \quad X'=\la X, \quad \Xi'= \la^2\Xi, \quad   H'=H, \quad \Hb'=\Hb,  \quad \Xib'=\la^{-2}\Xib,\\
   A'&=&\la^2A,\quad B'=\la B,    \quad P'=P,    \quad  \Bb'=\la^{-1} \Bb,\quad  \Ab'=\la^{-2} \Ab,
   \eeaa
   and
   \beaa
 \omb'&=& \la^{-1}\left(\omb +\frac{1}{2} e_3(\log \la)\right), \quad \om'= \la\left(\om -\frac{1}{2} e_4(\log \la)\right), \quad
 Z'= Z - \DD(\log \la).
\eeaa

\begin{definition}[$s$-conformally invariants]
\lab{def:sconformalinvariants}
We say that  a horizontal tensor $f$ is $s$-conformally invariant  if, under the  conformal  frame transformation above,  it changes as $f'=\la^s f $. 
\end{definition}

\begin{remark}
If $f$ $s$-conformal invariant, then  $\nab_3 f, \nab_4 f, \nab_a f$ are not conformal invariant.
\end{remark} 

 We correct the lack of being conformal invariant by making the following  definition.  

\begin{lemma}\label{lemma:definition-conformal-derivatives}
If $f$ is $s$-conformal invariant, then: 
 \begin{enumerate}
 \item $\nabc_3 f:= \nab_3f-2 s \omb f$ is $(s-1)$-conformally invariant.
 
 \item $\nabc_4 f:= \nab_4f+2 s \om f$ is $(s+1)$-conformally invariant.
 
 \item $\nabc_a f:= \nab_af+ s \ze_a f$ is $s$-conformally invariant. 
  \end{enumerate}
\end{lemma}

\begin{proof}
Immediate verification.
\end{proof}

\begin{remark} 
Note that $s$ is precisely what   in \cite{Ch-Kl} is called       the  signature of the tensor. In GHP formalism \cite{GHP}, the signature is related to the boost weights of the complex scalars.
\end{remark}

\begin{definition}
 We define the following conformal angular derivatives in the complex notation:
 \begin{itemize}
\item For $a+i b \in \sk_0(\mathbb{C}) $  we define
\beaa
\DDc(a+ib) &:=& \big(\nabc +i\dual \nabc\big)(a+ib).
\eeaa

\item For  $f+i \dual f \in\sk_1(\mathbb{C}) $ we define
\beaa
\DDc(f+i\dual f) &:=& \big(\nabc+i\dual \nabc\big) \c (f+i\dual f),
\\
\DDc \hot(f+i\dual f) &:=& (\nabc+i\dual\nabc )\hot(f+i\dual f).
\eeaa
\item For $u+ i \dual u \in \sk_2(\mathbb{C})$ we define
\beaa
\DDc \c(u+i\dual u) &:=& \big(\nabc +i\dual\nabc \big)\c(u+i\dual u).
\eeaa
\item In all the above cases we set
\beaa
\DDbc&:=&\nabc-i\nabc.
\eeaa
\end{itemize}
\end{definition}

Using these definitions we rewrite the main equations as follows.
\begin{proposition}
\label{prop-nullstr:complex-conf}
We have
\beaa
\nabc_3\tr\Xb +\frac{1}{2}(\tr\Xb)^2 &=& \DDc\c\ov{\Xib}+\Xib\c\ov{\Hb}+\ov{\Xib}\c H-\frac{1}{2}\Xbh\c\ov{\Xbh},\\
\nabc_3\Xbh+\Re(\tr\Xb) \Xbh&=&\frac 1 2  \DDc\hot \Xib+  \frac 1 2  \Xib\hot(H+\Hb)-\Ab,
\eeaa
\beaa
\nabc_3\tr X +\frac{1}{2}\tr\Xb\tr X &=& \DDc\c\ov{H}+H\c\ov{H}+2P+\Xib\c\ov{\Xi}-\frac{1}{2}\Xbh\c\ov{\Xh},\\
\nabc_3\widehat{X} +\frac{1}{2}\tr\Xb\, \widehat{X} &=&\frac 1 2 \DDc\hot H  +\frac 1 2 H\hot H -\frac{1}{2}\ov{\tr X} \widehat{\Xb}+\frac 1 4 \Xib\hot\Xi,
\eeaa
\beaa
\nabc_4\tr\Xb +\frac{1}{2}\tr X\tr\Xb &=& \DDc\c\ov{\Hb}+\Hb\c\ov{\Hb}+2\ov{P}+\Xi\c\ov{\Xib}-\frac{1}{2}\Xh\c\ov{\Xbh},\\
\nabc_4\widehat{\Xb} +\frac{1}{2}\tr X\, \widehat{\Xb} &=&\frac 1 2  \DDc\hot\Hb  +\frac 1 2 \Hb\hot\Hb -\frac{1}{2}\ov{\tr\Xb} \widehat{X}+\frac 1 4 \Xi\hot\Xib,
\eeaa
\beaa
\nabc_4\tr X +\frac{1}{2}(\tr X)^2 &=& \DDc\c\ov{\Xi}+\Xi\c\ov{H}+\ov{\Xi}\c\Hb-\frac{1}{2}\Xh\c\ov{\Xh},\\
\nabc_4\Xh+\Re(\tr X)\Xh&=&\frac 1 2  \DDc\hot \Xi+\frac 1 2   \Xi\hot(\Hb+H)-A,
\eeaa
\beaa
\nabc_3\Hb -\nabc_4\Xib &=&  -\frac{1}{2}\ov{\tr\Xb}(\Hb-H) -\frac{1}{2}\Xbh\c(\ov{\Hb}-\ov{H}) +\Bb,\\
\nabc_4H -\nabc_3\Xi &=&  -\frac{1}{2}\ov{\tr X}(H-\Hb) -\frac{1}{2}\Xh\c(\ov{H}-\ov{\Hb}) -B.
\eeaa
Also,
\beaa
\frac{1}{2}\ov{\DDc}\c\Xh &=& \frac{1}{2}\DDc\ov{\tr X}-i\Im(\tr X)H-i\Im(\tr \Xb)\Xi-B,\\
\frac{1}{2}\ov{\DDc}\c\Xbh &=& \frac{1}{2}\DDc\ov{\tr\Xb}-i\Im(\tr\Xb)\Hb-i\Im(\tr X)\Xib+\Bb.
\eeaa
\end{proposition}

\begin{proof}
See Proposition 2.4.14 in \cite{GKS22}.
\end{proof}
    
  \begin{proposition}\label{prop:bianchi:complex-conf} 
    We have
 \beaa
 \nabc_3A -\frac 1 2 \DDc\hot B &=& -\frac{1}{2}\tr\Xb A + 2 H   \hot B -3\ov{P}\Xh,\\
\nabc_4B -\frac{1}{2} \DDbc \c A &=& -2\ov{\tr X} B +\frac{1}{2}A\c \ov{\Hb}+3\ov{P} \,\Xi,\\
\nabc_3B-\DDc\ov{P} &=& -\tr\Xb B+\ov{\Bb}\c \Xh+3\ov{P}H +\frac{1}{2}A\c\ov{\Xib},\\
\nabc_4P -\frac{1}{2}\DDc\c \ov{B} &=& -\frac{3}{2}\tr X P + \Hb \c\ov{B} -\ov{\Xi}\c\Bb -\frac{1}{4}\Xbh\c \ov{A}, 
\eeaa
\beaa
\nabc_3P +\frac{1}{2}\DDbc \c\Bb &=& -\frac{3}{2}\ov{\tr\Xb} P - \ov{H} \c\Bb +\Xib\c \ov{B} -\frac{1}{4}\ov{\Xh}\c\Ab, \\
\nabc_4\Bb+\DDc P &=& -\tr X\Bb+\ov{B}\c \Xbh-3P\Hb -\frac{1}{2}\Ab\c\ov{\Xi},\\
\nabc_3\Bb +\frac{1}{2}\DDbc \c\Ab &=& -2\ov{\tr\Xb}\,\Bb  -\frac 1 2  \Ab\c \ov{H}-3P \,\Xib,\\
\nabc_4\Ab +\frac 1 2 \DDc\hot\Bb &=& -\frac{1}{2}\tr X \Ab - 2 \Hb\hot \Bb -3P\Xbh.
\eeaa
    \end{proposition} 
    
\begin{proof}
See Proposition 2.4.15 in \cite{GKS22}.
\end{proof}

%%%%%%%%%%%%%%%%%%%%%%%%%%%%%%%

\subsection{Kerr values}

%%%%%%%%%%%%%%%%%%%%%%%%%%%%%%%

%%%%%%%%%%%%%%%%%%%%%%%%%%%%%%%

\subsubsection{Normalized coordinates in Kerr spacetimes}
\label{subsect:normalizedcoords}

%%%%%%%%%%%%%%%%%%%%%%%%%%%%%%%

The Kerr metric in Boyer--Lindquist coordinates $(t,r,\th,\phi)$ is given by
\begin{align}
\gam={}& \g_{tt}dt^2 +\g_{rr}dr^2+(\g_{t\phi}+\g_{\phi t})dtd\phi +\g_{\phi\phi}d\phi^2 +\g_{\th\th}d\th^2,
\end{align}
where
\begin{equation}
\begin{split}
\g_{tt}={}&-\frac{\Delta-a^2\sin^2\theta}{|q|^2}, \quad \g_{t\phi}={}\g_{\phi t}=-\frac{2amr\sin\theta}{|q|^2}, \quad \g_{rr}={}\frac{|q|^2}{\Delta},\\
\g_{\phi\phi}={}&\frac{(r^2+a^2)^2-a^2\sin^2\theta\Delta}{|q|^2}\sin^2\theta, \quad \g_{\th\th}={}|q|^2,
\end{split}
\end{equation}
with 
\bea
\Delta:=r^2-2mr+a^2,\qquad |q|^2:=r^2+a^2\cos^2\th.
\eea
{In particular, $\partial_{t}$ and $\partial_{\phi}$ are Killing vectorfields and the larger root 
\begin{align}
r_+:=m + \sqrt{m^2 -a^2}
\end{align}
of $\Delta=\De(r)$ corresponds to the location of the event horizon.

It is well-known that the metric is singular on the event horizon in both the Boyer--Lindquist and the tortoise coordinates. To extend the Kerr metric beyond the future event horizon, we define the ingoing Eddington--Finkelstein coordinates $(v_+, r,\th,\phi_+)$ by 
\bea\lab{eq:definitionofingoingEFcoordiantesvplusandphiplus}
dv_+=dt+\frac{r^2+a^2}{\De}dr, \quad d\phi_+=d\phi+\frac{a}{\Delta}dr \,\,\, \text{mod } 2\pi.
\eea
The Kerr metric in this coordinate system is 
\begin{align}\lab{eq:KerrmetriciningoingEddigtonFinkelstein}
\gam={}&-\bigg(1-\frac{2mr}{|q|^2}\bigg)dv_+^2 +2dr dv_+ -\frac{4amr\sin^2\th}{|q|^2}dv_+d\phi_+ -2a\sin^2\th dr d\phi_+\nn\\
&+|q|^2 d\th^2+\frac{(r^2+a^2)^2-a^2\sin^2\theta\Delta}{|q|^2}\sin^2\th d\phi_+^2.
\end{align}

In the following lemma, we introduce coordinates systems, referred to as \textit{normalized coordinates}, used throughout the paper. 

\begin{lemma}[Normalized coordinates]
\label{lem:specificchoice:normalizedcoord}
We fix constants $\dhor$ and $\dbl$ such that 
$$0<\dhor\ll \dbl\ll 1-\frac{|a|}{m}.$$ 
There exists a choice of smooth functions $\tmod=\tmod(r)$ and $\phimod=\phimod(r)$ such that the coordinate systems $(\tt, r, x^1_0, x^2_0)$ and $(\tt, r, x^1_p, x^2_p)$, defined respectively on $\th\neq 0, \pi$ and $\th\neq\frac{\pi}{2}$, with 
\bea\lab{eq:definitionof:specificchoice:normalizedcoord}
\tau=v_+-\tmod, \quad \tphi=\phi_+ -  \phimod, \quad x^1_0=\th, \quad x^2_0=\tphi, \quad x^1_p=\sin\th\cos\tphi, \quad x^2_p=\sin\th\sin\tphi,
\eea
satisfy the following properties:
\begin{enumerate}
\item defining the causal spacetime region $\MM$ and corresponding spacelike boundary $\AA$ by 
\beaa
\bsplit
{}\qquad\MM&:=\big(\{(\tt, r, x^1_0, x^2_0),\,\, \th\neq 0, \pi\}\cup\{(\tt, r, x^1_p, x^2_p),\,\, \th\neq \pi/2\}\big)\cap\{r\geq r_+(1-\dhor)\}, \\ 
{}\qquad\AA&:=\pr\MM=\big(\{(\tt, r, x^1_0, x^2_0),\,\, \th\neq 0, \pi\}\cup\{(\tt, r, x^1_p, x^2_p),\,\, \th\neq \pi/2\}\big)\cap\{r=r_+(1-\dhor)\},
\end{split}
\eeaa
$\MM$ is covered by $(\tt, r, x^1_0, x^2_0)$ and $(\tt, r, x^1_p, x^2_p)$ with the metric components and  inverse metric components being smooth on their respective coordinate patch, 

\item $(\tt, r, x^1_0, x^2_0)$ coincides with Boyer-Lindquist coordinates\footnote{In particular, we have
\beaa
\tmod'(r)=\frac{r^2+a^2}{\De}, \qquad \phimod'(r)=\frac{a}{\De}\quad\textrm{on}\quad r\in [r_+(1+2\dbl), 12m].
\eeaa}
 in $r\in [r_+(1+2\dbl),  12m]$, 

\item for $r\notin (r_+(1+\dbl), 13m)$, we choose 
\beaa
\begin{split}
\tmod'(r) &=\frac{m^2}{r^2}, \qquad \phimod'(r)=0\quad\textrm{on}\quad r\leq r_+(1+\dbl),\\
\tmod'(r)&=\frac{2(r^2+a^2)}{\De} -\frac{m^2}{r^2}, \qquad \phimod'(r)=\frac{2a}{\De} \quad\textrm{on}\quad r\geq 13m,
\end{split}
\eeaa

\item the level sets of $\tt$ in $\MM$ are globally spacelike, transverse to the future event horizon $\HH_+$ and the spacelike boundary $\AA$, and asymptotically null to future null infinity $\II_+$.
\end{enumerate}

Furthermore, the nontrivial inverse metric components in the coordinate system $(\tt, r, \th, \tphi)$ are
\begin{align}
\label{eq:inverse:hypercoord}
\gam^{\tt \tt}={}&\frac{a^2\sin^2\th}{|q|^2} -\frac{2(r^2+a^2)}{|q|^2}\tmod'+\frac{\Delta}{|q|^2}(\tmod')^2, \qquad \gam^{rr}=\frac{\Delta}{|q|^2},\nn\\ 
\gam^{\tt r}={}&\gam^{r \tt}=\frac{r^2+a^2}{|q|^2}\left(1-\frac{\De}{r^2+a^2}\tmod'\right), \qquad \gam^{r\tphi} =\gam^{\tphi r} =\frac{a}{|q|^2}-\frac{\Delta}{|q|^2}\phimod',\nn\\
\gam^{\tt\tphi}={}& \gam^{\tphi \tt}= \frac{a}{|q|^2} (1-\tmod')-\phimod'\frac{r^2+a^2}{|q|^2}\left(1-\frac{\De}{r^2+a^2}\tmod'\right),\nn\\
\gam^{\th\th}={}&\frac{1}{|q|^2}, \qquad \gam^{\tphi\tphi}=\frac{1}{|q|^2\sin^2\th}-\frac{2a}{|q|^2}\phimod' +\frac{\Delta}{|q|^2}(\phimod')^2.
\end{align}
\end{lemma}

\begin{remark}\lab{rmk:phimodprimeisproportionaltoa!!}
Additionally, we may choose $\phimod$ such that 
\beaa
\phimod'(r)=a\phi_{\textrm{mod},0}'(r), \qquad \phi_{\textrm{mod},0}'(r)\geq 0\quad \forall r\in(r_+(1-\dhor), +\infty),
\eeaa
so that $\phimod'(r)$ has the same sign as $a$. From now on, we will assume that our choice of $\phimod$ satisfies this property. In view of \eqref{eq:inverse:hypercoord}, it implies that the inverse metric coefficients $\gam^{\a\b}$ in the normalized coordinates system $(\tau, r, x^1, x^2)$ are invariant under the change $(a, \tphi)\to (-a, -\tphi)$. 
\end{remark}

\begin{proof}
See the proof of Lemma 2.1 in \cite{MaSz24}.
\end{proof}

%%%%%%%%%%%%%%%%%%%%%%%%%%%%%%%

\subsubsection{Principal null pair in Kerr}
\label{subsect:principalnullpairinKerr}

%%%%%%%%%%%%%%%%%%%%%%%%%%%%%%%

We consider the principal null pair of Kerr which is regular across the future event horizon, i.e., in Boyer-Lindquist coordinates, 
\bea
\lab{def:e3e4inKerr}
 e_4 = \frac{r^2+a^2}{|q|^2} \pr_t +\frac{\De}{|q|^2} \pr_r +\frac{a}{|q|^2} \pr_{\phi}, \qquad 
 e_3=\frac{r^2+a^2}{\De} \pr_t -\pr_r +\frac{a}{\De} \pr_{\phi}.
\eea
Also, we consider its associated horizontal bundle $\{e_3, e_4\}^\perp$, which, for $\th\neq 0, \pi$, is spanned by 
\bea
\lab{def:e1e2inKerr}
 e_1=\frac{1}{|q|}\pr_\th,\quad e_2=\frac{a\sin\th}{|q|}\pr_t+\frac{1}{|q|\sin\th}\pr_\phi,
\eea
and we define the complex-valued scalar $q$ and the complex horizontal  $1$-forms $\Jk$ and $\Jk_\pm$ as 
\bea\lab{eq:def:Jkandq}
\bsplit
q=& r+ia\cos\th, \qquad \Jk=j+i\dual j, \qquad \Jk_\pm=j_\pm+i\dual j_\pm, \qquad j_1=0, \quad j_2=\frac{\sin\th}{|q|}, \\
(j_+)_1=&\frac{1}{|q|} \cos\th\cos\tphi, \,\,\,\, (j_+)_2=-\frac{1}{|q|} \sin\tphi, \,\,\,\, (j_-)_1 =\frac{1}{|q|} \cos\th\sin\tphi, \,\,\,\, (j_-)_2=\frac{1}{|q|}  \cos\tphi,
\end{split}
 \eea
where $\Jk$ and $\Jk_\pm$ are regular (even at the axis) as well as anti-self dual, i.e., $\Jk,\, \Jk_\pm\in\sk_1(\mathbb{C})$, and where the coordinate $\tphi$ involved in the definition of $j_\pm$ has been introduced in \eqref{eq:definitionof:specificchoice:normalizedcoord}. Note in particular the following identities, with $(x^1_p, x^2_p)$ introduced in \eqref{eq:definitionof:specificchoice:normalizedcoord}, 
\bea\lab{eq:usefulalgebraicidentitiesinvolvingscalarproductsReJkReJkpm}
\bsplit
\Re(\Jk)\c\Re(\Jk) &=\frac{(\sin\th)^2}{|q|^2}, \qquad \Re(\Jk)\c\Re(\Jk_+)=-\frac{x^2_p}{|q|^2},\qquad \Re(\Jk)\c\Re(\Jk_-)=\frac{x^1_p}{|q|^2},\\
\dual(\Re(\Jk))\c\Re(\Jk_+) &=\frac{\cos\th x^1_p}{|q|^2},\qquad \dual(\Re(\Jk))\c\Re(\Jk_-)=\frac{\cos\th x^2_p}{|q|^2}.
\end{split}
\eea

The complexified Ricci coefficients w.r.t. this principal null pair are given by  
\bea\lab{eq:KerrvaluesofcomplexifiedRicci}
\begin{split}
&\Xh=\Xbh=\Xi=\Xib=\omb=0, \qquad  \tr X=\frac{2\De \ov{q}}{|q|^4}, \qquad \tr\Xb=-\frac{2}{\ov{q}}, \qquad \om=- \frac 12\pr_r\left(\frac{\Delta}{|q|^2}\right),\\
&H=Z=\frac{a}{\ov{q}}\Jk=\frac{aq}{|q|^2}\Jk, \qquad \Hb=-\frac{a}{q}\Jk= -\frac{a\ov{q}}{|q|^2}\Jk.
\end{split}
\eea

The complexified curvature components are given by
\bea
A=B=\Bb=\Ab=0,\quad P=-\frac{2m}{q^3}.
\eea

Also, the principal null frame acts on the the normalized coordinates of Lemma \ref{lem:specificchoice:normalizedcoord} as follows 
\bea\lab{eq:actionofingoingprincipalnullframeonnormalizedcoordinates}
\bsplit
e_3(r)&=-1, \qquad\quad e_4(r)=\frac{\De}{|q|^2}, \qquad\qquad\qquad\qquad\quad\, e_1(r)=0, \quad\,\,\,\,\, e_2(r)=0,\\
e_3(\tau)&=\tmod'(r), \quad e_4(\tau)=\frac{2(r^2+a^2) - \De\tmod'(r)}{|q|^2}, \quad e_1(\tau)=0, \quad\,\,\,\, e_2(\tau)=\frac{a\sin\th}{|q|},\\
e_3(\th)&=0, \qquad\quad\,\,\,\,\, e_4(\th)=0, \qquad\qquad\qquad\qquad\qquad\,\, e_1(\th)=\frac{1}{|q|}, \quad e_2(\th)=0,\\
e_3(\tphi)&=\phimod'(r), \quad e_4(\tphi)=\frac{2a -\De\phimod'(r)}{|q|^2},\qquad\quad\,\,\,\,\, e_1(\tphi)=0, \quad\,\,\,\, e_2(\tphi)=\frac{1}{|q|\sin\th}.
\end{split}
\eea
In particular, recalling that $x^1_p=\sin\th\cos\tphi$ and $x^2_p=\sin\th\sin\tphi$, we have
\bea\lab{eq:actionofingoingprincipalnullframeonnormalizedcoordinates:bis}
\bsplit
e_3(x^1_p)&=-\phimod'(r)x^2_p, \qquad\, e_4(x^1_p)=-\frac{2a -\De\phimod'(r)}{|q|^2}x^2_p,\\ 
e_3(x^2_p)&=\phimod'(r)x^1_p, \qquad\quad e_4(x^2_p)=\frac{2a -\De\phimod'(r)}{|q|^2}x^1_p,
\end{split}
\eea
and, in view of the definition of $\Jk$ and $\Jk_\pm$,  
\bea\lab{eq:actionofingoingprincipalnullframeonnormalizedcoordinates:ter}
\DD(\tau)=a\Jk, \qquad \DD(\cos\th)=i\Jk,\qquad \DD(x^1_p)=\Jk_+, \qquad \DD(x^2_p)=\Jk_-.
\eea

Moreover the derivatives of $\Jk$ and $\Jk_\pm$ w.r.t. the principal null frame satisfy the following
\bea
\bsplit
\nab_3\Jk &=\frac{1}{\ov{q}}\Jk, \qquad \nab_4\Jk =- \frac{\De \ov{q}}{|q|^4}\Jk, \qquad \nab_3\Jk_\pm =\frac{1}{\ov{q}}\Jk_\pm, \qquad \nab_4 \Jk_\pm =- \frac{\De \ov{q}}{|q|^4}\Jk_{\pm} \mp  \frac{2a}{|q|^2}\Jk_{\mp},\\
\ov{\DD}\c\Jk &=\frac{4i(r^2+a^2)\cos\th}{|q|^4},\qquad \DD\hot\Jk=0, \\ 
\ov{\DD}\c \Jk_+ &= - \frac{4r^2 }{|q|^4}x^1_p - \frac{4ia^2\cos\th}{|q|^4}x^2_p, \qquad \ov{\DD}\c \Jk_- = - \frac{4r^2 }{|q|^4}x^2_p + \frac{4ia^2\cos\th}{|q|^4}x^1_p, \qquad \DD\hot \Jk_{\pm} =0.
\end{split}
\eea

%%%%%%%%%%%%%%%%%%%%%%%%%%%%%%%%%%%

\section{Teukolsky wave-transport system}
\lab{sec:sectionontheteukolskywavetransportsystem}

%%%%%%%%%%%%%%%%%%%%%%%%%%%%%%%%%%%

The main goal of this section is to state our main result concerning the Teukolsky wave-transport system in perturbations of Kerr. To this end, we first need to introduce some notations and properties for perturbations of Kerr.

%%%%%%%%%%%%%%%%%%%%%%%%%%%%%%%%%%%

\subsection{Perturbations of Kerr}
\lab{sec:Kerrpertbasic}

%%%%%%%%%%%%%%%%%%%%%%%%%%%%%%%%%%%

We consider a given vacuum spacetime $(\MM, \g)$ together with a null pair $(e_3, e_4)$ and its corresponding horizontal structure as in Section \ref{subsection:review-horiz.structures}. We will use the complexified Ricci and curvature coefficients  of Definition \ref{def:complexRicciandcurvaturecoefficients}. Moreover, we assume that $\MM$ is endowed with a pair of constants $(a, m)$, scalar functions $(\tau, r, \th, \vphi)$ and complex horizontal 1-forms $\Jk$, $\Jk_\pm$. We start by defining linearized  quantities.

%%%%%%%%%%%%%%%%%%%%%%%%%%%%%%%%%%%%%%%%%%%%%

\subsubsection{Definition of linearized quantities}
\lab{sec:definitionoflinearizedquantities:chap4}

%%%%%%%%%%%%%%%%%%%%%%%%%%%%%%%%%%%%%%%%%%%%%

Recall from Section \ref{subsect:principalnullpairinKerr} that, w.r.t. the principal null pair regular across  the future event horizon, the following quantities vanish in Kerr
\beaa
\Xh, \quad \Xbh,\quad \Xi, \quad \Xib, \quad \omb, \quad A, \quad B, \quad \Bb, \quad \Ab,\quad \nab(r), \quad e_4(\th), \quad e_3(\th), \quad \DD\hot\Jk.
\eeaa

We renormalize below all other quantities, not vanishing in Kerr as follows.
\begin{definition}
\lab{def:renormalizationofallnonsmallquantitiesinPGstructurebyKerrvalue}
We  define  the following renormalizations.
\begin{enumerate}
\item Linearization of the complex-valued Ricci and curvature coefficients:
\beaa
\bsplit
\trXc &:= \tr X-\frac{2\ov{q}\De}{|q|^4}, \qquad\trXbc := \tr\Xb+\frac{2}{\ov{q}},\\ 
\Pc &:= P+\frac{2m}{q^3},\qquad\qquad\, \omc  := \om  + \frac{1}{2}\pr_r\left(\frac{\De}{|q|^2} \right),\\
 \Hc &:= H-\frac{aq}{|q|^2}\Jk, \qquad\quad \Hbc:=\Hb+\frac{a\ov{q}}{|q|^2}\Jk,\qquad\quad \Zc := Z-\frac{aq}{|q|^2}\Jk.
 \end{split}
\eeaa

\item Linearization of derivatives of the scalar functions\footnote{Note that in Kerr we have $e_3(q)=e_3(r)$ and $e_4(q)=e_4(r)$ so that it suffices to linearize $e_3(r)$ and $e_4(r)$.} $r$, $\cos\th$, $q$, $\tau$, $x^1_p$ and $x^2_p$:
\beaa
\bsplit
\widecheck{e_3(r)} :=& e_3(r)+1, \qquad\qquad\qquad\,\,\,\,\, \widecheck{e_4(r)} := e_4(r)-\frac{\Delta}{|q|^2},\\
\widecheck{e_3(\tau)}:=& e_3(\tau)-\tmod'(r), \qquad\qquad\, \widecheck{e_4(\tau)}:=e_4(\tau)-\frac{2(r^2+a^2) - \De\tmod'(r)}{|q|^2},\\
\widecheck{e_3(x^1_p)}:=& e_3(x^1_p)+\phimod'(r)x^2_p, \qquad \widecheck{e_4(x^1_p)}:=  e_4(x^1_p)+\frac{2a -\De\phimod'(r)}{|q|^2}x^2_p, \\ \widecheck{e_3(x^2_p)}:=& e_3(x^2_p)-\phimod'(r)x^1_p, \qquad \widecheck{e_4(x^2_p)}:=e_4(x^2_p)-\frac{2a -\De\phimod'(r)}{|q|^2}x^1_p,\\
\widecheck{\DD q} :=& \DD q+a\Jk, \qquad\qquad\, \widecheck{\DD \ov{q}} :=\DD \ov{q}-a\Jk,\qquad \,\,\,\widecheck{\DD(\cos\th)} := \DD(\cos\th) -i\Jk,\\
\widecheck{\DD(\tau)}:=& \DD(\tau)-a\Jk, \qquad  \widecheck{\DD(x^1_p)} := \DD(x^1_p) - \Jk_+, \qquad \widecheck{\DD(x^1_p)} := \DD(x^2_p) - \Jk_-.
\end{split}
\eeaa

\item Linearization of derivatives of the complex 1-forms $\Jk$ and $\Jk_{\pm}$:
\beaa
\bsplit
\widecheck{\nab_3\Jk}:=&\nab_3\Jk -\frac{1}{\ov{q}}\Jk, \qquad\qquad  \widecheck{\nab_4\Jk}:=\nab_4\Jk +\frac{\De \ov{q}}{|q|^4}\Jk,\qquad \widecheck{\ov{\DD}\c\Jk}:= \ov{\DD}\c\Jk-\frac{4i(r^2+a^2)\cos\th}{|q|^4},\\
\widecheck{\nab_3\Jk_\pm}:=&\nab_3\Jk_\pm - \frac{1}{\ov{q}}\Jk_\pm\pm\phimod'(r)\Jk_{\mp},\qquad \widecheck{\nab_4 \Jk_\pm}:=\nab_4 \Jk_\pm + \frac{\De \ov{q}}{|q|^4}\Jk_{\pm} \pm\frac{2a-\De\phimod'(r)}{|q|^2}\Jk_{\mp},\\
\widecheck{\ov{\DD}\c \Jk_+}:=&\ov{\DD}\c \Jk_+ + \frac{4r^2 }{|q|^4}x^1_p + \frac{4ia^2\cos\th}{|q|^4}x^2_p,\qquad \widecheck{\ov{\DD}\c \Jk_-}:=\ov{\DD}\c \Jk_- + \frac{4r^2 }{|q|^4}x^2_p - \frac{4ia^2\cos\th}{|q|^4}x^1_p.
\end{split}
\eeaa
 \end{enumerate}
\end{definition}

%%%%%%%%%%%%%%%%%%%%%%%%%%%%%%%%%%%%%%%%%%%%%

\subsubsection{Definition of the notations $\Ga_b$ and $\Ga_g$ for error terms}
\lab{sec:definitionofGabandGagfirsttime}

%%%%%%%%%%%%%%%%%%%%%%%%%%%%%%%%%%%%%%%%%%%%%

\begin{definition}
\lab{definition.Ga_gGa_b}
The set of all linearized quantities is of the form $\Ga_g\cup \Ga_b$ with  $\Ga_g,  \Ga_b$
 defined as follows.
 \begin{enumerate}
\item 
 The set $\Ga_g$ is given by $\Ga_g=\Ga_{g,1}\cup \Ga_{g, 2}\cup\Ga_{g,3}$   with
 \bea
 \bsplit
 \Ga_{g,1} &= \Big\{\Xi, \quad \omc, \quad\trXc,\quad  \Xh,\quad \Zc,\quad \Hbc, \quad \trXbc , \quad r\Pc, \quad  rB, \quad  rA\Big\},\\
 \Ga_{g,2} &= \Big\{\widecheck{e_4(r)}, \quad r^{-1}\nab(r), \quad \widecheck{e_4(\tau)},   \quad r^{-1}\widecheck{\DD(\tau)}, \quad e_4(\cos\th)\Big\},\\
  \Ga_{g,3} &= \Big\{r\widecheck{\nab_4\Jk}\Big\}.
 \end{split}
 \eea
 
 \item The set $\Ga_b$ is given by $\Ga_b=\Ga_{b,1}\cup \Ga_{b, 2}\cup \Ga_{b,3}$   with
 \bea
 \bsplit
 \Ga_{b,1}&= \Big\{\Hc, \quad \Xbh, \quad \omb, \quad \Xib,\quad  r\Bb, \quad \Ab\Big\},\\
  \Ga_{b, 2}&= \Big\{r^{-1}\widecheck{e_3(r)}, \quad r^{-1}\widecheck{e_3(\tau)}, \quad  \widecheck{\DD(\cos\th)}, \quad e_3(\cos\th)\Big\}, \\
   \Ga_{b,3}&=\bigg\{ r\,\widecheck{\ov{\DD}\c\Jk}, \quad r\,\DD\hot\Jk, \quad r\,\widecheck{\nab_3\Jk}\bigg\}. 
   \end{split}
 \eea
\end{enumerate}
\end{definition}

\begin{remark}
The justification for the above decompositions has to do with the expected  decay properties of the linearized  components in perturbations of Kerr. More precisely, we will consider perturbations of Kerr for which $\Ga_g$ and $\Ga_b$ satisfy the following estimates,    see Section 4.1.5 in \cite{Sze} for details, 
\bea\lab{eq:expectedbehaviorGabGag:chap2}
\bsplit
\big|\dk^{\leq s}\Ga_g|&\les \ep \min\Big\{ r^{-2 }\tau^{-1/2-\dec},  \, r^{-1}\tau^{-1-\dec} \Big\}, \\
\big|\nab_3\dk^{\leq s-1}\Ga_g| &\les \ep  r^{-2 }\tau^{-1-\dec},\\
\big|\dk^{\leq s}\Ga_b\big| &\les \ep  r^{-1 }\tau^{-1-\dec},
\end{split}
\eea
for a small constant $\dec>0$, where $\dk=\{\nab_3, r\nab_4, \dkb=r\nab \}$ denotes weighted derivatives.
\end{remark}

\begin{remark}
In view of \eqref{eq:expectedbehaviorGabGag:chap2}, we note that $\Ga_g$ satisfies the assumptions of $\Ga_b$ and that $r^{-1}\Ga_b$ satisfies the assumptions of $\Ga_g$. Thus, in the rest of the paper, we will systematically replace $\Ga_g+\Ga_b$ by $\Ga_b$ and $r^{-1}\Ga_b+\Ga_g$ by $\Ga_g$.  
\end{remark}

\begin{remark}\lab{rmk:linearizedquantitiesJkpmandxbpnotincludedindefintionGagGab}
Note that we do not include linearized quantities associated to $\Jk_{\pm}$ and first order derivatives of $x^1_p$ and $x^2_p$ in the definition of $(\Ga_b, \Ga_g)$. Indeed, these quantities will only be used in the context of energy-Morawetz estimates for few derivatives when applying the results of \cite{MaSz24} \cite{MaSz26}, see Section \ref{sec:assumptionsforsec:energyMorawetzesitmatesforTeukoslkyonMM:upto15derivatives}, while we will need to estimate $(\Ga_b, \Ga_g)$ up to top order derivatives in Section 6 of \cite{Sze}.
\end{remark}

%%%%%%%%%%%%%%%%%%%%%%%%%%

\subsubsection{Commutation formulas revisited}

%%%%%%%%%%%%%%%%%%%%%%%%%%

The following lemma extends the commutation formulas of Section \ref{sec:generalcommutationformulasrealcase} to complex derivatives and makes use of the notations $(\Ga_b, \Ga_g)$ of Definition \ref{definition.Ga_gGa_b}.

   \begin{lemma}\label{LEMMA:COMMUTATION-FORMULAS-1}
   The following commutation formulas hold true.
   \begin{enumerate}
   \item 
   Let $h \in \sk_0(\CCC)$ $s$-conformally invariant. Then 
   \bea\label{eq:comm-nab4-nab3-DD-h-precise}
   \begin{split}
 \, [\nab_4 , \DD]h  &=  -\frac{1}{2}\tr X\DD h+(\Hb+Z)\nab_4 h -\frac 1 2 \Xh \c\ov{\DD} h+\Xi \nab_3h , \\
 \, [\nab_3 , \DD]h  &=   -\frac{1}{2}\tr \Xb\DD h+(H-Z)\nab_3 h -\frac 1 2 \Xbh\c\ov{\DD} h+\Xib \nab_4h. 
\end{split}
 \eea    
    
   \item
    Let $F\in  \sk_1 (\mathbb{C})$.  Then
    \bea\label{eq:comm:nab4-nab3-DDhot-precise}
    \begin{split}
\, [\nab_4,  \DD \hot] F  &=-\frac 1 2 \tr X\left( \DD \hot  F+\Hb\hot F\right)+(\Hb+Z)\hot\nab_4 F+ \Xi \hot \nab_3 F \\
 &-B \hot F - \frac 1 2 \tr \Xb \Xi \hot  F-\frac 1 2\Xh \c \ov{\DD} F+\frac12\Xh (\ov{\Hb}\c F)+ (\Ga_b \c \Ga_g) F, \\
\, [\nab_3,  \DD \hot] F  &=-\frac 1 2 \tr \Xb\left( \DD \hot  F+H \hot F\right)+(H-Z)\hot\nab_3 F+ \Xib \hot \nab_4 F \\
 &+\Bb \hot F - \frac 1 2 \tr X \Xib \hot  F-\frac 1 2\Xbh \c \ov{\DD} F+\frac12\Xbh (\ov{H}\c F)+ (\Ga_b \c \Ga_g) F.
 \end{split}
\eea
Using the schematic structure of the error terms, the above can be written as
\bea
 \lab{commutator-nab-43-D-hot} 
\bsplit
\, [\nab_4, \mathcal{D}\hot ]F&=- \frac 1 2 \tr X( \mathcal{D}\hot F + \underline{H} \hot F)+ (\underline{H}+Z) \hot \nab_4 F+\Xi \c \nabc_3 F+ r^{-1} \Ga_g \c  \dk^{\leq 1} F, \\
\, [\nab_3, \mathcal{D}\hot] F&=- \frac 1 2 \tr \Xb( \mathcal{D}\hot F + H \hot F)+ (H-Z) \hot \nab_3 F+r^{-1}\Ga_b \c \dk^{\leq 1} F.
\end{split}
\eea

\item Let $U\in \sk_2(\mathbb{C})$. Then
\bea\label{eq:comm:nab4-nab3-DDc-precise}
\begin{split}
\,[\nab_4,  \ov{\DD} \c ] U &=-\frac 1 2\ov{\tr X} \big(  \ov{\DD}\c U - 2\ov{\Hb} \c U\big) +\ov{(\Hb + Z)}\c\nab_4 U+ \ov{\Xi} \c \nab_3 U \\
&+2\ov{B} \c U  -\frac 1 2 \ov{\tr \Xb} \ov{\Xi}\c  U  -\frac 1 2 \Xh \c \ov{\DD} U-\frac 1 2 (\ov{\Xh}\c U)\ov{\Hb}+ (\Ga_b \c \Ga_g) U, \\
\,[\nab_3,  \ov{\DD} \c ] U &=-\frac 1 2\ov{\tr \Xb} \big(  \ov{\DD}\c U - 2\ov{H} \c U\big) +\ov{(H- Z)}\c\nab_3 U+ \ov{\Xib} \c \nab_4 U \\
&-2\ov{\Bb} \c U  -\frac 1 2 \ov{\tr X} \ov{\Xib}\c  U  -\frac 1 2 \Xbh \c \ov{\DD} U-\frac 1 2 (\ov{\Xbh}\c U)\ov{H}+ (\Ga_b \c \Ga_g) U.
\end{split}
\eea
Using the schematic structure of the error terms, the above can be written as
 \bea\label{commutator-nabc-3-ov-DDc-U}
 \bsplit
\, [\nab_4, \ov{\DD}\c] U&=- \frac 1 2\ov{\tr X}\, ( \ov{\DD} \c U - 2 \ov{\Hb} \c U)+\ov{(\Hb+Z)} \c \nab_4 U+\Xi \c \nabc_3 U +r^{-1} \Ga_g  \c \dk^{\leq 1} U, \\
\, [\nab_3, \ov{\DD}\c] U&=- \frac 1 2\ov{\tr\Xb}\, ( \ov{\DD} \c U -  2 \ov{H} \c U)+\ov{(H-Z)} \c \nab_3 U + r^{-1}\Ga_b \c  \dk^{\leq 1} U.
\end{split}
\eea
Similarly, for $F\in  \sk_1 (\mathbb{C})$ we have
 \bea\label{commutator-nab4-ov-DDcF}
 \bsplit
\, [\nab_4, \ov{\DD}\c] F&=- \frac 1 2\ov{\tr X}\, ( \ov{\DD} \c F -  \ov{\Hb} \c F)+\ov{(\Hb+Z)} \c \nab_4 F +\ov{\Xi} \c \nabc_3 F+ r^{-1} \Ga_g  \c \dk^{\leq 1} F, \\
\, [\nab_3, \ov{\DD}\c] F&=- \frac 1 2\ov{\tr\Xb}\, ( \ov{\DD} \c F -   \ov{H} \c F)+\ov{(H-Z)} \c \nab_3 F +r^{-1}\Ga_b \c  \dk^{\leq 1} F.
\end{split}
\eea
\item Let $U\in \sk_2(\mathbb{C})$. Then
\bea\label{correct-commutator-1} \lab{commutator-nab-4-D-c} \label{commutator-nab-3-nab-4-Psi}\label{correct-commutator}
\begin{split}
\,[\nab_3, \nab_4]U &= - 2\om \nab_3 U+ 2\omb \nab_4 U  + 2 (\eta_c-\etab_c) \nab_c U +4i \left(- \rhod+ \eta \wedge \etab  \right) U+\left(\Ga_b  \c \Ga_g \right) U.
\end{split}
\eea
Also,
\bea\label{commutator-nab-3-nab-a-U}
\bsplit
\, [\nab_3, \nab_a] U_{bc}&=-\frac  1 2   \trchb\, \Big(\nab_a U_{bc}+\eta_bU_{ac}+\eta_c U_{ab}-\de_{a b}(\eta \c U)_c-\de_{a c}(\eta \c U)_b \Big)\\
&-\frac 1 2 \atrchb\, \Big(\dual \nab_a  U_{bc} +\eta_b \dual U_{ac}+\eta_c \dual U_{ab}- \in_{a b}(\eta \c  U)_c- \in_{a c}(\eta \c  U)_b \Big)\\
 &+(\eta_a-\ze_a)\nab_3 U_{bc}+r^{-1}\Ga_b \c \dk^{\leq 1} U,
 \end{split}
\eea
\end{enumerate}
where the above error terms may also contain terms which are 
quadratic   in  the perturbation and  enjoy  better  decay  properties,  or are higher order  and decay at least as good.
\end{lemma}

\begin{proof}
See  Lemma 4.2.1 in \cite{GKS22}. 
\end{proof}

%%%%%%%%%%%%%%%%%%%%%%%%%%%%%%%%%%

\subsubsection{Approximate Killing vectorfields $\T$ and $\Z$}

%%%%%%%%%%%%%%%%%%%%%%%%%%%%%%%%%%

\begin{definition}\lab{Definition:vfsTZ} 
In $\MM$, we define $\T$ and $\Z$ as follows:
\bea\lab{eq:definitionofTandPhithataretheapproximateKillingvectorifeldinKerrpert}
\T := \frac{1}{2}\left(e_4+\frac{\Delta}{|q|^2}e_3 -2a\Re(\Jk)^be_b\right), \quad \Z := \frac 1 2 \left(2(r^2+a^2)\Re(\Jk)^be_b -a(\sin\th)^2 e_4 -\frac{a(\sin\th)^2\De}{ |q|^2} e_3\right).
\eea
\end{definition}

The following lemma provides basic relations between $\Lieb_\T, \Lieb_\Z$ and $\nab_\T, \nab_\Z$.
 \begin{lemma}\lab{lemma:basicpropertiesLiebTfasdiuhakdisug:chap9}
For a horizontal covariant k-tensor $U$, we have
\beaa
\nab_\T U_{b_1\cdots b_k} &=& \Lieb_\T U_{b_1\cdots b_k} +\frac{2amr\cos\th}{|q|^4}\sum_{j=1}^k\in_{b_jc} U_{b_1\cdots c\cdots b_k}+\Ga_b \c U,\\
\nab_\Z U_{b_1\cdots b_k} &=& \Lieb_\Z U_{b_1\cdots b_k} -\frac{\cos\th((r^2+a^2)^2-a^2(\sin\th)^2\De)}{|q|^4}\sum_{j=1}^k\in_{b_jc} U_{b_1\cdots c\cdots b_k}+r\Ga_b \c U.
\eeaa
\end{lemma}

\begin{proof}
See Lemma 4.3.6 in \cite{GKS22}.
\end{proof}

\begin{lemma}\lab{lemma:commutatorbetweenLieTLieZandnabnab4nab3}
We have 
\bea\lab{commutatorbetweenLieTLieZandnabnab4nab3:1}
\bsplit
 [\nab, \Lieb_T]U=&  \Ga_b \c \nab_3U +r^{-1}\dk^{\leq 1}(\Ga_b \c U),\\ 
   [\nab_4, \Lieb_T]U  =& \Ga_g \c \nab_3U+r^{-1}\dk^{\leq 1}(\Ga_b \c U),\\
    [\nab_3, \Lieb_T]U=& \dk^{\leq 1}(\Ga_b \c U),
\end{split}
\eea
and {
\bea\lab{commutatorbetweenLieTLieZandnabnab4nab3:2}
\bsplit
 [\nab, \Lieb_Z]U =&  r\Ga_g \c \nab_3U +\dk^{\leq 1}(\Ga_b\c U),\\ 
   [\nab_4, \Lieb_Z]U=&  r \xi \c \nab_3U+\dk^{\leq 1}(\Ga_g \c U),\\
     [\nab_3, \Lieb_Z]U=     & r\Ga_b\nab_3 U +r\dk^{\leq 1}(\Ga_g \c  U).
\end{split}
\eea}
This yields in particular the following non-sharp commutators
\bea
[r\nab, \Lieb_\T]=r\Ga_b\dk^{\leq 1}, \qquad [r\nab, \Lieb_\Z]=r^2\Ga_g\dk^{\leq 1}.
\eea
\end{lemma} 

\begin{proof}
See Lemma C.5.2 and (4.3.3) in \cite{GKS22}.
\end{proof}

%%%%%%%%%%%%%%%%%%%%%%%%%%%%%%%%%%%

\subsubsection{Decomposition of the wave operator in null frames}

%%%%%%%%%%%%%%%%%%%%%%%%%%%%%%%%%%%

The following lemma provides the decomposition of $\squared_k$ in null frames. 
\begin{lemma}\label{lemma:expression-wave-operator}
The wave operator for $\psi\in {\sk_k(\mathbb{C})}$, $k=0,1,2$, is given by
\bea\label{eq:wave-squared}
\begin{split}
\squared_k \psi&=-\nab_4 \nab_3 \psi  -\frac 1 2 \trchb \nab_4\psi+\left(2\om -\frac 1 2 \trch\right) \nab_3\psi+\lap_k\psi+2\etab \c\nab \psi \\
&+ ki \left( \rhod- \eta \wedge \etab \right) \psi+(\Ga_b \c \Ga_g) \c \psi,\\
\squared_k \psi&=-\nab_3\nab_4\psi +\left(2\omb -\frac 1 2 \trchb\right)\nab_4\psi -\frac 1 2 \trch \nab_3\psi+\lap_k\psi+2\eta \c\nab \psi \\
&- ki \left( \rhod- \eta \wedge \etab \right) \psi+(\Ga_b \c \Ga_g) \c \psi,
\end{split}
\eea
where $\lap_k=\nab^a \nab_a$ denotes the horizontal Laplacian for $k$-tensors. 
\end{lemma}

\begin{proof}
See Lemma 4.7.5 in \cite{GKS22} for the first identity of \eqref{eq:wave-squared} in the case $k=2$. The proof of Lemma 4.7.5 in \cite{GKS22} immediately extends to $k=0,1,2$ and to the second identity of \eqref{eq:wave-squared}.
\end{proof}

We will also need to decompose the following analog of Lemma \ref{lemma:expression-wave-operator} using complex derivatives.
\begin{corollary}\label{corollary-wave-complex} 
We have, for $\psi\in \sk_2(\CCC)$,
\bea
\begin{split}
\squared_2 \psi&=-\nab_4 \nab_3 \psi +\frac 1 4  \DD\hot( \DDb \c \psi)+\left(2\om -\frac 1 2 \tr X\right) \nab_3\psi- \frac 1 2 \tr\Xb \nab_4\psi+2\etab \c\nab \psi \\
& +  \left( - \frac 1 4 \tr X \ov{\tr \Xb}- \frac 1 4 \tr \Xb  \ov{\tr X}  - 2\ov{P}\right) \psi- 2i \left(\eta \wedge \etab\right)  \psi+(\Ga_b \c \Ga_g) \c \psi.
\end{split}
\eea
\end{corollary}

\begin{proof}
See Corollary 4.7.9 in \cite{GKS22}.
\end{proof}

%%%%%%%%%%%%%%%%%%%%%%%%%%%%%%%%%%%%%

\subsection{Teukolsky wave-transport system in perturbations of Kerr}
\lab{sec:precisederivationTeukolskywave-transportsystem:Kerrpert}

%%%%%%%%%%%%%%%%%%%%%%%%%%%%%%%%%%%%%

%%%%%%%%%%%%%%%%%%%%%%%%%%%%%%%%%%%%%%%%%%%%%%%%%%

\subsubsection{Definition  of $\pmb\phi_{s}^{(p)}$, $s=\pm 2$, $p=0,1,2$, in perturbations of Kerr}
\lab{sec:definitionofpmbphisspinpertubrtionsofKerr}

%%%%%%%%%%%%%%%%%%%%%%%%%%%%%%%%%%%%%%%%%%%%%%%%%%

We start with the case $s=+2$. In perturbations of Kerr, we define $\pmb\phi_{+2}^{(p)}\in\sk_2(\mathbb{C})$, $p=0,1,2$ as follows
\bea\lab{eq:definitionofthephiplus2phierarchy:perturbationofKerr}
\bsplit
\pmb\phi_{+2}^{(0)} &:=\frac{\ov{q}}{q}A,\\
\pmb\phi_{+2}^{(1)} &:= \ov{q}^2\left(\nabc_3+\frac{1}{2}\trchb -\frac{3}{2}\frac{\atrchb^2}{\trchb} -2i\atrchb\right)A,\\
\pmb\phi_{+2}^{(2)} &:=|q|^2\ov{q}^2\left( \nabc_3\nabc_3 A + C_1  \nabc_3A + C_2   A\right),
\end{split}
\eea
where the scalar function $C_1$, $C_2$ are given by\footnote{Notice that we always work with a frame in ingoing normalization, so that $\trchb=-\frac{2r}{|q|^2}+O(\frac{\ep}{r^2})\neq 0$ which justifies the fact that we may divide by $\trchb$ in \eqref{eq:C1-C2-comparison-Ma}.} 
\bea\lab{eq:C1-C2-comparison-Ma}
\begin{split}
C_1&=2\trchb -\frac{2\atrchb^2}{\trchb}  -4 i \atrchb, \\
C_2  &= \frac 1 2 \trchb^2- 4\atrchb^2+\frac 3 2\frac{\atrchb^4}{\trchb^2} +  i \left(-2\trchb\atrchb +\frac{4\atrchb^3}{\trchb}\right).
\end{split}
\eea

\begin{remark}\lab{rmk:compasisionphiplus2p=0withMaSz26andphiplus2p=2withqfGKS22}
Note that the definition of $\pmb\phi_{+2}^{(0)}$ in \eqref{eq:definitionofthephiplus2phierarchy:perturbationofKerr} coincides with the one in (5.31) of \cite{MaSz26}. Also, comparing the definition of $\pmb\phi_{+2}^{(2)}$ in \eqref{eq:definitionofthephiplus2phierarchy:perturbationofKerr} with Definition 5.2.2 in \cite{GKS22} for $\qf$, we have $\qf=\pmb\phi_{+2}^{(2)}$. 
\end{remark}

Next, we consider the case $s=-2$. In perturbations of Kerr, we define $\pmb\phi_{-2}^{(p)}\in\sk_2(\mathbb{C})$, $p=0,1,2$ as follows
\bea\lab{eq:definitionofthephiminus2phierarchy:perturbationofKerr}
\bsplit
\pmb\phi_{-2}^{(0)} &:=\frac{q}{\ov{q}}\left(\frac{\De}{|q|^2}\right)^2\Ab,\\
\pmb\phi_{-2}^{(1)} &:=q^2\frac{\De}{|q|^2}\left(\nabc_4 +\frac{1}{2}\trch  -\frac{3}{2}\left(\frac{\atrch}{\trch}\right)_{\vartheta}\atrch -2i\atrch\right)\Ab,\\
\pmb\phi_{-2}^{(2)} &:=|q|^2q^2\left( \nabc_4\nabc_4\Ab + \underline{C}_1  \nabc_4\Ab + \underline{C}_2   \Ab\right),
\end{split}
\eea
where the scalar function $\underline{C}_1$, $\underline{C}_2$ are given by
\bea\label{eq:Cb1-Cb2-comparison-Ma}
\bsplit
\und{C}_1=& 2\trch -2\left(\frac{\atrch}{\trch}\right)_{\vartheta}\atrch -4i\atrch, \\ 
\und{C}_2  =& \frac 1 2 \trch^2 -4(\atrchb)^2 +\frac 3 2\left(\left(\frac{\atrch}{\trch}\right)_{\vartheta}\right)^2\atrch^2 +i\left(-2\trch\atrch +4\left(\frac{\atrch}{\trch}\right)_{\vartheta}\atrch^2\right),
\end{split}
\eea
and where we have introduced the following 0-conformally invariant scalar
\bea\lab{eq:defintionatrchovertrchvarthetamodificationtoavoiddividingbytrchnearHHplus}
\left(\frac{\atrch}{\trch}\right)_{\vartheta}:=\vartheta(r)\frac{a\cos\th}{r}+(1-\vartheta(r))\frac{\atrch}{\trch},
\eea
with $\vartheta$ supported in $r\leq 11m$ and $\vartheta=1$ for $r\leq 10m$.

\begin{remark}
Since we will always work with the frame in ingoing normalization, so that $\trch=\frac{2\De r}{|q|^4}+O(\frac{\ep}{r^2})$, we cannot divide by $\trch$ near $r=r_+$ which justifies the definition of the 0-conformally invariant scalar in   \eqref{eq:defintionatrchovertrchvarthetamodificationtoavoiddividingbytrchnearHHplus}.
\end{remark}

\begin{remark}\lab{rmk:compasisionphiminus2p=0withMaSz26andphiminus2p=2withqfGKS22}
Note that the definition of $\pmb\phi_{-2}^{(0)}$ in \eqref{eq:definitionofthephiminus2phierarchy:perturbationofKerr} coincides with the one in (5.31) of \cite{MaSz26}. Also, comparing the definition of $\pmb\phi_{-2}^{(2)}$ in \eqref{eq:definitionofthephiminus2phierarchy:perturbationofKerr} with Definition 5.2.2 in \cite{GKS22} for $\qfb$, we have $\qfb=\pmb\phi_{-2}^{(2)}$. 
\end{remark}

%%%%%%%%%%%%%%%%%%%%%%%%%%%%%%%%%%%%%%

\subsubsection{Teukolsky wave-transport system in perturbations of Kerr}

%%%%%%%%%%%%%%%%%%%%%%%%%%%%%%%%%%%%%%

The following theorem provides the tensorial wave equations satisfied by $\pmb\phi_{s}^{(p)}$, $s=\pm 2$, $p=0,1,2$.
\begin{theorem}\lab{thm:derivationoftheTeukolskytensorialwavesystemfors=plusminus2:kerrpert:alternateformnullframeinsteadcoordvectorfield}
Let $(\Ga_g, \Ga_b)$ be given by Definition \ref{definition.Ga_gGa_b}. Then, the horizontal tensors $\pmb\phi_{s}^{(p)}\in\sk_2(\mathbb{C})$, $s=\pm 2$, $p=0,1,2$, defined in \eqref{eq:definitionofthephiplus2phierarchy:perturbationofKerr} \eqref{eq:definitionofthephiminus2phierarchy:perturbationofKerr} satisfy the following  tensorial wave equations
\bsub
\lab{eq:TensorialTeuSysandlinearterms:rescaleRHScontaine2:general:Kerrperturbation:alternateformnullframeinsteadcoordvectorfield}
\bea
\lab{eq:TensorialTeuSys:rescaleRHScontaine2:general:Kerrperturbation:alternateformnullframeinsteadcoordvectorfield}
\bigg(\squared_2 -\frac{4ia\cos\th}{|q|^2}\nab_{\T}- \frac{4-2\de_{p0}}{|q|^2}\bigg){\phis{p}} = \widetilde{\L}_{s}^{(p)}[\pmb\phi_{s}]+\widetilde{\N}_{W,s}^{(p)}, \quad s=\pm2, \quad p=0,1,2,
\eea
where the linear coupling terms $\widetilde{\L}_{s}^{(p)}[{\pmb\phi_s}]$ have the following schematic forms
\bea
\lab{eq:tensor:Lsn:onlye_2present:general:Kerrperturbation:alternateformnullframeinsteadcoordvectorfield}
\bsplit
{\widetilde{\L}_{s}^{(0)}[\pmb\phi_{s}]}={}& (2sr^{-3} +O(mr^{-4}))\phis{1}+ O(mr^{-3}) \nab_{\widehat{\mathcal{X}}_s}^{\leq 1}\phis{0},\\
{\widetilde{\L}_{s}^{(1)}[\pmb\phi_{s}]}={}& (sr^{-3} +O(mr^{-4}))\phis{2}+ O(mr^{-3}) \nab_{\widehat{\mathcal{X}}_s}^{\leq 1}  \phis{1}+O(mr^{-2})\nab_{\widehat{Z}}^{\leq 1}\phis{0},\\
{\widetilde{\L}_{s}^{(2)}[\pmb\phi_{s}]}={}&O(mr^{-3})\phis{2}+O(mr^{-2})\nab_{\Z+a\T}^{\leq 1}\phis{1}+O(m^2 r^{-2})\phis{0},
\end{split}
\eea
with $\widehat{\mathcal{X}}_s$, $s=\pm 2$, and $\widehat{Z}$ being the regular vectorfields defined respectively by
\bea\lab{eq:formofregularhorizontalvectorfieldwidetildemathcalXs:Kerrperturbation}
\widehat{\mathcal{X}}_s:=|q|^2\left(s\Re(\Jk)\c\nab -\frac{2a\cos\th}{r}\dual\Re(\Jk)\c\nab\right), \quad s=\pm 2,
\eea
and 
\bea\lab{eq:formofregularhorizontalvectorfieldwidehatZ:Kerrperturbation}
\widehat{Z}:=|q|^2\left(\Re(\Jk)\c\nab +\frac{a(\cos\th)^2}{|q|^2}\nab_{\T}\right),
\eea
with all the coefficients in  \eqref{eq:tensor:Lsn:onlye_2present:general:Kerrperturbation:alternateformnullframeinsteadcoordvectorfield} being independent of coordinates $\tau$ and 
$\tphi$, and with the coefficients in front of the terms $\phis{2}$ and $\nab_{\Z+a\T}\phis{1}$ on the RHS of equation of ${\widetilde{\L}_{s}^{(2)}[\pmb\phi_{s}]}$ in \eqref{eq:tensor:Lsn:onlye_2present:general:Kerrperturbation:alternateformnullframeinsteadcoordvectorfield} being real functions, and where the nonlinear correction terms $\widetilde{\N}_{W,s}^{(p)}$ have the following schematic form
\bea\lab{eq:schematicformofNpWsplus2} 
\bsplit
\widetilde{\N}^{(0)}_{W,+2} =& r^{-1}  \dk^{\leq 1}\big( \Ga_g \c  B\big) + \nab_3\Xi  \c B {+\big( r^{-1}\dkb^{\leq1} \Pc+ r^{-3} \Ga_b \big)\c \Xi +(\Ga_g \c \Xi) \c \Bb}\\
 & +(\Ga_b \c \Ga_g)\c A +\Ga_b\nab_3A+r^{-1}\dk^{\leq 1}(\Ga_bA),\\
 \widetilde{\N}^{(1)}_{W,+2} =& r\dk^{\leq 2}\big( \Ga_g \c  B\big)+r^2\dk^{\leq 1}\big(\nab_3\Xi  \c B) +r\dk^{\leq 1}\big((\dkb^{\leq1} \Pc+ r^{-2} \Ga_b)\c \Xi\big)\\
& +r^2\dk^{\leq 1}((\Ga_g \c \Xi) \c \Bb) +r^2\dk^{\leq 1}(\Ga_b\c\nab_3A)+r\dk^{\leq 2}(\Ga_b\c A)+r^2\dk^{\leq 1}\big(\Ga_b \c \Ga_g\c A\big),\\
 \widetilde{\N}^{(2)}_{W,+2} =& r^2 \dk^{\leq 3} (\Ga_g \c (A, B))+ {\nab_3 (r^3 \dk^{\leq 2}( \Ga_g \c (A, B)))} + \nab_3 (r^3\dk^{\leq1}\big(\Xi \c (\dkb^{\leq1} \Pc,  r^{-2}\Ga_b )\big))\\
&+\dk^{\leq 1} (\Ga_g \c \qf) {+\nab_3(r^4(\dk^{\leq 1}(\Ga_g\c\Xi)\c\Bb))},
\end{split}
\eea
\bea\lab{eq:schematicformofNpWsminus2}
   \bsplit
  \widetilde{\N}^{(0)}_{W,-2} =& r^{-1}\dk^{\leq 1}\big(  \Ga_b \c \Bb \big) + \Ga_b \c \Ga_b \c \Ga_g +r^{-1}\Ga_g\pmb\phi_{-2}^{(1)}+\Ga_g\nab_3\Ab +r^{-1}\Ga_b\dk^{\leq 1}\Ab,\\
 \widetilde{\N}^{(1)}_{W,-2} =& \dk^{\leq 2}\big(\Ga_b\c \Ga_g\big)+r\dk^{\leq 2}(\Ga_b\c\xi),\\
  \widetilde{\N}^{(2)}_{W,-2} =&  r^2 \dk^{\leq 2}(\Ga_b \c (A, B))+ \dk^{\leq 3} (\Ga_g \c \Ga_b),
\end{split}
 \eea
\esub
under the following additional properties for the global null frame of $\MM$:
\begin{itemize}
\item For \eqref{eq:schematicformofNpWsplus2} to hold, we require in addition that 
\bea\label{eq:additional-conditions-Hc-Xi-gRW-eq:0}
\Hc \in \Ga_g, \qquad \nab_3\Xi\in r^{-1}\dk^{\leq 1}\Ga_g.
\eea

\item For \eqref{eq:schematicformofNpWsminus2} to hold, we require in addition that 
\bea
\lab{eq:Xi-Hb-chapter12:0}
\Xi\in r^{-2}\Ga_g, \qquad \Hbc\in r^{-1}\Ga_g.
\eea
\end{itemize}
\end{theorem}

The following proposition provides the transport equations satisfied by $\pmb\phi_{s}^{(p)}$, $s=\pm 2$, $p=0,1,2$.
\begin{proposition}\lab{prop:transportequationphis=plusminus2p=0and1}
The horizontal tensors $\pmb\phi_{s}^{(p)}\in\sk_2(\mathbb{C})$, $s=\pm 2$, $p=0,1,2$, defined in \eqref{eq:definitionofthephiplus2phierarchy:perturbationofKerr} \eqref{eq:definitionofthephiminus2phierarchy:perturbationofKerr} satisfy the following transport equations in ingoing normalization 
\bsub\lab{eq:transportequationsins=plus2andminus2caseforp=0and1}
\bea\lab{eq:transportequationsins=+2caseforp=0and1}
\nab_3\left(\frac{r\ov{q}}{q}\left(\frac{r^2}{|q|^2}\right)^{p-2}\pmb\phi_{+2}^{(p)}\right) &=& \frac{\ov{q}}{rq}\left(\frac{r^2}{|q|^2}\right)^{p-1}\pmb\phi_{+2}^{(p+1)}+\N_{T,+2}^{(p)},
\eea
and
\bea\lab{eq:transportequationsins=-2caseforp=0and1}
\nab_4\left(\frac{rq}{\ov{q}}\left(\frac{r^2}{|q|^2}\right)^{p-2}\pmb\phi_{-2}^{(p)}\right) = \frac{q}{r\ov{q}}\left(\frac{r^2}{|q|^2}\right)^{p-1}\frac{\De}{|q|^2}\pmb\phi_{-2}^{(p+1)}+\N_{T,-2}^{(p)},
\eea
where the error terms $\N_{T,+2}^{(p)}$, $p=0,1$, are given by 
\bea\lab{eq:transportequationsins=+2caseforp=0and1:RHSNT+2p=0and1}
\N_{T,+2}^{(0)}=r\Ga_bA, \qquad \N_{T,+2}^{(1)}=r^3\Ga_b\nab_3A+r^2\dk^{\leq 1}(\Ga_b)A,
\eea
and where the error terms $\N_{T,-2}^{(p)}$, $p=0,1$, are given by 
\bea\lab{eq:transportequationsins=-2caseforp=0and1:RHSNT-2p=0and1}
\N_{T,-2}^{(0)}=r\Ga_g\c\Ga_b, \qquad \N_{T,-2}^{(1)}=r^2\dk^{\leq 1}(\Xi)\c\Ga_b+r\Ga_g\dk^{\leq 1}\Ga_b.
\eea
\esub
\end{proposition}

\begin{remark}
Note that the transport equations \eqref{eq:transportequationsins=+2caseforp=0and1} \eqref{eq:transportequationsins=-2caseforp=0and1} coincide with the ones in (5.33a) (5.33b) of \cite{MaSz26}. 
\end{remark}

The proof of Theorem \ref{thm:derivationoftheTeukolskytensorialwavesystemfors=plusminus2:kerrpert:alternateformnullframeinsteadcoordvectorfield} is postponed to Sections \ref{sec:proofofthm:derivationoftheTeukolskytensorialwavesystemfors=plusminus2:kerrpert:alternateformnullframeinsteadcoordvectorfield:cases=+2} and \ref{sec:proofofthm:derivationoftheTeukolskytensorialwavesystemfors=plusminus2:kerrpert:alternateformnullframeinsteadcoordvectorfield:cases=-2} which deal respectively with the case $s=+2$ and $s=-2$. We now proceed with the proof of Proposition \ref{prop:transportequationphis=plusminus2p=0and1}.

%%%%%%%%%%%%%%%%%%%%%%%%%%%%%%%%%%%%%%%%%%%%%%

\subsubsection{Proof of Proposition \ref{prop:transportequationphis=plusminus2p=0and1}} 
\lab{sec:proofof:prop:transportequationphis=plusminus2p=0and1}

%%%%%%%%%%%%%%%%%%%%%%%%%%%%%%%%%%%%%%%%%%%%%%

The goal of this section is to prove Proposition \ref{prop:transportequationphis=plusminus2p=0and1}, i.e., the fact that $\pmb\phi_{s}^{(p)}$, $=\pm 2$, $p=0,1$, verifies the transport equations \eqref{eq:transportequationsins=plus2andminus2caseforp=0and1}. We start with the case $s=+2$.

%%%%%%%%%%%%%%%%%%%%%%%%%%%%%%%%%%%%%%%

\paragraph{\textit{Transport equations for $\pmb\phi_{+2}^{(p)}$, $p=0,1$}.}

%%%%%%%%%%%%%%%%%%%%%%%%%%%%%%%%%%%%%%%

We first prove Proposition \ref{prop:transportequationphis=plusminus2p=0and1} in the case $s=+2$, i.e., the fact that $\pmb\phi_{+2}^{(p)}$, $p=0,1$ satisfies \eqref{eq:transportequationsins=+2caseforp=0and1} \eqref{eq:transportequationsins=+2caseforp=0and1:RHSNT+2p=0and1}.

We start with the case $p=0$. In view of the definition for $\pmb\phi_{+2}^{(0)}$ in \eqref{eq:definitionofthephiplus2phierarchy:perturbationofKerr}, we have
\beaa
\frac{rq}{\ov{q}}\frac{r^2}{|q|^2}\nabc_3\left(\frac{r\ov{q}}{q}\frac{|q|^4}{r^4}\pmb\phi_{+2}^{(0)}\right) &=& \frac{r^3q}{|q|^2\ov{q}}\nabc_3\left(\frac{\ov{q}^4}{r^3}A\right)\\
&=& \ov{q}^2\left(\nabc_3A +\left(\frac{4}{\ov{q}}e_3(\ov{q}) -\frac{3}{r}e_3(r)\right)A\right).
\eeaa 
Now, we use the following conformally invariant identities
\beaa
\frac{e_3(\ov{q})}{\ov{q}}=\frac{1}{2}\tr\Xb+\Ga_b, \qquad \frac{e_3(r)}{r}=\frac{1}{2}\trchb+\frac{1}{2}\frac{\atrchb^2}{\trchb}+\Ga_b,
\eeaa
which implies
\beaa
\frac{r^3q}{|q|^2\ov{q}}\nabc_3\left(\frac{|q|^4\ov{q}}{qr^3}\pmb\phi_{+2}^{(0)}\right) &=& \ov{q}^2\left(\nabc_3A +\left(2\tr\Xb -\frac{3}{2}\trchb- \frac{3}{2}\frac{\atrchb^2}{\trchb}\right)A\right)+r^2\Ga_bA\\
&=& \ov{q}^2\left(\nabc_3+\frac{1}{2}\trchb - \frac{3}{2}\frac{\atrchb^2}{\trchb} -2i\atrchb\right)A+r^2\Ga_bA\\
&=& \pmb\phi_{+2}^{(1)}+r^2\Ga_bA
\eeaa
and hence, using also $\omb\in\Ga_b$ given our ingoing normalization, 
\beaa
\nab_3\left(\frac{r\ov{q}}{q}\left(\frac{|q|^2}{r^2}\right)^2\pmb\phi_{+2}^{(0)}\right) &=& \frac{\ov{q}}{rq}\frac{|q|^2}{r^2}\pmb\phi_{+2}^{(1)}+r\Ga_bA,
\eeaa
which is the stated estimate for $p=0$.

Next, we consider the case $p=1$. In view of the definition for $\pmb\phi_{+2}^{(1)}$ in \eqref{eq:definitionofthephiplus2phierarchy:perturbationofKerr}, we have
\beaa
&&\frac{rq}{\ov{q}}\nabc_3\left(\frac{r\ov{q}}{q}\frac{|q|^2}{r^2}\pmb\phi_{+2}^{(1)}\right)\\
&=& \frac{rq}{\ov{q}}\nabc_3\left(\frac{\ov{q}^4}{r}\left(\nabc_3+\frac{1}{2}\trchb -\frac{3}{2}\frac{\atrchb^2}{\trchb} -2i\atrchb\right)A\right)\\
&=& |q|^2\ov{q}^2\left(\nabc_3+\frac{4}{\ov{q}}e_3(\ov{q}) - \frac{1}{r}e_3(r)\right)\left(\nabc_3+\frac{1}{2}\trchb - \frac{3}{2}\frac{\atrchb^2}{\trchb} -2i\atrchb\right)A\\
&=& |q|^2\ov{q}^2\left(\nabc_3+2\tr\Xb -  \frac{1}{2}\trchb - \frac{1}{2}\frac{\atrchb^2}{\trchb}\right)\left(\nabc_3+\frac{1}{2}\trchb - \frac{3}{2}\frac{\atrchb^2}{\trchb} -2i\atrchb\right)A\\
&&+r^4\Ga_b\nabc_3A+r^3\Ga_bA, 
\eeaa
or
\beaa
&&\frac{rq}{\ov{q}}\nabc_3\left(\frac{r\ov{q}}{q}\frac{|q|^2}{r^2}\pmb\phi_{+2}^{(1)}\right)\\
&=& |q|^2\ov{q}^2\Bigg\{\nabc_3^2A+\left( \frac{3}{2}\trchb - \frac{1}{2}\frac{\atrchb^2}{\trchb} -2i\atrchb   +\frac{1}{2}\trchb - \frac{3}{2}\frac{\atrchb^2}{\trchb} -2i\atrchb\right)\nabc_3A\\
&&+\Bigg[\nabc_3\left(\frac{1}{2}\trchb - \frac{3}{2}\frac{\atrchb^2}{\trchb} -2i\atrchb\right)\\
&&+\left(\frac{3}{2}\trchb - \frac{1}{2}\frac{\atrchb^2}{\trchb} -2i\atrchb\right)\left(\frac{1}{2}\trchb -\frac{3}{2}\frac{\atrchb^2}{\trchb} -2i\atrchb\right)\Bigg]A\Bigg\}+r^4\Ga_b\nabc_3A+r^3\Ga_bA.
\eeaa

Next, we have, using the definition of $C_1$ in \eqref{eq:C1-C2-comparison-Ma},  
\beaa
\frac{3}{2}\trchb - \frac{1}{2}\frac{\atrchb^2}{\trchb} -2i\atrchb   +\frac{1}{2}\trchb - \frac{3}{2}\frac{\atrchb^2}{\trchb} -2i\atrchb &=& 2\trchb -\frac{2\atrchb^2}{\trchb} -4 i \atrchb =C_1
\eeaa
and hence
\beaa
\frac{rq}{\ov{q}}\nabc_3\left(\frac{r\ov{q}}{q}\frac{|q|^2}{r^2}\pmb\phi_{+2}^{(1)}\right) &=& |q|^2\ov{q}^2\Bigg\{\nabc_3^2A+C_1\nabc_3A +\Bigg[\nabc_3\left(\frac{1}{2}\trchb -\frac{3}{2}\frac{\atrchb^2}{\trchb} -2i\atrchb\right)\\
&&+\left(\frac{3}{2}\trchb - \frac{1}{2}\frac{\atrchb^2}{\trchb} -2i\atrchb\right)\left(\frac{1}{2}\trchb -\frac{3}{2}\frac{\atrchb^2}{\trchb} -2i\atrchb\right)\Bigg]A\Bigg\}\\
&&+r^4\Ga_b\nabc_3A+r^3\Ga_bA.
\eeaa

Next, recall from Proposition \ref{prop-nullstr:complex-conf} that we have the following null structure equation in complex form using conformal derivatives 
\beaa
\nabc_3\tr\Xb +\frac{1}{2}(\tr\Xb)^2 &=& \DDc\c\ov{\Xib}+\Xib\c\ov{\Hb}+\ov{\Xib}\c H-\frac{1}{2}\Xbh\c\ov{\Xbh},
\eeaa
which we rewrite as 
\beaa
\nabc_3\tr\Xb +\frac{1}{2}(\tr\Xb)^2 &=& r^{-1}\dk^{\leq 1}\Ga_b,
\eeaa
or, taking the real and imaginary part,
\beaa
\nabc_3\trchb +\frac{1}{2}(\trchb^2-\atrchb^2) &=& r^{-1}\dk^{\leq 1}\Ga_b,\\
\nabc_3\atrchb +\trchb\atrchb &=& r^{-1}\dk^{\leq 1}\Ga_b.
\eeaa
We infer
\beaa
\nabc_3\left(\frac{1}{2}\trchb -\frac{3}{2}\frac{\atrchb^2}{\trchb} -2i\atrchb\right) &=& -\frac{1}{4}(\trchb^2-\atrchb^2) +3\atrchb^2 \\
&& -\frac{3}{4}\left(\frac{\atrchb}{\trchb}\right)^2(\trchb^2-\atrchb^2) +2i\trchb\atrchb+r^{-1}\dk^{\leq 1}\Ga_b\\
&=& -\frac{1}{4}\trchb^2 +\frac{5}{2}\atrchb^2 +\frac{3}{4}\frac{\atrchb^4}{\trchb^2} +2i\trchb\atrchb+r^{-1}\dk^{\leq 1}\Ga_b
\eeaa
and hence
\begin{align*}
\frac{rq}{\ov{q}}\nabc_3\left(\frac{r\ov{q}}{q}\frac{|q|^2}{r^2}\pmb\phi_{+2}^{(1)}\right) =& |q|^2\ov{q}^2\Bigg\{\nabc_3^2A+C_1\nabc_3A +\Bigg[-\frac{1}{4}\trchb^2  +\frac{5}{2}\atrchb^2 +\frac{3}{4}\frac{\atrchb^4}{\trchb^2}  +2i\trchb\atrchb\\
&+\left(\frac{3}{2}\trchb - \frac{1}{2}\frac{\atrchb^2}{\trchb} -2i\atrchb\right)\left(\frac{1}{2}\trchb - \frac{3}{2}\frac{\atrchb^2}{\trchb} -2i\atrchb\right)\Bigg]A\Bigg\}\\
&+r^4\Ga_b\nabc_3A+r^3\dk^{\leq 1}(\Ga_b)A.
\end{align*}

Since, we have, using the definition of $C_2$ in \eqref{eq:C1-C2-comparison-Ma},  
\beaa
&&-\frac{1}{4}\trchb^2 +\frac{5}{2}\atrchb^2  +\frac{3}{4}\frac{\atrchb^4}{\trchb^2} +2i\trchb\atrchb\\
&&+\left(\frac{3}{2}\trchb - \frac{1}{2}\frac{\atrchb^2}{\trchb} -2i\atrchb\right)\left(\frac{1}{2}\trchb - \frac{3}{2}\frac{\atrchb^2}{\trchb} -2i\atrchb\right)\\
&=& \frac{1}{2}\trchb^2  -4\atrchb^2  +\frac{3}{2}\frac{\atrchb^4}{\trchb^2} +  i \left(-2\trchb\atrchb +4\frac{\atrchb^3}{\trchb}\right) = C_2
\eeaa
we deduce 
\beaa
\frac{rq}{\ov{q}}\nabc_3\left(\frac{r\ov{q}}{q}\frac{|q|^2}{r^2}\pmb\phi_{+2}^{(1)}\right) &=& |q|^2\ov{q}^2\Big\{\nabc_3^2A+C_1\nabc_3A +C_2A\Big\}+r^4\Ga_b\nabc_3A+r^3\dk^{\leq 1}(\Ga_b)A\\
&=& \pmb\phi_{+2}^{(2)} +r^4\Ga_b\nabc_3A+r^3\dk^{\leq 1}(\Ga_b)A
\eeaa
and hence, using also $\omb\in\Ga_b$ given our ingoing normalization, 
\beaa
\nab_3\left(\frac{r\ov{q}}{q}\frac{|q|^2}{r^2}\pmb\phi_{+2}^{(1)}\right) &=& \frac{\ov{q}}{rq}\pmb\phi_{+2}^{(2)} +r^3\Ga_b\nab_3A+r^2\dk^{\leq 1}(\Ga_b)A
\eeaa
which is the stated estimate for $p=1$. This concludes the proof of Proposition \ref{prop:transportequationphis=plusminus2p=0and1} in the case $s=+2$.

%%%%%%%%%%%%%%%%%%%%%%%%%%%%%%%%%%%%%%

\paragraph{\textit{Transport equations for $\pmb\phi_{-2}^{(p)}$, $p=0,1$}.}

%%%%%%%%%%%%%%%%%%%%%%%%%%%%%%%%%%%%%%%

Next, we prove Proposition \ref{prop:transportequationphis=plusminus2p=0and1} in the case $s=-2$, i.e., the fact that $\pmb\phi_{-2}^{(p)}$, $p=0,1$ satisfies \eqref{eq:transportequationsins=-2caseforp=0and1} \eqref{eq:transportequationsins=-2caseforp=0and1:RHSNT-2p=0and1}.

We start with the case $p=0$. In view of the definition of $\pmb\phi_{-2}^{(0)}$ in \eqref{eq:definitionofthephiminus2phierarchy:perturbationofKerr}, we have 
\beaa
\frac{r\ov{q}}{q}\frac{r^2}{|q|^2}\nabc_4\left(\frac{rq}{\ov{q}}\frac{|q|^4}{r^4}\pmb\phi_{-2}^{(0)}\right) &=& \frac{r^3}{q^2}\nabc_4\left(\frac{q^4}{r^3}\left(\frac{\De}{|q|^2}\right)^2\Ab\right)\\
&=& q^2\left(\frac{\De}{|q|^2}\right)^2\left(\nabc_4\Ab+\left(\frac{4}{q}e_4(q) -\frac{3}{r}e_4(r)\right)\Ab\right)\\
&&+ 2q^2\left(\frac{\De}{|q|^2}\right)\left[\pr_r\left(\frac{\De}{|q|^2}\right)e_4(r)+\pr_{\cos\th}\left(\frac{\De}{|q|^2}\right)e_4(\cos\th)\right]\Ab.
\eeaa  
Now, we use the following conformally invariant identities
\beaa
\frac{e_4(q)}{q}=\frac{1}{2}\tr X+\Ga_g, \qquad \frac{e_4(r)}{r}=\frac{1}{2}\trch+\frac{1}{2}\left(\frac{\atrch}{\trch}\right)_{\vartheta}\atrch+\Ga_g,
\eeaa
as well as $e_4(\cos\th)=\Ga_g$, and the following non conformally invariant identities, which hold true for an ingoing normalization,
\beaa
\om = -\frac{1}{2}\pr_r\left(\frac{\De}{|q|^2}\right)+\Ga_g, \qquad e_4(r)=\frac{\De}{|q|^2}+\Ga_g, 
\eeaa
which implies
\beaa
\frac{r\ov{q}}{q}\frac{r^2}{|q|^2}\nabc_4\left(\frac{rq}{\ov{q}}\frac{|q|^4}{r^4}\pmb\phi_{-2}^{(0)}\right) &=& q^2\left(\frac{\De}{|q|^2}\right)^2\left(\nabc_4\Ab+\left(2\tr X -\frac{3}{2}\trch-\frac{3}{2}\left(\frac{\atrch}{\trch}\right)_{\vartheta}\atrch+\Ga_g\right)\Ab\right)\\
&& -4\om q^2\left(\frac{\De}{|q|^2}\right)^2\Ab+r^2\Ga_g\Ab\\
&=& q^2\left(\frac{\De}{|q|^2}\right)^2\left(\nabc_4+\frac{1}{2}\trch -2i\atrch -\frac{3}{2}\left(\frac{\atrch}{\trch}\right)_{\vartheta}\atrch\right)\Ab\\
&& -4\om\frac{r\ov{q}}{q}\frac{r^2}{|q|^2}\left(\frac{rq}{\ov{q}}\frac{|q|^4}{r^4}\pmb\phi_{-2}^{(0)}\right)+r^2\Ga_g\Ab\\
&=& -4\om\frac{r\ov{q}}{q}\frac{r^2}{|q|^2}\left(\frac{rq}{\ov{q}}\frac{|q|^4}{r^4}\pmb\phi_{-2}^{(0)}\right)+ \left(\frac{\De}{|q|^2}\right)\pmb\phi_{-2}^{(1)}+r^2\Ga_g\c\Ga_b,
\eeaa
and hence, since $\nabc_4f = \nab_4f -4\om f$ if $f$ is $-2$-conformally invariant, we infer
\beaa
\nab_4\left(\frac{rq}{\ov{q}}\left(\frac{r^2}{|q|^2}\right)^{-2}\pmb\phi_{-2}^{(0)}\right) = \frac{q}{r\ov{q}}\left(\frac{r^2}{|q|^2}\right)^{-1}\frac{\De}{|q|^2}\pmb\phi_{-2}^{(1)}+r\Ga_g\c\Ga_b,
\eeaa
which is the stated estimate for $p=0$.

Next, we consider the case $p=1$. In view of the definition for $\pmb\phi_{-2}^{(1)}$ in \eqref{eq:definitionofthephiminus2phierarchy:perturbationofKerr}, we have
\begin{align*}
\frac{r\ov{q}}{q}\nabc_4\left(\frac{rq}{\ov{q}}\frac{|q|^2}{r^2}\pmb\phi_{-2}^{(1)}\right) =& \frac{\ov{q}}{q}r\nabc_4\left(r\frac{q}{\ov{q}}\frac{|q|^2q^2}{r^2}\frac{\De}{|q|^2}\left(\nabc_4 +\frac{1}{2}\trch  -\frac{3}{2}\left(\frac{\atrch}{\trch}\right)_{\vartheta}\atrch -2i\atrch\right)\Ab\right)\\
=& \frac{\ov{q}}{q}r\nabc_4\left(\frac{q^4}{r}\frac{\De}{|q|^2}\left(\nabc_4 +\frac{1}{2}\trch  -\frac{3}{2}\left(\frac{\atrch}{\trch}\right)_{\vartheta}\atrch -2i\atrch\right)\Ab\right)\\
=& |q|^2q^2\frac{\De}{|q|^2}\left(\nabc_4 +2\tr X -\frac{1}{2}\trch - \frac{1}{2}\left(\frac{\atrch}{\trch}\right)_{\vartheta}\atrch+\Ga_g\right)\\
&\times\left(\nabc_4 +\frac{1}{2}\trch  -\frac{3}{2}\left(\frac{\atrch}{\trch} -2i\atrch\right)_{\vartheta}\atrch\right)\Ab\\
& +|q|^2q^2 e_4\left(\frac{\De}{|q|^2}\right)\left(\nabc_4 +\frac{1}{2}\trch  -\frac{3}{2}\left(\frac{\atrch}{\trch}\right)_{\vartheta}\atrch -2i\atrch\right)\Ab.
\end{align*}
Now, using again $e_4(\cos\th)=\Ga_g$, and the following non conformally invariant identities, which hold true for an ingoing normalization,
\beaa
\om = -\frac{1}{2}\pr_r\left(\frac{\De}{|q|^2}\right)+\Ga_g, \qquad e_4(r)=\frac{\De}{|q|^2}+\Ga_g, 
\eeaa
we have
\beaa
&& |q|^2q^2 e_4\left(\frac{\De}{|q|^2}\right)\left(\nabc_4 +\frac{1}{2}\trch  -\frac{3}{2}\left(\frac{\atrch}{\trch}\right)_{\vartheta}\atrch -2i\atrch\right)\Ab\\
&=& -2\om\frac{\De}{|q|^2}|q|^2q^2\left(\nabc_4 +\frac{1}{2}\trch  -\frac{3}{2}\left(\frac{\atrch}{\trch}\right)_{\vartheta}\atrch -2i\atrch\right)\Ab +r\Ga_g\dk^{\leq 1}\Ab\\
&=& -2\om \frac{r\ov{q}}{q}\frac{rq}{\ov{q}}\frac{|q|^2}{r^2}\pmb\phi_{-2}^{(1)} +r\Ga_g\dk^{\leq 1}\Ab.
\eeaa
Thus, we infer
\beaa
&&\frac{r\ov{q}}{q}\nabc_4\left(\frac{rq}{\ov{q}}\frac{|q|^2}{r^2}\pmb\phi_{-2}^{(1)}\right) +2\om \frac{r\ov{q}}{q}\frac{rq}{\ov{q}}\frac{|q|^2}{r^2}\pmb\phi_{-2}^{(1)}\\
&=& |q|^2q^2\frac{\De}{|q|^2}\left(\nabc_4 +2\tr X -\frac{1}{2}\trch - \frac{1}{2}\left(\frac{\atrch}{\trch}\right)_{\vartheta}\atrch+\Ga_g\right)\\
&&\times\left(\nabc_4 +\frac{1}{2}\trch  -\frac{3}{2}\left(\frac{\atrch}{\trch} -2i\atrch\right)_{\vartheta}\atrch\right)\Ab  +r\Ga_g\dk^{\leq 1}\Ab,
\eeaa
or, since $\nabc_4f = \nab_4f -2\om f$ if $f$ is $-1$-conformally invariant, 
\beaa
\frac{r\ov{q}}{q}\nab_4\left(\frac{rq}{\ov{q}}\frac{|q|^2}{r^2}\pmb\phi_{-2}^{(1)}\right) &=& |q|^2q^2\frac{\De}{|q|^2}\left(\nabc_4 +2\tr X -\frac{1}{2}\trch - \frac{1}{2}\left(\frac{\atrch}{\trch}\right)_{\vartheta}\atrch+\Ga_g\right)\\
&&\times\left(\nabc_4 +\frac{1}{2}\trch  -\frac{3}{2}\left(\frac{\atrch}{\trch} -2i\atrch\right)_{\vartheta}\atrch\right)\Ab  +r\Ga_g\dk^{\leq 1}\Ab.
\eeaa
We deduce
\beaa
&&\frac{r\ov{q}}{q}\nab_4\left(\frac{rq}{\ov{q}}\frac{|q|^2}{r^2}\pmb\phi_{-2}^{(1)}\right)\\
&=& |q|^2q^2\frac{\De}{|q|^2}\left(\nabc_4 +\frac{3}{2}\trch - \frac{1}{2}\left(\frac{\atrch}{\trch}\right)_{\vartheta}\atrch -2i\atrch\right)\\
&&\times\left(\nabc_4 +\frac{1}{2}\trch  -\frac{3}{2}\left(\frac{\atrch}{\trch}\right)_{\vartheta}\atrch -2i\atrch\right)\Ab +r^4\Ga_g\left(\nabc_4\Ab+\frac{1}{r}\Ab\right)+r^2\Ga_g\Ab\\
&=& |q|^2q^2\frac{\De}{|q|^2}\Bigg\{\nabc_4^2\Ab+\left(2\trch  -2\left(\frac{\atrch}{\trch}\right)_{\vartheta}\atrch -4i\atrch\right)\nabc_4\Ab\\
&&+\Bigg[\left(\frac{3}{2}\trch - \frac{1}{2}\left(\frac{\atrch}{\trch}\right)_{\vartheta}\atrch -2i\atrch\right)\left(\frac{1}{2}\trch  -\frac{3}{2}\left(\frac{\atrch}{\trch}\right)_{\vartheta}\atrch -2i\atrch\right)\\
&&+\nabc_4\left(\frac{1}{2}\trch  -\frac{3}{2}\left(\frac{\atrch}{\trch}\right)_{\vartheta}\atrch -2i\atrch\right)\Bigg]\Ab\Bigg\} +r^2\Ga_g\dk^{\leq 1}\Ga_b,
\eeaa
where we have used in particular in the last estimate 
\beaa
\nabc_4\Ab+\frac{1}{r}\Ab &=& r^{-2}\dk^{\leq 1}\Ga_b.
\eeaa

Next, we have from Proposition \ref{prop-nullstr:complex-conf} 
\beaa
\nabc_4\tr X +\frac{1}{2}(\tr X)^2 &=& r^{-1}\dk^{\leq 1}\Xi+r^{-2}\Ga_g,
\eeaa
or, taking the real and imaginary part,
\beaa
\nabc_4\trch +\frac{1}{2}(\trch^2-\atrch^2) &=& r^{-1}\dk^{\leq 1}\Xi+r^{-2}\Ga_g,\\
\nabc_4\atrch +\trch\atrch &=& r^{-1}\dk^{\leq 1}\Xi+r^{-2}\Ga_g.
\eeaa
Using also
\beaa
\left(\frac{\atrch}{\trch}\right)_{\vartheta}\trch &=& \vartheta(r)\frac{a\cos\th}{r}\trch+(1-\vartheta(r))\atrch \\
&=& \atrch +\vartheta(r)\left(\frac{a\cos\th}{r}\trch - \atrch\right)\\
 &=& \atrch + \vartheta(r)\Gac,
\eeaa
we infer
\beaa
&&\frac{r\ov{q}}{q}\nab_4\left(\frac{rq}{\ov{q}}\frac{|q|^2}{r^2}\pmb\phi_{-2}^{(1)}\right)\\
&=& |q|^2q^2\frac{\De}{|q|^2}\Bigg\{\nabc_4^2\Ab+\left(2\trch  -2\left(\frac{\atrch}{\trch}\right)_{\vartheta}\atrch -4i\atrch\right)\nabc_4\Ab\\
&&+\Bigg[-\frac{1}{4}(\trch^2-\atrch^2)  +\frac{3}{2}\left(\frac{\atrch}{\trch}\right)_{\vartheta}\trch\atrch \\
&& +\frac{3}{2}\left(\frac{\atrch}{\trch}\right)_{\vartheta}\left(\frac{1}{2}\trch+\frac{1}{2}\left(\frac{\atrch}{\trch}\right)_{\vartheta}\atrch\right)\atrch +2i\trch\atrch+r^{-1}\dk^{\leq 1}\Xi+r^{-2}\Ga_g\\
&& +\frac{3}{4}\trch^2 - \frac{5}{2}\left(\frac{\atrch}{\trch}\right)_{\vartheta}\trch\atrch + \frac{3}{4}\left(\left(\frac{\atrch}{\trch}\right)_{\vartheta}\right)^2\atrch^2 -4\atrch^2\\
&&+i\left(-4\trch\atrch+4\left(\frac{\atrch}{\trch}\right)_{\vartheta}\atrch^2\right)\Bigg]\Ab\Bigg\}+r^2\Ga_g\dk^{\leq 1}\Ga_b\\
&=& |q|^2q^2\frac{\De}{|q|^2}\Bigg\{\nabc_4^2\Ab+\left(2\trch  -2\left(\frac{\atrch}{\trch}\right)_{\vartheta}\atrch -4i\atrch\right)\nabc_4\Ab\\
&&+\Bigg[\frac 1 2 \trch^2 -4(\atrchb)^2 +\frac 3 2\left(\left(\frac{\atrch}{\trch}\right)_{\vartheta}\right)^2\atrch^2\\
&& +i\left(-2\trch\atrch +4\left(\frac{\atrch}{\trch}\right)_{\vartheta}\atrch^2\right)\Bigg]\Ab\Bigg\} +r^3\dk^{\leq 1}(\Xi)\c\Ga_b+r^2\Ga_g\dk^{\leq 1}\Ga_b\\
&=& \frac{\De}{|q|^2}\pmb\phi_{-2}^{(2)}+r^3\dk^{\leq 1}(\Xi)\c\Ga_b+r^2\Ga_g\dk^{\leq 1}\Ga_b
\eeaa
and hence
\beaa
\nab_4\left(\frac{rq}{\ov{q}}\frac{|q|^2}{r^2}\pmb\phi_{-2}^{(1)}\right) &=& \frac{q}{r\ov{q}}\frac{\De}{|q|^2}\pmb\phi_{-2}^{(2)}+r^2\dk^{\leq 1}(\Xi)\c\Ga_b+r\Ga_g\dk^{\leq 1}\Ga_b
\eeaa
which is the stated estimate for $p=1$. This concludes the proof of Proposition \ref{prop:transportequationphis=plusminus2p=0and1}.

%%%%%%%%%%%%%%%%%%%%%%%%%%%%%%%%%%%%%%%%%%%%%%%%%%%%%%%%%%%%%%%%%%%%%%%%%%%%%%

\subsection{Wave equations for $\pmb\phi_{+2}^{(p)}$, $p=0,1,2$}
\lab{sec:proofofthm:derivationoftheTeukolskytensorialwavesystemfors=plusminus2:kerrpert:alternateformnullframeinsteadcoordvectorfield:cases=+2}

%%%%%%%%%%%%%%%%%%%%%%%%%%%%%%%%%%%%%%%%%%%%%%%%%%%%%%%%%%%%%%%%%%%%%%%%%%%%%%

In this section, we prove Theorem \ref{thm:derivationoftheTeukolskytensorialwavesystemfors=plusminus2:kerrpert:alternateformnullframeinsteadcoordvectorfield} in the case $s=+2$, i.e., the fact that $\pmb\phi_{+2}^{(p)}$, $p=0,1,2$, verifies the system of tensorial wave equations \eqref{eq:TensorialTeuSysandlinearterms:rescaleRHScontaine2:general:Kerrperturbation:alternateformnullframeinsteadcoordvectorfield}. Note that in order to exhibit the structure of the nonlinear correction terms $\widetilde{\N}_{W,+2}^{(p)}$, $p=0,1,2$, stated in \eqref{eq:schematicformofNpWsplus2}, we will use in particular the additional assumptions \eqref{eq:additional-conditions-Hc-Xi-gRW-eq:0} on the global null frame.

%%%%%%%%%%%%%%%%%%%%%%%%%%%%%%%%%%%%%%%%%%

\subsubsection{Derivation of the wave equation for $\pmb\phi_{+2}^{(0)}$ in perturbations of Kerr}

%%%%%%%%%%%%%%%%%%%%%%%%%%%%%%%%%%%%%%%%%%

Recall from Proposition 2.4.1 in \cite{Gu} that 
\bea\label{Teukolsky-equation-tens}
 \LL(A)&=& \err[\LL(A)]
 \eea
 where
\bea\label{Teukolsky-operator-ch5}
\begin{split}
\LL(A) &=-\nab_4\nab_3A+ \frac{1}{4}\DD\hot (\DDb\c A)+\left(- \frac 1 2 \tr X -2\ov{\tr X} -2\om\right)\nab_3A\\
&-\left(\frac{1}{2}\tr\Xb-4\omb\right)\nab_4A+\left( 4H+\Hb +\ov{\Hb} +2Z+2\ov{Z}\right)\c \nab A+ VA+  H   \hot (\ov{\Hb} \c A),
\end{split}
\eea
where the potential $V$ is given by 
\bea
V := -\ov{\tr X} \tr \Xb +2\omb\tr X -2\om\tr\Xb +8\omb(\om+\ov{\tr X}) +\DD\c\ov{Z} +2Z\c\ov{Z} +2\Re(Z\c\ov{\Hb}) +2\ov{P} +4e_4(\omb),
\eea
and with error term  expressed schematically
\bea\lab{eq:ErrLLATeukolskyKerrpert}
\err[\LL(A)]=r^{-1}  \dk^{\leq 1}\big( \Ga_g \c  B\big) + \nab_3\Xi  \c B {+\big( r^{-1}\dkb^{\leq1} \Pc+ r^{-3} \Ga_b \big)\c \Xi +(\Ga_g \c \Xi) \c \Bb}.
\eea

Next, Corollary \ref{corollary-wave-complex} and \eqref{Teukolsky-operator-ch5} yield 
\beaa
\begin{split}
\LL(A) &= \squared_2A  -\big(2\ov{\tr X} +4\om\big)\nab_3A +4\omb\nab_4A +\big(4H +2Z+2\ov{Z}\big)\c \nab A\\
& + \left(V +  \frac 1 4 \tr X \ov{\tr \Xb} + \frac 1 4 \tr \Xb  \ov{\tr X}  + 2\ov{P} +2(2\eta\c\etab- 3i\eta\wedge\etab)\right)A  +(\Ga_b \c \Ga_g)\c A,
\end{split}
\eeaa
where we also used Lemma \ref{dot-hot-complex} to rewrite $H   \hot (\ov{\Hb} \c A)$. Plugging in \eqref{Teukolsky-equation-tens}, we infer
\beaa
\squared_2A &=& \big(2\ov{\tr X} +4\om\big)\nab_3A - 4\omb\nab_4A -\big(4H +2Z+2\ov{Z}\big)\c \nab A +\widetilde{V}A +\err[\LL(A)]+(\Ga_b \c \Ga_g)\c A,
\eeaa
where 
\beaa
\widetilde{V} &=& - \left(V +  \frac 1 4 \tr X \ov{\tr \Xb} + \frac 1 4 \tr \Xb  \ov{\tr X}  + 2\ov{P} +2(2\eta\c\etab- 3i\eta\wedge\etab)\right)\\
&=&   \frac{3}{4}\ov{\tr X} \tr \Xb -  \frac 1 4 \tr X \ov{\tr \Xb}   - 4\ov{P}  -4e_4(\omb) -2\omb\tr X +2\om\tr\Xb -8\omb(\om+\ov{\tr X})\\
&& -\DD\c\ov{Z} -2Z\c\ov{Z} -2\Re(Z\c\ov{\Hb}) -2(2\eta\c\etab- 3i\eta\wedge\etab).
\eeaa

Next, using \eqref{eq:definitionofthephiplus2phierarchy:perturbationofKerr}, we have
\beaa
\squared_2\pmb\phi_{+2}^{(0)} &=& \squared_2\left(\frac{\ov{q}}{q}A\right)=\frac{\ov{q}}{q}\squared_2A+2\g^{\a\b}e_\a\left(\frac{\ov{q}}{q}\right)\Ddot_\b A+\square\left(\frac{\ov{q}}{q}\right)A\\
&=& \frac{\ov{q}}{q}\squared_2A -e_3\left(\frac{\ov{q}}{q}\right)\nab_4A -e_4\left(\frac{\ov{q}}{q}\right)\nab_3A+2\nab\left(\frac{\ov{q}}{q}\right)\c\nab A+\square\left(\frac{\ov{q}}{q}\right)A
\eeaa
which together with \eqref{eq:definitionofTandPhithataretheapproximateKillingvectorifeldinKerrpert} implies 
\beaa
\squared_2\pmb\phi_{+2}^{(0)} &=& \frac{\ov{q}}{q}\squared_2A -e_3\left(\frac{\ov{q}}{q}\right)\left(2\nab_{\T} -\frac{\Delta}{|q|^2}\nab_3 +2a\Re(\Jk)\c\nab\right)A\\
&& -e_4\left(\frac{\ov{q}}{q}\right)\nab_3A+2\nab\left(\frac{\ov{q}}{q}\right)\c\nab A+\square\left(\frac{\ov{q}}{q}\right)A.
\eeaa
Plugging in the above, this yields 
\beaa
\squared_2\pmb\phi_{+2}^{(0)} &=& \frac{\ov{q}}{q}\bigg\{\big(2\ov{\tr X} +4\om\big)\nab_3A - 4\omb\nab_4A -\big(4H +2Z+2\ov{Z}\big)\c \nab A +\widetilde{V}A\\ 
&& +\err[\LL(A)]+(\Ga_b \c \Ga_g)\c A\bigg\} -e_3\left(\frac{\ov{q}}{q}\right)\left(2\nab_{\T} -\frac{\Delta}{|q|^2}\nab_3 +2a\Re(\Jk)\c\nab\right)A\\
&& -e_4\left(\frac{\ov{q}}{q}\right)\nab_3A+2\nab\left(\frac{\ov{q}}{q}\right)\c\nab A+\square\left(\frac{\ov{q}}{q}\right)A.
\eeaa
or
\beaa
\squared_2\pmb\phi_{+2}^{(0)} +2e_3\left(\frac{\ov{q}}{q}\right)\nab_{\T}A &=& \left(\frac{\ov{q}}{q}\big(2\ov{\tr X} +4\om\big) +e_3\left(\frac{\ov{q}}{q}\right)\frac{\Delta}{|q|^2} -e_4\left(\frac{\ov{q}}{q}\right)\right)\nab_3A  - 4\frac{\ov{q}}{q}\omb\nab_4A \\
&& +\left(-\frac{\ov{q}}{q}\big(4H +2Z+2\ov{Z}\big) -e_3\left(\frac{\ov{q}}{q}\right)2a\Re(\Jk) +2\nab\left(\frac{\ov{q}}{q}\right)\right)\c \nab A \\
&&+\left(\frac{\ov{q}}{q}\widetilde{V}+\square\left(\frac{\ov{q}}{q}\right)\right)A +\frac{\ov{q}}{q}\err[\LL(A)]+(\Ga_b \c \Ga_g)\c A. 
\eeaa

Next, we have
\beaa
\frac{q}{\ov{q}}e_3\left(\frac{\ov{q}}{q}\right) &=& \frac{e_3(\ov{q})}{\ov{q}} - \frac{e_3(q)}{q} = \frac{1}{q} - \frac{1}{\ov{q}}+\Ga_b= -\frac{2ia\cos\th}{|q|^2}+\Ga_b 
\eeaa
so that 
\beaa
2e_3\left(\frac{\ov{q}}{q}\right)\nab_{\T}A &=& -\frac{4ia\cos\th}{|q|^2}\frac{\ov{q}}{q}\nab_{\T}A+\Ga_b\nab_{\T}A\\
&=& -\frac{4ia\cos\th}{|q|^2}\nab_{\T}\pmb\phi_{+2}^{(0)}+\Ga_b\nab_3A+r^{-1}\Ga_b\dk^{\leq 1}A.
\eeaa
Also, we have
\beaa
2\ov{\tr X} +4\om +\frac{q}{\ov{q}}e_3\left(\frac{\ov{q}}{q}\right)\frac{\Delta}{|q|^2} - \frac{q}{\ov{q}}e_4\left(\frac{\ov{q}}{q}\right) &=& \frac{4}{\ov{q}}\frac{\De}{|q|^2} -2\pr_r\left(\frac{\De}{|q|^2}\right)  -\frac{2\De}{|q|^2}\left(\frac{1}{\ov{q}}-\frac{1}{q}\right)+\Ga_b\\
&=& \frac{4r}{|q|^2}\frac{\De}{|q|^2} -2\pr_r\left(\frac{\De}{|q|^2}\right)  +\Ga_b,
\eeaa
\beaa
&& -\big(4H +2Z+2\ov{Z}\big) - \frac{q}{\ov{q}}e_3\left(\frac{\ov{q}}{q}\right)2a\Re(\Jk) +2\frac{q}{\ov{q}}\nab\left(\frac{\ov{q}}{q}\right)\\
&=& -\frac{6aq}{|q|^2}\Jk  -\frac{2a\ov{q}}{|q|^2}\ov{\Jk} +\left(\frac{1}{\ov{q}}-\frac{1}{q}\right)2a\Re(\Jk) +2\left(\frac{1}{\ov{q}} +\frac{1}{q}\right)ai\dual\Re(\Jk)+\Ga_b\\
&=& -\frac{8ar}{|q|^2}\Re(\Jk)  +\frac{8a^2\cos\th}{|q|^2}\dual\Re(\Jk)  +\Ga_b,
\eeaa
and
\beaa
\widetilde{V}+\frac{q}{\ov{q}}\square\left(\frac{\ov{q}}{q}\right) &=& -\frac{2}{r^2}+O(mr^{-3})+r^{-1}\dk^{\leq 1}\Ga_b.
\eeaa
Plugging in the above, we infer
\begin{align*}
\squared_2\pmb\phi_{+2}^{(0)} -\frac{4ia\cos\th}{|q|^2}\nab_{\T}\pmb\phi_{+2}^{(0)} =& \left(\frac{4r}{|q|^2}\frac{\De}{|q|^2} -2\pr_r\left(\frac{\De}{|q|^2}\right)\right)\frac{\ov{q}}{q}\nab_3A    +\left(-\frac{8ar}{|q|^2}\Re(\Jk)  +\frac{8a^2\cos\th}{|q|^2}\dual\Re(\Jk)\right)\c \frac{\ov{q}}{q}\nab A \\
&+\left(-\frac{2}{r^2}+O(mr^{-3})\right)\frac{\ov{q}}{q}A +\frac{\ov{q}}{q}\err[\LL(A)]+(\Ga_b \c \Ga_g)\c A\\
& +\Ga_b\nab_3A+r^{-1}\dk^{\leq 1}(\Ga_bA). 
\end{align*}
Since 
\beaa
\frac{\ov{q}}{q}\nab_3A &=& \frac{1}{|q|^2}\left(\pmb\phi_{+2}^{(1)} +r\pmb\phi_{+2}^{(0)}  +O(a)\pmb\phi_{+2}^{(0)} +r^2\Ga_bA\right),\\
\frac{\ov{q}}{q}\nab A &=& \nab\pmb\phi_{+2}^{(0)} +O(ar^{-2})\pmb\phi_{+2}^{(0)}+\Ga_gA,
\eeaa
we deduce
\begin{align*}
\nn\squared_2\pmb\phi_{+2}^{(0)} -\frac{4ia\cos\th}{|q|^2}\nab_{\T}\pmb\phi_{+2}^{(0)} =& \left(\frac{4r}{|q|^2}\frac{\De}{|q|^2} -2\pr_r\left(\frac{\De}{|q|^2}\right)\right)\frac{1}{|q|^2}\pmb\phi_{+2}^{(1)}  +\left(-\frac{8ar}{|q|^2}\Re(\Jk)  +\frac{8a^2\cos\th}{|q|^2}\dual\Re(\Jk)\right)\c\nab\pmb\phi_{+2}^{(0)}  \\
\nn&+\left(\frac{2}{r^2}+O(mr^{-3})\right)\pmb\phi_{+2}^{(0)} +\frac{\ov{q}}{q}\err[\LL(A)]+(\Ga_b \c \Ga_g)\c A\\
& +\Ga_b\nab_3A+r^{-1}\dk^{\leq 1}(\Ga_bA). 
\end{align*}
which we write in simpler form as follows, using also \eqref{eq:formofregularhorizontalvectorfieldwidetildemathcalXs:Kerrperturbation}, 
\beaa
\squared_2\pmb\phi_{+2}^{(0)} -\frac{4ia\cos\th}{|q|^2}\nab_{\T}\pmb\phi_{+2}^{(0)} -\frac{2}{|q|^2}\pmb\phi_{+2}^{(0)} = \left(\frac{4}{r^3}+O(mr^{-4})\right)\pmb\phi_{+2}^{(1)}    +O(mr^{-3})\nab_{\widetilde{\Xcal}_{+2}}^{\leq 1}\pmb\phi_{+2}^{(0)}   +\widetilde{\N}_{W,+2}^{(0)},
\eeaa
with all the coefficients on the RHS  being independent of coordinates $\tau$ and $\phi$, and with $\widetilde{\N}^{(0)}_{W,+2}$ being given by 
\beaa
\widetilde{\N}^{(0)}_{W,+2}=  \frac{\ov{q}}{q}\err[\LL(A)]+(\Ga_b \c \Ga_g)\c A +\Ga_b\nab_3A+r^{-1}\dk^{\leq 1}(\Ga_bA)
\eeaa
so that, in view of  \eqref{eq:ErrLLATeukolskyKerrpert}, we have
\beaa
   \bsplit
 \widetilde{\N}^{(0)}_{W,+2} =& r^{-1}  \dk^{\leq 1}\big( \Ga_g \c  B\big) + \nab_3\Xi  \c B {+\big( r^{-1}\dkb^{\leq1} \Pc+ r^{-3} \Ga_b \big)\c \Xi +(\Ga_g \c \Xi) \c \Bb}\\
 & +(\Ga_b \c \Ga_g)\c A +\Ga_b\nab_3A+r^{-1}\dk^{\leq 1}(\Ga_bA).
\end{split}
 \eeaa
 This concludes the proof of \eqref{eq:TensorialTeuSysandlinearterms:rescaleRHScontaine2:general:Kerrperturbation:alternateformnullframeinsteadcoordvectorfield} in the case $s=+2$ and $p=0$.

%%%%%%%%%%%%%%%%%%%%%%%%%%%%%%%%%%%%%%%%%%%%%%%%%%%

\subsubsection{Derivation of the wave equation for $\pmb\phi_{+2}^{(1)}$ in perturbations of Kerr}

%%%%%%%%%%%%%%%%%%%%%%%%%%%%%%%%%%%%%%%%%%%%%%%%%%%

Recall that we have obtained above 
\beaa
\nn\squared_2\pmb\phi_{+2}^{(0)} -\frac{4ia\cos\th}{|q|^2}\nab_{\T}\pmb\phi_{+2}^{(0)} &=& \left(\frac{4}{r^3}+O(mr^{-4})\right)\pmb\phi_{+2}^{(1)}  +\left(-\frac{8ar}{|q|^2}\Re(\Jk)  +\frac{8a^2\cos\th}{|q|^2}\dual\Re(\Jk)\right)\c\nab\pmb\phi_{+2}^{(0)}  \\
\nn&&+\left(\frac{2}{r^2}+O(mr^{-3})\right)\pmb\phi_{+2}^{(0)} +\widetilde{\N}^{(0)}_{W,+2}.  
\eeaa
Next, we introduce the scalar function $h_1$ given by 
\bea
h_1:=\frac{r\ov{q}}{q}\left(\frac{r^2}{|q|^2}\right)^{-2}=\frac{\ov{q}^3q}{r^3}
\eea
and we compute 
\beaa
\squared_2(h_1\pmb\phi_{+2}^{(0)}) &=& h_1\squared_2\pmb\phi_{+2}^{(0)} -e_3(h_1)\nab_4\pmb\phi_{+2}^{(0)} -e_4(h_1)\nab_3\pmb\phi_{+2}^{(0)} +2\nab(h_1)\c\nab\pmb\phi_{+2}^{(0)}+\square(h_1)\pmb\phi_{+2}^{(0)}\\
&=& h_1\squared_2\pmb\phi_{+2}^{(0)} -e_3(h_1)\left(2\nab_{\T} -\frac{\Delta}{|q|^2}\nab_3 +2a\Re(\Jk)\c\nab\right)\pmb\phi_{+2}^{(0)} -e_4(h_1)\nab_3\pmb\phi_{+2}^{(0)}\\
&& +2\nab(h_1)\c\nab\pmb\phi_{+2}^{(0)}+\square(h_1)\pmb\phi_{+2}^{(0)}   
\eeaa
where we used \eqref{eq:definitionofTandPhithataretheapproximateKillingvectorifeldinKerrpert} in the last equality.  We infer
\beaa
\squared_2(h_1\pmb\phi_{+2}^{(0)}) &=& h_1\squared_2\pmb\phi_{+2}^{(0)} -2e_3(h_1)\nab_{\T}\pmb\phi_{+2}^{(0)}  
+\left(\frac{\Delta}{|q|^2}e_3(h_1) -e_4(h_1)\right)\nab_3\pmb\phi_{+2}^{(0)}\\
&& +\Big(-2e_3(h_1)a\Re(\Jk) +2\nab(h_1)\Big)\c\nab\pmb\phi_{+2}^{(0)}+\square(h_1)\pmb\phi_{+2}^{(0)}.   
\eeaa
Also, since 
\beaa
&&\frac{e_3(h_1)}{h_1} = -\frac{1}{r}+O(ar^{-2})+\Ga_b, \qquad \frac{e_4(h_1)}{h_1} = \frac{1}{r}+O(mr^{-2})+r^{-1}\Ga_g,\\
&& \frac{1}{h_1}\nab(h_1) =O(ar^{-2}) +\Ga_g, \qquad \square(h_1)=\frac{2}{r}+O(mr^{-2})+\dk^{\leq 1}\Ga_b, \qquad \T(h_1)=r\Ga_b,
\eeaa
we obtain 
\beaa
\squared_2(h_1\pmb\phi_{+2}^{(0)}) &=& h_1\squared_2\pmb\phi_{+2}^{(0)} -2\frac{e_3(h_1)}{h_1}\nab_{\T}(h_1\pmb\phi_{+2}^{(0)})  
+\big(-2+O(mr^{-1})\big)\nab_3\pmb\phi_{+2}^{(0)}\\
&& +\left(-2\frac{e_3(h_1)}{h_1}a\Re(\Jk) +2\frac{\nab(h_1)}{h_1}\right)\c\nab(h_1\pmb\phi_{+2}^{(0)})+\left(\frac{2}{r}+O(mr^{-2})\right)\pmb\phi_{+2}^{(0)}\\
&&+r\Ga_b\c\nab_3A+\dk^{\leq 1}(\Ga_b)\c A. 
\eeaa
Since we have in view of \eqref{eq:definitionofthephiplus2phierarchy:perturbationofKerr}
\beaa
\frac{1}{|q|^2}\pmb\phi_{+2}^{(1)} &=& \frac{\ov{q}^2}{|q|^2}\left(\nab_3  - \frac{1}{r} +O(ar^{-2}) + \Ga_b\right)A\\
&=& \left(\nab_3  - \frac{1}{r} +O(ar^{-2}) + \Ga_b\right)\pmb\phi_{+2}^{(0)},
\eeaa
we deduce
\beaa
\squared_2(h_1\pmb\phi_{+2}^{(0)}) &=& h_1\squared_2\pmb\phi_{+2}^{(0)} -2\frac{e_3(h_1)}{h_1}\nab_{\T}(h_1\pmb\phi_{+2}^{(0)})  +\left( -\frac{2}{r^2}+O(mr^{-3})\right)\pmb\phi_{+2}^{(1)}\\
&& +\left(-2\frac{e_3(h_1)}{h_1}a\Re(\Jk) +2\frac{\nab(h_1)}{h_1}\right)\c\nab(h_1\pmb\phi_{+2}^{(0)})+O(mr^{-2})\pmb\phi_{+2}^{(0)}+r\Ga_b\c\nab_3A+\dk^{\leq 1}(\Ga_b)\c A. 
\eeaa
Plugging the above tensorial wave equation for $\pmb\phi_{+2}^{(0)}$, and using also the fact that $\T(h_1)=r\Ga_b$ and $\nab(h_1)=O(ar^{-1})+r\Ga_g$, we infer
\beaa
\squared_2(h_1\pmb\phi_{+2}^{(0)}) &=& \left(\frac{4ia\cos\th}{|q|^2} -2\frac{e_3(h_1)}{h_1}\right)\nab_{\T}(h_1\pmb\phi_{+2}^{(0)}) +\left(\frac{2}{r^2}+O(mr^{-3})\right)\pmb\phi_{+2}^{(1)}   \\
\nn&& +\left(\left(-\frac{8r}{|q|^2} -2\frac{e_3(h_1)}{h_1}\right)a\Re(\Jk) +\frac{8a^2\cos\th}{|q|^2}\dual\Re(\Jk)  +2\frac{\nab(h_1)}{h_1}\right)\c\nab(h_1\pmb\phi_{+2}^{(0)}) \\
&&  +\left(\frac{2}{r}+O(mr^{-2})\right)\pmb\phi_{+2}^{(0)}+h_1\widetilde{\N}^{(0)}_{W,+2} +r\Ga_b\c\nab_3A+\dk^{\leq 1}(\Ga_b)\c A. 
\eeaa

Next, we commute the above tensorial wave equation for $h_1\pmb\phi_{+2}^{(0)}$ with $\nab_3$. We obtain 
\beaa
\squared_2(\nab_3(h_1\pmb\phi_{+2}^{(0)})) &=& -[\nab_3, \squared_2](h_1\pmb\phi_{+2}^{(0)}) +I+J+K+\err
\eeaa
where 
\beaa
I &:=& \left(\frac{4ia\cos\th}{|q|^2} -2\frac{e_3(h_1)}{h_1}\right)\nab_{\T}(\nab_3(h_1\pmb\phi_{+2}^{(0)}))+ \left(\frac{2}{r^2}+O(mr^{-3})\right)\nab_3\pmb\phi_{+2}^{(1)}   \\
&& +\left(\left(-\frac{8r}{|q|^2} -2\frac{e_3(h_1)}{h_1}\right)a\Re(\Jk) +\frac{8a^2\cos\th}{|q|^2}\dual\Re(\Jk)  +2\frac{\nab(h_1)}{h_1}\right)\c\nab(\nab_3(h_1\pmb\phi_{+2}^{(0)}))\\
\nn&&+\left(\frac{2}{r}+O(mr^{-2})\right)\nab_3\pmb\phi_{+2}^{(0)},
\eeaa
\beaa
J &:=& \left(\frac{4ia\cos\th}{|q|^2} -2\frac{e_3(h_1)}{h_1}\right)[\nab_3, \nab_{\T}](h_1\pmb\phi_{+2}^{(0)})  \\
&& +\left(\left(-\frac{8r}{|q|^2} -2\frac{e_3(h_1)}{h_1}\right)a\Re(\Jk) +\frac{8a^2\cos\th}{|q|^2}\dual\Re(\Jk)  +2\frac{\nab(h_1)}{h_1}\right)\c[\nab_3, \nab](h_1\pmb\phi_{+2}^{(0)}),
\eeaa
\beaa
K &:=& \nab_3\left(\frac{4ia\cos\th}{|q|^2} -2\frac{e_3(h_1)}{h_1}\right)\nab_{\T}(h_1\pmb\phi_{+2}^{(0)})+ \nab_3\left(\frac{2}{r^2}+O(mr^{-3})\right)\pmb\phi_{+2}^{(1)}   \\
&& +\nab_3\left(\left(-\frac{8r}{|q|^2} -2\frac{e_3(h_1)}{h_1}\right)a\Re(\Jk) +\frac{8a^2\cos\th}{|q|^2}\dual\Re(\Jk)  +2\frac{\nab(h_1)}{h_1}\right)\c\nab(h_1\pmb\phi_{+2}^{(0)})\\
\nn&&+\nab_3\left(\frac{2}{r}+O(mr^{-2})\right)\pmb\phi_{+2}^{(0)},   
\eeaa
and 
\bea\lab{eq:structureoferr=errortermindecompositionofsquared_2nab3h1pmbphiplus20}
\err &:=& \nab_3(h_1\widetilde{\N}^{(0)}_{W,+2})+r\dk^{\leq 1}(\Ga_b\c\nab_3A)+\dk^{\leq 2}(\Ga_b\c A).
\eea

Next, we compute the terms $I$, $J$ and $K$ starting with $I$. Using Proposition \ref{prop:transportequationphis=plusminus2p=0and1} and the definition of $h_1$, we have
\beaa
\nab_3(h_1\pmb\phi_{+2}^{(0)}) &=& \frac{\ov{q}}{rq}\left(\frac{r^2}{|q|^2}\right)^{-1}\pmb\phi_{+2}^{(1)}+r\Ga_bA,\\
\nab_3\pmb\phi_{+2}^{(0)} &=& \left(\frac{1}{r^2}+O(ar^{-3})\right)\pmb\phi_{+2}^{(1)}+\left(\frac{1}{r}+O(ar^{-2})\right)\pmb\phi_{+2}^{(0)}+\Ga_bA,\\
\nab_3\pmb\phi_{+2}^{(1)} &=& \left(\frac{1}{r^2}+O(ar^{-3})\right)\pmb\phi_{+2}^{(2)}+\left(\frac{1}{r}+O(ar^{-2})\right)\pmb\phi_{+2}^{(1)}+r^2\Ga_b\nab_3A+r\dk^{\leq 1}(\Ga_b)A,
\eeaa
so that 
\beaa
I &=& \left(\frac{4ia\cos\th}{|q|^2} -2\frac{e_3(h_1)}{h_1}\right)\nab_{\T}\left(\frac{\ov{q}}{rq}\left(\frac{r^2}{|q|^2}\right)^{-1}\pmb\phi_{+2}^{(1)}\right) + \left(\frac{2}{r^4}+O(mr^{-5})\right)\pmb\phi_{+2}^{(2)}   \\
&& +\left(\left(-\frac{8r}{|q|^2} -2\frac{e_3(h_1)}{h_1}\right)a\Re(\Jk) +\frac{8a^2\cos\th}{|q|^2}\dual\Re(\Jk)  +2\frac{\nab(h_1)}{h_1}\right)\c\nab\left(\frac{\ov{q}}{rq}\left(\frac{r^2}{|q|^2}\right)^{-1}\pmb\phi_{+2}^{(1)}\right)\\
&& + \left(\frac{4}{r^3}+O(mr^{-4})\right)\pmb\phi_{+2}^{(1)}   +\left(\frac{2}{r^2}+O(mr^{-3})\right)\pmb\phi_{+2}^{(0)}+\dk^{\leq 1}(\Ga_b\c A),
\eeaa
and hence
\beaa
I &=& \frac{\ov{q}}{rq}{\left(\frac{r^2}{|q|^2}\right)^{-1}}\left(\frac{4ia\cos\th}{|q|^2} -2\frac{e_3(h_1)}{h_1}\right)\nab_{\T}\pmb\phi_{+2}^{(1)} + \left(\frac{2}{r^4}+O(mr^{-5})\right)\pmb\phi_{+2}^{(2)}   \\
&& +\frac{\ov{q}}{rq}\left(\frac{r^2}{|q|^2}\right)^{-1}\left(\left(-\frac{8r}{|q|^2} -2\frac{e_3(h_1)}{h_1}\right)a\Re(\Jk) +\frac{8a^2\cos\th}{|q|^2}\dual\Re(\Jk)  +2\frac{\nab(h_1)}{h_1}\right)\c\nab\pmb\phi_{+2}^{(1)}\\
&& + \left(\frac{4}{r^3}+O(mr^{-4})\right)\pmb\phi_{+2}^{(1)}   +\left(\frac{2}{r^2}+O(mr^{-3})\right)\pmb\phi_{+2}^{(0)}+\dk^{\leq 1}(\Ga_b\c A).
\eeaa

Next, to compute $J$, we use the commutator formula \eqref{commutator-nab-3-nab-a-U} which yields, using that $\eta=\ze+\Ga_b$ ingoing normalization, for $U\in\sk_2$, 
\beaa
\, [\nab_3, \nab] U &=& -\frac  1 2   \trchb\nab U -\frac 1 2 \atrchb\dual \nab U +O(ar^{-3})U+\Ga_b \c \dk^{\leq 1}U.
\eeaa
Also, we have in view of Lemma \ref{lemma:commutatorbetweenLieTLieZandnabnab4nab3} and Lemma \ref{lemma:basicpropertiesLiebTfasdiuhakdisug:chap9}, for $U\in\sk_2$,  
\beaa
[\nab_3, \Lieb_{\T}]U=\dk^{\leq 1}(\Ga_b\c U), \qquad \nab_{\T}U=\Lieb_{\T}U+\frac{4amr\cos\th}{|q|^4}\dual U+\Ga_bU,
\eeaa
so that 
\beaa
[\nab_3, \nab_{\T}]U=O(am r^{-4})U+\dk^{\leq 1}(\Ga_b\c U).
\eeaa
In view of the definition of $J$, we deduce 
\beaa
J &=& \left(\left(-\frac{8r}{|q|^2} -2\frac{e_3(h_1)}{h_1}\right)a\Re(\Jk) +\frac{8a^2\cos\th}{|q|^2}\dual\Re(\Jk)  +2\frac{\nab(h_1)}{h_1}\right)\\
&&\c\left(-\frac  1 2   \trchb\nab -\frac 1 2 \atrchb\dual \nab\right)(h_1\pmb\phi_{+2}^{(0)})+O(amr^{-4})\pmb\phi_{+2}^{(0)}+\dk^{\leq 1}(\Ga_b\c A),
\eeaa
which together with the fact that $\nab(h_1)=O(ar^{-1})+r\Ga_g$ implies
\beaa
J &=& -\frac  1 2   \trchb h_1\left(\left(-\frac{8r}{|q|^2} -2\frac{e_3(h_1)}{h_1}\right)a\Re(\Jk) +\frac{8a^2\cos\th}{|q|^2}\dual\Re(\Jk)  +2\frac{\nab(h_1)}{h_1}\right)\c\nab\pmb\phi_{+2}^{(0)}\\
&& -\frac 1 2 \atrchb h_1\left(-\left(-\frac{8r}{|q|^2} -2\frac{e_3(h_1)}{h_1}\right)a\dual\Re(\Jk) +\frac{8a^2\cos\th}{|q|^2}\Re(\Jk)  -2\frac{\dual\nab(h_1)}{h_1}\right)\c\nab\pmb\phi_{+2}^{(0)}\\
&&+O(amr^{-4})\pmb\phi_{+2}^{(0)}+\dk^{\leq 1}(\Ga_b\c A).
\eeaa

Next, we compute $K$. We have
\beaa
K &=& \nab_3\left(\frac{4ia\cos\th}{|q|^2} -2\frac{e_3(h_1)}{h_1}\right)\nab_{\T}(h_1\pmb\phi_{+2}^{(0)})+ \left(\frac{4}{r^3}+O(mr^{-4})\right)\pmb\phi_{+2}^{(1)}   \\
&& +\nab_3\left(\left(-\frac{8r}{|q|^2} -2\frac{e_3(h_1)}{h_1}\right)a\Re(\Jk) +\frac{8a^2\cos\th}{|q|^2}\dual\Re(\Jk)  +2\frac{\nab(h_1)}{h_1}\right)\c\nab(h_1\pmb\phi_{+2}^{(0)})\\
\nn&&+\left(\frac{2}{r^2}+O(mr^{-3})\right)\pmb\phi_{+2}^{(0)}+\Ga_b\c\nab_3A+r^{-1}\Ga_b\c A,   
\eeaa
or
\beaa
K &=& h_1\nab_3\left(\frac{4ia\cos\th}{|q|^2} -2\frac{e_3(h_1)}{h_1}\right)\nab_{\T}\pmb\phi_{+2}^{(0)}+ \left(\frac{4}{r^3}+O(mr^{-4})\right)\pmb\phi_{+2}^{(1)}   \\
&& +h_1\nab_3\left(\left(-\frac{8r}{|q|^2} -2\frac{e_3(h_1)}{h_1}\right)a\Re(\Jk) +\frac{8a^2\cos\th}{|q|^2}\dual\Re(\Jk)  +2\frac{\nab(h_1)}{h_1}\right)\c\nab\pmb\phi_{+2}^{(0)}\\
\nn&&+\left(\frac{2}{r^2}+O(mr^{-3})\right)\pmb\phi_{+2}^{(0)}+\Ga_b\c\nab_3A+r^{-1}\Ga_b\c A.
\eeaa

Next, recalling that 
\beaa
\squared_2(\nab_3(h_1\pmb\phi_{+2}^{(0)})) &=& -[\nab_3, \squared_2](h_1\pmb\phi_{+2}^{(0)}) +I+J+K+\err,
\eeaa
and relying on the above computations of $I$, $J$ and $K$, we obtain
\beaa
&&\squared_2(\nab_3(h_1\pmb\phi_{+2}^{(0)}))\\
&=& -[\nab_3, \squared_2](h_1\pmb\phi_{+2}^{(0)}) + \frac{\ov{q}}{rq}{\left(\frac{r^2}{|q|^2}\right)^{-1}}\left(\frac{4ia\cos\th}{|q|^2} -2\frac{e_3(h_1)}{h_1}\right)\nab_{\T}\pmb\phi_{+2}^{(1)} + \left(\frac{2}{r^4}+O(mr^{-5})\right)\pmb\phi_{+2}^{(2)}   \\
&& +\frac{\ov{q}}{rq}\left(\frac{r^2}{|q|^2}\right)^{-1}\left(\left(-\frac{8r}{|q|^2} -2\frac{e_3(h_1)}{h_1}\right)a\Re(\Jk) +\frac{8a^2\cos\th}{|q|^2}\dual\Re(\Jk)  +2\frac{\nab(h_1)}{h_1}\right)\c\nab\pmb\phi_{+2}^{(1)}\\
&&+ h_1\nab_3\left(\frac{4ia\cos\th}{|q|^2} -2\frac{e_3(h_1)}{h_1}\right)\nab_{\T}\pmb\phi_{+2}^{(0)}\\ 
&& +\bigg\{-\frac  1 2   \trchb h_1\left(\left(-\frac{8r}{|q|^2} -2\frac{e_3(h_1)}{h_1}\right)a\Re(\Jk) +\frac{8a^2\cos\th}{|q|^2}\dual\Re(\Jk)  +2\frac{\nab(h_1)}{h_1}\right)\\
&& -\frac 1 2 \atrchb h_1\left(-\left(-\frac{8r}{|q|^2} -2\frac{e_3(h_1)}{h_1}\right)a\dual\Re(\Jk) +\frac{8a^2\cos\th}{|q|^2}\Re(\Jk)  -2\frac{\dual\nab(h_1)}{h_1}\right)\\
&& +h_1\nab_3\left(\left(-\frac{8r}{|q|^2} -2\frac{e_3(h_1)}{h_1}\right)a\Re(\Jk) +\frac{8a^2\cos\th}{|q|^2}\dual\Re(\Jk)  +2\frac{\nab(h_1)}{h_1}\right)\bigg\}\c\nab\pmb\phi_{+2}^{(0)}\\
\nn&& + \left(\frac{8}{r^3}+O(mr^{-4})\right)\pmb\phi_{+2}^{(1)}  +\left(\frac{4}{r^2}+O(mr^{-3})\right)\pmb\phi_{+2}^{(0)}+\err+\dk^{\leq 1}(\Ga_b\c A).
\eeaa

Next, we compute $[\nab_3, \squared_2](h_1\pmb\phi_{+2}^{(0)})$. To this end, relying on Corollary \ref{corollary-wave-complex}, we have  
\beaa
\begin{split}
[\nab_3, \squared_2]\psi=& -[\nab_3, \nab_4] \nab_3 \psi +\frac 1 4  [\nab_3, \DD\hot]( \DDb \c \psi)+\frac 1 4  \DD\hot( [\nab_3, \DDb\c] \psi) +\nab_3\left(2\om -\frac 1 2 \tr X\right) \nab_3\psi\\
& - \frac 1 2 \tr\Xb[\nab_3, \nab_4]\psi - \frac 1 2 e_3(\tr\Xb)\nab_4\psi +2\etab \c[\nab_3, \nab]\psi +2\nab_3\etab \c\nab \psi \\
& +  \nab_3\left( - \frac 1 4 \tr X \ov{\tr \Xb}- \frac 1 4 \tr \Xb  \ov{\tr X}  - 2\ov{P}\right) \psi- 2i \nab_3\left(\eta \wedge \etab\right)  \psi+\dk^{\leq 1}\big(\Ga_b \c \Ga_g)\big) \c \psi,
\end{split}
\eeaa
and hence
\beaa
\begin{split}
[\nab_3, \squared_2]\psi=& -[\nab_3, \nab_4] \nab_3 \psi +\frac 1 4  [\nab_3, \DD\hot]( \DDb \c \psi)+\frac 1 4  \DD\hot( [\nab_3, \DDb\c] \psi) +\left(-\frac{1}{r^2}+O(mr^{-3})\right) \nab_3\psi\\
& - \frac 1 2 \tr\Xb[\nab_3, \nab_4]\psi - \frac 1 2 e_3(\tr\Xb)\nab_4\psi +2\etab \c[\nab_3, \nab]\psi +2\nab_3\etab \c\nab \psi \\
& +  \left(\frac{4}{r^3}+O(mr^{-4})\right)\psi  +\dk^{\leq 1}\Ga_g\c\nab_3\psi+r^{-1}\dk^{\leq 1}\Ga_g\c\psi+\dk^{\leq 1}\big(\Ga_b \c \Ga_g)\big) \c \psi.
\end{split}
\eeaa
Then, we have, in view of Lemma \ref{LEMMA:COMMUTATION-FORMULAS-1} and the fact that $H=Z+\Ga_b$ in ingoing normalization, the following commutation formulas
\beaa
\bsplit
\, [\nab_3, \mathcal{D}\hot]\DDb\c\psi =& - \frac 1 2 \tr \Xb( \mathcal{D}\hot \DDb\c\psi  + H \hot \DDb\c\psi )+ r^{-1}\Ga_b\dk^{\leq 2}\psi,\\
\, [\nab_3, \ov{\DD}\c]\psi  =& - \frac 1 2\ov{\tr\Xb}\, ( \ov{\DD} \c\psi  -  2 \ov{H} \c\psi )+\Ga_b\c\nab_3\psi  + r^{-1}\Ga_b \c  \dk^{\leq 1}\psi ,\\
\,[\nab_3, \nab_4]\nab_3\psi  =& O(mr^{-2})\nab_3^2\psi   + 2 (\eta -\etab)\c\nab\nab_3\psi  +O(ma r^{-4})\nab_3\psi \\
&+\Ga_g\c\dk^{\leq 1}\nab_3\psi+r^{-1}\Ga_b\c\dk^{\leq 2}\psi +\left(\Ga_b  \c \Ga_g \right)\dk^{\leq 1}\psi ,\\
\,[\nab_3, \nab_4]\psi  =& O(mr^{-2})\nab_3\psi   + 2 (\eta -\etab)\c\nab\psi  +O(ma r^{-4})\psi+\Ga_g\c\dk\psi +\left(\Ga_b  \c \Ga_g \right)\psi ,\\
\, [\nab_3, \nab]\psi =& -\frac  1 2   \trchb\, \nab\psi -\frac 1 2 \atrchb\, \dual \nab\psi +O(ar^{-3})\psi +\Ga_b\c\nab_3\psi +r^{-1}\Ga_b \c \dk^{\leq 1}\psi,
 \end{split}
\eeaa
which yields
\beaa
\begin{split}
[\nab_3, \squared_2]\psi=& O(mr^{-2})\nab_3^2\psi  -2 (\eta -\etab)\c\nab\nab_3\psi   - \frac{1}{8}\tr \Xb( \mathcal{D}\hot \DDb\c\psi  + H \hot \DDb\c\psi)\\
& -\frac{1}{8}\DD\hot\big(\ov{\tr\Xb}\, ( \ov{\DD} \c\psi  -  2 \ov{H} \c\psi)\big) - \tr\Xb(\eta -\etab)\c\nab\psi -\etab \c\left(\trchb\, \nab\psi +\atrchb\, \dual \nab\psi\right)\\
& +2\nab_3\etab \c\nab \psi  +\left(-\frac{1}{r^2}+O(mr^{-3})\right) \nab_3\psi - \frac 1 2 e_3(\tr\Xb)\nab_4\psi\\ 
& +  \left(\frac{4}{r^3}+O(mr^{-4})\right)\psi  +\dk^{\leq 1}(\Ga_g\c\nab_3\psi)+r^{-1}\dk^{\leq 2}(\Ga_b\c\psi)+\dk^{\leq 1}\big(\Ga_b \c \Ga_g)\c\psi\big).
\end{split}
\eeaa

Next, we have the following consequence of Lemma \ref{SIMPLIFICATION-ANGULAR}
\begin{align*}
\DD\hot\left(\ov{\tr\Xb}\, ( \ov{\DD} \c\psi -  2 \ov{H} \c\psi)\right) =& \ov{\tr\Xb}\DD\hot( \ov{\DD} \c\psi -  2 \ov{H} \c\psi) +\DD(\ov{\tr\Xb})\hot( \ov{\DD} \c\psi -  2 \ov{H} \c\psi)\\
=& \ov{\tr\Xb}\DD\hot(\ov{\DD} \c\psi) -4\ov{\tr\Xb}(\ov{H}\c\DD)\psi -4\ov{\tr\Xb}(\DD\c\ov{H})\psi +\DD(\ov{\tr\Xb})\hot( \ov{\DD} \c\psi -  2 \ov{H} \c\psi)\\
=& \ov{\tr\Xb}\DD\hot(\ov{\DD} \c\psi) -4\ov{\tr\Xb}(\ov{H}\c\DD)\psi +\DD(\ov{\tr\Xb})\hot(\ov{\DD} \c\psi)\\
&  +O(ar^{-4})\psi +r^{-2}\dk^{\leq 1}(\Ga_b)\c\psi
\end{align*}
which implies 
\beaa
\begin{split}
[\nab_3, \squared_2]\psi=& O(mr^{-2})\nab_3^2\psi  -2 (\eta -\etab)\c\nab\nab_3\psi   - \frac{1}{4}{\trchb}\mathcal{D}\hot(\DDb\c\psi)  
- \frac{1}{8}\Big(\DD(\ov{\tr\Xb})+\tr \Xb H\Big)\hot(\DDb\c\psi)\\
& +\frac{1}{2}\ov{\tr\Xb}(\ov{H}\c\DD)\psi  - \tr\Xb(\eta -\etab)\c\nab\psi -\etab \c\left(\trchb\, \nab\psi +\atrchb\, \dual \nab\psi\right)\\
& +2\nab_3\etab \c\nab \psi  +\left(-\frac{1}{r^2}+O(mr^{-3})\right) \nab_3\psi - \frac 1 2 e_3(\tr\Xb)\nab_4\psi\\ 
& +  \left(\frac{4}{r^3}+O(mr^{-4})\right)\psi  +\dk^{\leq 1}(\Ga_g\c\nab_3\psi)+r^{-1}\dk^{\leq 2}(\Ga_b\c\psi)+\dk^{\leq 1}\big(\Ga_b \c \Ga_g)\c\psi\big).
\end{split}
\eeaa
Since we have in view of Corollary \ref{corollary-wave-complex}  
\beaa
\begin{split}
-\frac 1 4  \DD\hot( \DDb \c \psi)=& -\squared_2 \psi -\nab_4 \nab_3 \psi +\left(-\frac{1}{r}+O(mr^{-2})\right)\nab_3\psi- \frac 1 2 \tr\Xb \nab_4\psi+2\etab \c\nab \psi \\
& +  \left(\frac{2}{r^2}+O(mr^{-3})\right)\psi +\Ga_g\c\nab_3\psi+r^{-1}\Ga_g\c\psi+(\Ga_b \c \Ga_g) \c \psi.
\end{split}
\eeaa
we infer
\beaa
\begin{split}
[\nab_3, \squared_2]\psi=& O(mr^{-2})\nab_3^2\psi  -2 (\eta -\etab)\c\nab\nab_3\psi - {\trchb}\squared_2 \psi - {\trchb}\nab_4 \nab_3 \psi - \frac 1 2{\trchb\tr\Xb}\nab_4\psi\\
&+2{\trchb}\,\etab\c\nab \psi  - \frac{1}{8}\Big(\DD(\ov{\tr\Xb})+\tr \Xb H\Big)\hot(\DDb\c\psi) +\frac{1}{2}\ov{\tr\Xb}(\ov{H}\c\DD)\psi  - \tr\Xb(\eta -\etab)\c\nab\psi\\
& -\etab \c\left(\trchb\, \nab\psi +\atrchb\, \dual \nab\psi\right) +2\nab_3\etab \c\nab \psi  +\left(\frac{1}{r^2}+O(mr^{-3})\right) \nab_3\psi - \frac 1 2 e_3(\tr\Xb)\nab_4\psi\\ 
& +  O(mr^{-4})\psi  +\dk^{\leq 1}(\Ga_g\c\nab_3\psi)+r^{-1}\dk^{\leq 2}(\Ga_b\c\psi)+\dk^{\leq 1}\big(\Ga_b \c \Ga_g\c\psi\big).
\end{split}
\eeaa
We now choose $\psi=h_1\pmb\phi_{+2}^{(0)}$ which, together with 
\beaa
\nab_3\psi &=& \nab_3(h_1\pmb\phi_{+2}^{(0)}) = \frac{\ov{q}}{rq}\left(\frac{r^2}{|q|^2}\right)^{-1}\pmb\phi_{+2}^{(1)}+r\Ga_bA,\\
\nab_3^2\psi &=& \nab_3^2(h_1\pmb\phi_{+2}^{(0)})=O(r^{-1})\nab_3\pmb\phi_{+2}^{(1)}+O(r^{-2})\pmb\phi_{+2}^{(1)}+r\dk^{\leq 1}(\Ga_bA)\\
&=& O(r^{-3})\pmb\phi_{+2}^{(2)}+O(r^{-2})\pmb\phi_{+2}^{(1)}+r\dk^{\leq 1}(\Ga_bA),\\
\nab\psi &=& \nab(h_1\pmb\phi_{+2}^{(0)})=h_1\nab\pmb\phi_{+2}^{(0)}+O(ar^{-1})\pmb\phi_{+2}^{(0)}+r\Ga_g\c A,\\
\nab\nab_3\psi &=& \nab\nab_3(h_1\pmb\phi_{+2}^{(0)})=\frac{\ov{q}}{rq}\left(\frac{r^2}{|q|^2}\right)^{-1}\nab\pmb\phi_{+2}^{(1)}+O(ar^{-3})\pmb\phi_{+2}^{(1)}+r\Ga_g\c\nab_3A+\dk^{\leq 1}(\Ga_b\c A),
\eeaa
yields 
\beaa
\begin{split}
[\nab_3, \squared_2](h_1\pmb\phi_{+2}^{(0)})=& O(mr^{-5})\pmb\phi_{+2}^{(2)}  - \frac{2\ov{q}}{rq}\left(\frac{r^2}{|q|^2}\right)^{-1}(\eta -\etab)\c\nab\pmb\phi_{+2}^{(1)} - {\trchb}\squared_2(h_1\pmb\phi_{+2}^{(0)})\\
& - {\trchb}\nab_4\left(\frac{\ov{q}}{rq}\left(\frac{r^2}{|q|^2}\right)^{-1}\pmb\phi_{+2}^{(1)}\right) - \frac 1 2\big(e_3(\tr\Xb)+{\trchb\tr\Xb}\big)\nab_4(h_1\pmb\phi_{+2}^{(0)})\\
&+2h_1{\trchb}\,\etab \c\nab\pmb\phi_{+2}^{(0)}  - \frac{1}{8}h_1\Big(\DD(\ov{\tr\Xb})+\tr \Xb H\Big)\hot(\DDb\c\pmb\phi_{+2}^{(0)}) +\frac{1}{2}h_1\ov{\tr\Xb}(\ov{H}\c\DD)\pmb\phi_{+2}^{(0)}\\
&  - h_1\tr\Xb(\eta -\etab)\c\nab\pmb\phi_{+2}^{(0)} -h_1\etab \c\left(\trchb\, \nab\pmb\phi_{+2}^{(0)} +\atrchb\, \dual \nab\pmb\phi_{+2}^{(0)}\right)\\
& +2h_1\nab_3\etab \c\nab\pmb\phi_{+2}^{(0)}  +\left(\frac{1}{r^3}+O(mr^{-4})\right)\pmb\phi_{+2}^{(1)}  +  O(mr^{-3})\pmb\phi_{+2}^{(0)}\\
&  +r\dk^{\leq 1}(\Ga_g\c\nab_3A)+\dk^{\leq 2}(\Ga_b\c A)+r\dk^{\leq 1}\big(\Ga_b \c \Ga_g\c A\big).
\end{split}
\eeaa

Next, recall that we have obtained
\beaa
\squared_2(h_1\pmb\phi_{+2}^{(0)}) &=& \left(\frac{4ia\cos\th}{|q|^2} -2\frac{e_3(h_1)}{h_1}\right)\nab_{\T}(h_1\pmb\phi_{+2}^{(0)}) +\left(\frac{2}{r^2}+O(mr^{-3})\right)\pmb\phi_{+2}^{(1)}   \\
\nn&& +\left(\left(-\frac{8r}{|q|^2} -2\frac{e_3(h_1)}{h_1}\right)a\Re(\Jk) +\frac{8a^2\cos\th}{|q|^2}\dual\Re(\Jk)  +2\frac{\nab(h_1)}{h_1}\right)\c\nab(h_1\pmb\phi_{+2}^{(0)}) \\
&&  +\left(\frac{2}{r}+O(mr^{-2})\right)\pmb\phi_{+2}^{(0)}+h_1\widetilde{\N}^{(0)}_{W,+2} +r\Ga_b\c\nab_3A+\dk^{\leq 1}(\Ga_b)\c A. 
\eeaa
Plugging in the above, we obtain 
\beaa
\begin{split}
[\nab_3, \squared_2](h_1\pmb\phi_{+2}^{(0)})=& O(mr^{-5})\pmb\phi_{+2}^{(2)}  - \frac{2\ov{q}}{rq}\left(\frac{r^2}{|q|^2}\right)^{-1}(\eta -\etab)\c\nab\pmb\phi_{+2}^{(1)}\\
& - {\trchb}\nab_4\left(\frac{\ov{q}}{rq}\left(\frac{r^2}{|q|^2}\right)^{-1}\pmb\phi_{+2}^{(1)}\right) - \frac 1 2\big(e_3(\tr\Xb)+{\trchb\tr\Xb}\big)\nab_4(h_1\pmb\phi_{+2}^{(0)})\\
& -h_1{\trchb}\left(\frac{4ia\cos\th}{|q|^2} -2\frac{e_3(h_1)}{h_1}\right)\nab_{\T}\pmb\phi_{+2}^{(0)}  \\
& -h_1{\trchb}\left(\left(-\frac{8r}{|q|^2} -2\frac{e_3(h_1)}{h_1}\right)a\Re(\Jk) +\frac{8a^2\cos\th}{|q|^2}\dual\Re(\Jk)  +2\frac{\nab(h_1)}{h_1}\right)\c\nab\pmb\phi_{+2}^{(0)} \\
&+2h_1{\trchb}\,\etab \c\nab\pmb\phi_{+2}^{(0)}  - \frac{1}{8}h_1\Big(\DD(\ov{\tr\Xb})+\tr \Xb H\Big)\hot(\DDb\c\pmb\phi_{+2}^{(0)}) +\frac{1}{2}h_1\ov{\tr\Xb}(\ov{H}\c\DD)\pmb\phi_{+2}^{(0)}\\
&  - h_1\tr\Xb(\eta -\etab)\c\nab\pmb\phi_{+2}^{(0)} -h_1\etab \c\left(\trchb\, \nab\pmb\phi_{+2}^{(0)} +\atrchb\, \dual \nab\pmb\phi_{+2}^{(0)}\right)\\
& +2h_1\nab_3\etab \c\nab\pmb\phi_{+2}^{(0)}  +\left(\frac{5}{r^3}+O(mr^{-4})\right)\pmb\phi_{+2}^{(1)}  +  \left(\frac{4}{r^2}+O(mr^{-3})\right)\pmb\phi_{+2}^{(0)}\\
& +O(1)\widetilde{\N}^{(0)}_{W,+2} +r\dk^{\leq 1}(\Ga_g\c\nab_3A)+\dk^{\leq 2}(\Ga_b\c A)+r\dk^{\leq 1}\big(\Ga_b \c \Ga_g\c A\big).
\end{split}
\eeaa
Now, using in particular Proposition \ref{prop:transportequationphis=plusminus2p=0and1}, as well as $e_3(\tr\Xb)=-\frac{1}{2}(\tr\Xb)^2+r^{-1}\dk^{\leq 1}\Ga_b$ in view of Proposition \ref{prop-nullstr:complex}, we have
\beaa
&& - {\trchb}\nab_4\left(\frac{\ov{q}}{rq}\left(\frac{r^2}{|q|^2}\right)^{-1}\pmb\phi_{+2}^{(1)}\right) - \frac 1 2\big(e_3(\tr\Xb)+{\trchb\tr\Xb}\big)\nab_4(h_1\pmb\phi_{+2}^{(0)})\\
&=& - {\trchb}\left(2\nab_T-\frac{\De}{|q|^2}\nab_3+2a\Re(\Jk)\c\nab\right)\left(\frac{\ov{q}}{rq}\left(\frac{r^2}{|q|^2}\right)^{-1}\pmb\phi_{+2}^{(1)}\right)\\
&& - \frac{1}{4}{|\tr\Xb|^2}\left(2\nab_T-\frac{\De}{|q|^2}\nab_3+2a\Re(\Jk)\c\nab\right)(h_1\pmb\phi_{+2}^{(0)})+r^{-1}\dk^{\leq 1}(\Ga_b)\dk A\\
&=& -2\frac{\ov{q}}{rq}\left(\frac{r^2}{|q|^2}\right)^{-1}{\trchb}\nab_T\pmb\phi_{+2}^{(1)}+\left(-\frac{2}{r^4}+O(mr^{-5})\right)\pmb\phi_{+2}^{(2)} - 2a\frac{\ov{q}}{rq}\left(\frac{r^2}{|q|^2}\right)^{-1}{\trchb}\Re(\Jk)\c\nab\pmb\phi_{+2}^{(1)}\\
&& +\left(-\frac{3}{r^3}+O(mr^{-4})\right)\pmb\phi_{+2}^{(1)}  - \frac{1}{2}h_1{|\tr\Xb|^2}\nab_T\pmb\phi_{+2}^{(0)} - \frac{1}{2}{|\tr\Xb|^2}ah_1\Re(\Jk)\c\nab\pmb\phi_{+2}^{(0)}+O(a^2r^{-4})\pmb\phi_{+2}^{(0)}\\
&&+r\dk^{\leq 1}(\Ga_b)\nab_3A+\dk^{\leq 1}(\Ga_b)\dk^{\leq 1}A
\eeaa
and hence
\beaa
\begin{split}
[\nab_3, \squared_2](h_1\pmb\phi_{+2}^{(0)})=& \left(-\frac{2}{r^4}+O(mr^{-5})\right)\pmb\phi_{+2}^{(2)}  - \frac{2\ov{q}}{rq}\left(\frac{r^2}{|q|^2}\right)^{-1}(\eta -\etab)\c\nab\pmb\phi_{+2}^{(1)}\\
&  -2\frac{\ov{q}}{rq}\left(\frac{r^2}{|q|^2}\right)^{-1}{\trchb}\nab_T\pmb\phi_{+2}^{(1)} - 2a\frac{\ov{q}}{rq}\left(\frac{r^2}{|q|^2}\right)^{-1}{\trchb}\Re(\Jk)\c\nab\pmb\phi_{+2}^{(1)}\\
& -h_1\left({\trchb}\left(\frac{4ia\cos\th}{|q|^2} -2\frac{e_3(h_1)}{h_1}\right) +\frac{1}{2}{|\tr\Xb|^2}\right)\nab_{\T}\pmb\phi_{+2}^{(0)}  - \frac{1}{2}{|\tr\Xb|^2}ah_1\Re(\Jk)\c\nab\pmb\phi_{+2}^{(0)} \\
& -h_1{\trchb}\left(\left(-\frac{8r}{|q|^2} -2\frac{e_3(h_1)}{h_1}\right)a\Re(\Jk) +\frac{8a^2\cos\th}{|q|^2}\dual\Re(\Jk)  +2\frac{\nab(h_1)}{h_1}\right)\c\nab\pmb\phi_{+2}^{(0)} \\
&+2h_1{\trchb}\,\etab \c\nab\pmb\phi_{+2}^{(0)}  - \frac{1}{8}h_1\Big(\DD(\ov{\tr\Xb})+\tr \Xb H\Big)\hot(\DDb\c\pmb\phi_{+2}^{(0)}) +\frac{1}{2}h_1\ov{\tr\Xb}(\ov{H}\c\DD)\pmb\phi_{+2}^{(0)}\\
&  - h_1\tr\Xb(\eta -\etab)\c\nab\pmb\phi_{+2}^{(0)} -h_1\etab \c\left(\trchb\, \nab\pmb\phi_{+2}^{(0)} +\atrchb\, \dual \nab\pmb\phi_{+2}^{(0)}\right)\\
& +2h_1\nab_3\etab \c\nab\pmb\phi_{+2}^{(0)}  +\left(\frac{2}{r^3}+O(mr^{-4})\right)\pmb\phi_{+2}^{(1)}  +  \left(\frac{4}{r^2}+O(mr^{-3})\right)\pmb\phi_{+2}^{(0)}\\
& +O(1)\widetilde{\N}^{(0)}_{W,+2} +r\dk^{\leq 1}(\Ga_b\c\nab_3A)+\dk^{\leq 2}(\Ga_b\c A)+r\dk^{\leq 1}\big(\Ga_b \c \Ga_g\c A\big).
\end{split}
\eeaa
Now, recalling from the above that we have obtained 
\beaa
&&\squared_2(\nab_3(h_1\pmb\phi_{+2}^{(0)}))\\
&=& -[\nab_3, \squared_2](h_1\pmb\phi_{+2}^{(0)}) + \frac{\ov{q}}{rq}{\left(\frac{r^2}{|q|^2}\right)^{-1}}\left(\frac{4ia\cos\th}{|q|^2} -2\frac{e_3(h_1)}{h_1}\right)\nab_{\T}\pmb\phi_{+2}^{(1)} + \left(\frac{2}{r^4}+O(mr^{-5})\right)\pmb\phi_{+2}^{(2)}   \\
&& +\frac{\ov{q}}{rq}\left(\frac{r^2}{|q|^2}\right)^{-1}\left(\left(-\frac{8r}{|q|^2} -2\frac{e_3(h_1)}{h_1}\right)a\Re(\Jk) +\frac{8a^2\cos\th}{|q|^2}\dual\Re(\Jk)  +2\frac{\nab(h_1)}{h_1}\right)\c\nab\pmb\phi_{+2}^{(1)}\\
&&+ h_1\nab_3\left(\frac{4ia\cos\th}{|q|^2} -2\frac{e_3(h_1)}{h_1}\right)\nab_{\T}\pmb\phi_{+2}^{(0)}\\ 
&& +\bigg\{-\frac  1 2   \trchb h_1\left(\left(-\frac{8r}{|q|^2} -2\frac{e_3(h_1)}{h_1}\right)a\Re(\Jk) +\frac{8a^2\cos\th}{|q|^2}\dual\Re(\Jk)  +2\frac{\nab(h_1)}{h_1}\right)\\
&& -\frac 1 2 \atrchb h_1\left(-\left(-\frac{8r}{|q|^2} -2\frac{e_3(h_1)}{h_1}\right)a\dual\Re(\Jk) +\frac{8a^2\cos\th}{|q|^2}\Re(\Jk)  -2\frac{\dual\nab(h_1)}{h_1}\right)\\
&& +h_1\nab_3\left(\left(-\frac{8r}{|q|^2} -2\frac{e_3(h_1)}{h_1}\right)a\Re(\Jk) +\frac{8a^2\cos\th}{|q|^2}\dual\Re(\Jk)  +2\frac{\nab(h_1)}{h_1}\right)\bigg\}\c\nab\pmb\phi_{+2}^{(0)}\\
\nn&& + \left(\frac{8}{r^3}+O(mr^{-4})\right)\pmb\phi_{+2}^{(1)}  +\left(\frac{4}{r^2}+O(mr^{-3})\right)\pmb\phi_{+2}^{(0)}+\err+\dk^{\leq 1}(\Ga_b\c A),
\eeaa
we deduce 
\beaa
&&\squared_2(\nab_3(h_1\pmb\phi_{+2}^{(0)}))\\
&=&  \frac{\ov{q}}{rq}{\left(\frac{r^2}{|q|^2}\right)^{-1}}\left(\frac{4ia\cos\th}{|q|^2} -2\frac{e_3(h_1)}{h_1}+2{\trchb}\right)\nab_{\T}\pmb\phi_{+2}^{(1)}  + \left(\frac{4}{r^4}+O(mr^{-5})\right)\pmb\phi_{+2}^{(2)}   \\
&& +\frac{\ov{q}}{rq}\left(\frac{r^2}{|q|^2}\right)^{-1}\bigg\{\left(-\frac{8r}{|q|^2} -2\frac{e_3(h_1)}{h_1}\right)a\Re(\Jk) +\frac{8a^2\cos\th}{|q|^2}\dual\Re(\Jk)  +2\frac{\nab(h_1)}{h_1}\\
&&  + 2(\eta -\etab) + 2a{\trchb}\Re(\Jk)\bigg\}\c\nab\pmb\phi_{+2}^{(1)} + h_1\bigg\{\nab_3\left(\frac{4ia\cos\th}{|q|^2} -2\frac{e_3(h_1)}{h_1}\right)\\
&& +{\trchb}\left(\frac{4ia\cos\th}{|q|^2} -2\frac{e_3(h_1)}{h_1}\right) +\frac{1}{2}{|\tr\Xb|^2}\bigg\}\nab_{\T}\pmb\phi_{+2}^{(0)}\\ 
&& +h_1\bigg\{{\frac  1 2}\trchb\left(\left(-\frac{8r}{|q|^2} -2\frac{e_3(h_1)}{h_1}\right)a\Re(\Jk) +\frac{8a^2\cos\th}{|q|^2}\dual\Re(\Jk)  +2\frac{\nab(h_1)}{h_1}\right)\\
&& -\frac 1 2 \atrchb\left(-\left(-\frac{8r}{|q|^2} -2\frac{e_3(h_1)}{h_1}\right)a\dual\Re(\Jk) +\frac{8a^2\cos\th}{|q|^2}\Re(\Jk)  -2\frac{\dual\nab(h_1)}{h_1}\right)\\
&& +\nab_3\left(\left(-\frac{8r}{|q|^2} -2\frac{e_3(h_1)}{h_1}\right)a\Re(\Jk) +\frac{8a^2\cos\th}{|q|^2}\dual\Re(\Jk)  +2\frac{\nab(h_1)}{h_1}\right)\\
&&  + \frac{1}{2}{|\tr\Xb|^2}a\Re(\Jk) -2{\trchb}\,\etab + \tr\Xb(\eta -\etab)  -2\nab_3\etab\bigg\}\c\nab\pmb\phi_{+2}^{(0)} + h_1\etab \c\left(\trchb\, \nab\pmb\phi_{+2}^{(0)} +\atrchb\, \dual \nab\pmb\phi_{+2}^{(0)}\right)\\
&&+ \frac{1}{8}h_1\Big(\DD(\ov{\tr\Xb})+\tr \Xb H\Big)\hot(\DDb\c\pmb\phi_{+2}^{(0)}) -\frac{1}{2}h_1\ov{\tr\Xb}(\ov{H}\c\DD)\pmb\phi_{+2}^{(0)}  + \left(\frac{6}{r^3}+O(mr^{-4})\right)\pmb\phi_{+2}^{(1)} \\
&& +O(mr^{-3})\pmb\phi_{+2}^{(0)} +\err +O(1)\widetilde{\N}^{(0)}_{W,+2} +r\dk^{\leq 1}(\Ga_b\c\nab_3A)+\dk^{\leq 2}(\Ga_b\c A)+r\dk^{\leq 1}\big(\Ga_b \c \Ga_g\c A\big).
\eeaa

Next, using in particular \eqref{Leib-eq-DDb} as well as the identity
\beaa
\ov{H}\c\DD=(\eta-i\dual\eta)\c(\nab+i\dual\nab)=2(\eta-i\dual\eta)\c\nab=2\ov{H}\c\nab,
\eeaa
we have 
\beaa
&& h_1\etab \c\left(\trchb\, \nab\pmb\phi_{+2}^{(0)} +\atrchb\, \dual \nab\pmb\phi_{+2}^{(0)}\right) + \frac{1}{8}h_1\Big(\DD(\ov{\tr\Xb})+\tr \Xb H\Big)\hot(\DDb\c\pmb\phi_{+2}^{(0)}) -\frac{1}{2}h_1\ov{\tr\Xb}(\ov{H}\c\DD)\pmb\phi_{+2}^{(0)}\\
&=& h_1\bigg\{\trchb\, \etab - \atrchb\dual\etab + \frac{1}{2}\Big(\DD(\ov{\tr\Xb})+\tr \Xb H\Big) -\ov{\tr\Xb}\,\ov{H}\bigg\}\c\nab\pmb\phi_{+2}^{(0)},
\eeaa
and, plugging in the above, we infer
\beaa
\squared_2(\nab_3(h_1\pmb\phi_{+2}^{(0)})) &=&  {\frac{\ov{q}}{rq}\left(\frac{r^2}{|q|^2}\right)^{-1}\left(\frac{4ia\cos\th}{|q|^2} -2\frac{e_3(h_1)}{h_1}+2\trchb\right)}\nab_{\T}\pmb\phi_{+2}^{(1)}  + \left(\frac{4}{r^4}+O(mr^{-5})\right)\pmb\phi_{+2}^{(2)}   \\
&& +\frac{\ov{q}}{rq}\left(\frac{r^2}{|q|^2}\right)^{-1}\bigg\{\left(-\frac{8r}{|q|^2} -2\frac{e_3(h_1)}{h_1}\right)a\Re(\Jk) +\frac{8a^2\cos\th}{|q|^2}\dual\Re(\Jk)  +2\frac{\nab(h_1)}{h_1}\\
&&  + 2(\eta -\etab) + 2a{\trchb}\Re(\Jk)\bigg\}\c\nab\pmb\phi_{+2}^{(1)} + h_1\bigg\{\nab_3\left(\frac{4ia\cos\th}{|q|^2} -2\frac{e_3(h_1)}{h_1}\right)\\
&& +{\trchb}\left(\frac{4ia\cos\th}{|q|^2} -2\frac{e_3(h_1)}{h_1}\right) {+\frac{1}{2}|\tr\Xb|^2}\bigg\}\nab_{\T}\pmb\phi_{+2}^{(0)}\\ 
&& +h_1\bigg\{{\frac  1 2}\trchb\left(\left(-\frac{8r}{|q|^2} -2\frac{e_3(h_1)}{h_1}\right)a\Re(\Jk) +\frac{8a^2\cos\th}{|q|^2}\dual\Re(\Jk)  +2\frac{\nab(h_1)}{h_1}\right)\\
&& -\frac 1 2 \atrchb\left(-\left(-\frac{8r}{|q|^2} -2\frac{e_3(h_1)}{h_1}\right)a\dual\Re(\Jk) +\frac{8a^2\cos\th}{|q|^2}\Re(\Jk)  -2\frac{\dual\nab(h_1)}{h_1}\right)\\
&& +\nab_3\left(\left(-\frac{8r}{|q|^2} -2\frac{e_3(h_1)}{h_1}\right)a\Re(\Jk) +\frac{8a^2\cos\th}{|q|^2}\dual\Re(\Jk)  +2\frac{\nab(h_1)}{h_1}\right)  + \frac{1}{2}{|\tr\Xb|^2}a\Re(\Jk)\\
&&-2{\trchb}\,\etab + \tr\Xb(\eta -\etab)  -2\nab_3\etab +\trchb\, \etab - \atrchb\dual\etab + \frac{1}{2}\Big(\DD(\ov{\tr\Xb})+\tr \Xb H\Big)\\
&& -\ov{\tr\Xb}\,\ov{H}\bigg\}\c\nab\pmb\phi_{+2}^{(0)}  + \left(\frac{6}{r^3}+O(mr^{-4})\right)\pmb\phi_{+2}^{(1)}  +O(mr^{-3})\pmb\phi_{+2}^{(0)}\\
&& +\err +O(1)\widetilde{\N}^{(0)}_{W,+2} +r\dk^{\leq 1}(\Ga_b\c\nab_3A)+\dk^{\leq 2}(\Ga_b\c A)+r\dk^{\leq 1}\big(\Ga_b \c \Ga_g\c A\big).
\eeaa

Next, we compute 
\beaa
&& h_1\bigg\{\nab_3\left(\frac{4ia\cos\th}{|q|^2} -2\frac{e_3(h_1)}{h_1}\right) +{\trchb}\left(\frac{4ia\cos\th}{|q|^2} -2\frac{e_3(h_1)}{h_1}\right) {+\frac{1}{2}|\tr\Xb|^2}\bigg\}\\
&=& \frac{|q|^2\ov{q}^2}{r^3}\bigg\{\nab_3\left(\frac{4ia\cos\th}{|q|^2} -2\left(\frac{3}{r}-\frac{3}{\ov{q}}-\frac{1}{q}\right)\right) -{\frac{2r}{|q|^2}}\left(\frac{4ia\cos\th}{|q|^2} -2\left(\frac{3}{r}-\frac{3}{\ov{q}}-\frac{1}{q}\right)\right) {+\frac{2}{|q|^2}}\bigg\}+r\dk^{\leq 1}\Ga_b\\
&=& \frac{|q|^2\ov{q}^2}{r^3}\bigg\{\frac{8ira\cos\th}{|q|^4} +2\left(-\frac{3}{r^2}+\frac{3}{\ov{q}^2}+\frac{1}{q^2}\right) - {\frac{2r}{|q|^2}}\left(\frac{4ia\cos\th}{|q|^2} -\frac{6}{r}+\frac{6}{\ov{q}}+\frac{2}{q}\right) {+\frac{2}{|q|^2}}\bigg\}+r\dk^{\leq 1}\Ga_b\\
&=& {-\frac{6a^2(\cos\th)^2}{r^2}\frac{(\ov{q})^2}{r^3}}+r\dk^{\leq 1}\Ga_b.
\eeaa
Also, in view of 
\beaa
&&\nab(q)=-ia\dual\Re(\Jk) +r\Ga_g, \qquad \nab(\ov{q})=ia\dual\Re(\Jk) +r\Ga_g, \\
&&\frac{e_3(h_1)}{h_1}=\frac{3}{r}-\frac{3}{\ov{q}}-\frac{1}{q}+\Ga_b, \qquad \frac{\nab(h_1)}{h_1}=\frac{3ia}{\ov{q}}\dual\Re(\Jk)-\frac{ia}{q}\dual\Re(\Jk)+\Ga_g,\\
&&\nab_3\Re(\Jk)=\frac{r}{|q|^2}\Re(\Jk) - \frac{a\cos\th}{|q|^2}\dual\Re(\Jk)+r^{-1}\Ga_g, \\ 
&&\nab_3\dual\Re(\Jk)=\frac{a\cos\th}{|q|^2}\Re(\Jk) + \frac{r}{|q|^2}\dual\Re(\Jk)+r^{-1}\Ga_g,\\
&&\eta = \frac{ar}{|q|^2}\Re(\Jk) - \frac{a^2\cos\th}{|q|^2}\dual\Re(\Jk)+\Ga_b,\quad \etab =-\frac{ar}{|q|^2}\Re(\Jk) - \frac{a^2\cos\th}{|q|^2}\dual\Re(\Jk)+\Ga_g,\\
&&\nab_3\etab = -\frac{1}{2}\trchb(\etab-\eta)+\frac{1}{2}\atrchb(\dual\etab-\dual\eta)+r^{-1}\dk^{\leq 1}\Ga_b,\\
&&\DD(\ov{\tr\Xb})=\ov{\tr\Xb}H+2i\Im(\tr\Xb)\Hb+r^{-1}\dk^{\leq 1}\Ga_b,
\eeaa
we have 
\beaa
&&h_1\bigg\{{\frac  1 2}\trchb\left(\left(-\frac{8r}{|q|^2} -2\frac{e_3(h_1)}{h_1}\right)a\Re(\Jk) +\frac{8a^2\cos\th}{|q|^2}\dual\Re(\Jk)  +2\frac{\nab(h_1)}{h_1}\right)\\
&& -\frac 1 2 \atrchb\left(-\left(-\frac{8r}{|q|^2} -2\frac{e_3(h_1)}{h_1}\right)a\dual\Re(\Jk) +\frac{8a^2\cos\th}{|q|^2}\Re(\Jk)  -2\frac{\dual\nab(h_1)}{h_1}\right)\\
&& +\nab_3\left(\left(-\frac{8r}{|q|^2} -2\frac{e_3(h_1)}{h_1}\right)a\Re(\Jk) +\frac{8a^2\cos\th}{|q|^2}\dual\Re(\Jk)  +2\frac{\nab(h_1)}{h_1}\right)  + \frac{1}{2}{|\tr\Xb|^2}a\Re(\Jk)\\
&& -2{\trchb}\,\etab + \tr\Xb(\eta -\etab)  -2\nab_3\etab +\trchb\, \etab - \atrchb\dual\etab + \frac{1}{2}\Big(\DD(\ov{\tr\Xb})+\tr \Xb H\Big) -\ov{\tr\Xb}\,\ov{H}\bigg\}\\
&=& -\frac{6a}{r^2}h_1\Re(\Jk)+\dk^{\leq 1}\Ga_b.
\eeaa
Plugging in the above, this yields 
\beaa
\squared_2(\nab_3(h_1\pmb\phi_{+2}^{(0)})) &=&  {\frac{\ov{q}}{rq}\left(\frac{r^2}{|q|^2}\right)^{-1}\left(\frac{4ia\cos\th}{|q|^2} -2\frac{e_3(h_1)}{h_1}+2\trchb\right)}\nab_{\T}\pmb\phi_{+2}^{(1)}  + \left(\frac{4}{r^4}+O(mr^{-5})\right)\pmb\phi_{+2}^{(2)}   \\
&& +\frac{\ov{q}}{rq}\left(\frac{r^2}{|q|^2}\right)^{-1}\bigg\{\left(-\frac{8r}{|q|^2} -2\frac{e_3(h_1)}{h_1}\right)a\Re(\Jk) +\frac{8a^2\cos\th}{|q|^2}\dual\Re(\Jk)  +2\frac{\nab(h_1)}{h_1}\\
&&  + 2(\eta -\etab) + 2a{\trchb}\Re(\Jk)\bigg\}\c\nab\pmb\phi_{+2}^{(1)} {-\frac{6a^2(\cos\th)^2}{r^2}\frac{(\ov{q})^2}{r^3}}\nab_{\T}\pmb\phi_{+2}^{(0)}\\ 
&& {-\frac{6a}{r^2}h_1\Re(\Jk)}\c\nab\pmb\phi_{+2}^{(0)}  + \left(\frac{6}{r^3}+O(mr^{-4})\right)\pmb\phi_{+2}^{(1)}  +O(mr^{-3})\pmb\phi_{+2}^{(0)}\\
&& +\err +O(1)\widetilde{\N}^{(0)}_{W,+2} +r\dk^{\leq 1}(\Ga_b\c\nab_3A)+\dk^{\leq 2}(\Ga_b\c A)+r\dk^{\leq 1}\big(\Ga_b \c \Ga_g\c A\big).
\eeaa

Next, we introduce the scalar function $h_2$ given by 
\bea
h_2:=\frac{rq}{\ov{q}}\frac{r^2}{|q|^2}=\frac{r^3}{\ov{q}^2}
\eea
and we compute 
\beaa
\squared_2\big(h_2\nab_3(h_1\pmb\phi_{+2}^{(0)})\big) &=& h_2\squared_2(\nab_3(h_1\pmb\phi_{+2}^{(0)})) -e_3(h_2)\nab_4(\nab_3(h_1\pmb\phi_{+2}^{(0)})) -e_4(h_2)\nab_3(\nab_3(h_1\pmb\phi_{+2}^{(0)}))\\ &&+2\nab(h_2)\c\nab(\nab_3(h_1\pmb\phi_{+2}^{(0)}))+\square(h_2)\nab_3(h_1\pmb\phi_{+2}^{(0)})\\
&=& h_2\squared_2(\nab_3(h_1\pmb\phi_{+2}^{(0)})) -e_3(h_2)\left(2\nab_{\T} -\frac{\Delta}{|q|^2}\nab_3 +2a\Re(\Jk)\c\nab\right)\nab_3(h_1\pmb\phi_{+2}^{(0)})\\
&& -e_4(h_2)\nab_3^2(h_1\pmb\phi_{+2}^{(0)}) +2\nab(h_2)\c\nab\nab_3(h_1\pmb\phi_{+2}^{(0)})+\square(h_2)\nab_3(h_1\pmb\phi_{+2}^{(0)})   
\eeaa
where we used \eqref{eq:definitionofTandPhithataretheapproximateKillingvectorifeldinKerrpert} in the last equality.  We infer
\beaa
\squared_2\big(h_2\nab_3(h_1\pmb\phi_{+2}^{(0)})\big) &=& h_2\squared_2(\nab_3(h_1\pmb\phi_{+2}^{(0)})) -2e_3(h_2)\nab_{\T}\nab_3(h_1\pmb\phi_{+2}^{(0)})+\left(\frac{\Delta}{|q|^2}e_3(h_2) -e_4(h_2)\right)\nab_3^2(h_1\pmb\phi_{+2}^{(0)})\\
&& +\Big(-2e_3(h_2)a\Re(\Jk) +2\nab(h_2)\Big)\c\nab\nab_3(h_1\pmb\phi_{+2}^{(0)})+\square(h_2)\nab_3(h_1\pmb\phi_{+2}^{(0)}).   
\eeaa
Also, since 
\beaa
&&\frac{e_3(h_2)}{h_2} = -\frac{1}{r}+O(ar^{-2})+\Ga_b, \qquad \frac{e_4(h_2)}{h_2} = \frac{1}{r}+O(mr^{-2})+r^{-1}\Ga_g,\\
&& \frac{1}{h_2}\nab(h_2) =O(ar^{-2}) +\Ga_g, \qquad \square(h_2)=\frac{2}{r}+O(mr^{-2})+\dk^{\leq 1}\Ga_b, \qquad \T(h_2)=r\Ga_b,
\eeaa
we obtain 
\beaa
\squared_2\big(h_2\nab_3(h_1\pmb\phi_{+2}^{(0)})\big) &=& h_2\squared_2(\nab_3(h_1\pmb\phi_{+2}^{(0)})) -2\frac{e_3(h_2)}{h_2}\nab_{\T}(h_2\nab_3(h_1\pmb\phi_{+2}^{(0)}))\\  
&&+\big(-2+O(mr^{-1})\big)\nab_3^2(h_1\pmb\phi_{+2}^{(0)}) +\left(\frac{2}{r}+O(mr^{-2})\right)\nab_3(h_1\pmb\phi_{+2}^{(0)})\\
&&+\left(-2\frac{e_3(h_2)}{h_2}a\Re(\Jk) +2\frac{\nab(h_2)}{h_2}\right)\c\nab(h_2\nab_3(h_1\pmb\phi_{+2}^{(0)}))\\
&&+r\Ga_b\c\nab_3^2(h_1\pmb\phi_{+2}^{(0)})+\dk^{\leq 1}(\Ga_b)\c\nab_3(h_1\pmb\phi_{+2}^{(0)}). 
\eeaa
Since we have in view of Proposition \ref{prop:transportequationphis=plusminus2p=0and1}
\beaa
h_2\nab_3(h_1\pmb\phi_{+2}^{(0)}) &=& \pmb\phi_{+2}^{(1)}+r^2\Ga_bA,\\
\nab_3^2(h_1\pmb\phi_{+2}^{(0)}) &=& \left(\frac{1}{r}+O(ar^{-2})\right)\nab_3\pmb\phi_{+2}^{(1)}+\left(\frac{1}{r^2}+O(ar^{-3})\right)\pmb\phi_{+2}^{(1)}+r\dk^{\leq 1}(\Ga_bA)\\
&=& \left(\frac{1}{r^3}+O(ar^{-4})\right)\pmb\phi_{+2}^{(2)}+\left({\frac{2}{r^2}}+O(ar^{-3})\right)\pmb\phi_{+2}^{(1)}+r\dk^{\leq 1}(\Ga_bA),
\eeaa
we deduce
\beaa
\squared_2\pmb\phi_{+2}^{(1)} &=& h_2\squared_2(\nab_3(h_1\pmb\phi_{+2}^{(0)})) -2\frac{e_3(h_2)}{h_2}\nab_{\T}\pmb\phi_{+2}^{(1)}+\left(-\frac{2}{r^3}+O(mr^{-4})\right)\pmb\phi_{+2}^{(2)}\\
&& +{\left(-\frac{2}{r^2}+O(mr^{-3})\right)}\pmb\phi_{+2}^{(1)}+\left(-2\frac{e_3(h_2)}{h_2}a\Re(\Jk) +2\frac{\nab(h_2)}{h_2}\right)\c\nab\pmb\phi_{+2}^{(1)}\\
&&+r^2\dk^{\leq 1}(\Ga_b\c\nab_3A)+r\dk^{\leq 2}(\Ga_b\c A).
\eeaa
Plugging in the above, we infer
\beaa
\squared_2\pmb\phi_{+2}^{(1)} &=& \left\{{\frac{4ia\cos\th}{|q|^2} -2\frac{e_3(h_1)}{h_1}+2\trchb} -2\frac{e_3(h_2)}{h_2}\right\}\nab_{\T}\pmb\phi_{+2}^{(1)}\\
&&+\left(\frac{2}{r^3}+O(mr^{-4})\right)\pmb\phi_{+2}^{(2)}  +\bigg\{\left(-\frac{8r}{|q|^2} -2\frac{e_3(h_1)}{h_1} -2\frac{e_3(h_2)}{h_2}\right)a\Re(\Jk)\\
&& +\frac{8a^2\cos\th}{|q|^2}\dual\Re(\Jk)  +2\frac{\nab(h_1)}{h_1} +2\frac{\nab(h_2)}{h_2}  + 2(\eta -\etab) + 2a{\trchb}\Re(\Jk)\bigg\}\c\nab\pmb\phi_{+2}^{(1)}\\
&&  {-\frac{6a^2(\cos\th)^2}{r^2}}\nab_{\T}\pmb\phi_{+2}^{(0)} {-\frac{6a}{r^2}|q|^2\Re(\Jk)}\c\nab\pmb\phi_{+2}^{(0)}  + \left({\frac{4}{r^2}}+O(mr^{-3})\right)\pmb\phi_{+2}^{(1)} \\
&& +O(mr^{-2})\pmb\phi_{+2}^{(0)} +O(r)\err +O(r)\widetilde{\N}^{(0)}_{W,+2}\\
&& +r^2\dk^{\leq 1}(\Ga_b\c\nab_3A)+r\dk^{\leq 2}(\Ga_b\c A)+r^2\dk^{\leq 1}\big(\Ga_b \c \Ga_g\c A\big).
\eeaa

Next, using in particular the fact that $h_1h_2=|q|^2$, we compute
\beaa
{\frac{4ia\cos\th}{|q|^2} -2\frac{e_3(h_1)}{h_1}+2\trchb} -2\frac{e_3(h_2)}{h_2} &=& \frac{4ia\cos\th}{|q|^2} -2\frac{e_3(h_1h_2)}{h_1h_2} -\frac{4r}{|q|^2}+\Ga_g\\
&=& \frac{4ia\cos\th}{|q|^2} -2\frac{e_3(|q|^2)}{|q|^2} -\frac{4r}{|q|^2}+\Ga_g\\
&=& \frac{4ia\cos\th}{|q|^2}+\Ga_b
\eeaa
and 
\begin{align*}
& \left(-\frac{8r}{|q|^2} -2\frac{e_3(h_1)}{h_1} -2\frac{e_3(h_2)}{h_2}\right)a\Re(\Jk) +\frac{8a^2\cos\th}{|q|^2}\dual\Re(\Jk)  +2\frac{\nab(h_1)}{h_1} +2\frac{\nab(h_2)}{h_2}  + 2(\eta -\etab) + 2a{\trchb}\Re(\Jk)\\
=& \left(-\frac{8r}{|q|^2} -2\frac{e_3(h_1h_2)}{h_1h_2} \right)a\Re(\Jk) +\frac{8a^2\cos\th}{|q|^2}\dual\Re(\Jk)  +2\frac{\nab(h_1h_2)}{h_1h_2} \\
&  + 2\left(\frac{2ar}{|q|^2}\Re(\Jk)+\Ga_b\right) + 2a\left(-\frac{2r}{|q|^2}+\Ga_g\right)\Re(\Jk)\\ 
=& \left(-\frac{8r}{|q|^2} -2\frac{e_3(|q|^2)}{|q|^2} \right)a\Re(\Jk) +\frac{8a^2\cos\th}{|q|^2}\dual\Re(\Jk)  +2\frac{\nab(|q|^2)}{|q|^2}+\Ga_b\\
=& -\frac{4r}{|q|^2}a\Re(\Jk) +\frac{4a^2\cos\th}{|q|^2}\dual\Re(\Jk) +\Ga_b.
\end{align*}
Plugging in the above, we deduce 
\beaa
\squared_2\pmb\phi_{+2}^{(1)} &=& \frac{4ia\cos\th}{|q|^2}\nab_{\T}\pmb\phi_{+2}^{(1)} +\left(\frac{2}{r^3}+O(mr^{-4})\right)\pmb\phi_{+2}^{(2)}  +\left(-\frac{4r}{|q|^2}a\Re(\Jk) +\frac{4a^2\cos\th}{|q|^2}\dual\Re(\Jk)\right)\c\nab\pmb\phi_{+2}^{(1)}\nn\\
&&  {-\frac{6a^2(\cos\th)^2}{r^2}}\nab_{\T}\pmb\phi_{+2}^{(0)} {-\frac{6a}{r^2}|q|^2\Re(\Jk)}\c\nab\pmb\phi_{+2}^{(0)}  + \left({\frac{4}{r^2}}+O(mr^{-3})\right)\pmb\phi_{+2}^{(1)} +O(mr^{-2})\pmb\phi_{+2}^{(0)}\\
&& +O(r)\err +O(r)\widetilde{\N}^{(0)}_{W,+2} +r^2\dk^{\leq 1}(\Ga_b\c\nab_3A)+r\dk^{\leq 2}(\Ga_b\c A)+r^2\dk^{\leq 1}\big(\Ga_b \c \Ga_g\c A\big),
\eeaa
which we write in simpler form as follows, using also \eqref{eq:formofregularhorizontalvectorfieldwidetildemathcalXs:Kerrperturbation} \eqref{eq:formofregularhorizontalvectorfieldwidehatZ:Kerrperturbation},
\beaa
\squared_2\pmb\phi_{+2}^{(1)} -\frac{4ia\cos\th}{|q|^2}\nab_{\T}\pmb\phi_{+2}^{(1)} -\frac{4}{|q|^2}\pmb\phi_{+2}^{(1)} &=& \left(\frac{2}{r^3}+O(mr^{-4})\right)\pmb\phi_{+2}^{(2)}  +O(mr^{-3})\nab_{\widehat{\mathcal{X}}_{+2}}^{\leq 1}\pmb\phi_{+2}^{(1)}\\
&&  +O(mr^{-2})\nab_{\widehat{\Z}}^{\leq 1}\pmb\phi_{+2}^{(0)} +\widetilde{\N}_{W,+2}^{(1)},
\eeaa
with all the coefficients on the RHS  being independent of coordinates $\tau$ and $\phi$, and with $\widetilde{\N}^{(1)}_{W,+2}$ being given by 
\beaa
\widetilde{\N}^{(1)}_{W,+2}=  O(r)\err +O(r)\widetilde{\N}^{(0)}_{W,+2} +r^2\dk^{\leq 1}(\Ga_b\c\nab_3A)+r\dk^{\leq 2}(\Ga_b\c A)+r^2\dk^{\leq 1}\big(\Ga_b \c \Ga_g\c A\big)
\eeaa
so that, in view of  the structure of $\widetilde{\N}^{(0)}_{W,+2}$ and \eqref{eq:structureoferr=errortermindecompositionofsquared_2nab3h1pmbphiplus20}, we have
\beaa
   \bsplit
 \widetilde{\N}^{(1)}_{W,+2} =& r\dk^{\leq 2}\big( \Ga_g \c  B\big)+r^2\dk^{\leq 1}\big(\nab_3\Xi  \c B) +r\dk^{\leq 1}\big((\dkb^{\leq1} \Pc+ r^{-2} \Ga_b)\c \Xi\big)\\
& +r^2\dk^{\leq 1}((\Ga_g \c \Xi) \c \Bb) +r^2\dk^{\leq 1}(\Ga_b\c\nab_3A)+r\dk^{\leq 2}(\Ga_b\c A)+r^2\dk^{\leq 1}\big(\Ga_b \c \Ga_g\c A\big).
\end{split}
 \eeaa
This concludes the proof of \eqref{eq:TensorialTeuSysandlinearterms:rescaleRHScontaine2:general:Kerrperturbation:alternateformnullframeinsteadcoordvectorfield} in the case $s=+2$ and $p=1$.

%%%%%%%%%%%%%%%%%%%%%%%%%%%%%%%%%%%%%%%%%%%%%%%%%%

\subsubsection{Derivation of the wave equation for $\pmb\phi_{+2}^{(2)}$ in perturbations of Kerr}

%%%%%%%%%%%%%%%%%%%%%%%%%%%%%%%%%%%%%%%%%%%%%%%%%%

Recall from Remark \ref{rmk:compasisionphiplus2p=0withMaSz26andphiplus2p=2withqfGKS22} that $\qf=\pmb\phi_{+2}^{(2)}$. We may thus use the derivation of the wave equation for $\qf$ in Theorem 5.2.9 of \cite{GKS22} which yields\footnote{Note that the equation in Theorem 5.2.9 of \cite{GKS22} is given in outgoing normalization which accounts for the additional prefactor $\frac{\De}{|q|^2}$ in front of the second order derivative term on the RHS in (5.2.11) of \cite{GKS22}. Note also that we have slightly changed the notation for the lower order terms compared to (5.2.11)  of \cite{GKS22} with the following correspondence $q\ov{q}^3W_4\to W_4$ and similarly for the other lower order terms.}, in an ingoing frame, 
\bea\lab{eq:maingRWqfinGKS22}
\bigg(\squared_2 -\frac{4ia\cos\th}{|q|^2}\nab_{\T}- \frac{4}{|q|^2}\bigg)\pmb\phi^{(2)}_{+2} &=& V\pmb\phi^{(2)}_{+2} -\ov{q}^2\frac{8}{r^2}(a^2\nab_{\T}+a\nab_{\Z})\nab_3A +W_4\nab_4A+W_3\nab_3A\nn\\
&&+W\c\nab A+W_0A+ \err[\squared_2 \qf],
\eea
where:
\begin{itemize}
\item The approximate Killing vectorfields $\T$ and $\Z$ are given by \eqref{eq:definitionofTandPhithataretheapproximateKillingvectorifeldinKerrpert}.

\item $V$, $W_4$, $W_3$ and $W_0$ are complex functions of $(r, \th)$, and $W$ is the product of a complex function of $(r,\th)$ with $\Jk$, with
\beaa
V=O(mr^{-3}), \qquad \Im(V)=0, \qquad W_4, \, W_3,\, W=O(a), \qquad W_0=O(ar^{-1}). 
\eeaa

\item $\err[\squared_2 \qf]$ is  the nonlinear correction term, which under the additional conditions\footnote{The additional conditions hold true in view of our assumption \eqref{eq:additional-conditions-Hc-Xi-gRW-eq:0}.}  
\bea\label{eq:additional-conditions-Hc-Xi-gRW-eq}
\Hc \in \Ga_g{, \qquad \nab_3\Xi\in r^{-1}\dk^{\leq 1}\Ga_g,}
\eea
is given schematically by  the expression
   \bea
   \lab{eq:MaiThmParq-err}
   \bsplit
 \err[\squared_2 \qf]=& r^2 \dk^{\leq 3} (\Ga_g \c (A, B))+ {\nab_3 (r^3 \dk^{\leq 2}( \Ga_g \c (A, B)))} {+ \nab_3 (r^3\dk^{\leq1}\big(\Xi \c (\dkb^{\leq1} \Pc,  r^{-2}\Ga_b )\big))}\\
 &+\dk^{\leq 1} (\Ga_g \c \qf) {+\nab_3(r^4(\dk^{\leq 1}(\Ga_g\c\Xi)\c\Bb))}.
\end{split}
 \eea
\end{itemize}

Next, notice that we have\footnote{\lab{footnotelabel:simpleobservaationReJkcnabcosth}Notice in particular that $\Re(\Jk)\c\nab(\cos\th)=-\Re(\Jk)\c\dual\Re(\Jk)+r^{-1}\Ga_b=r^{-1}\Ga_b$.}
\bea\lab{eq:actionofTandPhioncoordinatesrandcosth}
\T(r)=r\Ga_b, \qquad \T(\cos\th)=\Ga_b, \qquad \Z(r)=r^2\Ga_g, \qquad \Z(\cos\th)=r\Ga_b,
\eea
so that 
\beaa
\ov{q}^2\frac{8}{r^2}(a^2\nab_{\T}+a\nab_{\Z})\nab_3A = \frac{8}{r^2}(a^2\nab_{\T}+a\nab_{\Z})(\ov{q}^2\nab_3A) +\Ga_g\nab_3A.
\eeaa
Next, we have together with \eqref{eq:definitionofthephiplus2phierarchy:perturbationofKerr}, 
\beaa
\frac{8}{r^2}(a^2\nab_{\T}+a\nab_{\Z})\pmb\phi_{+2}^{(1)} &=& \frac{8}{r^2}(a^2\nab_{\T}+a\nab_{\Z})(\ov{q}^2\nab_3A)\\
&&+\frac{8\ov{q}^2}{r^2}\left(\frac{1}{2}\trchb -\frac{3}{2}\frac{\atrchb^2}{\trchb} -2i\atrchb\right)(a^2\nab_{\T}+a\nab_{\Z})A+\dk^{\leq 1}(\Ga_b)A\\
&=& \frac{8}{r^2}(a^2\nab_{\T}+a\nab_{\Z})(\ov{q}^2\nab_3A)+\frac{8\ov{q}^2}{r^2}\left(\frac{3}{r} -\frac{4}{\ov{q}}\right)(a^2\nab_{\T}+a\nab_{\Z})A\\
&&+\Ga_g\dk^{\leq 1}A+\dk^{\leq 1}(\Ga_b)A.
\eeaa
Plugging in \eqref{eq:maingRWqfinGKS22}, we infer
\beaa
\bigg(\squared_2 -\frac{4ia\cos\th}{|q|^2}\nab_{\T}- \frac{4}{|q|^2}\bigg)\pmb\phi^{(2)}_{+2} &=& V\pmb\phi^{(2)}_{+2} -\frac{8}{r^2}(a^2\nab_{\T}+a\nab_{\Z})\pmb\phi_{+2}^{(1)}\\
&& +\frac{8\ov{q}^2}{r^2}\left(\frac{3}{r} -\frac{4}{\ov{q}}\right)(a^2\nab_{\T}+a\nab_{\Z})A+W_4\nab_4A+W_3\nab_3A\nn\\
&&+W\c\nab A+W_0A+ \err[\squared_2 \qf] +\Ga_g\dk^{\leq 1}A+\dk^{\leq 1}(\Ga_b)A
\eeaa
and hence
\beaa
\bigg(\squared_2 -\frac{4ia\cos\th}{|q|^2}\nab_{\T}- \frac{4}{|q|^2}\bigg)\pmb\phi^{(2)}_{+2} &=& V\pmb\phi^{(2)}_{+2} -\frac{8}{r^2}(a^2\nab_{\T}+a\nab_{\Z})\pmb\phi_{+2}^{(1)}\\
&& +\left(W_3+\frac{8\ov{q}^2}{r^2}\left(\frac{3}{r} -\frac{4}{\ov{q}}\right)\left(\frac{a^2\De}{2|q|^2} -\frac{a^2(\sin\th)^2\De}{2|q|^2}\right)\right)\nab_3A\\
&& +\left(W_4 +\frac{8\ov{q}^2}{r^2}\left(\frac{3}{r} -\frac{4}{\ov{q}}\right)\left(\frac{a^2}{2} -\frac{a^2(\sin\th)^2}{2}\right)\right)\nab_4A\\
&& +\left(W +\frac{8\ov{q}^2}{r^2}\left(\frac{3}{r} -\frac{4}{\ov{q}}\right)\big(-a^3\Re(\Jk)+a(r^2+a^2)\Re(\Jk)\big)\right)\c\nab A\nn\\
&&+W_0A+ \err[\squared_2 \qf] +\Ga_g\dk^{\leq 1}A+\dk^{\leq 1}(\Ga_b)A,
\eeaa
or,
\beaa
\bigg(\squared_2 -\frac{4ia\cos\th}{|q|^2}\nab_{\T}- \frac{4}{|q|^2}\bigg)\pmb\phi^{(2)}_{+2} &=& V\pmb\phi^{(2)}_{+2} -\frac{8}{r^2}(a^2\nab_{\T}+a\nab_{\Z})\pmb\phi_{+2}^{(1)}\\
&& +\left(W_3+\frac{4a^2(\cos\th)^2\De\ov{q}^2}{r^2|q|^2}\left(\frac{3}{r} -\frac{4}{\ov{q}}\right)\right)\nab_3A\\
&& +\left(W_4 +\frac{4a^2(\cos\th)^2\ov{q}^2}{r^2}\left(\frac{3}{r} -\frac{4}{\ov{q}}\right)\right)\nab_4A\\
&& +\left(W +8a\ov{q}^2\left(\frac{3}{r} -\frac{4}{\ov{q}}\right)\Re(\Jk)\right)\c\nab A\nn\\
&&+W_0A+ \err[\squared_2 \qf] +\Ga_g\dk^{\leq 1}A+\dk^{\leq 1}(\Ga_b)A.
\eeaa
Next, we check  the value of $W_4$, $W_3$ and $W$ in Theorem D.4.7 of \cite{GKS22} verify 
\beaa
W_4 = -\frac{4a^2(\cos\th)^2\ov{q}^2}{r^2}\left(\frac{3}{r} -\frac{4}{\ov{q}}\right),\qquad W = -8a\ov{q}^2\left(\frac{3}{r} -\frac{4}{\ov{q}}\right)\Re(\Jk),
\eeaa
and
\beaa
W_3+\frac{4a^2(\cos\th)^2\De\ov{q}^2}{r^2|q|^2}\left(\frac{3}{r} -\frac{4}{\ov{q}}\right) &=& \ov{q}^2\left(\frac{16ai}{r^2}+O(a^3r^{-4})\right),
\eeaa
and hence
\beaa
\bigg(\squared_2 -\frac{4ia\cos\th}{|q|^2}\nab_{\T}- \frac{4}{|q|^2}\bigg)\pmb\phi^{(2)}_{+2} &=& V\pmb\phi^{(2)}_{+2} -\frac{8}{r^2}(a^2\nab_{\T}+a\nab_{\Z})\pmb\phi_{+2}^{(1)} +\ov{q}^2\left(\frac{16ai}{r^2}+O(a^3r^{-4})\right)\nab_3A\\
&&+W_0A+ \err[\squared_2 \qf] +\Ga_g\dk^{\leq 1}A+\dk^{\leq 1}(\Ga_b)A.
\eeaa
Now, in view of \eqref{eq:definitionofthephiplus2phierarchy:perturbationofKerr}, we have
\beaa
\left(\frac{16ai}{r^2}+O(a^3r^{-4})\right)\pmb\phi_{+2}^{(1)} &=& \ov{q}^2\left(\frac{16ai}{r^2}+O(a^3r^{-4})\right)\nab_3A +\ov{q}^2\left(\frac{16ai}{r^2}+O(a^3r^{-4})\right)\left(\frac{3}{r} -\frac{4}{\ov{q}}\right)A+\Ga_bA\\
&=& \ov{q}^2\left(\frac{16ai}{r^2}+O(a^3r^{-4})\right)\nab_3A +\left(-\frac{16ai}{r}+O(a^2r^{-2})\right)A+\Ga_bA 
\eeaa
which yields 
\beaa
\bigg(\squared_2 -\frac{4ia\cos\th}{|q|^2}\nab_{\T}- \frac{4}{|q|^2}\bigg)\pmb\phi^{(2)}_{+2} &=& V\pmb\phi^{(2)}_{+2} -\frac{8}{r^2}(a^2\nab_{\T}+a\nab_{\Z})\pmb\phi_{+2}^{(1)}\\
&& +\left(\frac{16ai}{r^2}+O(a^3r^{-4})\right)\pmb\phi_{+2}^{(1)} +\left(W_0 -\frac{16ai}{r}+O(a^2r^{-2})\right)A\\
&&+ \err[\squared_2 \qf] +\Ga_g\dk^{\leq 1}A+\dk^{\leq 1}(\Ga_b)A.
\eeaa
Next, we check that the value of $W_0$ in Theorem D.4.7 of \cite{GKS22} verifies 
\beaa
W_0 -\frac{16ai}{r} =O(a^2r^{-2})
\eeaa
and hence, we obtain 
\beaa
\bigg(\squared_2 -\frac{4ia\cos\th}{|q|^2}\nab_{\T}- \frac{4}{|q|^2}\bigg)\pmb\phi^{(2)}_{+2} &=& V\pmb\phi^{(2)}_{+2} -\frac{8}{r^2}(a^2\nab_{\T}+a\nab_{\Z})\pmb\phi_{+2}^{(1)}\nn\\
&& +\left(\frac{16ai}{r^2}+O(a^3r^{-4})\right)\pmb\phi_{+2}^{(1)} +O(a^2r^{-2})\pmb\phi_{+2}^{(0)}\nn\\
&&+ \err[\squared_2 \qf] +\Ga_g\dk^{\leq 1}A+\dk^{\leq 1}(\Ga_b)A,
\eeaa
which we write in simpler form as follows 
\beaa
\bigg(\squared_2 -\frac{4ia\cos\th}{|q|^2}\nab_{\T}- \frac{4}{|q|^2}\bigg)\pmb\phi^{(2)}_{+2} &=& O(mr^{-3})\pmb\phi^{(2)}_{+2}+O(mr^{-2})\nab_{\Z+a\T}^{\leq 1}\pmb\phi^{(1)}_{+2}+O(m^2 r^{-2})\pmb\phi^{(0)}_{+2} +\widetilde{\N}^{(2)}_{W,+2},
\eeaa
with all the coefficients on the RHS  being independent of coordinates $\tau$ and $\tphi$, with the coefficients in front of the terms $\pmb\phi^{(2)}_{+2}$ and $\nab_{\Z+a\T}\pmb\phi^{(1)}_{+2}$ on the RHS being real functions, and with $\widetilde{\N}^{(2)}_{W,+2}$ being given by 
\beaa
\widetilde{\N}^{(2)}_{W,+2}=  \err[\squared_2 \qf] +\Ga_g\dk^{\leq 1}A+\dk^{\leq 1}(\Ga_b)A 
\eeaa
so that, in view of  \eqref{eq:MaiThmParq-err}, we have
\beaa
   \bsplit
 \widetilde{\N}^{(2)}_{W,+2} &= r^2 \dk^{\leq 3} (\Ga_g \c (A, B))+ {\nab_3 (r^3 \dk^{\leq 2}( \Ga_g \c (A, B)))} \\
 &{+ \nab_3 (r^3\dk^{\leq1}\big(\Xi \c (\dkb^{\leq1} \Pc,  r^{-2}\Ga_b )\big))}+\dk^{\leq 1} (\Ga_g \c \qf) {+\nab_3(r^4(\dk^{\leq 1}(\Ga_g\c\Xi)\c\Bb))}.
\end{split}
 \eeaa
This concludes the proof of \eqref{eq:TensorialTeuSysandlinearterms:rescaleRHScontaine2:general:Kerrperturbation:alternateformnullframeinsteadcoordvectorfield} in the case $s=+2$ and $p=2$, and hence the one of Theorem \ref{thm:derivationoftheTeukolskytensorialwavesystemfors=plusminus2:kerrpert:alternateformnullframeinsteadcoordvectorfield} in the case $s=+2$.

%%%%%%%%%%%%%%%%%%%%%%%%%%%%%%%%%%%%%%%%%%%%%%%%%%%%%%%%%%%%%%%%%%%%%%%%%%%%

\subsection{Wave equations for $\pmb\phi_{-2}^{(p)}$, $p=0,1,2$}
\lab{sec:proofofthm:derivationoftheTeukolskytensorialwavesystemfors=plusminus2:kerrpert:alternateformnullframeinsteadcoordvectorfield:cases=-2}

%%%%%%%%%%%%%%%%%%%%%%%%%%%%%%%%%%%%%%%%%%%%%%%%%%%%%%%%%%%%%%%%%%%%%%%%%%%%%

In this section, we prove Theorem \ref{thm:derivationoftheTeukolskytensorialwavesystemfors=plusminus2:kerrpert:alternateformnullframeinsteadcoordvectorfield} in the case $s=-2$, i.e., the fact that $\pmb\phi_{-2}^{(p)}$, $p=0,1,2$, verifies the system of tensorial wave equations \eqref{eq:TensorialTeuSysandlinearterms:rescaleRHScontaine2:general:Kerrperturbation:alternateformnullframeinsteadcoordvectorfield}. Note that in order to exhibit the structure of the nonlinear correction terms $\widetilde{\N}_{W,-2}^{(p)}$, $p=0,1,2$, stated in \eqref{eq:schematicformofNpWsminus2}, we will use in particular the additional assumptions \eqref{eq:Xi-Hb-chapter12:0} on the global null frame.

%%%%%%%%%%%%%%%%%%%%%%%%%%%%%%%%%%%%%%%%%%%%%%%%%%

\subsubsection{Derivation of the wave equation for $\pmb\phi_{-2}^{(0)}$ in perturbations of Kerr}

%%%%%%%%%%%%%%%%%%%%%%%%%%%%%%%%%%%%%%%%%%%%%%%%%%

Recall from Lemma 5.3.3 in \cite{GKS22} that
\bea\lab{eq:squared2AbformofTeukolskyfromLemma533inGKS22}
\squared_2 \Ab&=& \left(2\ov{\tr \Xb} +4\omb \right) \nab_4 \Ab-4\om \nab_3 \Ab -\left( 4\Hb-4\ze \right)\c \nab \Ab +V \Ab+ \err
 \eea
 where 
 \beaa
 V&=& \frac 3 4 \tr X \ov{\tr \Xb} -\frac 1 4 \tr \Xb  \ov{\tr X}  - 4P-4\om\left( \frac 1 2 \tr \Xb +2\ov{\tr \Xb} \right)+2 \omb \tr X-4\nab_3\om  -8\omb\om\\
&&+  \DD\c \ov{Z}- 2Z \c \ov{Z}+2 \ze\c \left( 4\Hb+2\eta \right)- 4 \etab \c \eta  +2i\etab \wedge \eta
 \eeaa
 and the error terms are given by
 \bea\lab{eq:errorterminTeukforAbinwaveform}
 \err&=& r^{-1}\dk^{\leq 1}\big(  \Ga_b \c \Bb \big) + \Ga_b \c \Ga_b \c \Ga_g.
 \eea

Next, using \eqref{eq:definitionofthephiminus2phierarchy:perturbationofKerr}, we have
\beaa
\squared_2\pmb\phi_{-2}^{(0)} &=& \squared_2\left(\frac{q}{\ov{q}}\left(\frac{\De}{|q|^2}\right)^2\Ab\right)\\
&=& \frac{q}{\ov{q}}\left(\frac{\De}{|q|^2}\right)^2\squared_2\Ab+2\g^{\a\b}e_\a\left(\frac{q}{\ov{q}}\left(\frac{\De}{|q|^2}\right)^2\right)\Ddot_\b\Ab +\square\left(\frac{q}{\ov{q}}\left(\frac{\De}{|q|^2}\right)^2\right)\Ab\\
&=& \frac{q}{\ov{q}}\left(\frac{\De}{|q|^2}\right)^2\squared_2\Ab -e_3\left(\frac{q}{\ov{q}}\left(\frac{\De}{|q|^2}\right)^2\right)\nab_4\Ab -e_4\left(\frac{q}{\ov{q}}\left(\frac{\De}{|q|^2}\right)^2\right)\nab_3\Ab\\
&&+2\nab\left(\frac{q}{\ov{q}}\left(\frac{\De}{|q|^2}\right)^2\right)\c\nab\Ab+\square\left(\frac{q}{\ov{q}}\left(\frac{\De}{|q|^2}\right)^2\right)\Ab.
\eeaa
Plugging in the above, this yields 
\beaa
\squared_2\pmb\phi_{-2}^{(0)} &=& \frac{q}{\ov{q}}\left(\frac{\De}{|q|^2}\right)^2\bigg\{\left(2\ov{\tr \Xb} +4\omb \right) \nab_4 \Ab-4\om \nab_3 \Ab -\left( 4\Hb-4\ze \right)\c \nab \Ab +V \Ab+ \err\bigg\}\\
&& -e_3\left(\frac{q}{\ov{q}}\left(\frac{\De}{|q|^2}\right)^2\right)\nab_4\Ab -e_4\left(\frac{q}{\ov{q}}\left(\frac{\De}{|q|^2}\right)^2\right)\nab_3\Ab+2\nab\left(\frac{q}{\ov{q}}\left(\frac{\De}{|q|^2}\right)^2\right)\c\nab\Ab\\
&&+\square\left(\frac{q}{\ov{q}}\left(\frac{\De}{|q|^2}\right)^2\right)\Ab,
\eeaa
or
\beaa
\squared_2\pmb\phi_{-2}^{(0)} &=& \left(\frac{q}{\ov{q}}\left(\frac{\De}{|q|^2}\right)^2\left(2\ov{\tr \Xb} +4\omb \right)  -e_3\left(\frac{q}{\ov{q}}\left(\frac{\De}{|q|^2}\right)^2\right)\right)\nab_4 \Ab\\
&& +\left(-\frac{q}{\ov{q}}\left(\frac{\De}{|q|^2}\right)^24\om  -e_4\left(\frac{q}{\ov{q}}\left(\frac{\De}{|q|^2}\right)^2\right)\right)\nab_3\Ab\\
&&+\left(-\frac{q}{\ov{q}}\left(\frac{\De}{|q|^2}\right)^2\left( 4\Hb-4\ze \right)+2\nab\left(\frac{q}{\ov{q}}\left(\frac{\De}{|q|^2}\right)^2\right)\right)\c\nab A\\
&&+\left(\frac{q}{\ov{q}}\left(\frac{\De}{|q|^2}\right)^2V +\square\left(\frac{q}{\ov{q}}\left(\frac{\De}{|q|^2}\right)^2\right)\right)\Ab +\frac{q}{\ov{q}}\left(\frac{\De}{|q|^2}\right)^2\err.
\eeaa
Next, we compute 
\beaa
&& -\frac{q}{\ov{q}}\left(\frac{\De}{|q|^2}\right)^24\om  -e_4\left(\frac{q}{\ov{q}}\left(\frac{\De}{|q|^2}\right)^2\right)\\
&=&
-\frac{2q}{\ov{q}}\left(\frac{\De}{|q|^2}\right)\left(2\om\left(\frac{\De}{|q|^2}\right)+e_4\left(\frac{\De}{|q|^2}\right)\right) -e_4\left(\frac{q}{\ov{q}}\right)\left(\frac{\De}{|q|^2}\right)^2\\
&=& -e_4\left(\frac{q}{\ov{q}}\right)\left(\frac{\De}{|q|^2}\right)^2 + \Ga_g 
\eeaa
which together with \eqref{eq:definitionofTandPhithataretheapproximateKillingvectorifeldinKerrpert} implies 
\beaa
\squared_2\pmb\phi_{-2}^{(0)} &=& \left(\frac{q}{\ov{q}}\left(\frac{\De}{|q|^2}\right)^2\left(2\ov{\tr \Xb} +4\omb \right)  -e_3\left(\frac{q}{\ov{q}}\left(\frac{\De}{|q|^2}\right)^2\right)\right)\nab_4 \Ab\\
&& -e_4\left(\frac{q}{\ov{q}}\right)\frac{\De}{|q|^2}\Big(2\nab_{\T} -\nab_4 + 2a\Re(\Jk)\c\nab\Big)\Ab\\
&&+\left(-\frac{q}{\ov{q}}\left(\frac{\De}{|q|^2}\right)^2\left( 4\Hb-4\ze \right)+2\nab\left(\frac{q}{\ov{q}}\left(\frac{\De}{|q|^2}\right)^2\right)\right)\c\nab A\\
&&+\left(\frac{q}{\ov{q}}\left(\frac{\De}{|q|^2}\right)^2V +\square\left(\frac{q}{\ov{q}}\left(\frac{\De}{|q|^2}\right)^2\right)\right)\Ab +\frac{q}{\ov{q}}\left(\frac{\De}{|q|^2}\right)^2\err +\Ga_g\nab_3\Ab.
\eeaa

Next, we have
\beaa
\frac{\ov{q}}{q}e_4\left(\frac{q}{\ov{q}}\right) &=& \left(\frac{1}{q} - \frac{1}{\ov{q}}\right)\frac{\De}{|q|^2}+r^{-1}\Ga_g = -\frac{2ia\cos\th}{|q|^2}\frac{\De}{|q|^2}+r^{-1}\Ga_g
\eeaa
so that 
\beaa
\frac{2\De}{|q|^2}e_4\left(\frac{q}{\ov{q}}\right)\nab_{\T}\Ab &=& -\frac{4ia\cos\th}{|q|^2}\frac{q}{\ov{q}}\left(\frac{\De}{|q|^2}\right)^2\nab_{\T}\Ab+r^{-1}\Ga_g\nab_{\T}\Ab\\
&=& -\frac{4ia\cos\th}{|q|^2}\nab_{\T}\pmb\phi_{-2}^{(0)}+r^{-1}\Ga_g\dk^{\leq 1}\Ab.
\eeaa
Also, we have
\beaa
&& \left(\frac{\De}{|q|^2}\right)^2\left(2\ov{\tr \Xb} +4\omb \right)  - \frac{\ov{q}}{q}e_3\left(\frac{q}{\ov{q}}\left(\frac{\De}{|q|^2}\right)^2\right) + \frac{\ov{q}}{q}e_4\left(\frac{q}{\ov{q}}\right)\frac{\De}{|q|^2}\\
&=& -\frac{4}{q}\left(\frac{\De}{|q|^2}\right)^2 +2\left(\frac{\De}{|q|^2}\right)\pr_r\left(\frac{\De}{|q|^2}\right)+2\left(\frac{\De}{|q|^2}\right)^2\left(\frac{1}{q}-\frac{1}{\ov{q}}\right)+\Ga_b\\
&=& \left(-\frac{4r}{|q|^2}\frac{\De}{|q|^2} +2\pr_r\left(\frac{\De}{|q|^2}\right)\right)\frac{\De}{|q|^2} +\Ga_b,
\eeaa
and
\beaa
&& -\left(\frac{\De}{|q|^2}\right)^2\left( 4\Hb-4\ze \right)+2\frac{\ov{q}}{q}\nab\left(\frac{q}{\ov{q}}\left(\frac{\De}{|q|^2}\right)^2\right) -\frac{\ov{q}}{q}e_4\left(\frac{q}{\ov{q}}\right)\frac{\De}{|q|^2}2a\Re(\Jk)\\
&=& -\left(\frac{\De}{|q|^2}\right)^2\left(-\frac{4a\ov{q}}{|q|^2}\Jk -\Re\left(\frac{4aq}{|q|^2}\Jk\right)\right) -2\left(\frac{1}{q}+\frac{1}{\ov{q}}\right)\left(\frac{\De}{|q|^2}\right)^2ai\dual\Re(\Jk)\\
&&-4\left(\frac{\De}{|q|^2}\right)^2\frac{\nab(|q|^2)}{|q|^2} -\left(\frac{1}{q}-\frac{1}{\ov{q}}\right)\left(\frac{\De}{|q|^2}\right)^22a\Re(\Jk)+\Ga_g\\
&=& -\left(\frac{\De}{|q|^2}\right)^2\left(-\frac{4a\ov{q}}{|q|^2}\Jk - \frac{2aq}{|q|^2}\Jk - \frac{2a\ov{q}}{|q|^2}\ov{\Jk}\right) -\frac{4r}{|q|^2}\left(\frac{\De}{|q|^2}\right)^2ai\dual\Re(\Jk)\\
&&+\frac{8a^2\cos\th}{|q|^2}\left(\frac{\De}{|q|^2}\right)^2\dual\Re(\Jk)+\frac{4ia^2\cos\th}{|q|^2}\left(\frac{\De}{|q|^2}\right)^2\Re(\Jk)+\Ga_g\\
&=& \left(\frac{\De}{|q|^2}\right)^2\frac{8ar}{|q|^2}\Re(\Jk)+\left(\frac{\De}{|q|^2}\right)^2\frac{8a^2\cos\th}{|q|^2}\dual\Re(\Jk)+\Ga_g.
\eeaa
Plugging in the above, we infer
\beaa
\squared_2\pmb\phi_{-2}^{(0)}  -\frac{4ia\cos\th}{|q|^2}\nab_{\T}\pmb\phi_{-2}^{(0)} &=& \left(-\frac{4r}{|q|^2}\frac{\De}{|q|^2} +2\pr_r\left(\frac{\De}{|q|^2}\right)\right)\frac{q}{\ov{q}}\frac{\De}{|q|^2}\nab_4\Ab\\
&&+\left(\frac{8ar}{|q|^2}\Re(\Jk)+\frac{8a^2\cos\th}{|q|^2}\dual\Re(\Jk)\right)\c \frac{q}{\ov{q}}\left(\frac{\De}{|q|^2}\right)^2\nab A\\
&&+\left(\frac{q}{\ov{q}}\left(\frac{\De}{|q|^2}\right)^2V +\square\left(\frac{q}{\ov{q}}\left(\frac{\De}{|q|^2}\right)^2\right)\right)\Ab\\
&& +\frac{q}{\ov{q}}\left(\frac{\De}{|q|^2}\right)^2\err +\Ga_g\nab_3\Ab +r^{-1}\Ga_b\dk^{\leq 1}\Ab
\eeaa
which we rewrite as 
\beaa
\squared_2\pmb\phi_{-2}^{(0)}  -\frac{4ia\cos\th}{|q|^2}\nab_{\T}\pmb\phi_{-2}^{(0)} &=& \left(-\frac{4r}{|q|^2}\frac{\De}{|q|^2} +2\pr_r\left(\frac{\De}{|q|^2}\right)\right)\frac{q}{\ov{q}}\frac{\De}{|q|^2}\nabc_4\Ab\\
&&+\left(\frac{8ar}{|q|^2}\Re(\Jk)+\frac{8a^2\cos\th}{|q|^2}\dual\Re(\Jk)\right)\c \frac{q}{\ov{q}}\left(\frac{\De}{|q|^2}\right)^2\nab A\\
&&+\bigg\{\frac{q}{\ov{q}}\left(\frac{\De}{|q|^2}\right)^2V +\square\left(\frac{q}{\ov{q}}\left(\frac{\De}{|q|^2}\right)^2\right)\\
&&+4\om\left(-\frac{4r}{|q|^2}\frac{\De}{|q|^2} +2\pr_r\left(\frac{\De}{|q|^2}\right)\right)\frac{q}{\ov{q}}\frac{\De}{|q|^2}\bigg\}\Ab\\
&& +\frac{q}{\ov{q}}\left(\frac{\De}{|q|^2}\right)^2\err +\Ga_g\nab_3\Ab +r^{-1}\Ga_b\dk^{\leq 1}\Ab.
\eeaa
In view of
\beaa
&&\frac{q}{\ov{q}}\left(\frac{\De}{|q|^2}\right)^2V +\square\left(\frac{q}{\ov{q}}\left(\frac{\De}{|q|^2}\right)^2\right) +4\om\left(-\frac{4r}{|q|^2}\frac{\De}{|q|^2} +2\pr_r\left(\frac{\De}{|q|^2}\right)\right)\frac{q}{\ov{q}}\frac{\De}{|q|^2}\\
&=& \left(-\frac{2}{r^2}+O(mr^{-3})\right)\frac{q}{\ov{q}}\left(\frac{\De}{|q|^2}\right)^2+r^{-1}\dk^{\leq 1}\Ga_b,
\eeaa
we infer
\beaa
\squared_2\pmb\phi_{-2}^{(0)}  -\frac{4ia\cos\th}{|q|^2}\nab_{\T}\pmb\phi_{-2}^{(0)} &=& \left(-\frac{4r}{|q|^2}\frac{\De}{|q|^2} +2\pr_r\left(\frac{\De}{|q|^2}\right)\right)\frac{q}{\ov{q}}\frac{\De}{|q|^2}\nabc_4\Ab\\
&&+\left(\frac{8ar}{|q|^2}\Re(\Jk)+\frac{8a^2\cos\th}{|q|^2}\dual\Re(\Jk)\right)\c \frac{q}{\ov{q}}\left(\frac{\De}{|q|^2}\right)^2\nab A\\
&&+\left(-\frac{2}{r^2}+O(mr^{-3})\right)\frac{q}{\ov{q}}\left(\frac{\De}{|q|^2}\right)^2\Ab \\
&&+\frac{q}{\ov{q}}\left(\frac{\De}{|q|^2}\right)^2\err +\Ga_g\nab_3\Ab +r^{-1}\Ga_b\dk^{\leq 1}\Ab.
\eeaa
Since
\beaa
\frac{q}{\ov{q}}\frac{\De}{|q|^2}\nabc_4\Ab &=& \frac{1}{|q|^2}\Big(\pmb\phi_{-2}^{(1)} - r\pmb\phi_{-2}^{(0)} +O(a)\pmb\phi_{-2}^{(0)}+r^2\Ga_g\Ab\Big),\\
\frac{q}{\ov{q}}\left(\frac{\De}{|q|^2}\right)^2\nab A &=& \nab\pmb\phi_{-2}^{(0)}+O(mr^{-2})\pmb\phi_{-2}^{(0)}+\Ga_g\Ab,
\eeaa
we deduce
\beaa
\squared_2\pmb\phi_{-2}^{(0)}  -\frac{4ia\cos\th}{|q|^2}\nab_{\T}\pmb\phi_{-2}^{(0)} &=& \left(-\frac{4r}{|q|^2}\frac{\De}{|q|^2} +2\pr_r\left(\frac{\De}{|q|^2}\right)\right)\frac{1}{|q|^2}\pmb\phi_{-2}^{(1)}\nn\\
&&+\left(\frac{8ar}{|q|^2}\Re(\Jk)+\frac{8a^2\cos\th}{|q|^2}\dual\Re(\Jk)\right)\c\nab\pmb\phi_{-2}^{(0)}\nn\\
&&+\left(\frac{2}{r^2}+O(mr^{-3})\right)\pmb\phi_{-2}^{(0)} +\frac{q}{\ov{q}}\left(\frac{\De}{|q|^2}\right)^2\err \nn\\
&&+r^{-1}\Ga_g\pmb\phi_{-2}^{(1)}+\Ga_g\nab_3\Ab +r^{-1}\Ga_b\dk^{\leq 1}\Ab,
\eeaa
which we write in simpler form as follows, using also \eqref{eq:formofregularhorizontalvectorfieldwidetildemathcalXs:Kerrperturbation},
\beaa
\squared_2\pmb\phi_{-2}^{(0)} -\frac{4ia\cos\th}{|q|^2}\nab_{\T}\pmb\phi_{-2}^{(0)} -\frac{2}{|q|^2}\pmb\phi_{-2}^{(0)} = \left(-\frac{4}{r^3}+O(mr^{-4})\right)\pmb\phi_{-2}^{(1)}    +O(mr^{-3})\nab_{\widehat{\mathcal{X}}_{-2}}^{\leq 1}\pmb\phi_{-2}^{(0)}   +\widetilde{\N}_{W,-2}^{(0)},
\eeaa
with all the coefficients on the RHS  being independent of coordinates $\tau$ and $\phi$, and with $\widetilde{\N}^{(0)}_{W,-2}$ being given by 
\beaa
\widetilde{\N}^{(0)}_{W,-2}=  \frac{q}{\ov{q}}\left(\frac{\De}{|q|^2}\right)^2\err +r^{-1}\Ga_g\pmb\phi_{-2}^{(1)}+\Ga_g\nab_3\Ab +r^{-1}\Ga_b\dk^{\leq 1}\Ab
\eeaa
so that, in view of  \eqref{eq:errorterminTeukforAbinwaveform}, we have
\beaa
   \bsplit
 \widetilde{\N}^{(0)}_{W,-2} =& r^{-1}\dk^{\leq 1}\big(  \Ga_b \c \Bb \big) + \Ga_b \c \Ga_b \c \Ga_g +r^{-1}\Ga_g\pmb\phi_{-2}^{(1)}+\Ga_g\nab_3\Ab +r^{-1}\Ga_b\dk^{\leq 1}\Ab.
\end{split}
 \eeaa
 This concludes the proof of \eqref{eq:TensorialTeuSysandlinearterms:rescaleRHScontaine2:general:Kerrperturbation:alternateformnullframeinsteadcoordvectorfield} in the case $s=-2$ and $p=0$.

%%%%%%%%%%%%%%%%%%%%%%%%%%%%%%%%%%%%%%%%%%%%%%%%%%%

\subsubsection{Derivation of the wave equation for $\pmb\phi_{-2}^{(1)}$ in perturbations of Kerr}

%%%%%%%%%%%%%%%%%%%%%%%%%%%%%%%%%%%%%%%%%%%%%%%%%%%

Let $\Ab_4\in\sk_2(\mathbb{C})$ be given by
\beaa
\Ab_4:= \nabc_4\Ab+\frac{1}{2}\tr X\Ab=\nab_4\Ab +\left(-4\om+\frac{1}{2}\tr X\right)\Ab.
\eeaa
$\Ab_4$ verifies in view of Bianchi identities 
\bea\lab{eq:Ab4isthegoodcomboofnabcAbandtrXAbdecayingbetterinr}
\Ab_4=r^{-2}\dk^{\leq 1}\Ga_b.
\eea
Also, in view of \eqref{eq:definitionofthephiminus2phierarchy:perturbationofKerr}, we have
\bea\lab{eq:linkAb4withpmbphiminus2p=1}
\nn\pmb\phi_{-2}^{(1)} &=& q^2\frac{\De}{|q|^2}\left(\nabc_4 +\frac{1}{2}\trch  -\frac{3}{2}\left(\frac{\atrch}{\trch}\right)_{\vartheta}\atrch -2i\atrch\right)\Ab\\
&=& q^2\frac{\De}{|q|^2}\Ab_4+O(a)\pmb\phi_{-2}^{(0)}+r^2\Ga_g\c\Ga_b,
\eea
and, using in addition \eqref{eq:Cb1-Cb2-comparison-Ma}, 
\beaa
\nn\left(\nabc_4+\frac{3}{2}\tr X\right)\Ab_4 &=& \left(\nabc_4+\frac{3}{2}\tr X\right)\left(\nabc_4\Ab+\frac{1}{2}\tr X\Ab\right)\\
\nn&=& \nabc_4^2\Ab+2\tr X\nabc_4\Ab+\left(\frac{1}{2}\nabc_4\tr X+\frac{3}{4}(\tr X)^2\right)\Ab\\
\nn&=& \nabc_4^2\Ab+2\tr X\nabc_4\Ab+\left(\frac{1}{2}\nabc_4\tr X+\frac{3}{4}(\tr X)^2\right)\Ab\\
\nn&=& \nabc_4^2\Ab+2\tr X\nabc_4\Ab+\left(\frac{1}{2}(\tr X)^2+r^{-1}\dk^{\leq 1}\Xi+\Ga_g\c\Ga_g\right)\Ab\\
\nn&=&  \frac{1}{|q|^2q^2}\phi_{-2}^{(2)} +\Big(2\tr X - \und{C}_1\Big)\nabc_4\Ab\\
&&+\left(\frac{1}{2}(\tr X)^2-\und{C}_2+r^{-1}\dk^{\leq 1}\Xi+\Ga_g\c\Ga_g\right)\Ab\\
\nn&=& \frac{1}{|q|^2q^2}\phi_{-2}^{(2)} + O(ar^{-4}\De)\Ab_4+ O(a^2r^{-5}\De)\nabc_4\Ab+O(a^2r^{-8}\De^2)\Ab\\
&&+r^{-1}\dk^{\leq 1}\Xi\c\Ga_b+r^{-2}\dk^{\leq 1}(\Ga_g\c\Ga_b),
\eeaa
which together with \eqref{eq:linkAb4withpmbphiminus2p=1} yields
\bea\lab{eq:linknabc4Ab4withpmbphiminus2p=2}
\nn\left(\nabc_4+\frac{3}{2}\tr X\right)\Ab_4 &=&  \frac{1}{|q|^2q^2}\pmb\phi_{-2}^{(2)} + O(ar^{-4})\pmb\phi_{-2}^{(1)}+O(a^2r^{-4})\pmb\phi_{-2}^{(0)}\\
&&+r^{-1}\dk^{\leq 1}\Xi\c\Ga_b+r^{-2}\dk^{\leq 1}(\Ga_g\c\Ga_b).
\eea

Next, according to Proposition D.7.1 in \cite{GKS22}, $\Ab_4$ satisfies 
\beaa
\begin{split}
&\Big( \nabc_4+  \tr X+ \frac 1 2\ov{\tr X}\Big) \Big( \nabc_3+ 2 \ov{\tr \Xb}+\frac 1 2 \tr \Xb\Big) \Ab_4\\
=&    \frac 1 4\big(  \DDc+H+5\Hb \big) \hot \Big( \DDbc \c \Ab_4+(\ov{H}+\ov{\Hb}) \c \Ab_4\Big)+3P\Ab_4\\
& +3\Big( \frac 1 2\ov{\tr X}-\tr X\Big) P\Ab -\frac 3 4 \tr X \Hb\hot \big(\DDbc \c\Ab +  \ov{H} \c \Ab \big)+  r^{-2}\dk^{\leq 2}(\Ga_b\c\Ga_g).
\end{split}
\eeaa
Using the definition of conformally invariant derivatives, together with the fact that $\Ab$ is $-2$-conformally invariant and $\Ab_4$ is $-1$-conformally invariant, we infer 
\beaa
\begin{split}
&\Big( \nab_4 -4\om +  \tr X+ \frac 1 2\ov{\tr X}\Big) \Big( \nab_3 +2\omb+ 2 \ov{\tr \Xb}+\frac 1 2 \tr \Xb\Big) \Ab_4\\
=&    \frac 1 4\big(  \DD -Z+H+5\Hb \big) \hot \Big( \DDb\c\Ab_4+(-\ov{Z}+\ov{H}+\ov{\Hb}) \c \Ab_4\Big)+3P\Ab_4\\
& +3\Big( \frac 1 2\ov{\tr X}-\tr X\Big) P\Ab -\frac 3 4 \tr X \Hb\hot \big(\DDb\c\Ab + \big(-2\ov{Z}+ \ov{H}\big)\c \Ab \big)+  r^{-2}\dk^{\leq 2}(\Ga_b\c\Ga_g),
\end{split}
\eeaa
which, together with \eqref{eq:Ab4isthegoodcomboofnabcAbandtrXAbdecayingbetterinr} and the fact that $\omb\in\Ga_b$ and $Z=H+\Ga_b$, implies
\beaa
\begin{split}
&\Big( \nab_4 -4\om +  \tr X+ \frac 1 2\ov{\tr X}\Big) \Big( \nab_3 + 2 \ov{\tr \Xb}+\frac 1 2 \tr \Xb\Big) \Ab_4\\
=&    \frac 1 4\big(  \DD +5\Hb \big) \hot \Big( \DDb\c\Ab_4+\ov{\Hb}\c \Ab_4\Big)+3P\Ab_4\\
& +3\Big( \frac 1 2\ov{\tr X}-\tr X\Big) P\Ab -\frac 3 4 \tr X \Hb\hot \big(\DDb\c\Ab -\ov{H}\c \Ab \big)+  r^{-2}\dk^{\leq 2}(\Ga_b\c\Ga_g).
\end{split}
\eeaa
Now, using Lemmas \ref{dot-hot-complex}  and \ref{SIMPLIFICATION-ANGULAR}, we have
\beaa
\big(  \DD +5\Hb \big) \hot \Big( \DDb\c\Ab_4+\ov{\Hb}\c \Ab_4\Big) &=& \DD\hot(\DDb\c\Ab_4) +20\Hb\c\nab\Ab_4 + 2(\ov{\Hb}\c\DD)\Ab_4 + 2(\DD\c\ov{\Hb})\Ab_4\\
&& +10|\Hb|^2\Ab_4\\
&=& \DD\hot(\DDb\c\Ab_4) +20\Hb\c\nab\Ab_4 + 2(\etab-i\dual\etab)\c(\nab+i\dual\nab)\Ab_4 \\
&& +O(ar^{-3})\Ab_4+r^{-3}\Ga_g\dk^{\leq 1}\Ga_b\\
&=& \DD\hot(\DDb\c\Ab_4) +\big(20\Hb+4\ov{\Hb}\big)\c\nab\Ab_4 +O(ar^{-3})\Ab_4\\
&&+r^{-3}\Ga_g\dk^{\leq 1}\Ga_b
\eeaa
and 
\beaa
\Hb\hot \big(\DDb\c\Ab -\ov{H}\c \Ab \big) &=& 4\Hb\c\nab\Ab - 4 \big(\etab\c\eta - i \eta\wedge\etab\big)\Ab\\
&=&  4\Hb\c\nab\Ab+O(a^2\De r^{-6})\Ab+r^{-2}\Ga_b\c\Ga_b,
\eeaa
and hence
\beaa
\begin{split}
&\Big( \nab_4 -4\om +  \tr X+ \frac 1 2\ov{\tr X}\Big) \Big( \nab_3 + 2 \ov{\tr \Xb}+\frac 1 2 \tr \Xb\Big) \Ab_4\\
=&  \frac 1 4\DD\hot(\DDb\c\Ab_4) +\big(5\Hb+\ov{\Hb}\big)\c\nab\Ab_4 +O(mr^{-3})\Ab_4 - 3\tr X\Hb\c\nab\Ab+O(m\De r^{-6})\Ab +  r^{-2}\dk^{\leq 2}(\Ga_b\c\Ga_g).
\end{split}
\eeaa
Then, we compute, using in particular \eqref{eq:Ab4isthegoodcomboofnabcAbandtrXAbdecayingbetterinr}, 
\beaa
&&\Big( \nab_4 -4\om +  \tr X+ \frac 1 2\ov{\tr X}\Big) \Big( \nab_3 + 2 \ov{\tr \Xb}+\frac 1 2 \tr \Xb\Big) \Ab_4\\
&=& \nab_4\nab_3\Ab_4 +\left(-4\om +  \tr X+ \frac 1 2\ov{\tr X}\right)\nab_3\Ab_4 +\Big( \nab_4 -4\om +  \tr X+ \frac 1 2\ov{\tr X}\Big)\Big(2 \ov{\tr \Xb}+\frac 1 2 \tr \Xb\Big) \Ab_4\\
&=& \nab_4\nab_3\Ab_4 +\left(-4\om +  \tr X+ \frac 1 2\ov{\tr X}\right)\nab_3\Ab_4 +\left(-\frac{5}{r}+O(ar^{-2})+\Ga_g\right)\nab_4\Ab_4\\
&&+\left(\nab_4\Big(2 \ov{\tr \Xb}+\frac 1 2 \tr \Xb\Big) + \Big(-4\om +  \tr X+ \frac 1 2\ov{\tr X}\Big)\Big(2 \ov{\tr \Xb}+\frac 1 2 \tr \Xb\Big)\right)\Ab_4\\
&=& \nab_4\nab_3\Ab_4 +\left(-4\om +  \tr X+ \frac 1 2\ov{\tr X}\right)\nab_3\Ab_4 +\left(-\frac{5}{r}+O(ar^{-2})\right)\nab_4\Ab_4\\
&&+\left(-\frac{10}{r^2}+O(mr^{-3})\right)\Ab_4 +r^{-3}\dk^{\leq 2}(\Ga_g\c\Ga_b),
\eeaa
and hence
\beaa
\begin{split}
&-\nab_4\nab_3\Ab_4 +\frac 1 4\DD\hot(\DDb\c\Ab_4) +\left(4\om -  \tr X - \frac 1 2\ov{\tr X}\right)\nab_3\Ab_4 +\left(\frac{5}{r}+O(ar^{-2})\right)\nab_4\Ab_4\\
&+\big(5\Hb+\ov{\Hb}\big)\c\nab\Ab_4+\left(\frac{10}{r^2}+O(mr^{-3})\right)\Ab_4 - 3\tr X\Hb\c\nab\Ab+O(m\De r^{-6})\Ab\\
=&   r^{-2}\dk^{\leq 2}(\Ga_b\c\Ga_g).
\end{split}
\eeaa
In particular, we have
\beaa
\begin{split}
&-\nab_4\nab_3\Ab_4 +\frac 1 4\DD\hot(\DDb\c\Ab_4) +\left(2\om -  \frac{1}{2}\tr X\right)\nab_3\Ab_4 -\frac{1}{2}\tr\Xb\nab_4\Ab_4 \\
& +2\etab\c\nab\Ab_4 + \frac{2}{r^2}\Ab_4\\
=& \big(-2\om + \trch\big)\nab_3\Ab_4 +\left(-\frac{4}{r}+O(ar^{-2})\right)\nab_4\Ab_4 -4\Hb\c\nab\Ab_4 +\left(-\frac{8}{r^2}+O(mr^{-3})\right)\Ab_4\\
& + 3\tr X\Hb\c\nab\Ab +O(m\De r^{-6})\Ab+  r^{-2}\dk^{\leq 2}(\Ga_b\c\Ga_g),
\end{split}
\eeaa
which implies, in view of Corollary \ref{corollary-wave-complex},
\beaa
\begin{split}
\squared_2\Ab_4 =& \big(-2\om + \trch\big)\nab_3\Ab_4 +\left(-\frac{4}{r}+O(ar^{-2})\right)\nab_4\Ab_4 -4\Hb\c\nab\Ab_4 +\left(-\frac{8}{r^2}+O(mr^{-3})\right)\Ab_4\\
& + 3\tr X\Hb\c\nab\Ab +O(m\De r^{-6})\Ab+  r^{-2}\dk^{\leq 2}(\Ga_b\c\Ga_g).
\end{split}
\eeaa

Next, we introduce the scalar function $\und{h}_1$ defined by 
\bea
\und{h}_1:=q^2\frac{\De}{|q|^2}. 
\eea
We compute 
\beaa
\squared_2\big(\und{h}_1\Ab_4\big) &=& \und{h}_1\squared_2\Ab_4 -e_4(\und{h}_1)\nab_3\Ab_4 -e_3(\und{h}_1)\nab_4\Ab_4 +2\nab(\und{h}_1)\c\nab\Ab_4+\square(\und{h}_1)\Ab_4. 
\eeaa
Since 
\beaa
&& e_4(\und{h}_1) = \left(\frac{2}{q}\frac{\De}{|q|^2}+\pr_r\left(\frac{\De}{|q|^2}\right)\right)\und{h}_1+r\Ga_g, \qquad e_3(\und{h}_1)=-\left(\frac{2}{q}\frac{\De}{|q|^2}+\pr_r\left(\frac{\De}{|q|^2}\right)\right)q^2+r^2\Ga_b,\\
&& \nab(\und{h}_1) = -\frac{2r}{|q|^2}\und{h}_1ia\dual\Re(\Jk)+r^2\Ga_g, \qquad\square(\und{h}_1)=\frac{\De}{|q|^2}\big(6+O(mr^{-1})\big)+\left(\pr_r\left(\frac{\De}{|q|^2}\right)\right)^2q^2+r\dk^{\leq 1}\Ga_b,
\eeaa
we obtain, using also \eqref{eq:Ab4isthegoodcomboofnabcAbandtrXAbdecayingbetterinr}, 
\begin{align*}
\squared_2\big(\und{h}_1\Ab_4\big) =& \und{h}_1\squared_2\Ab_4 -\left(\frac{2}{q}\frac{\De}{|q|^2}+\pr_r\left(\frac{\De}{|q|^2}\right)\right)q^2\frac{\De}{|q|^2}\nab_3\Ab_4 +\left(\frac{2}{q}\frac{\De}{|q|^2}+\pr_r\left(\frac{\De}{|q|^2}\right)\right)q^2\nab_4\Ab_4\\
& -\frac{4r}{|q|^2}\und{h}_1ia\dual\Re(\Jk)\c\nab\Ab_4+\left(\frac{\De}{|q|^2}\big(6+O(mr^{-1})\big)+\left(\pr_r\left(\frac{\De}{|q|^2}\right)\right)^2q^2\right)\Ab_4+r^{-1}\dk^{\leq 2}(\Ga_b\c\Ga_b). 
\end{align*}
Plugging the above expression for $\squared_2\Ab_4$,  we deduce, using also  \eqref{eq:Ab4isthegoodcomboofnabcAbandtrXAbdecayingbetterinr}, 
\beaa
\squared_2\big(\und{h}_1\Ab_4\big) &=& \frac{2ia\cos\th}{|q|^2}\und{h}_1\frac{\De}{|q|^2}\nab_3\Ab_4  +\left(\big(-2r+O(a)\big)\frac{\De}{|q|^2}+q^2\pr_r\left(\frac{\De}{|q|^2}\right)\right)\nab_4\Ab_4 \\
&& -4\und{h}_1\left(\Hb +\frac{r}{|q|^2}ia\dual\Re(\Jk)\right)\c\nab\Ab_4 +\bigg\{\frac{\De}{|q|^2}\big(-2+O(mr^{-1})\big)+\left(\pr_r\left(\frac{\De}{|q|^2}\right)\right)^2q^2\bigg\}\Ab_4\\
&& + 3\tr X q^2\frac{\De}{|q|^2}\Hb\c\nab\Ab +O(m\De^2r^{-6})\Ab  +  \dk^{\leq 2}(\Ga_b\c\Ga_g).
\eeaa
Since 
\beaa
q^2\pr_r\left(\frac{\De}{|q|^2}\right)\nab_4\Ab_4+\left(\pr_r\left(\frac{\De}{|q|^2}\right)\right)^2q^2\Ab_4 &=& q^2\pr_r\left(\frac{\De}{|q|^2}\right)\left(\nab_4\Ab_4+\pr_r\left(\frac{\De}{|q|^2}\right)\Ab_4\right)\\
&=& q^2\pr_r\left(\frac{\De}{|q|^2}\right)\Big(\nab_4\Ab_4 -2\om\Ab_4+\Ga_g\Ab_4\Big)\\
&=& O(m)\nabc_4\Ab_4+r^{-2}\Ga_g\dk^{\leq 1}\Ga_b, 
\eeaa
we infer
\beaa
\squared_2\big(\und{h}_1\Ab_4\big) &=& \frac{2ia\cos\th}{|q|^2}\und{h}_1\frac{\De}{|q|^2}\nab_3\Ab_4  +\big(-2r+O(m)\big)\nabc_4\Ab_4  -4\und{h}_1\left(\Hb +\frac{r}{|q|^2}ia\dual\Re(\Jk)\right)\c\nab\Ab_4\\
&& +\frac{\De}{|q|^2}\big(-2+O(mr^{-1})\big)\Ab_4 + 3\tr X q^2\frac{\De}{|q|^2}\Hb\c\nab\Ab +O(m\De^2r^{-6})\Ab  +  \dk^{\leq 2}(\Ga_b\c\Ga_g).
\eeaa
Together with \eqref{eq:definitionofTandPhithataretheapproximateKillingvectorifeldinKerrpert}, this yields 
\beaa
\squared_2\big(\und{h}_1\Ab_4\big) &=& \frac{2ia\cos\th}{|q|^2}\und{h}_1\Big(2\nab_{\T}  - \nab_4 + 2a\Re(\Jk)\c\nab\Big)\Ab_4  +\big(-2r+O(m)\big)\nabc_4\Ab_4\\
&&  -4\und{h}_1\left(\Hb +\frac{r}{|q|^2}ia\dual\Re(\Jk)\right)\c\nab\Ab_4 +\frac{\De}{|q|^2}\big(-2+O(mr^{-1})\big)\Ab_4\\
&& + 3\tr X q^2\frac{\De}{|q|^2}\Hb\c\nab\Ab +O(m\De^2r^{-6})\Ab  +  \dk^{\leq 2}(\Ga_b\c\Ga_g).
\eeaa
or, using also \eqref{eq:Ab4isthegoodcomboofnabcAbandtrXAbdecayingbetterinr}, 
\beaa
\squared_2\big(\und{h}_1\Ab_4\big) &=& \frac{4ia\cos\th}{|q|^2}\nab_{\T}(\und{h}_1\Ab_4)  +\big(-2r+O(m)\big)\nabc_4\Ab_4 \\
&& +4\und{h}_1\left(-\Hb -\frac{r}{|q|^2}ia\dual\Re(\Jk) +\frac{ia\cos\th}{|q|^2}a\Re(\Jk)\right)\c\nab\Ab_4 +\frac{\De}{|q|^2}\big(-2+O(mr^{-1})\big)\Ab_4\\
&& + 3\tr X q^2\frac{\De}{|q|^2}\Hb\c\nab\Ab  +O(m\De^2r^{-6})\Ab  +  \dk^{\leq 2}(\Ga_b\c\Ga_g).
\eeaa
Also, in view of 
\beaa
\Hb &=& -\frac{ar}{|q|^2}\Re(\Jk) - \frac{a^2\cos\th}{|q|^2}\dual\Re(\Jk) -\frac{air}{|q|^2}\dual\Re(\Jk) + \frac{a^2i\cos\th}{|q|^2}\Re(\Jk) +\Ga_g,
\eeaa
we infer
\beaa
\squared_2\big(\und{h}_1\Ab_4\big) &=& \frac{4ia\cos\th}{|q|^2}\nab_{\T}(\und{h}_1\Ab_4)  +\big(-2r+O(m)\big)\nabc_4\Ab_4 \\
&& +4\und{h}_1\left(\frac{ar}{|q|^2}\Re(\Jk) + \frac{a^2\cos\th}{|q|^2}\dual\Re(\Jk)\right)\c\nab\Ab_4 + \frac{\De}{|q|^2}\big(-2+O(mr^{-1})\big)\Ab_4\\
&& + 3\tr X q^2\frac{\De}{|q|^2}\Hb\c\nab\Ab  +O(m\De^2r^{-6})\Ab  +  \dk^{\leq 2}(\Ga_b\c\Ga_g).
\eeaa
Together with \eqref{eq:definitionofthephiminus2phierarchy:perturbationofKerr}, \eqref{eq:linkAb4withpmbphiminus2p=1} and \eqref{eq:linknabc4Ab4withpmbphiminus2p=2}, this yields 
\beaa
\squared_2\big(\und{h}_1\Ab_4\big) &=& \frac{4ia\cos\th}{|q|^2}\nab_{\T}(\und{h}_1\Ab_4)  +\left(-\frac{2}{r^3}+O(mr^{-4})\right)\pmb\phi_{-2}^{(2)} \\
&& +4\und{h}_1\left(\frac{ar}{|q|^2}\Re(\Jk) + \frac{a^2\cos\th}{|q|^2}\dual\Re(\Jk)\right)\c\nab\Ab_4 + \left(\frac{4}{r^2}+O(mr^{-3})\right)\pmb\phi_{-2}^{(1)}\\
&& + 3\tr X q^2\frac{\De}{|q|^2}\Hb\c\nab\Ab  +O(mr^{-2})\pmb\phi_{-2}^{(0)} +  \dk^{\leq 2}(\Ga_b\c\Ga_g).
\eeaa
Now, we have 
\beaa
4\left(\frac{ar}{|q|^2}\Re(\Jk) + \frac{a^2\cos\th}{|q|^2}\dual\Re(\Jk)\right)\c\nab(\und{h}_1) &=&  O(a^2r^{-4}\De)+\Ga_b,
\eeaa
where we use Footnote \ref{footnotelabel:simpleobservaationReJkcnabcosth}, for the computation of $\Re(\Jk)\c\nab(\und{h}_1)$, and hence, using also \eqref{eq:linkAb4withpmbphiminus2p=1} and \eqref{eq:Ab4isthegoodcomboofnabcAbandtrXAbdecayingbetterinr}, we obtain
\beaa
\squared_2\big(\und{h}_1\Ab_4\big) &=& \frac{4ia\cos\th}{|q|^2}\nab_{\T}(\und{h}_1\Ab_4)  +\left(-\frac{2}{r^3}+O(mr^{-4})\right)\pmb\phi_{-2}^{(2)} \\
&&+4\left(\frac{ar}{|q|^2}\Re(\Jk) + \frac{a^2\cos\th}{|q|^2}\dual\Re(\Jk)\right)\c\nab(\und{h}_1\Ab_4) + \left(\frac{4}{r^2}+O(mr^{-3})\right)\pmb\phi_{-2}^{(1)}\\
&& + 3\tr X q^2\frac{\De}{|q|^2}\Hb\c\nab\Ab  +O(mr^{-2})\pmb\phi_{-2}^{(0)} +  \dk^{\leq 2}(\Ga_b\c\Ga_g),
\eeaa
or
\bea\lab{eq:intermediarytensorialwaveeqationforundh1Ab4}
\nn&&\squared_2\big(\und{h}_1\Ab_4\big) - \frac{4ia\cos\th}{|q|^2}\nab_{\T}(\und{h}_1\Ab_4) - \frac{4}{|q|^2}\pmb\phi_{-2}^{(1)}  \\
\nn&=&   \left(-\frac{2}{r^3}+O(mr^{-4})\right)\pmb\phi_{-2}^{(2)} +4\left(\frac{ar}{|q|^2}\Re(\Jk) + \frac{a^2\cos\th}{|q|^2}\dual\Re(\Jk)\right)\c\nab(\und{h}_1\Ab_4)  + O(mr^{-3})\pmb\phi_{-2}^{(1)}\\
&& + 3\tr X q^2\frac{\De}{|q|^2}\Hb\c\nab\Ab +O(mr^{-2})\pmb\phi_{-2}^{(0)} +  \dk^{\leq 2}(\Ga_b\c\Ga_g).
\eea

Next, we derive a more precise version of \eqref{eq:linkAb4withpmbphiminus2p=1}. In view of \eqref{eq:definitionofthephiminus2phierarchy:perturbationofKerr}, we have
\bea\lab{eq:linkAb4withpmbphiminus2p=1:moreprecise}
\pmb\phi_{-2}^{(1)} &=& q^2\frac{\De}{|q|^2}\left(\nabc_4 +\frac{1}{2}\trch  -\frac{3}{2}\left(\frac{\atrch}{\trch}\right)_{\vartheta}\atrch -2i\atrch\right)\Ab = \und{h}_1\Ab_4 - \frac{3i}{2}j\Ab,
\eea
where we have introduced the 1-conformally invariant scalar function $j$ defined by 
 \bea\lab{eq:scalarfuctionjsuchthatphiplus21isscalarfucntiontimeaAb4plusjAb}
 j &:=&  q^2\frac{\De}{|q|^2}\left(1 -i\left(\frac{\atrch}{\trch}\right)_{\vartheta}\right)\atrch.
 \eea 
We now focus on deriving a wave equation for $j\Ab$. We have
\beaa
\squared_2(j\Ab) &=& j\squared_2\Ab-e_3(j)\nab_4\Ab -e_4(j)\nab_3\Ab+2\nab(j)\c\nab\Ab+\square(j)\Ab.
\eeaa
Also, we have from Lemma 5.3.3 in \cite{GKS22} 
\beaa
\squared_2 \Ab = \left(2\ov{\tr \Xb} +4\omb \right) \nab_4 \Ab-4\om \nab_3 \Ab -\left( 4\Hb-4\ze \right)\c \nab \Ab +V \Ab+ r^{-1}\dk^{\leq 1}\big(  \Ga_b \c \Bb \big) + \Ga_b \c \Ga_b \c \Ga_g
 \eeaa
 where 
 \beaa
 V&=& \frac 3 4 \tr X \ov{\tr \Xb} -\frac 1 4 \tr \Xb  \ov{\tr X}  - 4P-4\om\left( \frac 1 2 \tr \Xb +2\ov{\tr \Xb} \right)+2 \omb \tr X-4\nab_3\om  -8\omb\om\\
&&+  \DD\c \ov{Z}- 2Z \c \ov{Z}+2 \ze\c \left( 4\Hb+2\eta \right)- 4 \etab \c \eta  +2i\etab \wedge \eta
 \eeaa
so that
\beaa
\squared_2 \Ab = O(r^{-1})\nab_4 \Ab-4\om \nab_3 \Ab -\left( 4\Hb-4\eta \right)\c \nab \Ab +O(r^{-2})\Ab + r^{-1}\dk^{\leq 1}\big(  \Ga_b \c \Ga_b\big)
\eeaa
or
\beaa
\squared_2 \Ab &=& O(r^{-1})\Ab_4 -4\om \nab_3 \Ab -\left( 4\Hb-4\eta \right)\c \nab \Ab +O(r^{-2})\Ab + r^{-1}\dk^{\leq 1}\big(  \Ga_b \c \Ga_b\big)
\eeaa
which, since $j=O(ar^{-4}\De^2)+r^2\Ga_g$, yields 
\beaa
j\squared_2 \Ab &=& O(ar^{-3}\De)\Ab_4 -4j\om \nab_3 \Ab - j\left( 4\Hb-4\eta \right)\c \nab \Ab +O(ar^{-6}\De^2)\Ab + \dk^{\leq 1}\big(  \Ga_g \c \Ga_b\big).
\eeaa
Together with \eqref{eq:definitionofthephiminus2phierarchy:perturbationofKerr} and \eqref{eq:linkAb4withpmbphiminus2p=1}, this implies
\beaa
j\squared_2 \Ab &=& O(ar^{-3})\pmb\phi_{-2}^{(1)} -4j\om \nab_3 \Ab - j\left( 4\Hb-4\eta \right)\c \nab \Ab +O(ar^{-2})\pmb\phi_{-2}^{(0)} + \dk^{\leq 1}\big(  \Ga_g \c \Ga_b\big).
\eeaa
Plugging in the above equation for $\squared_2(j\Ab)$, we infer
\beaa
\squared_2(j\Ab) &=& O(ar^{-3})\pmb\phi_{-2}^{(1)} -e_3(j)\nab_4\Ab  -\big(e_4(j)+4j\om\big)\nab_3 \Ab +\Big(2\nab(j) - j\big(4\Hb-4\eta\big)\Big)\c \nab \Ab\\
&& +O(ar^{-2})\pmb\phi_{-2}^{(0)} +\square(j)\Ab + \dk^{\leq 1}\big(  \Ga_g \c \Ga_b\big).
\eeaa

Next, we have
\beaa
e_4(j)+4j\om &=& \left(\frac{\De}{|q|^2}\right)^2\pr_r\left(\frac{\De}{|q|^2}\right)q^2\left(1-\frac{ia\cos\th}{r}\right)\frac{2a\cos\th}{|q|^2}+\frac{\De}{|q|^2}\left(1 -i\left(\frac{\atrch}{\trch}\right)_{\vartheta}\right)e_4(q^2\atrch)\\
&& -q^2\frac{\De}{|q|^2}ie_4\left(\left(\frac{\atrch}{\trch}\right)_{\vartheta}\right)\atrch -2q^2\left(\frac{\De}{|q|^2}\right)^2\left(1-\frac{ia\cos\th}{r}\right)\frac{2a\cos\th}{|q|^2}\pr_r\left(\frac{\De}{|q|^2}\right)+\Ga_g\\
&=&  -q^2\left(\frac{\De}{|q|^2}\right)^2\left(1-\frac{ia\cos\th}{r}\right)\frac{2a\cos\th}{|q|^2}\pr_r\left(\frac{\De}{|q|^2}\right) +iq^2\left(\frac{\De}{|q|^2}\right)^3\frac{2a^2(\cos\th)^2}{r^2|q|^2}\\
&&+\frac{\De}{|q|^2}\left(1 -i\left(\frac{\atrch}{\trch}\right)_{\vartheta}\right)q^2\big(e_4(\atrch)+\tr X\atrch\big)+\Ga_g\\
&=&  -\left(\frac{\De}{|q|^2}\right)^3\frac{2ia^2(\cos\th)^2}{r^2} +\Ga_g+r\dk^{\leq 1}\xi
\eeaa
where we used in particular the null structure equation for $e_4(\atrch)$ that can be deduced from taking the imaginary part of the one for $\nab_4\tr X$ in Proposition \ref{prop-nullstr:complex}, and hence
\beaa
\squared_2(j\Ab) &=& O(ar^{-3})\pmb\phi_{-2}^{(1)} -e_3(j)\nab_4\Ab +\Big(2\nab(j) - j\big(4\Hb-4\eta\big)\Big)\c \nab \Ab + \frac{2ia^2(\cos\th)^2}{r^2}\left(\frac{\De}{|q|^2}\right)^3\nab_3 \Ab\\
&& +O(ar^{-2})\pmb\phi_{-2}^{(0)} +\square(j)\Ab + \dk^{\leq 1}\big(  \Ga_g \c \Ga_b\big)+r\dk^{\leq 1}\Ga_b\c\dk^{\leq 1}\xi.
\eeaa
Together with \eqref{eq:definitionofTandPhithataretheapproximateKillingvectorifeldinKerrpert}, this yields 
\beaa
\squared_2(j\Ab) &=& O(ar^{-3})\pmb\phi_{-2}^{(1)} -e_3(j)\nab_4\Ab +\Big(2\nab(j) - j\big(4\Hb-4\eta\big)\Big)\c \nab \Ab\\
&& +\left(\frac{\De}{|q|^2}\right)^2\frac{2ia^2(\cos\th)^2}{r^2}\Big(2\nab_{\T}  - \nab_4 + 2a\Re(\Jk)\c\nab\Big)\Ab\\
&& +O(ar^{-2})\pmb\phi_{-2}^{(0)} +\square(j)\Ab + \dk^{\leq 1}\big(  \Ga_g \c \Ga_b\big)+r\dk^{\leq 1}\Ga_b\c\dk^{\leq 1}\xi,
\eeaa
or
\beaa
\squared_2(j\Ab) &=& O(ar^{-3})\pmb\phi_{-2}^{(1)} -\left(e_3(j)+\left(\frac{\De}{|q|^2}\right)^2\frac{2ia^2(\cos\th)^2}{r^2}\right)\nab_4\Ab\\
&& +\left(2\nab(j) - j\big(4\Hb-4\eta\big)+\left(\frac{\De}{|q|^2}\right)^2\frac{4ia^2(\cos\th)^2}{r^2}a\Re(\Jk)\right)\c \nab \Ab\\
&& +\left(\frac{\De}{|q|^2}\right)^2\frac{4ia^2(\cos\th)^2}{r^2}\nab_{\T}\Ab +O(ar^{-2})\pmb\phi_{-2}^{(0)} +\square(j)\Ab + \dk^{\leq 1}\big(  \Ga_g \c \Ga_b\big)+r\dk^{\leq 1}\Ga_b\c\dk^{\leq 1}\xi,
\eeaa
which we rewrite as 
\beaa
\squared_2(j\Ab) &=& O(ar^{-3})\pmb\phi_{-2}^{(1)} -\left(e_3(j)+\left(\frac{\De}{|q|^2}\right)^2\frac{2a^2i(\cos\th)^2}{r^2}\right)\Ab_4 +\left(\frac{\De}{|q|^2}\right)^2\frac{4ia^2(\cos\th)^2}{r^2}\nab_{\T}\Ab\\
&& +\left(2\nab(j) - j\big(4\Hb-4\eta\big)+\left(\frac{\De}{|q|^2}\right)^2\frac{4ia^2(\cos\th)^2}{r^2}a\Re(\Jk)\right)\c \nab \Ab\\
&& +\left(\square(j)+e_3(j)\left(\frac{1}{2}\tr X -4\om\right)\right)\Ab  +O(ar^{-2})\pmb\phi_{-2}^{(0)} + \dk^{\leq 1}\big(  \Ga_g \c \Ga_b\big)+r\dk^{\leq 1}\Ga_b\c\dk^{\leq 1}\xi.
\eeaa

Next, we have
\beaa
e_3(j) &=& -\pr_r\left(q^2\frac{\De}{|q|^2}\left(1 -\frac{ai\cos\th}{r}\right)\right)\frac{\De}{|q|^2}\frac{2a\cos\th}{|q|^2}+q^2\frac{\De}{|q|^2}\left(1 -\frac{ai\cos\th}{r}\right)e_3(\atrch)+\dk^{\leq 1}\Ga_b\\
&=& -\pr_r\left(q^2\left(\frac{\De}{|q|^2}\right)^2\left(1 -\frac{ai\cos\th}{r}\right)\frac{2a\cos\th}{|q|^2}\right)+r^2\widecheck{\curl\eta}+r\dk^{\leq 1}\Ga_g\\
&=& -4q^2\frac{\De}{|q|^2}\left(1 -\frac{ai\cos\th}{r}\right)\frac{a\cos\th}{|q|^2}\pr_r\left(\frac{\De}{|q|^2}\right) +\frac{2ia^2(\cos\th)^2}{r^2}\left(\frac{\De}{|q|^2}\right)^2+r^2\widecheck{\curl\eta} +r\dk^{\leq 1}\Ga_g.
\eeaa
Also, we have
\beaa
&&\square(j) +e_3(j)\left(\frac{1}{2}\tr X -4\om\right)\\
&=& -e_4(e_3(j)) +\left(2\om -\frac{1}{2}\tr X\right)e_3(j) +e_3(j)\left(\frac{1}{2}\tr X -4\om\right)+O(a\De^2r^{-6})+\dk^{\leq 2}\Ga_g\\
&=& -e_4(e_3(j)) - 2\om e_3(j) +O(a\De^2r^{-6})+\dk^{\leq 2}\Ga_g,
\eeaa
which together with the above computation of $e_3(j)$ yields
\beaa
\square(j) +e_3(j)\left(\frac{1}{2}\tr X -4\om\right) &=& -e_4(r^2\widecheck{\curl\eta}) +O(a\De^2r^{-6})+\dk^{\leq 2}\Ga_g\\
&=& -r^2\Big(e_4(\widecheck{\curl\eta})+\trch\widecheck{\curl\eta}\Big)+O(a\De^2r^{-6})+\dk^{\leq 2}\Ga_g\\
&=& O(a\De^2r^{-6})+r\dk^{\leq 2}\xi+\dk^{\leq 2}\Ga_g,
\eeaa
where we used in the last equality a commutator formula for $[\nab_4, \curl]$ that can be deduced from \eqref{commutator-nab4-ov-DDcF}, and the null structure equation for $\nab_4\eta$ that can be deduced from taking the real part of the one for $\nab_4H$ in Proposition \ref{prop-nullstr:complex}. Plugging in the above and using \eqref{eq:Ab4isthegoodcomboofnabcAbandtrXAbdecayingbetterinr}, we infer
\beaa
\squared_2(j\Ab) &=& O(ar^{-3})\pmb\phi_{-2}^{(1)} +O(am\De r^{-4})\Ab_4 +\left(\frac{\De}{|q|^2}\right)^2\frac{4ia^2(\cos\th)^2}{r^2}\nab_{\T}\Ab\\
&& +\left(2\nab(j) - j\big(4\Hb-4\eta\big)+\left(\frac{\De}{|q|^2}\right)^2\frac{4ia^2(\cos\th)^2}{r^2}a\Re(\Jk)\right)\c \nab \Ab\\
&&   +O(ar^{-2})\pmb\phi_{-2}^{(0)} + \dk^{\leq 1}\big(  \Ga_g \c \Ga_b\big)+r\dk^{\leq 2}(\Ga_b\c\xi),
\eeaa
which together with \eqref{eq:linkAb4withpmbphiminus2p=1} yields 
\beaa
\squared_2(j\Ab) &=& O(ar^{-3})\pmb\phi_{-2}^{(1)}  +\left(\frac{\De}{|q|^2}\right)^2\frac{4ia^2(\cos\th)^2}{r^2}\nab_{\T}\Ab\\
&& +\left(2\nab(j) - j\big(4\Hb-4\eta\big)+\left(\frac{\De}{|q|^2}\right)^2\frac{4ia^2(\cos\th)^2}{r^2}a\Re(\Jk)\right)\c \nab \Ab\\
&&   +O(ar^{-2})\pmb\phi_{-2}^{(0)} + \dk^{\leq 1}\big(  \Ga_g \c \Ga_b\big)+r\dk^{\leq 2}(\Ga_b\c\xi)
\eeaa
and hence
\beaa
\squared_2(j\Ab) - \frac{4ia\cos\th}{|q|^2}\nab_{\T}(j\Ab) &=& O(ar^{-3})\pmb\phi_{-2}^{(1)} + \left(- \frac{4ia\cos\th}{|q|^2}j+\left(\frac{\De}{|q|^2}\right)^2\frac{4ia^2(\cos\th)^2}{r^2}\right)\nab_{\T}\Ab\\
&& +\left(2\nab(j) - j\big(4\Hb-4\eta\big)+\left(\frac{\De}{|q|^2}\right)^2\frac{4ia^2(\cos\th)^2}{r^2}a\Re(\Jk)\right)\c \nab \Ab\\
&&   +O(ar^{-2})\pmb\phi_{-2}^{(0)} + \dk^{\leq 1}\big(  \Ga_g \c \Ga_b\big)+r\dk^{\leq 2}(\Ga_b\c\xi).
\eeaa

Next, we have 
\beaa
&& - \frac{4ia\cos\th}{|q|^2}j+\left(\frac{\De}{|q|^2}\right)^2\frac{4ia^2(\cos\th)^2}{r^2}\\
&=& - \frac{4ia\cos\th}{|q|^2}q^2\frac{\De}{|q|^2}\left(1 -i\frac{a\cos\th}{r}\right)\frac{\De}{|q|^2}\frac{2a\cos\th}{|q|^2} +\left(\frac{\De}{|q|^2}\right)^2\frac{4ia^2(\cos\th)^2}{r^2}+\Ga_g\\
&=& \frac{4ia^2(\cos\th)^2}{r^2|q|^2}\left(\frac{\De}{|q|^2}\right)^2\Big(-2qr+|q|^2\Big)+\Ga_g\\
&=& -\frac{4ia^2(\cos\th)^2}{r^2}\left(\frac{\De}{|q|^2}\right)^2\frac{q}{\ov{q}}+\Ga_g,
\eeaa
so that 
\beaa
\squared_2(j\Ab) - \frac{4ia\cos\th}{|q|^2}\nab_{\T}(j\Ab) &=& O(ar^{-3})\pmb\phi_{-2}^{(1)}  -\frac{4ia^2(\cos\th)^2}{r^2}\nab_{\T}\pmb\phi_{-2}^{(0)}\\
&& +\left(2\nab(j) - j\big(4\Hb-4\eta\big)+\left(\frac{\De}{|q|^2}\right)^2\frac{4ia^2(\cos\th)^2}{r^2}a\Re(\Jk)\right)\c \nab \Ab\\
&&   +O(ar^{-2})\pmb\phi_{-2}^{(0)} + \dk^{\leq 1}\big(  \Ga_g \c \Ga_b\big)+r\dk^{\leq 2}(\Ga_b\c\xi).
\eeaa
Together with \eqref{eq:intermediarytensorialwaveeqationforundh1Ab4} and \eqref{eq:linkAb4withpmbphiminus2p=1:moreprecise}, we infer
\beaa
&&\squared_2\pmb\phi_{-2}^{(1)} - \frac{4ia\cos\th}{|q|^2}\nab_{\T}\pmb\phi_{-2}^{(1)} -\frac{4}{|q|^2}\pmb\phi_{-2}^{(1)}\\
&=& \squared_2\big(\und{h}_1\Ab_4\big) - \frac{4ia\cos\th}{|q|^2}\nab_{\T}(\und{h}_1\Ab_4) - \frac{4}{|q|^2}\pmb\phi_{-2}^{(1)}  - \frac{3i}{2}\left(\squared_2(j\Ab) - \frac{4ia\cos\th}{|q|^2}\nab_{\T}(j\Ab)\right)\\
&=&    \left(-\frac{2}{r^3}+O(mr^{-4})\right)\pmb\phi_{-2}^{(2)}  +4\left(\frac{ar}{|q|^2}\Re(\Jk) + \frac{a^2\cos\th}{|q|^2}\dual\Re(\Jk)\right)\c\nab(\und{h}_1\Ab_4)  + O(mr^{-3})\pmb\phi_{-2}^{(1)} \\
&& +3\left(\tr X q^2\frac{\De}{|q|^2}\Hb -i\nab(j) + 2ij\big(\Hb-\eta\big)+\left(\frac{\De}{|q|^2}\right)^2\frac{2a^2(\cos\th)^2}{r^2}a\Re(\Jk)\right)\c \nab \Ab\\
&&  - \frac{6a^2(\cos\th)^2}{r^2}\nab_{\T}\pmb\phi_{-2}^{(0)} +O(mr^{-2})\pmb\phi_{-2}^{(0)} + \dk^{\leq 2}\big(\Ga_b\c \Ga_g\big)+r\dk^{\leq 2}(\Ga_b\c\xi).
\eeaa
Also, using again \eqref{eq:linkAb4withpmbphiminus2p=1:moreprecise}, we have
\beaa
&&4\left(\frac{ar}{|q|^2}\Re(\Jk) + \frac{a^2\cos\th}{|q|^2}\dual\Re(\Jk)\right)\c\nab(\und{h}_1\Ab_4)\\ 
&=& 4\left(\frac{ar}{|q|^2}\Re(\Jk) + \frac{a^2\cos\th}{|q|^2}\dual\Re(\Jk)\right)\c\nab\pmb\phi_{-2}^{(1)}+6i\left(\frac{ar}{|q|^2}\Re(\Jk) + \frac{a^2\cos\th}{|q|^2}\dual\Re(\Jk)\right)\c\nab(j\Ab)\\
&=& 4\left(\frac{ar}{|q|^2}\Re(\Jk) + \frac{a^2\cos\th}{|q|^2}\dual\Re(\Jk)\right)\c\nab\pmb\phi_{-2}^{(1)} + \frac{6iarj}{|q|^2}\left(\Re(\Jk) +\frac{a\cos\th}{r}\dual\Re(\Jk)\right)\c\nab\Ab\\
&& +O(a^3r^{-4})\pmb\phi_{-2}^{(0)} +r^{-1}\Ga_b\c\dk^{\leq 1}\Ga_g,
\eeaa
where we used in particular the fact that $\Re(\Jk)\c\nab(\cos\th)=r^{-1}\Ga_b$, see Footnote \ref{footnotelabel:simpleobservaationReJkcnabcosth}, for the computation of $\Re(\Jk)\c\nab(j)$. We deduce 
\beaa
&&\squared_2\pmb\phi_{-2}^{(1)} - \frac{4ia\cos\th}{|q|^2}\nab_{\T}\pmb\phi_{-2}^{(1)} -\frac{4}{|q|^2}\pmb\phi_{-2}^{(1)}\\
&=&    \left(-\frac{2}{r^3}+O(mr^{-4})\right)\pmb\phi_{-2}^{(2)} +4\left(\frac{ar}{|q|^2}\Re(\Jk) + \frac{a^2\cos\th}{|q|^2}\dual\Re(\Jk)\right)\c\nab\pmb\phi_{-2}^{(1)}  + O(mr^{-3})\pmb\phi_{-2}^{(1)} \\
&& +3\bigg\{\tr X q^2\frac{\De}{|q|^2}\Hb -i\nab(j) + 2ij\big(\Hb-\eta\big)+\left(\frac{\De}{|q|^2}\right)^2\frac{2a^2(\cos\th)^2}{r^2}a\Re(\Jk)\\
&& +\frac{2iarj}{|q|^2}\left(\Re(\Jk) +\frac{a\cos\th}{r}\dual\Re(\Jk)\right)\bigg\}\c \nab \Ab  - \frac{6a^2(\cos\th)^2}{r^2}\nab_{\T}\pmb\phi_{-2}^{(0)} +O(mr^{-2})\pmb\phi_{-2}^{(0)}\\
&& + \dk^{\leq 2}\big(\Ga_b\c \Ga_g\big)+r\dk^{\leq 2}(\Ga_b\c\xi).
\eeaa

Next, we compute
\beaa
\nab(j) &=& \nab\left(\frac{2qa\cos\th}{r}\left(\frac{\De}{|q|^2}\right)^2\right)+r\dk^{\leq 1}\Ga_g = \frac{2a\De^2}{r}\nab\left(\frac{q\cos\th}{|q|^4}\right)+r\dk^{\leq 1}\Ga_g\\
&=& \frac{2a\De^2}{r}\left(\frac{\cos\th}{|q|^4}\nab(q)+\frac{q}{|q|^4}\nab(\cos\th) -\frac{2q\cos\th}{|q|^6}\nab(|q|^2)\right)+r\dk^{\leq 1}\Ga_g\\
&=& \frac{2a}{r}\left(\frac{\De}{|q|^2}\right)^2\left(-ia\cos\th -q +\frac{4qa^2(\cos\th)^2}{|q|^2}\right)\dual\Re(\Jk)+r\dk^{\leq 1}\Ga_g,
\eeaa
\beaa
\tr X q^2\frac{\De}{|q|^2}\Hb = 2q\left(\frac{\De}{|q|^2}\right)^2\left(-\frac{ar}{|q|^2}\Re(\Jk) - \frac{a^2\cos\th}{|q|^2}\dual\Re(\Jk) -\frac{air}{|q|^2}\dual\Re(\Jk) + \frac{ia^2\cos\th}{|q|^2}\Re(\Jk)\right)+r\Ga_g,
\eeaa
and
\beaa
2j\big(\Hb-\eta\big) &=& \frac{4qa\cos\th}{r}\left(\frac{\De}{|q|^2}\right)^2\bigg\{-\frac{ar}{|q|^2}\Re(\Jk) - \frac{a^2\cos\th}{|q|^2}\dual\Re(\Jk) -\frac{air}{|q|^2}\dual\Re(\Jk) + \frac{ia^2\cos\th}{|q|^2}\Re(\Jk)\\
&& - \frac{ar}{|q|^2}\Re(\Jk) + \frac{a^2\cos\th}{|q|^2}\dual\Re(\Jk)\bigg\}+\Ga_b\\
&=& \frac{4qa\cos\th}{r}\left(\frac{\De}{|q|^2}\right)^2\bigg\{-\frac{2ar}{|q|^2}\Re(\Jk)  -\frac{air}{|q|^2}\dual\Re(\Jk) + \frac{ia^2\cos\th}{|q|^2}\Re(\Jk)\bigg\}+\Ga_b.
\eeaa
This yields
\beaa
&& \tr X q^2\frac{\De}{|q|^2}\Hb -i\nab(j) + 2ij\big(\Hb-\eta\big)+\left(\frac{\De}{|q|^2}\right)^2\frac{2a^2(\cos\th)^2}{r^2}a\Re(\Jk) +\frac{2iarj}{|q|^2}\left(\Re(\Jk) +\frac{a\cos\th}{r}\dual\Re(\Jk)\right)\\
&=& \left(\frac{\De}{|q|^2}\right)^2\bigg\{2q\left(-\frac{ar}{|q|^2}   + \frac{ia^2\cos\th}{|q|^2}\right) +\frac{4qa\cos\th}{r}\left(-\frac{2iar}{|q|^2}  - \frac{a^2\cos\th}{|q|^2}\right)\\
&&+\frac{2a^3(\cos\th)^2}{r^2}+\frac{4ia^2\cos\th q}{|q|^2}\bigg\}\Re(\Jk) +\left(\frac{\De}{|q|^2}\right)^2\bigg\{2q\left( - \frac{a^2\cos\th}{|q|^2} -\frac{air}{|q|^2}\right)\\
&& +\frac{2a}{r}\left(-a\cos\th +iq -\frac{4iqa^2(\cos\th)^2}{|q|^2}\right) +\frac{4qa^2\cos\th}{|q|^2}+\frac{4ia^3(\cos\th)^2q}{r|q|^2}\bigg\}\dual\Re(\Jk) +r\dk^{\leq 1}\Ga_g\\
&=& -\frac{2aq^2}{r^2}\left(\frac{\De}{|q|^2}\right)^2\Re(\Jk) +r\dk^{\leq 1}\Ga_g
\eeaa
and hence
\beaa
&&\squared_2\pmb\phi_{-2}^{(1)} - \frac{4ia\cos\th}{|q|^2}\nab_{\T}\pmb\phi_{-2}^{(1)} -\frac{4}{|q|^2}\pmb\phi_{-2}^{(1)}\nn\\
&=&    \left(-\frac{2}{r^3}+O(mr^{-4})\right)\pmb\phi_{-2}^{(2)} +4\left(\frac{ar}{|q|^2}\Re(\Jk) + \frac{a^2\cos\th}{|q|^2}\dual\Re(\Jk)\right)\c\nab\pmb\phi_{-2}^{(1)} + O(mr^{-3})\pmb\phi_{-2}^{(1)}\nn\\
&&  -\frac{6aq^2}{r^2}\left(\frac{\De}{|q|^2}\right)^2\Re(\Jk)\c \nab \Ab  - \frac{6a^2(\cos\th)^2}{r^2}\nab_{\T}\pmb\phi_{-2}^{(0)} +O(mr^{-2})\pmb\phi_{-2}^{(0)} + \dk^{\leq 2}\big(\Ga_b\c \Ga_g\big)+r\dk^{\leq 2}(\Ga_b\c\xi),
\eeaa
or, using the definition of $\pmb\phi_{-2}^{(0)}$ in \eqref{eq:definitionofthephiminus2phierarchy:perturbationofKerr}
\begin{align}\lab{eq:almostfinalformoftensorialwaveequationforpmbphiminus2p=1perturbationofKerrusingnullframesinsteadofcoordinatesvectorfields}
&\squared_2\pmb\phi_{-2}^{(1)} - \frac{4ia\cos\th}{|q|^2}\nab_{\T}\pmb\phi_{-2}^{(1)} -\frac{4}{|q|^2}\pmb\phi_{-2}^{(1)}\nn\\
=&    \left(-\frac{2}{r^3}+O(mr^{-4})\right)\pmb\phi_{-2}^{(2)} +4\left(\frac{ar}{|q|^2}\Re(\Jk) + \frac{a^2\cos\th}{|q|^2}\dual\Re(\Jk)\right)\c\nab\pmb\phi_{-2}^{(1)} + O(mr^{-3})\pmb\phi_{-2}^{(1)}\nn\\
&  -\frac{6a|q|^2}{r^2}\Re(\Jk)\c\nab\pmb\phi_{-2}^{(0)}  - \frac{6a^2(\cos\th)^2}{r^2}\nab_{\T}\pmb\phi_{-2}^{(0)} +O(mr^{-2})\pmb\phi_{-2}^{(0)} + \dk^{\leq 2}\big(\Ga_b\c \Ga_g\big)+r\dk^{\leq 2}(\Ga_b\c\xi),
\end{align}
which we write in simpler form as follows, \eqref{eq:formofregularhorizontalvectorfieldwidetildemathcalXs:Kerrperturbation} \eqref{eq:formofregularhorizontalvectorfieldwidehatZ:Kerrperturbation},
\beaa
\squared_2\pmb\phi_{-2}^{(1)} - \frac{4ia\cos\th}{|q|^2}\nab_{\T}\pmb\phi_{-2}^{(1)} -\frac{4}{|q|^2}\pmb\phi_{-2}^{(1)} &=&    \left(-\frac{2}{r^3}+O(mr^{-4})\right)\pmb\phi_{-2}^{(2)} +O(mr^{-3})\nab_{\widehat{\mathcal{X}}_{-2}}^{\leq 1}\pmb\phi_{-2}^{(1)}\\
&&   +O(mr^{-2})\nab_{\widehat{\Z}}^{\leq 1}\pmb\phi_{-2}^{(0)}   +\widetilde{\N}^{(1)}_{W,-2},
\eeaa
with all the coefficients on the RHS  being independent of coordinates $\tau$ and $\tphi$, and with $\widetilde{\N}^{(1)}_{W,-2}$ being given by 
\beaa
\widetilde{\N}^{(1)}_{W,-2} &=& \dk^{\leq 2}\big(\Ga_b\c \Ga_g\big)+r\dk^{\leq 2}(\Ga_b\c\xi).
\eeaa
This concludes the proof of \eqref{eq:TensorialTeuSysandlinearterms:rescaleRHScontaine2:general:Kerrperturbation:alternateformnullframeinsteadcoordvectorfield} in the case $s=-2$ and $p=1$.

%%%%%%%%%%%%%%%%%%%%%%%%%%%%%%%%%%%%%%%%%%%%%%%%%%

\subsubsection{Derivation of the wave equation for $\pmb\phi_{-2}^{(2)}$ in perturbations of Kerr}

%%%%%%%%%%%%%%%%%%%%%%%%%%%%%%%%%%%%%%%%%%%%%%%%%%

Recall from Remark \ref{rmk:compasisionphiminus2p=0withMaSz26andphiminus2p=2withqfGKS22} that $\qfb=\pmb\phi_{-2}^{(2)}$. We may thus use the derivation of the wave equation for $\qfb$ in Theorem 5.3.9 of \cite{GKS22} which yields\footnote{Note that the equation in Theorem 5.3.9 of \cite{GKS22} has two typos concerning the sign in front of the second term on the LHS and the coefficient $q^2$ in front of the second term on the RHS (the conjugate quantity appears in \cite{GKS22}). Note also that we have slightly changed the notation for the lower order terms compared to (5.3.8) of \cite{GKS22} with the following correspondence $q^3\ov{q}\und{W}_4\to \und{W}_4$ and similarly for the other lower order terms.}, in an ingoing frame, 
\bea\lab{eq:maingRWqfbinGKS22}
\bigg(\squared_2 - \frac{4ia\cos\th}{|q|^2}\nab_{\T}- \frac{4}{|q|^2}\bigg)\pmb\phi^{(2)}_{-2} &=& V\pmb\phi^{(2)}_{-2} - q^2\frac{\De}{|q|^2}\frac{8}{r^2}(a^2\nab_{\T}+a\nab_{\Z})\Ab_4 +\und{W}_4\Ab_4+\und{W}_3\nab_3\Ab\nn\\
&&+\und{W}\c\nab\Ab+\und{W}_0\Ab+ \err[\squared_2 \qfb],
\eea
where:
\begin{itemize}
\item $\Ab_4\in\sk_2(\mathbb{C})$ is given by 
\bea
\Ab_4:= \nabc_4\Ab+\frac{1}{2}\tr X\Ab.
\eea

\item The approximate Killing vectorfields $\T$ and $\Z$ are given by \eqref{eq:definitionofTandPhithataretheapproximateKillingvectorifeldinKerrpert}.

\item $V$, $\und{W}_4$, $\und{W}_3$ and $\und{W}_0$ are complex functions of $(r, \th)$, and $\und{W}$ is the product of a complex function of $(r,\th)$ with $\Jk$, with
\beaa
V=O(mr^{-3}), \qquad \Im(V)=0, \qquad \und{W}_4, \, \und{W}=O(a^2r^{-1}), \qquad \und{W}_3,\, \und{W}_0=O(a^2r^{-2}). 
\eeaa

\item $\err[\squared_2 \qfb]$ is  the nonlinear correction term, which under the additional conditions\footnote{In fact, it suffices to assume that $\Xi\in r^{-2}\Ga_g$ and $\Hbc\in r^{-1}\Ga_g$. These additional conditions make the structure of $\err[\squared_2 \qfb]$ in \eqref{eq:MaiThmParq-err-bar} possible and hold true in view of our assumption \eqref{eq:Xi-Hb-chapter12:0}.}  
\beaa
\Xi=0, \qquad \Hbc=0, \quad\textrm{for}\quad r\geq r_0,
\eeaa
is given schematically by  the expression
 \bea\lab{eq:MaiThmParq-err-bar}
\err[{\squared_2} \qfb]&=&  r^2 \dk^{\leq 2}(\Ga_b \c (A, B))+ \dk^{\leq 3} (\Ga_g \c \Ga_b).
 \eea
\end{itemize}

Next, using \eqref{eq:actionofTandPhioncoordinatesrandcosth},  we have
\beaa
q^2\frac{\De}{|q|^2}\frac{8}{r^2}(a^2\nab_{\T}+a\nab_{\Z})\Ab_4 &=& \frac{8}{r^2}(a^2\nab_{\T}+a\nab_{\Z})\left(q^2\frac{\De}{|q|^2}\Ab_4\right)+\Ga_g\Ab_4.
\eeaa
Also, we have together with \eqref{eq:definitionofthephiminus2phierarchy:perturbationofKerr}, 
\beaa
\frac{8}{r^2}(a^2\nab_{\T}+a\nab_{\Z})\pmb\phi_{-2}^{(1)} &=& \frac{8}{r^2}(a^2\nab_{\T}+a\nab_{\Z})\left(q^2\frac{\De}{|q|^2}\Ab_4\right)\\
&&+\frac{8q^2}{r^2}\frac{\De}{|q|^2}\left(-\frac{3}{2}\left(\frac{\atrch}{\trch}\right)_{\vartheta}\atrch -\frac{3}{2}i\atrch\right)(a^2\nab_{\T}+a\nab_{\Z})\Ab+\dk^{\leq 1}(\Ga_g)\Ab\\
&=& \frac{8}{r^2}(a^2\nab_{\T}+a\nab_{\Z})\left(q^2\frac{\De}{|q|^2}\Ab_4\right)\\
&& -\frac{24q^2a\cos\th}{r^2|q|^2}\left(\frac{\De}{|q|^2}\right)^2\left(\frac{a\cos\th}{r}+i\right)(a^2\nab_{\T}+a\nab_{\Z})\Ab +\dk^{\leq 1}(\Ga_g\Ab).
\eeaa
Plugging in \eqref{eq:maingRWqfbinGKS22}, we infer
\begin{align*}
\bigg(\squared_2 - \frac{4ia\cos\th}{|q|^2}\nab_{\T}- \frac{4}{|q|^2}\bigg)\pmb\phi^{(2)}_{-2} =& V\pmb\phi^{(2)}_{-2} -\frac{8}{r^2}(a^2\nab_{\T}+a\nab_{\Z})\pmb\phi_{-2}^{(1)}\\
& -\frac{24q^2a\cos\th}{r^2|q|^2}\left(\frac{\De}{|q|^2}\right)^2\left(\frac{a\cos\th}{r}+i\right)(a^2\nab_{\T}+a\nab_{\Z})\Ab\\
&  +\und{W}_4\Ab_4+\und{W}_3\nab_3\Ab +\und{W}\c\nab\Ab+\und{W}_0\Ab + \err[\squared_2 \qfb] +\Ga_g\Ab_4 +\dk^{\leq 1}(\Ga_g\Ab)
\end{align*}
and hence
\beaa
\bigg(\squared_2 - \frac{4ia\cos\th}{|q|^2}\nab_{\T}- \frac{4}{|q|^2}\bigg)\pmb\phi^{(2)}_{-2} &=& V\pmb\phi^{(2)}_{-2} -\frac{8}{r^2}(a^2\nab_{\T}+a\nab_{\Z})\pmb\phi_{-2}^{(1)}\\
&& +\left(\und{W}_3 -\frac{24q^2a\cos\th}{r^2|q|^2}\left(\frac{\De}{|q|^2}\right)^2\left(\frac{a\cos\th}{r}+i\right)\frac{a^2(\cos\th)^2\De}{2|q|^2}\right)\nab_3\Ab\\
&& +\left(\und{W}_4 -\frac{24q^2a\cos\th}{r^2|q|^2}\left(\frac{\De}{|q|^2}\right)^2\left(\frac{a\cos\th}{r}+i\right)\frac{a^2(\cos\th)^2}{2}\right)\Ab_4\\
&&   +\left(\und{W} -\frac{24q^2a\cos\th}{r^2|q|^2}\left(\frac{\De}{|q|^2}\right)^2\left(\frac{a\cos\th}{r}+i\right)ar^2\Re(\Jk)\right)\c\nab\Ab\\
&&+\und{W}_0\Ab +O(a^3r^{-3})\Ab + \err[\squared_2 \qfb] +\Ga_g\Ab_4 +\dk^{\leq 1}(\Ga_g\Ab),
\eeaa
where we used the fact that $\Ab_4=\nab_4\Ab+O(r^{-1})\Ab+\Ga_g\Ab$. We simplify the above as
\beaa
\bigg(\squared_2 - \frac{4ia\cos\th}{|q|^2}\nab_{\T}- \frac{4}{|q|^2}\bigg)\pmb\phi^{(2)}_{-2} &=& V\pmb\phi^{(2)}_{-2} -\frac{8}{r^2}(a^2\nab_{\T}+a\nab_{\Z})\pmb\phi_{-2}^{(1)}\\
&& +\left(\und{W}_3 -\frac{12q^2a^3(\cos\th)^3}{r^2|q|^2}\left(\frac{\De}{|q|^2}\right)^3\left(\frac{a\cos\th}{r}+i\right)\right)\nab_3\Ab\\
&& +\left(\und{W}_4 -\frac{12q^2a^3(\cos\th)^3}{r^2|q|^2}\left(\frac{\De}{|q|^2}\right)^2\left(\frac{a\cos\th}{r}+i\right)\right)\Ab_4\\
&&   +\left(\und{W} -\frac{24q^2a^2\cos\th}{|q|^2}\left(\frac{\De}{|q|^2}\right)^2\left(\frac{a\cos\th}{r}+i\right)\Re(\Jk)\right)\c\nab\Ab\\
&&+\und{W}_0\Ab +O(a^3r^{-3})\Ab + \err[\squared_2 \qfb] +\Ga_g\Ab_4 +\dk^{\leq 1}(\Ga_g\Ab).
\eeaa
Next, we check that the values of $\und{W}_3$, $\und{W}_4$ and $\und{W}$ in Appendix D7 of \cite{GKS22} verify 
\beaa
\und{W}_3 = \frac{12q^2a^3(\cos\th)^3}{r^2|q|^2}\left(\frac{\De}{|q|^2}\right)^3\left(\frac{a\cos\th}{r}+i\right),\qquad \und{W} = \frac{24q^2a^2\cos\th}{|q|^2}\left(\frac{\De}{|q|^2}\right)^2\left(\frac{a\cos\th}{r}+i\right)\Re(\Jk),
\eeaa
and 
\beaa
\und{W}_4 -\frac{12q^2a^3(\cos\th)^3}{r^2|q|^2}\left(\frac{\De}{|q|^2}\right)^2\left(\frac{a\cos\th}{r}+i\right) &=& q^2\frac{\De}{|q|^2}\left(\frac{16ai}{r^2}+O(a^3r^{-4})\right)
\eeaa
and hence
\beaa
\bigg(\squared_2 - \frac{4ia\cos\th}{|q|^2}\nab_{\T}- \frac{4}{|q|^2}\bigg)\pmb\phi^{(2)}_{-2} &=& V\pmb\phi^{(2)}_{-2} -\frac{8}{r^2}(a^2\nab_{\T}+a\nab_{\Z})\pmb\phi_{-2}^{(1)} +q^2\frac{\De}{|q|^2}\left(\frac{16ai}{r^2}+O(a^3r^{-4})\right)\Ab_4\\
&&+\und{W}_0\Ab +O(a^3r^{-3})\Ab + \err[\squared_2 \qfb] +\Ga_g\Ab_4 +\dk^{\leq 1}(\Ga_g\Ab).
\eeaa
Now, in view of \eqref{eq:definitionofthephiminus2phierarchy:perturbationofKerr}, we have
\beaa
\left(\frac{16ai}{r^2}+O(a^3r^{-4})\right)\pmb\phi_{-2}^{(1)} &=& q^2\frac{\De}{|q|^2}\left(\frac{16ai}{r^2}+O(a^3r^{-4})\right)\Ab_4 +O(a^2r^{-2})\Ab+\Ga_g\Ab
\eeaa
which yields 
\beaa
\bigg(\squared_2 - \frac{4ia\cos\th}{|q|^2}\nab_{\T}- \frac{4}{|q|^2}\bigg)\pmb\phi^{(2)}_{-2} &=& V\pmb\phi^{(2)}_{-2} -\frac{8}{r^2}(a^2\nab_{\T}+a\nab_{\Z})\pmb\phi_{-2}^{(1)} +\left(\frac{16ai}{r^2}+O(a^3r^{-4})\right)\pmb\phi_{-2}^{(1)}\\
&&+\und{W}_0\Ab +O(a^2r^{-2})\Ab + \err[\squared_2 \qfb] +\Ga_g\Ab_4 +\dk^{\leq 1}(\Ga_g\Ab).
\eeaa
and hence
\beaa
\bigg(\squared_2 - \frac{4ia\cos\th}{|q|^2}\nab_{\T}- \frac{4}{|q|^2}\bigg)\pmb\phi^{(2)}_{-2} &=& V\pmb\phi^{(2)}_{-2} -\frac{8}{r^2}(a^2\nab_{\T}+a\nab_{\Z})\pmb\phi_{-2}^{(1)} +\left(\frac{16ai}{r^2}+O(a^3r^{-4})\right)\pmb\phi_{-2}^{(1)}\nn\\
&&+O(a^2r^{-2})\pmb\phi_{-2}^{(0)} + \err[\squared_2 \qfb] +\Ga_g\Ab_4 +\dk^{\leq 1}(\Ga_g\Ab).
\eeaa
which we write in simpler form as follows 
\beaa
\bigg(\squared_2 -\frac{4ia\cos\th}{|q|^2}\nab_{\T}- \frac{4}{|q|^2}\bigg)\pmb\phi^{(2)}_{-2} &=& O(mr^{-3})\pmb\phi^{(2)}_{-2}+O(mr^{-2})\nab_{\Z+a\T}^{\leq 1}\pmb\phi^{(1)}_{-2}+O(m^2 r^{-2})\pmb\phi^{(0)}_{-2} +\widetilde{\N}^{(2)}_{W,-2},
\eeaa
with all the coefficients on the RHS  being independent of coordinates $\tau$ and $\tphi$, with the coefficients in front of the terms $\pmb\phi^{(2)}_{-2}$ and $\nab_{\Z+a\T}\pmb\phi^{(1)}_{-2}$ on the RHS being real functions, and with $\widetilde{\N}^{(2)}_{W,-2}$ being given by 
\beaa
\widetilde{\N}^{(2)}_{W,-2}=  \err[\squared_2 \qfb] +\Ga_g\Ab_4 +\dk^{\leq 1}(\Ga_g\Ab)
\eeaa
so that, in view of  \eqref{eq:MaiThmParq-err-bar}, we have
\beaa
 \widetilde{\N}^{(2)}_{W,+2}  = r^2 \dk^{\leq 2}(\Ga_b \c (A, B))+ \dk^{\leq 3} (\Ga_g \c \Ga_b),
 \eeaa 
which concludes the proof of \eqref{eq:TensorialTeuSysandlinearterms:rescaleRHScontaine2:general:Kerrperturbation:alternateformnullframeinsteadcoordvectorfield} in the case $s=-2$ and $p=2$, and hence the one of Theorem \ref{thm:derivationoftheTeukolskytensorialwavesystemfors=plusminus2:kerrpert:alternateformnullframeinsteadcoordvectorfield}.

%%%%%%%%%%%%%%%%%%%%%%%%%%%%%%%%%%%%%%%%%%%%%%%%%%

\section{Energy-Morawetz estimates for Teukolsky equations in perturbations of Kerr}
\lab{sec:energyMorawetzesitmatesforTeukoslkyonMM:upto15derivatives}

%%%%%%%%%%%%%%%%%%%%%%%%%%%%%%%%%%%%%%%%%%%%%%%%%%

In this section, we rely on \cite{MaSz26} to derive energy-Morawetz estimates for up to 14 derivatives for solutions to Teukolsky equation. The proofs in \cite{MaSz26} need to be slightly adapted given that the spacetime in \cite{MaSz26} extends to $\II_+$ while $\MM$ only extends to the spacelike hypersuface $\Si_*$. We start by stating the main properties of $(\MM, \g)$ needed for Section \ref{sec:energyMorawetzesitmatesforTeukoslkyonMM:upto15derivatives}.

%%%%%%%%%%%%%%%%%%%%%%%%%%%%%%%%%%%%%%%%%

\subsection{Regular scalarization of tensorial wave equations in \cite{MaSz26}}
\lab{sec:regularscalarizationinMaSz26}
  
%%%%%%%%%%%%%%%%%%%%%%%%%%%%%%%%%%%%%%%%%

In this section, we recall the regular scalarization procedure for tensorial wave equations introduced in \cite{MaSz26}.

%%%%%%%%%%%%%%%%%%%%%%%%%%%%%%%%%%%%%%%%%%

\subsubsection{Scalarization using a regular triplet $\Om_i$, $i=1,2,3$}

%%%%%%%%%%%%%%%%%%%%%%%%%%%%%%%%%%%%%%%%%%

In order to scalarize horizontal tensors, we will rely on the following definition, see Definition 3.1 in \cite{MaSz26}. 

\begin{definition}[Regular triplet]
\lab{def:definitionofregulartripletOmii=123}
Let $(\MM, \g)$ be a spacetime, $(e_3, e_4)$ be a null pair, and consider the corresponding horizontal structure $\O(\MM)$ introduced in Section \ref{subsection:review-horiz.structures}. We say that vectorfields $\Om_i$, $i=1,2,3$, identified with elements of $\sk_1$, form a regular triplet if they are regular and satisfy the following identities 
\bea\lab{eq:fundamentalpropertiesof1formsOmi}
x^i\Om_i=0, \qquad (\Om^i)_a(\Om_i)_b=\de_{ab}, \qquad \Om_i\c\Om_j=\de_{ij}-x^ix^j, \qquad \Om_i\c\dual\Om_j=\in_{ijk}x^k,
\eea
where, by convention, we denote $\Om^i=\Om_i$, i.e., the $i$-index is lowered or raised using $\de_{ij}$ or $\de^{ij}$.
\end{definition}

\begin{remark}\lab{rmk:generalcontructionofregulartripletsingivenspacetime}
For a specific choice of a regular triplet in Kerr, see Definition \ref{def:regulartripletinKerrOmii=123}.
\end{remark}

Next, we introduce the following 1-forms on $\MM$.
\begin{definition}
\lab{def:Mialphaj:Kerr}
Let $(\MM, \g)$ be a spacetime, $(e_3, e_4)$ be a null pair, and consider the corresponding horizontal structure $\O(\MM)$ introduced in Section \ref{subsection:review-horiz.structures}. Let $\Om_i$, $1,2,3$ be a regular triplet in the sense of Definition \ref{def:definitionofregulartripletOmii=123}. We  define the following 1-forms on $\MM$ 
\bea\lab{eq:definitionofMalphaijwithoutambiguity}
M_{i\a}^j:=(\Ddot_\a\Om_i)\c\Om^j, \quad \forall \a,i,j.
\eea
Further, we define $M_{i}^{j \a}:=\g^{\a\b} M_{i\b}^j$.
\end{definition}

\begin{lemma}\lab{lemma:introductionandpropertiesoftheMalphaij}
Let $M_{i\a}^j$ be the 1-forms on $(\MM, \g)$ as defined in Definition \ref{def:Mialphaj:Kerr}.  Then we have
\bea
\label{def:Mialphaj}
\Ddot_\a\Om_i=M_{i\a}^j\Om_j.
\eea
\end{lemma}

\begin{proof}
See Lemma 3.4 in \cite{MaSz26}.
\end{proof}

The following lemma provides a useful property for the symmetric part of $M_{i\a}^j$.
\begin{lemma}
\label{lem:property:Omi}
The symmetric part $(M_{S})_{i\a}^{j}$ of $M_{i\a}^{j}$ w.r.t. $(i,j)$ satisfies 
\bea
\lab{formula:symmetricpartofMmatrices}
(M_{S})_{i\a}^{j}:=\frac{1}{2}(M_{i\a}^{j}+M_{j\a}^{i}), \qquad (M_{S})_{i\a}^{j}=-\frac{1}{2}\pr_\a(x^ix^j).
\eea
\end{lemma}

\begin{proof}
See Lemma 3.5 in \cite{MaSz26}.
\end{proof}

%%%%%%%%%%%%%%%%%%%%%%%%%%%%%%%%

\subsubsection{From tensors to regular scalars and back}

%%%%%%%%%%%%%%%%%%%%%%%%%%%%%%%%

The following lemma allows to pass from horizontal tensors in $\sk_2(\mathbb{C})$ to scalars and reciprocally.
\begin{lemma}\lab{lemma:backandforthbetweenhorizontaltensorsk2andscalars:complex}
Let $(\MM, \g)$ be a spacetime, $(e_3, e_4)$ be a null pair, and consider the corresponding horizontal structure $\O(\MM)$ introduced in Section \ref{subsection:review-horiz.structures}. Assume that $\Om_i$, $1,2,3$ is a regular triplet in the sense of Definition \ref{def:definitionofregulartripletOmii=123}. Then, the following holds:
\begin{enumerate}
\item\lab{item1:sk2Csatisfyconditions} Let $\pmb\psi\in\sk_2(\mathbb{C})$ and define the complex-valued scalars $\psi_{ij}:=\pmb\psi(\Om_i, \Om_j)$, $i,j=1,2,3$. Then:
\begin{itemize}
\item The complex-valued scalars $\psi_{ij}$ satisfy 
\bea\lab{eq:fundamentalidentitiestoderivefromscalarizationoftensor:complexcase}
\psi_{ij}=\psi_{ji}, \qquad x^i\psi_{ij}=0, \qquad (\de^{ij}-x^ix^j)\psi_{ij}=0, \qquad \in_{ikl}x^l\psi_{kj}+i\psi_{ij}=0.
\eea
\item We may recover the tensor $\pmb\psi$ from the scalars $\psi_{ij}$ by the formula 
\beaa
\pmb\psi_{ab}=\psi_{ij}(\Om^i)_a(\Om^j)_b.
\eeaa
\end{itemize}
\item Reciprocally, let $\psi_{ij}$ be complex-valued scalars satisfying the identities \eqref{eq:fundamentalidentitiestoderivefromscalarizationoftensor:complexcase}, and introduce the complex-valued horizontal 2-tensor $\pmb\psi$ by $\pmb\psi_{ab}:=\psi_{ij}(\Om^i)_a(\Om^j)_b$, $a,b=1,2$. Then, we have $\pmb\psi\in\sk_2(\mathbb{C})$ and $\pmb\psi(\Om_i, \Om_j)=\psi_{ij}$ for all $i,j=1,2,3$.
\end{enumerate}
\end{lemma}

\begin{proof}
See Lemma 3.8 in \cite{MaSz26}.
\end{proof}

%%%%%%%%%%%%%%%%%%%%%%%%%%%%%%%%%%%%%%%%

\subsubsection{Scalarization of the tensorial wave operator $\squared_2$}

%%%%%%%%%%%%%%%%%%%%%%%%%%%%%%%%%%%%%%%%

The following lemma provides the scalarization of the tensorial wave operator $\squared_2$.
\begin{lemma}\lab{lemma:formoffirstordertermsinscalarazationtensorialwaveeq}
Let $(\MM, \g)$ be a spacetime, $(e_3, e_4)$ be a null pair, and consider the corresponding horizontal structure $\O(\MM)$ introduced in Section \ref{subsection:review-horiz.structures}. Assume that $\Om_i$, $1,2,3$ is a regular triplet in the sense of Definition \ref{def:definitionofregulartripletOmii=123}. Also, let $\pmb\psi\in\sk_2(\mathbb{C})$ and let $\psi_{ij}$ be the scalars associated to it in view of Lemma \ref{lemma:backandforthbetweenhorizontaltensorsk2andscalars:complex}. Then, we have
\bea
\squared_2\pmb\psi(\Om_i, \Om_j) &=& \square_\g(\psi_{ij}) -S(\psi)_{ij} - (Q\psi)_{ij}
\eea
where 
\bsub
\label{SandV}
\begin{align}
S(\psi)_{ij} ={}& 2M_{i}^{k\a}\pr_\a(\psi_{kj}) +2M_{j}^{k\a}\pr_\a(\psi_{ik}),\\
(Q\psi)_{ij} ={}& (\Ddot^\a M_{i\a}^k)\psi_{kj}+(\Ddot^\a M_{j\a}^k)\psi_{ik} -M_{i\a}^kM_k^{l\a}\psi_{lj}-2M_{i\a}^kM_{j}^{l\a}\psi_{kl}-M_{j\a}^kM_k^{l\a}\psi_{il},
\end{align}
\esub
with the 1-forms $M_{i\a}^j$ defined by \eqref{eq:definitionofMalphaijwithoutambiguity}.
\end{lemma}

\begin{proof}
See Lemma 3.9 in \cite{MaSz26}.
\end{proof}

%%%%%%%%%%%%%%%%%%%%%%

\subsubsection{Tensorization defect}

%%%%%%%%%%%%%%%%%%%%%%

We now consider more general families of scalars $\psi_{ij}$ that are not necessarily derived from the scalarization of a tensor in $\sk_2(\mathbb{C})$, and we aim at estimating their tensorization defect which is defined as follows. 

\begin{definition}[Tensorization defect]
\lab{def:definitionofthenotationerrforthescalarizationdefect}
For a general family of complex-valued scalars $\psi_{ij}$, define the following error term which estimates the corresponding tensorization defect 
\bea\lab{eq:definitionofthenotationerrforthescalarizationdefect}
\bsplit
\err_{\textrm{TDefect}}[\psi]:=&\Big(\err_{\textrm{TDefect},n}[\psi],\,n=1,2,3,4,5\Big),\\
(\err_{\textrm{TDefect},1}[\psi])_{ij}:=&\psi_{ij}-\psi_{ji},\qquad (\err_{\textrm{TDefect},2}[\psi])_j:=x^i\psi_{ij}, \\ 
(\err_{\textrm{TDefect},3}[\psi])_j:=&x^i\psi_{ji},\qquad \err_{\textrm{TDefect},4}[\psi]:=(\de^{ij}-x^ix^j)\psi_{ij},\\
(\err_{\textrm{TDefect},5}[\psi])_{ij}:=& \in_{ikl}x^l\psi_{kj}+i\psi_{ij}.
\end{split}
\eea
\end{definition}

\begin{remark}
In view of Lemma \ref{lemma:backandforthbetweenhorizontaltensorsk2andscalars:complex}, a family of complex-valued scalars $\psi_{ij}$ comes from the scalarization of a tensor in $\sk_2(\mathbb{C})$ if and only if $\err_{\textrm{TDefect}}[\psi]=0$.
\end{remark}

We also approximate general families of scalars by families that are generated by scalarization of a tensor in $\sk_2(\mathbb{C})$. 
\begin{lemma}\lab{lemma:computationerrorscalarizationdeffect}
Let $\psi_{ij}$ be a general family of complex-valued scalars and let $\Pi_2[\psi]$ be defined by
\bea\lab{eq:computationerrorscalarizationdeffect:defPi2}
(\Pi_2[\psi])_{ij}:=\frac{1}{2}\Big((\pi_2[\psi])_{ij}+i\in_{ikl}x^l(\pi_2[\psi])_{kj}\Big),
\eea
where 
\bea
\nn (\pi_2[\psi])_{ij} &:=& \frac{1}{2}(\psi_{ij}+\psi_{ji}) -\frac{1}{2}\Big(x^k(\psi_{kj}+\psi_{jk})x^i+x^k(\psi_{ki}+\psi_{ik})x^j\Big) +x^kx^l\psi_{kl}x^ix^j\\
&&-\frac{1}{2}(\de^{kl}-x^kx^l)\psi_{kl}(\de^{ij}-x^ix^j).
\eea
Then, $\Pi_2[\psi]$ is generated by the scalarization of a tensor in $\sk_2(\mathbb{C})$, i.e., we have 
\beaa
(\Pi_2[\psi])_{ij}=\pmb\Pi_2[\psi](\Om_i, \Om_j), \quad i,j=1,2,3, \quad \pmb\Pi_2[\psi]\in\sk_2(\mathbb{C}).
\eeaa
Furthermore, we have
\beaa
\psi_{ij}=(\Pi_2[\psi])_{ij} +(\widetilde{\err}_0[\psi])_{ij}
\eeaa
where 
\bea
\bsplit
(\widetilde{\err}_0[\psi])_{ij}:=& -\frac{i}{2}(\err_{\textrm{TDefect},5}[\psi])_{ij}+\frac{1}{2}\Big((\err_0[\psi])_{ij}+i\in_{ikl}x^l(\err_0[\psi])_{kj}\Big),\\
(\err_0[\psi])_{ij} :=& \frac{1}{2}(\err_{\textrm{TDefect},1}[\psi])_{ij} +\frac{1}{2}\Big((\err_{\textrm{TDefect},2}[\psi])_j+(\err_{\textrm{TDefect},3}[\psi])_j\Big)x^i\\
& +\frac{1}{2}\Big((\err_{\textrm{TDefect},2}[\psi])_i+(\err_{\textrm{TDefect},3}[\psi])_i\Big)x^j\\
& -x^k(\err_{\textrm{TDefect},2}[\psi])_kx^ix^j +\frac{1}{2}\err_{\textrm{TDefect},4}[\psi](\de^{ij}-x^ix^j).
\end{split}
\eea
\end{lemma}

\begin{proof}
See Lemma 3.13 in \cite{MaSz26}.
\end{proof}

Finally, we derive wave equations for the tensorization defect. 
\begin{lemma}\lab{lemma:waveequationsfortensordeffects}
Let $\psi_{ij}$ be a family of complex-valued scalars satisfying 
\beaa
\big(\square_\g+V\big)\psi_{ij} &=& S(\psi)_{ij} + (Q\psi)_{ij},
\eeaa
where $S$ and $Q$ are defined in \eqref{SandV} and where $V$ is a real-valued function, and let $\err_{\textrm{TDefect}}[\psi]$ be associated to the family $\psi_{ij}$ as in \eqref{eq:definitionofthenotationerrforthescalarizationdefect}. Then, $\err_{\textrm{TDefect},1}[\psi]$ satisfies 
\bea\lab{eq:waveeqpsiijminuspsiji}
(\square_\g+V)\big((\err_{\textrm{TDefect},1}[\psi])_{ij}\big) &=& S(\err_{\textrm{TDefect},1}[\psi])_{ij} + (Q\err_{\textrm{TDefect},1}[\psi])_{ij},
\eea
$\err_{\textrm{TDefect},2}[\psi]$ satisfies 
\bea\lab{eq:waveeqxipsiij}
\nn(\square_\g+V)((\err_{\textrm{TDefect},2}[\psi])_j) &=& 2M_j^{k\a}\pr_\a((\err_{\textrm{TDefect},2}[\psi])_k) +(\Ddot^\a M_{j\a}^k)(\err_{\textrm{TDefect},2}[\psi])_k\\
&& -M_{j\a}^kM_k^{l\a}(\err_{\textrm{TDefect},2}[\psi])_l,
\eea
$\err_{\textrm{TDefect},4}[\psi]$ satisfies 
\bea\lab{eq:waveeqfortraceofpsi:deltaijpsiij}
\nn &&(\square_\g+V)(\err_{\textrm{TDefect},4}[\psi])\\ 
\nn &=& -2\pr^{\a} (x^i)\pr_{\a}((\err_{\textrm{TDefect},2}[\psi])_i)  -2\pr^{\a} (x^i)\pr_{\a}((\err_{\textrm{TDefect},3}[\psi])_i)\\
&& +\Big(-2\square_\g(x^i) - x^j(\Ddot^\a M_{j\a}^i) -\pr_\a(x^k)M_k^{i\a}\Big)(\err_{\textrm{TDefect},2}[\psi])_i\nn\\
&& +\Big(- 2\square_\g(x^i) -  x^j(\Ddot^\a M_{j\a}^i)  + \pr_\a(x^k)M_k^{i\a}\Big)(\err_{\textrm{TDefect},3}[\psi])_i
\eea
and $\err_{\textrm{TDefect},5}[\psi]$ satisfies 
\bea\lab{eq:waveeqforlastdefecttensor:antiselfdual}
\nn&& (\square_\g+V)((\err_{\textrm{TDefect},5}[\psi])_{ij})\\ 
\nn&=& S(\err_{\textrm{TDefect},5}[\psi])_{ij} +(Q\err_{\textrm{TDefect},5}[\psi])_{ij} + 2\in_{ikl}x^k\g^{\a\b}\pr_\a(x^l)\pr_\b((\err_{\textrm{TDefect},2}[\psi])_j)\\ 
\nn&& +\Big(\in_{ikl}x^k\square_\g(x^l)  +\in_{ikl}\g^{\a\b}\pr_\a(x^l)\pr_\b(x^k) -\in_{knl}x^nM_{i\a}^k\pr^\a(x^l)\Big)(\err_{\textrm{TDefect},2}[\psi])_j\\
&&  -2\in_{ikl}x^kM_{j}^{m\a}\pr_\a(x^l)(\err_{\textrm{TDefect},2}[\psi])_m.
\eea
\end{lemma}

\begin{proof}
See Lemma 3.14 in \cite{MaSz26}.
\end{proof}

%%%%%%%%%%%%%%%%%%%%%%%%%%%%%%%%%%%%%%%%%%%%%%%%%%%%%%%%%%%%%%%%%%%%%%%%%%%%%%%%%%%%%%%%%%%%%%%%

\subsubsection{Differentiation with respect to $\pr_\tau$ and $\widehat{\pr}_{\tphi}$ preserving identities \eqref{eq:fundamentalidentitiestoderivefromscalarizationoftensor:complexcase}}

%%%%%%%%%%%%%%%%%%%%%%%%%%%%%%%%%%%%%%%%%%%%%%%%%%%%%%%%%%%%%%%%%%%%%%%%%%%%%%%%%%%%%%%%%%%%%%%%

We start by noticing that differentiation w.r.t. $\pr_\tau$ preserves the identities  \eqref{eq:fundamentalidentitiestoderivefromscalarizationoftensor:complexcase}. 
\begin{lemma}\lab{lemma:differentiatingwrtprtaupreservetheidentitiesscaloftensors}
Let $\psi_{ij}$ be a family of complex-valued scalars satisfying the identities \eqref{eq:fundamentalidentitiestoderivefromscalarizationoftensor:complexcase}. Then, $\pr_\tau(\psi_{ij})$ satisfies the identities \eqref{eq:fundamentalidentitiestoderivefromscalarizationoftensor:complexcase} as well.
\end{lemma}

\begin{proof}
See Lemma 3.15 in \cite{MaSz26}.
\end{proof}

While $\pr_{\tphi}(x^3)=0$, we have $\pr_{\tphi}(x^1)=-x^2$ and $\pr_{\tphi}(x^2)=x^1$. Hence, differentiation w.r.t.  $\pr_{\tphi}$ does not preserve the identities \eqref{eq:fundamentalidentitiestoderivefromscalarizationoftensor:complexcase} and we will instead use the following modification. 

\begin{definition}\lab{def:widehatprtphi}
Let $\widehat{\pr}_{\tphi}$ denote the first-order operator acting on families of complex-valued scalars $\psi_{ij}$ as follows 
\beaa
\widehat{\pr}_{\tphi}(\psi)_{ij}:=\pr_{\tphi}(\psi_{ij}) +\in_{ik3}\psi_{kj} +\in_{jk3}\psi_{ki}.
\eeaa
\end{definition}

\begin{remark}
In Kerr, if $\psi_{ij}=\pmb\psi(\Om_i, \Om_j)$ with $\pmb\psi\in\sk_2$, then we have $\widehat{\pr}_{\tphi}(\psi)_{ij}=\Lieb_{\pr_{\tphi}}\pmb\psi(\Om_i, \Om_j)$, see Lemma \ref{lemma:InKerrlinkbetweenprtauwidehatprtphiandLieb[rtauLiebprtphi}, where  the horizontal Lie derivative $\Lieb$ has been introduced in Definition \ref{definition:hor-Lie-derivative}. This motivates Definition \ref{def:widehatprtphi}.
\end{remark}

The following lemma proves that differentiation w.r.t. $\widehat{\pr}_{\tphi}$ preserves the identities \eqref{eq:fundamentalidentitiestoderivefromscalarizationoftensor:complexcase}. 
\begin{lemma}\lab{lemma:differentiatingwrtwidehatprtphipreservetheidentitiesscaloftensors}
Let $\psi_{ij}$ be a family of complex-valued scalars satisfying the identities \eqref{eq:fundamentalidentitiestoderivefromscalarizationoftensor:complexcase} and let $\widehat{\pr}_{\tphi}$ be as in Definition \ref{def:widehatprtphi}. Then, $\widehat{\pr}_{\tphi}(\psi)_{ij}$ satisfies the identities \eqref{eq:fundamentalidentitiestoderivefromscalarizationoftensor:complexcase} as well.
\end{lemma}

\begin{proof}
See Lemma 3.18 in \cite{MaSz26}.
\end{proof}
 
\begin{remark}
\lab{rem:prtauandprtphihatpreservesk2C}
In view of Lemmas \ref{lemma:backandforthbetweenhorizontaltensorsk2andscalars:complex},  \ref{lemma:differentiatingwrtprtaupreservetheidentitiesscaloftensors} and   \ref{lemma:differentiatingwrtwidehatprtphipreservetheidentitiesscaloftensors}, we immediately infer the fact that 
if $\psi_{ij}=\pmb\psi(\Om_i, \Om_j)$ for $\pmb\psi\in \sk_2(\mathbb{C})$, then for any $k,l\in\mathbb{N}$, there exists $\pmb\psi_{(k,l)}\in\sk_2(\mathbb{C})$ such that $\pr_\tau^k\widehat{\pr}_{\tphi}^l(\psi)_{ij}=\pmb\psi_{(k,l)}(\Om_i,\Om_j)$.
\end{remark}

%%%%%%%%%%%%%%%%%%%%%%%%%%%%%

\subsubsection{Regular triplet in Kerr}

%%%%%%%%%%%%%%%%%%%%%%%%%%%%%

\begin{definition}[Regular triplet in Kerr]
\lab{def:regulartripletinKerrOmii=123}
Let 
\bea\lab{eq:definitionofxii=123usedintheformulaoftheregulartripletsinKerr}
x^1:=\cos\tphi\sin\th, \qquad x^2:=\sin\tphi\sin\th, \qquad x^3:=\cos\th.
\eea
Then, we define the following horizontal vectorfields $\Om^i$ in Kerr by 
\bea
\Om^i :=|q|\dual\nab(x^i), \quad i=1,2,3.
\eea
\end{definition}

\begin{lemma}\lab{lemma:fundamentalpropertiesof1formsOmi}
The horizontal vectorfields $\Om_i$ in Kerr introduced in Definition \ref{def:regulartripletinKerrOmii=123} satisfy \eqref{eq:fundamentalpropertiesof1formsOmi}. In particular, they form a regular triplet in Kerr in the sense of Definition \ref{def:definitionofregulartripletOmii=123}.
\end{lemma}

\begin{proof}
See Lemma 3.21 in \cite{MaSz26}.
\end{proof}

We also derive the following properties of the 1-forms $M_{i\a}^j$ in Kerr.
\begin{lemma}\lab{lemma:computationoftheMialphajinKerr}
Let $\Om_i$, $i=1,2,3$, be the regular triplet in Kerr of Definition \ref{def:regulartripletinKerrOmii=123}, and let $M_{i\a}^j$ be the corresponding 1-forms in Kerr given by \eqref{eq:definitionofMalphaijwithoutambiguity}. Then, we have, for $i,j=1,2,3$,  
\beaa
M_{i3}^j &=& \phimod'(r)\in_{ki3}(\de^{kj}-x^kx^j)+\frac{a\cos\th}{|q|^2}\in_{ijk}x^k,\\
M_{i4}^j &=& \frac{2a -\De\phimod'(r)}{|q|^2}\in_{ki3}(\de^{kj}-x^kx^j)+\frac{a\cos\th\De}{|q|^4}\in_{ijk}x^k,\\
M_{i\a}^j(\pr_\tau)^\a &=& -\frac{2amr\cos\th}{|q|^4}\in_{ijk}x^k,
\eeaa
and the following asymptotic holds, for $r$ large and $i,j=1,2,3$,
\beaa
M_{ia}^j=O(r^{-1}), \quad a=1,2.
\eeaa
Also, we have for $r\in [r_+(1+2\dbl),  12m]$
\beaa
M_{i\a}^j(\pr_r)^\a=0,\qquad \g^{r\a}M_{i\a}^j=0, \quad i,j=1,2,3.
\eeaa
\end{lemma}

\begin{proof}
See Lemma 3.22 in \cite{MaSz26}.
\end{proof}

\begin{lemma}\lab{lemma:InKerrlinkbetweenprtauwidehatprtphiandLieb[rtauLiebprtphi}
In Kerr, if $\psi_{ij}=\pmb\psi(\Om_i, \Om_j)$ with $\pmb\psi\in\sk_2$, then 
\beaa
\Lieb_{\pr_\tau}\Om_i=0, \qquad \Lieb_{\pr_{\tphi}}\Om_i = -\in_{ij3}\Om_j,
\eeaa
and
\beaa
\pr_\tau(\psi_{ij})=\Lieb_{\pr_\tau}\pmb\psi(\Om_i, \Om_j), \qquad \widehat{\pr}_{\tphi}(\psi)_{ij}=\Lieb_{\pr_{\tphi}}\pmb\psi(\Om_i, \Om_j),
\eeaa
where the horizontal Lie derivative $\Lieb$ has been introduced in Definition \ref{definition:hor-Lie-derivative}, and where $\widehat{\pr}_{\tphi}$ has been introduced in Definition \ref{def:widehatprtphi}.
\end{lemma}

\begin{proof}
See Lemma 3.23 in \cite{MaSz26}.
\end{proof}

Finally, for convenience, we compute the explicit formula of $M_{ia}^j$, $i,j=1,2,3$, $a=1,2$, in Kerr which is not given by Lemma \ref{lemma:computationoftheMialphajinKerr}.
\begin{lemma}\lab{lemma:explicitcomputationofMiaja=1or2}
Let $\Om_i$, $i=1,2,3$, be the regular triplet in Kerr of Definition \ref{def:regulartripletinKerrOmii=123}, and let $M_{i\a}^j$ be the corresponding 1-forms in Kerr given by \eqref{eq:definitionofMalphaijwithoutambiguity}. Then, we have, for $a=1,2$, and $j=1,2,3$,  
\beaa
M_{1a}^j &=& -\frac{a^2x^3}{|q|^3}(\de_{1j}-x^1x^j)(\dual\Om_3)_a +\frac{r^2}{|q|^3}x^1(\dual\Om_j)_{a}-\frac{2a^2x^3}{|q|^3}x^2(\Om_j)_{a},\\
M_{2a}^j &=& -\frac{a^2x^3}{|q|^3}(\de_{2j}-x^2x^j)(\dual\Om_3)_a +\frac{r^2}{|q|^3}x^2(\dual\Om_j)_{a}+\frac{2a^2x^3}{|q|^3}x^1(\Om_j)_{a},\\
M_{3a}^j &=& -\frac{a^2x^3}{|q|^3}(\de_{3j}-x^3x^j)(\dual\Om_3)_a +\frac{r^2+a^2}{|q|^3}x^3(\dual\Om_j)_{a}.
\eeaa
\end{lemma}

\begin{proof}
We compute 
\beaa
(\nab_a\Om_i)_b &=& \nab_a(|q|\dual\nab_b(x^i)) = \frac{1}{|q|}\nab_a(|q|)(\Om_i)_b+|q|\in_{bc}\nab_a\nab_c(x^i)\\
&=& \frac{1}{2|q|^2}\nab_a(|q|^2)(\Om_i)_b+|q|\in_{bc}\nab_a\nab_c(x^i)\\
&=& \frac{a^2\cos\th}{|q|^2}e_a(\cos\th)(\Om_i)_b+|q|\in_{bc}\nab_a\nab_c(x^i)\\
&=& \frac{a^2x^3}{|q|^2}e_a(x^3)(\Om_i)_b+|q|\in_{bc}\nab_a\nab_c(x^i)\\
&=& -\frac{a^2x^3}{|q|^3}(\dual\Om_3)_a(\Om_i)_b+|q|\in_{bc}\nab_a\nab_c(x^i).
\eeaa
Since we have in view of section 2.4.6 in \cite{KS:Kerr}
\beaa
\nab\hot\nab(x^i)=0, \quad i=1,2,3,
\eeaa
we infer
\beaa
\nab_a\nab_b(x^i) &=& \frac{1}{2}\De(x^i)\de_{ab}+\frac{1}{2}\in_{cd}\nab_c\nab_d(x^i)\in_{ab}+\frac{1}{2}(\nab\hot\nab(x^i))_{ab}\\
&=& \frac{1}{2}\De(x^i)\de_{ab}+\frac{1}{2}\in_{cd}\nab_c\nab_d(x^i)\in_{ab}.
\eeaa
Also, we have in view of (2.1.27) in \cite{GKS22}
\beaa
\in_{cd}\nab_c\nab_d(x^i) &=& \frac{1}{2}(\atrch e_3(x^i)+\frac{1}{2}\atrchb e_4(x^i))\in_{cd}\in^{cd}\\
&=& \frac{2a\cos\th}{|q|^2}\left(\frac{\De}{|q|^2}e_3(x^i)+e_4(x^i)\right)
\eeaa
so that 
\beaa
\nab_a\nab_b(x^i) &=& \frac{1}{2}\De(x^i)\de_{ab}+\frac{a\cos\th}{|q|^2}\left(\frac{\De}{|q|^2}e_3(x^i)+e_4(x^i)\right)\in_{ab},
\eeaa
and hence
\beaa
(\nab_a\Om_i)_b &=& -\frac{a^2x^3}{|q|^3}(\dual\Om_3)_a(\Om_i)_b -\frac{|q|}{2}\De(x^i)\in_{ab}+\frac{ax^3}{|q|}\left(\frac{\De}{|q|^2}e_3(x^i)+e_4(x^i)\right)\de_{ab}.
\eeaa
We deduce 
\beaa
M_{ia}^j &=& \nab_a\Om_i\c\Om_j=(\nab_a\Om_i)_b(\Om_j)^b\\
&=& -\frac{a^2x^3}{|q|^3}(\dual\Om_3)_a\Om_i\c\Om_j-\frac{|q|}{2}\De(x^i)(\dual\Om_j)_{a}+\frac{ax^3}{|q|}\left(\frac{\De}{|q|^2}e_3(x^i)+e_4(x^i)\right)(\Om_j)_{a}\\
&=& -\frac{a^2x^3}{|q|^3}(\de_{ij}-x^ix^j)(\dual\Om_3)_a-\frac{|q|}{2}\De(x^i)(\dual\Om_j)_{a}+\frac{ax^3}{|q|}\left(\frac{\De}{|q|^2}e_3(x^i)+e_4(x^i)\right)(\Om_j)_{a}.
\eeaa
Since 
\beaa
\frac{\De}{|q|^2}e_3(x^1)+e_4(x^1) =-\frac{2a}{|q|^2}x^2, \qquad \frac{\De}{|q|^2}e_3(x^2)+e_4(x^2)=\frac{2a}{|q|^2}x^1, \qquad \frac{\De}{|q|^2}e_3(x^3)+e_4(x^3)=0,
\eeaa
and since we have in view of (2.4.18) in \cite{KS:Kerr}
\beaa
\De(x^j)=-\frac{2r^2}{|q|^4}x^j, \quad j=1,2, \qquad \De(x^3)=-\frac{2(r^2+a^2)}{|q|^4}x^3,
\eeaa
we deduce
\beaa
M_{1a}^j &=& -\frac{a^2x^3}{|q|^3}(\de_{1j}-x^1x^j)(\dual\Om_3)_a +\frac{r^2}{|q|^3}x^1(\dual\Om_j)_{a}-\frac{2a^2x^3}{|q|^3}x^2(\Om_j)_{a},\\
M_{2a}^j &=& -\frac{a^2x^3}{|q|^3}(\de_{2j}-x^2x^j)(\dual\Om_3)_a +\frac{r^2}{|q|^3}x^2(\dual\Om_j)_{a}+\frac{2a^2x^3}{|q|^3}x^1(\Om_j)_{a},\\
M_{3a}^j &=& -\frac{a^2x^3}{|q|^3}(\de_{3j}-x^3x^j)(\dual\Om_3)_a +\frac{r^2+a^2}{|q|^3}x^3(\dual\Om_j)_{a},
\eeaa
as stated. This concludes the proof of Lemma \ref{lemma:explicitcomputationofMiaja=1or2}.
\end{proof}

%%%%%%%%%%%%%%%%%%%%%%%%%%%%%%%%%%%%%%%%%%%%%%%%%%%%%%%%%%%%%%%%

\subsection{Set-up for Section \ref{sec:energyMorawetzesitmatesforTeukoslkyonMM:upto15derivatives}}

%%%%%%%%%%%%%%%%%%%%%%%%%%%%%%%%%%%%%%%%%%%%%%%%%%%%%%%%%%%%%%%%

%%%%%%%%%%%%%%%%%%%%%%%

\subsubsection{Choices of constants}
\lab{sec:smallnesconstants}

%%%%%%%%%%%%%%%%%%%%%%%

The following constants  are  involved in the derivation of the energy-Morawetz estimates for Teukolsky:
\begin{itemize}
\item The constants $m>0$ and $a$, with $|a|<m$, are the mass and the angular momentum per unit mass of the Kerr solution relative to which the perturbation to Kerr is measured. 

\item The size of the perturbation to Kerr is measured by $\ep>0$. 

\item $r_0>0$ is tied to  $\Mint\cap\Mext=\{r=r_0\}$. 

\item The constant $\dhor$ is tied to the boundary of $\MM$ given by $\pr\MM=\AA=\{r=r_+(1-\dhor)\}$. 

\item The constant $\dred$ measures the width of the redshift region.

\item The constant $\dbl$ appears in the construction of normalized coordinates, see Lemma \ref{lem:specificchoice:normalizedcoord}.

\item The constant $\dec$ is tied to decay estimates in $(r, \tau)$ of Kerr perturbations, see Section \ref{sec:assumptionsforsec:energyMorawetzesitmatesforTeukoslkyonMM:upto15derivatives}. 

\item The large integer $\Nmic$ is tied to the choice of a contant\footnote{The contant $\Rmic$ is used to separate the proof of energy-Morawetz estimates in two regions, $\MM_{r\leq\Rmic}$ and $\MM_{r\geq\Rmic}$, see beginning of Section \ref{sect:microlocalenergyMorawetztensorialwaveequation}.} $\Rmic\in [\Nmic m, (\Nmic+1)m]$, see Remark \ref{rmk:choiceofconstantRbymeanvalue}.
\end{itemize}

These  constants are chosen such that 
\bea\lab{eq:constraintsonthemainsmallconstantsepanddelta}
\bsplit
&0<\ep_0,\,\ep\ll\dhor\ll\dred\ll \dbl \ll 1-\frac{|a|}{m}, \qquad \ep_0,\,\ep\ll \de\ll\dec, \\ 
&\ep_0,\,\ep\ll \frac{1}{N_0},\,\frac{1}{r_0},
\end{split}
\eea

From now on, in the rest of the paper, $\lesssim$ means bounded by a positive constant multiple, with this positive constant depending only on universal constants (such as constants arising from Sobolev embeddings, elliptic estimates,...) as well as the constants 
$$m,\,\, a, \,\, \dhor,\,\, \dred,\,\, \dbl, \,\, \dec,\,\, r_0,$$
\textit{but not on} $\ep$. {Also, note that the} constants $\dhor, \dred$ and $\dbl$ can be {chosen} to be only dependent on $m$ and $a$.

%%%%%%%%%%%%%%%%%%%%%%%%%%%%%%%%%%     
   
     \subsubsection{Subregions and hypersurfaces of $\MM$}
     \lab{section:SpacetimeMM-chap6}

%%%%%%%%%%%%%%%%%%%%%%%%%%%%%%%%%%

As in Section \ref{sec:Kerrpertbasic}, we consider a given vacuum spacetime $(\MM, \g)$ together with a null pair $(e_3, e_4)$ and its corresponding horizontal structure as in Section \ref{subsection:review-horiz.structures}. We will use the complexified Ricci and curvature coefficients  of Definition \ref{def:complexRicciandcurvaturecoefficients}. Moreover, we assume that $\MM$ is endowed with a pair of constants $(a, m)$, scalar functions $(\tau, r, \th, \tphi)$ and complex horizontal 1-forms $\Jk$, $\Jk_{\pm}$. 

In addition, we assume the following:
\begin{enumerate}
\item The level sets $\Si(\tau)$ of $\tau$ are spacelike, and $\tau\in[1,\tau_*]$ on $\MM$ for some arbitrary large constant $\tau_*$. 

\item The boundary of $\MM$ is given by 
\bea\lab{eq:prMM=AAcupSistarcupSi1cupSitaustar}
\pr\MM=\AA\cup\Si_*\cup\Si(1)\cup\Si(\tau_*)
\eea
where 
\bea\lab{eq:defAA=r=rplusoneminusedeh}
\AA:=\Big\{r=r_+(1-\deh), \, 1\leq\tau\leq\tau_*\Big\}, 
\eea
and $\Si_*$ is a spacelike hypersurface such that 
\bea\lab{eq:rangetauandronSigmastar}
1\leq\tau\leq\tau_*,\qquad r\geq r_*, \qquad \tau+r=c_*+h_*(r), \qquad h_*(r)=O(m^2r^{-1}),\qquad\textrm{on}\quad\Si_*,
\eea
with $r_*$, $c_*$ constants, and with $r_*$ satisfying the following dominance condition
\bea\lab{eq:dominantconditionforrstarcomparedtotaustartonS*}
r_*\simeq \ep_0^{-1}\tau_*^{1+\dec}.
\eea

\item Let $r_0$ a large enough fixed constant. We decompose $\MM$ as follows
\bea\lab{eq:definitionofMextandMint}
\Mint:=\MM\cap\{r\leq r_0\}, \qquad \Mext:=\MM\cap\{r\geq r_0\}. 
\eea
\end{enumerate} 

\begin{remark}
In view of \eqref{eq:rangetauandronSigmastar} and \eqref{eq:dominantconditionforrstarcomparedtotaustartonS*}, we have on $\Si_*$
\bea\lab{eq:controlofsizeofc*assimeqr*:onSi*}
c_*=\big(1+O(\ep_0)\big)r_*,
\eea
as well as 
\beaa
|r-r_*|\les |\tau-\tau_*|+\frac{m^2}{r_*}\les \tau_*+\frac{m^2}{r_*}\les \ep_0r_*+\frac{m^2}{r_*}\les \ep_0r_*\quad\Longrightarrow\quad r\simeq r_*
\eeaa
so that \eqref{eq:dominantconditionforrstarcomparedtotaustartonS*} in fact holds for any $r$ on $\Si_*$, i.e., 
\bea\lab{eq:dominantconditionforrstarcomparedtotaustartonS*:holdsforallronSi*}
r\simeq \ep_0^{-1}\tau_*^{1+\dec}\quad\textrm{on}\quad\Si_*.
\eea
\end{remark}

Finally, we introduce the following subregions of $\MM$
\bsub\lab{eq:defofsubregionsofMM}
 \begin{align}
 \MM(\tt_1,\tt_2):={}&\MM\cap\{\tt_1\leq \tt\leq \tt_2\}, \quad \forall\tau_1<\tau_2,\\
 \MM_{red}:={}&\MM\cap\{r\leq r_+(1+\dred)\},\\
 \Mtrap:={}&\MM_{r_+(1+2\dbl), 10m},\\ 
 \Mntrap:={}&\MM\setminus\Mtrap.
 \end{align}
 \esub

\begin{remark}
The conditions \eqref{eq:rangetauandronSigmastar} and \eqref{eq:dominantconditionforrstarcomparedtotaustartonS*} on $\Si_*$ are consistent with (3.2.7) and (3.4.5) in \cite{KS:Kerr}.
\end{remark}

%%%%%%%%%%%%%%%%%%%%%%%%%%%%%%%%%%%%%%%%%%%%%%%%%%%%%%%%%%%%%%%%%%%%%%

\subsubsection{Main properties of $(\MM, \g)$}
\lab{sec:assumptionsforsec:energyMorawetzesitmatesforTeukoslkyonMM:upto15derivatives}

%%%%%%%%%%%%%%%%%%%%%%%%%%%%%%%%%%%%%%%%%%%%%%%%%%%%%%%%%%%%%%%%%%%%%%

We assume that  the following identities hold on $\MM$ 
\bea\lab{eq:usefulalgebraicidentitiesinvolvingscalarproductsReJkReJkpm:Kerrpert}
\bsplit
\Re(\Jk)\c\Re(\Jk) &=\frac{(\sin\th)^2}{|q|^2}, \qquad \Re(\Jk)\c\Re(\Jk_+)=-\frac{x^2_p}{|q|^2},\qquad \Re(\Jk)\c\Re(\Jk_-)=\frac{x^1_p}{|q|^2},\\
\dual(\Re(\Jk))\c\Re(\Jk_+) &=\frac{\cos\th x^1_p}{|q|^2},\qquad \dual(\Re(\Jk))\c\Re(\Jk_-)=\frac{\cos\th x^2_p}{|q|^2}.
\end{split}
\eea

%%%%%%%%%%%%%%%%%%%%%%%%%%%%%%%%%%%%%%%%%%%%%%%%%%%%%%%%%%%%%%%%%%%%%%%%%%%%%
 
 \paragraph{\textit{Definition and estimates for $(\Ga_g, \Ga_b)$}.}
 
%%%%%%%%%%%%%%%%%%%%%%%%%%%%%%%%%%%%%%%%%%%%%%%%%%%%%%%%%%%%%%%%%%%%%%%%%%%%%

Using the linearized quantities introduced in Definition \ref{def:renormalizationofallnonsmallquantitiesinPGstructurebyKerrvalue}, we introduce the following definition for $(\Ga_g, \Ga_b)$ valid throughout Section \ref{sec:energyMorawetzesitmatesforTeukoslkyonMM:upto15derivatives}.   
\begin{definition}
\lab{definition.Ga_gGa_b:forMaSz26}
The set of all linearized quantities is of the form $\Ga_g\cup \Ga_b$ with  $\Ga_g,  \Ga_b$
 defined as follows.
 \begin{enumerate}
\item 
 The set $\Ga_g$ is given by $\Ga_g=\Ga_{g,1}\cup \Ga_{g, 2}\cup\Ga_{g,3}$   with
 \bea
 \bsplit
 \Ga_{g,1} &= \Big\{\Xi, \quad \omc, \quad\trXc,\quad  \Xh,\quad \Zc,\quad \Hbc, \quad \trXbc , \quad r\Pc, \quad  rB, \quad  rA\Big\},\\
 \Ga_{g,2} &= \Big\{\widecheck{e_4(r)}, \,\,\,\, r^{-1}\nab(r), \,\,\,\, r^{-1}\widecheck{e_4(\tau)}, \,\,\,\, r^{-1}\widecheck{\DD(\tau)}, \,\,\,\, r^{-1}\widecheck{e_3(\tau)},  \,\,\,\, e_4(\cos\th), \,\,\,\, \widecheck{e_4(x^1_p)}, \,\,\,\, \widecheck{e_4(x^2_p)}\Big\},\\
  \Ga_{g,3} &= \Big\{r\widecheck{\nab_4\Jk}, \quad r\widecheck{\nab_4\Jk_\pm}\Big\}.
 \end{split}
 \eea
 
 \item The set $\Ga_b$ is given by $\Ga_b=\Ga_{b,1}\cup \Ga_{b, 2}\cup \Ga_{b,3}$   with
 \bea
 \bsplit
 \Ga_{b,1}&= \Big\{\Hc, \quad \Xbh, \quad \omb, \quad \Xib,\quad  r\Bb, \quad \Ab\Big\},\\
  \Ga_{b, 2}&= \Big\{r^{-1}\widecheck{e_3(r)}, \quad  \widecheck{\DD(\cos\th)}, \quad e_3(\cos\th),  \quad  \widecheck{\DD(x^1_p)}, \quad {\widecheck{e_3(x^1_p)}}, \quad  \widecheck{\DD(x^2_p)}, \quad {\widecheck{e_3(x^2_p)}}\Big\}, \\
   \Ga_{b,3}&=\bigg\{ r\,\widecheck{\ov{\DD}\c\Jk}, \quad r\,\DD\hot\Jk, \quad r\,\widecheck{\nab_3\Jk}, \quad r\,\widecheck{\ov{\DD}\c\Jk_\pm}, \quad r\,\DD\hot\Jk_\pm, \quad r\,\widecheck{\nab_3\Jk_\pm}\bigg\}. 
   \end{split}
 \eea
\end{enumerate}
\end{definition}

\begin{remark}\lab{rmk:infactGagGabinthssectioniswidetildeGagwidetildeGab}
Note that Definition \ref{definition.Ga_gGa_b:forMaSz26} for $(\Ga_g, \Ga_b)$ differs slightly from Definition \ref{definition.Ga_gGa_b} and coincides in fact with Definition 4.24 in \cite{Sze}.
\end{remark}

We assume that  $(\Ga_g, \Ga_b)$ introduced in Definition \ref{definition.Ga_gGa_b} satisfy the following estimates 
\bea\lab{eq:decaypropertiesofGabGag}
\bsplit
&|\dk^{\leq {15}}\Ga_g|\les \min\left\{\frac{\ep}{r^2\tau^{\frac{1+\dec}{2}}}, \, \frac{\ep}{r\tau^{1+\frac{3\dec}{4}}}\right\}, \qquad |\dk^{\leq {15}}\Ga_b|\les \frac{\ep}{r\tau^{1+\frac{3\dec}{4}}},\\ 
&|\dk^{\leq {15}}\xi| \les \frac{\ep}{r^3},\qquad |\dk^{\leq {15}}\widecheck{e_4(\tau)}|\les \frac{\ep}{r^2},
\end{split}
\eea
where we recall that the weighted derivatives $\dk$ are defined by
\bea\lab{eq:defweightedderivative}
\dk:=\{\nab_3,\, r\nab_4,\, r\nabla\}.
\eea
For convenience, we also define the unweighted derivatives $\pr$ as follows
\bea\lab{eq:defunweightedderivative}
\pr:=\{\nab_3,\, \nab_4,\, \nabla\}.
\eea

%%%%%%%%%%%%%%%%%%%%%%%%%%%%%%%%%%%%%%%%%%%%%%%%%%%%%%%%%%%%%%%%%%%%%%%%%%%%

\paragraph{\textit{Assumptions on the inverse metric perturbation}.}
 
%%%%%%%%%%%%%%%%%%%%%%%%%%%%%%%%%%%%%%%%%%%%%%%%%%%%%%%%%%%%%%%%%%%%%%%%%%%%

We define, in the normalized coordinates $(\tt,r,x^1_0, x^2_0)$ and {$(\tt,r,x^1_p, x^2_p)$}, the inverse metric difference 
\bea
\gcheck^{\a\b}:=\g^{\a\b}-\gam^{\a\b}.
\eea
Then, with $(\Ga_g, \Ga_b)$ verifying \eqref{eq:decaypropertiesofGabGag}, we assume that $\gcheck^{\a\b}$ satisfies the following\footnote{These estimates hold in fact on each coordinate patch, i.e., in the coordinates $(\tt,r,x^1_0, x^2_0)$ for $\th\in[\frac{\pi}{4}, \frac{3\pi}{4}]$, and in the coordinates {$(\tt,r,x^1_p, x^2_p)$} for $\th\in [0,\pi]\setminus(\frac{\pi}{3}, \frac{2\pi}{3})$.}:
\bea\lab{eq:controloflinearizedinversemetriccoefficients}
\widecheck{\g}^{rr}=r\Ga_b, \quad \widecheck{\g}^{r\tau}=r\Ga_g,\quad \widecheck{\g}^{\tau\tau}=r\Ga_g, \quad \widecheck{\g}^{ra}=\Ga_b, \quad \widecheck{\g}^{\tau a}=\Ga_g,\quad \widecheck{\g}^{ab}=r^{-1}\Ga_g,
\eea 
and
\bea\lab{eq:controloflinearizedinversemetriccoefficients:inversegtautau}
|\dk^{\leq 15}\widecheck{\g}^{\tau\tau}|\les \frac{\ep}{r^2},
\eea 
with $(\Ga_b, \Ga_g)$ verifying \eqref{eq:decaypropertiesofGabGag}.

The following immediate non-sharp consequence of Lemma \ref{lem:specificchoice:normalizedcoord}, \eqref{eq:controloflinearizedinversemetriccoefficients}, \eqref{eq:controloflinearizedinversemetriccoefficients:inversegtautau} and  \eqref{eq:decaypropertiesofGabGag} will be useful
\bea\lab{eq:consequenceasymptoticKerrandassumptionsinverselinearizedmetric}
\begin{split}
\g^{rr}&=O(1), \qquad \g^{r\tau}=O(1), \qquad \g^{ra}=O(r^{-1}),\\
\g^{\tau\tau}&=O(m^2r^{-2}), \qquad \g^{\tau a}=O(mr^{-2}), \qquad \g^{ab}=O(r^{-2}).
\end{split}
\eea

%%%%%%%%%%%%%%%%%%%%%%%%%%%%%%%%%%%%%%%%%%%%%%%%%%%%%%%%%%%%%%%%%%%%%%%%%%%%%%%%%

\paragraph{\textit{Further properties of the metric coefficients}.}

%%%%%%%%%%%%%%%%%%%%%%%%%%%%%%%%%%%%%%%%%%%%%%%%%%%%%%%%%%%%%%%%%%%%%%%%%%%%%%%%%

We now provide further properties of the metric coefficients. All the statements and estimates are from Section 2.4 in \cite{MaSz24}. We start with the control of the perturbed metric coefficients. 
\begin{lemma}\lab{lemma:controlofmetriccoefficients:bis}
Assume that $\widecheck{\g}^{\a\b}$ verifies \eqref{eq:controloflinearizedinversemetriccoefficients} \eqref{eq:controloflinearizedinversemetriccoefficients:inversegtautau}. Then, $\widecheck{\g}_{\a\b}:=\g_{\a\b}-(\gam)_{\a\b}$ verifies  
\bea\lab{eq:controloflinearizedmetriccoefficients}
\begin{split}
\widecheck{\g}_{rr} &=r\Ga_g, \qquad\quad \widecheck{\g}_{r\tau}=r\Ga_g, \qquad\quad \widecheck{\g}_{\tau\tau}=r\Ga_b, \\
\widecheck{\g}_{\tau a} &=r^2\Ga_b, \qquad\, \widecheck{\g}_{ra}=r^2\Ga_g, \qquad\,\,\,\, \widecheck{\g}_{ab}=r^3\Ga_g,
\end{split}
\eea
and
\bea\lab{eq:controloflinearizedmetriccoefficients:grr}
|\dk^{\leq 15}\widecheck{\g}_{rr}|\les \frac{\ep}{r^2}.
\eea 
Also, we have
\bea\lab{eq:controloflinearizedmetriccoefficients:det}
\widecheck{\det(\g)}=\det(\g_{a,m})r^2\Ga_g, \qquad \widecheck{\det(\g)}:=\det(\g)-\det(\g_{a,m}).
\eea
\end{lemma}

The following immediate non-sharp consequence of Lemma \ref{lem:specificchoice:normalizedcoord}, \eqref{eq:controloflinearizedmetriccoefficients}, \eqref{eq:controloflinearizedmetriccoefficients:grr} and  \eqref{eq:decaypropertiesofGabGag} will be useful
\bea\lab{eq:consequenceasymptoticKerrandassumptionsinverselinearizedmetric:bis}
\bsplit
\g_{rr}&=O(m^2r^{-2}), \qquad \g_{r \tt}=O(1),\qquad \g_{ra}=O(m),\\
\g_{\tt \tt}&= O(1),\qquad\qquad\,\, \g_{\tau a}=O(r),\qquad \, \g_{ab}=O(r^2).
 \end{split}
\eea

\begin{lemma}\lab{lemma:computationofthederiveativeofsrqtg}
Let the 1-form $N_{det}$ be defined by  
\beaa
(N_{det})_\mu:=\frac{1}{\sqrt{|\g|}}\pr_\mu\sqrt{|\g|} - \frac{1}{\sqrt{|\g_{a,m}|}}\pr_\mu\sqrt{|\g_{a,m}|}. 
\eeaa
Then, we have 
\beaa
(N_{det})_r=\dk^{\leq 1}\Ga_g, \qquad (N_{det})_\tau=r\dk^{\leq 1}\Ga_g, \qquad (N_{det})_{x^a}=r\dk^{\leq 1}\Ga_g,
\eeaa 
and 
\beaa
(N_{det})^r=r\dk^{\leq 1}\Ga_g, \qquad (N_{det})^\tau=\dk^{\leq 1}\Ga_g, \qquad (N_{det})^{x^a}=r^{-1}\dk^{\leq 1}\Ga_g.
\eeaa
\end{lemma}

We have the following corollary of Lemma \ref{lemma:computationofthederiveativeofsrqtg}. 
\begin{corollary}\lab{cor:controloflinearizeddivergencecoordvectorfields}
We have
\beaa
\widecheck{\textbf{\textrm{Div}}(\pr_r)}=\dk^{\leq 1}\Ga_g, \qquad \widecheck{\textbf{\textrm{Div}}(\pr_\tau)}=r\dk^{\leq 1}\Ga_g, \qquad \widecheck{\textbf{\textrm{Div}}(\pr_{x^a})}=r\dk^{\leq 1}\Ga_g, \quad a=1,2.
\eeaa
\end{corollary}

Next, we provide the control of deformation tensors involved in energy-Morawetz estimates{, where the deformation tensor of a vectorfield $X$ is given by
\bea
\label{def:deformationtensor:lastsect}
{}^{(X)}\pi_{\a\b} :=  \D_{\a}X_{\b} + \D_{\b}X_{\a}=\LL_X\g_{\a\b}.
\eea}

\begin{lemma}\lab{lemma:controlofdeformationtensorsforenergyMorawetz}
The deformation {tensors of $\pr_\tau$ and $\pr_{\tphi}$ satisfy}
\beaa
\begin{split}
{}^{(\pr_\tau)}\pi_{rr}{, {}^{(\pr_{\tphi})}\pi_{rr}} &=r\dk^{\leq 1}\Ga_g, \qquad {}^{(\pr_\tau)}\pi_{r\tau}{, {}^{(\pr_{\tphi})}\pi_{r\tau}}=r\dk^{\leq 1}\Ga_g, \qquad {}^{(\pr_\tau)}\pi_{\tau\tau}{, {}^{(\pr_{\tphi})}\pi_{\tau\tau}}=r\dk^{\leq 1}\Ga_b, \\
{}^{(\pr_\tau)}\pi_{\tau a}{, {}^{(\pr_{\tphi})}\pi_{\tau a}} &=r^2\dk^{\leq 1}\Ga_b, \qquad {}^{(\pr_\tau)}\pi_{ra}{, {}^{(\pr_{\tphi})}\pi_{ra}}=r^2\dk^{\leq 1}\Ga_g, \qquad {}^{(\pr_\tau)}\pi_{ab}{, {}^{(\pr_{\tphi})}\pi_{ab}}=r^3\dk^{\leq 1}\Ga_g,
\end{split}
\eeaa
\beaa
\begin{split}
{}^{(\pr_\tau)}\pi^{rr}{, {}^{(\pr_{\tphi})}\pi^{rr}} &=r\dk^{\leq 1}\Ga_b, \qquad {}^{(\pr_\tau)}\pi^{r\tau}{, {}^{(\pr_{\tphi})}\pi^{r\tau}}=r\dk^{\leq 1}\Ga_g, \qquad {}^{(\pr_\tau)}\pi^{\tau\tau}{, {}^{(\pr_{\tphi})}\pi^{\tau\tau}}=r\dk^{\leq 1}\Ga_g, \\
{}^{(\pr_\tau)}\pi^{\tau a}{, {}^{(\pr_{\tphi})}\pi^{\tau a}} &=\dk^{\leq 1}\Ga_g, \qquad {}^{(\pr_\tau)}\pi^{ra}{, {}^{(\pr_{\tphi})}\pi^{ra}}=\dk^{\leq 1}\Ga_b, \qquad {}^{(\pr_\tau)}\pi^{ab}{, {}^{(\pr_{\tphi})}\pi^{ab}}=r^{-1}\dk^{\leq 1}\Ga_g,
\end{split}
\eeaa
and
\beaa
|\dk^{\leq 14}{}^{(\pr_\tau)}\pi_{rr}|+|\dk^{\leq 14}{}^{(\pr_{\tphi})}\pi_{rr}|\les \frac{\ep}{r^2}, \quad |\dk^{\leq 14}{}^{(\pr_\tau)}\pi^{\tau\tau}|+|\dk^{\leq 14}{}^{(\pr_{\tphi})}\pi^{\tau\tau}|\les \frac{\ep}{r^2}.
\eeaa

Also, the perturbed deformation tensor of $\pr_r$ satisfies 
\beaa
\begin{split}
\widecheck{{}^{(\pr_r)}\pi}_{rr} &=\dk^{\leq 1}\Ga_g, \qquad \widecheck{{}^{(\pr_r)}\pi}_{r\tau}=\dk^{\leq 1}\Ga_g, \qquad \widecheck{{}^{(\pr_r)}\pi}_{\tau\tau}=\dk^{\leq 1}\Ga_b, \\
\widecheck{{}^{(\pr_r)}\pi}_{\tau a} &=r\dk^{\leq 1}\Ga_b, \qquad \widecheck{{}^{(\pr_r)}\pi}_{ra}=r\dk^{\leq 1}\Ga_g, \qquad \widecheck{{}^{(\pr_r)}\pi}_{ab}=r^2\dk^{\leq 1}\Ga_g,
\end{split}
\eeaa
and
\beaa
\begin{split}
\widecheck{{}^{(\pr_r)}\pi}^{rr} &=\dk^{\leq 1}\Ga_b, \qquad \widecheck{{}^{(\pr_r)}\pi}^{r\tau}=\dk^{\leq 1}\Ga_g, \qquad \widecheck{{}^{(\pr_r)}\pi}^{\tau\tau}=\dk^{\leq 1}\Ga_g, \\
\widecheck{{}^{(\pr_r)}\pi}^{\tau a} &=r^{-1}\dk^{\leq 1}\Ga_g, \qquad \widecheck{{}^{(\pr_r)}\pi}^{ra}=r^{-1}\dk^{\leq 1}\Ga_b, \qquad \widecheck{{}^{(\pr_r)}\pi}^{ab}=r^{-2}\dk^{\leq 1}\Ga_g.
\end{split}
\eeaa
\end{lemma}

%%%%%%%%%%%%%%%%%%%%%%%%%%%%%%%%%%%%%%%%%%%%%%%%%%%%%%%%%%%%%%%%%%%%%%%%%%%%%%%%%

\paragraph{\textit{Assumptions on the regular triplet $\Om_i$, $i=1,2,3$ in Section \ref{sec:energyMorawetzesitmatesforTeukoslkyonMM:upto15derivatives}}.}

%%%%%%%%%%%%%%%%%%%%%%%%%%%%%%%%%%%%%%%%%%%%%%%%%%%%%%%%%%%%%%%%%%%%%%%%%%%%%%%%%

Recall from Definition \ref{def:Mialphaj:Kerr} that we associate to the regular triplet $\Om_i$, $i=1,2,3$, in Section \ref{sec:assumptionsforsec:energyMorawetzesitmatesforTeukoslkyonMM:upto15derivatives} the following 1-forms on $\MM$ 
\beaa
M_{i\a}^j:=(\Ddot_\a\Om_i)\c\Om^j, \quad \forall \a,i,j.
\eeaa
Our assumptions on the regular triplet $\Om_i$, $i=1,2,3$, are the following 
\bsub\lab{eq:assumptionsonregulartripletinperturbationsofKerr}
\bea
&&\lab{eq:assumptionsonregulartripletinperturbationsofKerr:0}
\widecheck{M_{i4}^j}=\Ga_g, \qquad \widecheck{M_{i3}^j}=\Ga_b,  \qquad \widecheck{M_{ia}^j}=\Ga_b, \quad \forall\, i,j,a,\\
&&\lab{eq:assumptionforLiebprtphiOmiinKerrperturbation}
\Lieb_{\pr_{\tphi}}\Om_i +\in_{ij3}\Om^j=r\Ga_b, \quad\textrm{for}\quad i=1,2,3.
\eea
\esub 

In view of Lemma \ref{lemma:computationoftheMialphajinKerr} for $(M_{i\a}^j)_K$ and the assumption \eqref{eq:assumptionsonregulartripletinperturbationsofKerr:0}
 for $\widecheck{M_{i\a}^j}$, we have the following estimates for $M_{i\a}^j$.
\begin{lemma}
\lab{lem:estimatesforMialphaj:Kerrpert}
Under the assumption \eqref{eq:assumptionsonregulartripletinperturbationsofKerr:0}  for $\widecheck{M_{i\a}^j}$, we have
 \begin{equation}
 \lab{estimates:Mialphaj:Kerrperturbations}
 M_{i4}^j={O(r^{-2})}, \quad M_{ia}^j = {O(r^{-1})}, \quad  M_{i3}^j=O(r^{-2})+\Ga_b,\quad M_{i\a}^j(\pr_{\tau})^{\a}=O(r^{-3}) + \Ga_b, 
 \end{equation}
 and 
 \bea
M_{i\a}^j( \pr_{r})^{\a}=\Ga_b, \qquad \g^{r\a}M_{i\a}^{j}=\Ga_b, \quad r\in [r_+(1+2\dbl), 12m].
 \eea
\end{lemma}

%%%%%%%%%%%%%%%%%%%%%%%%%%%%%%%%%%%%%%%%%%%%

\subsubsection{Teukolsky wave/transport systems in perturbations of Kerr}
\lab{sec:TeukolskyWavesysteminperturbationsofKerr}

%%%%%%%%%%%%%%%%%%%%%%%%%%%%%%%%%%%%%%%%%%%%

In this section, we first recall the form of the tensorial Teukolsky wave/transport systems in perturbations of Kerr derived in Section \ref{sec:precisederivationTeukolskywave-transportsystem:Kerrpert}. We then provide the corresponding scalarized form using the regular triplet introduced in Section \ref{sec:assumptionsforsec:energyMorawetzesitmatesforTeukoslkyonMM:upto15derivatives}.

%%%%%%%%%%%%%%%%%%%%%%%%%%%%%%%%%%%%%%%%%%%%%%

\paragraph{\textit{Tensorial Teukolsky wave/transport systems in perturbations of Kerr}.}

%%%%%%%%%%%%%%%%%%%%%%%%%%%%%%%%%%%%%%%%%%%%%%

We consider $\pmb\phi_s^{(p)}\in\sk_2(\mathbb{C})$, $s=\pm 2$, $p=0,1,2$, given by \eqref{eq:definitionofthephiplus2phierarchy:perturbationofKerr} \eqref{eq:definitionofthephiminus2phierarchy:perturbationofKerr}. 
Then, the tensorial Teukolsky wave equations in perturbations of Kerr are given by 
\bsub
\lab{eq:TensorialTeuSysandlinearterms:rescaleRHScontaine2:general:Kerrperturbation}
\bea
\lab{eq:TensorialTeuSys:rescaleRHScontaine2:general:Kerrperturbation}
\bigg(\squared_2 -\frac{4ia\cos\th}{|q|^2}\nab_{\pr_{\tt}}- \frac{4-2\de_{p0}}{\qs}\bigg){\phis{p}} = \L_{s}^{(p)}[\pmb\phi_{s}]+\N_{W,s}^{(p)}, \quad s=\pm2, \quad p=0,1,2,
\eea
where nonlinear terms $\N_{W,s}^{(p)}$ are given by \eqref{eq:schematicformofNpWsplus2} \eqref{eq:schematicformofNpWsminus2}, and where the linear coupling terms $\L_{s}^{(p)}[{\pmb\phi_s}]$ have the following schematic forms
\bea
\lab{eq:tensor:Lsn:onlye_2present:general:Kerrperturbation:bis}
\bsplit
{\L_{s}^{(0)}[\pmb\phi_{s}]}={}& (2sr^{-3} +O(mr^{-4}))\phis{1}+ O(mr^{-3}) \nab_{\Xcal_s}^{\leq 1}\phis{0},\\
{\L_{s}^{(1)}[\pmb\phi_{s}]}={}& (sr^{-3} +O(mr^{-4}))\phis{2}+ O(mr^{-3}) \nab_{\Xcal_s}^{\leq 1}  \phis{1}+O(mr^{-2})\nab_{\pr_{\tphi}+a\pr_{\tt}}^{\leq 1}\phis{0},\\
{\L_{s}^{(2)}[\pmb\phi_{s}]}={}&O(mr^{-3})\phis{2}+O(mr^{-2})\nab_{\pr_{\tphi}+a\pr_{\tt}}^{\leq 1}\phis{1}+O(m^2 r^{-2})\phis{0},
\end{split}
\eea
\esub
with $\Xcal_s$, $s=\pm 2$, being the regular vectorfields defined by
\bea\lab{eq:formofregularhorizontalvectorfieldmathcalXs:Kerrperturbation}
\mathcal{X}_{s} &:=& s\big(\pr_{\tphi}+a(\sin\th)^2\pr_\tau\big) - \frac{2a\cos\th}{r}\sin\th\pr_\th, \quad s=\pm 2, 
\eea
with all the coefficients in  \eqref{eq:tensor:Lsn:onlye_2present:general:Kerrperturbation:bis} being independent of coordinates $\tau$ and 
$\tphi$, and with the coefficients in front of the terms $\phis{2}$ and $\nab_{\pr_{\tphi}+a\pr_{\tt}}\phis{1}$ on the RHS of equation of ${\L_{s}^{(2)}[\pmb\phi_{s}]}$ in \eqref{eq:tensor:Lsn:onlye_2present:general:Kerrperturbation:bis} being real functions. 

Moreover, the tensorial Teukolsky transport equations in perturbations of Kerr are given by
\bsub\lab{def:TensorialTeuScalars:wavesystem:Kerrperturbation}
 \bea
\lab{def:TensorialTeuScalars:wavesystem:Kerrperturbation:+2}
\nab_3 \left(\frac{r\bar{q}}{q}\left(\frac{r^2}{|q|^2}\right)^{p-2}\pmb\phi_{+2}^{(p)}\right)=\frac{\bar{q}}{rq}\left(\frac{r^2}{|q|^2}\right)^{p-1}\pmb\phi_{+2}^{(p+1)}+\N_{T,+2}^{(p)}, \quad p=0,1,
\eea
and
\bea
\lab{def:TensorialTeuScalars:wavesystem:Kerrperturbation:-2}
\nab_4\left(\frac{rq}{\bar{q}}\left(\frac{{r^2}}{|q|^2}\right)^{p-2}\pmb\phi_{-2}^{(p)}\right)=\frac{q}{r\bar{q}}\left(\frac{r^2}{|q|^2}\right)^{p-1}\frac{\De}{\qs}\pmb\phi_{-2}^{(p+1)}+\N_{T,-2}^{(p)}, \,\,\,\, p=0,1,
\eea
\esub
where the nonlinear terms $\N_{T,s}^{(p)}$ are given by \eqref{eq:transportequationsins=+2caseforp=0and1:RHSNT+2p=0and1} \eqref{eq:transportequationsins=-2caseforp=0and1:RHSNT-2p=0and1}. The equations \eqref{eq:TensorialTeuSysandlinearterms:rescaleRHScontaine2:general:Kerrperturbation}  and \eqref{def:TensorialTeuScalars:wavesystem:Kerrperturbation} then correspond to the tensorial Teukolsky wave/transport systems in perturbations of Kerr.

In view of \eqref{eq:definitionofthephiminus2phierarchy:perturbationofKerr}, $\pmb\phi_{-2}^{(0)}$ degenerates at $r=r_+$ and does thus not allow to recover estimates for $\Ab$ near $r=r_+$. To remedy this problem, we will rely on the following lemma.
\begin{lemma}\lab{lemma:pmbphiminus2p=0isdegenerateatr=rplus}
The horizontal tensor $\Ab\in\sk_2(\mathbb{C})$ satisfies in the redshift region $r\leq r_+(1+2\dred)$:
\bea\lab{eq:waveequationpmbphip=0sminus2nodeginredshiftregion}
\squared_2\nab_4^p\Ab = (2-p)\pr_r\left(\frac{\De}{|q|^2}\right)\nab_3\nab_4^p\Ab
+ O(1)\big(\nab_{4}\nab_4^{\leq p}\Ab,\nab\nab_4^{\leq p}\Ab, \nab_4^{\leq p}\Ab\big)+\N_{\nab_4^p\Ab}, \quad p=0,1,2,
\eea
where we have in $r\leq r_+(1+2\dred)$
\bea
\N_{\Ab}=\dk^{\leq 1}(\Ga_g\c\Ga_b), \qquad \N_{\nab_4^p\Ab}=\nab_4^p\N_{\Ab}+\dk^{\leq p+1}(\Ga_g\c\Ab), \quad p=1,2.
\eea
\end{lemma}

\begin{proof}
See Lemma 5.3 in \cite{Sze}.
\end{proof}

%%%%%%%%%%%%%%%%%%%%%%%%%%%%%%%%%%%%%%%%%%%%%%

\paragraph{\textit{Scalarized Teukolsky wave/transport systems in perturbations of Kerr}.}

%%%%%%%%%%%%%%%%%%%%%%%%%%%%%%%%%%%%%%%%%%%%%%

Applying Lemma \ref{lemma:formoffirstordertermsinscalarazationtensorialwaveeq} to the tensorial Teukolsky wave/transport systems \eqref{eq:TensorialTeuSysandlinearterms:rescaleRHScontaine2:general:Kerrperturbation}  \eqref{def:TensorialTeuScalars:wavesystem:Kerrperturbation} by using the regular triplet $\Om_i$, $i=1,2,3$ introduced in Section \ref{sec:assumptionsforsec:energyMorawetzesitmatesforTeukoslkyonMM:upto15derivatives}, we obtain  the following scalarized Teukolsky wave/transport systems in perturbations of Kerr.  

\begin{lemma}[Scalarized Teukolsky wave/transport systems in perturbations of Kerr]
\lab{lem:scalarizedTeukolskywavetransportsysteminKerrperturbation:Omi}
Let $\Om_i$, $i=1,2,3$ be the regular triplet introduced in Section \ref{sec:assumptionsforsec:energyMorawetzesitmatesforTeukoslkyonMM:upto15derivatives}, and define the complex-valued scalars $\phiss{ij}{p}$ by 
\bea
\phiss{ij}{p}:=\pmb\phi_s^{(p)}(\Om_i,\Om_j), \quad s=\pm 2,\quad i,j=1,2,3, \quad p=0,1,2, 
\eea
where $\pmb\phi_s^{(p)}\in\sk_2(\mathbb{C})$, $s=\pm 2$, $p=0,1,2$, satisfy the tensorial Teukolsky wave/transport systems \eqref{eq:TensorialTeuSysandlinearterms:rescaleRHScontaine2:general:Kerrperturbation}  \eqref{def:TensorialTeuScalars:wavesystem:Kerrperturbation}. Then, the tensorial Teukolsky wave systems  \eqref{eq:TensorialTeuSysandlinearterms:rescaleRHScontaine2:general:Kerrperturbation} scalarized using $(\Om_i)_{i=1,2,3}$ take the following form 
\bea
\lab{eq:ScalarizedTeuSys:general:Kerrperturbation} 
\widehat\square_\g(\phi_s^{(p)})_{ij} -\frac{4-2\de_{p0}}{\qs} \phiss{ij}{p} = {L_{s,ij}^{(p)}}+N_{W,s,ij}^{(p)}, \quad s=\pm2, \quad i,j=1,2,3, \quad p=0,1,2,
\eea
where we have defined 
\bea
\lab{eq:defofwidehatsquaregoperator}
\widehat\square_\g(\phi_s^{(p)})_{ij}&:=&\square_\g\phiss{ij}{p}- \widehat{S}(\phi_s^{(p)})_{ij} -(\widehat{Q}\phi_s^{(p)})_{ij}
\eea
with
\bsub\lab{eq:definitionwidehatSandwidehatQperturbationsofKerr}
\begin{align}
\widehat{S}(\phi_s^{(p)})_{ij} ={}&S(\phi_s^{(p)})_{ij} +\frac{4ia\cos\th}{|q|^2} \pr_\tau\phiss{ij}{p},\\
(\widehat{Q}\phi_s^{(p)})_{ij} ={}&(Q\phi_s^{(p)})_{ij} 
-\frac{4ia\cos\th}{|q|^2}\big(M_{i\tau}^l \phiss{lj}{p}+M_{j\tau}^l \phiss{il}{p}\big),
\end{align}
\esub
where the linear coupling terms ${L_{s,ij}^{(p)}}:=(\L_{s}^{(p)}[\pmb\phi_s])_{ij}$ are given by
\bsub
\lab{eq:linearterms:ScalarizedTeuSys:general:Kerrperturbation}
\bea
{L_{s,ij}^{(0)}}&=& (2sr^{-3} +O(mr^{-4}))\phiss{ij}{1} + O(mr^{-3}){\Xcal_s}\phiss{ij}{0}+\sum_{k,l=1,2,3}O(mr^{-3})\phiss{kl}{0},\\
{L_{s,ij}^{(1)}}&=&(sr^{-3} +O(mr^{-4}))\phiss{ij}{2} + O(mr^{-3}){\Xcal_s}\phiss{ij}{1}+O(mr^{-2})(\pr_{\tphi} +a\pr_\tau)\phiss{ij}{0}\nn\\
  &&+\sum_{k,l=1,2,3}\Big(O(mr^{-3})\phiss{kl}{1}+O(mr^{-2})\phiss{kl}{0}\Big),\\
{L_{s,ij}^{(2)}}&=&O(mr^{-3}) \phiss{ij}{2}+O(mr^{-2})(\pr_{\tphi} +a\pr_\tau)\phiss{ij}{1} +O(m^2r^{-2})\phiss{ij}{0}\nn\\
&&+\sum_{k,l=1,2,3}O(mr^{-2})\phiss{kl}{1}
\eea
\esub
with $\Xcal_s$, $s=\pm 2$, being the regular vectorfields introduced in \eqref{eq:formofregularhorizontalvectorfieldmathcalXs:Kerrperturbation}, with the coefficients in front of the terms $\phiss{ij}{2}$ and $(\pr_{\tphi}+a\pr_{\tau})\phiss{ij}{1}$ on the RHS of equation of $L_{s,ij}^{(2)}$ in \eqref{eq:linearterms:ScalarizedTeuSys:general:Kerrperturbation} being real functions, and with all the coefficients in the first line of the three equations in  \eqref{eq:linearterms:ScalarizedTeuSys:general:Kerrperturbation} being independent of coordinates $\tau$ and 
$\tphi$, where $S(\phi_s^{(p)})_{ij}$ and $(Q\phi_s^{(p)})_{ij}$ are given as in \eqref{SandV}, and where the complex-valued scalars $N_{W,s,ij}^{(p)}$ are given by 
\bea
N_{W,s,ij}^{(p)}:=\N_{W,s}^{(p)}(\Om_i, \Om_j),\quad s=\pm 2, \quad i,j=1,2,3, \quad p=0,1,2.
\eea

Moreover, the tensorial Teukolsky transport equations \eqref{def:TensorialTeuScalars:wavesystem:Kerrperturbation} scalarized using $(\Om_i)_{i=1,2,3}$ take the following form 
\bsub
 \lab{eq:ScalarizedQuantitiesinTeuSystem:Kerrperturbation}
\bea
\bsplit
& e_3\bigg(\frac{r\bar{q}}{q}\bigg({\frac{r^2}{|q|^2}}\bigg)^{p-2}\phipluss{ij}{p} \bigg) - \frac{r\bar{q}}{q}\bigg({\frac{r^2}{|q|^2}}\bigg)^{p-2}\Big(M_{i3}^k \phipluss{kj}{p} + M_{j3}^k \phipluss{ik}{p}\Big)\\
=&\frac{\bar{q}}{rq}\bigg({\frac{r^2}{|q|^2}}\bigg)^{p-1}\phipluss{ij}{p+1}+N_{T,+2,ij}^{(p)}
\end{split}
\eea
and 
\bea
\bsplit
& e_4\bigg(\frac{rq}{\bar{q}}\bigg(\frac{r^2}{|q|^2}\bigg)^{p-2}\phiminuss{ij}{p}\bigg) - \frac{rq}{\bar{q}}\bigg(\frac{r^2}{|q|^2}\bigg)^{p-2}\Big(M_{i4}^k \phiminuss{kj}{p} + M_{j4}^k \phiminuss{ik}{p}\Big)\\
=& \frac{q}{r\bar{q}}\left(\frac{r^2}{|q|^2}\right)^{p-1}\frac{\De}{\qs}\phiminuss{ij}{p+1} +N_{T,-2,ij}^{(p)},
\end{split}
\eea
\esub
where the complex-valued scalars $N_{T,s,ij}^{(p)}$ are given by 
\bea
N_{T,s,ij}^{(p)}:=\N_{T,s}^{(p)}(\Om_i, \Om_j),\quad s=\pm 2, \quad i,j=1,2,3, \quad p=0,1.
\eea
The equations \eqref{eq:ScalarizedTeuSys:general:Kerrperturbation} \eqref{eq:ScalarizedQuantitiesinTeuSystem:Kerrperturbation} then correspond to the scalarized Teukolsky wave/transport systems in perturbations of Kerr.
\end{lemma}

\begin{proof}
See Lemma 5.26 in \cite{MaSz26}.
\end{proof}

%%%%%%%%%%%%%%%%%%%%%%%%%%%%%%%%%%%%%%%%%%%%%%%%%%%%%%%%%

\subsubsection{The spacetime region $\RR_*$ and the expanded spacetime $\MMh$}
\lab{sec:derfintionofRRstarandMMh}

%%%%%%%%%%%%%%%%%%%%%%%%%%%%%%%%%%%%%%%%%%%%%%%%%%%%%%%%%

In \cite{MaSz26}, we derive energy-Morawetz estimates on $\MM(\tau_1, \tau_2)$ where $\MM$ extends all the way to $\II_+$. In this section, we adapt the results in that paper to the case where $\MM$ only extends to $\Si_*$. We start by introducing the spacetime region where the modifications will be needed.

%%%%%%%%%%%%%%%%%%%%%%%%%%%%%%%%%%%%%%%%%%%%%%%%%%%%%%%%%

\paragraph{\textit{The spacetime region $\RR_*$}.}

%%%%%%%%%%%%%%%%%%%%%%%%%%%%%%%%%%%%%%%%%%%%%%%%%%%%%%%%%

Recall from \eqref{eq:rangetauandronSigmastar} that 
\beaa
\Si_*=\{\tau+r=c_*+h_*(r)\}\cap\{1\leq\tau\leq\tau_*\}, \qquad h_*(r)=O(m^2r^{-1}),
\eeaa
with $r_*$ satisfying the dominance condition \eqref{eq:dominantconditionforrstarcomparedtotaustartonS*}. Then, we consider the following causal spacetime regions of $\MM$
\bea\lab{eq:definitionofthecausalspacetimeregionsRR*andRR*tau1tau2}
\bsplit
\RR_*:=&\MM\cap\left\{\frac{c_*}{2}\leq \tau+r-h_*(r)\leq c_*\right\}, \\ 
\RR_*(\tau_1, \tau_2):=&\RR_*\cap\{\tau_1\leq\tau\leq\tau_2\}, \quad 1\leq \tau_1<\tau_2\leq \tau_*,
\end{split}
\eea
so that $\RR_*=\RR_*(1,\tau_*)$ and $\RR_*(\tau_1, \tau_2)$ has the following causal boundaries 
\beaa
\bsplit
\pr^-\RR_*(\tau_1, \tau_2) &= \left(\Si(\tau_1)\cap\left\{\frac{c_*}{2}\leq \tau+r-h_*(r)\leq c_*\right\}\right)\cup\left(\left\{\tau+r=\frac{c_*}{2}+h_*(r)\right\}\cap\MM(\tau_1, \tau_2)\right),\\
\pr^+\RR_*(\tau_1, \tau_2) &= \left(\Si(\tau_2)\cap\left\{\frac{c_*}{2}\leq \tau+r-h_*(r)\leq c_*\right\}\right)\cup\Si_*(\tau_1, \tau_2),
\end{split}
\eeaa
where, arguing as for \eqref{eq:dominantconditionforrstarcomparedtotaustartonS*:holdsforallronSi*}, using also \eqref{eq:controlofsizeofc*assimeqr*:onSi*}, we have
\bea\lab{eq:dominantconditionforrstarcomparedtotaustartonS*:holdsforallronRR*}
r\simeq r_*, \qquad r\simeq \ep_0^{-1}\tau_*^{1+\dec},\quad\textrm{on}\quad\RR_*.
\eea

Next, we introduce the following 1-parameter family of hypersurfaces of $\RR_*$
\bea\lab{eq:defofSigmastar:c}
\Si_{*,c}:=\{\tau+r=c+h_*(r)\}\cap\{1\leq\tau\leq\tau_*\}, \quad \frac{c_*}{2}\leq c\leq c_*, 
\eea
which are spacelike in view of \eqref{eq:dominantconditionforrstarcomparedtotaustartonS*:holdsforallronRR*}, and $\Si_{*,c}(\tau_1, \tau_2):=\Si_{*,c}\cap\{\tau_1\leq\tau\leq\tau_2\}$, so that 
\beaa
\Si_{*,c_*}=\Si_*, \qquad \RR_*=\bigcup_{\frac{c_*}{2}\leq c\leq c_*}\Si_{*,c}, \qquad \RR_*(\tau_1, \tau_2)=\bigcup_{\frac{c_*}{2}\leq c\leq c_*}\Si_{*,c}(\tau_1, \tau_2).
\eeaa
Note also from \eqref{eq:defofSigmastar:c}, \eqref{eq:dominantconditionforrstarcomparedtotaustartonS*} and \eqref{eq:controlofsizeofc*assimeqr*:onSi*} that 
\bea\lab{eq:rsimcuptoOep0inSi*cforcbetweenc*on2andc*}
|r-c|\les \ep_0r_*\les \ep_0 c_*\les \ep_0c\quad\Longrightarrow\quad r=c\big(1+O(\ep_0)\big)\quad\textrm{on}\quad\Si_{*,c}, \quad \frac{c_*}{2}\leq c\leq c_*.
\eea
Finally, we introduce for convenience the decreasing function $\varsigma_*$ defined by
\bea\lab{eq:defdecreasingfunctionvarsigma*definingSi*c}
r+\tau=c+h_*(r), \quad 1\leq\tau\leq\tau_*, \quad \frac{c_*}{2}\leq c\leq c_* \quad\Longleftrightarrow\quad r=\varsigma_*(\tau-c),
\eea
so that 
\beaa
\bsplit
\RR_*&=\{\varsigma_*(\tau-c_*/2)\leq r\leq\varsigma_*(\tau-c_*)\}\cap\{1\leq\tau\leq\tau_*\}, \qquad \MM\setminus\RR_*=\MM\cap\{r\leq \varsigma_*(\tau-c_*/2)\},\\
\Si_{*,c}&=\{r=\varsigma_*(\tau-c)\}\cap\{1\leq\tau\leq\tau_*\}, \quad \frac{c_*}{2}\leq c\leq c_*.
\end{split}
\eeaa

%%%%%%%%%%%%%%%%%%%%%%%%%%%%%%%%%%%%%%%%%%%%%%%%%%%%%%%%%

\paragraph{\textit{The expanded spacetime $(\MMh, \gh)$}.}

%%%%%%%%%%%%%%%%%%%%%%%%%%%%%%%%%%%%%%%%%%%%%%%%%%%%%%%%%

Next, we introduce the extended spacetime region $(\MMh, \gh)$ as follows: 
\begin{itemize}
\item The four dimensional manifold $\MMh$ is covered by coordinate systems $(\tt, r, x^1_0, x^2_0)$ and $(\tt, r, x^1_p, x^2_p)$, defined respectively on $\th\neq 0, \pi$ and $\th\neq\frac{\pi}{2}$, with 
\beaa
\tau\in\mathbb{R}, \quad r_+(1-\dhor)\leq r<+\infty, \quad x^1_0=\th, \quad x^2_0=\tphi, \quad {x^1_p=\sin\th\cos\tphi, \quad x^2_p=\sin\th\sin\tphi}.
\eeaa

\item Let $\chi_{\tau_*}=\chi_{\tau_*}(\tau)$ be a smooth cut-off function satisfying 
\bea\lab{eq:propertieschitau*forwidehatg}
\chi_{\tau_*}(\tau)=0\,\,\,\textrm{on}\,\,\, \mathbb{R}\setminus(1, \tau_*), \quad\chi_{\tau_*}(\tau)=1\,\,\,\textrm{on}\,\,\, [2, \tau_*-1], \quad\|\chi_{\tau_*}\|_{W^{{15},+\infty}(\mathbb{R})}\les 1,
\eea
and let $\chi_{r_*}=\chi_{r_*}(\tau, r)$ be a smooth cut-off function satisfying 
\bea\lab{eq:supportpropertiesofcutoffchistar}
\bsplit
\chi_{r_*}&=1\quad\textrm{for}\quad r\leq\varsigma_*\left(\tau-\frac{3c_*}{4}\right),\\ 
\chi_{r_*}&=0\quad\textrm{for}\quad r\geq\varsigma_*\left(\tau-\frac{7c_*}{8}\right), \qquad \|(r\pr_\tau, r\pr_r)^{\leq 15}\chi_{r_*}\|_{L^\infty}\lesssim 1.
\end{split}
\eea
The Lorentzian metric $\gh$ is given on $\MMh$ by 
\bea\lab{eq:definitionofLorentzianmetricwidehatgonwidehatMM}
\gh:=\chi_{r_*}\chi_{\tau_*}\g+(1-\chi_{r_*}\chi_{\tau_*})\gam, 
\eea
so that 
\bea\lab{eq:wherewidehatgisequaltogandisequaltogam}
\bsplit
\gh &=\g\quad\textrm{on}\quad\MM(2, \tau_*-1)\cap\left\{r\leq\varsigma_*\left(\tau-\frac{3c_*}{4}\right)\right\}, \\
\gh &=\gam\quad\textrm{on}\quad\MMh\setminus\left(\MM\cap\left\{r\geq\varsigma_*\left(\tau-\frac{7c_*}{8}\right)\right\}\right).
\end{split}
\eea
\end{itemize}

We introduce the following subregions of $\MMh$
\bsub\lab{eq:defofsubregionsofwidehatMM}
 \begin{align}
 \MMh(\tt_1,\tt_2):={}&\MMh\cap\{\tt_1\leq \tt\leq \tt_2\}, \quad \forall\tau_1<\tau_2,\\
 \MMh_{red}:={}&\MMh\cap\{r\leq r_+(1+\dred)\},\\
 \Mtraph:={}&\MMh_{r_+(1+2\dbl), 10m},\\ 
\Mntraph:={}&\MMh\setminus\Mtraph,
 \end{align}
 \esub
 and denote by $\Sih(\tau)$ the level sets of $\tau$ in $\MMh$ on which $r$ takes all values in $(r_+(1-\dhor), +\infty)$.
 
Also, in view of \eqref{eq:definitionofLorentzianmetricwidehatgonwidehatMM}, we have
\beaa
\widecheck{\gh}^{\a\b}&:=&\gh^{\a\b} -\gam\\
&=& \chi_{r_*}\chi_{\tau_*}\g^{\a\b}+(1-\chi_{r_*}\chi_{\tau_*})\gam^{\a\b} -\gam^{\a\b}=\chi_{r_*}\chi_{\tau_*}\left(\g^{\a\b}-\gam^{\a\b}\right)=\chi_{r_*}\chi_{\tau_*}\widecheck{\g}^{\a\b}
\eeaa
and hence, using also \eqref{eq:propertieschitau*forwidehatg} \eqref{eq:supportpropertiesofcutoffchistar}, we infer
\beaa
|\dk^{\leq {15}}\widecheck{\gh}^{\a\b}|\les |\dk^{\leq {15}}(\chi_{r_*}\chi_{\tau_*})||\dk^{\leq {15}}\widecheck{\g}^{\a\b}|\les |\pr_\tau^{\leq {15}}\chi_{\tau_*}||(\pr_\tau, r_*\pr_r)^{\leq 15}\chi_{r_*}||\dk^{\leq {15}}\widecheck{\g}^{\a\b}| \les |\dk^{\leq {15}}\widecheck{\g}^{\a\b}|,
\eeaa
where we used the fact that $|\dk^{\leq {15}}(\chi_{r_*}\chi_{\tau_*})|\les |(\pr_\tau, r\pr_r)^{\leq {15}}(\chi_{r_*}\chi_{\tau_*})|$. Since $\g$ satisfies \eqref{eq:controloflinearizedinversemetriccoefficients} \eqref{eq:controloflinearizedinversemetriccoefficients:inversegtautau}, and in view of the properties \eqref{eq:propertieschitau*forwidehatg} \eqref{eq:supportpropertiesofcutoffchistar} respectively of $\chi_{\tau_*}$ and $\chi_{r_*}$, we deduce that $\gh$ satisfies
\bea\lab{eq:controloflinearizedinversemetriccoefficients:widehat}
\widecheck{\gh}^{rr}=r\Ga_b, \quad \widecheck{\gh}^{r\tau}=r\Ga_g,\quad \widecheck{\gh}^{\tau\tau}=r\Ga_g, \quad \widecheck{\gh}^{ra}=\Ga_b, \quad \widecheck{\gh}^{\tau a}=\Ga_g,\quad \widecheck{\gh}^{ab}=r^{-1}\Ga_g,
\eea 
and
\bea\lab{eq:controloflinearizedinversemetriccoefficients:inversegtautau:widehat}
|\dk^{\leq 15}\widecheck{\gh}^{\tau\tau}|\les \frac{\ep}{r^2}.
\eea

%%%%%%%%%%%%%%%%%%%%%%%%%%%%%%%%%%%%%

\paragraph{\textit{Future null infinity of $\MMh$}.}

%%%%%%%%%%%%%%%%%%%%%%%%%%%%%%%%%%%%%

We start by constructing an auxiliary ingoing optical function $\tauu$ in a subregion of $(\MMh, \gh)$. 
\begin{lemma}\lab{lemma:constructionoftheingoingopticalfunctiontauu}
There exists an ingoing optical function $\tauu$ defined in $\MMh\cap\{r\geq c_*\}$ by 
\bea
\tauu:=\tauu_0 +\tauut, \qquad \tauu_0:=\tau+2r+4m\log\left(\frac{r}{2m}\right),
\eea
where $\tauut$ satisfies 
\bea
|\dk^{\leq 2}\tauut|\les r^{-1}\quad\textrm{in}\quad\MMh\cap\{r\geq c_*\}.
\eea
\end{lemma}

\begin{proof}
This follows from the fact that the existence of $\tauu$ is immediate in Kerr and the fact that $\gh=\gam$ for $r\geq c_*$ in view of \eqref{eq:wherewidehatgisequaltogandisequaltogam}.
\end{proof}

Making use of the ingoing optical function $\tauu$, we may now define $\II_+$. 
\begin{definition}[Definition of $\II_+$]
\lab{def:howmathcalIplusisdefinedinMM} 
Consider the coordinates $(\tauu, \tau, x^1, x^2)$ covering the spacetime region $\MMh\cap\{r\geq c_*\}$, where $\tauu$ is the ingoing optical function constructed in Lemma \ref{lemma:constructionoftheingoingopticalfunctiontauu}. Then, the future null infinity of $(\MMh, \gh)$ is defined as
\bea
\II_+:=\MMh\cap\{\tauu=+\infty\}.
\eea
\end{definition}

The following lemma provides the control of the induced geometry on $\II_+$ in the perturbed spacetime $(\MMh, \gh)$. 
\begin{lemma}\lab{lemma:controllinearizedmetric:inducedmetricII+}
Let $\II_+$ be given by Definition \ref{def:howmathcalIplusisdefinedinMM}. Consider the coordinates system $(\tau, x^1, x^2)$ covering $\II_+$, and denote by $(\pr_\tau^{\II_+}, \pr_{x^1}^{\II_+}, \pr_{x^2}^{\II_+})$ the corresponding coordinate vectorfields. Then, 
\begin{enumerate}
\item the coordinate vectorfields $\pr_{x^a}^{\II_+}$, $a=1,2,$ satisfy 
\bea\label{expression:prxaIIplus:nullinf}
\pr_{x^a}^{\II_+}=\pr_{x^a},\,\,\, a=1,2,
\eea

\item the spheres $S^{\II_+}(\tau_1):=\II_+\cap\{\tau=\tau_1\}$ foliating $\II_+$ are round,

\item $\pr_\tau^{\II_+}$ is ingoing null and 
\bea\label{expression:prtauIIplus:nullinf}
\pr_\tau^{\II_+}=\pr_\tau -\frac{1}{2}\pr_r, 
\eea

\item $\pr_r$ is an outgoing null vectorfield on $\II_+$ and satisfies 
\bea\label{expression:prrIIplus:nullinf}
\g(\pr_\tau^{\II_+}, \pr_r)=-1, \qquad \g(r^{-1}\pr_{x^a}^{\II_+}, \pr_r)=0.
\eea
\end{enumerate}
\end{lemma}

\begin{proof}
The proof is immediate since $\gh=\gam$ for $r\geq c_*$ in view of \eqref{eq:wherewidehatgisequaltogandisequaltogam}.
\end{proof}

%%%%%%%%%%%%%%%%%%%%%%%%%%%%%%%%%%%%%%%%%%%%%%%%%%%%%%%
   
\subsubsection{Main norms energy-Morawetz norms on $\MMh$ and $\MM$}
        \lab{subsection:basicnormsforpsi:MMh}
        
%%%%%%%%%%%%%%%%%%%%%%%%%%%%%%%%%%%%%%%%%%%%%%%%%%%%%%%

We introduce in this section the energy, Morawetz and flux norms needed to state our main result. First, given any $(\tau, r)$, and for any scalar function $F$ on the spheres $S(\tau, r)$ of constant $\tau$ and $r$, we introduce the following notation 
\beaa
\int_{\mathbb{S}^2}F(\tau, r, \om)d\mathring{\ga} := \int F(\tau, r, x^1, x^2)\sqrt{\det(\mathring{\ga})}dx^1dx^2,
\eeaa 
as well as the corresponding notation for the spheres $S^{\II_+}(\tau)$ of constant $\tau$ on $\II_+$. We start by introducing energy, Morawetz and flux norms for horizontal tensors in $\sk_k(\mathbb{C})$, $k=1,2$ on $\MMh$.

%%%%%%%%%%%%%%%%%%%%%%%%%%%%%%%%%%%%%%%%%%%%%%%%%%%%%%%%

\paragraph{\textit{Energy, Morawetz and flux norms for $\pmb\psi\in\sk_k(\mathbb{C})$, $k=1,2$ on $\MMh$}.}

%%%%%%%%%%%%%%%%%%%%%%%%%%%%%%%%%%%%%%%%%%%%%%%%%%%%%%%%

For $\pmb\psi\in\sk_k(\mathbb{C})$, $k=1,2$, and $\tau_1<\tau_2$, we define flux norms\footnote{For $\F_{\II_+}[\psi](\tau_1,\tau_2)$, recall that $\II_+=\MM\cap\{\tauu=+\infty\}$ where the ingoing optical function $\tauu$ has been constructed in Lemma \ref{lemma:constructionoftheingoingopticalfunctiontauu}, and recall that the notations $\pr_\tau^{\II_+}$ and $\pr_{x^a}^{\II_+}$ on $\II_+$ have been introduced in Lemma \ref{lemma:controllinearizedmetric:inducedmetricII+}.}
\bsub
\label{def:variousMorawetzIntegrals}
\bea
\bsplit
\Fh_{\AA}[\pmb\psi](\tau_1,\tau_2):=& \int_{\tau_1}^{\tau_2}\int_{\mathbb{S}^2}\big(|\mu| |\nab_3\pmb\psi|^2 +|\nab_4\pmb\psi|^2+|\nabla\pmb\psi|^2+|\pmb\psi|^2\big)(\tau, r=r_+(1-\dhor), \om){d\mathring{\ga}d\tau},\\
\Fh_{\II_+}[\pmb\psi](\tau_1,\tau_2):=& \int_{\tau_1}^{\tau_2}\int_{\mathbb{S}^2}\Bigg(|\nab_{\pr_\tau^{\II_+}}\pmb\psi|^2\\
&\qquad\qquad\quad +r^{-2}\left(|\nabla_{\pr_{x^1}^{\II_+}}\pmb\psi|^2+|\nabla_{\pr_{x^2}^{\II_+}}\pmb\psi|^2+|\pmb\psi|^2\right)\Bigg)(\underline{\tau}=+\infty, \tau, \om)r^2d\mathring{\ga}d\tau,\\
\Fh[\pmb\psi](\tau_1,\tau_2):=& \Fh_{\II_+}[\pmb\psi](\tau_1,\tau_2)+\Fh_{\AA}[\pmb\psi](\tau_1,\tau_2),
\end{split}
\eea
the energy norm
\bea
\Eh[\pmb\psi](\tau):= \int_{r_+(1-\dhor)}^{+\infty}\int_{\mathbb{S}^2}\Big(|\nab_4\pmb\psi|^2+|\nab\pmb\psi|^2+r^{-2}|\nab_3\pmb\psi|^2+r^{-2}|\pmb\psi|^2\Big)r^2{d\mathring{\ga}dr},
\eea
and the Morawetz norms
\bea
\bsplit
\Mh_\de[\pmb\psi](\tau_1,\tau_2):=&\int_{\Mntraph(\tau_1,\tau_2)}\left(\frac{|\nab_{\partial_{\tt}}\pmb\psi|^2}{r^{1+\de}} +\frac{|\nabla\pmb\psi|^2}{r}\right)
+\int_{\MMh(\tau_1,\tau_2)}\left(\frac{|\nab_{\partial_r}\pmb\psi|^2}{r^{1+\de}} +\frac{|\pmb\psi|^2}{r^3}\right),\\
\Mh[\pmb\psi](\tau_1,\tau_2):=&\int_{\Mntraph(\tau_1,\tau_2)}\left(\frac{|\nab_{\partial_{\tt}}\pmb\psi|^2}{r^2} +\frac{|\nabla\pmb\psi|^2}{r}\right)
+\int_{\MMh(\tau_1,\tau_2)}\left(\frac{|\nab_{\partial_r}\pmb\psi|^2}{r^2} +\frac{|\pmb\psi|^2}{r^3}\right),
\end{split}
\eea
\esub
for any given $0\leq\de\leq 1$, where $\Mh[\pmb\psi](\tau_1,\tau_2)=\Mh_1[\pmb\psi](\tau_1,\tau_2)$. We also define, for $\H\in\sk_k(\mathbb{C})$, $k=1,2$,
\bea\lab{eq:defmathcalNpsif}
&&\NNh_\de[\pmb\psi, \H](\tau_1, \tau_2)\nn\\ 
&:=& \!\!\int_{\MMh(\tau_1, \tau_2)}r^{1+\de}|\H|^2\nn\\
&+&\!\!\!\!\min\left[\left(\int_{\Mtraph(\tau_1, \tau_2)}|\H|^2\right)^{\frac{1}{2}} \left(\int_{\Mtraph(\tau_1, \tau_2)}|\dk{^{\leq 1}}\pmb\psi|^2\right)^{\frac{1}{2}}, 
\int_{\Mtraph(\tau_1, \tau_2)}\tau^{1+\de}|\H|^2\right]\!.
\eea

Next, for any nonnegative integer $\reg$, let
\bea\lab{eq:definitionofhigherordernormsFregEregMdeltaregMregNNdeltaregwidehatNreg}
\bsplit
\Fh^{(\reg)}[\pmb\psi](\tau_1,\tau_2)&:=\Fh[\dk^{\leq\reg}\pmb\psi](\tau_1,\tau_2),\qquad\qquad\qquad\quad \Eh^{(\reg)}[\pmb\psi](\tau):=\Eh[\dk^{\leq\reg}\pmb\psi](\tau), \\ 
\Mh^{(\reg)}_\de[\pmb\psi](\tau_1,\tau_2)&:=\Mh_\de[\dk^{\leq\reg}\pmb\psi](\tau_1,\tau_2),\qquad\qquad\, \Mh^{(\reg)}[\pmb\psi](\tau_1,\tau_2):=\Mh[\dk^{\leq\reg}\pmb\psi](\tau_1,\tau_2),\\
\NNh^{(\reg)}_\de[\pmb\psi, \H](\tau_1, \tau_2)&:=\NNh_\de[\dk^{\leq\reg}\pmb\psi, \dk^{\leq \reg}\H](\tau_1, \tau_2).
\end{split}
\eea

Finally, we define for any nonnegative integer $\reg$ the following combined norms 
\bea\lab{eq:definitionofhigherordernormsFregEregMdeltaregMregNNdeltaregwidehatNreg:combinednorms}
\bsplit
\EMFh^{(\reg)}_\de[\pmb\psi](\tau_1,\tau_2) := \sup_{\tt\in [\tau_1, \tau_2]} \Eh^{(\reg)}[\pmb\psi](\tt) + \Mh^{(\reg)}_\de[\pmb\psi](\tau_1,\tau_2)+\Fh^{(\reg)}[\pmb\psi](\tau_1,\tau_2),\\
\EMFh^{(\reg)}[\pmb\psi](\tau_1,\tau_2) := \sup_{\tt\in [\tau_1, \tau_2]} \Eh^{(\reg)}[\pmb\psi](\tt) + \Mh^{(\reg)}[\pmb\psi](\tau_1,\tau_2)+\Fh^{(\reg)}[\pmb\psi](\tau_1,\tau_2),
\end{split}
\eea
with $\EMh^{(\reg)}_\de[\pmb\psi](\tau_1,\tau_2)$, $\EMh^{(\reg)}[\pmb\psi](\tau_1,\tau_2)$, $\MFh^{(\reg)}[\pmb\psi](\tau_1,\tau_2)$ and $\EFh^{(\reg)}[\pmb\psi](\tau_1,\tau_2)$ being defined in a similar way.

%%%%%%%%%%%%%%%%%%%%%%%%%%%%%%%%%%%%%%%%%%%%%%%%%%%%%%%%%

\paragraph{\textit{Energy, Morawetz and flux norms for $\pmb\psi\in\sk_k(\mathbb{C})$, $k=1,2$, on $\MM$}.}

%%%%%%%%%%%%%%%%%%%%%%%%%%%%%%%%%%%%%%%%%%%%%%%%%%%%%%%%%

On $\MM$, we define the norms $\F_\AA[\pmb\psi]$, $\F_{\Si_*}[\psi]$, $\F[\pmb\psi]$, $\E[\pmb\psi]$, $\M_\de[\pmb\psi]$, $\M[\pmb\psi]$, $\NNh_\de[\pmb\psi, \H]$, as well as the corresponding higher order derivative norms and combined quantities as above by replacing:
\begin{itemize}
\item $\MMh$, $\Mtraph$, $\Mntraph$ and $\Sih(\tau)$ respectively with $\MM$, $\Mtrap$, $\Mntrap$ and $\Si(\tau)$,

\item $r<+\infty$ with $r\leq\varsigma_*(\tau-c_*)$ and $\tau_1<\tau_2$ with $1\leq\tau_1<\tau_2\leq\tau_*$, 

\item $\II_+$ with $\Si_*$ and $\Fh_{\II_+}[\pmb\psi](\tau_1, \tau_2)$ with $\F_{\Si_*}[\pmb\psi](\tau_1, \tau_2)$ where 
\bea
\F_{\Si_*}[\pmb\psi](\tau_1,\tau_2) :=  \int_{\tau_1}^{\tau_2}\int_{\mathbb{S}^2}\Big( |\nab_4\pmb\psi |^2+|\nab_3\pmb\psi|^2+|\nab\pmb\psi|^2+ r^{-2} |\pmb\psi|^2\Big)(\tau, r=\varsigma_*(\tau-c_*), \om)r^2d\mathring{\ga}d\tau.
\eea
\end{itemize}

In addition,  we introduce the following generalized fluxes
\bea\lab{eq:definitionofgeneralizedfluxFSi*commac}
\nn\F_{\Si_{*,c}}[\pmb\psi](\tau_1,\tau_2) &:=&  \int_{\tau_1}^{\tau_2}\int_{\mathbb{S}^2}\Big( |\nab_4\pmb\psi |^2+|\nab_3\pmb\psi|^2+|\nab\pmb\psi|^2+ r^{-2} |\pmb\psi|^2\Big)\\
&&\qquad\qquad\qquad\qquad\qquad\qquad (\tau, r=\varsigma_*(\tau-c), \om)r^2d\mathring{\ga}d\tau, \quad \frac{c_*}{2}\leq c\leq c_*,
\eea
where we have in particular 
\beaa
\F_{\Si_*}[\pmb\psi](\tau_1,\tau_2)=\F_{\Si_{*,c_*}}[\pmb\psi](\tau_1,\tau_2).
\eeaa

%%%%%%%%%%%%%%%%%%%%%%%%%%%%%%%%%%%%%%%%%%%%%%

\paragraph{\textit{Energy, Morawetz and flux norms for scalars on $\MM$ and $\MMh$}.}

%%%%%%%%%%%%%%%%%%%%%%%%%%%%%%%%%%%%%%%%%%%%%%

For any scalars $\psi$ and $H$, we define the norms $\Fh_\AA[\psi]$, $\Fh_{\II_+}[\psi]$, $\Fh[\psi]$, $\Eh[\psi]$, $\Mh_\de[\psi]$, $\Mh[\psi]$, $\NNh_\de[\psi, H]$, $\F_\AA[\psi]$, $\F_{\Si_*}[\psi]$, $\F_{\Si_{*,c}}[\psi]$, $\E[\psi]$, $\M_\de[\psi]$, $\M[\psi]$, $\NN_\de[\psi, H]$, as well as the corresponding higher order derivative norms and combined quantities as above by replacing $\pmb\psi, \,\H\in\sk_k(\mathbb{C})$, $k=1,2$, with scalars $\psi$, $H$ in the formulas.

%%%%%%%%%%%%%%%%%%%%%%%%%%%%%%%%%%%%%%%%%%%%%%%%%

\paragraph{\textit{Energy, Morawetz and flux norms for scalarized tensors using regular triplets}.}

%%%%%%%%%%%%%%%%%%%%%%%%%%%%%%%%%%%%%%%%%%%%%%%%%

We start with the following lemma.
\begin{lemma}\lab{lemma:equivalentofmodulussquarepmbpsiandsumijpisijalsowithdkreg}
Let $\pmb\psi\in\sk_2(\mathbb{C})$ on $\MM$ and let $\psi_{ij}$ be the corresponding scalars given by $\psi_{ij}:=\pmb\psi(\Om_i, \Om_j)$ for $i,j=1,2,3$ where $\Om_i$, $i=1,2,3$ is the regular triplet introduced in Section \ref{sec:assumptionsforsec:energyMorawetzesitmatesforTeukoslkyonMM:upto15derivatives}. Then, for any integer $\reg\leq {15}$, we have, for $\a=a, 3, 4$, $a=1,2$, 
\beaa
|\dk^{\leq \reg}\pmb\psi|^2\simeq \sum_{i,j=1}^3|\dk^{\leq\reg}(\psi_{ij})|^2, \qquad |\Ddot_\a\dk^{\leq\reg}\pmb\psi|^2\simeq \sum_{i,j=1}^3|e_\a(\dk^{\leq\reg}(\psi_{ij}))|^2+O(r^{-2})\sum_{i,j=1}^3|\dk^{\leq\reg}(\psi_{ij})|^2.
\eeaa
\end{lemma}

\begin{proof}
See Lemma 5.31 in \cite{MaSz26}.
\end{proof}

Then, we infer  immediately from Lemma \ref{lemma:equivalentofmodulussquarepmbpsiandsumijpisijalsowithdkreg} the following equivalence relations between the norms for tensors and the norms for the scalars obtained by scalarization using regular triplets.
\begin{corollary}
[Equivalence of the norms for tensors and for scalars using regular triplets]
\lab{coro:equivalenceofnormsfortensorsandscalars}
Let $\pmb\psi, \,\H\in\sk_2(\mathbb{C})$ on $\MM$, and let $\psi_{ij}, H_{ij}$ be the corresponding scalars given by $\psi_{ij}:=\pmb\psi(\Om_i, \Om_j)$ and $H_{ij}:=\H(\Om_i, \Om_j)$, for $i,j=1,2,3$, respectively, where $\Om_i$, $i=1,2,3$ is the regular triplet introduced in Section \ref{sec:assumptionsforsec:energyMorawetzesitmatesforTeukoslkyonMM:upto15derivatives}. 
Then we have the equivalence relations, for any $\reg\leq 14$, $\de\in[0,1]$ and $\tau_1<\tau_2$,
\bea\lab{eq:equivalenceofnormsfortensorsandscalars}
\E^{(\reg)}[\pmb\psi](\tau)\simeq \sum_{i,j=1}^3\E^{(\reg)}[\psi_{ij}](\tau), \qquad \mathcal{N}^{(\reg)}_\de[\pmb\psi, \H](\tau_1, \tau_2)\simeq \sum_{i,j=1}^3\mathcal{N}^{(\reg)}_\de[\psi_{ij}, H_{ij}](\tau_1, \tau_2),
\eea
and similarly for all the other norms for $\pmb\psi, \,\H\in\sk_2(\mathbb{C})$ appearing in this section.
\end{corollary}

\begin{remark}\lab{rmk:abusenotationbetweentensorandscalarizedversioninEMFnorms}
More generally, we will also consider family of scalars $\psi_{ij}$, $i,j=1,2,3,$ on $\MM$ that are not generated\footnote{Recall from Lemma \ref{lemma:backandforthbetweenhorizontaltensorsk2andscalars:complex} that complex-valued scalars $\psi_{ij}$, $i,j=1,2,3,$ come from the scalarization of a tensor in $\sk_2(\mathbb{C})$ if and only if they satisfy the identities stated in the second item of Lemma \ref{lemma:backandforthbetweenhorizontaltensorsk2andscalars:complex}.} by the scalarization of a tensor in $\sk_2(\mathbb{C})$. In that case\footnote{This will be the case in $\MM(\tau_2-2, \tau_2)$ and in $\RR_*(\tau_1, \tau_2)$ due to the semi global extension procedure of Proposition  \ref{prop:extensionprocedureoftheTeukolskywaveequations} which does not preserve the identities stated in the second item of Lemma \ref{lemma:backandforthbetweenhorizontaltensorsk2andscalars:complex}.}, by a slight abuse of notation, we will still denote the norms appearing on the RHS of the identities in  \eqref{eq:equivalenceofnormsfortensorsandscalars} by the ones appearing on the LHS (even though the corresponding tensors $\pmb\psi$, $\H$, do not exist).
\end{remark}

%%%%%%%%%%%%%%%%%%%%%%%%%%%%%%%%%%%%%%%%%%%%

\subsubsection{Basic estimates for wave equations in perturbations of Kerr}
\lab{sect:basicestimatesforwaveequations}

%%%%%%%%%%%%%%%%%%%%%%%%%%%%%%%%%%%%%%%%%%%%

%%%%%%%%%%%%%%%%%%%%%%%%%%%%%%%%%%

\paragraph{\textit{Standard calculation for generalized currents}.}

%%%%%%%%%%%%%%%%%%%%%%%%%%%%%%%%%% 

Consider  variational  wave equations  for  tensors  $\pmb\psi\in \sk_k(\mathbb{C})$, $k=0,1,2$, of the form 
\bea\lab{eq:Gen.RW-general}
\squared_k \pmb\psi-V\pmb\psi=\pmb{N},
\eea
where $V$ is a real potential. The variational wave equation \eqref{eq:Gen.RW-general} has Lagrangian
\beaa
 \LL[\pmb\psi]:= \g^{\mu\nu}\Re\left(\Db_\mu \pmb\psi\c\ov{\Db_\nu  \pmb\psi }\right)+ V |\pmb\psi |^2,
 \eeaa
  where  the dot product   here denotes full contraction with respect to the  horizontal indices.
The  corresponding   energy-momentum tensor associated to \eqref{eq:Gen.RW-general} is given by 
 \bea\label{eq:definition-QQ-mu-nu}
\nn  \QQ_{\mu\nu}[\pmb\psi] &:=& \Re\Big(\Db_\mu  \pmb\psi \c \ov{\Db _\nu  \pmb\psi }\Big)
          -\frac 12 \g_{\mu\nu} \left( \Re\left(\Db_\la  \pmb\psi\c\ov{\Db^\la  \pmb\psi }\right)+ V |\pmb\psi |^2\right)\\
&=& \Re\Big(\Db_\mu  \pmb\psi \c \ov{\Db _\nu  \pmb\psi }\Big)  -\frac 12 \g_{\mu\nu}  \LL[\pmb\psi]. 
 \eea
Also, recall from \eqref{def:deformationtensor:lastsect} that the deformation tensor of a vectorfield $X$ is defined by
\beaa
{}^{(X)}\pi_{\a\b} =  \D_{\a}X_{\b} + \D_{\b}X_{\a}.
\eeaa

We collect here some general calculations for generalized currents associated to equation \eqref{eq:Gen.RW-general}.

\begin{proposition}\lab{prop-app:stadard-comp-Psi}
 Let   $\pmb\psi\in \mathfrak{s}_k(\mathbb{C})$, $k=0,1,2$, and  let $X$ be  a real-valued vectorfield and $w$ a real scalar function.  Define the $1$-form $\PP_\mu[\pmb\psi](X, w)$ by
 \bea
\lab{definitionofcurrentPPmuXw:generaltensor}
\PP_\mu[\pmb\psi](X, w):=\QQ_{\mu\nu}[\pmb\psi] X^\nu +\frac 1 2  w \Re\Big(\pmb\psi \c \ov{\Db_\mu \pmb\psi }\Big)-\frac 1 4|\pmb\psi|^2   \pr_\mu w,
  \eea
and define the 1-form ${}^{(X)}A_\nu$  by
\bea\lab{eq:thespacetime1formXA}
{}^{(X)}A_\nu &:=& X^\mu \in^{ab}\Rdot_{ ab   \nu\mu}.
\eea
Then, we have
 \bea\lab{eq:DivofPPmu:tensor:RW:prop}
  \D^\mu  \PP_\mu[\pmb\psi](X, w)&=& \frac 1 2 \QQ[\pmb\psi]  \c\piX - \frac 1 2 X( V ) |\pmb\psi|^2 +\frac{k}{2}{}^{(X)}A_\nu\Im\Big(\pmb\psi\c\ov{\Ddot^{\nu}\pmb\psi}\Big)\nn\\
  &&+\frac 12  w \LL[\pmb\psi] -\frac 1 4|\pmb\psi|^2   \square_\g  w   +  \Re\bigg(\ov{\bigg(\nab_X\pmb\psi +\frac 1 2   w \pmb\psi\bigg)}\c \left(\squared_k \pmb\psi- V\pmb\psi\right)\bigg).
 \eea
\end{proposition}

\begin{proof}
See Proposition 6.2 in \cite{MaSz26}.
\end{proof}

In order to control the term of the type $\pmb\psi\c\Ddot_\nu\pmb\psi$ in \eqref{eq:DivofPPmu:tensor:RW:prop}, we will rely on the following corollaries of Lemma 6.3 in \cite{MaSz26}.

\begin{corollary}\lab{cor:computationofthecomponentsofthetensorAforexactenergyconservationKerr}
Let $(\MM, \g)$ satisfy the assumptions of Section \ref{sec:assumptionsforsec:energyMorawetzesitmatesforTeukoslkyonMM:upto15derivatives}. Then, we have
\bea
{}^{({\pr_\tau})}A_\mu=-\D_\mu\left(\Im\left(\frac{2m}{q^2}\right)\right) {+r^{-1}\Ga_b}.
\eea
\end{corollary}

\begin{proof}
See Corollary 6.4 in \cite{MaSz26}.
\end{proof}

\begin{corollary}\lab{cor:asymtpticbehavioroftheRdottermintensorialenergyidentity:larger}
Let $(\MM, \g)$ satisfy the assumptions of Section \ref{sec:assumptionsforsec:energyMorawetzesitmatesforTeukoslkyonMM:upto15derivatives}. For $X$ such that 
\beaa
X^4=O(1), \qquad X^3=O(1), \qquad X^b=O(r^{-1}), \,\, b=1,2,
\eeaa
we have 
\beaa
{}^{(X)}A_\nu\Im\Big(\pmb\psi\c\ov{\Ddot^{\nu}\pmb\psi}\Big) = O(r^{-3})\Im\Big(\pmb\psi\c\ov{\nab_4\pmb\psi }\Big)+ O(r^{-3})\Im\Big(\pmb\psi\c\ov{\nab_3\pmb\psi }\Big) + \Big(O(r^{-3})+r^{-1}\Ga_b\Big)\Im\Big(\pmb\psi\c\ov{\nab\pmb\psi }\Big).  
\eeaa
\end{corollary}

\begin{proof}
See Corollary 6.5 in \cite{MaSz26}.
\end{proof}

Next, we consider solutions to the following tensorial wave equation for $\pmb\psi\in \sk_k(\mathbb{C})$, $k=0,1,2$,
\bea
\lab{eq:tensorialwaveRW:withprttderivative}
 \squared_k \pmb\psi-\frac{4ia\cos\th}{\qs}\nab_{\pr_{\tt}}\pmb\psi-V\pmb\psi=\pmb{N}
 \eea
and, as in Section  7.3 of \cite{GKS22}, we make use of Corollary \ref{cor:computationofthecomponentsofthetensorAforexactenergyconservationKerr} to derive an energy identity.

\begin{lemma}\lab{cor:modifiedcurrentsforprtandprvphi}
Let $(\MM, \g)$ satisfy the assumptions of Section \ref{sec:assumptionsforsec:energyMorawetzesitmatesforTeukoslkyonMM:upto15derivatives}. Let $\pmb\psi\in\sk_k(\mathbb{C})$, $k=0,1,2$, be a solution to the tensorial wave equation \eqref{eq:tensorialwaveRW:withprttderivative}, with the real potential $V$ satisfying $\pr_{\tt}V=0$.
Define 
\bea
\widetilde{w} := \Im\left(\frac{m}{q^2}\right)= -\frac{2amr\cos\th}{|q|^4},
\eea
and define the following modified current associated to the vectorfield $\pr_{\tt}$:
\bea
{}^{(\pr_{\tt})}\widetilde{\PP}_\mu[\pmb\psi] :=  {\PP}_\mu[\pmb\psi](\pr_{\tt}, 0) + k\widetilde{w}\Im\left( \pmb\psi\c\ov{\Ddot_\mu\pmb\psi}\right)+(\pr_{\tt})_{\mu}\frac{2a\cos\th}{\qs}k\widetilde{w}|\pmb\psi|^2.
\eea
Then, we have
\bea
\D^\mu{}^{(\pr_{\tt})}\widetilde{\PP}_\mu [\pmb\psi]&=&  \Re\bigg(\ov{\big(\nab_{\pr_{\tt}}\pmb\psi-ik\widetilde{w}\pmb\psi\big)}\c \left(\squared_k\pmb\psi-\frac{4ia\cos\th}{\qs}\nab_{\pr_{\tt}}\pmb\psi-V\pmb\psi\right)\bigg)\nn\\
&&+\frac 1 2 \QQ[\pmb\psi]  \c {^{(\pr_{\tt})}}\pi+\Div (\pr_{\tt})\frac{2a\cos\th}{\qs}k\widetilde{w}|\pmb\psi|^2+r^{-1}\Ga_b\Im\Big(\pmb\psi\c\ov{\Ddot^{\nu}\pmb\psi}\Big).
\eea
\end{lemma}

\begin{remark}
In the case $\g=\gam$ and $\pmb{N}=0$, this induces a conservation of energy.
\end{remark}

\begin{proof}
See Lemma 6.6 in \cite{MaSz26}.
\end{proof}

We have the following corollary of Lemma \ref{cor:modifiedcurrentsforprtandprvphi}.
\begin{corollary}\lab{cor:energyestimatetensoriallevel:scalarizedversion}
Let $(\MM, \g)$ satisfy the assumptions of Section \ref{sec:assumptionsforsec:energyMorawetzesitmatesforTeukoslkyonMM:upto15derivatives}. Let $\pmb\psi\in\sk_2(\mathbb{C})$, and let $\psi_{ij}$ be the corresponding scalars given by $\psi_{ij}:=\pmb\psi(\Om_i, \Om_j)$, for $i,j=1,2,3$, where $\Om_i$, $i=1,2,3$ is the regular triplet introduced in Section \ref{sec:assumptionsforsec:energyMorawetzesitmatesforTeukoslkyonMM:upto15derivatives}. Under the assumptions of Lemma \ref{cor:modifiedcurrentsforprtandprvphi}, and given ${}^{(\pr_{\tt})}\widetilde{\PP}_\mu[\pmb\psi]$ and $\widetilde{w}$ defined as in Lemma \ref{cor:modifiedcurrentsforprtandprvphi}, we have
\bea
\D^\mu{}^{(\pr_{\tt})}\widetilde{\PP}_\mu [\pmb\psi] &=& \Re\left(\Big(\square_{\g}(\psi^{ij})+\widehat{S}(\psi)^{ij}+(\widehat{Q}\psi)^{ij}-V\psi^{ij}\Big)\ov{\Big(\pr_{\tt}(\psi_{ij}) - M_{i\tt}^k\psi_{kj} - M_{j\tt}^k\psi_{ik} -i2\widetilde{w}\psi_{ij}\Big)}\right)\nn\\
&&+\frac 1 2 \QQ[\pmb\psi]  \c {^{(\pr_{\tt})}}\pi+\Div (\pr_{\tt})\frac{4a\cos\th}{\qs}\widetilde{w}|\pmb\psi|^2+r^{-1}\Ga_b\Im\Big(\pmb\psi\c\ov{\Ddot^{\nu}\pmb\psi}\Big).
\eea
\end{corollary}

\begin{proof}
See Corollary 6.8 in \cite{MaSz26}.
\end{proof}

%%%%%%%%%%%%%%%%%%%%%%%%%%%

\paragraph{\textit{Control of error terms}.}

%%%%%%%%%%%%%%%%%%%%%%%%%%% 

The following two lemmas, taken respectively from \cite[Lemma 3.3]{MaSz24} and \cite[Lemma 3.5]{MaSz24}, will allow us to control the error terms arising in the derivation of energy-Morawetz estimates in $\MM(\tau_1,\tau_2)$.  

\begin{lemma}\lab{lemma:basiclemmaforcontrolNLterms:ter}
Let $h\in r^{-1}\dk^{\leq 1}\Ga_b$ be a scalar function and let 
$M^{\a\b}$ be symmetric and satisfy
\beaa
&& M^{rr}\in r\dk^{\leq 1}\Ga_b, \qquad M^{r\tau}\in r\dk^{\leq 1}\Ga_g, \qquad M^{\tau\tau}\in r\dk^{\leq 1}\Ga_g
, \qquad |M^{\tau\tau}|\les\frac{\ep}{r^2},\\
&& M^{rx^a}\in \dk^{\leq 1}\Ga_b, \qquad M^{\tau x^a}\in \dk^{\leq 1}\Ga_g,  \qquad M^{x^ax^b}\in r^{-1}\dk^{\leq 1}\Ga_g,
\eeaa
where $a,b=1,2$. Then, the following estimate holds
\beaa
\int_{\MM(\tau_1, \tau_2)}\Big(\big|M^{\a\b}\pr_\a\psi\pr_\b\psi\big|+h|\psi|^2\Big) &\les& \ep\EM[\psi](\tau_1, \tau_2).
\eeaa
\end{lemma}

\begin{lemma}\lab{lemma:basiclemmaforcontrolNLterms:bis}
Let $M^{\a\b}$ be symmetric and satisfy
\beaa
&& M^{rr}\in r\dk^{\leq 1}\Ga_b, \qquad M^{r\tau}\in r\dk^{\leq 1}\Ga_g, \qquad M^{\tau\tau}\in r\dk^{\leq 1}\Ga_g, \qquad |M^{\tau\tau}|\les\frac{\ep}{r^2},\\
&& M^{rx^a}\in \dk^{\leq 1}\Ga_b, \qquad M^{\tau x^a}\in \dk^{\leq 1}\Ga_g,  \qquad M^{x^ax^b}\in r^{-1}\dk^{\leq 1}\Ga_g,
\eeaa
where $a,b=1,2$. Then, the following estimate holds
\beaa
\int_{\MM(\tau_1, \tau_2)}\big|M^{\a\b}\pr_\a\pr_\b\psi\big|^2
+\bigg|\int_{\MM(\tau_1, \tau_2)}M^{\a\b}\pr_\a\pr_\b\psi \pr_\tau(\pr^{\leq 1}\psi)\bigg|\\
+\int_{\MM(\tau_1, \tau_2)}r^{-1}\big|M^{\a\b}\pr_\a\pr_\b\psi\big| \big|\dk^{\leq 1}\pr^{\leq 1}\psi\big| &\les& \ep\EM[\pr^{\leq 1}\psi](\tau_1, \tau_2). 
\eeaa

Also, let $N$ be a spacetime vectorfield such that we have 
\beaa
N^r\in r\dk^{\leq 2}\Ga_g, \qquad N^\tau\in r\dk^{\leq 2}\Ga_g, \qquad |N^\tau|\les\frac{\ep}{r^2}, \qquad N^{x^a}\in \dk^{\leq 2}\Ga_g.
\eeaa
Then, the following holds
\beaa
\int_{\MM(\tau_1, \tau_2)}\big|N^\a\pr_\a\psi\big|^2+ \int_{\MM(\tau_1, \tau_2)}\big|N^\a\pr_\a\psi\big| \Big|\big(\pr_\tau, \pr_r, r^{-1}\pr_{x^a}, r^{-1}\big)\pr^{\leq 1}\psi\Big| \les \ep\EM[\pr^{\leq 1}\psi](\tau_1, \tau_2). 
\eeaa
\end{lemma}

\begin{remark}
In practice, concerning the quantities estimated in Lemma \ref{lemma:basiclemmaforcontrolNLterms:bis}:
\begin{itemize}
\item $\big(\pr_\tau, \pr_r, r^{-1}\pr_{x^a}, r^{-1}\big)\pr^{\leq 1}\psi$ will be due to  energy-Morawetz multipliers, 

\item  $M^{\a\b}\pr_\a\pr_\b\psi$ and $N^\a\pr_\a\psi$ will come from the RHS of the wave equation, in particular after commutation with {various vectorfields such as} $\pr_\tau$.
\end{itemize}
\end{remark}

%%%%%%%%%%%%%%%%%%%%%%%%%%%%%

\paragraph{\textit{Commutators with the D'Alembertian}.}

%%%%%%%%%%%%%%%%%%%%%%%%%%%%%

The following two lemmas provide the structure of commutators between first-order derivatives and the scalar wave operator, respectively for unweighted and weighted derivatives.

\begin{lemma}
\lab{lem:commutatorwithwave:firstorderderis:0}
Let $(\MM, \g)$ satisfy the assumptions of Section \ref{sec:assumptionsforsec:energyMorawetzesitmatesforTeukoslkyonMM:upto15derivatives}. Then, the  commutator between $\square_{\g}$ and $\pr_{\tau}$ satisfies 
\bea
\label{esti:commutatorBoxgandT:general:0}
\, [ \pr_{\tau}, \square_{\g}]\psi = \pr_{\tau}(\gcheck^{\a\b})\pr_{\a}\pr_{\b}\psi  +\pr_{\tau}^2(\gcheck^{\tau\tau})\pr_\tau\psi  +
\dk^{\leq 2}\Ga_g\c\dk\psi,
\eea
and the commutator between $\square_\g$ and $(\pr_r,r^{-1}\pr_{x^a})$ satisfies 
\bea
\lab{eq:localwavecommutators:withfirstordergoodderis:0}
{[(\pr_r,r^{-1}\pr_{x^a}), \square_{\g}]\psi} = O(r^{-2}) (\pr_{\tt}, \pr_{r}, \pr_{x^a})^{\leq 1}\pr\psi 
+\dk^{\leq 1}\Ga_g\dk\pr\psi+r^{-1}\dk^{\leq 2}\Ga_g\dk \psi.
\eea
\end{lemma}

\begin{proof}
See Lemma 6.14 in \cite{MaSz26}.
\end{proof}

\begin{lemma}
\lab{lem:commutatorwithwave:firstorderderis}
Let $(\MM, \g)$ satisfy the assumptions of Section \ref{sec:assumptionsforsec:energyMorawetzesitmatesforTeukoslkyonMM:upto15derivatives}. Then, the commutator between $\square_{\g}$ and $\pr_{\tau}$ satisfies
\bea
\label{esti:commutatorBoxgandT:general}
\, [ \pr_{\tau}, \square_{\g}]\psi =\pr_{\tau}(\dk(\gcheck^{\tau\tau})\pr_{\tau}\psi)+  \dk^{\leq 1}(\dk^{\leq 1}\Ga_g\c\dk\psi),
\eea
and the commutator between $\square_\g$ and $(r\pr_r,\pr_{x^a})$ satisfies 
\bea
\lab{eq:localwavecommutators:withfirstordergoodderis}
{[(r\pr_r, \pr_{x^a}), \square_{\g}]\psi} = \big(-\square_{\g}\psi, 0\big)+O(r^{-2})\dk^{\leq 1}\dk\psi +\pr_{\tau}(\gcheck^{\tau\tau}\pr_{\tau}\psi) + \dk^{\leq 1}(\dk^{\leq 1}\Ga_g\c\dk\psi).
\eea
\end{lemma}

\begin{proof}
See Lemma 6.15 in \cite{MaSz26}.
\end{proof}

%%%%%%%%%%%%%%%%%%%%%%%%%%%

\paragraph{\textit{Local energy estimates}.}

%%%%%%%%%%%%%%%%%%%%%%%%%%%

We have the following basic local (in time) energy estimates for systems of wave equations.

\begin{lemma}[Local energy estimate]
\lab{lemma:localenergyestimate}
Let $\g$ satisfy the assumptions of Section \ref{sec:assumptionsforsec:energyMorawetzesitmatesforTeukoslkyonMM:upto15derivatives}. Let $(\psi)_{ij}$, $i,j=1,2,3$,  satisfy the following coupled system of scalar wave equations
\bea
\lab{eq:eqsforlocalenergyestimatelemma:general}
\square_{\g}\psi_{ij}={N_{ij}, \qquad N_{ij}:=D_1r^{-1}\pr_{\tt} \psi_{ij}+O(r^{-2})\dk^{\leq 1}\psi_{kl} +F_{ij},}
\eea
with the constant $D_1\geq 0$.
For any $\tau_0\in\mathbb{R}$, $0\leq\reg\leq 14$, $\de\in (0,1]$ and $q>0$, we have the following future directed local energy estimates
\bsub
\bea
\label{eq:localenergyestimate:future}
\EMF_\de^{(\reg)}[\psi](\tau_0, \tau_0+q) &\les_q & E^{(\reg)}[\psi](\tau_0)  +\NNtlede^{(\reg)}[\psi, \F](\tau_0, \tau_0+q),\\
\EMF_\de^{(\reg)}[\psi](\tau_0, \tau_0+q) &\les_q & E^{(\reg)}[\psi](\tau_0) +\sum_{i,j}\int_{\MM(\tau_0, \tau_0+q)}r^{1+\de}|\dk^{\leq \reg}F_{ij}|^2,\label{eq:localenergyestimate:future:bis}
\eea
\esub
where for any $\tau'<\tau''$, 
\begin{align}\lab{def:NNtleinlocalenergyestimate}
\NNtlede^{(\reg)}[\psi, \F](\tau', \tau'')&:=\sum_{i,j}\bigg(\int_{\MM_{r\geq 10m}(\tau', \tau'')}r^{-1}|\dk^{\leq \reg}F_{ij}| |\dk^{\leq \reg+1}\psi_{ij}| 
+\int_{\MM(\tau', \tau'')}|\dk^{\leq \reg}F_{ij}|^2\bigg)  \nn\\
+& \sup_{\tau'''\in[\tau', \tau'']}\bigg|\sum_{i,j}\int_{\MM_{r\geq 10m}(\tau', \tau''')} \Re\Big(\ov{\dk^{\leq \reg}F_{ij}}\big(1+O(r^{-\de})\big)\pr_{\tt} \dk^{\leq \reg}\psi_{ij}\Big)\bigg| .
\end{align}
\end{lemma}

\begin{proof}
See Lemma 6.19 in \cite{MaSz26}.
\end{proof}

The following lemma is the analog of Lemma \ref{lemma:localenergyestimate} upon replacing weighted derivatives $\dk$ with unweighted derivatives $\pr$.
\begin{lemma}[Local energy estimate with unweighted derivatives]
\lab{lemma:localenergyestimate:unweigthedderivatives}
Let $\g$ satisfy the assumptions of Section \ref{sec:assumptionsforsec:energyMorawetzesitmatesforTeukoslkyonMM:upto15derivatives}. Let $(\psi)_{ij}$, $i,j=1,2,3$,  satisfy the coupled system of scalar wave equations \eqref{eq:eqsforlocalenergyestimatelemma:general} with the constant $D_1\geq 0$.
For any $\tau_0\in\mathbb{R}$, $0\leq\reg\leq 14$, $\de\in (0,1]$ and $q>0$, we have the following future directed local energy estimates
\bsub
\bea
\label{eq:localenergyestimate:future:unweigthedderivatives}
\EMF_\de[\pr^{\leq\reg}\psi](\tau_0, \tau_0+q) &\les_q & E[\pr^{\leq\reg}\psi](\tau_0)  +\NNtlede[\pr^{\leq\reg}\psi, \pr^{\leq\reg}\F](\tau_0, \tau_0+q),\\
\EMF_\de[\pr^{\leq\reg}\psi](\tau_0, \tau_0+q) &\les_q & E[\pr^{\leq\reg}\psi](\tau_0) +\sum_{i,j}\int_{\MM(\tau_0, \tau_0+q)}r^{1+\de}|\pr^{\leq\reg}F_{ij}|^2,\label{eq:localenergyestimate:future:unweigthedderivatives:bis}
\eea
\esub
and the following past directed local energy estimates 
\bsub
\begin{align}
\label{eq:localenergyestimate:past:unweigthedderivatives}
\EMF_\de[\pr^{\leq\reg}\psi](\tau_0-q, \tau_0) \les_q & \,E[\pr^{\leq\reg}\psi](\tau_0) +F[\pr^{\leq\reg}\psi](\tau_0-q, \tau_0)\nn\\
&+ \NNtlede[\pr^{\leq\reg}\psi, \pr^{\leq\reg}\F](\tau_0-q, \tau_0),\\
\EMF_\de[\pr^{\leq\reg}\psi](\tau_0-q, \tau_0) \les_q & \,E[\pr^{\leq\reg}\psi](\tau_0) +F[\pr^{\leq\reg}\psi](\tau_0-q, \tau_0)\nn\\
&+ \sum_{i,j}\int_{\MM(\tau_0-q, \tau_0)}r^{1+\de}|\pr^{\leq \reg}F_{ij}|^2,\label{eq:localenergyestimate:past:bis:unweigthedderivatives}
\end{align}
\esub
where $\NNtlede$ has been introduced in \eqref{def:NNtleinlocalenergyestimate}.
\end{lemma}

\begin{proof}
See Lemma 6.20 in \cite{MaSz26}.
\end{proof}

The following lemma focuses on local energy estimate in the region $\RR_*$.
\begin{lemma}[Local energy estimate in the region $\RR_*$]
\lab{lemma:localenergyestimate:weightedandunweigthedderivatives:RR*}
Let $\g$ satisfy the assumptions of Section \ref{sec:assumptionsforsec:energyMorawetzesitmatesforTeukoslkyonMM:upto15derivatives}. Let $(\psi)_{ij}$, $i,j=1,2,3$,  satisfy the following coupled system of scalar wave equations
\bea
\lab{eq:eqsforlocalenergyestimatelemma:general:versionmoregeneralforRR*}
\square_{\g}\psi_{ij}=N_{ij}, \qquad N_{ij}:=O(r^{-1})\pr\psi_{kl}+O(r^{-2})\psi_{kl} +F_{ij},
\eea
For any $1\leq\tau_1<\tau_2\leq\tau_*$, any $c_*/2\leq c_0\leq c_*$ and $0\leq\reg\leq 14$, we have the following future directed local energy estimates for weighted derivatives 
\bea
\label{eq:localenergyestimate:future:RR*}
\nn&&\EM^{(\reg)}_{0,r\geq\varsigma_*(\tau-c_0)}[\psi](\tau_1, \tau_2)+\sup_{c_0\leq c\leq c_*}\F_{\Si_{*,c}}^{(\reg)}[\psi](\tau_1, \tau_2)\\
 &\les& \E^{(\reg)}[\psi](\tau_1) +\F_{\Si_{*,c_0}}^{(\reg)}[\psi](\tau_1, \tau_2)+ \int_{\MM_{r\geq\varsigma_*(\tau-c_0)}(\tau_1, \tau_2)}r|\dk^{\leq\reg}F|^2,
\eea
and the corresponding ones for unweighted derivatives
\bea
\label{eq:localenergyestimate:future:unweigthedderivatives:RR*}
\nn&&\EM_{0,r\geq\varsigma_*(\tau-c_0)}[\pr^{\leq\reg}\psi](\tau_1, \tau_2)+\sup_{c_0\leq c\leq c_*}\F_{\Si_{*,c}}[\pr^{\leq\reg}\psi](\tau_1, \tau_2)\\
& \les &\E[\pr^{\leq\reg}\psi](\tau_1) +\F_{\Si_{*,c_0}}[\pr^{\leq\reg}\psi](\tau_1, \tau_2)+ \int_{\MM_{r\geq\varsigma_*(\tau-c_0)}(\tau_1, \tau_2)}r|\pr^{\leq\reg}F|^2.
\eea
\end{lemma}

\begin{proof}
We start with the proof of \eqref{eq:localenergyestimate:future:RR*}. We commute \eqref{eq:eqsforlocalenergyestimatelemma:general:versionmoregeneralforRR*} by $(\pr_\tau, r\pr_r, \pr_{{x^a}})^{\leq\reg}$ and obtain in view of Lemma \ref{lem:commutatorwithwave:firstorderderis}, for $\reg\leq 14$, 
\beaa
\square_{\g}\dk^{\leq\reg}\psi_{ij}= O(r^{-1})\pr\dk^{\leq\reg}\psi_{kl}+O(r^{-2})\dk^{\leq\reg}\psi_{kl} +\dk^{\leq\reg}F_{ij},
\eeaa
and then, integrating the current based on the timelike vectorfield $\pr_\tau$ on the spacetime region $\MM_{\varsigma_*(\tau-c_0)\leq r\leq \varsigma_*(\tau-c)}(\tau_1, \tau_2)$ for $c_0\leq c\leq c_*$, we obtain the following analog of (6.34) in \cite{MaSz26}
\bea
\lab{eq:localenergyestimate:zeroorder:general}
\nn&&\E^{(\reg)}_{\varsigma_*(\tau_2-c_0)\leq r\leq \varsigma_*(\tau_2-c)}[\psi](\tau_2)+\F_{\Si_{*,c}}^{(\reg)}[\psi](\tau_1, \tau_2)\nn\\ 
&\les& \E^{(\reg)}[\psi](\tau_1) +\F_{\Si_{*,c_0}}^{(\reg)}[\psi](\tau_1, \tau_2)+\frac{1}{2} \bigg|\int_{\MM_{\varsigma_*(\tau-c_0)\leq r\leq \varsigma_*(\tau-c)}(\tau_1, \tau_2)}{}^{(\pr_\tau)} \pi \cdot \QQ[\dk^{\leq\reg}\psi]\bigg|\\
\nn&&+\bigg| \int_{\MM_{\varsigma_*(\tau-c_0)\leq r\leq \varsigma_*(\tau-c)}(\tau_1, \tau_2)} \Re\Big(\pr_\tau(\dk^{\leq\reg}\psi)\ov{O(r^{-1})\pr\dk^{\leq\reg}\psi_{kl}+O(r^{-2})\dk^{\leq\reg}\psi_{kl} +\dk^{\leq\reg}F_{ij}\big)}\Big) \bigg|.
\eea
Next, we estimate the before to last integral in \eqref{eq:localenergyestimate:zeroorder:general}. Since we have in view of Lemma \ref{lemma:controlofdeformationtensorsforenergyMorawetz}, 
\beaa
&& \big({}^{(\pr_{\tau})} \pi\big)^{rr}\in r\dk^{\leq 1}\Ga_b, \qquad \big({}^{(\pr_{\tau})} \pi\big)^{r\tau}\in r\dk^{\leq 1}\Ga_g, \qquad |\big({}^{(\pr_{\tau})} \pi\big)^{\tau\tau}|\les \frac{\ep}{r^2},\\  
&& \big({}^{(\pr_{\tau})} \pi\big)^{rx^a}\in \dk^{\leq 1}\Ga_b, \qquad \big({}^{(\pr_{\tau})} \pi\big)^{\tau x^a}\in \dk^{\leq 1}\Ga_g,  \qquad \big({}^{(\pr_{\tau})} \pi\big)^{x^ax^b}\in r^{-1}\dk^{\leq 1}\Ga_g,
\eeaa
we deduce from Lemma \ref{lemma:basiclemmaforcontrolNLterms:ter}
\beaa
\bigg|\int_{\MM_{\varsigma_*(\tau-c_0)\leq r\leq \varsigma_*(\tau-c)}(\tau_1, \tau_2)}{}^{(\pr_\tau)} \pi \cdot \QQ[\dk^{\leq\reg}\psi]\bigg| &\les& \ep\EM^{(\reg)}_{r\geq\varsigma_*(\tau-c_0)}[\psi](\tau_1, \tau_2)
\eeaa
which thus yields
\beaa
&&\E^{(\reg)}_{\varsigma_*(\tau_2-c_0)\leq r\leq \varsigma_*(\tau_2-c)}[\psi](\tau_2)+\F_{\Si_{*,c}}^{(\reg)}[\psi](\tau_1, \tau_2)\nn\\
& \les &\E^{(\reg)}[\psi](\tau_1) +\F_{\Si_{*,c_0}}^{(\reg)}[\psi](\tau_1, \tau_2) +\int_{\MM_{\varsigma_*(\tau-c_0)\leq r\leq \varsigma_*(\tau-c)}(\tau_1, \tau_2)}r^{-1}\big(|\pr\dk^{\leq \reg}\psi|^2+r^{-2}|\dk^{\leq \reg}\psi|^2|\big)\\
&&+ \int_{\MM_{\varsigma_*(\tau-c_0)\leq r\leq \varsigma_*(\tau-c)}(\tau_1, \tau_2)}|\pr_{\tau}(\dk^{\leq\reg}\psi)||\dk^{\leq\reg}F|+\ep\EM^{(\reg)}_{r\geq\varsigma_*(\tau-c_0)}[\psi](\tau_1, \tau_2)
\eeaa
and hence, using in particular the definition of $\F_{\Si_{*,c}}^{(\reg)}[\psi](\tau_1, \tau_2)$, 
\beaa
&&\E^{(\reg)}_{\varsigma_*(\tau_2-c_0)\leq r\leq \varsigma_*(\tau_2-c)}[\psi](\tau_2)+\F_{\Si_{*,c}}^{(\reg)}[\psi](\tau_1, \tau_2)\nn\\
& \les &\E^{(\reg)}[\psi](\tau_1) +\F_{\Si_{*,c_0}}^{(\reg)}[\psi](\tau_1, \tau_2)+\int_{c_0}^c\frac{1}{r}\F_{\Si_{*,c'}}^{(\reg)}[\psi](\tau_1, \tau_2)dc'+ \int_{\MM_{r\geq \varsigma_*(\tau-c_0)}(\tau_1, \tau_2)}r|\dk^{\leq\reg}F|^2\\
&&+\ep\EM^{(\reg)}_{r\geq\varsigma_*(\tau-c_0)}[\psi](\tau_1, \tau_2).
\eeaa
Together with Gronwall, and using the fact that $r\sim c_*$ on $\RR_*(\tau_1, \tau_2)$ in view of \eqref{eq:rsimcuptoOep0inSi*cforcbetweenc*on2andc*}, we deduce
\beaa
&&\E^{(\reg)}_{r\geq\varsigma_*(\tau_2-c_0)}[\psi](\tau_2)+\sup_{c_0\leq c\leq c_*}\F_{\Si_{*,c}}^{(\reg)}[\psi](\tau_1, \tau_2)\nn\\
& \les &\E^{(\reg)}[\psi](\tau_1) +\F_{\Si_{*,c_0}}^{(\reg)}[\psi](\tau_1, \tau_2) + \int_{\MM_{r\geq\varsigma_*(\tau-c_0)}(\tau_1, \tau_2)}r|\dk^{\leq\reg}F|^2+\ep\EM^{(\reg)}_{r\geq\varsigma_*(\tau-c_0)}[\psi](\tau_1, \tau_2).
\eeaa
Then, running again the above argument, this time integrating the current based on the timelike vectorfield $\pr_\tau$ on the spacetime region $\MM_{\varsigma_*(\tau-c_0)\leq r\leq \varsigma_*(\tau-c)}(\tau_1, \tau)$ for $\tau_1\leq\tau\leq \tau_2$, we upgrade the above estimate to 
\beaa
&&\sup_{\tau_1\leq\tau\leq\tau_2}\E^{(\reg)}_{r\geq\varsigma_*(\tau-c_0)}[\psi](\tau)+\sup_{c_0\leq c\leq c_*}\F_{\Si_{*,c}}^{(\reg)}[\psi](\tau_1, \tau_2)\nn\\
& \les &\E^{(\reg)}[\psi](\tau_1) +\F_{\Si_{*,c_0}}^{(\reg)}[\psi](\tau_1, \tau_2) + \int_{\MM_{r\geq\varsigma_*(\tau-c_0)}(\tau_1, \tau_2)}r|\dk^{\leq\reg}F|^2+\ep\EM^{(\reg)}_{r\geq\varsigma_*(\tau-c_0)}[\psi](\tau_1, \tau_2).
\eeaa
Finally, since 
\beaa
\M^{(\reg)}_{0,r\geq\varsigma_*(\tau-c_0)}[\psi](\tau_1, \tau_2)\les\int_{c_0}^{c_*}\frac{1}{r}\F_{\Si_{*,c}}^{(\reg)}[\psi](\tau_1, \tau_2)dc\les \sup_{c_0\leq c\leq c_*}\F_{\Si_{*,c}}^{(\reg)}[\psi](\tau_1, \tau_2),
\eeaa
we deduce, for $\ep>0$ small enough,
\beaa
&&\EM^{(\reg)}_{0,r\geq\varsigma_*(\tau-c_0)}[\psi](\tau_1, \tau_2)+\sup_{c_0\leq c\leq c_*}\F_{\Si_{*,c}}^{(\reg)}[\psi](\tau_1, \tau_2)\\
&\les& \E^{(\reg)}[\psi](\tau_1) +\F_{\Si_{*,c_0}}^{(\reg)}[\psi](\tau_1, \tau_2)+ \int_{\MM_{r\geq\varsigma_*(\tau-c_0)}(\tau_1, \tau_2)}r|\dk^{\leq\reg}F|^2
\eeaa
as stated in \eqref{eq:localenergyestimate:future:RR*}.

Next, we consider \eqref{eq:localenergyestimate:future:unweigthedderivatives:RR*}. First, we commute \eqref{eq:eqsforlocalenergyestimatelemma:general:versionmoregeneralforRR*} by $(\pr_\tau, \pr_r, r^{-1}\pr_{{x^a}})^{\leq\reg}$ and obtain in view of Lemma \ref{lem:commutatorwithwave:firstorderderis:0}, for $\reg\leq 14$, 
\beaa
\square_{\g}\pr^{\leq\reg}\psi_{ij}= O(r^{-1})\pr\pr^{\leq\reg}\psi_{kl} +O(r^{-2})\pr^{\leq \reg}\psi_{kl} +\pr^{\leq\reg}F_{ij}+\pr^{\leq\reg}(\gcheck^{\a\b})\pr_{\a}\pr_{\b}\pr^{\leq\reg -1}\psi,
\eeaa
and then, integrating the current based on the timelike vectorfield $\pr_\tau$ on the spacetime region $\MM_{\varsigma_*(\tau-c_0)\leq r\leq \varsigma_*(\tau-c)}(\tau_1, \tau_2)$ for $c_0\leq c\leq c_*$, we obtain, arguing as above, 
\beaa
&&\E_{\varsigma_*(\tau_2-c_0)\leq r\leq \varsigma_*(\tau_2-c)}[\pr^{\leq\reg}\psi](\tau_2)+\F_{\Si_{*,c}}[\pr^{\leq\reg}\psi](\tau_1, \tau_2)\nn\\
& \les &\E[\pr^{\leq\reg}\psi](\tau_1) +\F_{\Si_{*,c_0}}[\pr^{\leq\reg}\psi](\tau_1, \tau_2) +\int_{c_0}^c\frac{1}{r}\F_{\Si_{*,c'}}[\pr^{\leq\reg}\psi](\tau_1, \tau_2)dc'\\
&&+ \int_{\MM_{r\geq\varsigma_*(\tau-c_0)}(\tau_1, \tau_2)}r|\pr^{\leq\reg}F|^2+\ep\EM_{r\geq\varsigma_*(\tau-c_0)}[\pr^{\leq\reg}\psi](\tau_1, \tau_2)\\
&&+\left|\int_{\MM_{\varsigma_*(\tau-c_0)\leq r\leq \varsigma_*(\tau-c)}(\tau_1, \tau_2)}\pr^{\leq\reg}(\gcheck^{\a\b})\pr_{\a}\pr_{\b}\pr^{\leq\reg -1}\psi\ov{\pr_\tau\pr^{\leq\reg}\psi}\right|.
\eeaa
Estimating the last term on the RHS of the above inequality using Lemma \ref{lemma:basiclemmaforcontrolNLterms:bis}, we infer, using also Gronwall, and the fact that $r\sim c_*$ on $\RR_*(\tau_1, \tau_2)$,
\beaa
&&\E_{r\geq\varsigma_*(\tau_2-c_0)}[\pr^{\leq\reg}\psi](\tau_2)+\sup_{c_0\leq c\leq c_*}\F_{\Si_{*,c}}[\pr^{\leq\reg}\psi](\tau_1, \tau_2)\nn\\
 &\les &\E[\pr^{\leq\reg}\psi](\tau_1) +\F_{\Si_{*,c_0}}[\pr^{\leq\reg}\psi](\tau_1, \tau_2) + \int_{\MM_{r\geq\varsigma_*(\tau-c_0)}(\tau_1, \tau_2)}r|\pr^{\leq\reg}F|^2+\ep\EM_{r\geq\varsigma_*(\tau-c_0)}[\pr^{\leq\reg}\psi](\tau_1, \tau_2).
\eeaa
We then proceed as above to recover the control of $\EM_{0,r\geq\varsigma_*(\tau-c_0)}[\pr^{\leq\reg}\psi](\tau_1, \tau_2)$ and to absorb the last term on the RHS for $\ep>0$ which then yields the stated estimate  \eqref{eq:localenergyestimate:future:unweigthedderivatives:RR*}. This concludes the proof of Lemma \ref{lemma:localenergyestimate:weightedandunweigthedderivatives:RR*}.
\end{proof}

%%%%%%%%%%%%%%%%%%%%%%%%%%%%%%%%%%%%%%%

\subsection{Statement and proof the main result of Section \ref{sec:energyMorawetzesitmatesforTeukoslkyonMM:upto15derivatives}}
\lab{sect:maintheoremandproof}

%%%%%%%%%%%%%%%%%%%%%%%%%%%%%%%%%%%%%%%

 In this section, we state and prove Energy-Morawetz estimates for up to 14 derivatives for solutions to Teukolsky equations,  which is the main result of Section \ref{sec:energyMorawetzesitmatesforTeukoslkyonMM:upto15derivatives}.

%%%%%%%%%%%%%%%%%%%%%%%%%%%%%%%%%%%%%%%%%%%%%%%%%%%%%%%%%%%%%%%%

\subsubsection{Statement of the main result of Section \ref{sec:energyMorawetzesitmatesforTeukoslkyonMM:upto15derivatives}}

%%%%%%%%%%%%%%%%%%%%%%%%%%%%%%%%%%%%%%%%%%%%%%%%%%%%%%%%%%%%%%%%

We now state Energy-Morawetz estimates for up to 14 derivatives for solutions to Teukolsky equations.

\begin{theorem}
\label{thm:main:MaSz26}
Let $(\MM, \g)$ satisfy the assumptions of Section \ref{sec:assumptionsforsec:energyMorawetzesitmatesforTeukoslkyonMM:upto15derivatives}. Then, for $\ep>0$ small enough, we have for solutions $\pmb\phi_s^{(p)}$, $s=\pm 2$, $p=0,1,2$, to the tensorial Teukolsky wave/transport systems \eqref{eq:TensorialTeuSysandlinearterms:rescaleRHScontaine2:general:Kerrperturbation}  \eqref{def:TensorialTeuScalars:wavesystem:Kerrperturbation}  on spacetimes $(\MM,\g)$ perturbations of Kerr the following energy-Morawetz-flux estimates, for any $1\leq\tau_1<\tau_2\leq\tau_*$ and any $0<\de\leq \frac{1}{3}$, 
\bea
\label{MainEnerMora:psi:plus2case}
\sum_{p=0}^2\EMF^{({11})}_{\de}[\pmb\phi_{+2}^{(p)}](\tau_1, \tau_2)
&\les&  \sum_{p=0}^2\E^{({11})}[\pmb\phi_{+2}^{(p)}](\tau_1)
+\sum_{p=0}^2\NN^{({11})}_\de[\pmb\phi_{+2}^{(p)}, \N_{W,+2}^{(p)}](\tt_1, \tt_2)\nn\\ 
&&
+\sum_{p=0}^1\int_{\MM(\tt_1, \tt_2)}r^{-1+\de}|\dk^{\leq 12}\N_{T,+2}^{(p)}|^2
\eea
and, assuming also that $\Ab$ satisfies \eqref{eq:waveequationpmbphip=0sminus2nodeginredshiftregion},
\bea
\label{MainEnerMora:psi:minus2case}
\nn&&\sum_{p=0}^2\EMF^{({14})}_\de[\pmb\phi_{-2}^{(p)}](\tau_1, \tau_2) + \sum_{p=0}^2\EMF^{({14})}_{r\leq r_+(1+\dred)}[\nab_4^p\Ab](\tau_1, \tau_2)\\
\nn&\les& \sum_{p=0}^2\E^{({14})}[\pmb\phi_{-2}^{(p)}](\tau_1)+\sum_{p=0}^2\E^{({14})}_{r\leq r_+(1+\dred)}[\nab_4^p\Ab](\tau_1) +\int_{\Si(\tau_1)}|\dk^{\leq 12}\N_{T,-2}^{(0)}|^2\\
\nn&&+\sum_{p=0}^1\int_{\MM(\tt_1, \tt_2)}r^{-1+\de}|\dk^{\leq {15}}\N_{T,-2}^{(p)}|^2+\sum_{p=0}^2\NN_\de^{({14})}[\pmb\phi_{-2}^{(p)}, \N_{W,-2}^{(p)}](\tt_1, \tt_2)\\ 
&&+\sum_{p=0}^2\int_{\MM_{r\leq r_+(1+2\dred)}(\tau_1, \tau_2)}|\dk^{\leq {14}}\N_{\nab_4^p\Ab}|^2,
\eea
where the norms $\EMF^{(\reg)}_\de[\c](\tau_1,\tau_2)$, $\E^{(\reg)}[\c](\tau_1)$ and $\mathcal{N}^{(\reg)}_\de[\c, \c](\tau_1, \tau_2)$ have been introduced in Section \ref{subsection:basicnormsforpsi:MMh}, where $\N_{W,s}^{(p)}$, $\N_{T,s}^{(p)}$ and $\N_{\nab_4^p\Ab}$ are introduced in equations \eqref{eq:TensorialTeuSysandlinearterms:rescaleRHScontaine2:general:Kerrperturbation},   \eqref{def:TensorialTeuScalars:wavesystem:Kerrperturbation} and \eqref{eq:waveequationpmbphip=0sminus2nodeginredshiftregion}, and where the implicit constant in $\lesssim$ only depends on $a$, $m$, $\de$, $\Nmic$ and $\dec$.
\end{theorem}

Theorem \ref{thm:main:MaSz26} is the analog of Theorem 7.1 in \cite{MaSz26}. In order to prove Theorem \ref{thm:main:MaSz26}, we will need to slightly adapt the proof of Theorem 7.1 in \cite{MaSz26} given that the spacetime in \cite{MaSz26} extends to $\II_+$ while $\MM$ only extends to the spacelike hypersuface $\Si_*$. We start by extending the solutions to Teukolsky on $\MM$ to a semi-global problem.

%%%%%%%%%%%%%%%%%%%%%%%%%%%%%%%%

\subsubsection{Extension to a semi-global problem}
\lab{sect:extensiontoglobalproblem:Teu}

%%%%%%%%%%%%%%%%%%%%%%%%%%%%%%%%

We start with the following definition.
\begin{definition}\lab{def:definitionofthetimetauR}
Let $\Nmic$ be the large integer introduced in  Section \ref{sec:smallnesconstants}. For $\tau_1\geq 1$, we define $\tmic<\tau_1$ as the smallest real number such that\footnote{Note that $\{\tau=\tmic\}$ is possibly included in the extension $(\MMh, \gh)$ of $(\MM, \g)$ since $\tau\geq 1$ on $\MM$.}
\bea
\Sih(\tmic)\cap\{r\leq (\Nmic+1)m\}\subset D^-\big(\Si(\tau_1)\cap\{r\leq 2\Nmic m\}\big),
\eea
i.e., $\Sih(\tmic)\cap\{r\leq (\Nmic+1)m\}$ is included in the past domain of dependence of $\Si(\tau_1)\cap\{r\leq 2\Nmic m\}$. We also define the interval $\Iti:=(\tmic, +\infty)$.
\end{definition}

\begin{remark}
The choice of $\tmic$ in Definition \ref{def:definitionofthetimetauR} satisfies  
\bea
-\frac{m}{\Nmic}\les \tmic -\tau_1\les -\frac{m}{\Nmic},
\eea
see Remark 7.4 in \cite{MaSz26}.
\end{remark}

Recall from Lemma \ref{lem:scalarizedTeukolskywavetransportsysteminKerrperturbation:Omi} that the tensorial Teukolsky wave/transport systems \eqref{eq:TensorialTeuSysandlinearterms:rescaleRHScontaine2:general:Kerrperturbation}  \eqref{def:TensorialTeuScalars:wavesystem:Kerrperturbation} are equivalent to the scalarized Teukolsky wave/transport systems  \eqref{eq:ScalarizedTeuSys:general:Kerrperturbation} \eqref{eq:ScalarizedQuantitiesinTeuSystem:Kerrperturbation}. The proof of Theorem \ref{thm:main:MaSz26} will follow from global energy-Morawetz estimates for an extension to $\tau\in\Iti$, with $\Iti$ as in Definition \ref{def:definitionofthetimetauR}, of the scalarized system of wave equations \eqref{eq:ScalarizedTeuSys:general:Kerrperturbation}. The goal of this section is to construct this extended solution. The following  proposition is an adaptation of Proposition 7.5 in \cite{MaSz26}.

\begin{proposition}\lab{prop:extensionprocedureoftheTeukolskywaveequations}
Let $(\MM, \g)$ satisfy the assumptions of Section \ref{sec:assumptionsforsec:energyMorawetzesitmatesforTeukoslkyonMM:upto15derivatives}. Assume that $\phi^{(p)}_{s,ij}$ satisfies the scalarized system of wave equations \eqref{eq:ScalarizedTeuSys:general:Kerrperturbation} for $\tau\in(\tau_1, \tau_2)$ with $\tau_1\geq 1$ and $\tau_2$ satisfying\footnote{Since the proof of Theorem \ref{thm:main:MaSz26} in the case $\tau_1<\tau_2\leq\tau_1+10$ follows immediately from local energy type arguments, see Step 0 in Section \ref{subsect:proofofThm4.1}, we focus here on the case where $\tt_2$ satisfies  \eqref{eq:tau2geqtau1plus10conditionisthemaincase}.} 
\bea\lab{eq:tau2geqtau1plus10conditionisthemaincase}
\tau_1+10\leq \tau_2\leq\tau_*.
\eea
Also, let $\tmic<\tau_1$ and $\Iti$ be as in Definition \ref{def:definitionofthetimetauR}, let $\chi_{\tau_1, \tau_2}=\chi_{\tau_1, \tau_2}(\tau)$ be a smooth cut-off function satisfying 
\bea\lab{eq:propertieschitoextendmetricg}
\chi_{\tau_1, \tau_2}(\tau)=0\,\,\,\textrm{on}\,\,\, \mathbb{R}\setminus(\tau_1, \tau_2-1), \,\,\,\,\chi_{\tau_1, \tau_2}(\tau)=1\,\,\,\textrm{on}\,\,\, [\tau_1+1, \tau_2-2], \,\,\,\, \|\chi_{\tau_1, \tau_2}\|_{W^{{15},+\infty}(\mathbb{R})}\les 1,
\eea
let $\chi^{(1)}_{\tau_1, \tau_2}=\chi^{(1)}_{\tau_1, \tau_2}(\tau)$ be a smooth cut-off function satisfying 
\bea\lab{eq:propertieschi:thisonetodealwithRHSofTeukolsky}
\chi^{(1)}_{\tau_1, \tau_2}(\tau)=0\,\,\,\textrm{on}\,\,\, \mathbb{R}\setminus(\tau_1, \tau_2-2), \,\,\,\,\chi^{(1)}_{\tau_1, \tau_2}(\tau)=1\,\,\,\textrm{on}\,\,\, [\tau_1+1, \tau_2-3], \,\,\,\, \|\chi^{(1)}_{\tau_1, \tau_2}\|_{W^{{15},+\infty}(\mathbb{R})}\les 1,
\eea
let $\chi_{r_*}(\tau, r)$ be a smooth cut-off introduced in \eqref{eq:supportpropertiesofcutoffchistar}, and let $\chi^{(1)}_{r_*}=\chi^{(1)}_{r_*}(\tau, r)$ be a smooth cut-off function satisfying 
\bea\lab{eq:supportpropertiesofcutoffchistar:(1)}
\bsplit
\chi^{(1)}_{r_*}&=1\quad\textrm{for}\quad r\leq\varsigma_*\left(\tau-\frac{5c_*}{8}\right),\\ 
\chi^{(1)}_{r_*}&=0\quad\textrm{for}\quad r\geq\varsigma_*\left(\tau-\frac{3c_*}{4}\right), \qquad \|(r\pr_\tau, r\pr_r)^{\leq 15}\chi^{(1)}_{r_*}\|_{L^\infty}\lesssim 1,
\end{split}
\eea
and define the extended metric on $\MMh$
\bea\lab{eq:extendedmetricgchitau1tau2}
\g_{\tau_1, \tau_2,*}^{\a\b}=\chi_{\tau_1, \tau_2}\chi_{r_*}\g^{\a\b}+(1-\chi_{\tau_1, \tau_2}\chi_{r_*})\g_{a,m}^{\a\b}.
\eea
Then, there exists $\psi^{(p)}_{s,ij}$ satisfying the following system of scalar wave equations defined by
\bea\lab{eq:waveeqwidetildepsi1}
\big({\square}_{\g_{\tau_1, \tau_2,*}}- (4-2\de_{p0})|q|^{-2}\big)\psi^{(p)}_{s,ij} &=F_{total,s,ij}^{(p)} \quad\textrm{on}\quad\MMh(\Iti),
\eea
with\footnote{Below, $\widehat{S}_K$ and $\widehat{Q}_K$ denote $\widehat{S}$ and $\widehat{Q}$ of \eqref{eq:definitionwidehatSandwidehatQperturbationsofKerr} computed using the regular triplet in Kerr of Definition \ref{def:regulartripletinKerrOmii=123}.}
\bea\lab{def:tildef}
\bsplit
F_{total,s,ij}^{(p)}:=&\widehat{F}^{(p)}_{s,ij} +\widetilde{F}^{(p)}_{s,ij}+\underline{F}^{(p)}_{s,ij}+\breve{F}^{(p)}_{s,ij},\\
\widehat{F}^{(p)}_{s,ij}:=&\chi_{\tau_1, \tau_2}\chi_{r_*}\big( \widehat{S}(\pmb\psi^{(p)}_s)_{ij} +(\widehat{Q}\pmb\psi^{(p)}_s)_{ij}\big)\\
& +(1-\chi_{\tau_1, \tau_2}\chi_{r_*})\big(\widehat{S}_K(\pmb\psi^{(p)}_s)_{ij} +(\widehat{Q}_K\pmb\psi^{(p)}_s)_{ij} +f_p\psi^{(p)}_{s,ij}\big),\\
f_p:=& f_p(r, \cos\th)=\frac{2\de_{p0}}{|q|^2}- \frac{4a^2\cos^2\th(|q|^2+6mr)}{|q|^6},\\
\widetilde{F}^{(p)}_{s,ij}:=&\chi^{(1)}_{\tau_1, \tau_2}\chi^{(1)}_{r_*}F^{(p)}_{s,ij}, \qquad F^{(p)}_{s,ij}=L^{(p)}_{s,ij}+N^{(p)}_{W,s,ij},  \quad i,j=1,2,3,\,\,\, p=0,1,2, \,\,\, s=\pm 2,
\end{split}
\eea
${L_{s,ij}^{(p)}}=(\L_{s}^{(p)}[\pmb\phi_s])_{ij}$ given by \eqref{eq:linearterms:ScalarizedTeuSys:general:Kerrperturbation}, 
and $\underline{F}^{(p)}_{s,ij}$ and $\breve{F}^{(p)}_{s,ij}$ given respectively by \eqref{def:tildef0} and \eqref{def:breveFpsij},
such that the following properties hold:
\begin{itemize}
\item $\g_{\tau_1, \tau_2,*}$ satisfies \eqref{eq:controloflinearizedinversemetriccoefficients} \eqref{eq:controloflinearizedinversemetriccoefficients:inversegtautau} and coincides with Kerr in $\MMh\setminus(\MM_{r\leq\varsigma_*(\tau-\frac{7c_*}{8})}(\tau_1, \tau_2-1))$, 

\item $\psi^{(p)}_{s,ij}$ satisfies the following identities 
\bea\lab{eq:causlityrelationsforwidetildepsi1}
\bsplit
\psi^{(p)}_{s,ij}&=\phi^{(p)}_{s,ij}\quad\textrm{on}\quad\MM_{r\leq\varsigma_*(\tau-\frac{5c_*}{8})}(\tau_1+1, \tau_2-3), \\
\dk^{\leq\reg}(\psi^{(p)}_{s,ij}) &= 0\quad\textrm{on}\quad \Sih(\tmic)\cap 
\{r\leq(\Nmic+1)m\}\,\,\,\,\forall\,\reg\in\mathbb{N},
\end{split}
\eea
so that $\psi^{(p)}_{s,ij}$ can be smoothly extended\footnote{This property is used to derive microlocal energy-Morawetz estimates on $\MM_{r\leq\Rmic}$ in section 8 of \cite{MaSz26}, where $\Rmic\in [\Nmic m, (\Nmic+1)m]$ is introduced in Remark \ref{rmk:choiceofconstantRbymeanvalue}.} by $0$ in $\MMh(-\infty, \tmic)\cap\{r\leq (\Nmic+1)m\}$.

\item $\psi^{(p)}_{s,ij}$ satisfies the following local energy estimate on $\tau\in[\tmic, \tau_1+3]$, for all $\reg\leq {14}$, 
\bea\lab{eq:localenergyforpsiontau1minus1tau1plus1}
&&\sum_{p=0}^2\sum_{p=0}^2\EMFh_\de[\pr^{\leq \reg}\psi^{(p)}_{s}](\tmic, \tau_1+3)\nn\\
&\les& \sum_{i,j=1}^3\sum_{p=0}^2\bigg(\E[\pr^{\leq \reg}\phi^{(p)}_{s,ij}](\tau_1)+\int_{\MM(\tau_1, \tau_1+3)}r^{1+\de} |\pr^{\leq \reg}N_{W,s,ij}^{(p)}|^2\bigg).
\eea

\item $\psi^{(p)}_{s,ij}$ satisfies the following local energy estimates on $\RR_*(\tau_1, \tau_2-3)\cup\MM(\tau_2-3, \tau_2)$, for all $\reg\leq {14}$, 
\bea\lab{eq:localenergyestimateforpsisijponSigmatau2minus3tau2:partialderivatives}
\nn&&\sum_{i,j=1}^3\sum_{p=0}^2\bigg(\EM_{0,r\geq\varsigma_*(\tau-c_*/2)}[\pr^{\leq \reg}\psi^{(p)}_{s,ij}](\tau_1, \tau_2-1)\\
\nn&&+\sup_{\frac{c_*}{2}\leq c\leq c_*}\F_{\Si_{*,c}}[\pr^{\leq\reg}\psi^{(p)}_{s,ij}](\tau_1, \tau_2-1)+\EMF_\de[\pr^{\leq \reg}\psi^{(p)}_{s,ij}](\tau_2-3, \tau_2)\bigg)\\ 
\nn&\les& \sum_{i,j=1}^3\sum_{p=0}^2\bigg(\E[\pr^{\leq \reg}\phi^{(p)}_{s,ij}](\tau_1)+\E[\pr^{\leq \reg}\phi^{(p)}_{s,ij}](\tau_2-3)+\F_{\Si_{*,c_*/2}}[\pr^{\leq\reg}\phi^{(p)}_{s,ij}](\tau_1+1, \tau_2-1)\\
&&+\int_{\MM(\tau_1, \tau_1+1)\cup\RR_*(\tau_1+1, \tau_2-3)\cup\MM(\tau_2-3, \tau_2-2)}r^{1+\de}|\pr^{\leq \reg}N_{W,s,ij}^{(p)}|^2\bigg).
\eea

\item The tensorization defect $\err_{\textrm{TDefect}}[\psi^{(p)}_{s}]$ corresponding to the family of complex-valued scalars $\psi_{ij}$ as introduced in  Definition \ref{def:definitionofthenotationerrforthescalarizationdefect} satisfies
\begin{align}\lab{eq:propertyofscalarizationdefect:prop1}
\nn&\err_{\textrm{TDefect}}[\psi^{(p)}_{s}]=0\quad\textrm{on}\\
&\qquad\,\,\MMh(\tau_1+1, +\infty)\setminus(\MM(\tau_2-2, \tau_2)\cup\MM_{\varsigma_*(\tau-\frac{5c_*}{8})\leq r\leq\varsigma_*(\tau-\frac{31c_*}{32})}(\tau_1+1, \tau_2-2)),
\end{align}
and, for all $\reg\leq {14}$,
\begin{align}\lab{eq:propertyofscalarizationdefect:prop2}
\nn&\sum_{p=0}^2\Big(\EM_{0,r\geq\varsigma_*(\tau-c_*/2)}[\pr^{\leq \reg}\err_{\textrm{TDefect}}[\psi^{(p)}_{s}]](\tau_1+1, \tau_2-2)\\
\nn&+\sup_{c_*/2\leq c\leq c_*}\F_{\Si_{*,c}}[\pr^{\leq \reg}\err_{\textrm{TDefect}}[\psi^{(p)}_{s}]](\tau_1+1, \tau_2-2)\\
\nn&+\EMF_\de[\pr^{\leq \reg}\err_{\textrm{TDefect}}[\psi^{(p)}_{s}]](\tau_2-2, \tau_2)\Big)\\ 
\nn\les& \ep^2\sum_{i,j=1}^3\sum_{p=0}^2\bigg(\E[\pr^{\leq \reg}\phi^{(p)}_{s,ij}](\tau_1)+\E[\pr^{\leq \reg}\phi^{(p)}_{s,ij}](\tau_2-3)+\F_{\Si_{*,c_*/2}}[\pr^{\leq\reg}\phi^{(p)}_{s,ij}](\tau_1+1, \tau_2-1)\\
&+\int_{\MM(\tau_1, \tau_1+1)\cup\RR_*(\tau_1+1, \tau_2-3)\cup\MM(\tau_2-3, \tau_2-2)}r^{1+\de}|\pr^{\leq \reg}N_{W,s,ij}^{(p)}|^2\bigg).
\end{align}
\end{itemize}
\end{proposition}

\begin{proof}
We adapt the proof of Proposition 7.5 in \cite{MaSz26}.

We proceed along the following steps. 

\noindent{\bf Step 1.} Since $\g_{\tau_1, \tau_2,*}$ is defined by \eqref{eq:extendedmetricgchitau1tau2}, we have
\beaa
\widecheck{\g}_{\tau_1, \tau_2,*}^{\a\b}=\chi_{\tau_1, \tau_2}\chi_{r_*}\g^{\a\b}+(1-\chi_{\tau_1, \tau_2}\chi_{r_*})\gam^{\a\b} -\gam^{\a\b}=\chi_{\tau_1, \tau_2}\chi_{r_*}\left(\g^{\a\b}-\gam^{\a\b}\right)=\chi_{\tau_1, \tau_2}\chi_{r_*}\widecheck{\g}^{\a\b}
\eeaa
and hence
\beaa
|\dk^{\leq {15}}\widecheck{\g}_{\tau_1, \tau_2,*}^{\a\b}|\les |\dk^{\leq {15}}(\chi_{\tau_1, \tau_2}\chi_{r*})||\dk^{\leq {15}}\widecheck{\g}^{\a\b}|\les |\pr_\tau^{\leq 15}\chi_{\tau_1, \tau_2}||(\pr_\tau, r_*\pr_r)^{\leq 15}\chi_{r_*}||\dk^{\leq {15}}\widecheck{\g}^{\a\b}| \les |\dk^{\leq {15}}\widecheck{\g}^{\a\b}|,
\eeaa
where we used the fact that $|\dk^{\leq {15}}(\chi_{\tau_1, \tau_2}\chi_{r_*})|\les |(\pr_\tau, r\pr_r)^{\leq {15}}(\chi_{\tau_1, \tau_2}\chi_{r_*})|$. Since $\g$ satisfies \eqref{eq:controloflinearizedinversemetriccoefficients} \eqref{eq:controloflinearizedinversemetriccoefficients:inversegtautau}, and in view of the properties \eqref{eq:propertieschitoextendmetricg} \eqref{eq:supportpropertiesofcutoffchistar} respectively of $\chi_{\tau_1,\tau_2}$ and $\chi_{r_*}$, we deduce that $\g_{\tau_1, \tau_2,*}$ also satisfies \eqref{eq:controloflinearizedinversemetriccoefficients} \eqref{eq:controloflinearizedinversemetriccoefficients:inversegtautau}, and in addition coincides with $\gam$ in the spacetime region $\MMh\setminus(\MM_{r\leq\varsigma_*(\tau-\frac{7c_*}{8})}(\tau_1, \tau_2-1))$. 

\noindent{\bf Step 2.}  Next, we introduce the solutions $\widetilde{\phi}^{(p)}_{s,ij}$ to the following auxiliary system of scalar wave equations, for $p=0,1,2$, $s=\pm 2$, and $i,j=1,2,3$,
\bea\lab{eq:waveequationdefiningfirstextensionwidetildepsi}
\bsplit
\big({\square}_{\g_{\tau_1, \tau_2,*}}-(4-2\de_{p0})|q|^{-2}\big)\widetilde{\phi}^{(p)}_{s,ij} =& \chi_{\tau_1, \tau_2}\chi_{r_*}\big( \widehat{S}(\widetilde{\pmb\phi}^{(p)}_s)_{ij} +(\widehat{Q}\widetilde{\pmb\phi}^{(p)}_s)_{ij}\big)\\
& +(1-\chi_{\tau_1, \tau_2}\chi_{r_*})\big( \widehat{S}_K(\widetilde{\pmb\phi}^{(p)}_s)_{ij} +(\widehat{Q}_K\widetilde{\pmb\phi}^{(p)}_s)_{ij}+f_p\widetilde{\phi}^{(p)}_{s,ij}\big)\\
&+\widetilde{F}^{(p)}_{s,ij}\,\,\,\,\textrm{on}\,\,\,\MMh(\tau_1, \tau_2), \\
 \widetilde{\phi}^{(p)}_{s,ij} =& \chi_{r_*}\phi^{(p)}_{s,ij}, \,\,\,\, N_{\Sih(\tau_1+1)}\widetilde{\phi}^{(p)}_{s,ij}=\chi_{r_*}N_{\Sigma(\tau_1+1)}\phi^{(p)}_{s,ij}\,\,\,\textrm{on}\,\,\,\Sih(\tau_1+1),\\
\widetilde{\phi}^{(p)}_{s,ij} =& \phi^{(p)}_{s,ij}\quad\textrm{on}\quad\AA(\tau_1, \tau_1+1), \quad\widetilde{\phi}^{(p)}_{s,ij} =0\quad\textrm{on}\quad\II_+(\tau_1, \tau_1+1).
\end{split}
\eea 
Then, in view of \eqref{eq:waveequationdefiningfirstextensionwidetildepsi} and the support properties of $\chi_{r_*}$ in \eqref{eq:supportpropertiesofcutoffchistar}, we have by causality
\bea\lab{eq:bycausalityanddefintionwidetildephipsij=0inMMtau1tau1plus1}
\widetilde{\phi}^{(p)}_{s,ij}=0\quad\textrm{on}\quad(\MMh\setminus\MM_{r\leq\varsigma_*(\tau-\frac{7c_*}{8})})(\tau_1, \tau_1+1),
\eea
and we have in view of the local energy estimate\footnote{In view of \eqref{eq:bycausalityanddefintionwidetildephipsij=0inMMtau1tau1plus1}, it suffices to apply \eqref{eq:localenergyestimate:past:bis:unweigthedderivatives} on $\MM(\tau_1, \tau_1+1)$ with trivial initial data on $\Si_*(\tau_1, \tau_1+1)$.} \eqref{eq:localenergyestimate:past:bis:unweigthedderivatives}  with $\reg\leq {14}$ for \eqref{eq:waveequationdefiningfirstextensionwidetildepsi}
\beaa
&&\sum_{i,j=1}^3\sum_{p=0}^2\EMFh_\de[\pr^{\leq \reg}\widetilde{\phi}^{(p)}_{s,ij}](\tau_1, \tau_1+1)\nn\\
&\les& \sum_{i,j=1}^3\sum_{p=0}^2\bigg(\EF[\pr^{\leq \reg}\phi^{(p)}_{s,ij}](\tau_1, \tau_1+1)+\int_{\MM(\tau_1, \tau_1+1)}r^{1+\de} |\pr^{\leq \reg}N_{W,s,ij}^{(p)}|^2\bigg).
\eeaa
Together with the local energy estimate \eqref{eq:localenergyestimate:future:unweigthedderivatives:bis} with $\reg\leq {14}$ for $\phi^{(p)}_{s,ij}$, we infer
\bea\lab{eq:localenergyestimateforphisijponSigmatau1tau1plus1}
\nn&&\sum_{i,j=1}^3\sum_{p=0}^2\EMF_\de[\pr^{\leq \reg}\phi^{(p)}_{s,ij}](\tau_1, \tau_1+1)\\ 
&\les& \sum_{i,j=1}^3\sum_{p=0}^2\bigg(\E[\pr^{\leq \reg}\phi^{(p)}_{s,ij}](\tau_1)+\int_{\MM(\tau_1, \tau_1+1)}r^{1+\de} |\pr^{\leq \reg}N_{W,s,ij}^{(p)}|^2\bigg)
\eea
and
\bea\lab{eq:localenergyestimateforwidetildepsionSigma}
\nn&&\sum_{i,j=1}^3\sum_{p=0}^2\EMFh_\de[\pr^{\leq \reg}\widetilde{\phi}^{(p)}_{s,ij}](\tau_1, \tau_1+1)\\ 
&\les& \sum_{i,j=1}^3\sum_{p=0}^2\bigg(\E[\pr^{\leq \reg}\phi^{(p)}_{s,ij}](\tau_1)+\int_{\MM(\tau_1, \tau_1+1)}r^{1+\de} |\pr^{\leq \reg}N_{W,s,ij}^{(p)}|^2\bigg).
\eea

Also, let $\chi_{\Nmic}(r)$ be a smooth cut-off function such that $\chi_{\Nmic}(r)=1$ for $r\leq 2\Nmic m$ and $\chi_{\Nmic}(r)=0$ for $r\geq 4\Nmic m$. Then, we introduce the solution $(\phi_{aux})^{(p)}_{s,ij}$ to the following auxiliary system of wave equations
\bea\lab{eq:waveequationdefiningpsiaux}
\bsplit
\big({\square}_{\g_{\tau_1, \tau_2,*}}-(4-2\de_{p0})|q|^{-2}\big)(\phi_{aux})^{(p)}_{s,ij} =&\widehat{S}_K((\pmb\phi_{aux})^{(p)}_s)_{ij} +(\widehat{Q}_K(\pmb\phi_{aux})^{(p)}_s)_{ij}\\
& +f_p(\phi_{aux})^{(p)}_{s,ij}\,\,\,\textrm{on}\,\,\,\MMh(\tmic, \tau_1),\\
(\phi_{aux})^{(p)}_{ij}=&\chi_{\Nmic}(r)\widetilde{\phi}^{(p)}_{s,ij}\quad\textrm{on}\quad\Sih(\tau_1),\\
N_{\Sigma(\tau_1)}(\phi_{aux})^{(p)}_{s,ij}=&\chi_{\Nmic}(r)N_{\Sigma(\tau_1)}(\widetilde{\phi}^{(p)}_{s,ij})\quad\textrm{on}\quad\Sih(\tau_1),\\
(\phi_{aux})^{(p)}_{s,ij}=& (\phi_\AA)^{(p)}_{s,ij}\quad\textrm{on}\quad\AA_+\cap\{\tmic\leq\tau\leq\tau_1\},\\ 
(\phi_{aux})^{(p)}_{s,ij}=&0\quad\textrm{on}\quad\II_+\cap\{\tmic\leq\tau\leq\tau_1\},
\end{split}
\eea 
where $(\phi_\AA)^{(p)}_{s,ij}$ is a smooth extension of $\phi^{(p)}_{s,ij}$ from $\AA\cap\{\tau\geq \tau_1\}$ to $\AA\cap\{\tmic\leq\tau\leq\tau_1\}$ satisfying, for $\reg\leq {14}$,
\bea\lab{eq:propertyextensionpsiAAofpsitotauleqtau1}
\F_\AA[\pr^{\leq\reg}(\phi_\AA)^{(p)}_{s,ij}](\tmic, \tau_1) \les \F_\AA[\pr^{\leq\reg}\phi^{(p)}_{s,ij}](\tau_1, \tau_1+1).
\eea 
In view of \eqref{eq:waveequationdefiningpsiaux}, the support properties of $\chi_{\Nmic}(r)$ and the fact that  $4\Nmic m\ll \ep^{-1}<r_*$, we have by causality
\bea\lab{eq:bycausalityanddefintionphiauxpsij=0inMMhsetminusMMtmictau1}
(\phi_{aux})^{(p)}_{s,ij}=0\quad\textrm{on}\quad(\MMh\setminus\MM)(\tmic, \tau_1),
\eea
and the local energy estimate\footnote{In view of \eqref{eq:bycausalityanddefintionphiauxpsij=0inMMhsetminusMMtmictau1}, it suffices to apply \eqref{eq:localenergyestimate:past:bis:unweigthedderivatives} on $\MM(\tmic, \tau_1)$ with trivial initial data on $\Si_*(\tmic, \tau_1)$.}  \eqref{eq:localenergyestimate:past:bis:unweigthedderivatives} for  \eqref{eq:waveequationdefiningpsiaux} yields, for $\reg\leq {14}$, 
\beaa
&&\sum_{i,j=1}^3\sum_{p=0}^2\EMFh_\de[\pr^{\leq \reg}(\phi_{aux})^{(p)}_{s,ij}](\tmic, \tau_1)\\ 
&\les& \sum_{i,j=1}^3\sum_{p=0}^2\Big(\E[\pr^{\leq \reg}(\phi_{aux})^{(p)}_{s,ij}](\tau_1)+\F[\pr^{\leq \reg}(\phi_{aux})^{(p)}_{s,ij}](\tmic, \tau_1)\Big)\\
&\les& \sum_{i,j=1}^3\sum_{p=0}^2\Big(\E[\pr^{\leq \reg}\widetilde{\phi}^{(p)}_{s,ij}](\tau_1)+\F_\AA[\pr^{\leq \reg}(\phi_\AA)^{(p)}_{s,ij}](\tmic, \tau_1)\Big)
\eeaa
which together with  \eqref{eq:localenergyestimateforwidetildepsionSigma} and \eqref{eq:propertyextensionpsiAAofpsitotauleqtau1} implies
\beaa
&&\sum_{i,j=1}^3\sum_{p=0}^2\EMFh_\de[\pr^{\leq \reg}(\phi_{aux})^{(p)}_{s,ij}](\tmic, \tau_1)\\
&\les& \sum_{i,j=1}^3\sum_{p=0}^2\bigg(\EF[\pr^{\leq \reg}\phi^{(p)}_{s,ij}](\tau_1, \tau_1+1)+\int_{\MM(\tau_1, \tau_1+1)}r^{1+\de} |\pr^{\leq \reg}N_{W,s,ij}^{(p)}|^2\bigg).
\eeaa
Using \eqref{eq:localenergyestimateforphisijponSigmatau1tau1plus1}, we deduce, for $\reg\leq {14}$,
\bea\lab{eq:localenergyestimateforpsiauxonSigma}
\nn&&\sum_{i,j=1}^3\sum_{p=0}^2\EMFh_\de[\pr^{\leq \reg}(\phi_{aux})^{(p)}_{s,ij}](\tmic, \tau_1)\\ 
&\les& \sum_{i,j=1}^3\sum_{p=0}^2\bigg(\E[\pr^{\leq \reg}\phi^{(p)}_{s,ij}](\tau_1)+\int_{\MM(\tau_1, \tau_1+1)}r^{1+\de} |\pr^{\leq \reg}N_{W,s,ij}^{(p)}|^2\bigg).
\eea

\noindent{\bf Step 3.} Next, we define
\bea
\label{def:tildef0}
\underline{F}^{(p)}_{s,ij} = \left\{\begin{array}{l}
\Big({\square}_{\g_{\tau_1, \tau_2,*}}-(4-2\de_{p0})|q|^{-2} -\big( \widehat{S}_K+\widehat{Q}_K +f_p\big)\Big)(\chi_{\tau_1}(\phi_{aux})^{(p)}_{s,ij})\,\,\,\textrm{on}\,\,\,\MMh(\tmic, \tau_1),\\[2mm]
0\quad\textrm{on}\quad \MMh\setminus\MMh(\tmic, \tau_1),
\end{array}\right.
\eea
where the smooth cut-off $\chi_{\tau_1}=\chi_{\tau_1}(\tau)$ is such that $\chi_{\tau_1}=1$ for $\tau\geq \tau_1$ and 
$\chi_{\tau_1}=0$ for $\tau\leq \tmic$. In particular, \eqref{eq:waveequationdefiningpsiaux} and \eqref{def:tildef0} imply that, for all $\tau\in\mathbb{R}$, 
\beaa
\underline{F}^{(p)}_{s,ij} &=&\big({\square}_{\g_{\tau_1, \tau_2,*}}-(4-2\de_{p0})|q|^{-2} -\big( \widehat{S}_K +\widehat{Q}_K+f_p\big)\big)(\chi_{\tau_1}(\phi_{aux})^{(p)}_{s,ij}) \\
&=& 2\g_{\tau_1, \tau_2,*}^{\a\b}\pr_\a(\chi_{\tau_1})\pr_\b((\phi_{aux})^{(p)}_{s,ij})+\square_{\g_{\tau_1, \tau_2,*}}(\chi_{\tau_1})(\phi_{aux})^{(p)}_{s,ij} -[\widehat{S}_K, \chi_{\tau_1}](\phi_{aux})^{(p)}_{s,ij}.
\eeaa
Now, since $\g_{\tau_1, \tau_2,*}$ satisfies the assumptions \eqref{eq:controloflinearizedinversemetriccoefficients} \eqref{eq:controloflinearizedinversemetriccoefficients:inversegtautau} in view of Step 1, we easily infer, as in (7.23) in \cite{MaSz26},
\bea\lab{eq:structureoftildef0}
\underline{F}^{(p)}_{s,ij} = O(r^{-1})\Big(\chi_{\tau_1}''(\tau), \chi_{\tau_1}'(\tau)\Big)\dk^{\leq 1}(\phi_{aux})^{(p)}_{s,ij}.
\eea
Also, notice from \eqref{eq:waveequationdefiningpsiaux} and finite speed of propagation that $(\phi_{aux})^{(p)}_{s,ij}$ vanishes in the past domain of dependence of $\Sih(\tau_1)\cap\{r\geq 4\Nmic m\}$ which clearly includes $\{r\geq 5\Nmic m\}\cap\{\tmic\leq\tau\leq\tau_1\}$, and hence
\bea\lab{eq:structureoftildef0:support}
\textrm{Supp}\left(\underline{F}^{(p)}_{s,ij}\right)\subset\{\tmic\leq\tau\leq\tau_1\}\cap\{r\leq 5\Nmic m\}.
\eea
By using the control of the energy of $(\phi_{aux})^{(p)}_{s,ij}$ provided by  \eqref{eq:localenergyestimateforpsiauxonSigma} together with \eqref{eq:structureoftildef0} \eqref{eq:structureoftildef0:support}, we obtain
\bea\lab{eq:localenergyestimateforunderlineFsijponSigma}
\nn&& \sum_{i,j=1}^3\sum_{p=0}^2\int_{\MM(\tmic, \tau_1)}r^{1+\de}|\pr^{\leq\reg}\underline{F}_{s,ij}^{(p)}|^2\\ 
&\les& \sum_{i,j=1}^3\sum_{p=0}^2\bigg(\E[\pr^{\leq \reg}\phi^{(p)}_{s,ij}](\tau_1)+\int_{\MM(\tau_1, \tau_1+1)}r^{1+\de}{|\pr^{\leq \reg}N_{W,s,ij}^{(p)}|^2}\bigg).
\eea

\noindent{\bf Step 4.} Next, we introduce the solutions $\breve{\phi}^{(p)}_{s,ij}$ to the following auxiliary system of scalar wave equations, for $p=0,1,2$, $s=\pm 2$, and $i,j=1,2,3$,
\bea\lab{eq:waveequationdefiningfirstextensionbrevepsi}
\bsplit
\big({\square}_{\gam}-(4-2\de_{p0})|q|^{-2}\big)\breve{\phi}^{(p)}_{s,ij} =&  \widehat{S}_K(\breve{\phi}^{(p)}_s)_{ij} +(\widehat{Q}_K\breve{\phi}^{(p)}_s)_{ij}\\
&+f_p\breve{\phi}^{(p)}_{s,ij},\,\,\,\textrm{on}\,\,\,\MMh(\tau_2-1, +\infty)\cup\MMh_{r\geq\varsigma_*(\tau-\frac{7c_*}{8})}(\tau_1,\tau_2-1), \\
 \breve{\phi}^{(p)}_{s,ij} =& (\Pi_2[\widetilde{\phi}_{s}^{(p)}])_{ij}\quad\textrm{on}\quad\Si_{r\leq\varsigma_*(\tau_2-1-\frac{7c_*}{8})}(\tau_2-1),\\ 
  N_{\Sigma(\tau_2-1)}\breve{\phi}^{(p)}_{s,ij}=&N_{\Sigma(\tau_2-1)}(\Pi_2[\widetilde{\phi}^{(p)}_{s}])_{ij}\quad\textrm{on}\quad\Si_{r\leq\varsigma_*(\tau_2-1-\frac{7c_*}{8})}(\tau_2-1),\\
  \breve{\phi}^{(p)}_{s,ij} =& (\Pi_2[\widetilde{\phi}_{s}^{(p)}])_{ij}\quad\textrm{on}\quad\Si_{*,\frac{7c_*}{8}}(\tau_1, \tau_2-1),\\ 
 N_{\Si_{*,\frac{7c_*}{8}}}\breve{\phi}^{(p)}_{s,ij}=&N_{\Si_{*,\frac{7c_*}{8}}}(\Pi_2[\widetilde{\phi}^{(p)}_{s}])_{ij}\quad\textrm{on}\quad\Si_{*,\frac{7c_*}{8}}(\tau_1, \tau_2-1),\\
 \breve{\phi}^{(p)}_{s,ij} =& 0\quad\textrm{on}\quad\Sih_{r\geq\varsigma_*(\tau_1-\frac{7c_*}{8})}(\tau_1),\\
  N_{\Sih(\tau_1)}\breve{\phi}^{(p)}_{s,ij}=&0\quad\textrm{on}\quad\Sih_{r\geq\varsigma_*(\tau_1-\frac{7c_*}{8})}(\tau_1),\\
\end{split}
\eea
where $\widetilde{\phi}_{s,ij}^{(p)}$ is the solution of \eqref{eq:waveequationdefiningfirstextensionwidetildepsi}, and where $\Pi_2$ has been  introduced in \eqref{eq:computationerrorscalarizationdeffect:defPi2}. Then, since we have $\err_{\textrm{TDefect}}[\Pi_2[\breve{\phi}^{(p)}_{s}]]=0$ in view of Lemma \ref{lemma:computationerrorscalarizationdeffect}, we infer from \eqref{eq:waveequationdefiningfirstextensionbrevepsi} that 
\beaa
&&\err_{\textrm{TDefect}}[\breve{\phi}^{(p)}_{s}] =0, \qquad N_{\Sigma}\err_{\textrm{TDefect}}[\breve{\phi}^{(p)}_{s}]= 0,\\
&&\qquad\qquad\textrm{on}\quad\Si_{r\leq\varsigma_*(\tau_2-1-\frac{7c_*}{8})}(\tau_1-1)\cup\Si_{*,\frac{7c_*}{8}}(\tau_1, \tau_2-1)\cup\Sih_{r\geq\varsigma_*(\tau_1-\frac{7c_*}{8})}(\tau_1),
\eeaa
which together with uniqueness for the system of wave equations for $\err_{\textrm{TDefect}}[\breve{\phi}^{(p)}_{s}]$ of Lemma \ref{lemma:waveequationsfortensordeffects} implies 
\bea\lab{eq:thetensordeffectofbrevephivanishesidentically}
\err_{\textrm{TDefect}}[\breve{\phi}^{(p)}_{s}] = 0\quad\textrm{on}\quad\MMh(\tau_2-1, +\infty)\cup\MMh_{r\geq\varsigma_*(\tau-\frac{7c_*}{8})}(\tau_1,\tau_2-1).
\eea
Additionally, note that we have by causality, in view \eqref{eq:waveequationdefiningfirstextensionbrevepsi} and \eqref{eq:bycausalityanddefintionwidetildephipsij=0inMMtau1tau1plus1}, 
\bea\lab{eq:bycausalityanddefintionbrevephipsij=0inMMtau1tau1plus1}
\breve{\phi}^{(p)}_{s,ij}=0\quad\textrm{on}\quad(\MMh\setminus\MM_{r\leq\varsigma_*(\tau-\frac{7c_*}{8})})(\tau_1, \tau_1+1).
\eea

Also, we define $\breve{F}^{(p)}_{s,ij}$ as follows 
\bea
\label{def:breveFpsij}
\breve{F}^{(p)}_{s,ij} = \left\{\begin{array}{l}
\Big({\square}_{\gam}-(4-2\de_{p0})|q|^{-2} -\big( \widehat{S}_K+\widehat{Q}_K +f_p\big)\Big)(\chi_{\tau_2}\chi^{(2)}_{r_*}\widetilde{\phi}^{(p)}_{s,ij}+(1-\chi_{\tau_2}\chi^{(2)}_{r_*})\breve{\phi}^{(p)}_{s,ij})\\
\qquad\qquad\qquad\qquad\qquad\textrm{on}\,\,\,\MM(\tau_2-1, \tau_2)\cup(\MM_{r\geq\varsigma_*(\tau-\frac{7c_*}{8})}(\tau_1+1,\tau_2-1)),\\[2mm]
0\quad\textrm{on}\quad \MMh\setminus(\MM(\tau_2-1, \tau_2)\cup(\MM_{r\geq\varsigma_*(\tau-\frac{7c_*}{8})}(\tau_1+1,\tau_2-1))),
\end{array}\right.
\eea
where the smooth cut-off $\chi_{\tau_2}=\chi_{\tau_2}(\tau)$ is such that $\chi_{\tau_2}=1$ for $\tau\leq \tau_2-2/3$ and 
$\chi_{\tau_2}=0$ for $\tau\geq \tau_2-1/3$, and where the smooth cut-off $\chi^{(2)}_{r_*}=\chi^{(2)}_{r_*}(\tau, r)$ is such that
\beaa
\bsplit
\chi^{(2)}_{r_*}&=1\quad\textrm{for}\quad r\leq\varsigma_*\left(\tau-\frac{29c_*}{32}\right),\\ 
\chi^{(2)}_{r_*}&=0\quad\textrm{for}\quad r\geq\varsigma_*\left(\tau-\frac{31c_*}{32}\right), \qquad \|(r\pr_\tau, r\pr_r)^{\leq 15}\chi^{(2)}_{r_*}\|_{L^\infty}\lesssim 1.
\end{split}
\eeaa
In particular, introducing the smooth cut-off $\chi^{(1)}_{\tau_1}=\chi^{(1)}_{\tau_1}(\tau)$ such that $\chi_{\tau_1}=1$ for $\tau\geq \tau_1+2/3$ and $\chi_{\tau_1}^{(1)}=0$ for $\tau\leq \tau_1+1/3$, \eqref{eq:waveequationdefiningfirstextensionwidetildepsi}, \eqref{eq:waveequationdefiningfirstextensionbrevepsi} and \eqref{def:breveFpsij} imply that, on $\MMh$, 
\beaa
\breve{F}^{(p)}_{s,ij} = \chi^{(1)}_{\tau_1}\Big(2\gam^{\a\b}\pr_\a(\chi_{\tau_2}\chi^{(2)}_{r_*})\pr_\b((\widetilde{\phi}-\breve{\phi})^{(p)}_{s,ij})+\square_{\gam}(\chi_{\tau_2}\chi^{(2)}_{r_*})(\widetilde{\phi}-\breve{\phi})^{(p)}_{s,ij} - [\widehat{S}_K, \chi_{\tau_2}\chi^{(2)}_{r_*}](\widetilde{\phi}-\breve{\phi})^{(p)}_{s,ij}\Big).
\eeaa
We infer from Lemma \ref{lem:specificchoice:normalizedcoord}, Lemma  \ref{lemma:computationoftheMialphajinKerr} and \eqref{eq:bycausalityanddefintionwidetildephipsij=0inMMtau1tau1plus1} \eqref{eq:bycausalityanddefintionbrevephipsij=0inMMtau1tau1plus1} that 
\bea\lab{eq:structureofbreveFpsij}
\breve{F}^{(p)}_{s,ij} &=&  \chi^{(2)}_{r_*}\bigg(-2\chi_{\tau_2}'(\tau)r^{-1}\pr_r(r(\widetilde{\phi}-\breve{\phi})^{(p)}_{s,ij}) +\sum_{k,l}O(r^{-2})\Big(\chi_{\tau_2}''(\tau), \chi_{\tau_2}'(\tau)\Big)\dk^{\leq 1}(\widetilde{\phi}-\breve{\phi})^{(p)}_{s,kl}\bigg)\nn\\
\nn&&+\chi^{(1)}_{\tau_1}\chi_{\tau_2}\sum_{k,l}O(r^{-1})(r\pr_r, r\pr_\tau)\chi_{r_*}\pr(\widetilde{\phi}-\breve{\phi})^{(p)}_{s,kl}\\
&&+\chi^{(1)}_{\tau_1}\chi_{\tau_2}\sum_{k,l}O(r^{-2})\big((r\pr_r, r\pr_\tau)\chi_{r_*}, (r\pr_r, r\pr_\tau)^2\chi_{r_*}\big)(\widetilde{\phi}-\breve{\phi})^{(p)}_{s,kl}.
\eea

\noindent{\bf Step 5.} Next, we introduce the solution $\psi^{(p)}_{s,ij}$ to the following scalar wave equation 
\bea\lab{eq:definitionofthescalarizedsystemwaveequatationpsisijp}
\bsplit
\big({\square}_{\g_{\tau_1, \tau_2,*}}-(4-2\de_{p0})|q|^{-2}\big)\psi^{(p)}_{s,ij} &= \widehat{F}^{(p)}_{s,ij}+\widetilde{F}^{(p)}_{s,ij}+\underline{F}^{(p)}_{s,ij} +\breve{F}^{(p)}_{s,ij}\quad\textrm{on}\quad\MM(\Iti),\\
\psi^{(p)}_{s,ij} &= \widetilde{\phi}^{(p)}_{s,ij}, \quad N_{\Sih(\tau_1)}\psi^{(p)}_{s,ij}=N_{\Sih(\tau_1)}\widetilde{\phi}^{(p)}_{s,ij}\quad\textrm{on}\quad\Sih(\tau_1),\\
\psi^{(p)}_{s,ij} &= \chi_{\tau_1}(\phi_\AA)^{(p)}_{s,ij}\quad\textrm{on}\quad\AA\cap\{\tau\leq\tau_1\},\\
\psi^{(p)}_{s,ij} &= 0\quad\textrm{on}\quad\II_+\cap\{\tau\leq\tau_1\},
\end{split}
\eea
where we recall that $(\phi_\AA)^{(p)}_{s,ij}$ is a smooth extension of $\phi^{(p)}_{s,ij}$  from $\AA\cap\{\tau\geq \tau_1\}$ to $\AA\cap\{\tmic\leq\tau\leq\tau_1\}$ satisfying \eqref{eq:propertyextensionpsiAAofpsitotauleqtau1}. In particular, note by causality, using in particular Definition \ref{def:definitionofthetimetauR} of $\tmic$, that we have
\bsub
\lab{eq:psispij:intermsoftildephibrevephiandphiaux}
\bea
\psi^{(p)}_{s,ij} &=& \chi_{\tau_2}\chi^{(2)}_{r_*}\widetilde{\phi}^{(p)}_{s,ij}+(1-\chi_{\tau_2}\chi^{(2)}_{r_*})\breve{\phi}^{(p)}_{s,ij}\quad\textrm{on}\quad\MMh(\tau_1, +\infty), \\ 
\psi^{(p)}_{s,ij} &=& \chi_{\tau_1}(\phi_{aux})^{(p)}_{s,ij}\quad\textrm{on}\quad\MM(\tmic, \tau_1)\cap\{r\leq (\Nmic+1)m\}.
\eea
\esub
On the other hand, we have 
\beaa
\widetilde{\phi}^{(p)}_{s,ij}=\phi^{(p)}_{s,ij}\quad\textrm{on}\quad \MM_{r\leq\varsigma_*(\tau-\frac{5c_*}{8})}(\tau_1+1, \tau_2-3)
\eeaa
by causality in view of \eqref{eq:waveequationdefiningfirstextensionwidetildepsi}, and we thus deduce
\beaa
\bsplit
\psi^{(p)}_{s,ij}&=\phi^{(p)}_{s,ij}\quad\textrm{on}\quad \MM_{r\leq\varsigma_*(\tau-\frac{5c_*}{8})}(\tau_1+1, \tau_2-3), \\
\dk^{\leq\reg}(\psi^{(p)}_{s,ij}) &= 0\quad\textrm{on}\quad \Si(\tmic)\cap 
\{r\leq(\Nmic+1)m\}\,\,\,\,\forall\,\reg\in\mathbb{N}.
\end{split}
\eeaa

\noindent{\bf Step 6.} We now derive local energy estimates for $\psi^{(p)}_{s,ij}$ on  $\tau\in[\tmic, \tau_1+3]$ using the system of scalar wave equations in \eqref{eq:definitionofthescalarizedsystemwaveequatationpsisijp}. Using \eqref{eq:localenergyestimate:future:unweigthedderivatives:bis} for the wave system\footnote{In fact, we also need to derive a local energy estimate for $\psiss{ij}{p}$ on $(\MMh\setminus\MM)(\tau_1,\tau_1+3)$ with data on $\Si_*(\tau_1,\tau_1+3)$ and $(\Sih\setminus\Si)(\tau_1)$ which is a straightforward adaptation of \eqref{eq:localenergyestimate:future:unweigthedderivatives:bis} left to the reader. Moreover, the local energy estimate for $\psiss{ij}{p}$ on $\MMh(\tmic,\tau_1)$ requires an extension of \eqref{eq:localenergyestimate:past:bis:unweigthedderivatives} to $\MMh$ which can be found in (6.37b) of \cite{MaSz26}.} consisting of the wave equations for $\psiss{ij}{p}$ in \eqref{eq:definitionofthescalarizedsystemwaveequatationpsisijp} and the wave equations of $\phiss{ij}{p}$ in $\MM(\tau_1,\tau_1+3)$, using \eqref{eq:localenergyestimate:past:bis:unweigthedderivatives} for the wave equations for $\psiss{ij}{p}$ in \eqref{eq:definitionofthescalarizedsystemwaveequatationpsisijp} in $\MMh(\tmic,\tau_1)$,  and using the initialization of $\psi^{(p)}_{s,ij}$ on $\Si(\tau_1)$, see \eqref{eq:definitionofthescalarizedsystemwaveequatationpsisijp}, we infer, for all $\reg\leq {14}$, 
\beaa
&&\sum_{i,j=1}^3\sum_{p=0}^2\EMFh_\de[\pr^{\leq \reg}\psi^{(p)}_{s,ij}](\tmic, \tau_1+3)+\sum_{i,j=1}^3\sum_{p=0}^2\EMF_\de[\pr^{\leq \reg}\phi^{(p)}_{s,ij}](\tt_1, \tau_1+3)\\ 
&\les& \sum_{i,j=1}^3\sum_{p=0}^2\Big(\Eh[\pr^{\leq \reg}\widetilde{\phi}^{(p)}_{s,ij}](\tau_1)
+\E[\pr^{\leq \reg}{\phi}^{(p)}_{s,ij}](\tau_1)+\F_\AA[\pr^{\leq \reg}(\phi_\AA)^{(p)}_{s,ij}](\tmic, \tau_1)\Big)\\
&&+\sum_{i,j=1}^3\sum_{p=0}^2\int_{\MM(\tau_1, \tau_1+3)}r^{1+\de}|\pr^{\leq\reg}N_{W,s,ij}^{(p)}|^2+\sum_{i,j=1}^3\sum_{p=0}^2\int_{\MM(\tmic, \tau_1)}r^{1+\de}|\pr^{\leq\reg}\underline{F}_{s,ij}^{(p)}|^2,
\eeaa
which together with  \eqref{eq:localenergyestimateforwidetildepsionSigma},  \eqref{eq:propertyextensionpsiAAofpsitotauleqtau1}  and \eqref{eq:localenergyestimateforunderlineFsijponSigma}  yields, for all $\reg\leq {14}$, 
\beaa
&&\sum_{p=0}^2\sum_{p=0}^2\EMFh_\de[\pr^{\leq \reg}\psi^{(p)}_{s}](\tmic, \tau_1+3)\nn\\
&\les& \sum_{i,j=1}^3\sum_{p=0}^2\bigg(\E[\pr^{\leq \reg}\phi^{(p)}_{s,ij}](\tau_1)+\int_{\MM(\tau_1, \tau_1+3)}r^{1+\de} |\pr^{\leq \reg}N_{W,s,ij}^{(p)}|^2\bigg).
\eeaa

\noindent{\bf Step 7.} Next, we derive local energy estimates for $\psi^{(p)}_{s,ij}$ on $\RR_*(\tau_1, \tau_2-1)$ and $\MM(\tau_2-3, \tau_2)$. To this end, we first derive local energy estimates for $\widetilde{\phi}^{(p)}_{s,ij}$. Applying the local energy estimate \eqref{eq:localenergyestimate:future:unweigthedderivatives:RR*} with $\reg\leq {14}$ to \eqref{eq:waveequationdefiningfirstextensionwidetildepsi}, we have
\beaa
\nn&&\sum_{i,j=1}^3\sum_{p=0}^2\bigg(\EM_{0,r\geq\varsigma_*(\tau-c_*/2)}[\pr^{\leq \reg}\widetilde{\phi}^{(p)}_{s,ij}](\tau_1, \tau_2-1)+\sup_{\frac{c_*}{2}\leq c\leq c_*}\F_{\Si_{*,c}}[\pr^{\leq\reg}\widetilde{\phi}^{(p)}_{s,ij}](\tau_1, \tau_2-1)\bigg)\\ 
&\les& \sum_{i,j=1}^3\sum_{p=0}^2\bigg(\E[\pr^{\leq \reg}\widetilde{\phi}^{(p)}_{s,ij}](\tau_1)+\F_{\Si_{*,c_*/2}}[\pr^{\leq\reg}\widetilde{\phi}^{(p)}_{s,ij}](\tau_1, \tau_2-1)+\int_{\RR_*(\tau_1, \tau_2-2)}r|\pr^{\leq \reg}N_{W,s,ij}^{(p)}|^2\bigg),
\eeaa
which, together with \eqref{eq:localenergyestimateforwidetildepsionSigma}, \eqref{eq:bycausalityanddefintionwidetildephipsij=0inMMtau1tau1plus1}, and the fact that $\widetilde{\phi}^{(p)}_{s,ij}=\phi^{(p)}_{s,ij}$ on $(\MM\setminus\RR_*)(\tau_1+1, \tau_2-3)$  in view of Step 5, yields, for $\reg\leq {14}$,
\bea\lab{eq:controllingwidetildephipsijlocalenergyestimateRR*tau1plus3tau2minus3}
\nn&&\sum_{i,j=1}^3\sum_{p=0}^2\bigg(\EM_{0,r\geq\varsigma_*(\tau-c_*/2)}[\pr^{\leq \reg}\widetilde{\phi}^{(p)}_{s,ij}](\tau_1, \tau_2-1)+\sup_{\frac{c_*}{2}\leq c\leq c_*}\F_{\Si_{*,c}}[\pr^{\leq\reg}\widetilde{\phi}^{(p)}_{s,ij}](\tau_1, \tau_2-1)\bigg)\\ 
\nn&\les& \sum_{i,j=1}^3\sum_{p=0}^2\bigg(\E[\pr^{\leq \reg}\phi^{(p)}_{s,ij}](\tau_1)+\F_{\Si_{*,c_*/2}}[\pr^{\leq\reg}\phi^{(p)}_{s,ij}](\tau_1+1, \tau_2-1)\\
&&+\int_{\MM(\tau_1, \tau_1+1)\cup\RR_*(\tau_1+1, \tau_2-2)}r^{1+\de}|\pr^{\leq \reg}N_{W,s,ij}^{(p)}|^2\bigg).
\eea
Next, applying the local energy estimate \eqref{eq:localenergyestimate:future:unweigthedderivatives:bis} with $\reg\leq {14}$ to \eqref{eq:waveequationdefiningfirstextensionwidetildepsi}, we have
\beaa
\nn&&\sum_{i,j=1}^3\sum_{p=0}^2\EMF_\de[\pr^{\leq \reg}\widetilde{\phi}^{(p)}_{s,ij}](\tau_2-3, \tau_2)\\ 
&\les& \sum_{i,j=1}^3\sum_{p=0}^2\bigg(\E[\pr^{\leq \reg}\widetilde{\phi}^{(p)}_{s,ij}](\tau_2-3)+\E[\pr^{\leq \reg}\phi^{(p)}_{s,ij}](\tau_2-3)+\int_{\MM(\tau_2-3, \tau_2-2)}r^{1+\de} |\pr^{\leq \reg}N_{W,s,ij}^{(p)}|^2\bigg)
\eeaa
which, together with \eqref{eq:controllingwidetildephipsijlocalenergyestimateRR*tau1plus3tau2minus3} and the fact that $\widetilde{\phi}^{(p)}_{s,ij}=\phi^{(p)}_{s,ij}$ on $(\MM\setminus\RR_*)(\tau_1+1, \tau_2-3)$  in view of Step 5, yields, for $\reg\leq {14}$,
\bea\lab{eq:localenergyestimateforwidetildephisijponSigmatau2minus3tau2}
\nn&&\sum_{i,j=1}^3\sum_{p=0}^2\bigg(\EM_{0,r\geq\varsigma_*(\tau-c_*/2)}[\pr^{\leq \reg}\widetilde{\phi}^{(p)}_{s,ij}](\tau_1, \tau_2-1)\\
\nn&&+\sup_{\frac{c_*}{2}\leq c\leq c_*}\F_{\Si_{*,c}}[\pr^{\leq\reg}\widetilde{\phi}^{(p)}_{s,ij}](\tau_1, \tau_2-1)+\EMF_\de[\pr^{\leq \reg}\widetilde{\phi}^{(p)}_{s,ij}](\tau_2-3, \tau_2)\bigg)\\ 
\nn&\les& \sum_{i,j=1}^3\sum_{p=0}^2\bigg(\E[\pr^{\leq \reg}\phi^{(p)}_{s,ij}](\tau_1)+\E[\pr^{\leq \reg}\phi^{(p)}_{s,ij}](\tau_2-3)+\F_{\Si_{*,c_*/2}}[\pr^{\leq\reg}\phi^{(p)}_{s,ij}](\tau_1+1, \tau_2-1)\\
&&+\int_{\MM(\tau_1, \tau_1+1)\cup\RR_*(\tau_1+1, \tau_2-3)\cup\MM(\tau_2-3, \tau_2-2)}r^{1+\de}|\pr^{\leq \reg}N_{W,s,ij}^{(p)}|^2\bigg).
\eea

Next, applying the energy estimate \eqref{eq:localenergyestimate:future:unweigthedderivatives:RR*} with $\reg\leq {14}$ to $\breve{\phi}^{(p)}_{s,ij}$ solution of \eqref{eq:waveequationdefiningfirstextensionbrevepsi}, we have
\beaa
\nn&&\sum_{i,j=1}^3\sum_{p=0}^2\bigg(\EM_{0,r\geq\varsigma_*(\tau-7c_*/8)}[\pr^{\leq \reg}\breve{\phi}^{(p)}_{s,ij}](\tau_1, \tau_2-1)+\sup_{\frac{7c_*}{8}\leq c\leq c_*}\F_{\Si_{*,c}}[\pr^{\leq\reg}\breve{\phi}^{(p)}_{s,ij}](\tau_1, \tau_2-1)\bigg)\\ 
&\les& \sum_{i,j=1}^3\sum_{p=0}^2\F_{\Si_{*,\frac{7c_*}{8}}}[\pr^{\leq\reg}\Pi_2[\widetilde{\phi}_{s}^{(p)}]](\tau_1, \tau_2-1)\les \sum_{i,j=1}^3\sum_{p=0}^2\F_{\Si_{*,\frac{7c_*}{8}}}[\pr^{\leq\reg}\widetilde{\phi}_{s}^{(p)}](\tau_1, \tau_2-1),
\eeaa
and applying the local energy estimate \eqref{eq:localenergyestimate:future:unweigthedderivatives:bis} with $\reg\leq {14}$ to $\breve{\phi}^{(p)}_{s,ij}$, we have
\beaa
&&\sum_{i,j=1}^3\sum_{p=0}^2\EMF_\de[\pr^{\leq \reg}\breve{\phi}^{(p)}_{s,ij}](\tau_2-1, \tau_2)\\
 &\les& \sum_{i,j=1}^3\sum_{p=0}^2\Big(\E_{r\leq\varsigma_*(\tau_2-1-7c_*/2)}[\pr^{\leq \reg}(\Pi_2[\widetilde{\phi}_{s}^{(p)}])](\tau_2-1)+\E_{r\geq\varsigma_*(\tau_2-1-7c_*/2)}[\pr^{\leq \reg}\breve{\phi}_{s}^{(p)}](\tau_2-1)\Big)\\
&\les& \sum_{i,j=1}^3\sum_{p=0}^2\Big(\E_{r\geq\varsigma_*(\tau_2-1-7c_*/2)}[\pr^{\leq \reg}\widetilde{\phi}_{s}^{(p)}](\tau_2-1)+\E_{r\geq\varsigma_*(\tau_2-1-7c_*/2)}[\pr^{\leq \reg}\breve{\phi}_{s}^{(p)}](\tau_2-1)\Big)
\eeaa
which together with \eqref{eq:localenergyestimateforwidetildephisijponSigmatau2minus3tau2} yields, for $\reg\leq {14}$,
\bea\lab{eq:localenergyestimateforbrevephisijponSigmatau2minus1tau2}
\nn&&\sum_{i,j=1}^3\sum_{p=0}^2\bigg(\EM_{0,r\geq\varsigma_*(\tau-7c_*/8)}[\pr^{\leq \reg}\breve{\phi}^{(p)}_{s,ij}](\tau_1, \tau_2-1)+\sup_{\frac{7c_*}{8}\leq c\leq c_*}\F_{\Si_{*,c}}[\pr^{\leq\reg}\breve{\phi}^{(p)}_{s,ij}](\tau_1, \tau_2-1)\bigg)\\ 
\nn&&+\sum_{i,j=1}^3\sum_{p=0}^2\EMF_\de[\pr^{\leq \reg}\breve{\phi}^{(p)}_{s,ij}](\tau_2-1, \tau_2)\\ 
\nn&\les& \sum_{i,j=1}^3\sum_{p=0}^2\bigg(\E[\pr^{\leq \reg}\phi^{(p)}_{s,ij}](\tau_1)+\E[\pr^{\leq \reg}\phi^{(p)}_{s,ij}](\tau_2-3)+\F_{\Si_{*,c_*/2}}[\pr^{\leq\reg}\phi^{(p)}_{s,ij}](\tau_1+1, \tau_2-1)\\
&&+\int_{\MM(\tau_1, \tau_1+1)\cup\RR_*(\tau_1+1, \tau_2-3)\cup\MM(\tau_2-3, \tau_2-2)}r^{1+\de}|\pr^{\leq \reg}N_{W,s,ij}^{(p)}|^2\bigg).
\eea
Since we have, in view of Step 5,
\beaa
\psi^{(p)}_{s,ij} &=& \chi_{\tau_2}\chi^{(2)}_{r_*}\widetilde{\phi}^{(p)}_{s,ij}+(1-\chi_{\tau_2}\chi^{(2)}_{r_*})\breve{\phi}^{(p)}_{s,ij}\quad\textrm{on}\quad\MM(\tau_1, +\infty), 
\eeaa
we immediately infer from \eqref{eq:localenergyestimateforwidetildephisijponSigmatau2minus3tau2} and 
\eqref{eq:localenergyestimateforbrevephisijponSigmatau2minus1tau2}, for $\reg\leq {14}$,
\beaa
\nn&&\sum_{i,j=1}^3\sum_{p=0}^2\bigg(\EM_{0,r\geq\varsigma_*(\tau-c_*/2)}[\pr^{\leq \reg}\psi^{(p)}_{s,ij}](\tau_1, \tau_2-1)\\
\nn&&+\sup_{\frac{c_*}{2}\leq c\leq c_*}\F_{\Si_{*,c}}[\pr^{\leq\reg}\psi^{(p)}_{s,ij}](\tau_1, \tau_2-1)+\EMF_\de[\pr^{\leq \reg}\psi^{(p)}_{s,ij}](\tau_2-3, \tau_2)\bigg)\\ 
\nn&\les& \sum_{i,j=1}^3\sum_{p=0}^2\bigg(\E[\pr^{\leq \reg}\phi^{(p)}_{s,ij}](\tau_1)+\E[\pr^{\leq \reg}\phi^{(p)}_{s,ij}](\tau_2-3)+\F_{\Si_{*,c_*/2}}[\pr^{\leq\reg}\phi^{(p)}_{s,ij}](\tau_1+1, \tau_2-1)\\
&&+\int_{\MM(\tau_1, \tau_1+1)\cup\RR_*(\tau_1+1, \tau_2-3)\cup\MM(\tau_2-3, \tau_2-2)}r^{1+\de}|\pr^{\leq \reg}N_{W,s,ij}^{(p)}|^2\bigg).
\eeaa
which is the desired estimate \eqref{eq:localenergyestimateforpsisijponSigmatau2minus3tau2:partialderivatives}. 

\noindent{\bf Step 8.} Next, we estimate the tensorization defect $\err_{\textrm{TDefect}}[\psi^{(p)}_{s}]$ corresponding to the family of complex-valued scalars $\psi_{ij}$ as introduced in  Definition \ref{def:definitionofthenotationerrforthescalarizationdefect}. First, recall from Step 5 that 
\beaa
\psi^{(p)}_{s,ij} &=& \phi^{(p)}_{s,ij}\quad\textrm{on}\quad\MM_{r\leq\varsigma_*(\tau-\frac{5c_*}{8})}(\tau_1+1, \tau_2-3), \\ 
\psi^{(p)}_{s,ij} &=& \breve{\phi}^{(p)}_{s,ij}\quad\textrm{on}\quad\MMh(\tau_2, +\infty)\cup\MMh_{r\geq\varsigma_*(\tau-\frac{31c_*}{32})}(\tau_1, \tau_2),
\eeaa
which together with \eqref{eq:thetensordeffectofbrevephivanishesidentically}, and the fact that $\phi^{(p)}_{s,ij}=\pmb\phi^{(p)}_{s}(\Om_i, \Om_j)$ by definition, yields
\beaa
\err_{\textrm{TDefect}}[\psi^{(p)}_{s}]=0\quad\textrm{on}\quad\MM_{r\leq\varsigma_*(\tau-\frac{5c_*}{8})}(\tau_1+1, \tau_2-3)\cup\MMh(\tau_2, +\infty)\cup\MMh_{r\geq\varsigma_*(\tau-\frac{31c_*}{32})}(\tau_1, \tau_2).
\eeaa

Then, we consider the region $\MM_{r\leq\varsigma_*(\tau-\frac{5c_*}{8})}(\tau_2-3, \tau_2-2)$ and, to this end, we introduce the following auxiliary coupled system of tensorial wave equations, for $s=\pm 2$, $p=0,1,2$,
\beaa
&&\bigg(\squared_2 -\frac{4ia\cos\th}{|q|^2}\nab_{\pr_{\tt}}- \frac{4-2\de_{p0}}{\qs}\bigg){(\pmb\phi_{aux,1})_s^{(p)}}\\
&=& \chi^{(1)}_{\tau_1, \tau_2}\big(\L_{s}^{(p)}[(\pmb\phi_{aux,1})_s]+\N_{W,s}^{(p)}\big),  \quad\textrm{on}\quad\MM_{r\leq\varsigma_*(\tau-\frac{5c_*}{8})}(\tau_2-3, \tau_2),\\
&&(\pmb\phi_{aux,1})^{(p)}_s =\pmb\phi^{(p)}_{s}, \quad \nab_{N_{\Sigma(\tau_2-3)}}(\pmb\phi_{aux,1})^{(p)}_s =\nab_{N_{\Sigma(\tau_2-3)}}\pmb\phi^{(p)}_{s} \quad\textrm{on}\quad\Sigma_{r\leq\varsigma_*(\tau_2-3-\frac{5c_*}{8})}(\tau_2-3).
\eeaa
Defining $(\phi_{aux,1})_{s,ij}^{(p)}:=(\pmb\phi_{aux,1})_s^{(p)}(\Om_i, \Om_j)$, one easily checks that $(\phi_{aux,1})_{s,ij}^{(p)}$ satisfies the same system of wave equations as $\psi_{s,ij}^{(p)}$ on $\MM_{r\leq\varsigma_*(\tau-\frac{5c_*}{8})}(\tau_2-3, \tau_2-2)$ and the same initial data on $\Si_{r\leq\varsigma_*(\tau_2-3-\frac{5c_*}{8})}(\tau_2-3)$. By uniqueness for the wave equation, we infer 
\beaa
\psi_{s,ij}^{(p)}=(\phi_{aux,1})_{s,ij}^{(p)}=(\pmb\phi_{aux,1})_s^{(p)}(\Om_i, \Om_j)\quad\textrm{on}\quad\MM_{r\leq\varsigma_*(\tau-\frac{5c_*}{8})}(\tau_2-3, \tau_2-2)
\eeaa
which implies that $\err_{\textrm{TDefect}}[\psi^{(p)}_{s}]=0$ on $\MM_{r\leq\varsigma_*(\tau-\frac{5c_*}{8})}(\tau_2-3, \tau_2-2)$. In view of the above, we deduce 
\beaa
\err_{\textrm{TDefect}}[\psi^{(p)}_{s}]=0\quad\textrm{on}\quad\MMh(\tau_1+1, +\infty)\setminus(\MM(\tau_2-2, \tau_2)\cup\MM_{\varsigma_*(\tau-\frac{5c_*}{8})\leq r\leq\varsigma_*(\tau-\frac{31c_*}{32})}(\tau_1+1, \tau_2-2)).
\eeaa

Next, we derive energy estimates for $\err_{\textrm{TDefect}}[\psi^{(p)}_{s}]$ on $\MM(\tau_2-2, \tau_2)\cup\RR_*(\tau_1+1, \tau_2-2)$. To this end, we first derive a system of wave equations for $\err_{\textrm{TDefect}}[\widetilde{\phi}^{(p)}_{s}]$. In view of \eqref{eq:waveequationdefiningfirstextensionwidetildepsi}, we have on $\MM(\tau_2-2, \tau_2)\cup\RR_*(\tau_1+1, \tau_2-2)$
\beaa
\big({\square}_{\g_{\tau_1, \tau_2,*}}-(4-2\de_{p0})|q|^{-2}\big)\widetilde{\phi}^{(p)}_{s,ij} &=& \chi_{\tau_1, \tau_2}\chi_{r_*}\big( \widehat{S}(\widetilde{\pmb\phi}^{(p)}_s)_{ij} +(\widehat{Q}\widetilde{\pmb\phi}^{(p)}_s)_{ij}\big)\\
&& +(1-\chi_{\tau_1, \tau_2}\chi_{r_*})\big( \widehat{S}_K(\widetilde{\pmb\phi}^{(p)}_s)_{ij} +(\widehat{Q}_K\widetilde{\pmb\phi}^{(p)}_s)_{ij}+f_p\widetilde{\phi}^{(p)}_{s,ij}\big)
\eeaa 
and hence, together with \eqref{eq:assumptionsonregulartripletinperturbationsofKerr:0}, we infer, on $\MM(\tau_2-2, \tau_2)\cup\RR_*(\tau_1+1, \tau_2-2)$,
\bea\lab{eq:waveequationdefiningfirstextensionwidetildepsi:tau2minus2totau2}
\bsplit
\Big({\square}_{\g_{\tau_1, \tau_2,*}}+V\Big)\widetilde{\phi}^{(p)}_{s,ij} =&  \widehat{S}(\widetilde{\pmb\phi}^{(p)}_s)_{ij} +(\widehat{Q}\widetilde{\pmb\phi}^{(p)}_s)_{ij} +\Ga_g\dk^{\leq 1}\widetilde{\phi}^{(p)}_{s},\\
V:=& -(4-2\de_{p0})|q|^{-2}-(1-\chi_{\tau_1, \tau_2})f_p.
\end{split}
\eea 
Also, in view of Lemma \ref{lemma:computationoftheMialphajinKerr} and \eqref{eq:assumptionsonregulartripletinperturbationsofKerr:0}, and \eqref{eq:definitionofthenotationerrforthescalarizationdefect}, we have
\beaa
M_{i\tau}^l\widetilde{\phi}^{(p)}_{s,lj}+M_{j\tau}^l\widetilde{\phi}^{(p)}_{s,il} &=& -\frac{2amr\cos\th}{|q|^4}\Big(\in_{ilk}x^k\widetilde{\phi}^{(p)}_{s,lj}+\in_{jlk}x^k\widetilde{\phi}^{(p)}_{s,il}\Big)+\Ga_g\widetilde{\phi}^{(p)}_{s}\\
&=& -\frac{2amr\cos\th}{|q|^4}\Big(-2i\widetilde{\phi}^{(p)}_{s,ij}+(\err_{\textrm{TDefect},5}[\psi])_{ij}+(\err_{\textrm{TDefect},5}[\psi])_{ji}\\
&&+i(\err_{\textrm{TDefect},1}[\psi])_{ij}+\in_{jlk}x^k(\err_{\textrm{TDefect},1}[\psi])_{il}\Big)+\Ga_g\widetilde{\phi}^{(p)}_{s}.
\eeaa
Together with \eqref{eq:waveequationdefiningfirstextensionwidetildepsi:tau2minus2totau2} and Lemma \ref{lemma:waveequationsfortensordeffects}, using also Lemma \ref{lemma:computationoftheMialphajinKerr}, \eqref{eq:assumptionsonregulartripletinperturbationsofKerr:0}, and the properties of $\g_{\tau_1, \tau_2,*}$, we infer on $\MM(\tau_2-2, \tau_2)\cup\RR_*(\tau_1+1, \tau_2-2)$
\beaa
\square_{\g_{\tau_1, \tau_2,*}}\err_{\textrm{TDefect}}[\widetilde{\phi}^{(p)}_{s}] =
O(r^{-2})\dk^{\leq 1}\err_{\textrm{TDefect}}[\widetilde{\phi}^{(p)}_{s}]+\Ga_g\dk^{\leq 1}\widetilde{\phi}^{(p)}_{s}.
\eeaa
Since the initial data for $\err_{\textrm{TDefect}}[\widetilde{\phi}^{(p)}_{s}]$ is trivial on $\Si_{*,c_*/2}(\tau_1, \tau_2-3)$, $\Si_{r\leq\varsigma_*(\tau_2-3-c_*/2)}(\tau_2-3)$ and $\Si(\tau_1+1)$, using the local energy estimates \eqref{eq:localenergyestimate:future:unweigthedderivatives:RR*} on $\RR_*(\tau_1+1, \tau_2-2)$ and \eqref{eq:localenergyestimate:future:unweigthedderivatives:bis} on $\MM(\tau_2-2, \tau_2)$, we infer, for $\reg\leq 14$,
\beaa
\bsplit
&\EM_{0,r\geq\varsigma_*(\tau-c_*/2)}[\pr^{\leq \reg}\err_{\textrm{TDefect}}[\widetilde{\phi}^{(p)}_{s}]](\tau_1+1, \tau_2-2)\\
&+\sup_{c_*/2\leq c\leq c_*}\F_{\Si_{*,c}}[\pr^{\leq \reg}\err_{\textrm{TDefect}}[\widetilde{\phi}^{(p)}_{s}]](\tau_1+1, \tau_2-2)+\EMF_\de[\pr^{\leq \reg}\err_{\textrm{TDefect}}[\widetilde{\phi}^{(p)}_{s}]](\tau_2-2, \tau_2)\\
\les& \ep^2\M_{0, r\geq\varsigma_*(\tau-c_*/2)}[\pr^{\leq \reg}\widetilde{\phi}^{(p)}_{s}](\tau_2+1, \tau_2-2)+\ep^2\EM_\de[\pr^{\leq \reg}\widetilde{\phi}^{(p)}_{s}](\tau_2-2, \tau_2).
\end{split}
\eeaa
Together with \eqref{eq:localenergyestimateforwidetildephisijponSigmatau2minus3tau2}, we deduce, for $\reg\leq 14$,
\begin{align}
\nn&\sum_{p=0}^2\Big(\EM_{0,r\geq\varsigma_*(\tau-c_*/2)}[\pr^{\leq \reg}\err_{\textrm{TDefect}}[\widetilde{\phi}^{(p)}_{s}]](\tau_1+1, \tau_2-2)\\
\nn&+\sup_{c_*/2\leq c\leq c_*}\F_{\Si_{*,c}}[\pr^{\leq \reg}\err_{\textrm{TDefect}}[\widetilde{\phi}^{(p)}_{s}]](\tau_1+1, \tau_2-2)+\EMF_\de[\pr^{\leq \reg}\err_{\textrm{TDefect}}[\widetilde{\phi}^{(p)}_{s}]](\tau_2-2, \tau_2)\Big)\\ 
\nn\les& \ep^2\sum_{i,j=1}^3\sum_{p=0}^2\bigg(\E[\pr^{\leq \reg}\phi^{(p)}_{s,ij}](\tau_1)+\E[\pr^{\leq \reg}\phi^{(p)}_{s,ij}](\tau_2-3)+\F_{\Si_{*,c_*/2}}[\pr^{\leq\reg}\phi^{(p)}_{s,ij}](\tau_1+1, \tau_2-1)\\
&+\int_{\MM(\tau_1, \tau_1+1)\cup\RR_*(\tau_1+1, \tau_2-3)\cup\MM(\tau_2-3, \tau_2-2)}r^{1+\de}|\pr^{\leq \reg}N_{W,s,ij}^{(p)}|^2\bigg).
\end{align}
Now, recalling from Step 5 that 
\beaa
\psi^{(p)}_{s,ij} = \chi_{\tau_2}\chi^{(2)}_{r_*}\widetilde{\phi}^{(p)}_{s,ij}+(1-\chi_{\tau_2}\chi^{(2)}_{r_*})\breve{\phi}^{(p)}_{s,ij}\quad\textrm{on}\quad\MMh(\tau_1, +\infty), 
\eeaa
we infer
\beaa
\err_{\textrm{TDefect}}[\psi^{(p)}_{s}] = \chi_{\tau_2}\chi^{(2)}_{r_*}\err_{\textrm{TDefect}}[\widetilde{\phi}^{(p)}_{s}]+(1-\chi_{\tau_2}\chi^{(2)}_{r_*})\err_{\textrm{TDefect}}[\breve{\phi}^{(p)}_{s}]\quad\textrm{on}\quad\MMh(\tau_1, +\infty), 
\eeaa
which together with \eqref{eq:thetensordeffectofbrevephivanishesidentically} yields
\beaa
\err_{\textrm{TDefect}}[\psi^{(p)}_{s}] = \chi_{\tau_2}\chi^{(2)}_{r_*}\err_{\textrm{TDefect}}[\widetilde{\phi}^{(p)}_{s}]\quad\textrm{on}\quad\MMh(\tau_1, +\infty).
\eeaa
In view of the above control for $\err_{\textrm{TDefect}}[\widetilde{\phi}^{(p)}_{s}]$, this yields, for $\reg\leq 14$,
\begin{align}\lab{eq:controlscalarizationdefectforwidetildephiontau2minus2tau2}
\nn&\sum_{p=0}^2\Big(\EM_{0,r\geq\varsigma_*(\tau-c_*/2)}[\pr^{\leq \reg}\err_{\textrm{TDefect}}[\psi^{(p)}_{s}]](\tau_1+1, \tau_2-2)\\
\nn&+\sup_{c_*/2\leq c\leq c_*}\F_{\Si_{*,c}}[\pr^{\leq \reg}\err_{\textrm{TDefect}}[\psi^{(p)}_{s}]](\tau_1+1, \tau_2-2)+\EMF_\de[\pr^{\leq \reg}\err_{\textrm{TDefect}}[\psi^{(p)}_{s}]](\tau_2-2, \tau_2)\Big)\\ 
\nn\les& \ep^2\sum_{i,j=1}^3\sum_{p=0}^2\bigg(\E[\pr^{\leq \reg}\phi^{(p)}_{s,ij}](\tau_1)+\E[\pr^{\leq \reg}\phi^{(p)}_{s,ij}](\tau_2-3)+\F_{\Si_{*,c_*/2}}[\pr^{\leq\reg}\phi^{(p)}_{s,ij}](\tau_1+1, \tau_2-1)\\
&+\int_{\MM(\tau_1, \tau_1+1)\cup\RR_*(\tau_1+1, \tau_2-3)\cup\MM(\tau_2-3, \tau_2-2)}r^{1+\de}|\pr^{\leq \reg}N_{W,s,ij}^{(p)}|^2\bigg).
\end{align}

In addition, we also control $\widetilde{\phi}^{(p)}_{s,ij}- \breve{\phi}^{(p)}_{s,ij}$ on $\MM(\tau_2-1, \tau_2)\cup\MM_{r\geq\varsigma_*(\tau-\frac{7c_*}{8})}(\tau_1,\tau_2-1)$. In view of \eqref{eq:waveequationdefiningfirstextensionwidetildepsi} and \eqref{eq:waveequationdefiningfirstextensionbrevepsi}, $\widetilde{\phi}^{(p)}_{s,ij}- \breve{\phi}^{(p)}_{s,ij}$ satisfies
\beaa
\bsplit
\big({\square}_{\gam}-(4-2\de_{p0})|q|^{-2}\big)(\widetilde{\phi}-\breve{\phi})^{(p)}_{s,ij} =&  \big(\widehat{S}_K +\widehat{Q}_K+f_p\big)(\widetilde{\phi}-\breve{\phi})^{(p)}_{s,ij}\\
&\qquad\qquad\textrm{on}\quad\MM(\tau_2-1,\tau_2)\cup\MM_{r\geq\varsigma_*(\tau-\frac{7c_*}{8})}(\tau_1,\tau_2-1), \\
(\widetilde{\phi}-\breve{\phi})^{(p)}_{s,ij} =& (\widetilde{\phi}-\Pi_2[\widetilde{\phi}_{s}^{(p)}])_{ij}\quad\textrm{on}\quad\Si_{r\leq\varsigma_*(\tau_2-1-\frac{7c_*}{8})}(\tau_2-1),\\ 
  N_{\Sigma(\tau_2-1)}(\widetilde{\phi}-\breve{\phi})^{(p)}_{s,ij}=&N_{\Sigma(\tau_2-1)}(\widetilde{\phi}-\Pi_2[\widetilde{\phi}^{(p)}_{s}])_{ij}\quad\textrm{on}\quad\Si_{r\leq\varsigma_*(\tau_2-1-\frac{7c_*}{8})}(\tau_2-1),\\
  (\widetilde{\phi}-\breve{\phi})^{(p)}_{s,ij} =& (\widetilde{\phi}-\Pi_2[\widetilde{\phi}_{s}^{(p)}])_{ij}\quad\textrm{on}\quad\Si_{*,\frac{7c_*}{8}}(\tau_1, \tau_2-1),\\ 
 N_{\Si_{*,\frac{7c_*}{8}}}(\widetilde{\phi}-\breve{\phi})^{(p)}_{s,ij}=&N_{\Si_{*,\frac{7c_*}{8}}}(\widetilde{\phi}-\Pi_2[\widetilde{\phi}^{(p)}_{s}])_{ij}\quad\textrm{on}\quad\Si_{*,\frac{7c_*}{8}}(\tau_1, \tau_2-1),\\
 (\widetilde{\phi}-\breve{\phi})^{(p)}_{s,ij} =& 0\quad\textrm{on}\quad\Si_{r\geq\varsigma_*(\tau_1-\frac{7c_*}{8})}(\tau_1),\\
  N_{\Si(\tau_1)}(\widetilde{\phi}-\breve{\phi})^{(p)}_{s,ij}=&0\quad\textrm{on}\quad\Si_{r\geq\varsigma_*(\tau_1-\frac{7c_*}{8})}(\tau_1),\\
\end{split}
\eeaa
Thus, applying the local energy estimates \eqref{eq:localenergyestimate:future:unweigthedderivatives:RR*} on $\MM_{r\geq\varsigma_*(\tau-\frac{7c_*}{8})}(\tau_1, \tau_2-1)$ and \eqref{eq:localenergyestimate:future:unweigthedderivatives:bis} on $\MM(\tau_2-1, \tau_2)$, we infer, for $\reg\leq {14}$, 
\beaa
\bsplit
&\sum_{i,j=1}^3\sum_{p=0}^2\bigg(\EM_{0,r\geq\varsigma_*(\tau-7c_*/8)}[\pr^{\leq\reg}(\widetilde{\phi}-\breve{\phi})](\tau_1, \tau_2-1)+\sup_{\frac{7c_*}{8}\leq c\leq c_*}\F_{\Si_{*,c}}[\pr^{\leq\reg}(\widetilde{\phi}-\breve{\phi})](\tau_1, \tau_2-1)\\
&+\EMF_\de[\pr^{\leq \reg}(\widetilde{\phi}-\breve{\phi})^{(p)}_{s,ij}](\tau_2-1, \tau_2)\bigg)\\
 \les& \sum_{i,j=1}^3\sum_{p=0}^2\bigg(\F_{\Si_{*,\frac{7c_*}{8}}}[\pr^{\leq \reg}(\widetilde{\phi}^{(p)}_{s}-\Pi_2[\widetilde{\phi}^{(p)}_{s}])_{ij}](\tau_1, \tau_2-1)+\E[\pr^{\leq \reg}(\widetilde{\phi}^{(p)}_{s}-\Pi_2[\widetilde{\phi}^{(p)}_{s}])_{ij}](\tau_2-1)\bigg).
\end{split}
\eeaa
 Together with Lemma \ref{lemma:computationerrorscalarizationdeffect}, we deduce, for 
  $\reg\leq {14}$, 
\beaa
\bsplit
&\sum_{i,j=1}^3\sum_{p=0}^2\bigg(\EM_{0,r\geq\varsigma_*(\tau-7c_*/8)}[\pr^{\leq\reg}(\widetilde{\phi}-\breve{\phi})](\tau_1, \tau_2-1)+\sup_{\frac{7c_*}{8}\leq c\leq c_*}\F_{\Si_{*,c}}[\pr^{\leq\reg}(\widetilde{\phi}-\breve{\phi})](\tau_1, \tau_2-1)\\
&+\EMF_\de[\pr^{\leq \reg}(\widetilde{\phi}-\breve{\phi})^{(p)}_{s,ij}](\tau_2-1, \tau_2)\bigg)\\
 \les& \sum_{i,j=1}^3\sum_{p=0}^2\bigg(\F_{\Si_{*,\frac{7c_*}{8}}}[\pr^{\leq \reg}\err_{\textrm{TDefect}}[\widetilde{\phi}^{(p)}_{s}]](\tau_1, \tau_2-1)+\E[\pr^{\leq \reg}\err_{\textrm{TDefect}}[\widetilde{\phi}^{(p)}_{s}]](\tau_2-1)\bigg).
\end{split}
\eeaa
Plugging \eqref{eq:controlscalarizationdefectforwidetildephiontau2minus2tau2}, and using \eqref{eq:bycausalityanddefintionwidetildephipsij=0inMMtau1tau1plus1}, we infer, for $\reg\leq {14}$, 
\begin{align}\lab{eq:controlofwidetildephiminusbreevephiontau2minus1tau2}
\nn&\sum_{i,j=1}^3\sum_{p=0}^2\bigg(\EM_{0,r\geq\varsigma_*(\tau-7c_*/8)}[\pr^{\leq\reg}(\widetilde{\phi}-\breve{\phi})](\tau_1, \tau_2-1)+\sup_{\frac{7c_*}{8}\leq c\leq c_*}\F_{\Si_{*,c}}[\pr^{\leq\reg}(\widetilde{\phi}-\breve{\phi})](\tau_1, \tau_2-1)\\
\nn&+\EMF_\de[\pr^{\leq \reg}(\widetilde{\phi}-\breve{\phi})^{(p)}_{s,ij}](\tau_2-1, \tau_2)\bigg)\\
\nn\les& \ep^2\sum_{i,j=1}^3\sum_{p=0}^2\bigg(\E[\pr^{\leq \reg}\phi^{(p)}_{s,ij}](\tau_1)+\E[\pr^{\leq \reg}\phi^{(p)}_{s,ij}](\tau_2-3)+\F_{\Si_{*,c_*/2}}[\pr^{\leq\reg}\phi^{(p)}_{s,ij}](\tau_1+1, \tau_2-1)\\
&+\int_{\MM(\tau_1, \tau_1+1)\cup\RR_*(\tau_1+1, \tau_2-3)\cup\MM(\tau_2-3, \tau_2-2)}r^{1+\de}|\pr^{\leq \reg}N_{W,s,ij}^{(p)}|^2\bigg).
\end{align}

\noindent{\bf Step 9.} Finally, we have obtained the following:
\begin{itemize}
\item in view of \eqref{eq:extendedmetricgchitau1tau2} and Step 1, $\g_{\tau_1, \tau_2,*}$ satisfies \eqref{eq:controloflinearizedinversemetriccoefficients} \eqref{eq:controloflinearizedinversemetriccoefficients:inversegtautau} and coincides with Kerr in $\MMh\setminus(\MM_{r\leq\varsigma_*(\tau-\frac{7c_*}{8})}(\tau_1, \tau_2-1))$, 

\item in view of Step 5 and the fact that $F_{total,s,ij}^{(p)}=\widehat{F}^{(p)}_{s,ij}+\widetilde{F}^{(p)}_{s,ij}+\underline{F}^{(p)}_{s,ij}+\breve{F}^{(p)}_{s,ij}$, $\psi^{(p)}_{s,ij}$ satisfies \eqref{eq:waveeqwidetildepsi1} \eqref{eq:causlityrelationsforwidetildepsi1},

\item in view of Step 6, $\psi^{(p)}_{s,ij}$ satisfies \eqref{eq:localenergyforpsiontau1minus1tau1plus1},

\item in view of Step 7, $\psi^{(p)}_{s,ij}$ satisfies \eqref{eq:localenergyestimateforpsisijponSigmatau2minus3tau2:partialderivatives},

\item and in view of Step 8, $\psi^{(p)}_{s,ij}$ satisfies \eqref{eq:propertyofscalarizationdefect:prop1} and \eqref{eq:propertyofscalarizationdefect:prop2}.
\end{itemize}
This concludes the proof of Proposition \ref{prop:extensionprocedureoftheTeukolskywaveequations}.
\end{proof}

%%%%%%%%%%%%%%%%%%%%%%%%%%%%%%%%%%%%%%%%%%%%%%%%%%%%%%%%%%%%%

\subsubsection{Remaining main intermediary results for the proof of Theorem \ref{thm:main:MaSz26}}

%%%%%%%%%%%%%%%%%%%%%%%%%%%%%%%%%%%%%%%%%%%%%%%%%%%%%%%%%%%%%

%%%%%%%%%%%%%%%%%%%%%%%%%%%%%%%%%%%%%%%%%%%%%%%%%%%%%%%%%%%%%%%%%%

\paragraph{\textit{Global energy-Morawetz estimates for unweighted derivatives of solutions to \eqref{eq:waveeqwidetildepsi1}}.}

%%%%%%%%%%%%%%%%%%%%%%%%%%%%%%%%%%%%%%%%%%%%%%%%%%%%%%%%%%%%%%%%%%
 
In order to prove Theorem \ref{thm:main:MaSz26}, i.e., the derivation of energy-Morawetz estimates for $\tau$ in $(\tau_1, \tau_2)$, we first state in this section global energy-Morawetz estimates for \eqref{eq:waveeqwidetildepsi1}, i.e., energy-Morawetz estimates for $\tau$ in $\mathbb{R}$, that hold for unweighted derivatives $\pr$ introduced in \eqref{eq:defunweightedderivative}.

\begin{theorem}\lab{th:main:intermediary}
Let $(\MM, \g)$ satisfy the assumptions of Section \ref{sec:assumptionsforsec:energyMorawetzesitmatesforTeukoslkyonMM:upto15derivatives}, and assume that $\tt_1$ and $\tt_2$ satisfy \eqref{eq:tau2geqtau1plus10conditionisthemaincase}. Let $\pmb\phi_s^{(p)}$, $s=\pm 2$, $p=0,1,2$, be a solution to the tensorial Teukolsky wave/transport systems \eqref{eq:TensorialTeuSysandlinearterms:rescaleRHScontaine2:general:Kerrperturbation}  \eqref{def:TensorialTeuScalars:wavesystem:Kerrperturbation}  in perturbations of Kerr, 
and let $\psis{p}$, $s=\pm 2$, $p=0,1,2$, be a solution to \eqref{eq:waveeqwidetildepsi1} satisfying \eqref{eq:causlityrelationsforwidetildepsi1}. Then, we have, for any $\reg\leq 14$ and $0<\de\leq \frac{1}{3}$,
\begin{align}
\lab{eq:mainhighorderunweightedEMF:maintheorem:pm2}
&\sum_{p=0}^2\EMF_{\de}[\pr^{\leq \reg}\phis{p}](\tau_1,\tau_2)
+ \sum_{p=0}^2{\EMFh}_{\de}[\pr^{\leq \reg}\psis{p}](\Iti)\nn\\
\les&\sum_{p=0}^2\E[\pr^{\leq \reg}\phis{p}](\tau_1)+\sum_{p=0}^2\widetilde{\NN}_\de[\pr^{\leq \reg}\phis{p}, \pr^{\leq \reg}\N_{W,s}^{(p)}](\tt_1, \tt_2)+\sum_{p=0}^1\int_{\MM(\tau_1,\tau_2)}r^{-3+\de}|\pr^{\leq \reg+1}\N_{T,s}^{(p)}|^2\nn\\
&
+\sum_{p=0}^1\int_{\MM_{r\geq 10m}(\tau_1,\tau_2)} \Big(r^{-1+\de}|\pr^{\leq \reg}\N^{(p)}_{W,s}|+r^{-2+\de}|\pr^{\leq \reg+1}\N^{(p)}_{T,s}|\Big)|\pr^{\leq \reg}\phis{p}|\nn\\
& +\sum_{p=0}^2\left(\EMFh[\psis{p}](\Iti)+\sup_{\frac{c_*}{2}\leq c\leq c_*}\F_{\Si_{*,c}}[\psis{p}](\tau_1+1, \tau_2-2)
+\EMF[\phis{p}](\tau_1, \tau_2)\right).
\end{align}
\end{theorem}

\begin{remark}
Theorem \ref{th:main:intermediary} is the analog of Theorem 7.6 in \cite{MaSz26}. Note the following two differences with the statement of Theorem 7.6 in \cite{MaSz26}:
\begin{itemize}
\item The norms for $\phis{p}$ and $\psis{p}$ in \eqref{eq:mainhighorderunweightedEMF:maintheorem:pm2} are respectively on $\MM$ and $\MMh$, while they are both on $\MM$ in Theorem 7.6 in \cite{MaSz26}.

\item The presence of the additional before to last terms involving the norms $\F_{\Si_{*,c}}[\psis{p}]$ on the RHS of \eqref{eq:mainhighorderunweightedEMF:maintheorem:pm2} compared to the corresponding estimate in Theorem 7.6 of  \cite{MaSz26}.
\end{itemize}
\end{remark}
 
The proof of Theorem \ref{th:main:intermediary} is postponed to Section \ref{sec:proofofth:main:intermediary:00}.

%%%%%%%%%%%%%%%%%%%%%%%%%%%%%%%%%%%%%%%%%%%%

\paragraph{\textit{Energy-Morawetz estimates near infinity in perturbations of Kerr}.}

%%%%%%%%%%%%%%%%%%%%%%%%%%%%%%%%%%%%%%%%%%%%

The following proposition provides energy-Morawetz estimates for Teukolsky equations in perturbations of Kerr for $r\geq R$ with $R$ large enough.

\begin{proposition}\lab{prop:EnergyMorawetznearinfinitytensorialTeuk}
Let $(\MM, \g)$ satisfy the assumptions of Section \ref{sec:assumptionsforsec:energyMorawetzesitmatesforTeukoslkyonMM:upto15derivatives}.  We have for solutions $\pmb\phi_s^{(p)}$, $s=\pm 2$, $p=0,1,2$, to the tensorial Teukolsky wave/transport systems \eqref{eq:TensorialTeuSysandlinearterms:rescaleRHScontaine2:general:Kerrperturbation}  \eqref{def:TensorialTeuScalars:wavesystem:Kerrperturbation}  in perturbations of Kerr the following energy-Morawetz-flux estimates, for any $1\leq\tau_1<\tau_2 <+\infty$ and any $0<\de\leq \frac{1}{3}$, for $\reg\leq 14$, and for $20m\leq R\leq c_*/3$ large enough, 
\bea
\lab{eq:EMnearinfinity:highorderweightedderivatives:Teu:pm2}
\nn&&\sum_{p=0}^2\EMF^{(\reg)}_{\de, \geq R}[\pmb\phi_{s}^{(p)}](\tau_1, \tau_2) +\sum_{p=0}^1\int_{\MM_{r\geq R}(\tau_1,\tau_2)}r^{-3+\de} \big|(r\nab)^{\leq 1} \dk^{\leq \reg}\phis{p}\big|^2\\
\nn&&+\sup_{\frac{c_*}{2}\leq c\leq c_*}\F^{(\reg)}_{\Si_{*,c}}[\phis{p}](\tau_1, \tau_2)\\
\nn&\les_R& \sum_{p=0}^2\M^{(\reg)}_{R/2, R}[\pmb\phi_{s}^{(p)}](\tau_1, \tau_2)+\sum_{p=0}^2\E^{(\reg)}_{r\geq R/2}[\pmb\phi_{s}^{(p)}](\tau_1)+\sum_{p=0}^2\int_{\MM_{r\geq R/2}(\tt_1, \tt_2)}r^{1+\de}|\dk^{\leq \reg}\N_{W,s}^{(p)}|^2\\
&&+\sum_{p=0}^1\int_{\MM_{r\geq R/2}(\tt_1, \tt_2)}r^{-1+\de}|\dk^{\leq \reg+1} \N^{(p)}_{T,s}|^2.
\eea
\end{proposition}

\begin{proof}
Proposition \ref{prop:EnergyMorawetznearinfinitytensorialTeuk} is the analog of Proposition 7.8 in \cite{MaSz26} whose proof in section 11 of \cite{MaSz26} immediately extends to the case where $\MM$ has the spacelike hypersurface $\Si_*$, instead of $\II_+$, as part of its future boundary.
\end{proof}

%%%%%%%%%%%%%%%%%%%%%%%%%%%%%%%%%%%%%%%%%%%%%%%%%%

\paragraph{\textit{Energy-Morawetz estimates for tensorial wave equations in subextremal Kerr}.}

%%%%%%%%%%%%%%%%%%%%%%%%%%%%%%%%%%%%%%%%%%%%%%%%%%

In order to control the lower order terms appearing last on the RHS of \eqref{eq:mainhighorderunweightedEMF:maintheorem:pm2}, we will rely on the two energy-Morawetz estimates in Kerr stated below. First, we consider solutions $\pmb\phi_s\in\sk_2(\mathbb{C})$, $s=\pm 2$, to the following tensorial wave equations on a subextremal Kerr background 
\bea\lab{eq:TeukolskyequationforAandAbintensorialforminKerr:inhomogenouscase}
\nn\left(\squared_2 -\frac{4ia\cos\th}{|q|^2}\nab_{\pr_t} -  \frac{s}{|q|^2}\right)\pmb\phi_s +\frac{2s}{|q|^2}(r-m)\nab_3\pmb\phi_s  -\frac{4sr}{|q|^2}\nab_{\pr_t}\pmb\phi_s\\
+\frac{4a\cos\th}{|q|^6}\Big(a\cos\th\big(|q|^2+6mr\big) - is\big((r-m)|q|^2+4mr^2\big)\Big)\pmb\phi_s &=& \N_s,
\eea
with $\N_s\in\sk_2(\mathbb{C})$, $s=\pm 2$, where \eqref{eq:TeukolskyequationforAandAbintensorialforminKerr:inhomogenouscase} is the inhomogeneous version of the tensorial Teukolsky equations in Kerr.

\begin{theorem}[Weak-Morawetz for Teukolsky in subextremal Kerr]
\lab{cor:weakMorawetzforTeukolskyfromMillet:bis}
Let $\tau_0\geq 1$ and assume that $\pmb\phi_s\in\sk_2(\mathbb{C})$, $s=\pm 2$, satisfies the inhomogeneous Teukosky equation 
\eqref{eq:TeukolskyequationforAandAbintensorialforminKerr:inhomogenouscase}  on a subextremal Kerr background.  Then, for any $\de\in (0,\frac{1}{3}]$, we have
\bea\lab{eq:weakMorawetzforTeukolskyfromMillet:bis:minus2}
\nn\int_{\MMh(\tau\geq\tau_0)}r^{-3+\de}|\dk^{\leq 3}\pmb\phi_{-2}|^2 &\les& \Eh^{(13)}[\pmb\phi_{-2}](\tau_0)+\Eh^{(11)}[r^{\frac{1+\de}{2}}\nab_{\pr_r}(r\pmb\phi_{-2})](\tau_0)\\
&&+\int_{\MMh(\tau\geq\tau_0)}r^{1+\de}|\dk^{\leq 13}\N_{-2}|^2,
\eea
and, for any $\de\in(0,\frac{1}{3}]$, any $10m\leq R_0\leq r_*/6$ and any $\tau_0+1<\tau_1\leq\tau_*$, we have
\bea\lab{eq:weakMorawetzforTeukolskyfromMillet:plus2:nolossinr}
&&\int_{(\MM\setminus\RR_*)(\tau_0,\tau_1-1)}r^{-11+\frac{\de}{2}}|\dk^{\leq 3}\pmb\phi_{+2}|^2\nn\\
\nn &\les&R_0^{1+\frac{\de}{2}}\E_{r\leq 2R_0}^{(9)}[r^{-4}\pmb\phi_{+2}](\tau_0) +\int_{\MM(\tau_0,\tau_1)}r^{-7+\de}|\dk^{\leq 10}\N_{+2}|^2 \\
\nn&&+R_0^{-\frac{\de}{2}}\bigg({\EMF}^{(11)}_{\de, r\leq\varsigma_*(\tau-c_*/2)}[r^{-4}\pmb\phi_{+2}](\tau_0,\tau_1)+\F^{(10)}_{\Si_{*,c_*/2}}[r^{-4}\pmb\phi_{+2}](\tau_0,\tau_1)\\
&&\qquad\qquad+\int_{(\MM\setminus\RR_*)(\tau_0,\tau_1)} r^{-3+\de}|\dk^{\leq 11} \nab_{3}(r^{-3}\pmb\phi_{+2})|^2\bigg).
\eea
\end{theorem}

\begin{proof}
The estimate \eqref{eq:weakMorawetzforTeukolskyfromMillet:bis:minus2} is proved in Theorem 7.9 in \cite{MaSz26}, while the proof of the estimate \eqref{eq:weakMorawetzforTeukolskyfromMillet:plus2:nolossinr} proceeds by a slight modification of the proof of the corresponding estimate in Theorem 7.9 in \cite{MaSz26}. Indeed, we have from (12.33) in \cite{MaSz26}, for any $10m\leq R\leq r_*/6$,
\bea
\lab{eq:weakmorawetz:TeuKerr:compactregion:+2}
\int_{\MM_{r\leq R}(\tau_0,\tau_1-1)}r^{-11+\frac{\de}{2}}|\dk^{\leq 3}\pmb\phi_{+2}|^2 
&\les& R^{1+\frac{\de}{2}}\E_{r\leq 2R}^{(9)}[r^{-4}\pmb\phi_{+2}](\tau_0) +\int_{\MM_{r\leq 2R}(\tau_0,\tau_1)}r^{-7+\frac{\de}{2}}|\dk^{\leq 9}\N_{+2}|^2 \nn\\
&&+\int_{\MM_{R, 2R}(\tau_0,\tau_1)}r^{-3+\frac{\de}{2}}|\dk^{\leq 9}\nab_3(r^{-3}\pmb\phi_{+2})|^2\nn\\
&&+\int_{\MM_{R, 2R}(\tau_0,\tau_1)}r^{-11+\frac{\de}{2}}|\dk^{\leq 10}\pmb\phi_{+2}|^2.
\eea
Next, we control the last term on the RHS of  \eqref{eq:weakmorawetz:TeuKerr:compactregion:+2}. Recall the identities (12.34) and (12.35) in \cite{MaSz26}, i.e.,
\bea
\lab{choiceoftensorsandsource:weakMora:Kerr:pf}
\phiplus{0}=|q|^{-4}\pmb\phi_{+2}, \quad \phiplus{1}=\frac{r^2}{|q|^2}\bigg({\frac{q}{\bar{q}}}r \nab_{3} r{\frac{\bar{q}}{q}}\bigg)\bigg(\frac{r^2}{|q|^2}\bigg)^{-2}\pmb\phi_{+2}^{(0)}, \quad \N_{T,+2}^{(0)}=0, \quad \N_{W,+2}^{(0)}=|q|^{-4}\N_{+2},
\eea
and
\bea
\lab{eq:Laplacianofphiprewrittenbyderivativeofphip+1:pm2:weightedderi:inKerr}
&&\qs \De_2\dk^{(\reg)}\phiplus{0}- 2\dk^{(\reg)}\phiplus{0}\nn\\
&=& O(r^{-1})(r\nab_4)^{\leq 1}\dk^{\leq \reg}\phiplus{1}
+O(1)\big(\nab_3\dk^{\leq \reg}\phiplus{0}, r^{-1}\dk^{\leq \reg+1}\phiplus{0}\big)
+\de_{\reg\geq 1}O(1)(r\nab)^{\leq 1}\dk^{\leq \reg-1}\phiplus{0}
\nn\\
&&+ O(r^2)\dk^{\leq \reg}\L_{+2}^{(0)}[\pmb\phi_{+2}] + O(r^2)\dk^{\leq \reg}\N_{W,+2}^{(0)}.
\eea 
Multiplying on both sides of \eqref{eq:Laplacianofphiprewrittenbyderivativeofphip+1:pm2:weightedderi:inKerr} by $r^{-3+\de}\ov{\dk^{(\reg)}\phiplus{0}}$, taking the real part, integrating over $(\MM\setminus\RR_*)(\tau_0,\tau_1)$, and in view of the expression of $\L_{+2}^{(0)}[\pmb\phi_{+2}]$ in \eqref{eq:tensor:Lsn:onlye_2present:general:Kerrperturbation:bis}, we deduce, for $\reg\leq 10$ and $0<\de\leq \frac{1}{3}$, 
\beaa
&&\int_{(\MM\setminus\RR_*)(\tau_0,\tau_1)}r^{-3+\de}|(r\nab)^{\leq 1}\dk^{\leq \reg}\phiplus{0}|^2\\
&\les&
\int_{(\MM\setminus\RR_*)(\tau_0,\tau_1)}\Big(r^{-1+\de}|\dk^{\leq \reg} \N_{W,{+2}}^{(0)}|+r^{-4+\de}|\dk^{\leq \reg+1}\phiplus{1}|\Big)|\dk^{\leq \reg} \phiplus{0}| \nn\\
&&
+\EMF^{(\reg+1)}_{\de,r\leq\varsigma_*(\tau-c_*/2)}[\phiplus{0}](\tau_0,\tau_1) +\F^{(\reg)}_{\Si_{*,c_*/2}}[r^{-4}\pmb\phi_{+2}](\tau_0,\tau_1)\\
&&+\de_{\reg\geq 1}\int_{(\MM\setminus\RR_*)(\tau_0,\tau_1)}r^{-3+\de}|(r\nab)^{\leq 1}\dk^{\leq \reg-1}\phiplus{0}|^2.
\eeaa
Summing over $\reg\leq 10$, and in view of \eqref{choiceoftensorsandsource:weakMora:Kerr:pf}, this yields the following analog of (12.36) in \cite{MaSz26}
\bea\lab{eq:weakmorawetz:TeuKerr:tocontrolcompactregion:+2}
\nn&&\int_{(\MM\setminus\RR_*)(\tau_0,\tau_1)} r^{-11+\de}|\dk^{\leq 10}\pmb\phi_{+2}|^2\\
&\les&{\EMF}^{(11)}_{\de,r\leq\varsigma_*(\tau-c_*/2)}[r^{-4}\pmb\phi_{+2}](\tau_0,\tau_1)+\F^{(10)}_{\Si_{*,c_*/2}}[r^{-4}\pmb\phi_{+2}](\tau_0,\tau_1)
+\int_{\MM(\tau_0,\tau_1)}r^{-7+\de}|\dk^{\leq 10} \N_{{+2}}|^2\nn\\
&&+ \int_{(\MM\setminus\RR_*)(\tau_0,\tau_1)} r^{-3+\de}|\dk^{\leq 11} \nab_{3}(r^{-3}\pmb\phi_{+2})|^2 .
\eea
Finally, combining the estimates\footnote{Note also from \eqref{eq:controlofsizeofc*assimeqr*:onSi*},  \eqref{eq:rsimcuptoOep0inSi*cforcbetweenc*on2andc*}, \eqref{eq:defdecreasingfunctionvarsigma*definingSi*c} and the choice of $R$ that 
\beaa
2R\leq \frac{r_*}{3}\leq\frac{c_*}{3}(1+O(\ep_0))\leq \varsigma_*(\tau-c_*/2), \quad \forall\tau\in[1,\tau_*]\quad\Longrightarrow\quad \MM_{R, 2R}(\tau_0,\tau_1)\subset(\MM\setminus\RR_*)(\tau_0,\tau_1).
\eeaa} \eqref{eq:weakmorawetz:TeuKerr:compactregion:+2} and  \eqref{eq:weakmorawetz:TeuKerr:tocontrolcompactregion:+2} yields the desired estimate \eqref{eq:weakMorawetzforTeukolskyfromMillet:plus2:nolossinr}. This concludes the proof of Theorem  \ref{cor:weakMorawetzforTeukolskyfromMillet:bis}.
\end{proof}

Next, we also state an energy-Morawetz estimates for $\breve{\pmb\phi}\in\sk_2(\mathbb{C})$ satisfying the following tensorial wave equation in a subextremal Kerr background
\bea\lab{eq:basictensorialwaveequationinKerr:forbrevephi}
\squared_2\breve{\pmb\phi} -\frac{4ia\cos\th}{|q|^2}\nab_{\pr_t}\breve{\pmb\phi} -\bigg(\frac{4}{\qs}-\frac{4a^2\cos^2\th}{|q|^6}\Big(|q|^2+6mr\Big)\bigg)\breve{\pmb\phi}=0.
\eea

\begin{theorem}[Energy-Morawetz for \eqref{eq:basictensorialwaveequationinKerr:forbrevephi} in subextremal Kerr]
\lab{prop:weakMorawetzfortensorialwaveeqfromDRSR}
Let $\breve{\pmb\phi}\in\sk_2(\mathbb{C})$ be a solution in a subextremal Kerr background to the tensorial wave equation \eqref{eq:basictensorialwaveequationinKerr:forbrevephi} in $\MMh(\tau\geq \tau_0)$, 
where $\tau_0\geq 1$ is a constant. Then, we have 
\bea
\lab{eq:EMFestimates:tensorialwaveinsubextremalKerr:weakMorathm}
\EMFh[\breve{\pmb\phi}](\tau_0,+\infty) \les \Eh[\breve{\pmb\phi}](\tau_0).
\eea
Furthermore, for $1\leq\tau_0<\tau_1\leq\tau_*$, we have 
\bea
\lab{eq:EMFestimates:tensorialwaveinsubextremalKerr:weakMorathm:MMhsetminusMMregion}
\EMFh_{r\geq \varsigma_*(\tau-c_*)}[\breve{\pmb\phi}](\tau_0,\tau_1) \les \Eh_{r\geq \varsigma_*(\tau_0-c_*)}[\breve{\pmb\phi}](\tau_0)+\F_{\Si_*}[\breve{\pmb\phi}](\tau_0,\tau_1).
\eea
\end{theorem}

\begin{proof}
The proof of \eqref{eq:EMFestimates:tensorialwaveinsubextremalKerr:weakMorathm} is given in Theorem 7.10 in \cite{MaSz26}. Concerning  \eqref{eq:EMFestimates:tensorialwaveinsubextremalKerr:weakMorathm:MMhsetminusMMregion}, note that \eqref{eq:basictensorialwaveequationinKerr:forbrevephi} corresponds to the tensorial wave equation \eqref{eq:tensorialwaveRW:withprttderivative} on Kerr in the particular case 
\beaa
k=2, \qquad \N=0, \qquad V=\frac{4}{\qs}-\frac{4a^2\cos^2\th}{|q|^6}\Big(|q|^2+6mr\Big).
\eeaa
In particular, integrating the divergence identity of Lemma \ref{cor:modifiedcurrentsforprtandprvphi} on $(\MMh\setminus\MM)(\tau_0, \tau)$ for $1\leq\tau_0<\tau\leq\tau_1\leq\tau_*$ immediately yields, using also the fact that $r\geq r_*\gg m$ on $(\MMh\setminus\MM)(1,\tau_*)$, 
\bea\lab{eq:EMFestimates:tensorialwaveinsubextremalKerr:weakMorathm:MMhsetminusMMregion:onlyenergypart}
\EFh_{r\geq \varsigma_*(\tau-c_*)}[\breve{\pmb\phi}](\tau_0,\tau_1) \les \Eh_{r\geq \varsigma_*(\tau_0-c_*)}[\breve{\pmb\phi}](\tau_0)+\F_{\Si_*}[\breve{\pmb\phi}](\tau_0,\tau_1).
\eea
Then, we integrate the divergence identity \eqref{eq:DivofPPmu:tensor:RW:prop} on $(\MMh\setminus\MM)(\tau_0, \tau_1)$ with $1\leq\tau_0<\tau_1\leq\tau_*$ using the vectorfield $X_1$ and the scalar function $w_1$ given by  
\beaa
\lab{eq:choicesofX0w:Morawetzde=1:Kerrpert}
X_1=2\mu(1-mr^{-1})\pr_{r}^{\text{BL}},\qquad w_1=4\mu r^{-1}(1-mr^{-1}). 
\eeaa
Using the following estimate which holds for $r$ large enough, see below (8.152) in \cite{MaSz26}, 
\beaa
\bigg(\QQ[\breve{\pmb\phi}]  \c {^{(X_1)}}\pi +w_1 \LL[\breve{\pmb\phi}]  -X_1(V) |\breve{\pmb\phi}|^2-\frac{1}{2}|\breve{\pmb\phi}|^2   \square_\g  w_1\bigg) &\gtrsim& \frac{|\nab_{\pr_{\tt}}\breve{\pmb\phi}|^2}{r^2}
+\frac{|\nab_{\pr_{r}}\breve{\pmb\phi}|^2}{r^2} +\frac{|\nab\breve{\pmb\phi}|^2}{r}+\frac{|\breve{\pmb\phi}|^2}{r^3},
\eeaa
 together with Corollary \ref{cor:asymtpticbehavioroftheRdottermintensorialenergyidentity:larger}, we obtain, using also the fact that $r\geq r_*\gg m$ on $(\MMh\setminus\MM)(1,\tau_*)$,  
\beaa
\Mh_{r\geq \varsigma_*(\tau-c_*)}[\breve{\pmb\phi}](\tau_0,\tau_1) \les \EFh_{r\geq \varsigma_*(\tau-c_*)}[\breve{\pmb\phi}](\tau_0,\tau_1),
\eeaa
which together with \eqref{eq:EMFestimates:tensorialwaveinsubextremalKerr:weakMorathm:MMhsetminusMMregion:onlyenergypart} yields the stated estimate  \eqref{eq:EMFestimates:tensorialwaveinsubextremalKerr:weakMorathm:MMhsetminusMMregion}. This concludes the proof of Theorem \ref{prop:weakMorawetzfortensorialwaveeqfromDRSR}.
\end{proof}

%%%%%%%%%%%%%%%%%%%%%%%%%%%

\subsubsection{Proof of Theorem \ref{thm:main:MaSz26}}
\lab{subsect:proofofThm4.1}

%%%%%%%%%%%%%%%%%%%%%%%%%%%

Let $(\MM, \g)$ satisfy the assumptions of Section \ref{sec:assumptionsforsec:energyMorawetzesitmatesforTeukoslkyonMM:upto15derivatives}, let $\tau_1$ and $\tau_2$ be such that $1\leq\tau_1<\tau_2 <+\infty$, let $\phis{p}$, $s=\pm 2$, $p=0,1,2$, be solutions to the tensorial Teukolsky wave/transport systems \eqref{eq:TensorialTeuSysandlinearterms:rescaleRHScontaine2:general:Kerrperturbation}  \eqref{def:TensorialTeuScalars:wavesystem:Kerrperturbation}  in perturbations of Kerr, and  assume also that $\Ab$ satisfies \eqref{eq:waveequationpmbphip=0sminus2nodeginredshiftregion}. We proceed in the following steps, following closely section 7.3 in \cite{MaSz26}.

\noindent{\bf Step 0.} We first prove Theorem \ref{thm:main:MaSz26} in the case $\tau_1<\tau_2\leq \tau_1+10$. This follows easily by local existence and redshift estimates, see Step 0 in section 7.3 in \cite{MaSz26} which immediately applies here.

\noindent{\bf Step 1.} We now assume that \eqref{eq:tau2geqtau1plus10conditionisthemaincase} holds, i.e., we consider the case $\tau_2\geq\tau_1+10$, and we begin by upgrading the estimate \eqref{eq:mainhighorderunweightedEMF:maintheorem:pm2} to energy-Morawetz estimates controlling in addition weighted derivatives of $\pmb\phi_s^{(p)}$. Proceeding as in Step 1 in \cite{MaSz26}, which relies on combining the energy-Morawetz estimates \eqref{eq:mainhighorderunweightedEMF:maintheorem:pm2} and \eqref{eq:EMnearinfinity:highorderweightedderivatives:Teu:pm2}, we obtain, for $s=\pm 2$, and for any $\reg\leq 14$ and $0<\de\leq \frac{1}{3}$,
\begin{align}\lab{eq:allderivativesrecoveredandonlylowerordertermonRHS:plus2caseandminus2caseatthesametime}
&\sum_{p=0}^2\EMF_{\de}^{(\reg)}[\phis{p}](\tau_1,\tau_2) + \sum_{p=0}^2{\EMFh}_{\de}[\pr^{\leq \reg}\psis{p}](\Iti)\nn\\
\les&\sum_{p=0}^2\E^{(\reg)}[\phis{p}](\tau_1)+\sum_{p=0}^2\NN_\de^{(\reg)}[\phis{p}, \N_{W,s}^{(p)}](\tt_1, \tt_2)+\sum_{p=0}^1\int_{\MM(\tau_1,\tau_2)}r^{-1+\de}|\dk^{\leq \reg+1}\N_{T,s}^{(p)}|^2\nn\\
& +\sum_{p=0}^2\bigg(\EMFh[\psis{p}](\Iti) +\sup_{\frac{c_*}{2}\leq c\leq c_*}\F_{\Si_{*,c}}[\psis{p}](\tau_1+1, \tau_2-2)
+\EMF[\phis{p}](\tau_1, \tau_2)\bigg).
\end{align}

\noindent{\bf Step 2.} In this step, we remove the lower order terms which appear on the last line of the RHS of \eqref{eq:allderivativesrecoveredandonlylowerordertermonRHS:plus2caseandminus2caseatthesametime} by relying on the energy-Morawetz estimates in Kerr of Theorems \ref{cor:weakMorawetzforTeukolskyfromMillet:bis} and \ref{prop:weakMorawetzfortensorialwaveeqfromDRSR}. To this end, using the regular triplet $(\Om_K)_i$, $i=1,2,3$, in Kerr introduced in Definition \ref{def:regulartripletinKerrOmii=123}, we define the horizontal tensors $\pmb\phi_s$, $s=\pm 2$, on $r\leq\varsigma_*(\tau-c_*)$ by 
\bea\lab{eq:definitionofphiplusminus2inKerrfromphi0plusminus2Kerrpert}
(\pmb\phi_{+2})_{ab}=|q|^4\phi_{+2,ij}^{(0)}(\Om_K^i)_a(\Om_K^j)_b, \qquad (\pmb\phi_{-2})_{ab}=\frac{|q|^4}{\De^2}\phi_{-2,ij}^{(0)}(\Om_K^i)_a(\Om_K^j)_b,
\eea
where the definition of $\pmb\phi_{-2}$ will only be used in $r\geq r_+(1+2\dred )$. Proceeding as in the first part Step 2 of \cite{MaSz26}, we first obtain the following inhomogeneous Teukolsky equations in Kerr, for $s=\pm 2$ and $r\leq\varsigma_*(\tau-c_*)$, 
\bea\lab{eq:TeukolskyequationforAandAbintensorialforminKerr:inhomogeneous} 
\nn\left(\squared_{2,K} -\frac{4ia\cos\th}{|q|^2}\nab_{\pr_t} -  \frac{s}{|q|^2}\right)\pmb\phi_s +\frac{2s}{|q|^2}(r-m)\nab_{(e_K)_3}\pmb\phi_s  -\frac{4sr}{|q|^2}\nab_{\pr_t}\pmb\phi_s\\
+\frac{4a\cos\th}{|q|^6}\Big(a\cos\th\big(|q|^2+6mr\big) - is\big((r-m)|q|^2+4mr^2\big)\Big)\pmb\phi_s &=& \N_s,  
\eea
where the inhomogeneous RHS $\N_s$, $s=\pm 2$, are given by  
\bsub\lab{eq:structureofinhomogenoustermsNplusminus2inTeukolskyforphisinKerr}
\bea
\N_{+2} &=& |q|^4\Big[\N_{W,+2}^{(0)}+O(r^{-2})\N_{T,+2}^{(0)}+\dk^{\leq 2}(\Ga_g\c\pmb\phi_{+2}^{(0)})\Big],\\
\N_{-2} &=& \frac{|q|^4}{\De^2}\Big[\N_{W,-2}^{(0)}+O(\De^{-1})\N_{T,-2}^{(0)}+\frac{r^2}{\De}\dk^{\leq 2}(\Ga_g\c\pmb\phi_{-2}^{(0)})\Big], 
\eea
\esub
with the definition for $\N_{-2}$ in \eqref{eq:structureofinhomogenoustermsNplusminus2inTeukolskyforphisinKerr} being used only for $r\geq r_+(1+2\dred)$, see (7.57) (7.58) in \cite{MaSz26}. For, $r\leq r_+(1+2\dred)$, similarly to \eqref{eq:definitionofphiplusminus2inKerrfromphi0plusminus2Kerrpert}, we define $\pmb\phi_{-2}\in\sk_2(\mathbb{C})$ in Kerr for $r\leq r_+(1+2\dred)$ by\footnote{The definitions of $\pmb\phi_{-2}$ in \eqref{eq:definitionofphiplusminus2inKerrfromphi0plusminus2Kerrpert} and \eqref{eq:definitionofphiplusminus2inKerrfromphi0plusminus2Kerrpert:bis} agree since $
\phi_{-2,ij}^{(0)}=q(\ov{q})^{-1}\De^2|q|^{-4}\Ab(\Om_i, \Om_j)$ in view of \eqref{eq:definitionofthephiminus2phierarchy:perturbationofKerr}.}
\bea\lab{eq:definitionofphiplusminus2inKerrfromphi0plusminus2Kerrpert:bis}
(\pmb\phi_{-2})_{ab}=|q|^{-2}q^2\Ab(\Om_i, \Om_j)(\Om_K^i)_a(\Om_K^j)_b,
\eea
and proceeding again as in the first part Step 2 of \cite{MaSz26}, see (5.102) in \cite{MaSz26}, we have
\bea\lab{eq:structureofinhomogenoustermsNplusminus2inTeukolskyforphisinKerr:HHm2}
\N_{-2}=|q|^{-2}q^2\Big[\N_{\Ab}+\dk^{\leq 2}(\Ga_g\c\Ab)\Big].
\eea

Next, we introduce $\pmb\varphi_s\in\sk_2(\mathbb{C})$, $s=\pm 2$, in Kerr which are the solutions of the following modification of \eqref{eq:TeukolskyequationforAandAbintensorialforminKerr:inhomogeneous}
\bea\lab{eq:TeukolskyequationforAandAbintensorialforminKerr:inhomogeneous:tilde} 
\nn\left(\squared_{2,K} -\frac{4ia\cos\th}{|q|^2}\nab_{\pr_t} -  \frac{s}{|q|^2}\right)\pmb\varphi_s +\frac{2s}{|q|^2}(r-m)\nab_{(e_K)_3}\pmb\varphi_s  -\frac{4sr}{|q|^2}\nab_{\pr_t}\pmb\varphi_s\\
+\frac{4a\cos\th}{|q|^6}\Big(a\cos\th\big(|q|^2+6mr\big) - is\big((r-m)|q|^2+4mr^2\big)\Big)\pmb\varphi_s &=& \widetilde{\N}_s,  
\eea
with 
\bea\lab{eq:definitiionofwidetildeNsextendingNsinKerr}
\widetilde{\N}_s=\widetilde\chi_{\tau_2}(\tau)\chi_{r_*}(\tau, r)\N_s, \qquad s=\pm 2,
\eea
where $\chi_{r_*}$ is the smooth cut-off function satisfying \eqref{eq:supportpropertiesofcutoffchistar}, and where $\widetilde\chi_{\tau_2}$ is a smooth cut-off function satisfying $\widetilde\chi_{\tau_2}(\tau)=1$ for $\tau\leq \tau_2-1$ and $\widetilde\chi_{\tau_2}=0$ for $\tau\geq \tau_2$, and with initial data 
\bea\lab{eq:choiceofinitialdataforpmbvarphisisequalchir*initialdatapmbphis}
\pmb\varphi_s=\chi_{r_*}(\tau_1, r_*)\pmb\phi_s, \qquad \nab_{\pr_{\tau}}\pmb\varphi_s=\chi_{r_*}(\tau_1, r_*)\nab_{\pr_{\tau}}\pmb\phi_s, \quad\textrm{on}\quad\{\tau=\tau_1\},
\eea
so that we have by causality
\bea\lab{eq:causalityrelationbetweenpmbphisandpmbpsisinKerr}
\pmb\varphi_s=\pmb\phi_s, \quad s=\pm 2, \quad \textrm{on}\quad (\MM\setminus\RR_*)(\tau_1, \tau_2-1).
\eea
Now, in view of \eqref{eq:TeukolskyequationforAandAbintensorialforminKerr:inhomogeneous:tilde}, $\pmb\varphi_s$ satisfies \eqref{eq:TeukolskyequationforAandAbintensorialforminKerr:inhomogenouscase}, and we may thus apply the weak Morawetz estimates in Kerr of Theorem \ref{cor:weakMorawetzforTeukolskyfromMillet:bis}, with $\tau_0=\tau_1$,  which implies, for any $\de\in(0,\frac{1}{3}]$,
\bea\lab{eq:weakMorawetzestimatesforphiminus2:consequenceThoerem}
\nn\int_{\MMh(\tau_1,+\infty)}r^{-3+\de}|\dk^{\leq 3}\pmb\varphi_{-2}|^2 &\les& \Eh^{({13})}[\pmb\varphi_{-2}](\tau_1)+\Eh^{({11})}[r\nab_{\pr_r}(r\pmb\varphi_{-2})](\tau_1)\\
&&+\int_{\MMh(\tau_1,+\infty)}r^{1+\de}|\dk^{\leq {13}}\widetilde{\N}_{-2}|^2,
\eea
and with $(\tau_0,\tau_1)=(\tau_1,\tau_2-1)$, which implies, for any $\de\in(0,\frac{1}{3}]$ and any $10m\leq R_0\leq r_*/6$, 
\begin{align}\lab{eq:weakMorawetzestimatesforphiplus2:consequenceThoerem}
&\int_{(\MM\setminus\RR_*)(\tau_1,\tau_2-2)}r^{-11+\frac{\de}{2}}|\dk^{\leq 3}\pmb\varphi_{+2}|^2\nn\\
 \les{}& R_0^{1+\frac{\de}{2}}\E^{({9})}[r^{-4}\pmb\varphi_{+2}](\tau_1)+\int_{\MM(\tau_1,\tau_2-1)}r^{-7+\de}|\dk^{\leq {10}}\widetilde{\N}_{+2}|^2 \nn\\
\nn&+R_0^{-\frac{\de}{2}}\bigg({\EMF}^{(11)}_{\de, r\leq\varsigma_*(\tau-c_*/2)}[r^{-4}\pmb\varphi_{+2}](\tau_1,\tau_2-1)+\F^{(10)}_{\Si_{*,c_*/2}}[r^{-4}\pmb\varphi_{+2}](\tau_1,\tau_2-1)\\
&\qquad\qquad+\int_{(\MM\setminus\RR_*)(\tau_1,\tau_2-1)} r^{-3+\de}|\dk^{\leq 11} \nab_{3,\K}(r^{-3}\pmb\varphi_{+2})|^2\bigg).
\end{align}

We first derive a  weak Morawetz estimate for $\phiminus{0}$. In view of  \eqref{eq:definitionofphiplusminus2inKerrfromphi0plusminus2Kerrpert}, \eqref{eq:definitionofphiplusminus2inKerrfromphi0plusminus2Kerrpert:bis} and \eqref{eq:causalityrelationbetweenpmbphisandpmbpsisinKerr}, we have
\beaa
\int_{(\MM\setminus\RR_*)(\tau_1, \tau_2-1)}r^{-3+\de}|\dk^{\leq 3}\pmb\phi_{-2}^{(0)}|^2 &\les&   \int_{(\MM\setminus\RR_*)(\tau_1, \tau_2-1)}r^{-3+\de}|\dk^{\leq 3}\pmb\phi_{-2}|^2\\
&\les&\int_{(\MM\setminus\RR_*)(\tau_1, \tau_2-1)}r^{-3+\de}|\dk^{\leq 3}\pmb\varphi_{-2}|^2
\eeaa
and, using in addition \eqref{eq:choiceofinitialdataforpmbvarphisisequalchir*initialdatapmbphis} and the properties \eqref{eq:supportpropertiesofcutoffchistar} of the cut-off $\chi_{r_*}$, 
\beaa
\Eh^{({13})}[\pmb\varphi_{-2}](\tau_1)+\Eh^{({11})}[r\nab_{\pr_r}(r\pmb\varphi_{-2})](\tau_1) &\les& \E^{({13})}[\pmb\phi_{-2}](\tau_1)+\E^{({11})}[r\nab_{\pr_r}(r\pmb\phi_{-2})](\tau_1)\\
&\les& \E^{({13})}[\pmb\phi_{-2}^{(0)}](\tau_1)+\E^{({11})}[r\nab_{\pr_r}(r\pmb\phi_{-2}^{(0)})](\tau_1)\\
&&+\E^{({13})}_{r\leq r_+(1+\dred )}[\Ab](\tau_1).
\eeaa
Also, we have in view of \eqref{eq:definitiionofwidetildeNsextendingNsinKerr} and the properties of the cut-off functions $\widetilde\chi_{\tau_2}$ and $\chi_{r_*}$
\beaa
\int_{\MMh(\tau_1,+\infty)}r^{1+\de}|\dk^{\leq {13}}\widetilde{\N}_{-2}|^2\les \int_{\MM(\tau_1,\tau_2)}r^{1+\de}|\dk^{\leq {13}}\N_{-2}|^2.
\eeaa 
Plugging the above estimates in \eqref{eq:weakMorawetzestimatesforphiminus2:consequenceThoerem}, we infer
\beaa
\nn\int_{(\MM\setminus\RR_*)(\tau_1, \tau_2-1)}r^{-3+\de}|\dk^{\leq 3}\pmb\phi_{-2}^{(0)}|^2 &\les& \E^{({13})}[\pmb\phi_{-2}^{(0)}](\tau_1)+\E^{({11})}[r\nab_{\pr_r}(r\pmb\phi_{-2}^{(0)})](\tau_1)\\
&&+\E^{({13})}_{r\leq r_+(1+\dred )}[\Ab](\tau_1)+\int_{\MM(\tau_1,\tau_2)}r^{1+\de}|\dk^{\leq {13}}\N_{-2}|^2.
\eeaa
Finally, together with \eqref{eq:structureofinhomogenoustermsNplusminus2inTeukolskyforphisinKerr} and \eqref{eq:structureofinhomogenoustermsNplusminus2inTeukolskyforphisinKerr:HHm2}, (5.9) in \cite{MaSz26} which states
\beaa
e_4 = \left(1+O(mr^{-1})\right)\pr_r + O(m^2r^{-2})\pr_\tau+O(\ep r^{-2})\pr_{x^a},
\eeaa
and the fact that, in view of \eqref{def:TensorialTeuScalars:wavesystem:Kerrperturbation},  
\beaa
\E^{({11})}[r\nab_{\pr_r}(r\pmb\phi_{-2}^{(0)})](\tau_1) &\les& \E^{({11})}[\pmb\phi_{-2}^{(1)}](\tau_1)+\E^{(12)}[\pmb\phi_{-2}^{(0)}](\tau_1)+\int_{\Si(\tau_1)}|\dk^{\leq {12}}\N_{T,-2}^{(0)}|^2,
\eeaa
we deduce the following weak Morawetz estimate for $\phiminus{0}$  in perturbations of Kerr\footnote{Again, recall, from our convention introduced in Section \ref{sec:smallnesconstants}, that we do not need to track the dependence of $\les$ on $\dred $ in the two first lines on the RHS of \eqref{eq:weakMorawetzestimatesforphiminus2:consequenceThoerem:1}. Additionally, recall that $\ep\ll\dred $ in view of \eqref{eq:constraintsonthemainsmallconstantsepanddelta} which will allow us, below, to absorb the nonlinear terms on the last line of \eqref{eq:weakMorawetzestimatesforphiminus2:consequenceThoerem:1} regardless of the dependence on $\dred $ so that we do not track this dependence in the last line either.} 
\bea\lab{eq:weakMorawetzestimatesforphiminus2:consequenceThoerem:1}
\nn&&\int_{(\MM\setminus\RR_*)(\tau_1, \tau_2-1)}r^{-3+\de}|\dk^{\leq 3}\pmb\phi_{-2}^{(0)}|^2\\ 
\nn&\les& \E^{({13})}[\pmb\phi_{-2}^{(0)}](\tau_1)+ \E^{({11})}[\pmb\phi_{-2}^{(1)}](\tau_1)+\int_{\Si(\tau_1)}|\dk^{\leq {12}}\N_{T,-2}^{(0)}|^2 +\E^{({13})}_{r\leq r_+(1+\dred )}[\Ab](\tau_1)\\
\nn&+&\int_{\MM(\tau_1,\tau_2)}r^{1+\de}|\dk^{\leq {13}}\N_{W,-2}^{(0)}|^2+\int_{\MM(\tau_1,\tau_2)}r^{-3+\de}|\dk^{\leq {13}}\N_{T,-2}^{(0)}|^2+\int_{\MM_{r\leq r_+(1+\dred )}(\tau_1,\tau_2)}|\dk^{\leq {13}}\N_{\Ab}|^2\\
&&+\ep^2\sup_{\tau\in[\tau_1, \tau_2]}\E^{({14})}[\pmb\phi_{-2}^{(0)}](\tau)+\ep^2\sup_{\tau\in[\tau_1, \tau_2]}\E^{({14})}_{r\leq r_+(1+\dred )}[\Ab](\tau).
\eea

Next, we derive a  weak Morawetz estimate for $\phiplus{0}$. In view of \eqref{eq:definitionofphiplusminus2inKerrfromphi0plusminus2Kerrpert}, \eqref{eq:definitiionofwidetildeNsextendingNsinKerr}, \eqref{eq:choiceofinitialdataforpmbvarphisisequalchir*initialdatapmbphis}, \eqref{eq:causalityrelationbetweenpmbphisandpmbpsisinKerr} and the properties of the cut-off functions $\widetilde\chi_{\tau_2}$ and $\chi_{r_*}$, the estimate \eqref{eq:weakMorawetzestimatesforphiplus2:consequenceThoerem} implies
\beaa
\int_{(\MM\setminus\RR_*)(\tau_1,\tau_2-2)}r^{-3+\frac{\de}{2}}|\dk^{\leq 3}\phiplus{0}|^2 &\les& R_0^{1+\frac{\de}{2}}\E^{({9})}[\phiplus{0}](\tau_1)+\int_{\MM(\tau_1,\tau_2-1)}r^{-7+\de}|\dk^{\leq {10}}{\N}_{+2}|^2 \nn\\
&&+R_0^{-\frac{\de}{2}}\bigg({\EMF}^{(11)}_{\de}[\phiplus{0}](\tau_1,\tau_2-1)+\F^{(10)}_{\Si_{*,c_*/2}}[\phiplus{0}](\tau_1,\tau_2-1)\\
&&\qquad\qquad+\int_{\MM(\tau_1,\tau_2-1)} r^{-3+\de}|\dk^{\leq 11} \nab_{3,\K}(r\phiplus{0})|^2\bigg).
\eeaa
Together with (5.11) in \cite{MaSz26} which implies
\beaa
e_{3,\K} =e_3 +r\Ga_b \pr_r + r\Ga_g \pr_{\tau} + \Ga_b \pr_{x^a}=e_3 + r\Ga_g \dk,
\eeaa
and the fact that, in view of \eqref{def:TensorialTeuScalars:wavesystem:Kerrperturbation},  
\beaa
\int_{\MM(\tau_1,\tau_2-1)} r^{-3+\de}|\dk^{\leq 11}\nab_3(r\phiplus{0})|^2\les 
\sum_{p=0}^1\int_{\MM(\tau_1,\tau_2-1)} r^{-5+\de}|\dk^{\leq 11} \phiplus{p}|^2+\int_{\MM(\tau_1,\tau_2-1)} r^{-3+\de}|\dk^{\leq 11} \N_{T, +2}^{(0)}|^2,
\eeaa
we deduce the following weak Morawetz estimates for $\phiplus{0}$  in perturbations of Kerr
\bea\lab{eq:weakMorawetzestimatesforphiplus2:consequenceThoerem:1}
&&\int_{(\MM\setminus\RR_*)(\tau_1,\tau_2-2)}r^{-3+\frac{\de}{2}}|\dk^{\leq 3}\phiplus{0}|^2\nn\\
 &\les& R_0^{1+\frac{\de}{2}}\E^{({9})}[\phiplus{0}](\tau_1)+\int_{\MM(\tau_1,\tau_2-1)}\Big(r^{-7+\de}|\dk^{\leq {10}}{\N}_{+2}|^2+r^{-3+\de}|\dk^{\leq 11} \N_{T, +2}^{(0)}|^2\Big) \nn\\
&&+R_0^{-\frac{\de}{2}}\bigg({\EMF}^{(11)}_{\de}[\phiplus{0}](\tau_1,\tau_2-1)+\F^{(10)}_{\Si_{*,c_*/2}}[\phiplus{0}](\tau_1,\tau_2-1)\nn\\
&&\qquad\qquad+\int_{\MM(\tau_1,\tau_2-1)} r^{-5+\de}|\dk^{\leq 11} \phiplus{1}|^2\bigg)\nn\\
&\les& R_0^{1+\frac{\de}{2}}\E^{({9})}[\phiplus{0}](\tau_1)+\int_{\MM(\tau_1,\tau_2-1)}\Big(r^{1+\de}|\dk^{\leq {10}}{\N}_{W,+2}^{(0)}|^2+r^{-3+\de}|\dk^{\leq 11} \N_{T, +2}^{(0)}|^2\Big) \nn\\
&&+\Big(R_0^{-\frac{\de}{2}}+\ep^2\Big){\EMF}^{(11)}_{\de}[\phiplus{0}](\tau_1,\tau_2-1)+R_0^{-\frac{\de}{2}}\F^{(10)}_{\Si_{*,c_*/2}}[\phiplus{0}](\tau_1,\tau_2-1)\nn\\
&&+R_0^{-\frac{\de}{2}}\int_{\MM(\tau_1,\tau_2-1)} r^{-5+\de}|\dk^{\leq 11} \phiplus{1}|^2,
\eea
where we have used \eqref{eq:structureofinhomogenoustermsNplusminus2inTeukolskyforphisinKerr} in the last step.

Next, notice that  we have, for any $R\leq c_*/3$,
\beaa
\sum_{p=0}^2\EMF_{r\leq R}[\pmb\phi_{+2}^{(p)}](\tau_1, \tau_2-2)\les_{R}\int_{(\MM\setminus\RR_*)(\tau_1, \tau_2-2)}r^{-3}|\dk^{\leq 3}\pmb\phi_{+2}^{(0)}|^2
+\sum_{p=0}^2\int_{\MM_{r\leq R+m}(\tau_1, \tau_2-2)}|\N_{W,+2}^{(p)}|^2,\\
\sum_{p=0}^2\EMF_{r\leq R}[\pmb\phi_{-2}^{(p)}](\tau_1, \tau_2-1)\les_{R} \int_{(\MM\setminus\RR_*)(\tau_1, \tau_2-1)}r^{-3}|\dk^{\leq 3}\pmb\phi_{-2}^{(0)}|^2+\sum_{p=0}^2\int_{\MM_{r\leq R+m}(\tau_1, \tau_2-1)}|\N_{W,-2}^{(p)}|^2,
\eeaa
where
\begin{itemize} 
\item the control of the Morawetz norm is immediate, 
\item the control of the flux norm on $\AA$ follows from redshift estimates and the control of Morawetz, 
\item and the control of the energy norm follows from the mean value argument which yields the control of the energy for $\tau=\tau_n$ with $\tau_n\in[n,n+1)$ for any $[n,n+1)\subset(\tau_1, \tau_2)$ and local energy estimates. 
\end{itemize}
In view of the above and \eqref{eq:weakMorawetzestimatesforphiminus2:consequenceThoerem:1} \eqref{eq:weakMorawetzestimatesforphiplus2:consequenceThoerem:1}, we may apply Proposition \ref{prop:EnergyMorawetznearinfinitytensorialTeuk} with $\reg=0$ to deduce
\beaa
\nn&&\sum_{p=0}^2\Big(\EMF[\pmb\phi_{+2}^{(p)}](\tau_1, \tau_2-1)+\F_{\Si{*,c_*/2}}[\pmb\phi_{+2}^{(p)}](\tau_1, \tau_2-1)\Big)\\ 
\nn&\les& R_0^{1+\de}\sum_{p=0}^2\E^{({9})}[\pmb\phi_{+2}^{(p)}](\tau_1)
+\sum_{p=0}^2\int_{\MM(\tau_1,\tau_2-1)}r^{1+\de}|\dk^{\leq {10}}{\N}_{W,+2}^{(p)}|^2\\
\nn&&+\sum_{p=0}^1\int_{\MM(\tau_1,\tau_2-1)}r^{-1+\de}|\dk^{\leq 11} \N_{T, +2}^{(p)}|^2
+\Big(R_0^{-\frac{\de}{2}}+\ep^2\Big)\sum_{p=0}^1{\EMF}^{(11)}_{\de}[\phiplus{p}](\tau_1,\tau_2-1)\\
\nn&&+R_0^{-\frac{\de}{2}}\F^{(10)}_{\Si_{*,c_*/2}}[\phiplus{0}](\tau_1,\tau_2-1),
\eeaa
where we have also used the local energy estimate \eqref{eq:localenergyestimate:future:bis} applied with $\reg=0$, $\tau_0=\tau_2-2$,  $q=1$ and $\psi=(\phiplus{p})_{p=0,1,2}$, and Proposition \ref{prop:EnergyMorawetznearinfinitytensorialTeuk} with $\reg=0$ on $\MM(\tau_2-2, \tau_2-1)$, and 
\beaa
\nn&&\sum_{p=0}^2\Big(\EMF[\pmb\phi_{-2}^{(p)}](\tau_1, \tau_2-1)+\F_{\Si{*,c_*/2}}[\pmb\phi_{-2}^{(p)}](\tau_1, \tau_2-1)\Big)\\ 
\nn&\les& \sum_{p=0}^2\E^{({13})}[\pmb\phi_{-2}^{(p)}](\tau_1)+\int_{\Si(\tau_1)}|\dk^{\leq {12}}\N_{T,-2}^{(0)}|^2 +\E^{({13})}_{r\leq r_+(1+\dred )}[\Ab](\tau_1)\\
\nn&+&\sum_{p=0}^2\int_{\MM(\tau_1,\tau_2)}r^{1+\de}|\dk^{\leq {13}}\N_{W,-2}^{(p)}|^2+\sum_{p=0}^1\int_{\MM(\tau_1,\tau_2)}r^{-1+\de}|\dk^{\leq {13}}\N_{T,-2}^{(p)}|^2\\
&+&\int_{\MM_{r\leq r_+(1+\dred )}(\tau_1,\tau_2)}|\dk^{\leq {13}}\N_{\Ab}|^2+\ep^2\sup_{\tau\in[\tau_1, \tau_2]}\E^{({14})}[\pmb\phi_{-2}^{(0)}](\tau)+\ep^2\sup_{\tau\in[\tau_1, \tau_2]}\E^{({14})}_{r\leq r_+(1+\dred )}[\Ab](\tau).
\eeaa
Together with \eqref{eq:causlityrelationsforwidetildepsi1}, this implies for the scalars $\psi^{(p)}_{s,ij}$
\beaa
\nn&&\sum_{p=0}^2\EMF_{r\leq\varsigma_*(\tau-c_*/2)}[\psiplus{p}](\tau_1+1, \tau_2-3)+
\sum_{p=0}^2\EMF[\pmb\phi_{+2}^{(p)}](\tau_1, \tau_2-1)\\
&& +\sum_{p=0}^2\F_{\Si{*,c_*/2}}[\pmb\phi_{+2}^{(p)}](\tau_1, \tau_2-1)\\ 
\nn&\les& R_0^{1+\de}\sum_{p=0}^2\E^{({9})}[\pmb\phi_{+2}^{(p)}](\tau_1)
+\sum_{p=0}^2\int_{\MM(\tau_1,\tau_2-1)}r^{1+\de}|\dk^{\leq {10}}{\N}_{W,+2}^{(p)}|^2\\
\nn&&+\sum_{p=0}^1\int_{\MM(\tau_1,\tau_2-1)}r^{-1+\de}|\dk^{\leq 11} \N_{T, +2}^{(p)}|^2
+\Big(R_0^{-\frac{\de}{2}}+\ep^2\Big)\sum_{p=0}^1{\EMF}^{(11)}_{\de}[\phiplus{p}](\tau_1,\tau_2-1)\\
\nn&&+R_0^{-\frac{\de}{2}}\F^{(10)}_{\Si_{*,c_*/2}}[\phiplus{0}](\tau_1,\tau_2-1)
\eeaa
and
\beaa
\nn&&\sum_{p=0}^2\EMF_{r\leq\varsigma_*(\tau-c_*/2)}[\psiminus{p}](\tau_1+1, \tau_2-3)
+\sum_{p=0}^2\EMF[\pmb\phi_{-2}^{(p)}](\tau_1, \tau_2-1)\\
&&+\sum_{p=0}^2\F_{\Si{*,c_*/2}}[\pmb\phi_{-2}^{(p)}](\tau_1, \tau_2-1)\\
\nn&\les& \sum_{p=0}^2\E^{({13})}[\pmb\phi_{-2}^{(p)}](\tau_1)+\int_{\Si(\tau_1)}|\dk^{\leq {12}}\N_{T,-2}^{(0)}|^2 +\E^{({13})}_{r\leq r_+(1+\dred )}[\Ab](\tau_1)\\
\nn&&+\sum_{p=0}^2\int_{\MM(\tau_1,\tau_2)}r^{1+\de}|\dk^{\leq {13}}\N_{W,-2}^{(p)}|^2+\sum_{p=0}^1\int_{\MM(\tau_1,\tau_2)}r^{-1+\de}|\dk^{\leq {13}}\N_{T,-2}^{(0)}|^2\\
&+&\int_{\MM_{r\leq r_+(1+\dred )}(\tau_1,\tau_2)}|\dk^{\leq {13}}\N_{\Ab}|^2+\ep^2\sup_{\tau\in[\tau_1, \tau_2]}\E^{({14})}[\pmb\phi_{-2}^{(0)}](\tau)+\ep^2\sup_{\tau\in[\tau_1, \tau_2]}\E^{({14})}_{r\leq r_+(1+\dred )}[\Ab](\tau).
\eeaa
Together with the local energy estimates in  \eqref{eq:localenergyforpsiontau1minus1tau1plus1} with $\reg=0$, the local energy estimates \eqref{eq:localenergyestimateforpsisijponSigmatau2minus3tau2:partialderivatives} with $\reg=0$, the one in  \eqref{eq:localenergyestimate:future:bis} applied with $\reg=0$, $\tau_0=\tau_2-1$,  $q=1$ and $\psi=(\phiplus{p})_{p=0,1,2}$, the one in \eqref{eq:localenergyestimate:future:bis} applied with $\reg=0$, $\tau_0=\tau_2-1$, $q=1$ and $\psi=(\phiminus{p})_{p=0,1,2}$, and applying \eqref{eq:EMFestimates:tensorialwaveinsubextremalKerr:weakMorathm:MMhsetminusMMregion} with $(\tau_0, \tau_1)=(\tau_1, \tau_2)$ and then \eqref{eq:EMFestimates:tensorialwaveinsubextremalKerr:weakMorathm} with $\tau_0=\tau_2$, we infer
\bea\lab{eq:weakMorawetzestimatesforphiplus2:consequenceThoerem:2}
\nn&&\sum_{p=0}^2\left(\EMFh[\psiplus{p}](\Iti)+\sup_{\frac{c_*}{2}\leq c\leq c_*}\F_{\Si_{*,c}}[\psis{p}](\tau_1, \tau_2-1)
+\EMF[\pmb\phi_{+2}^{(p)}](\tau_1, \tau_2)\right)\\ 
\nn&\les& R_0^{1+\de}\sum_{p=0}^2\E^{({9})}[\pmb\phi_{+2}^{(p)}](\tau_1)
+\sum_{p=0}^1\int_{\MM(\tau_1,\tau_2-1)}r^{-1+\de}|\dk^{\leq 11}\N_{T, +2}^{(p)}|^2\\
&&
+\sum_{p=0}^2\int_{\MM(\tau_1,\tau_2)}r^{1+\de}|\dk^{\leq {10}}{\N}_{W,+2}^{(p)}|^2+\Big(R_0^{-\frac{\de}{2}}+\ep^2\Big)\sum_{p=0}^1{\EMF}^{(11)}_{\de}[\phiplus{p}](\tau_1,\tau_2-1)\nn\\
&&+R_0^{-\frac{\de}{2}}\F^{(10)}_{\Si_{*,c_*/2}}[\phiplus{0}](\tau_1,\tau_2-1)
\eea
and  
\bea\lab{eq:weakMorawetzestimatesforphiminus2:consequenceThoerem:2}
\nn&&\sum_{p=0}^2\left(\EMFh[\psiminus{p}](\Iti)+\sup_{\frac{c_*}{2}\leq c\leq c_*}\F_{\Si_{*,c}}[\psis{p}](\tau_1, \tau_2-1)
+\EMF[\pmb\phi_{-2}^{(p)}](\tau_1, \tau_2)\right)\\ 
\nn&\les& \sum_{p=0}^2\E^{({13})}[\pmb\phi_{-2}^{(p)}](\tau_1)+\int_{\Si(\tau_1)}|\dk^{\leq 12}\N_{T,-2}^{(0)}|^2 +\E^{({13})}_{r\leq r_+(1+\dred )}[\Ab](\tau_1)\\
\nn&&+\sum_{p=0}^2\int_{\MM(\tau_1,\tau_2)}r^{1+\de}|\dk^{\leq {13}}\N_{W,-2}^{(p)}|^2+\sum_{p=0}^1\int_{\MM(\tau_1,\tau_2)}r^{-1+\de}|\dk^{\leq {13}}\N_{T,-2}^{(p)}|^2\\
\nn&&+\int_{\MM_{r\leq r_+(1+\dred )}(\tau_1,\tau_2)}|\dk^{\leq {13}}\N_{\Ab}|^2\\
&&+\ep^2\sup_{\tau\in[\tau_1, \tau_2]}\E^{({14})}[\pmb\phi_{-2}^{(0)}](\tau)+\ep^2\sup_{\tau\in[\tau_1, \tau_2]}\E^{({14})}_{r\leq r_+(1+\dred )}[\Ab](\tau).
\eea

Finally, using \eqref{eq:weakMorawetzestimatesforphiplus2:consequenceThoerem:2} and \eqref{eq:weakMorawetzestimatesforphiminus2:consequenceThoerem:2} to control the lower order terms which appear on the last line of the RHS of \eqref{eq:allderivativesrecoveredandonlylowerordertermonRHS:plus2caseandminus2caseatthesametime}, we infer
\bea\lab{eq:allderivativesrecoveredevenlowerordertermonRHS:plus2case:00000000}
\nn&& \sum_{p=0}^2\EMF^{({11})}_{\de}[\pmb\phi_{+2}^{(p)}](\tau_1, \tau_2)\\ 
\nn&\les& R_0^{1+\de}\sum_{p=0}^2\E^{({11})}[\pmb\phi_{+2}^{(p)}](\tau_1)+\sum_{p=0}^1\int_{\MM(\tt_1, \tt_2)}r^{-1+\de}|\dk^{\leq {12}}\N_{T,+2}^{(p)}|^2\\
\nn&&+\sum_{p=0}^2\NN^{({11})}_\de[\pmb\phi_{+2}^{(p)}, \N_{W,+2}^{(p)}](\tt_1, \tt_2)
+\Big(R_0^{-\frac{\de}{2}}+\ep^2\Big)\sum_{p=0}^1{\EMF}^{(11)}_{\de}[\phiplus{p}](\tau_1,\tau_2-1)\\
&&+R_0^{-\frac{\de}{2}}\F^{(10)}_{\Si_{*,c_*/2}}[\phiplus{0}](\tau_1,\tau_2-1)
\eea
and
\bea\lab{eq:allderivativesrecoveredevenlowerordertermonRHS:minus2case}
\nn \sum_{p=0}^2\EMF^{({14})}_{\de}[\pmb\phi_{-2}^{(p)}](\tau_1, \tau_2)&\les& \sum_{p=0}^2\E^{({14})}[\pmb\phi_{-2}^{(p)}](\tau_1)+\E^{({13})}_{r\leq r_+(1+\dred )}[\Ab](\tau_1)+\int_{\Si(\tau_1)}|\dk^{\leq {12}}\N_{T,-2}^{(0)}|^2\\
\nn&&+\sum_{p=0}^2\NN_\de^{({14})}[\pmb\phi_{-2}^{(p)}, \N_{W,-2}^{(p)}](\tt_1, \tt_2)+\int_{\MM_{r\leq r_+(1+\dred )}(\tau_1,\tau_2)}|\dk^{\leq {13}}\N_{\Ab}|^2\\
\nn&&+\sum_{p=0}^1\int_{\MM(\tt_1, \tt_2)}r^{-1+\de}|\dk^{\leq {15}}\N_{T,-2}^{(p)}|^2\\
&&+\ep^2\sup_{\tau\in[\tau_1, \tau_2]}\E^{({14})}[\pmb\phi_{-2}^{(0)}](\tau)
+\ep^2\sup_{\tau\in[\tau_1, \tau_2]}\E^{({14})}_{r\leq r_+(1+\dred )}[\Ab](\tau),
\eea
where we used in the derivation of \eqref{eq:allderivativesrecoveredevenlowerordertermonRHS:plus2case} \eqref{eq:allderivativesrecoveredevenlowerordertermonRHS:minus2case} the fact that, for $s=\pm 2$ and $\reg\leq 14$,
\beaa
\sum_{p=0}^2\int_{\MM(\tau_1,\tau_2)}r^{1+\de}|\dk^{\leq \reg}{\N}_{W,s}^{(p)}|^2
\les\sum_{p=0}^2\NN^{(\reg)}_\de[\pmb\phi_{s}^{(p)}, \N_{W,s}^{(p)}](\tt_1, \tt_2)
\eeaa
in view of the definition of $\NN^{(\reg)}_\de[\c, \c](\tau_1, \tau_2)$ in Section \ref{subsection:basicnormsforpsi:MMh}. Also, combining \eqref{eq:allderivativesrecoveredevenlowerordertermonRHS:plus2case:00000000} with Proposition \ref{prop:EnergyMorawetznearinfinitytensorialTeuk}, we obtain
\bea\lab{eq:allderivativesrecoveredevenlowerordertermonRHS:plus2case}
\nn&& \sum_{p=0}^2\EMF^{({11})}_{\de}[\pmb\phi_{+2}^{(p)}](\tau_1, \tau_2)+\F^{(10)}_{\Si_{*,c_*/2}}[\phiplus{0}](\tau_1,\tau_2-1)\\ 
\nn&\les& R_0^{1+\de}\sum_{p=0}^2\E^{({11})}[\pmb\phi_{+2}^{(p)}](\tau_1)+\sum_{p=0}^1\int_{\MM(\tt_1, \tt_2)}r^{-1+\de}|\dk^{\leq {12}}\N_{T,+2}^{(p)}|^2\\
\nn&&+\sum_{p=0}^2\NN^{({11})}_\de[\pmb\phi_{+2}^{(p)}, \N_{W,+2}^{(p)}](\tt_1, \tt_2)
+\Big(R_0^{-\frac{\de}{2}}+\ep^2\Big)\sum_{p=0}^1{\EMF}^{(11)}_{\de}[\phiplus{p}](\tau_1,\tau_2-1)\\
&&+R_0^{-\frac{\de}{2}}\F^{(10)}_{\Si_{*,c_*/2}}[\phiplus{0}](\tau_1,\tau_2-1),
\eea
and combining \eqref{eq:allderivativesrecoveredevenlowerordertermonRHS:minus2case} with the redshift estimates of Corollary 6.24 in \cite{MaSz26} implies 
\bea\lab{eq:allderivativesrecoveredevenlowerordertermonRHS:minus2case:withredshift}
\nn &&\sum_{p=0}^2\EMF^{({14})}_{\de}[\pmb\phi_{-2}^{(p)}](\tau_1, \tau_2)+\sum_{p=0}^2\EMF^{({14})}_{r\leq r_+(1+\dred)}[\nab_4^p\Ab](\tau_1, \tau_2)\\
\nn&\les& \sum_{p=0}^2\E^{({14})}[\pmb\phi_{-2}^{(p)}](\tau_1)+\sum_{p=0}^2\E_{r\leq r_+(1+2\dred )}^{({14})}[\nab_4^p\Ab](\tau_1)+\int_{\Si(\tau_1)}|\dk^{\leq 12}\N_{T,-2}^{(0)}|^2\\
\nn&&+\sum_{p=0}^2\NN_\de^{({14})}[\pmb\phi_{-2}^{(p)}, \N_{W,-2}^{(p)}](\tt_1, \tt_2)+\sum_{p=0}^1\int_{\MM(\tt_1, \tt_2)}r^{-1+\de}|\dk^{\leq {15}}\N_{T,-2}^{(p)}|^2\\
\nn&&+\sum_{p=0}^2\int_{\MM_{r\leq r_+(1+2\dred)}(\tau_1, \tau_2)}|\dk^{\leq {14}}\N_{\nab_4^p\Ab}|^2\\
&&+\ep^2\sup_{\tau\in[\tau_1, \tau_2]}\E^{({14})}[\pmb\phi_{-2}^{(0)}](\tau)
+\ep^2\sup_{\tau\in[\tau_1, \tau_2]}\E^{({14})}_{r\leq r_+(1+\dred )}[\Ab](\tau).
\eea
Now, for $R_0$ large enough and $\ep$ small enough, we may absorb the last terms with $\ep^2$ coefficients on the RHS of \eqref{eq:allderivativesrecoveredevenlowerordertermonRHS:minus2case:withredshift} and the last terms respectively with $R_0^{-\frac{\de}{2}}+\ep^2$ and $R_0^{-\frac{\de}{2}}$ coefficients on the RHS of \eqref{eq:allderivativesrecoveredevenlowerordertermonRHS:plus2case} which yields 
\beaa
\nn&&\sum_{p=0}^2\EMF^{({11})}_{\de}[\pmb\phi_{+2}^{(p)}](\tau_1, \tau_2)\\ 
\nn&\les&  \sum_{p=0}^2\E^{({11})}[\pmb\phi_{+2}^{(p)}](\tau_1)
+\sum_{p=0}^1\int_{\MM(\tt_1, \tt_2)}r^{-1+\de}|\dk^{\leq {12}}\N_{T,+2}^{(p)}|^2+\sum_{p=0}^2\NN^{({11})}_\de[\pmb\phi_{+2}^{(p)}, \N_{W,+2}^{(p)}](\tt_1, \tt_2)
\eeaa
and
\beaa
\nn &&\sum_{p=0}^2\EMF^{({14})}_{\de}[\pmb\phi_{-2}^{(p)}](\tau_1, \tau_2)+\sum_{p=0}^2\EMF^{({14})}_{r\leq r_+(1+\dred)}[\nab_4^p\Ab](\tau_1, \tau_2)\\
\nn&\les& \sum_{p=0}^2\E^{({14})}[\pmb\phi_{-2}^{(p)}](\tau_1)+\sum_{p=0}^2\E_{r\leq r_+(1+2\dred )}^{({14})}[\nab_4^p\Ab](\tau_1)+\int_{\Si(\tau_1)}|\dk^{\leq 12}\N_{T,-2}^{(0)}|^2\\
\nn&&+\sum_{p=0}^2\NN_\de^{({14})}[\pmb\phi_{-2}^{(p)}, \N_{W,-2}^{(p)}](\tt_1, \tt_2)+\sum_{p=0}^1\int_{\MM(\tt_1, \tt_2)}r^{-1+\de}|\dk^{\leq {15}}\N_{T,-2}^{(p)}|^2\\
&&+\sum_{p=0}^2\int_{\MM_{r\leq r_+(1+2\dred)}(\tau_1, \tau_2)}|\dk^{\leq {14}}\N_{\nab_4^p\Ab}|^2,
\eeaa
as stated. This concludes the proof of Theorem \ref{thm:main:MaSz26}.

%%%%%%%%%%%%%%%%%%%%%%%%%%%%%%%%

\subsection{Proof of Theorem \ref{th:main:intermediary}} 
\lab{sec:proofofth:main:intermediary:00}

%%%%%%%%%%%%%%%%%%%%%%%%%%%%%%%%

We adapt the proof of Theorem 7.6 in \cite{MaSz26}, see sections 8--10 in \cite{MaSz26}, where Theorem 7.6 in \cite{MaSz26} is the analog of Theorem \ref{th:main:intermediary} in the case where $\MM$ extends to $\II_+$ rather than $\Si_*$. 

\begin{remark}\lab{rmk:wesimplyhavetoextendtheproofofsections8to10ofMaSz26toregionMrgeqR0}
Note that the proof of Theorem 8.19 in \cite{MaSz26} treats separately the spacetime regions $\MMh_{r\leq R_0}$ and $\MMh_{r\geq R_0}$, where the large constant $R_0$ is chosen in Remark \ref{rmk:choiceofconstantRbymeanvalue} below. Since $\Rmic$ is such that $\chi_{r_*}=\chi_{r_*}^{(1)}=1$ on $\MMh_{r\leq R_0}$, where $\chi_{r_*}$ and $\chi^{(1)}_{r_*}$ are introduced respectively in \eqref{eq:supportpropertiesofcutoffchistar} and \eqref{eq:supportpropertiesofcutoffchistar:(1)}, the system of scalar wave equations \eqref{eq:waveeqwidetildepsi1} coincides in $\MMh_{r\leq R_0}$ with (7.9) in \cite{MaSz26}. In particular, all proofs of \cite{MaSz26} in $\MMh_{r\leq R_0}$ apply, and it thus suffices here to show how to extend to our setting the proofs of \cite{MaSz26} in $\MMh_{r\geq R_0}$.  
\end{remark}

Taking into account Remark \ref{rmk:wesimplyhavetoextendtheproofofsections8to10ofMaSz26toregionMrgeqR0}, we extend all results of Sections 8--10 in \cite{MaSz26} concerning the region $\MMh_{r\geq R_0}$ to our setting as follows:
\begin{itemize}
\item The results of Sections 8 in \cite{MaSz26} concerning the region $\MMh_{r\geq R_0}$ are extended to our setting in Section \ref{sect:microlocalenergyMorawetztensorialwaveequation}. 

\item The results of Sections 9 in \cite{MaSz26} concerning the region $\MMh_{r\geq R_0}$ are extended to our setting in Section \ref{sec:proofofth:main:intermediary:0order}.

\item The results of Sections 10 in \cite{MaSz26} concerning the region $\MMh_{r\geq R_0}$ are extended to our setting in Section \ref{sec:proofofth:main:intermediary}.
\end{itemize}

\begin{remark}
Note that $\MMh_{r\geq R_0}$ lies within the large $r$ region where $\gh$ is close to Minkowski and satisfies in particular $\MMh_{r\geq R_0}\subset\Mntraph$. In particular, the treatment of $\{r\geq R_0\}$  corresponds with the easy part of Sections 8--10 in \cite{MaSz26}. 
\end{remark}

%%%%%%%%%%%%%%%%%%%%%%%%%%%%%%%%%%%%%%%%%%%%%%%

\subsubsection{Global energy-Morawetz estimates for a system of scalar wave equations}
\lab{sect:microlocalenergyMorawetztensorialwaveequation}

%%%%%%%%%%%%%%%%%%%%%%%%%%%%%%%%%%%%%%%%%%%%%%%

In this section, we consider a system of scalar wave equations for complex-valued scalars $\psi_{ij}$
\begin{align}
\lab{eq:ScalarizedWaveeq:general:Kerrpert} 
\big({\square}_{\g}-D_0|q|^{-2}\big)\psi_{ij} =\widehat{F}_{ij}:=&\chi_{\tau_1, \tau_2}\chi_{r_*}\big(\widehat{S}(\psi)_{ij} +(\widehat{Q}\psi)_{ij})+(1-\chi_{\tau_1, \tau_2}\chi_{r_*})\big( \widehat{S}_K(\psi)_{ij} +(\widehat{Q}_K\psi)_{ij}\big)\nn\\
&+(1-\chi_{\tau_1, \tau_2}\chi_{r_*})f_{D_0}\psi_{ij}+F_{ij} ,\quad i,j=1,2,3, \quad\textrm{on}\quad\MMh,
\end{align}
where $\g$ satisfy satisfies \eqref{eq:controloflinearizedinversemetriccoefficients} \eqref{eq:controloflinearizedinversemetriccoefficients:inversegtautau} and coincides with Kerr in $\MMh\setminus(\MM_{r\leq\varsigma_*(\tau-\frac{7c_*}{8})}(\tau_1, \tau_2-1))$, where the cut-offs $\chi_{\tau_1, \tau_2}$ and $\chi_{r_*}$ satisfy respectively \eqref{eq:propertieschitoextendmetricg} and \eqref{eq:supportpropertiesofcutoffchistar}, where  $D_0>0$ is a constant, where
\bea
f_{D_0}:=\frac{4-D_0}{|q|^2}- \frac{4a^2\cos^2\th(|q|^2+6mr)}{|q|^6},
\eea
where
\bsub
\label{hatSandV:generalwave:Kerrpert}
\bea
\widehat{S}(\psi)_{ij} &=& {S}(\psi)_{ij}+\frac{4ia\cos\th}{|q|^2} \pr_{\tt}(\psi_{ij})\nn\\
&=&2M_{i}^{k\a}\pr_\a(\psi_{kj}) +2M_{j}^{k\a}\pr_\a(\psi_{ik})+\frac{4ia\cos\th}{|q|^2} \pr_{\tt}(\psi_{ij}),\\
(\widehat{Q}\psi)_{ij} &=&({Q}\psi)_{ij} -\frac{4ia\cos\th}{|q|^2}\big(M_{i\tau}^l \psi_{jl}+M_{j\tau}^l \psi_{il}\big)\nn\\
&=& (\Ddot^\a M_{i\a}^k)\psi_{kj}+(\Ddot^\a M_{j\a}^k)\psi_{ik} -M_{i\a}^kM_k^{l\a}\psi_{lj}-2M_{i\a}^kM_{j}^{l\a}\psi_{kl}\nn\\
&&-M_{j\a}^kM_k^{l\a}\psi_{il}-\frac{4ia\cos\th}{|q|^2}\big(M_{i\tau}^l \psi_{jl}+M_{j\tau}^l \psi_{il}\big)
\eea
\esub
with the 1-forms $M_{i\a}^j$ defined by \eqref{eq:definitionofMalphaijwithoutambiguity}.

\begin{remark}
Note that the solution $\psi^{(p)}_{s,ij}$ of \eqref{eq:waveeqwidetildepsi1} satisfies \eqref{eq:ScalarizedWaveeq:general:Kerrpert} with the choices $D_0=4-2\de_{p0}$ and $F_{ij}=\widetilde{F}^{(p)}_{s,ij}+\underline{F}^{(p)}_{s,ij}+\breve{F}^{(p)}_{s,ij}$, where $\widetilde{F}^{(p)}_{s,ij}$, $\underline{F}^{(p)}_{s,ij}$ and $\breve{F}^{(p)}_{s,ij}$ are defined respectively in \eqref{def:tildef}, \eqref{def:tildef0} and \eqref{def:breveFpsij}. This justifies the introduction of the model problem \eqref{eq:ScalarizedWaveeq:general:Kerrpert}.
\end{remark}

The aim of this section is to prove the following theorem which extends the microlocal energy-Morawetz estimates stated in Theorem 8.19 of \cite{MaSz26} to the above coupled system of inhomogenous scalar wave equations \eqref{eq:ScalarizedWaveeq:general:Kerrpert}. 

\begin{theorem}[Global energy-Morawetz estimates for \eqref{eq:ScalarizedWaveeq:general:Kerrpert}]\lab{th:mainenergymorawetzmicrolocal}
Let $(\MM, \g)$ satisfy the assumptions of Section \ref{sec:assumptionsforsec:energyMorawetzesitmatesforTeukoslkyonMM:upto15derivatives} with  $\g=\gam$ in the region $\MMh\setminus(\MM_{r\leq\varsigma_*(\tau-\frac{7c_*}{8})}(\tau_1, \tau_2-1))$, and let $\psi_{ij}$ be a solution to the system of scalarized wave equations \eqref{eq:ScalarizedWaveeq:general:Kerrpert}  with RHS $F_{ij}$, where the cut-offs $\chi_{\tau_1, \tau_2}$ and $\chi_{r_*}$ appearing in \eqref{eq:ScalarizedWaveeq:general:Kerrpert} satisfy respectively \eqref{eq:propertieschitoextendmetricg} and \eqref{eq:supportpropertiesofcutoffchistar}. Assume that $\psi_{ij}$ vanishes in $\MM(-\infty, \tmic)\cap\{r\leq (\Nmic+1)m\}$, with $\tmic$ and $N_0$ introduced in Definition \ref{def:definitionofthetimetauR}, and satisfies $\psi_{ij}=\pmb\psi(\Om_i, \Om_j)$ in $\MM(\tau_1+1, \tau_2-2)$ for a tensor $\pmb\psi\in\sk_2(\mathbb{C})$ and $\psi_{ij}=\pmb\psi((\Om_K)_i, (\Om_K)_j)$ in $\MM(\tau_2, +\infty)$ for a tensor $\pmb\psi\in\sk_{2,K}(\mathbb{C})$. Finally, assume that $F_{ij}$ are supported in $\MMh(\tmic, \tau_1)\cup\MM(\tau_1, \tau_2)$. Then, we have 
\begin{align}\lab{th:eq:mainenergymorawetzmicrolocal:tensorialwave:scalarized:eachpsisp}
&\widetilde{\EMF}[\pmb \psi] +\sup_{\frac{c_*}{2}\leq c\leq c_*}\F_{\Si_{*,c}}[\pmb\psi](\tau_1+1, \tau_2-2)\nn\\
\les& \sup_{\tau\in[\tmic,\tau_1+2]}\E[\pmb\psi](\tau)+
\Errdefect[\pmb\psi]+\Errdefects[\pmb\psi] +\A[\pmb \psi](\Iti)+\A_*[\pmb \psi]+\NNt[\pmb \psi,\pmb F]+\NNt_*[\pmb \psi,\pmb F],
\end{align}
where we have defined
\bea
\widetilde{\EMF}[\pmb \psi]=\EMFh[\pmb \psi]+\widetilde{\M}[\pmb\psi],
\eea
with $\widetilde{\M}[\pmb\psi]$ being the microlocal Morawetz norm on $\MMh_{r\leq 10m}$ introduced in (8.23) of \cite{MaSz26}, 
\bea
\lab{def:NNtintermsofNNtMora:NNtEner:NNtaux:wavesystem:EMF} 
\NNt[\pmb \psi,\pmb F]:=\NNtmora[\pmb \psi, \pmb F]+\NNtener[\pmb \psi, \pmb F]
+\NNtaux[\pmb F],
\eea
with $\NNtmora[\pmb \psi, \pmb F]$, $\NNtener[\pmb \psi, \pmb F]$ and $\NNtaux[\pmb F]$ introduced in (8.26) of \cite{MaSz26}, 
\bea\lab{eq:definitionofNNt*pmbpsipmbF:additionalterm*}
\NNt_*[\pmb \psi,\pmb F]:=\sup_{\frac{c_*}{2}\leq c\leq c_*}\bigg|\int_{\MM_{11m\leq r\leq\varsigma_*(\tau-c)}(\tau_1+1,\tau_2-2)}\sum_{i,j}\Re\big(F_{ij}\ov{T_{\tau_2}(\psi)_{ij}}\big)\bigg|,
\eea
with $T_{\tau_2}(\psi)_{ij}$ given by \eqref{def:widehatTpsiij:equalsnabTpsiij}, and
\bsub\lab{def:EM-1norms:Reals}
\begin{align}
\A[\pmb \psi](\Iti) :=&\sum_{i,j}\bigg(\int_{\MMh(\Iti)}r^{-3}|\psi_{ij}|^2 +\sup_{\tau\in\Iti}\int_{\Sih(\tau)}r^{-2}|\psi_{ij}|^2+\int_{\II_+(\Iti)}r^{-2}|\psi_{ij}|^2\bigg),\\
\A_*[\pmb \psi] :=& \sum_{i,j}\sup_{\frac{c_*}{2}\leq c\leq c_*}\int_{\Si_{*,c}(\tau_1+1, \tau_2-2)}r^{-2}|\psi_{ij}|^2,
\end{align}
\esub
\bsub\lab{def:Errdefectofpsi}
\begin{align}
\Errdefect[\pmb\psi] :=&\sup_{\tau\in[\tmic,\tau_1+1]\cup[\tau_2-2,\tau_2]}\sum_{j}\Big(\E[x^i\psi_{ij}](\tau)+\E[x^i\psi_{ji}](\tau)\Big),\\
\Errdefects[\pmb\psi] :=&\sup_{\frac{c_*}{2}\leq c\leq c_*}\sum_{j}\Big(\F_{\Si_{*,c}}[x^i\psi_{ij}](\tau_1+1, \tau_2-2)+\F_{\Si_{*,c}}[x^i\psi_{ji}](\tau_1+1, \tau_2-2)\Big).
\end{align}
\esub 
\end{theorem}

\begin{remark}
Comparing \eqref{th:eq:mainenergymorawetzmicrolocal:tensorialwave:scalarized:eachpsisp} with the corresponding estimate (8.31) in Theorem 8.19 of \cite{MaSz26}, we emphasize the following differences:
\begin{itemize}
\item the LHS of \eqref{th:eq:mainenergymorawetzmicrolocal:tensorialwave:scalarized:eachpsisp} controls in addition the quantity $\sup_{\frac{c_*}{2}\leq c\leq c_*}\F_{\Si_{*,c}}[\pmb\psi](\tau_1+1, \tau_2-2)$,

\item the new terms on the RHS of \eqref{th:eq:mainenergymorawetzmicrolocal:tensorialwave:scalarized:eachpsisp} compared to (8.31) in \cite{MaSz26} are $\Errdefects[\pmb\psi]$, $\A_*[\pmb \psi]$ and $\NNt_*[\pmb \psi,\pmb F]$.
\end{itemize}
\end{remark}

As emphasized in Remark \ref{rmk:wesimplyhavetoextendtheproofofsections8to10ofMaSz26toregionMrgeqR0}, the proof of Theorem 8.19 of \cite{MaSz26} treats separately the spacetime regions $\MMh_{r\leq R_0}$ and $\MMh_{r\geq R_0}$, where the large constant $R_0$ is chosen in Remark \ref{rmk:choiceofconstantRbymeanvalue} below.

\begin{remark}[Choice of the constant $\Rmic$]\label{rmk:choiceofconstantRbymeanvalue}
The constant $\Rmic\in [{\Nmic}m, ({\Nmic}+1)m]$, with $\Nmic\geq 20$ a large enough integer, is chosen to verify 
\bea
\label{eq:choiceofRvalue:Kerr}
\nn&&\sum_{i,j=1}^3\int_{H_{\Rmic}(\Iti)}\big(|{\pr^{\leq 1}}\psi_{ij}|^2+|\square_\g\psi_{ij}|^2\big) d\tt dx^1dx^2\\
&\leq& \frac{1}{m}\sum_{i,j=1}^3\int_{\MM_{\Nmic m, (\Nmic +1)m}(\Iti)}\big(|{\pr^{\leq 1}}\psi_{ij}|^2+|\square_\g\psi_{ij}|^2\big)d\tau d r dx^1 dx^2.
\eea
This implies that for solutions to the coupled system of scalar wave equations \eqref{eq:ScalarizedWaveeq:general:Kerrpert}, we have
\bea
\label{eq:choiceofRvalue:Kerr:extendedRWsystem}
 \sum_{i,j=1}^3\int_{H_{\Rmic}(\Reals)}\big(|{\pr^{\leq 1}}\psi_{ij}|^2+|\square_\g\psi_{ij}|^2\big) d\tt dx^1dx^2\les \frac{1}{m}\sum_{i,j=1}^3\int_{\MM_{{\Nmic}m, ({\Nmic}+1)m}(\Iti)}\big(|{\pr^{\leq 1}}\psi_{ij}|^2+|F_{ij}|^2\big).
\eea
\end{remark}

Since $\Rmic\leq (\Nmic+1)m\ll \ep^{-1}\ll c_*$, and in view of \eqref{eq:supportpropertiesofcutoffchistar}, we have $\chi_{r_*}=1$ on $\MMh_{r\leq R_0}$. Therefore, in the region $\MMh_{r\leq R_0}$,  \eqref{eq:ScalarizedWaveeq:general:Kerrpert} coincides with (8.1) in \cite{MaSz26}. In particular, all proofs of \cite{MaSz26} in $\MMh_{r\leq R_0}$ regarding energy-Morawetz estimates for \eqref{eq:ScalarizedWaveeq:general:Kerrpert} apply, and, to prove Theorem \ref{th:mainenergymorawetzmicrolocal}, it thus suffices to show how to extend the proofs of \cite{MaSz26} in $\MMh_{r\geq R_0}$ to \eqref{eq:ScalarizedWaveeq:general:Kerrpert}. In particular, the results and proofs in section 8.3 of \cite{MaSz26} immediately apply to \eqref{eq:ScalarizedWaveeq:general:Kerrpert}, and we initiate the proof of Theorem \ref{th:mainenergymorawetzmicrolocal} by showing how to extend the proofs of section 8.4 in \cite{MaSz26} to \eqref{eq:ScalarizedWaveeq:general:Kerrpert}.

We start by introducing the following quantity 
\bea\lab{def:NNtlocalinr:NNtEner:NNtaux:wavesystem:EMF:proof}
\NNtlocal[\pmb \psi](\tau',\tau'') &:=&\bigg(\sup_{\tau\in (\tau',\tau'')}\int_{\Sih(\tau)}r^{-2}|\pmb \psi|^2
 +\int_{\MMh(\tau',\tau'')}r^{-3}|\pmb\psi|^2+\int_{\II_+(\tau',\tau'')}r^{-2}|\pmb \psi|^2\nn\\
 &&+\sup_{\frac{c_*}{2}\leq c\leq c_*}\int_{\Si_{*,c}(\tau_*', \tau_*'')}r^{-2}|\pmb \psi|^2+\Errdefect[\pmb\psi]+\Errdefects[\pmb\psi]\bigg)^{\frac{1}{2}}\nn\\
 &&\times\bigg(\widetilde{\EMF}[\pmb \psi](\tau',\tau'')+\sup_{\frac{c_*}{2}\leq c\leq c_*}\F_{\Si_{*,c}}[\pmb\psi](\tau_*', \tau_*'')\bigg)^{\frac{1}{2}},
 \eea
where $(\tau_*', \tau_*'')$ is such that 
\bea\lab{eq:definitionoftau*primetau*doubleprimefromtauprimtaudoubleprime}
[\tau_*', \tau_*'']=[\tau', \tau'']\cap[\tau_1+1, \tau_2-2].
\eea

Next, recall from Proposition \ref{prop-app:stadard-comp-Psi} that the $1$-form   $\PP_\mu[\pmb\psi](X, w)$, for a vectorfield $X$, a real scalar function $w$ and a tensor $\pmb\psi\in \mathfrak{s}_k(\mathbb{C})$, $k=0,1,2$, is given by 
 \beaa
\PP_\mu[\pmb\psi](X, w)=\QQ_{\mu\nu}[\pmb\psi] X^\nu +\frac 1 2  w \Re\big(\pmb\psi \c \ov{\Db_\mu \pmb\psi }\big)-\frac 1 4|\pmb\psi|^2   \pr_\mu w
  \eeaa
and satisfies
  \bea\lab{eq:DivofPPmu:tensor:RW}
\D^\mu  \PP_\mu[\pmb\psi](X, w)&=& \frac 1 2 \QQ[\pmb\psi]  \c\piX - \frac 1 2 X( V) |\pmb\psi|^2 +\frac{k}{2}{}^{(X)}A_\nu\Im\Big(\pmb\psi\c\ov{\Ddot^{\nu}\pmb\psi}\Big)+\frac 12  w \LL[\pmb\psi]\nn\\
  && -\frac 1 4|\pmb\psi|^2   \square_\g  w   +  \Re\bigg(\ov{\bigg(\nab_X\pmb\psi +\frac 1 2   w \pmb\psi\bigg)}\c \left(\squared_k \pmb\psi- V\pmb\psi\right)\bigg),
\eea
with ${}^{(X)}A_\nu$ the 1-form introduced in \eqref{eq:thespacetime1formXA}.
  
The following basic estimate for the $1$-form $\PP_\mu[\pmb\psi](X, w)$ is the analog of Lemma 8.28 in \cite{MaSz26}.

\begin{lemma}
\lab{lem:estimatesforcurrentforprttandprrBL}
Let $\pmb\psi\in\sk_2(\mathbb{C})$  and let $\psi_{ij}$ be the corresponding scalars given by $\psi_{ij}:=\pmb\psi(\Om_i, \Om_j)$ for $i,j=1,2,3$ where $\Om_i$, $i=1,2,3$ is the regular triplet introduced in Section \ref{sec:assumptionsforsec:energyMorawetzesitmatesforTeukoslkyonMM:upto15derivatives}. For a real-valued vectorfield $X=O(1)\pr_{\tt}+O(1)\pr_{r} + O(r^{-2})\pr_{x^a}$ and  a real-valued scalar function $w=O(r^{-1})$, we have
\bsub
\lab{eq:difference:PPmupmbpsiandpsiij}
\bea
\lab{eq:difference:PPmupmbpsiandpsiij:energy}
&&\int_{\Sih(\tau)}\bigg|\bigg({\PP}_{\mu}[\pmb\psi](X, w)- \sum_{i,j}\PP_\mu[\psi_{ij}](X, w)\bigg)N_{\Sih_{\tau}}^{\mu}\bigg|\nn\\
 &\les&\left(\int_{\Sih(\tau)}r^{-2}|\pmb\psi|^2\right)^{\frac{1}{2}}\Big(\E[\pmb\psi](\tau)\Big)^{\frac{1}{2}},\quad \forall \tau\geq \tmic, \\
\lab{eq:difference:PPmupmbpsiandpsiij:scri}
&&\int_{\II_+(\tau',\tau'')}\bigg|\bigg({\PP}_{\mu}[\pmb\psi](X, w)- \sum_{i,j}\PP_\mu[\psi_{ij}](X, w)\bigg)N_{\II_+}^{\mu}\bigg| \nn\\
&\les& \left(\int_{\II_+(\tau', \tau'')}r^{-2}|\pmb\psi|^2\right)^{\frac{1}{2}}\Big(\F_{\II_+}[\pmb \psi](\tau',\tau'')\Big)^{\frac{1}{2}},\quad \forall \tmic\leq \tau'<\tau'',\\
\lab{eq:difference:PPmupmbpsiandpsiij:Hr}
&&\int_{\Hh_{r}(\tau',\tau'')}\bigg|\bigg({\PP}_{\mu}[\pmb\psi](X, w)- \sum_{i,j}\PP_\mu[\psi_{ij}](X, w)\bigg)N_{H_r}^{\mu}\bigg| \nn\\
&\les& \int_{\Hh_{r}(\tau',\tau'')}r^{-1}|\pmb\psi| |\pr^{\leq 1}\pmb\psi|, \quad \forall r_+(1-\dhor)\leq r <+\infty, \,\, \tmic\leq \tau'<\tau'',\\
\lab{eq:difference:PPmupmbpsiandpsiij:RR*}
&&\int_{\RR_*(\tau',\tau'')}r^{-1}\bigg|\bigg({\PP}_{\mu}[\pmb\psi](X, w)- \sum_{i,j}\PP_\mu[\psi_{ij}](X, w)\bigg)(\DD(\tt)^{\mu}, \DD(r)^{\mu})\bigg| \nn\\
&\les& \left(\int_{\MM(\tau', \tau'')}r^{-3}|\pmb\psi|^2\right)^{\frac{1}{2}}\Big(\EM[\pmb\psi](\tau', \tau'')\Big)^{\frac{1}{2}}, \quad \forall \,\, 1\leq \tau'<\tau''\leq\tau_*,
\eea
\esub
where $\Hh_{r_1}=\{r=r_1\}$. 
\end{lemma}

\begin{proof}
The estimates \eqref{eq:difference:PPmupmbpsiandpsiij:energy}, \eqref{eq:difference:PPmupmbpsiandpsiij:scri} and \eqref{eq:difference:PPmupmbpsiandpsiij:Hr} are proved in Lemma 8.28 in \cite{MaSz26}. It has remains to prove \eqref{eq:difference:PPmupmbpsiandpsiij:RR*}. To thus end, recall from (8.94) in \cite{MaSz26} that we have the following identity for a real-valued vectorfield $X=O(1)\pr_{\tt}+O(1)\pr_{r} + O(r^{-2})\pr_{x^a}$ and a real-valued scalar function $w$
\beaa
&&\PP_\mu[\pmb\psi](X, w)- \sum_{i,j}\PP_\mu[\psi_{ij}](X, w)\nn\\
&=&\Re\Big(M_{i\mu}^j \psi\ov{X\psi}  + (O(r^{-3})+\Ga_b) \psi\ov{\pr_{\mu}\psi} +M_{i\mu}^j (O(r^{-3})+\Ga_b) \psi\ov{\psi} \Big)+w\Re\big(M_{i\mu}^j \psi \ov{\psi}\big)\nn\\
&&+\g_{\mu \a}X^{\a} \Re\Big(O(r^{-2})\psi \ov{(e_3)^{\leq 1}\psi }+O(r^{-1})\psi\ov{e_a\psi}+(O(r^{-2})+\Ga_b) \psi \ov{e_4 \psi}\Big),
\eeaa
where the schematic notations $\Re(\psi\ov{e_{\a}\psi})$ and $\Re(\psi\ov{\psi})$ denote respectively any term of the form $\Re(\psi_{ij}\ov{e_\a(\psi_{kl})})$ and $\Re(\psi_{ij}\ov{\psi_{kl}})$. Together with Lemma \ref{lem:estimatesforMialphaj:Kerrpert} according to which we have
\beaa
M_{i4}^j=O(r^{-2}), \quad M_{ia}^j = O(r^{-1}), \quad  M_{i3}^j=O(r^{-2})+\Ga_b,
\eeaa
we infer
\beaa
&&\int_{\RR_*(\tau',\tau'')}r^{-1}\bigg|\bigg({\PP}_{\mu}[\pmb\psi](X, w)- \sum_{i,j}\PP_\mu[\psi_{ij}](X, w)\bigg)(\DD(\tt)^{\mu}, \DD(r)^{\mu})\bigg|\\
&\les& \int_{\RR_*(\tau',\tau'')}r^{-2}|\psi|\Big(|e_a\psi|+r^{-1}(|e_4\psi|+|e_3\psi|+|\psi|\big)+\tau^{-1-\frac{3\dec}{4}}|e_4\psi|\Big)\\
&\les& \left(\int_{\MM(\tau', \tau'')}r^{-3}|\pmb\psi|^2\right)^{\frac{1}{2}}\Big(\EM[\pmb\psi](\tau', \tau'')\Big)^{\frac{1}{2}},
\eeaa
as stated in \eqref{eq:difference:PPmupmbpsiandpsiij:RR*}. This concludes the proof of Lemma \ref{lem:estimatesforcurrentforprttandprrBL}.
\end{proof}

Next, we introduce cut-offs that will appear in the statement of Proposition \ref{prop-app:stadard-comp-Psi:extendedscalarizedRW} below. Let $\chi_n(\tau)$, $n=1,2,3,4$, be smooth nonnegative cut-off functions satisfying
\begin{equation}
\lab{def:cutoffsintime1234}
\begin{split}
&\sum_{n=1}^4\chi_n(\tau)=1\,\,\forall\tau\in\Reals, \qquad\mathrm{supp}(\chi_1)\subset(-\infty, \tau_1+2), \qquad \mathrm{supp}(\chi_2)\subset(\tau_1+1, \tau_2-2),\\ 
&\mathrm{supp}(\chi_3)\subset(\tau_2-3, \tau_2+1), \qquad \mathrm{supp}(\chi_4)\subset(\tau_2, +\infty),
\end{split}
\end{equation}
let $\chi_{n,*}(\tau,r)$, $n=1,2,3$, be smooth nonnegative cut-off functions satisfying
\begin{equation}
\lab{def:cutoffsintime123*}
\begin{split}
&\sum_{n=1}^3\chi_{n,*}(\tau,r)=1\,\,\forall\tau\in\Reals,\,\, r\geq r_+(1-\dhor), \qquad\mathrm{supp}(\chi_{1,*})\subset\left\{r\leq\varsigma_*\left(\tau-\frac{5c_*}{8}\right)\right\},\\ 
&\mathrm{supp}(\chi_{2,*})\subset\RR_*, \qquad \mathrm{supp}(\chi_{3,*})\subset\left\{r\geq\varsigma_*\left(\tau-\frac{31c_*}{32}\right)\right\},\\
& \|(r\pr_\tau, r\pr_r)^{\leq 1}\chi_{r_*,n}\|_{L^\infty}\lesssim 1, \quad n=1,2,3.
\end{split}
\end{equation}
Finally, using the cut-offs in \eqref{def:cutoffsintime1234} and \eqref{def:cutoffsintime123*}, let $\widetilde{\chi}_{n}(\tau,r)$, $n=1,2,3$, be smooth nonnegative cut-off functions given by
\bea\lab{def:cutoffsintime123:widetilde}
\widetilde{\chi}_{1}:= \chi_1+\chi_3(\chi_{1,*}+\chi_{2,*})+\chi_2\chi_{2,*},\qquad \widetilde{\chi}_{2}:=&\chi_2\chi_{1,*}, \qquad \widetilde{\chi}_{3}:=\big(\chi_2+\chi_3\big)\chi_{3,*}+\chi_4,
\eea
which satisfy in view of  \eqref{def:cutoffsintime1234} \eqref{def:cutoffsintime123*}
\bea\lab{eq:mainpropertiesofwidetildechinn=123}
\bsplit
&\sum_{n=1}^3\widetilde{\chi}_{n}(\tau, r)=1\,\,\forall\tau\in\Reals,\,\, r\geq r_+(1-\dhor), \\ 
&\textrm{supp}(\widetilde{\chi}_1)\subset\MMh(-\infty, \tau_1+2)\cup\RR_*(\tau_1+2,\tau_2-3)\cup\MMh(\tau_2-3, \tau_2+1),\\
&\textrm{supp}(\widetilde{\chi}_2)\subset\MM_{r\leq\varsigma_*(\tau-\frac{5c_*}{8})}(\tau_1+1, \tau_2-2),\\
&\textrm{supp}(\widetilde{\chi}_3)\subset\MMh_{r\geq\varsigma_*(\tau-\frac{31c_*}{32})}(\tau_1+1, \tau_2)\cup\MMh(\tau_2, +\infty).  
\end{split}
\eea
The following divergence identity is the analog of Proposition 8.31 in \cite{MaSz26}.

\begin{proposition}\lab{prop-app:stadard-comp-Psi:extendedscalarizedRW}
Let $\psi_{ij}$ be a solution to the system of scalarized wave equations \eqref{eq:ScalarizedWaveeq:general:Kerrpert}. Under the same assumptions for the scalars $\psi_{ij}$, $F_{ij}$ and the spacetime $(\MM,\g)$ as in Theorem \ref{th:mainenergymorawetzmicrolocal}, let $X$ be  a real-valued vectorfield satisfying
 \bea
\lab{eq:assump:generalvectorfieldX:EMnearinf}
X=O(1)\pr_r + O(1)\pr_{\tau} + O(r^{-2})\pr_{x^a},
\eea
let $w$ be a real scalar function satisfying $w=O(r^{-1})$, let $\widetilde{\chi}_{n}(\tau,r)$, $n=1,2,3$, be the smooth nonnegative cut-off functions given by \eqref{def:cutoffsintime123:widetilde}, and define the following modified current
\begin{align}
\lab{eq:DefofwidetilePPmu:extendedscalarizedRW}
{\PP}_{\mu,\tau_2}[\pmb\psi](X, w):=& \widetilde{\chi}_{1}(\tau,r)\sum_{i,j}{\PP}_\mu[\psi_{ij}](X,w)+\widetilde{\chi}_{2}(\tau,r){\PP}_{\mu}[\pmb\psi](X, w)+\widetilde{\chi}_{3}(\tau,r)\big({\PP}_{\mu}[\pmb\psi](X, w)\big)_{K}, 
\end{align}
with ${\PP}_\mu[\pmb\psi](X, w)$ given in \eqref{definitionofcurrentPPmuXw:generaltensor} for $\pmb\psi\in\sk_2(\mathbb{C})$ and $V=D_0|q|^{-2}$, $\big({\PP}_\mu[\pmb\psi](X, w)\big)_{K}$ being the corresponding quantity in Kerr for $\pmb\psi\in\sk_2(\mathbb{C})$ and $V=\frac{4}{|q|^2}- \frac{4a^2\cos^2\th(|q|^2+6mr)}{|q|^6}$, and ${\PP}_\mu[\psi_{ij}](X,w)$ given in 
\eqref{definitionofcurrentPPmuXw:generaltensor} for a scalar $\psi_{ij}\in\sk_0(\mathbb{C})$ and $V=D_0|q|^{-2}$. 
 Then, ${\PP}_{\mu,\tau_2}[\pmb\psi](X, w)$ satisfies
 \bea
\lab{eq:1formPP_mu:equality}
{\PP}_{\mu,\tau_2}[\pmb\psi](X, w)=\sum_{i,j}{\PP}_\mu[\psi_{ij}](X,w)
+H_{\mu}[\psi,\pr\psi]
\eea
where $H_{\mu}[\psi,\pr\psi]$ satisfies for any $\tmic\leq\tau'<\tau''$, 
\bsub
\lab{eq:1formBB_mu:equality:errortermestimates}
\begin{align}
&\int_{\Sigma(\tau')}\big|H_{\mu}[\psi,\pr\psi]N_{\Sigma(\tau)}^{\mu}\big| \les \left(\int_{\Sigma(\tau')}r^{-2}|\pmb\psi|^2\right)^{\frac{1}{2}}\big(\E[\pmb\psi](\tau')\big)^{\frac{1}{2}},  \\
&\int_{\II_+(\tau',\tau'')}\big|H_{\mu}[\psi,\pr\psi]N_{\II_+}^{\mu} \big|\les  \left(\int_{\II_+(\tau', \tau'')}r^{-2}|\pmb\psi|^2\right)^{\frac{1}{2}}\Big(\F_{\II_+}[\pmb \psi](\tau',\tau'')\Big)^{\frac{1}{2}},\\
&\int_{H_{r}(\tau',\tau'')}\big|H_{\mu}[\psi,\pr\psi]N_{H_r}^{\mu}\big| \les \int_{H_{r}(\tau',\tau'')}r^{-1}|\pmb\psi| |\pr^{\leq 1}\pmb\psi|, \quad\forall r_+(1-\dhor)\leq r<+\infty,
\end{align}
\esub
and its divergence equals
\bea
\lab{eq:DivwidetilePPmu:extendedscalarizedRW}
\D^\mu {\PP}_{\mu,\tau_2}[\pmb\psi](X, w) &=&\sum_{i,j}\Re\bigg(F_{ij}\ov{\bigg(X_{\tau_2}(\psi)_{ij} + \frac{1}{2}w\psi_{ij}  \bigg)}\bigg) \nn\\
&&+\widetilde{\chi}_{1}(\tau,r)\frac{1}{2}\sum_{i,j}\Big(\QQ[\psi_{ij}]  \c\piX +w \LL[\psi_{ij}]  -X( V) |\psi_{ij}|^2-\frac{1}{2}|\psi_{ij}|^2   \square_\g  w\Big) \nn \\
&&+\widetilde{\chi}_{2}(\tau,r)\frac 1 2 \Big(\QQ[\pmb\psi]  \c\piX +w \LL[\pmb\psi] -X( V) |\pmb\psi|^2-\frac{1}{2}|\pmb\psi|^2   \square_\g  w\Big)\nn\\
&&+\widetilde{\chi}_{3}(\tau,r)\frac 1 2 \Big(\QQ[\pmb\psi]  \c\piX +w \LL[\pmb\psi] -X( V) |\pmb\psi|^2-\frac{1}{2}|\pmb\psi|^2   \square_\g  w\Big)_{K}\nn\\
&&+\Err_{\text{l.o.t.}}, \quad\textrm{for}\quad \tau\geq\tmic,
\eea
 where 
\bea
\lab{def:Xtau2psiij:generalforEM}
X_{\tau_2}(\psi)_{ij} := 
X(\psi_{ij}) - \widetilde{\chi}_{2}\big(X^{\a}M_{i\a}^k\psi_{kj}+ X^{\a}M_{j\a}^k\psi_{ik}\big) -\widetilde{\chi}_{3}\big(X^{\a} (M_K)_{i\a}^k\psi_{kj}+ X^{\a} (M_K)_{j\a}^k\psi_{ik}\big),
\eea
and where $\Err_{\text{l.o.t.}}$ denotes terms satisfying the following bound for any $r_1 \geq 10m$ and $(\tau',\tau'')\subset (\tmic,+\infty)$, 
\bea
\lab{def:Errlowerorderterm:DivPPmuchi}
\bigg|\int_{\MMh_{r_1, +\infty}(\tau',\tau'')}\Err_{\text{l.o.t.}}\bigg|
&\les&\sup_{\tau\in[\tmic,\tau_1+2]}\Eh[\pmb \psi](\tau)+\ep\sup_{\tt\in[\tau',\tau'']}\Eh[\pmb\psi](\tau)\nn\\
&& +\ep\sup_{\frac{c_*}{2}\leq c\leq c_*}\F_{\Si_{*,c}}[\pmb\psi](\tau_*', \tau_*'')+ \NNtlocal[\pmb \psi](\tau',\tau'')\nn\\
&&+\int_{H_{r_1}(\tau',\tau'')}r^{-1} |\pmb\psi| |\pr^{\leq 1}\pmb\psi|
\eea
with $\NNtlocal[\pmb \psi](\tau',\tau'')$ as given in \eqref{def:NNtlocalinr:NNtEner:NNtaux:wavesystem:EMF:proof} and $(\tau_*', \tau_*'')$ as given by \eqref{eq:definitionoftau*primetau*doubleprimefromtauprimtaudoubleprime}.
\end{proposition}

\begin{proof}
The identity \eqref{eq:1formPP_mu:equality} with $H_{\mu}[\psi,\pr\psi]$ satisfying \eqref{eq:1formBB_mu:equality:errortermestimates} follows from the definition \eqref{eq:DefofwidetilePPmu:extendedscalarizedRW} for ${\PP}_{\mu,\tau_2}[\pmb\psi](X, w)$ together with the first identity in \eqref{eq:mainpropertiesofwidetildechinn=123} and the estimates \eqref{eq:difference:PPmupmbpsiandpsiij} for a real-valued vectorfield $X=O(1)\pr_{\tt}+O(1)\pr_{r} + O(r^{-2})\pr_{x^a}$ and   a real-valued scalar function $w=O(r^{-1})$.

Next, we compute the divergence of ${\PP}_{\mu,\tau_2}[\pmb\psi](X, w)$.
We have, from the formula \eqref{eq:DefofwidetilePPmu:extendedscalarizedRW} for ${\PP}_{\mu,\tau_2}[\pmb\psi](X, w)$,
\bea
\lab{eq:computeDivergenceofBB:intermsofBulk:000000}
 \D^\mu{\PP}_{\mu,\tau_2}[\pmb\psi](X, w)
 &=&\sum_{n=1}^3\Bulkxw{n}+\bigg(\D^\mu(\widetilde{\chi}_{1})\sum_{i,j}{\PP}_\mu[\psi_{ij}](X,w)  \nn\\
 &&\qquad\qquad\qquad\qquad+\D^\mu(\widetilde{\chi}_{2}){\PP}_{\mu}[\pmb\psi](X, w)+\D^\mu(\widetilde{\chi}_{3})\big({\PP}_{\mu}[\pmb\psi](X, w)\big)_{K}\bigg)\nn\\
&=&\sum_{n=1}^3\Bulkxw{n}+\D^\mu(\widetilde{\chi}_{2})\bigg({\PP}_{\mu}[\pmb\psi](X, w)- \sum_{i,j}{\PP}_\mu[\psi_{ij}](X, w)\bigg)\nn\\
&&+\D^\mu(\widetilde{\chi}_{3})\bigg(\big({\PP}_{\mu}[\pmb\psi](X, w)\big)_{K}- \sum_{i,j}{\PP}_\mu[\psi_{ij}](X, w)\bigg),
\eea
where we have defined in the first equality of \eqref{eq:computeDivergenceofBB:intermsofBulk}
\bsub
\lab{def:Bulkxw:1234}
\begin{align}
\Bulkxw{1} :=& \widetilde{\chi}_{1}(\tau,r)\sum_{i,j}\D^\mu{\PP}_\mu[\psi_{ij}](X,w), \\
\Bulkxw{2} :=&\widetilde{\chi}_{2}(\tau,r)\D^\mu{\PP}_{\mu}[\pmb\psi](X, w),  \\
\Bulkxw{3} :=&\widetilde{\chi}_{3}(\tau,r)\D_K^\mu\big({\PP}_{\mu}[\pmb\psi](X, w)\big)_{K},
\end{align}
\esub
and used in the second equality of \eqref{eq:computeDivergenceofBB:intermsofBulk} the fact $\sum_{n=1}^3\D^\mu(\widetilde{\chi}_n)=0$ which follows from the first identity in \eqref{eq:mainpropertiesofwidetildechinn=123}. Then, in view of \eqref{def:cutoffsintime123:widetilde}, we have
\begin{align}\lab{eq:moredetialsontermsofdivergenceDDmuPPmupmspsionwhichderivativefallsoncutoffswidetildechinn123}
&\D^\mu(\widetilde{\chi}_{2})\bigg({\PP}_{\mu}[\pmb\psi](X, w)- \sum_{i,j}{\PP}_\mu[\psi_{ij}](X, w)\bigg)+\D^\mu(\widetilde{\chi}_{3})\bigg(\big({\PP}_{\mu}[\pmb\psi](X, w)\big)_{K}- \sum_{i,j}{\PP}_\mu[\psi_{ij}](X, w)\bigg)\nn\\
=& \chi_2'(\tau)\chi_{1,*}(\tau, r)\DD^\mu(\tau)\bigg({\PP}_{\mu}[\pmb\psi](X, w)- \sum_{i,j}{\PP}_\mu[\psi_{ij}](X, w)\bigg)\nn\\
&+\Big(\big(\chi_2'(\tau)+\chi_3'(\tau)\big)\chi_{3,*}(\tau, r)+\chi_4'(\tau)\Big)\DD^\mu(\tau)\bigg(\big({\PP}_{\mu}[\pmb\psi](X, w)\big)_{K}- \sum_{i,j}{\PP}_\mu[\psi_{ij}](X, w)\bigg)\nn\\
&+\chi_2(\tau)\DD^\mu(\chi_{1,*})\bigg({\PP}_{\mu}[\pmb\psi](X, w)- \sum_{i,j}{\PP}_\mu[\psi_{ij}](X, w)\bigg)\nn\\
&+\big(\chi_2(\tau)+\chi_3(\tau)\big)\DD^\mu(\chi_{3,*})\bigg(\big({\PP}_{\mu}[\pmb\psi](X, w)\big)_{K}- \sum_{i,j}{\PP}_\mu[\psi_{ij}](X, w)\bigg)
\end{align}
which together with \eqref{eq:computeDivergenceofBB:intermsofBulk:000000} implies
\bea
\lab{eq:computeDivergenceofBB:intermsofBulk}
 \D^\mu{\PP}_{\mu,\tau_2}[\pmb\psi](X, w) &=& \sum_{n=1}^3\Bulkxw{n}+\Err_{\text{l.o.t.}},
\eea
where we used in \eqref{eq:computeDivergenceofBB:intermsofBulk}:
\begin{itemize} 
\item the first estimate in \eqref{eq:difference:PPmupmbpsiandpsiij} and the support properties of $\chi_2'(\tau)$, $\chi_3'(\tau)$ and $\chi_4'(\tau)$ in \eqref{def:cutoffsintime1234} for the first two terms on the RHS of \eqref{eq:moredetialsontermsofdivergenceDDmuPPmupmspsionwhichderivativefallsoncutoffswidetildechinn123},

\item the last estimate in \eqref{eq:difference:PPmupmbpsiandpsiij} and the support properties and estimates of $\DD(\chi_{1,*})$ and $\DD(\chi_{3,*})$ in \eqref{def:cutoffsintime123*} for the last two terms on the RHS of \eqref{eq:moredetialsontermsofdivergenceDDmuPPmupmspsionwhichderivativefallsoncutoffswidetildechinn123}.
\end{itemize}

Next, we estimate the terms $\Bulkxw{n}$, $n=1,2,3$. The estimates for $\Bulkxw{2}$ and $\Bulkxw{3}$ are identical respectively to the ones for $\Bulkxw{2}$ and $\Bulkxw{4}$ in \cite{MaSz26} so that (8.141) (8.142) in \cite{MaSz26} holds, i.e.,
\begin{align}
\lab{eq:Bulkxw:general:proof:part2}
\Bulkxw{2}
={}&\sum_{i,j}\Re\bigg(\ov{\widetilde{\chi}_2\bigg(X(\psi_{ij}) - \big(X^{\a}M_{i\a}^k\psi_{kj}+ X^{\a}M_{j\a}^k\psi_{ik}\big) +\frac{1}{2}w\psi_{ij}\bigg)}F_{ij}\bigg)
\nn\\
&+\widetilde{\chi}_2\frac 1 2 \Big(\QQ[\pmb\psi]  \c\piX +w \LL[\pmb\psi] -X( V ) |\pmb\psi|^2-\frac{1}{2}|\pmb\psi|^2   \square_\g  w\Big)+\Err_{\text{l.o.t.}},
\end{align}
and
\begin{align}
\lab{eq:Bulkxw:general:proof:part4}
\Bulkxw{3}
={}&\sum_{i,j}\Re\bigg(\widetilde{\chi}_3\ov{\bigg(X(\psi_{ij}) - \big(X^{\a}(M_K)_{i\a}^k\psi_{kj}+ X^{\a}(M_K)_{j\a}^k\psi_{ik}\big) +\frac{1}{2}w\psi_{ij}\bigg)}F_{ij}\bigg)
\nn\\
&+\widetilde{\chi}_3\frac 1 2 \Big(\QQ[\pmb\psi]  \c\piX +w \LL[\pmb\psi] -X( V) |\pmb\psi|^2-\frac{1}{2}|\pmb\psi|^2   \square_\g  w\Big)_{K}+\Err_{\text{l.o.t.}} .
\end{align} 

Next, we consider the term $\Bulkxw{1}$. Using \eqref{def:cutoffsintime123:widetilde}, we decompose it as follows 
\beaa
\Bulkxw{1} &=& \Bulkxw{1,1}+\Bulkxw{1,2},\\
\Bulkxw{1,1} &:=& \Big(\chi_1+\chi_3(\chi_{1,*}+\chi_{2,*})\Big)\sum_{i,j}\D^\mu{\PP}_\mu[\psi_{ij}](X,w),\\
\Bulkxw{1,2} &:=& \chi_2\chi_{2,*}\sum_{i,j}\D^\mu{\PP}_\mu[\psi_{ij}](X,w).
\eeaa
As $\chi_1+\chi_3(\chi_{1,*}+\chi_{2,*})$ is supported in $\MMh(-\infty, \tau_1+2)\cup\MMh(\tau_2-3, \tau_2+1)$, the estimate for $\Bulkxw{1,1}$ is identical to the sum of the ones for $\Bulkxw{1}$ and $\Bulkxw{3}$ in \cite{MaSz26} so that we obtain in view of (8.144) and (8.145) in \cite{MaSz26}
\begin{align}
\lab{eq:Bulkxw:general:proof:part3:00}
\Bulkxw{1,1}
=& \Big(\chi_1+\chi_3(\chi_{1,*}+\chi_{2,*})\Big)\frac{1}{2}\sum_{i,j}\Big(\QQ[\psi_{ij}]  \c\piX +w \LL[\psi_{ij}]  -X( V) |\psi_{ij}|^2-\frac{1}{2}|\psi_{ij}|^2   \square_\g  w\Big)
\nn\\
&+\sum_{i,j}\Re\bigg(F_{ij}\Big(\chi_1+\chi_3(\chi_{1,*}+\chi_{2,*})\Big)\ov{\bigg(X(\psi_{ij}) + \frac{1}{2}w\psi_{ij}  \bigg)}\bigg) +\Err_{\text{l.o.t.}}.
\end{align}

It remains to estimate the term $\Bulkxw{1,2}$ involving the cut-off $\chi_2\chi_{2,*}$ which is supported in $\RR_*(\tau_1+1,\tau_2-2)$. To this end, analogously to (8.115) (8.116) in \cite{MaSz26}, we rewrite the coupled system of wave equations \eqref{eq:ScalarizedWaveeq:general:Kerrpert} as follows
\bea
\lab{eq:scalarizedextendedeqs:psiij:moraandener}
\square_{\g}(\psi_{ij})-D_0|q|^{-2}\psi_{ij}=F_{ij} +\sum_{n=1}^4G_{n,ij} ,
\eea
where 
\bsub
\lab{def:Gnij:energyesti:withpotential}
\bea
G_{1,ij}&=& S_K(\psi)_{ij} ,\\
 G_{2,ij}&=&\frac{4ia\cos\th}{|q|^2} \pr_{\tt}\psi_{ij},\\
  G_{3,ij}&=&\chi_{\tau_1, \tau_2}\chi_{r_*}(\widehat{Q}\psi)_{ij}+(1-\chi_{\tau_1, \tau_2}\chi_{r_*})\big((\widehat{Q}_K\psi)_{ij}+f_{D_0}\psi_{ij}\big) ,\\
    G_{4,ij}&=& \chi_{\tau_1, \tau_2}\chi_{r_*}\big({S}(\psi)_{ij}-{S}_K(\psi)_{ij}\big),
\eea
\esub
Then, applying \eqref{eq:DivofPPmu:tensor:RW:prop} to the wave equations \eqref{eq:scalarizedextendedeqs:psiij:moraandener}, multiplying on both sides by $\chi_2\chi_{2,*}$, and summing over $i,j$, we deduce, as in (8.143) in \cite{MaSz26},
\bea
\lab{eq:energyidentity:generalvectorfieldX:witherrortermsexplicit}
\Bulkxw{1,2} 
&=&\sum_{i,j}\chi_2\chi_{2,*}\D^\mu{\PP}_\mu[\psi_{ij}](X,w)\nn\\
&=&\sum_{i,j}\Re\bigg(F_{ij}\chi_2\chi_{2,*}\ov{\bigg(X(\psi_{ij}) + \frac{1}{2}w\psi_{ij}  \bigg)}\bigg) \nn\\
&&
+ \chi_2\chi_{2,*}\frac{1}{2}\sum_{i,j}\Big(\QQ[\psi_{ij}]  \c\piX +w \LL[\psi_{ij}]  -X( V) |\psi_{ij}|^2-\frac{1}{2}|\psi_{ij}|^2   \square_\g  w\Big) \nn \\
&&+\sum_{i,j}\sum_{n=1}^4\Re\bigg(\chi_2\chi_{2,*}G_{n,ij}\ov{\bigg(X(\psi_{ij}) + \frac{1}{2}w\psi_{ij}  \bigg)}\bigg).
\eea
It remains to  estimate the term in the last line of  the above equation \eqref{eq:energyidentity:generalvectorfieldX:witherrortermsexplicit}. In view of the estimates \eqref{eq:assumptionsonregulartripletinperturbationsofKerr:0}  \eqref{estimates:Mialphaj:Kerrperturbations} for ${M_{i\a}^j}$, the form of ${S}(\psi)_{ij}$ and $({Q}\psi)_{ij}$ in \eqref{SandV}, and the form of $(\widehat{Q}\psi)_{ij}$ in \eqref{hatSandV:generalwave:Kerrpert}, we have
\bea\lab{estimates:widehatSandQ:perturbationsofKerr}
\bsplit
& S(\psi)_{ij}= O(r^{-2})\pr\psi+O(r^{-1})\nab\psi +\Ga_g\dk\psi, \qquad (\widehat{Q}\psi)_{ij}=O(r^{-2})\psi +\dk^{\leq 1}\Ga_g \psi, \\ 
& S(\psi)_{ij}-S_K(\psi)_{ij}= \Ga_g\dk\psi.
\end{split}
 \eea
Then, using \eqref{estimates:widehatSandQ:perturbationsofKerr} and the fact that $f_{D_0}=O(r^{-2})$ and $w=O(r^{-1})$, we have
\beaa
\sum_{i,j}\sum_{n=1}^4\Re\big(\chi_2\chi_{2,*}G_{n,ij}\ov{w\psi_{ij}  }\big)=\chi_2\chi_{2,*}O(r^{-2})|\psi|\big(|\pr\psi|+r^{-1}|\psi|\big),
\eeaa
and further, in view of the assumption \eqref{eq:assump:generalvectorfieldX:EMnearinf} for the vectorfield $X$, we have
\beaa
\sum_{i,j}\sum_{n=3}^4\Re\big(\chi_2\chi_{2,*}G_{n,ij}\ov{X\psi_{ij} }\big)=\chi_2\chi_{2,*}\Big(O(r^{-2})|\psi| + |\Ga_g||\dk\psi|\Big)\big(|\pr\psi|+r^{-1}|\psi|\big).
\eeaa

Next, we recall the following differential identities, see (8.43) (8.44) in \cite{MaSz26},
\bea
\lab{eq:integrationbypartsforitimesfirstordertimesfirstorder:general}
2\Re(i f \pr_{\a} \psi \ov{\pr_{\b}\psi}) &=& \pr_{\a}\big(\Re(i f \psi \ov{\pr_{\b}\psi})\big)-\pr_{\b}\big(\Re(i f \psi \ov{\pr_{\a}\psi})\big) \nn\\
&-& \Re(i \pr_{\a}f \psi \ov{\pr_{\b}\psi})+\Re(i \pr_{\b}f \psi \ov{\pr_{\a}\psi}), \text{for any real-valued function } f,
\eea
\bea
\lab{eq:integrationbypartsforitimesfirstordertimesfirstorder:general:caseantisymmatrix}
2\Re(A^{ij}\pr_{\a} \psi_i \ov{\pr_{\b}\psi_j}) &=& \pr_{\a}\big(\Re(A^{ij}\psi_i \ov{\pr_{\b}\psi_j})\big)-\pr_{\b}\big(\Re(A^{ij} \psi_i \ov{\pr_{\a}\psi_j})\big) \nn\\
&-& \Re(\pr_{\a}(A^{ij}) \psi_i \ov{\pr_{\b}\psi_j})+\Re(\pr_{\b}(A^{ij})\psi_i\ov{\pr_{\a}\psi_j}), \,\,\,\,\text{for any antisymmetric}\nn\\
&&\qquad\qquad\qquad\qquad\qquad\quad\textrm{family of real-valued functions }A^{ij}.
\eea
For the term with $G_{2,ij}=\frac{4ia\cos\th}{|q|^2} \pr_{\tt}\psi_{ij}$, we integrate the differential identity \eqref{eq:integrationbypartsforitimesfirstordertimesfirstorder:general} with the choice $f=\chi_3\frac{2a\cos\th}{|q|^2}\sqrt{|\textrm{det}(\gh)|}X^\a$ and $x^\b=\tau$ which yields, for any $10m\leq r_1\ll r_*$, 
\beaa
\bigg|\int_{\MM_{r_1, +\infty}(\tau',\tau'')}\Re\big(\chi_2\chi_{2,*}G_{2,ij}\ov{X\psi_{ij} }\big)\bigg| &\les& \int_{\RR_*(\tau_*', \tau_*'')}r^{-2} |\psi| |\pr_{\tt}\psi|+
\sup_{\tau\in [\tau',\tau'']}\int_{\Si(\tau)}r^{-2} |\psi| |\pr^{\leq 1}\psi|\\
&\les& \NNtlocal[\pmb \psi](\tau',\tau''),
\eeaa
where we used in the last inequality the following estimate
\begin{align}\lab{eq:maininequalityrminus2psiprpsiplusrminuspsionRR*}
\int_{\RR_*(\tau_*', \tau_*'')}r^{-2}|\psi|\big(|\pr\psi|+r^{-1}|\psi|\big) \les& \int_{\frac{c_*}{2}}^{c_*}\int_{\Si_{*,c}(\tau_*', \tau_*'')}r^{-2}|\psi|\big(|\pr\psi|+r^{-1}|\psi|\big)\nn\\
\les& \sup_{\frac{c_*}{2}\leq c\leq c_*}\int_{\Si_{*,c}(\tau_*', \tau_*'')}r^{-1}|\psi|\big(|\pr\psi|+r^{-1}|\psi|\big)\nn\\
\les& \sup_{\frac{c_*}{2}\leq c\leq c_*}\bigg[\bigg(\int_{\Si_{*,c}(\tau_*', \tau_*'')}r^{-2}|\psi|^2\bigg)^{\frac{1}{2}}\Big(\F_{\Si_{*,c}}[\pmb\psi](\tau_*', \tau_*'')\Big)^{\frac{1}{2}}\bigg].
\end{align}

For the term with $G_{1,ij}$, we proceed as for the proof of (8.121) in \cite{MaSz26}. In particular, we estimate only the integral of the part involving $(M_K)_{i}^{k\a}\pr_\a(\psi_{kj})$, the control of the integral of the other part involving $(M_K)_{j}^{k\a}\pr_\a(\psi_{ik})$ being estimated in the exact same way.  
Using \eqref{formula:symmetricpartofMmatrices}, we have
\beaa
(M_K)_{i}^{k\a}&=&(M_{K,S})_{i}^{k\a} +(M_{K,A})_{i}^{k\a}= -\frac{1}{2}\pr^{\a} (x^i x^k) +(M_{K,A})_{i}^{k\a}\nn\\
&=& -\frac{1}{2}\gam^{\a\b}\pr_{\b}(x^i x^k) +\gam^{\a\b}(M_{K,A})_{i\b}^{k},
\eeaa
and in view of the estimate for ${(M_K)_{i\a}^j}$ in Lemma \ref{lemma:computationoftheMialphajinKerr}, we obtain the following:
\begin{itemize}
\item By integrating the differential identity \eqref{eq:integrationbypartsforitimesfirstordertimesfirstorder:general:caseantisymmatrix} with $A_i^k=\chi_2\chi_{2,*}\gam^{\a\b}(M_{K,A})_{i\b}^{k}X^\mu$, we deduce that the integral involving $\gam^{\a\b}(M_{K,A})_{i\b}^{k}$ is bounded by\footnote{The weight $r^{-2}$ in the term $r^{-2}|\psi||\pr\psi|$ integrated on $\RR_*(\tau', \tau'')$ comes from \eqref{eq:assump:generalvectorfieldX:EMnearinf} and the fact that $\pr_\tau(A_i^k)=0$ while $(\pr_r(A_i^k), r^{-1}\pr_{x^a})=O(r^{-2})$.}
\beaa
\sup_{\tau\in[\tau',\tau'']}\int_{\Sigma(\tau)} r^{-2} |\psi| |\dk^{\leq 1}\psi| +\int_{\RR_*(\tau_*', \tau_*'')}r^{-2}|\psi||\pr\psi|\les\NNtlocal[\pmb \psi](\tau',\tau'')
\eeaa
in view of the estimates \eqref{eq:maininequalityrminus2psiprpsiplusrminuspsionRR*}.

\item Since we have 
\begin{align*}
\gam^{\a\b}\pr_\b(x^ix^k)\pr_\a(\psi_{kj})\ov{\pr_\tau(\psi_{ij})} 
={}&\gam^{\a\b}\pr_\b(x^k)\pr_\a(\psi_{kj})\ov{\pr_\tau(x^i\psi_{ij})}+\gam^{\a\b}\pr_\b(x^i)\pr_\a(x^k\psi_{ik})\ov{\pr_\tau(\psi_{ij})}\\
& - \gam^{\a\b}\pr_\b(x^i)\pr_\a(x^k)\psi_{ik}\ov{\pr_\tau(\psi_{ij})},
\end{align*}
and in view of the fact that
\beaa
\gam^{\a\b}\pr_\b(x^k)\pr_\a=r^{-1}\pr, \quad \gam^{\a\b}\pr_\b(x^i)\pr_\a(x^k)=O(r^{-2}),
\eeaa
the integral involving $-\frac{1}{2}\gam^{\a\b}\pr_\b (x^i x^k)$ is, in view of \eqref{def:Errdefectofpsi}, \eqref{eq:maininequalityrminus2psiprpsiplusrminuspsionRR*} and the support properties of $\chi_2\chi_{2,*}$, bounded by 
\beaa
\sup_{\frac{c_*}{2}\leq c\leq c_*}\bigg[\Big(\F_{\Si_{*,c}}[\pmb\psi](\tau_*', \tau_*'')\Big)^{\frac{1}{2}}\bigg(\Errdefects[\pmb\psi]+\int_{\Si_{*,c}(\tau_*', \tau_*'')}r^{-2}|\psi|^2\bigg)^{\frac{1}{2}}\bigg]\les\NNtlocal[\pmb \psi](\tau',\tau'').
\eeaa
\end{itemize}
Therefore, we deduce,
\beaa
\bigg|\int_{\MM_{r_1, +\infty}(\tau',\tau'')}\Re\big(\chi_2\chi_{2,*}G_{1,ij}\ov{X\psi_{ij}}\big)\bigg|
&\les& \NNtlocal[\pmb \psi](\tau',\tau'').
\eeaa

In view of the above estimates, we infer for any $r_1 \geq 10m$, 
\beaa
&&\bigg|\int_{\MM_{r_1, +\infty}(\tau',\tau'')}  \sum_{i,j}\sum_{n=1}^4\Re\bigg(\chi_2\chi_{2,*}G_{n,ij}\ov{\bigg(X(\psi_{ij}) + \frac{1}{2}w\psi_{ij}  \bigg)}\bigg)\bigg|\nn\\
&\les&\ep\sup_{\frac{c_*}{2}\leq c\leq c_*}\F_{\Si_{*,c}}[\pmb\psi](\tau_*', \tau_*'') +  \NNtlocal[\pmb \psi](\tau',\tau'')
\eeaa
which implies
\beaa
\sum_{i,j}\sum_{n=1}^4\Re\bigg(\chi_2\chi_{2,*}G_{n,ij}\ov{\bigg(X(\psi_{ij}) + \frac{1}{2}w\psi_{ij}  \bigg)}\bigg)=\Err_{\text{l.o.t.}}
\eeaa
and hence, together with \eqref{eq:energyidentity:generalvectorfieldX:witherrortermsexplicit},
\bea
\lab{eq:Bulkxw:general:proof:part3}
\Bulkxw{1,2}
&=&\chi_2\chi_{2,*}\frac{1}{2}\sum_{i,j}\Big(\QQ[\psi_{ij}]  \c\piX +w \LL[\psi_{ij}]  -X( V) |\psi_{ij}|^2-\frac{1}{2}|\psi_{ij}|^2   \square_\g  w\Big)
\nn\\
&&
+\sum_{i,j}\Re\bigg(F_{ij}\chi_2\chi_{2,*}\ov{\bigg(X(\psi_{ij}) + \frac{1}{2}w\psi_{ij}  \bigg)}\bigg) +\Err_{\text{l.o.t.}} .
\eea

Plugging  the estimates \eqref{eq:Bulkxw:general:proof:part2}, \eqref{eq:Bulkxw:general:proof:part4}, \eqref{eq:Bulkxw:general:proof:part3:00} and \eqref{eq:Bulkxw:general:proof:part3} into \eqref{eq:computeDivergenceofBB:intermsofBulk} and using $\sum_{n=1}^3\widetilde{\chi}_n=1$ and $\widetilde{\chi}_1=\chi_1+\chi_3(\chi_{1,*}+\chi_{2,*})+\chi_2\chi_{2,*}$, we infer the identity  \eqref{eq:DivwidetilePPmu:extendedscalarizedRW} with $X_{\tau_2}$ defined as in \eqref{def:Xtau2psiij:generalforEM} and with $\Err_{\text{l.o.t.}}$ satisfying \eqref{def:Errlowerorderterm:DivPPmuchi}, which then concludes the proof of Proposition \ref{prop-app:stadard-comp-Psi:extendedscalarizedRW}.
\end{proof}

The following proposition is the analog of Proposition 8.32 in \cite{MaSz26}. It provides energy and Morawetz estimates for the globally extended coupled system of wave equations \eqref{eq:ScalarizedWaveeq:general:Kerrpert} in a large radius region away from the trapping region. 

\begin{proposition}
\lab{prop:energymorawetznearinfinity:extendedRWsystem:kerrpert}
Let $\psi_{ij}$ be a solution to the system of wave equations \eqref{eq:ScalarizedWaveeq:general:Kerrpert}, and let $X_{\tau_2}(\psi)_{ij}$ and $\NNtlocal[\pmb \psi](\tau',\tau'')$ be given as in \eqref{def:Xtau2psiij:generalforEM} and 
\eqref{def:NNtlocalinr:NNtEner:NNtaux:wavesystem:EMF:proof}, respectively. 
 Then, under the assumptions of Theorem \ref{th:mainenergymorawetzmicrolocal} for the scalars $\psi_{ij}$, $F_{ij}$ and the spacetime $(\MM,\g)$:
 \begin{itemize}
 \item there exists a constant $c>0$ such that the following energy estimate holds, for  any  constant $r_1\geq 10m$,
\bea
\label{eq:energynearinf:extendedsystem:Kerrpert:1}
&&c\Big(\Eh_{r\geq r_1}[\pmb\psi](\tau'')
+\Fh_{\II_+}[\pmb\psi](\tau',\tau'')\Big)
-\int_{\Hh_{r_1}(\tau',\tau'')}\sum_{i,j}\PP_\mu[\psi_{ij}](\pr_{\tt}, 0)N_{\Hh_r}^{\mu}\nn\\
&&+\int_{\MMh_{r\geq r_1}(\tau',\tau'')}\sum_{i,j}\Re\big(F_{ij}\ov{T_{\tau_2}(\psi)_{ij}}\big) 
\nn \\
&\les& \sup_{\tau\in[\tmic,\tau_1+2]}\Eh[\pmb \psi](\tau)+\Eh_{r\geq r_1}[\pmb\psi](\tau') +\int_{\Hh_{r_1}(\tau',\tau'')} (r_1)^{-1}|\psi| |\pr^{\leq 1}\psi|  \nn\\
&&+ \NNtlocal[\pmb \psi](\tau',\tau'') + \ep\EMh[\pmb\psi](\tau',\tau'')+\ep\sup_{\frac{c_*}{2}\leq c\leq c_*}\F_{\Si_{*,c}}[\pmb\psi](\tau_*', \tau_*''),
\eea
where $(\tau_*', \tau_*'')$ is given by \eqref{eq:definitionoftau*primetau*doubleprimefromtauprimtaudoubleprime}, and where $T_{\tau_2}(\psi)_{ij}$ is given by \eqref{def:Xtau2psiij:generalforEM} in the particular case where $X=\pr_\tau$, i.e.,
\bea
 \lab{def:widehatTpsiij:equalsnabTpsiij}
T_{\tau_2}(\psi)_{ij} &:=& \pr_{\tt}(\psi_{ij})  -\widetilde{\chi}_2\big(M_{i\tt}^k\psi_{kj} + M_{j\tt}^k\psi_{ik} + 2i\widetilde{w}\psi_{ij}\big)\nn\\
&& - \widetilde{\chi}_3\big((M_K)_{i\tt}^k\psi_{kj} + (M_K)_{j\tt}^k\psi_{ik} + 2i\widetilde{w}\psi_{ij}\big);
\eea

 \item there exists a constant $c>0$ such that the following energy estimate holds, for  any  constant $r_1\geq 11m$,
\bea
\label{eq:energynearinf:extendedsystem:Kerrpert:2:cutoff}
&&c\Big(\Eh_{r\geq r_1}[\pmb\psi](\tau'')
+\Fh_{\II_+}[\pmb\psi](\tau',\tau'')\Big)
+\int_{\MMh_{r\geq r_1}(\tau',\tau'')}\sum_{i,j}\Re\big(F_{ij}\ov{T_{\tau_2}(\psi)_{ij}}\big) 
\nn \\
&\les& \sup_{\tau\in[\tmic,\tau_1+2]}\Eh[\pmb \psi](\tau)+\Eh_{r\geq r_1-m}[\pmb\psi](\tau')  + \NNtlocal[\pmb \psi](\tau',\tau'') + \ep\EMh[\pmb\psi](\tau',\tau'') \nn\\
&&+\ep\sup_{\frac{c_*}{2}\leq c\leq c_*}\F_{\Si_{*,c}}[\pmb\psi](\tau_*', \tau_*'')
+\Mh_{r_1-m, r_1}[\pmb\psi](\tau',\tau'')
+\int_{\MMh_{r_1-m, r_1}(\tau',\tau'')}|\pmb F|^2;
\eea

\item there exists a constant $c>0$ such that the following energy estimate holds, for  any  constant $r_1\geq 11m$ and any $c_*/2\leq c_0\leq c_*$, assuming in addition that $\tau_1+1\leq\tau'<\tau''\leq\tau_2-2$,
\begin{align}
\label{eq:energynearinf:extendedsystem:Kerrpert:2:cutoff:stopinRR*onSi*c0}
&c\Big(\E_{r_1\leq r\leq\varsigma_*(\tau-c_0)}[\pmb\psi](\tau'')
+\F_{\Si_{*,c_0}}[\pmb\psi](\tau',\tau'')\Big)
+\int_{\MM_{r_1\leq r\leq\varsigma_*(\tau-c_0)}(\tau',\tau'')}\sum_{i,j}\Re\big(F_{ij}\ov{T_{\tau_2}(\psi)_{ij}}\big) 
\nn \\
\les& \sup_{\tau\in[\tmic,\tau_1+2]}\Eh[\pmb \psi](\tau)+\Eh_{r\geq r_1-m}[\pmb\psi](\tau')  + \NNtlocal[\pmb \psi](\tau',\tau'') + \ep\EMh[\pmb\psi](\tau',\tau'')\nn\\
& +\ep\sup_{\frac{c_*}{2}\leq c\leq c_*}\F_{\Si_{*,c}}[\pmb\psi](\tau', \tau'')  +\Mh_{r_1-m, r_1}[\pmb\psi](\tau',\tau'')
+\int_{\MMh_{r_1-m, r_1}(\tau',\tau'')}|\pmb F|^2;
\end{align}

\item there exists  a constant $c>0$ such that the following Morawetz estimate holds for  a suitably large constant $R_1\gg 12m$:
\bea
\label{eq:Morawetznearinf:extendedsystem:Kerrpert:1:00}
\nn&& c\M_{r\geq R_1}[\pmb\psi](\tau',\tau'') -\int_{H_{R_1}(\tau',\tau'')}\sum_{i,j}\PP_\mu[\psi_{ij}](X_1, w_1)N_{H_r}^{\mu}\\
&&+\int_{\MM_{r\geq R_1}(\tau',\tau'')}\sum_{i,j}\Re\left(F_{ij}\ov{\left((X_{1})_{\tau_2}(\psi)_{ij} +\frac{1}{2}w_1\psi_{ij}\right)}\right) 
\nn \\
&\les&\sup_{\tau\in[\tmic,\tau_1+2]}\E[\pmb \psi](\tau)+\E_{r\geq R_1}[\pmb\psi](\tau'')
+\E_{r\geq R_1}[\pmb\psi](\tau') 
+\F_{\II_+}[\pmb\psi](\tau',\tau'')
 \nn\\
&&
+\ep\sup_{\tt\in[\tau',\tau'']}\E[\pmb\psi](\tau) +\ep\sup_{\frac{c_*}{2}\leq c\leq c_*}\F_{\Si_{*,c}}[\pmb\psi](\tau_*', \tau_*'')+\int_{H_{R_1}(\tau',\tau'')}(R_1)^{-1}|\pmb\psi| |\pr^{\leq 1}\pmb\psi|\nn\\
&&+ \NNtlocal[\pmb \psi](\tau',\tau''),
\eea
where\footnote{Note that the value of the function $w$ here is twice the value of the function $w$ chosen in \cite[Lemma 3.10]{MaSz24} which is due to a different normalization in the definition of the current $\PP_{\mu}[\pmb\psi](X,w)$. This is also the case for the choice of the function $w_{\de}$ in \eqref{def:Xandw:improvedMorawetz:extendedRW:Kerrpert}.}
\beaa
X_1=2\mu(1-mr^{-1})\pr_{r}^{\text{BL}},\qquad w_1=4\mu r^{-1}(1-mr^{-1});
\eeaa

\item for any $\de\in(0,1]$, there exists  a constant $c>0$ such that the following Morawetz estimate holds for  a suitably large constant $R_1\gg 12m$:
\begin{align}
\label{eq:Morawetznearinf:extendedRWsystem:Kerrpert:delta}
\nn& c\de\M_{\de, r\geq R_1}[\pmb\psi](\tau',\tau'')+\int_{\MM_{r\geq R_1-m}(\tau',\tau'')}\sum_{i,j}\Re\left(F_{ij}\ov{\left((X_{\de})_{\tau_2}(\psi)_{ij} +\frac{1}{2}w_{\de}\psi_{ij}\right)}\right) 
\nn \\
\les&\sup_{\tau\in[\tmic,\tau_1+2]}\E[\pmb \psi](\tau)+\E_{r\geq R_1-m}[\pmb\psi](\tau'')
+\E_{r\geq R_1-m}[\pmb\psi](\tau') 
+\F_{\II_+}[\pmb\psi](\tau',\tau'')
\nn\\
&+\ep\sup_{\tt\in[\tau',\tau'']}\E[\pmb\psi](\tau) +\ep\sup_{\frac{c_*}{2}\leq c\leq c_*}\F_{\Si_{*,c}}[\pmb\psi](\tau_*', \tau_*'') +\M_{R_1-m,R_1}[\pmb\psi](\tau',\tau'')+ \NNtlocal[\pmb \psi] (\tau',\tau''),
\end{align}
with
\bea
\lab{def:Xandw:improvedMorawetz:extendedRW:Kerrpert}
X_{\de}=2\mu f_{\de}{\pr}_r^{\text{BL}}, \quad w_{\de}=4\mu h_{\de}, \quad f_{\de}=\chi_{R_1}(1-m^{\de}r^{-\de}),\quad 
h_{\de}=\chi_{R_1} r^{-1} (1-m^{\de}r^{-\de}),
\eea
where $\chi_{R_1}=\chi_{R_1}(r)$ is a smooth cutoff function that equals $1$ for $r\geq R_1$ and vanishes for $r\leq R_1-m$.
 \end{itemize}
\end{proposition}

\begin{proof}
The proof of the estimates \eqref{eq:energynearinf:extendedsystem:Kerrpert:1}, \eqref{eq:energynearinf:extendedsystem:Kerrpert:2:cutoff}, \eqref{eq:Morawetznearinf:extendedsystem:Kerrpert:1:00} and \eqref{eq:Morawetznearinf:extendedRWsystem:Kerrpert:delta} is an immediate adaptation of the proof of the corresponding estimates in Proposition 8.32 of \cite{MaSz26}, and relies on the integration on $\MMh_{r_1, +\infty}(\tau',\tau'')$, where $r_1\geq 10m$, of the divergence identity \eqref{eq:DivwidetilePPmu:extendedscalarizedRW} for various choices of multipliers $(X, w)$, using \eqref{eq:1formPP_mu:equality} \eqref{eq:1formBB_mu:equality:errortermestimates} to control the boundary terms and \eqref{def:Errlowerorderterm:DivPPmuchi} to control the error terms in the bulk. Finally, \eqref{eq:energynearinf:extendedsystem:Kerrpert:2:cutoff:stopinRR*onSi*c0} is an immediate extension of \eqref{eq:energynearinf:extendedsystem:Kerrpert:2:cutoff} where the integration of the divergence identity \eqref{eq:DivwidetilePPmu:extendedscalarizedRW} is done with the same multiplier but on the spacetime region $\MM_{r_1\leq r\leq\varsigma_*(\tau-c_0)}(\tau',\tau'')$ with $\tau_1+1\leq\tau'<\tau''\leq\tau_2-2$.
\end{proof}

We are now in position to prove Theorem \ref{th:mainenergymorawetzmicrolocal}.

\begin{proof}[Proof of Theorem \ref{th:mainenergymorawetzmicrolocal}]
We proceed as in section 8.5 of \cite{MaSz26}:
\begin{enumerate}
\item In the region $r\geq\Rmic$, with $\Rmic$ as in Remark \ref{rmk:choiceofconstantRbymeanvalue}, we multiply the energy estimate \eqref{eq:energynearinf:extendedsystem:Kerrpert:1} by a large constant $A$ and add it to the Morawetz estimate \eqref{eq:Morawetznearinf:extendedsystem:Kerrpert:1:00}. We deduce the following energy-Morawetz estimate near infinity for solutions to \eqref{eq:ScalarizedWaveeq:general:Kerrpert}, for any $\tmic\leq \tau'<\tau''$,
\begin{align}
\lab{eq:largerMoraesti:scalarizedwave:withFij}
&c\MFh_{r\geq {\Rmic}}[\pmb\psi ](\tau', \tau'')+c\Eh_{r\geq {\Rmic}}[\pmb\psi](\tau'')+\sum_{i,j}\textbf{BDR}[\psi_{ij} ]\vert_{\Hh_{\Rmic}(\tau',\tau'')}\nn\\
&
+\sum_{i,j}\int_{\MMh_{{\Rmic},+\infty}(\tau',\tau'')}\Re\Big(F_{ij}\ov{ ({X_1} +A{\pr_{\tt}}+ w_1)\psi_{ij}}\Big)\nn\\
&
+\sum_{i,j}\int_{\MMh_{{\Rmic},+\infty}(\tau',\tau'')}\Re\bigg(F_{ij}\ov{ \Big( (X_1)_{\tau_2}(\psi)_{ij} - X_1 (\psi_{ij} ) + AT_{\tau_2}(\psi)_{ij} - A \pr_{\tau}(\psi_{ij})\Big)}\bigg)\nn\\
\les&\Eh[\pmb\psi](\tau')+\sup_{\tau\in[\tmic,\tau_1+2]}\Eh[\pmb \psi](\tau)+\ep\sup_{\tau\in[\tau_2-3,\tau_2+1]}\Eh[\pmb\psi](\tau)+\ep\EMh[\pmb\psi](\tau',\tau'')\nn\\
&+\ep\sup_{\frac{c_*}{2}\leq c\leq c_*}\F_{\Si_{*,c}}[\pmb\psi](\tau', \tau'') + \NNtlocal[\pmb \psi](\tau',\tau'')
+\int_{\Hh_{{\Rmic}}(\tau',\tau'')}(R_0)^{-1}|\pmb\psi| |\pr^{\leq 1}\pmb\psi|.
\end{align}

\item Then, we take $\tau'=\tmic$ and $\tau''\to +\infty$ in \eqref{eq:largerMoraesti:scalarizedwave:withFij} and sum the resulting estimate on $\MMh_{r\geq\Rmic}(\Iti)$ with the microlocal Morawetz estimate (8.129) in \cite{MaSz26} on $\MMh_{r_+(1+\dhor),\Rmic}(\Iti)$ which holds true in our setting. In particular, as in \cite{MaSz26}, the boundary terms at $r=\Rmic$, which are identical in our case, cancel to each other up to lower order terms. 

\item Next, following the proof of Proposition 8.34 in \cite{MaSz26}, we add the redshift estimate (8.156) in \cite{MaSz26} which holds true in our setting, and we obtain the analog of the global microlocal Morawetz estimate (8.159) in \cite{MaSz26} on $\MMh(\Iti)$.  

\item Next, we control the energy in $r\geq 11m$ by using the energy estimate \eqref{eq:energynearinf:extendedsystem:Kerrpert:2:cutoff} with $r_1=11m$ and $\tau'=\tmic$, and taking the supremum in $\tau''$ which yields the following analog of (8.178) in \cite{MaSz26}
\bea
\lab{eq:energyestiinrgeqRregion:extendedRWsystem:Kerrpert:proof}
\nn\EFh_{r\geq 11m}[\pmb\psi](\Iti) &\les&  \sup_{\tau\in[\tmic,\tau_1+2]}\Eh[\pmb \psi](\tau)+\Mh[\pmb \psi](\Iti) +  \NNtlocal[\pmb \psi](\Iti)+\int_{\Mntraph}|\pmb F|^2\\
&&+\ep\sup_{\tau\in\Iti}\Eh[\pmb\psi](\tau)+\ep\sup_{\frac{c_*}{2}\leq c\leq c_*}\F_{\Si_{*,c}}[\pmb\psi](\tau_1+1, \tau_2-2)\nn\\
&&+\sum_{i,j}\sup_{\tau\geq\tmic}\bigg|\int_{\Mntraph(\tmic, \tau)}{\Re\Big(F_{ij}\ov{\pr_{\tau}\psi_{ij}}\Big)}\bigg|.
\eea

\item In view of \eqref{eq:energyestiinrgeqRregion:extendedRWsystem:Kerrpert:proof} and the control of the energy in the redshift region provided by the redshift estimate in (8.156) in \cite{MaSz26} which holds true in our setting, it remains to control the energy in $\MMh_{r_+(1+\dhor),11m}(\Iti)$. This is done in Proposition 8.37 of \cite{MaSz26} which holds true in our setting. 

\item In addition to the estimates in section 8 of \cite{MaSz26}, we apply \eqref{eq:energynearinf:extendedsystem:Kerrpert:2:cutoff:stopinRR*onSi*c0} with $r_1=11m$, $\tau'=\tau_1+1$, $\tau''=\tau_2-2$ and any $\frac{c_*}{2}\leq c_0\leq c_*$, and then take the supremum in $c_0$ to obtain, for $\ep>0$ small enough, 
\begin{align}
\label{eq:energynearinf:extendedsystem:Kerrpert:2:cutoff:stopinRR*onSi*c0:controlofFluxtermsinRR*inaddtiontoestimatesMaSz26}
& \sup_{\frac{c_*}{2}\leq c\leq c_*}\F_{\Si_{*,c}}[\pmb\psi](\tau_1+1, \tau_2-2)\nn \\
\les& \sup_{\tau\in[\tmic,\tau_1+2]}\Eh[\pmb \psi](\tau) + \NNtlocal[\pmb \psi](\tau_1+1,\tau_2-2) + \ep\EMh[\pmb\psi](\tau_1+1,\tau_2-2)\nn\\
&   +\Mh_{10m, 11m}[\pmb\psi](\tau',\tau'')
+\int_{\MMh_{10m, 11m}(\tau',\tau'')}|\pmb F|^2\nn\\
& +\sup_{\frac{c_*}{2}\leq c\leq c_*}\bigg|\int_{\MM_{11m\leq r\leq\varsigma_*(\tau-c)}(\tau_1+1,\tau_2-2)}\sum_{i,j}\Re\big(F_{ij}\ov{T_{\tau_2}(\psi)_{ij}}\big)\bigg|. 
\end{align}

\item Finally, as in section 8.5.3 in \cite{MaSz26}, we conclude the proof of the extension of the proof of \cite{MaSz26} in $\MMh_{r\geq R_0}$ to \eqref{eq:ScalarizedWaveeq:general:Kerrpert} by combing the above estimates and taking $\dhor$ and $\ep$ suitably small, and applying Cauchy-Schwarz to the term $\NNtlocal[\pmb \psi](\Iti)$ defined as in \eqref{def:NNtlocalinr:NNtEner:NNtaux:wavesystem:EMF:proof}. This yields the following analog of the main microlocal energy-Morawetz estimate (8.31) in \cite{MaSz26} in our setting
\beaa
&&\widetilde{\EMF}[\pmb \psi] +\sup_{\frac{c_*}{2}\leq c\leq c_*}\F_{\Si_{*,c}}[\pmb\psi](\tau_1+1, \tau_2-2)\nn\\
&\les& \sup_{\tau\in[\tmic,\tau_1+2]}\E[\pmb\psi](\tau)+
\Errdefect[\pmb\psi]+\Errdefects[\pmb\psi] +\A[\pmb \psi](\Iti)+\A_*[\pmb \psi]+\NNt[\pmb \psi,\pmb F]\nn\\
&&+\sup_{\frac{c_*}{2}\leq c\leq c_*}\bigg|\int_{\MM_{11m\leq r\leq\varsigma_*(\tau-c)}(\tau_1+1,\tau_2-2)}\sum_{i,j}\Re\big(F_{ij}\ov{T_{\tau_2}(\psi)_{ij}}\big)\bigg|.
\eeaa
as stated in \eqref{th:eq:mainenergymorawetzmicrolocal:tensorialwave:scalarized:eachpsisp}, where:
\begin{enumerate}
\item the additional control of $\sup_{\frac{c_*}{2}\leq c\leq c_*}\F_{\Si_{*,c}}[\pmb\psi](\tau_1+1, \tau_2-2)$ on the LHS of \eqref{th:eq:mainenergymorawetzmicrolocal:tensorialwave:scalarized:eachpsisp} compared to (8.31) in \cite{MaSz26} comes from the LHS of \eqref{eq:energynearinf:extendedsystem:Kerrpert:2:cutoff:stopinRR*onSi*c0:controlofFluxtermsinRR*inaddtiontoestimatesMaSz26},

\item the new terms on the RHS of \eqref{th:eq:mainenergymorawetzmicrolocal:tensorialwave:scalarized:eachpsisp} compared to (8.31) in \cite{MaSz26} are $\Errdefects[\pmb\psi]$ and $\A_*[\pmb \psi]$, as well as the last term which is due to the last term on the RHS of \eqref{eq:energynearinf:extendedsystem:Kerrpert:2:cutoff:stopinRR*onSi*c0:controlofFluxtermsinRR*inaddtiontoestimatesMaSz26}.
\end{enumerate}
\end{enumerate}
This concludes the proof of Theorem \ref{th:mainenergymorawetzmicrolocal}.
\end{proof}

%%%%%%%%%%%%%%%%%%%%%%%%%%%%%%%%%%%%%%%%

\subsubsection{EMF estimates for $\phis{p}$ and $\psis{p}$ at zeroth-order}
\lab{sec:proofofth:main:intermediary:0order}

%%%%%%%%%%%%%%%%%%%%%%%%%%%%%%%%%%%%%%%%

Let $(\MM, \g)$ satisfy the assumptions of Section \ref{sec:assumptionsforsec:energyMorawetzesitmatesforTeukoslkyonMM:upto15derivatives}. By a slight abuse of notation, we denote throughout this section the metric $\g_{\chi_{\tau_1, \tau_2},*}$ in \eqref{eq:extendedmetricgchitau1tau2} by $\g$. Also, let $\{\pmb\phi_s^{(p)}\}_{s=\pm 2, p=0,1,2}$ be a solution to the tensorial Teukolsky wave/transport systems \eqref{eq:TensorialTeuSysandlinearterms:rescaleRHScontaine2:general:Kerrperturbation}  \eqref{def:TensorialTeuScalars:wavesystem:Kerrperturbation}  in perturbations of Kerr, and let $\{\psis{p}\}_{s=\pm 2, p=0,1,2}$ be a solution to the following system of wave equations\footnote{In practice, $\{\psis{p}\}_{s=\pm 2, p=0,1,2}$ is a solution to \eqref{eq:waveeqwidetildepsi1} \eqref{def:tildef} so that $\widehat{F}_{total,s,ij}^{(p)}$ is  given by 
\bea\lab{eq:definitionofwidehatFtotalsij}
\widehat{F}_{total,s,ij}^{(p)}:=\chi_{\tau_1, \tau_2}^{(1)}\chi_{r_*}^{(1)}N^{(p)}_{W,s,ij}+\underline{F}^{(p)}_{s,ij} +\breve{F}^{(p)}_{s,ij},
\eea
with $\underline{F}^{(p)}_{s,ij}$ and $\breve{F}^{(p)}_{s,ij}$ given respectively by \eqref{def:tildef0} and  \eqref{def:breveFpsij}. In this section, we do not specify the form of $\widehat{F}_{total,s,ij}^{(p)}$ appearing in \eqref{eq:definitionofwidehatFtotalsij} as we need a version that is stable under commutation w.r.t. higher order derivatives in view of our applications in Section \ref{sec:proofofth:main:intermediary}.}
\begin{align}
\lab{eq:waveeqwidetildepsi1:gtilde:Teu}
\bigg({\square}_{{\g}} -\frac{4-2\de_{p0}}{|q|^{2}}\bigg)\psi^{(p)}_{s,ij}
={}&\chi_{\tau_1, \tau_2}\chi_{r_*}\big( \widehat{S}(\psi^{(p)}_s)_{ij} +(\widehat{Q}\psi^{(p)}_s)_{ij}\big) +(1-\chi_{\tau_1, \tau_2}\chi_{r_*})\big(\widehat{S}_K(\psi^{(p)}_s)_{ij} +(\widehat{Q}_K\psi^{(p)}_s)_{ij} \big)\nn\\
&+(1-\chi_{\tau_1, \tau_2}\chi_{r_*})f_p\psi^{(p)}_{s,ij}+\chi_{\tau_1, \tau_2}^{(1)}\chi_{r_*}^{(1)}L^{(p)}_{s,ij}+\widehat{F}_{total,s,ij}^{(p)} 
\end{align}
on $\MM$ satisfying \eqref{eq:causlityrelationsforwidetildepsi1}. 

The goal of this section is to initiate the proof of the EMF estimates for unweighted derivatives of $\phis{p}$ and $\psis{p}$ stated in Theorem \ref{th:main:intermediary} by first proving the EMF estimates for $\phis{p}$ and $\psis{p}$ at zeroth-order stated below in Theorem \ref{thm:EMF:systemofTeuscalarized:order0:final}. Relying on Theorem \ref{thm:EMF:systemofTeuscalarized:order0:final}, Theorem \ref{th:main:intermediary} will then be proved in Section \ref{sec:proofofth:main:intermediary}.

We define the following analogs of the notations in \eqref{def:EM-1norms:Reals}
\bsub
\lab{def:AandAonorms:tau1tau2:phisandpsis}
\begin{align}
\A[\pmb\psi_s](\Iti):=&\sum_{i,j=1}^3\sum_{p=0}^2\bigg(\int_{\MMh(\Iti)}\frac{|\psiss{ij}{p}|^2}{r^3}  +\sup_{\tau\in\Iti}\int_{\Sih_{\tau}}\frac{|\psiss{ij}{p}|^2}{r^2}+\int_{\II_+(\Iti)}\frac{|\psiss{ij}{p}|^2}{r^2}\bigg),\\
\A_*[\pmb\psi_s]:=&\sum_{i,j=1}^3\sum_{p=0}^2\sup_{\frac{c_*}{2}\leq c\leq c_*}\int_{\Si_{*,c}(\tau_1+1, \tau_2-2)}r^{-2}|\psiss{ij}{p}|^2,\\
\A[\pmb\phi_s](\tau_1, \tau_2):=&\sum_{p=0}^2\bigg(\int_{\MM(\tau_1, \tau_2)}\frac{|\pmb\phi_s^{(p)}|^2}{r^3}  +\sup_{\tau\in[\tau_1, \tau_2]}\int_{\Sih_{\tau}}\frac{|\pmb\phi_s^{(p)}|^2}{r^2}+\int_{\Si_*(\tau_1, \tau_2)}\frac{|\pmb\phi_s^{(p)}|^2}{r^2}\bigg),\\
\A_*[\pmb\phi_s]:=& \sum_{p=0}^2\sup_{\frac{c_*}{2}\leq c\leq c_*}\int_{\Si_{*,c}(\tau_1+1, \tau_2-2)}r^{-2}|\pmb\phi_s^{(p)}|^2.
\end{align}
\esub

Next, we define, for any $\tau'<\tau''$, $\reg\in\mathbb{N}$ and $\de\in [0,1]$, the following EMF norm 
\bea
\lab{def:widehatEMFdenorm}
\bsplit
\widehat{\EMFh}_{\de}^{(\reg)}[ \pmb\psi](\tau',\tau''):={}& \EFh^{(\reg)}[ \pmb\psi](\tau',\tau'')+\widehat{\Mh}_{\de}^{(\reg)}[ \pmb\psi](\tau',\tau''),\\
\widehat{\Mh}_{\de}^{(\reg)}[ \pmb\psi](\tau',\tau''):={}& \Mh_{\de}^{(\reg)}[ \pmb\psi](\tau',\tau'')+\int_{\Mtraph(\tau',\tau'')} |\pr^{\leq \reg+1}\pmb\psi|^2\\
\simeq{}&\int_{\MMh(\tau',\tau'')} \bigg(\frac{|\nab_{\pr_r} \dk^{\leq \reg}\pmb\psi|^2}{r^{1+\de}}
+\frac{|\nab_{\pr_\tt} \dk^{\leq \reg}\pmb\psi|^2}{r^{1+\de}}
+\frac{|(r\nab)^{\leq 1}\dk^{\leq \reg}\pmb\psi|^2}{r^{3}}\bigg),
\end{split}
\eea
in which the spacetime integrand is non-degenerate w.r.t. all derivatives in the trapping region. Also, denote for convenience  $\widehat{\EMFh}_{1}^{(\reg)}[\pmb\psi](\tau',\tau'')$ by $\widehat{\EMFh}^{(\reg)}[\pmb\psi](\tau',\tau'')$. Finally, for $1\leq\tau'<\tau''\leq\tau_*$, we denote by $\widehat{\EMF}_{\de}^{(\reg)}[ \pmb\phi](\tau',\tau'')$, $\widehat{\M}_{\de}^{(\reg)}[ \pmb\phi](\tau',\tau'')$ and $\widehat{\EMF}^{(\reg)}[\pmb\phi](\tau',\tau'')$ the analog of \eqref{def:widehatEMFdenorm} where the integration are on $\MM$ instead of $\MMh$.

We define the following early time energy norm for $s=\pm 2$ 
\begin{align}
\lab{eq:definitionofinitialenergyofphiandpsi:IEterm}
\IE{\pmb\phi_s, \pmb\psi_s}:={}&\sum_{p=0,1,2}\left(\E[\phis{p}](\tau_1)+\sup_{\tau\in[\tmic, \tau_1+2]}\Eh[\psis{p}](\tau)\right).
\end{align}
Next, we define the following EMF norm for $s=\pm 2$: 
\begin{align}
\lab{def:EMFtotalps:pm2}
\EMFtotalp{s}:={}&\sum_{p=0,1,2}\widetilde{\EMF}[\psis{p}]
+\sum_{p=0,1}\widehat{\Mh}_{\de}[\psis{p}](\tau_1+1, \tau_2-3)
\nn\\
&+\sum_{p=0,1}\widehat{\EMF}_{\de}[\phis{p}](\tau_1,\tau_2)+ {\EMF}_{\de}[\phis{2}](\tau_1,\tau_2).
\end{align}
Also, we define, for $\pmb\psi\in\sk_k(\mathbb{C})$ and $\pmb H\in\sk_k(\mathbb{C})$, $k=1,2$, 
\bea
\lab{eq:defofNNhat'}
\widehat{\mathcal{N}}'[\pmb\psi, \H](\tau_1, \tau_2)
&:=& \sup_{\tau_1< \tau'<\tau''< \tau_2}\bigg|\int_{\Mntrap(\tau', \tau'')}\big(1+O(r^{-\de})\big)\nab_{\pr_\tau}\pmb\psi\c\ov{\H}\bigg| \nn\\
&&+\int_{\Mntrap(\tau_1, \tau_2)}r^{-1}|\dk^{\leq 1}\pmb\psi||\H|+\int_{\MM(\tau_1, \tau_2)}|\H|^2,
\eea
where the coefficient $1+O(r^{-\de})$ appearing on the RHS of \eqref{eq:defofNNhat'} is in practice either equal to $1$ or to the smooth function $f_\de=f_\de(r)$ introduced in \eqref{def:Xandw:improvedMorawetz:extendedRW:Kerrpert}. Finally, we define, for $s=\pm 2$,
\bea\lab{expression:NNttotalps:pm2}
\NNttotalp{s}
&:=&\sum_{p=0,1,2}\widetilde{\mathcal{N}}[\psis{p},\widehat{\F}_{total,s}^{(p)}]+\sum_{p=0,1,2}\widehat{\mathcal{N}}'[\phis{p}, \N_{W, s}^{(p)}](\tau_1, \tau_2)\nn\\
&&+\sum_{p=0,1}\widehat{\mathcal{N}}'[ \phis{p}, r^{-2}\N_{T, s}^{(p)}](\tau_1, \tau_2)+\sum_{p=0,1}\int_{\MM_{r\leq 12m}(\tau_1,\tau_2)}|\pr^{\leq 1}\N_{T, s}^{(p)}|^2\nn\\
&&+\sum_{p=0,1}\int_{\MM(\tau_1,\tau_2)} \Big(r^{-1+\de}|\N^{(p)}_{W,s}|+r^{-2+\de}|\pr\N^{(p)}_{T,s}|\Big)|\phis{p}|,
\eea
where $\widetilde{\mathcal{N}}[\c, \c]$ and $\widehat{\NN}'[\c, \c]$ are defined in \eqref{def:NNtintermsofNNtMora:NNtEner:NNtaux:wavesystem:EMF} and \eqref{eq:defofNNhat'}. Also,  we define $\EMFtotalhps{s}{\pr^{\leq \reg}}$ and $\NNttotalph{s}{\pr^{\leq \reg}}$ accordingly by making the replacements 
\beaa
(\phis{p}, \psis{p}, \widehat{\F}_{total,s}^{(p)}, \N_{W, s}^{(p)}, \N_{T, s}^{(p)}) \to  (\pr^{\leq \reg}\phis{p}, \pr^{\leq \reg}\psis{p}, \pr^{\leq \reg}\widehat{\F}_{total,s}^{(p)}, \pr^{\leq \reg}\N_{W, s}^{(p)}, \pr^{\leq \reg}\N_{T, s}^{(p)}) .
\eeaa

We are now ready to state our main EMF estimates for $\phis{p}$ and $\psis{p}$ at zeroth-order which are the analog of the ones in Theorem 9.1 of \cite{MaSz26}.

\begin{theorem}[EMF estimates for  $\phis{p}$ and $\psis{p}$ at zeroth-order]
\lab{thm:EMF:systemofTeuscalarized:order0:final}
Let $(\MM, \g)$ satisfy the assumptions of Section \ref{sec:assumptionsforsec:energyMorawetzesitmatesforTeukoslkyonMM:upto15derivatives}. Let $\{\pmb\phi_s^{(p)}\}_{s=\pm 2, p=0,1,2}$ be a solution to the tensorial Teukolsky wave/transport systems \eqref{eq:TensorialTeuSysandlinearterms:rescaleRHScontaine2:general:Kerrperturbation}  \eqref{def:TensorialTeuScalars:wavesystem:Kerrperturbation}  in perturbations of Kerr, 
and let $\{\psis{p}\}_{s=\pm 2, p=0,1,2}$ be a solution to \eqref{eq:waveeqwidetildepsi1:gtilde:Teu} satisfying \eqref{eq:causlityrelationsforwidetildepsi1}. Then, we have, for $s=\pm 2$ and any $\de\in(0,\frac{1}{3}]$, 
\bea\lab{eq:EMFtotalp:sumup:rweightscontrolled:pm2}
&&\EMFtotalp{s} +\sum_{p=0,1,2}\sup_{\frac{c_*}{2}\leq c\leq c_*}\F_{\Si_{*,c}}[\pmb\psi_p](\tau_1+1, \tau_2-2)\nn\\ 
&\les&  \IE{\pmb\phi_s, \pmb\psi_s}  + \NNttotalp{s}+\sum_{p=0,1,2}\widetilde{\mathcal{N}}_*[\psis{p},\widehat{\F}_{total,s}^{(p)}]+\A[\pmb\psi_{s}](\Iti)+\A[\pmb\phi_{s}](\tau_1,\tau_2)+\A_*[\pmb\psi_s]\nn\\
&&+\A_*[\pmb\phi_s]+\sum_{p=0}^2\Errdefect[\psis{p}]+\sum_{p=0}^2\Errdefects[\psis{p}],
\eea
where $\EMFtotalp{s}$, $\IE{\pmb\phi_s, \pmb\psi_s}$, $\NNttotalp{s}$, $\NNt_*[\c,\c]$, $\Errdefect[\c]$, $\Errdefects[\c]$, $\A[\c](\c,\c)$ and $\A_*[\c]$  are given respectively as in \eqref{def:EMFtotalps:pm2}, \eqref{eq:definitionofinitialenergyofphiandpsi:IEterm},  \eqref{expression:NNttotalps:pm2}, \eqref{eq:definitionofNNt*pmbpsipmbF:additionalterm*}, \eqref{def:Errdefectofpsi} and \eqref{def:AandAonorms:tau1tau2:phisandpsis}.
\end{theorem}

\begin{remark}
Comparing \eqref{eq:EMFtotalp:sumup:rweightscontrolled:pm2} with the corresponding estimate (9.10) in Theorem 9.1 of \cite{MaSz26}, we emphasize the following differences:
\begin{itemize}
\item the LHS of \eqref{eq:EMFtotalp:sumup:rweightscontrolled:pm2} controls in addition the quantities $\sup_{\frac{c_*}{2}\leq c\leq c_*}\F_{\Si_{*,c}}[\pmb\psi_p](\tau_1+1, \tau_2-2)$, for $p=0,1,2$,

\item the new terms on the RHS of \eqref{eq:EMFtotalp:sumup:rweightscontrolled:pm2} compared to (9.10) in \cite{MaSz26} are $\widetilde{\mathcal{N}}_*[\psis{p},\widehat{\F}_{total,s}^{(p)}]$ and $\Errdefects[\psis{p}]$ for $p=0,1,2$, as well as $\A_*[\pmb\psi_s]$ and $\A_*[\pmb\phi_s]$. 
\end{itemize}
\end{remark}

We initiate the proof of Theorem \ref{thm:EMF:systemofTeuscalarized:order0:final} by deriving local energy estimates for $\phis{p}$ in the regions $\MM(\tau_1,\tau_1+1)$, $\MM_{r\geq\varsigma_*(\tau-\frac{5c_*}{8})}(\tau_1+1, \tau_2-3)$ and $\MM(\tau_2-3,\tau_2)$, where $\psis{p}$ differs from $\phis{p}$. First, the local in time energy estimates (9.12) on $\MM(\tau_2-3,\tau_2)$ and $\MM(\tau_1,\tau_1+1)$ on (9.13) in \cite{MaSz26} apply immediately to our setting, i.e.,
\bsub
\label{eq:localenergyestimate:future:Teu:phis01}
\begin{align}
\label{eq:localenergyestimate:future:Teu:phis0}
\EMF_{\de}[\phis{0}](\tau_2-3, \tau_2)  \les& \E[\phis{0}](\tau_2-3) + \NNtlede[\phis{0}, \N_{W,s}^{(0)}](\tau_2-3, \tau_2)+\A[\phis{1}](\tau_2-3, \tau_2),\\
\label{eq:localenergyestimate:future:Teu:phis1}
\EMF_{\de}[\phis{1}](\tau_2-3, \tau_2)  \les& \E[\phis{1}](\tau_2-3) + \NNtlede[\phis{1}, \N_{W,s}^{(1)}](\tau_2-3, \tau_2)\nn\\
&+\A[\phis{2}](\tau_2-3, \tau_2)
+\sup_{\tau\in[\tau_2-3,\tau_2]} \E[\phis{0}](\tau),\\
\label{eq:localenergyestimate:future:Teu:phis2}
\EMF_{\de}[\phis{2}](\tau_2-3, \tau_2)  \les& \E[\phis{2}](\tau_2-3) + \NNtlede[\phis{2}, \N_{W,s}^{(2)}](\tau_2-3, \tau_2)\nn\\
&+\sum_{p=0,1}\sup_{\tau\in[\tau_2-3,\tau_2]} \E[\phis{p}](\tau)
\end{align}
\esub
and
\bsub
\label{eq:localenergyestimate:past:Teu:phis01}
\begin{align}
\label{eq:localenergyestimate:past:Teu:phis0}
\EMF_{\de}[\phis{0}](\tau_1, \tau_1+1) \les& \E[\phis{0}](\tau_1) + \NNtlede[\phis{0}, \N_{W,s}^{(0)}](\tau_1, \tau_1+1)+\A[\phis{1}](\tau_1, \tau_1+1),\\
\label{eq:localenergyestimate:past:Teu:phis1}
\EMF_{\de}[\phis{1}](\tau_1, \tau_1+1)  \les& \E[\phis{1}](\tau_1)+ \NNtlede[\phis{1}, \N_{W,s}^{(1)}](\tau_1, \tau_1+1)\nn\\
&
+\A[\phis{2}](\tau_1, \tau_1+1)+\sup_{\tau\in[\tau_1,\tau_1+1]} \E[\phis{0}](\tau),\\
\label{eq:localenergyestimate:past:Teu:phis2}
\EMF_{\de}[\phis{2}](\tau_1, \tau_1+1)  \les& \E[\phis{2}](\tau_1)
+ \NNtlede[\phis{2}, \N_{W,s}^{(2)}](\tau_1, \tau_1+1)\nn\\
&
+\sum_{p=0,1}\sup_{\tau\in[\tau_1,\tau_1+1]} \E[\phis{p}](\tau),
\end{align}
\esub
where $\NNtlede$ is given by \eqref{def:NNtleinlocalenergyestimate}. Also, the Teukolsky equations \eqref{eq:ScalarizedTeuSys:general:Kerrperturbation} satisfied by $\phiss{ij}{p}$ have the following form
\beaa
\square_{\g}\phiss{ij}{0}&=&\sum_{k,l}\big(O(r^{-1})\pr\phiss{kl}{0}+O(r^{-2})\phiss{kl}{0}\big) + O(r^{-3}) \phiss{ij}{1} + N_{W,s,ij}^{(0)},\\
\square_{\g}\phiss{ij}{1}&=&\sum_{k,l}\big(O(r^{-1})\pr\phiss{kl}{1}+O(r^{-2})\phiss{kl}{1}\big) + O(r^{-3}) \phiss{ij}{2} +\sum_{k,l}O(r^{-2})\dk^{\leq 1}\phiss{kl}{0}+ N_{W,s,ij}^{(1)},\\
\square_{\g}\phiss{ij}{2}&=&\sum_{k,l}\sum_{p=0,1,2}\big(O(r^{-1})\pr\phiss{kl}{p}+O(r^{-2})\phiss{kl}{p}\big) + N_{W,s,ij}^{(2)},
\eeaa
so that they all solve a coupled system of scalar wave equations of the form \eqref{eq:eqsforlocalenergyestimatelemma:general:versionmoregeneralforRR*}. Hence, we may apply \eqref{eq:localenergyestimate:future:unweigthedderivatives:RR*} with $\reg=0$ on $\MM_{r\geq\varsigma_*(\tau-\frac{5c_*}{8})}(\tau_1+1, \tau_2-2)$ which yields 
\bsub
\label{eq:localenergyestimate:Teu:phis01:RR*regioncase}
\bea
\nn&&\EM_{0,r\geq\varsigma_*(\tau-\frac{5c_*}{8})}[\phis{0}](\tau_1+1, \tau_2-2)\\
& \les &\E[\pr^{\leq\reg}\psi](\tau_1+1) +\F_{\Si_{*,\frac{5c_*}{8}}}[\phis{0}](\tau_1+1, \tau_2-2)+ \int_{\RR_*(\tau_1+1, \tau_2-2)}r|\N_{W,s}^{(0)}|^2,\\
\nn&&\EM_{0,r\geq\varsigma_*(\tau-\frac{5c_*}{8})}[\phis{1}](\tau_1+1, \tau_2-2)\\
& \les &\E[\pr^{\leq\reg}\psi](\tau_1+1) +\F_{\Si_{*,\frac{5c_*}{8}}}[\phis{1}](\tau_1+1, \tau_2-2)+ \int_{\RR_*(\tau_1+1, \tau_2-2)}r|\N_{W,s}^{(1)}|^2,\\
\nn&&\EM_{0,r\geq\varsigma_*(\tau-\frac{5c_*}{8})}[\phis{2}](\tau_1+1, \tau_2-2)\\
& \les &\E[\pr^{\leq\reg}\psi](\tau_1+1) +\F_{\Si_{*,\frac{5c_*}{8}}}[\phis{2}](\tau_1+1, \tau_2-2)+ \int_{\RR_*(\tau_1+1, \tau_2-2)}r|\N_{W,s}^{(2)}|^2.
\eea
\esub

Next, an application of Theorem \ref{th:mainenergymorawetzmicrolocal} yields the following analog of Lemma 9.2 in \cite{MaSz26}.

\begin{lemma}
For $s=\pm 2$ and $p=0,1,2$, we have the following energy-Morawetz estimate 
\begin{align}\lab{th:eq:mainenergymorawetzmicrolocal:tensorialwave:scalarized}
&\widetilde{\EMF}[\psis{p}] +\sup_{\frac{c_*}{2}\leq c\leq c_*}\F_{\Si_{*,c}}[\psis{p}](\tau_1+1, \tau_2-2)\nn\\
\les& \sup_{\tau\in[\tmic, \tau_1+2]}\E[\psis{p}](\tau)+\A[\psis{p}](\Iti)+\A_*[\psis{p}]+\Errdefect[\psis{p}]+\Errdefects[\psis{p}]\nn\\
&+\widetilde{\mathcal{N}}[\psis{p},\chi_{\tau_1, \tau_2}^{(1)}\chi_{r_*}^{(1)}\L^{(p)}_{s}]+\widetilde{\mathcal{N}}[\psis{p},\widehat{\F}_{total,s}^{(p)}]+\NNt_*[\psis{p}, \chi_{\tau_1, \tau_2}^{(1)}\chi_{r_*}^{(1)}\L^{(p)}_{s}]+\NNt_*[\psis{p}, \widehat{\F}_{total,s}^{(p)}].
\end{align}
\end{lemma}

\begin{proof}
Noticing that, for each $p=0,1,2$, the system of wave equations \eqref{eq:waveeqwidetildepsi1:gtilde:Teu} for $\psiss{ij}{p}$ corresponds to \eqref{eq:ScalarizedWaveeq:general:Kerrpert} with $D_0=4-2\de_{p0}$  and $F_{ij}=\chi_{\tau_1, \tau_2}^{(1)}\chi_{r_*}^{(1)}L^{(p)}_{s,ij}+\widehat{F}_{total,s,ij}^{(p)}$, we may apply Theorem \ref{th:mainenergymorawetzmicrolocal}
to the system of equations \eqref{eq:waveeqwidetildepsi1:gtilde:Teu} of $\psiss{ij}{p}$ to infer, for $s=\pm 2$ and $p=0,1,2$,
\beaa
&&\widetilde{\EMF}[\psis{p}] +\sup_{\frac{c_*}{2}\leq c\leq c_*}\F_{\Si_{*,c}}[\psis{p}](\tau_1+1, \tau_2-2)\\
&\les&\sup_{\tau\in[\tmic, \tau_1+2]}\E[\psis{p}](\tau)+\A[\psis{p}](\Iti)+\A_*[\psis{p}]+\Errdefect[\psis{p}]+\Errdefects[\psis{p}]\nn\\
&&+ \widetilde{\mathcal{N}}[\psis{p}, \chi_{\tau_1, \tau_2}^{(1)}\chi_{r_*}^{(1)}\L^{(p)}_{s}+\widehat{\F}_{total,s}^{(p)}]+\NNt_*[\psis{p}, \chi_{\tau_1, \tau_2}^{(1)}\chi_{r_*}^{(1)}\L^{(p)}_{s}+\widehat{\F}_{total,s}^{(p)}]\nn\\
&\les& \sup_{\tau\in[\tmic, \tau_1+2]}\E[\psis{p}](\tau)+\A[\psis{p}](\Iti)+\A_*[\psis{p}]+\Errdefect[\psis{p}]+\Errdefects[\psis{p}]\nn\\
&&+\widetilde{\mathcal{N}}[\psis{p},\chi_{\tau_1, \tau_2}^{(1)}\chi_{r_*}^{(1)}\L^{(p)}_{s}]+\widetilde{\mathcal{N}}[\psis{p},\widehat{\F}_{total,s}^{(p)}]+\NNt_*[\psis{p}, \chi_{\tau_1, \tau_2}^{(1)}\chi_{r_*}^{(1)}\L^{(p)}_{s}]+\NNt_*[\psis{p}, \widehat{\F}_{total,s}^{(p)}]
\eeaa
as desired. This concludes the proof of the lemma.
\end{proof}

We are now in position to prove Theorem \ref{thm:EMF:systemofTeuscalarized:order0:final}.

\begin{proof}[Proof of Theorem \ref{thm:EMF:systemofTeuscalarized:order0:final}]
We proceed as in Sections 9.3 to 9.6 of \cite{MaSz26}:
\begin{enumerate}
\item First, we obtain the analog of the preliminary EMF estimates for $\psis{p}$ of Proposition 9.3 in \cite{MaSz26} on the control of 
$$\widetilde{\EMF}[\psis{p}]+\widehat{\M}_{\de}[\psis{p}](\tau_1+1, \tau_2-3)+\widehat{\EMF}_{\de}[\phis{p}](\tau_1,\tau_2), \quad p=0,1,$$
and 
$$\widetilde{\EMF}[\psis{2}]+\M_{\de}[\psis{2}](\tau_1+1, \tau_2-3)+{\EMF}_{\de}[\phis{2}](\tau_1,\tau_2).$$
The extension of the proof of Proposition 9.3 in \cite{MaSz26} to our setting is immediate and modifies the main estimates (9.17) in \cite{MaSz26} as follows
\begin{itemize}
\item on the LHS, we control in addition $\sup_{\frac{c_*}{2}\leq c\leq c_*}\F_{\Si_{*,c}}[\psis{p}](\tau_1+1, \tau_2-2)$, $p=0,1,2$, in view of \eqref{th:eq:mainenergymorawetzmicrolocal:tensorialwave:scalarized},

\item on the RHS, we have the additional terms $\A_*[\psis{p}]$, $\Errdefects[\psis{p}]$, $\NNt_*[\psis{p}, \widehat{\F}_{total,s}^{(p)}]$ and $\NNt_*[\psis{p}, \chi_{\tau_1, \tau_2}^{(1)}\chi_{r_*}^{(1)}\L^{(p)}_{s}]$, $p=0,1,2$, in view of \eqref{th:eq:mainenergymorawetzmicrolocal:tensorialwave:scalarized}, as well as the additional terms $\int_{\RR_*(\tau_1+1, \tau_2-2)}r|\N_{W,s}^{(p)}|^2$, $p=0,1,2$ in view of \eqref{eq:localenergyestimate:Teu:phis01:RR*regioncase}. In particular, the terms involving the linear coupling terms in Teukolsky, i.e., $\L_{s}^{(p)}$, $p=0,1,2$, on the RHS of the analog of (9.17) in \cite{MaSz26} are
\bea\lab{eq:termsinvolvinglinearcouplingterminTeukLsptobeabsorbed}
\widetilde{\mathcal{N}}[\psis{0},\chi_{\tau_1, \tau_2}^{(1)}\chi_{r_*}^{(1)}\L^{(0)}_{s}], \quad\NNtdemora[\psis{0},\L_{s}^{(0)}](\tau_1+1,\tau_2-3), \quad\NNt_*[\psis{p},\chi_{\tau_1, \tau_2}^{(1)}\chi_{r_*}^{(1)}\L^{(p)}_{s}],
\eea
with the quantity $\NNtdemora[\pmb\psi,\F](\tau_1+1,\tau_2-3)$ being defined as follows
\begin{align}
\lab{def:NNtdemora}
\NNtdemora[\pmb\psi,\F](\tau_1+1,\tau_2-3):=&\Bigg|\int_{\MM_{r\geq 12m}(\tau_1+1,\tau_2-3)}\sum_{i,j}\Re\Bigg(F_{ij}\ov{\bigg(X_{\de}(\psi_{ij})  - \chi_2(X_\de)^{\a}M_{i\a}^k\psi_{kj}}\nn\\
&\qquad\qquad\qquad\qquad\qquad\qquad\ov{- \chi_2(X_\de)^{\a}M_{j\a}^k\psi_{ik} + \frac{1}{2}w_{\de}\psi_{ij}   \bigg)}\Bigg) \Bigg|,
\end{align}
where the vectorfield $X_\de$ and the scalar function $w_{\de}$ are given as in \eqref{def:Xandw:improvedMorawetz:extendedRW:Kerrpert} and where the cut-off function $\chi_2$,  introduced in \eqref{def:cutoffsintime1234}, is supported in $\tau\in(\tau_1+1, \tau_2-2)$ and satisfies $\chi_2=1$ on $\tau\in(\tau_1+2, \tau_2-3)$. 
\end{itemize}
More precisely:
\begin{enumerate}
\item The first part of the proof of Proposition 9.3 in \cite{MaSz26} concerns the removal of degeneracy in $\Mtrap$ for $p=0,1$, which immediately extends to our setting. 

\item Then, the proof uses the fact that $\phis{p}=\psis{p}$ on $\MM_{r\leq\varsigma_*(\tau-\frac{5c_*}{8})}(\tau_1+1, \tau_2-3)$ in view of \eqref{eq:causlityrelationsforwidetildepsi1} to control $\phis{p}$ there thanks to the control of $\psis{p}$ on $\MMh(\Iti)$ from \eqref{th:eq:mainenergymorawetzmicrolocal:tensorialwave:scalarized}, and then the local energy estimates \eqref{eq:localenergyestimate:future:Teu:phis01} \eqref{eq:localenergyestimate:past:Teu:phis01}
\eqref{eq:localenergyestimate:Teu:phis01:RR*regioncase} to deduce the control of $\phis{p}$ on $\MM(\tau_1, \tau_2)$. 

\item Finally, the last part of Proposition 9.3 in \cite{MaSz26} consist in upgrading the control of $\Mh[\psis{p}](\tau_1+1, \tau_2-3)$, $p=0,1,2$, provided by \eqref{th:eq:mainenergymorawetzmicrolocal:tensorialwave:scalarized} to a control of $\Mh_\de[\psis{p}](\tau_1+1, \tau_2-3)$, $p=0,1,2$, based on a direct application of the improved Morawetz estimate \eqref{eq:Morawetznearinf:extendedRWsystem:Kerrpert:delta} which is immediate\footnote{Note that the additional term $\ep\sup_{\frac{c_*}{2}\leq c\leq c_*}\F_{\Si_{*,c}}[\pmb\psi](\tau_1+1, \tau_2-3)$ on the RHS of \eqref{eq:Morawetznearinf:extendedRWsystem:Kerrpert:delta} compared to (8.149) in \cite{MaSz26} is absorbed for $\ep>0$ small enough since we control $\sup_{\frac{c_*}{2}\leq c\leq c_*}\F_{\Si_{*,c}}[\psis{p}](\tau_1+1, \tau_2-2)$ on the LHS in addition to the LHS of (9.17) in \cite{MaSz26}.}.
\end{enumerate}

\item Next, we extend Proposition  9.4 in \cite{MaSz26} to our setting, which corresponds to the control of the linear coupling terms in Teukolsky appearing on the RHS of the analog of (9.17) in \cite{MaSz26}, i.e., the terms listed in \eqref{eq:termsinvolvinglinearcouplingterminTeukLsptobeabsorbed}. The first two terms in \eqref{eq:termsinvolvinglinearcouplingterminTeukLsptobeabsorbed} are estimated exactly as in the proof of Proposition  9.4 in \cite{MaSz26}, while the last term in \eqref{eq:termsinvolvinglinearcouplingterminTeukLsptobeabsorbed} is estimated in view of \eqref{eq:definitionofNNt*pmbpsipmbF:additionalterm*}
\beaa
\NNt_*[\psis{p}, \chi_{\tau_1, \tau_2}^{(1)}\chi_{r_*}^{(1)}\L^{(p)}_{s}] &\les& \bigg(\int_{\MM_{r\geq 11m}(\tau_1+1,\tau_2-3)}r^{1+\de}\big|\L^{(p)}_{s}\big|^2 \bigg)^{\frac{1}{2}}\Big({\M}_{\de}[\psis{p}](\tau_1+1,\tau_2-3)\Big)^{\frac{1}{2}}\\
&&+\bigg(\int_{\MM_{r\geq 11m}(\tau_2-3,\tau_2-2)}r^2\big|\L^{(p)}_{s}\big|^2 \bigg)^{\frac{1}{2}}\bigg(\sup_{\tau\in (\tau_2-3,\tau_2-2)}\E[\psis{p}](\tau)\bigg)^{\frac{1}{2}},
\eeaa
where the terms involving $\L^{(p)}_{s}$ on the RHS are estimated in (9.33) (9.34) of \cite{MaSz26} so that the third term in \eqref{eq:termsinvolvinglinearcouplingterminTeukLsptobeabsorbed} is estimated exactly as the estimate in Lemma 9.6 of \cite{MaSz26} for the second term in \eqref{eq:termsinvolvinglinearcouplingterminTeukLsptobeabsorbed}. 

\item Next, we consider the proof of Proposition 9.8 in \cite{MaSz26}. It concerns the control of $\pmb\phi_s^{(p)}$, $p=0,1$, on $\MM(\tau_1, \tau_2)$ and its proof in \cite{MaSz26} extends immediately to our setting. Note also that while Proposition 9.9 in \cite{MaSz26}, a higher order derivative version of Proposition 9.8 in \cite{MaSz26}, is not relevant for the proof of  Theorem \ref{thm:EMF:systemofTeuscalarized:order0:final}, it nevertheless also immediately extends to our setting. 

\item Then, Proposition 9.10 in \cite{MaSz26} is a simple consequence of Propositions 9.3, 9.4 and 9.8 in \cite{MaSz26}, so that it immediately extends to our setting. Finally, the end of the proof of Theorem 9.1 is  a simple consequence of Proposition 9.10 in \cite{MaSz26} and hence again immediately extends to our setting.
\end{enumerate}
This concludes the proof of Theorem \ref{thm:EMF:systemofTeuscalarized:order0:final}.
\end{proof}

%%%%%%%%%%%%%%%%%%%%%%%%%%%%%%%%%%%%%

\subsubsection{Proof of Theorem \ref{th:main:intermediary}}
\lab{sec:proofofth:main:intermediary}

%%%%%%%%%%%%%%%%%%%%%%%%%%%%%%%%%%%%%

We start by defining, for $\pmb\psi\in\sk_k(\mathbb{C})$ and $\pmb H\in\sk_k(\mathbb{C})$, $k=1,2$, 
\bea\lab{eq:defintionwidehatmathcalNfpsinormRHS}
&&\widehat{\NN}[\pmb\psi, \H](\tau_1, \tau_2) \nn\\
\nn&:=&\sup_{\tau_1< \tau'<\tau''< \tau_2}\bigg|\int_{\Mntraph(\tau', \tau'')}\nab_{\pr_\tau}\pmb\psi\c\ov{\H}\bigg|+\int_{\Mntraph(\tau_1, \tau_2)}r^{-1}|\dk^{\leq 1}\pmb\psi||\H|+\int_{\MMh(\tau_1, \tau_2)}|\H|^2\nn\\
&+&\!\!\!\!\min\left[\left(\int_{\Mtraph(\tau_1, \tau_2)}|\H|^2\right)^{\frac{1}{2}} \left(\int_{\Mtraph(\tau_1, \tau_2)}|\dk{^{\leq 1}}\pmb\psi|^2\right)^{\frac{1}{2}}, 
\int_{\Mtraph(\tau_1, \tau_2)}\tau^{1+\de}|\H|^2\right],
\eea
and
\bea
\lab{eq:defofNNhat''}
\widehat{\mathcal{N}}''[\pmb\psi, \H](\tau_1, \tau_2)
&:=&\sup_{\tau_1< \tau'<\tau''< \tau_2}\bigg|\int_{\Mntraph(\tau', \tau'')}\nab_{\pr_\tau}\pmb\psi\c\ov{\H}\bigg| \nn\\
&&+\int_{\Mntraph(\tau_1, \tau_2)}r^{-1}|\dk^{\leq 1}\pmb\psi||\H|+\int_{\MMh(\tau_1, \tau_2)}|\H|^2.
\eea
The following lemma allows to derive conditional EMF estimates for all high order (unweighted) derivatives conditional on the control of high order $(\pr_\tau, \chi_0(r)\pr_{\tphi})$ derivatives. 

\begin{lemma}\lab{lemma:higherorderenergyMorawetzestimates}
Let $\g$ satisfy the assumptions of Section \ref{sec:assumptionsforsec:energyMorawetzesitmatesforTeukoslkyonMM:upto15derivatives}, let $\tau_2>\tau_1$ satisfy \eqref{eq:tau2geqtau1plus10conditionisthemaincase}, and let $\reg$ be a positive integer satisfying $1\leq \reg \leq 14$. Let $(\psi)_{ij}$ and $(F)_{ij}$ be families of complex-valued scalars satisfying the following coupled system of wave equations
\bea
\lab{eq:eqsforconditionalhighorderestimatelemma:general}
\square_{\g}\psi_{ij}=N_{ij}, \qquad N_{ij}:=D_1r^{-1}\pr_{\tt}\psi_{ij}+O(r^{-2})\dk^{\leq 1}\psi_{kl} +F_{ij},
\eea
where $D_1\geq0$ is a constant and the coefficients of the terms in $N_{ij}$ may depend on $\dk^{\leq 1}\Ga_b$ and $\dk^{\leq 1}\Ga_g$.
Assume that $(\psi)_{ij}$ satisfy the following redshift estimate:
\bea
\lab{eq:assum:redshift:highorderunweightedEMF:generalwave}
\nn\EMFh_{r\leq r_+(1+\dred)}[\pr^{\leq \reg}\pmb\psi](\tau_1, \tau_2)&\les& \Eh[\pr^{\leq \reg}\pmb\psi](\tau_1)+\dred^{-1}\Mh_{r_+(1+\dred), r_+(1+2\dred)}[\pr^{\leq \reg}\pmb\psi](\tau_1, \tau_2)\\
&&+\int_{\MMh_{r\leq r_+(1+2\dred)}(\tau_1, \tau_2)}|\pr^{\leq \reg}{\F}|^2.
\eea
Let $\chi_0=\chi_0(r)$ be a smooth cut-off function such that
\bea
\lab{def:cutoffcuntionchi0}
0\leq\chi_0\leq 1, \qquad \chi_0=1 \qquad \text{for } \,\, r\leq 11m, \quad \chi_0=0 \quad \text{for }\,\, r\geq 12m.
\eea
Then, for  any $1\leq\tau_1<\tau_2 <+\infty$, we have  the improved estimate
\bea\lab{eq:higherorderenergymorawetzforlargerconditionalonsprtau}
\nn&&\EMFh_{r\geq 11m}[\pr^{\leq \reg}\pmb\psi](\tau_1, \tau_2) \nn\\
&\les& \EMFh[\pr_\tau^{\leq \reg}\pmb\psi](\tau_1, \tau_2)+\big(\EMFh[\pr^{\leq \reg-1}\pmb\psi](\tau_1, \tau_2)\big)^{\frac{1}{2}}\left(\EMFh[\pr^{\leq \reg}\pmb\psi](\tau_1, \tau_2)\right)^{\frac{1}{2}}\nn\\
&&+\sup_{\tau\in[\tau_1,\tau_2-2]}\widehat{\mathcal{N}}_{r\geq {10m}}''[\pr^{\leq \reg}\pmb\psi, \pr^{\leq \reg}\F](\tau, \tau+2)+\int_{\MMh_{r\geq 10m}(\tau_1,\tau_2)}|\pr^{\leq  \reg}\F|^2,
\eea
and
\bea\lab{eq:higherorderenergymorawetzconditionalonsprtausprphis}
\nn\EMFh[\pr^{\leq \reg}\pmb\psi](\tau_1, \tau_2) &\les& \Eh[\pr^{\leq \reg}\pmb\psi](\tau_1)+\EMFh[(\pr_\tau, \chi_0\pr_{\tphi})^{\leq \reg}\pmb\psi](\tau_1, \tau_2)+\int_{\MMh(\tau_1, \tau_2)}|\pr^{\leq \reg}\F|^2\\
&&+\sup_{\tau\in[\tau_1,\tau_2-2]}\widehat{\mathcal{N}}_{r\geq {10m}}''[\pr^{\leq \reg}\pmb\psi, \pr^{\leq \reg}\F](\tau, \tau+2)
\eea
as well as
\bea\lab{eq:higherorderenergymorawetzconditionalonsprtausprphis:nontrapped}
\nn\widehat{\EMFh}[\pr^{\leq \reg}\pmb\psi](\tau_1, \tau_2) &\les& \Eh[\pr^{\leq \reg}\pmb\psi](\tau_1)+\widehat{\EMFh}[(\pr_\tau, \chi_0\pr_{\tphi})^{\leq \reg}\pmb\psi](\tau_1, \tau_2)+\int_{\MMh(\tau_1, \tau_2)}|\pr^{\leq \reg}\F|^2\\
&&+\sup_{\tau\in[\tau_1,\tau_2-2]}\widehat{\mathcal{N}}_{r\geq {10m}}''[\pr^{\leq \reg}\pmb\psi, \pr^{\leq \reg}\F](\tau, \tau+2),
\eea
where $\widehat{\EMFh}$ and $\widehat{\NN}$ are given as in \eqref{def:widehatEMFdenorm} and \eqref{eq:defintionwidehatmathcalNfpsinormRHS}, respectively.
Furthermore, we have, for any $0\leq\de\leq 1$, 
\bea\lab{eq:improvedhigherorderenergymorawetzforlargerconditionalonsprtau}
\nn&&\EMFh_{\de, r\geq 11m}[\pr^{\leq \reg}\pmb\psi](\tau_1, \tau_2) \nn\\
&\les& \EMFh_\de [\pr_\tau^{\leq \reg}\pmb\psi](\tau_1, \tau_2)+\big(\EMFh[\pr^{\leq \reg-1}\pmb\psi](\tau_1, \tau_2)\big)^{\frac{1}{2}}\left(\EMFh[\pr^{\leq \reg}\pmb\psi](\tau_1, \tau_2)\right)^{\frac{1}{2}}\nn\\
&&+\sup_{\tau\in[\tau_1,\tau_2-2]}\widehat{\mathcal{N}}_{r\geq {10m}}''[\pr^{\leq \reg}\pmb\psi, \pr^{\leq \reg}\F](\tau, \tau+2)+\int_{\MMh_{r\geq 10m}(\tau_1,\tau_2)}|\pr^{\leq  \reg}\F|^2
\eea
and
\bea\lab{eq:imporvedhigherorderenergymorawetzconditionalonlowerrweigth}
\nn\EMFh_{\de}[\pr^{\leq \reg}\pmb\psi](\tau_1, \tau_2) &\les& \Eh[\pr^{\leq \reg}\pmb\psi](\tau_1)+\EMFh_{\de}[(\pr_\tau, \chi_0\pr_{\tphi})^{\leq \reg}\pmb\psi](\tau_1, \tau_2)+\int_{\MMh(\tau_1, \tau_2)}|\pr^{\leq \reg}\F|^2\\
&&+\sup_{\tau\in[\tau_1,\tau_2-2]}\widehat{\mathcal{N}}_{r\geq {10m}}''[\pr^{\leq \reg}\pmb\psi, \pr^{\leq \reg}\Fh](\tau, \tau+2)
\eea
as well as
\bea\lab{eq:imporvedhigherorderenergymorawetzconditionalonlowerrweigth:nontrapped}
\widehat{\EMFh}_{\de}[\pr^{\leq \reg}\pmb\psi](\tau_1, \tau_2) &\les& \Eh[\pr^{\leq \reg}\pmb\psi](\tau_1)+\widehat{\EMFh}_{\de}[(\pr_\tau, \chi_0\pr_{\tphi})^{\leq \reg}\pmb\psi](\tau_1, \tau_2)+\int_{\MM(\tau_1, \tau_2)}|\pr^{\leq \reg}\F|^2\nn\\
&&+\sup_{\tau\in[\tau_1,\tau_2-2]}\widehat{\mathcal{N}}_{r\geq {10m}}''[\pr^{\leq \reg}\pmb\psi, \pr^{\leq \reg}\F](\tau, \tau+2).
\eea
Finally, for any $1\leq\tau_1<\tau_2\leq\tau_*$, we have  
\bea\lab{eq:recoveringsupFSi*cforcbetweenc*overc*forallderivativesfromprtaudervativesonly}
\sup_{\frac{c_*}{2}\leq c\leq c_*}\F_{\Si_{*,c}}[\pr^{\leq \reg}\pmb\psi](\tau_1, \tau_2) &\les& \sup_{\frac{c_*}{2}\leq c\leq c_*}\F_{\Si_{*,c}}[\pr_\tau^{\leq \reg}\pmb\psi](\tau_1, \tau_2)+\sup_{\tau\in[\tau_1, \tau_2]}\E_{r\geq 11m}[\pr^{\leq \reg}\pmb\psi](\tau)\nn\\
&& +\int_{\RR_*(\tau_1, \tau_2)}|\pr^{\leq\reg}\F|^2.
\eea
\end{lemma}

\begin{proof}
As the estimates \eqref{def:cutoffcuntionchi0}--\eqref{eq:imporvedhigherorderenergymorawetzconditionalonlowerrweigth:nontrapped} are proved in Lemma 10.1 \cite{MaSz26}, it remains to prove  \eqref{eq:recoveringsupFSi*cforcbetweenc*overc*forallderivativesfromprtaudervativesonly}. To this end, we consider the following non-sharp consequence of (10.15) in \cite{MaSz26}, for $\reg_1+\reg_2\leq 13$, 
\begin{equation}
\lab{eq:eqsforconditionalhighorderestimatelemma:general:prreg}
\bsplit
\square_{\g}(\pr_r, r^{-1}\pr_{x^a})^{\reg_1}\pr_{\tt}^{\reg_2}\psi_{ij}={}& O(r^{-1})(\pr, r^{-1})\pr^{\leq\reg_1+\reg_2}\psi_{ij}+{F}_{ij}^{(\reg_1,\reg_2)},\\
{F}_{ij}^{(\reg_1,\reg_2)}={}&
(\pr_r, r^{-1}\pr_{x^a})^{\reg_1}\pr_{\tt}^{\reg_2}F_{ij}
+O(\ep)\pr\pr^{\leq \reg_1+\reg_2}\psi_{kl}.
\end{split}
\end{equation}
Together with (3.24) in \cite{MaSz24}, which follows from the assumptions for $\g$ satisfy the in  Section \ref{sec:assumptionsforsec:energyMorawetzesitmatesforTeukoslkyonMM:upto15derivatives}, we infer
\bea
\lab{eq:eqsforconditionalhighorderestimatelemma:general:prreg:consequence}
\gam^{ij}\pr_i\pr_j(\pr_r, r^{-1}\pr_{x^a})^{\reg_1}\pr_{\tt}^{\reg_2}\psi_{ij}&=& O(1)\pr\pr_\tau(\pr_r, r^{-1}\pr_{x^a})^{\reg_1}\pr_{\tt}^{\reg_2}\psi_{kl}\nn\\
&&+O(r^{-1})(\pr, r^{-1})\pr^{\leq\reg_1+\reg_2}\psi_{ij}+{F}_{ij}^{(\reg_1,\reg_2)}.
\eea
Next, squaring both sides of \eqref{eq:eqsforconditionalhighorderestimatelemma:general:prreg:consequence}, taking the definition of ${F}_{ij}^{(\reg_1,\reg_2)}$ in \eqref{eq:eqsforconditionalhighorderestimatelemma:general:prreg} into account and integrating on $\MM_{\varsigma_*(\tau-n)\leq r\leq\varsigma_*(\tau-n-1)}(\tau_1, \tau_2)$ for any integer $n$ such that $\frac{c_*}{2}\leq n<n+1\leq c_*$, we obtain, using also integration by parts,  
\beaa
&&\int_{\MM_{\varsigma_*(\tau-n)\leq r\leq\varsigma_*(\tau-n-1)}(\tau_1, \tau_2)}|(\pr_r, r^{-1}\pr_{x^a})^{\reg_1+2}\pr_{\tt}^{\reg_2}\psi_{ij}|^2\\ 
&\les& \int_{\MM_{\varsigma_*(\tau-n)\leq r\leq\varsigma_*(\tau-n-1)}(\tau_1, \tau_2)}|(\pr_r, r^{-1}\pr_{x^a})^{\leq\reg_1+1}\pr_{\tt}\pr_{\tt}^{\leq\reg_2}\psi_{kl}|^2\\
&&+(\ep+r_*^{-1})\int_{\MM_{\varsigma_*(\tau-n)\leq r\leq\varsigma_*(\tau-n-1)}(\tau_1, \tau_2)}|\pr^{\leq\reg_1+\reg_2+1}\pr\psi_{kl}|^2+\int_{\RR_*(\tau_1, \tau_2)}|\pr^{\leq\reg_1+\reg_2}F_{ij}|^2\\
&&+\left(\sup_{\frac{c_*}{2}\leq c\leq c_*}\F_{\Si_{*,c}}[\pr^{\leq \reg_1+\reg_2+1}\pmb\psi](\tau_1, \tau_2)\right)^{\frac{1}{2}}\left(\sup_{\frac{c_*}{2}\leq c\leq c_*}\F_{\Si_{*,c}}[\pr^{\leq \reg_1+\reg_2}\pmb\psi](\tau_1, \tau_2)\right)^{\frac{1}{2}}.
\eeaa
By induction we infer, for $\ep>0$ small enough and $r_*$ large enough, and for $1\leq\reg\leq 14$, 
\beaa
&&\int_{\MM_{\varsigma_*(\tau-n)\leq r\leq\varsigma_*(\tau-n-1)}(\tau_1, \tau_2)}|\pr^{\leq\reg}\pr\pmb\psi|^2\\
&\les& \sup_{\frac{c_*}{2}\leq c\leq c_*}\F_{\Si_{*,c}}[\pr_\tau^{\leq \reg}\pmb\psi](\tau_1, \tau_2)+\int_{\RR_*(\tau_1, \tau_2)}|\pr^{\leq\reg}\F|^2\\
&&+\left(\sup_{\frac{c_*}{2}\leq c\leq c_*}\F_{\Si_{*,c}}[\pr^{\leq \reg}\pmb\psi](\tau_1, \tau_2)\right)^{\frac{1}{2}}\left(\sup_{\frac{c_*}{2}\leq c\leq c_*}\F_{\Si_{*,c}}[\pr^{\leq \reg-1}\pmb\psi](\tau_1, \tau_2)\right)^{\frac{1}{2}}.
\eeaa
In particular, for any integer $n$ such that $\frac{c_*}{2}\leq n<n+1\leq c_*$, there exists $c_n\in[n, n+1]$ such that, for $1\leq\reg\leq 14$, 
\beaa
\int_{\Si_{*,c_n}(\tau_1, \tau_2)}|\pr^{\leq\reg}\pr\pmb\psi|^2 &\les& \sup_{\frac{c_*}{2}\leq c\leq c_*}\F_{\Si_{*,c}}[\pr_\tau^{\leq \reg}\pmb\psi](\tau_1, \tau_2)+\int_{\RR_*(\tau_1, \tau_2)}|\pr^{\leq\reg}\F|^2\\
&&+\left(\sup_{\frac{c_*}{2}\leq c\leq c_*}\F_{\Si_{*,c}}[\pr^{\leq \reg}\pmb\psi](\tau_1, \tau_2)\right)^{\frac{1}{2}}\left(\sup_{\frac{c_*}{2}\leq c\leq c_*}\F_{\Si_{*,c}}[\pr^{\leq \reg-1}\pmb\psi](\tau_1, \tau_2)\right)^{\frac{1}{2}}.
\eeaa
Then, arguing by local existence on $\MM_{\varsigma_*(\tau-c_n+1)\leq r\leq\varsigma_*(\tau-c_n-1)}(\tau_1, \tau_2)$ similarly as in the proof of \eqref{eq:localenergyestimate:future:unweigthedderivatives:RR*}, and then taking the supremum on integers $n$ such that $\frac{c_*}{2}\leq n<n+1\leq c_*$, we deduce, for $1\leq\reg\leq 14$,  
\beaa
&&\sup_{\frac{c_*}{2}\leq c\leq c_*}\F_{\Si_{*,c}}[\pr^{\leq \reg}\pmb\psi](\tau_1, \tau_2)\\ 
&\les& \sup_{\frac{c_*}{2}\leq c\leq c_*}\F_{\Si_{*,c}}[\pr_\tau^{\leq \reg}\pmb\psi](\tau_1, \tau_2)+\sup_{\tau\in[\tau_1, \tau_2]}\E_{r\geq 11m}[\pr^{\leq \reg}\pmb\psi](\tau) +\int_{\RR_*(\tau_1, \tau_2)}|\pr^{\leq\reg}\F|^2\\
&&+\left(\sup_{\frac{c_*}{2}\leq c\leq c_*}\F_{\Si_{*,c}}[\pr^{\leq \reg}\pmb\psi](\tau_1, \tau_2)\right)^{\frac{1}{2}}\left(\sup_{\frac{c_*}{2}\leq c\leq c_*}\F_{\Si_{*,c}}[\pr^{\leq \reg-1}\pmb\psi](\tau_1, \tau_2)\right)^{\frac{1}{2}},
\eeaa
from which \eqref{eq:recoveringsupFSi*cforcbetweenc*overc*forallderivativesfromprtaudervativesonly} follows by induction. This concludes the proof of Lemma \ref{lemma:higherorderenergyMorawetzestimates}.
\end{proof}

Next, we introduce the following EMF norms for $s=\pm 2$ and any $\de\in(0,\frac{1}{3}]$
\begin{align}
\lab{def:EMFtotalpsnotilde:pm2}
\EMFtotalps{s}:={}&\sum_{p=0,1,2}{\EMFh}[\psis{p}](\Iti)
+\sum_{p=0,1}\widehat{\M}_{\de}[\psis{p}](\tau_1+1, \tau_2-3)
\nn\\
&+\sum_{p=0,1}\widehat{\EMF}_{\de}[\phis{p}](\tau_1,\tau_2)+ {\EMF}_{\de}[\phis{2}](\tau_1,\tau_2),
\end{align}
and the following early time EMF norms, for $s=\pm 2$ and any $\de\in[0,1]$,
\bea\lab{eq:definitionofinitialenergyofphiandpsi:IEtotalterm:delta}
\IEde{\pr^{\leq \reg}\pmb\phi_s, \pr^{\leq \reg}\pmb\psi_s} :=\IE{\pr^{\leq \reg}\pmb\phi_s, \pr^{\leq \reg}\pmb\psi_s} + \sum_{p=0}^2{\MFh}_{\de}[\pr^{\leq \reg}\psis{p}](\tmic,\tau_1+1),
\eea
where $\IE{\pr^{\leq \reg}\pmb\phi_s, \pr^{\leq \reg}\pmb\psi_s}$ is given by \eqref{eq:definitionofinitialenergyofphiandpsi:IEterm}. Then, we apply Lemma \ref{lemma:higherorderenergyMorawetzestimates} to the coupled Teukolsky wave system and deduce the following high order EMF estimates which recover the control of all unweighted derivatives from the control of high order $(\pr_{\tt}, \chi_0\pr_{\tphi})$ derivatives.

\begin{lemma}\lab{lemma:higherorderenergyMorawetzestimates:Teu}
Under the same assumptions as in Theorem \ref{th:main:intermediary}, we have, for  $s=\pm 2$, $1\leq \reg\leq 14$ and $\de\in(0,\frac{1}{3}]$, 
\begin{align}\lab{eq:higherorderenergymorawetzforlargerconditionalonsprtau:Teu}
&\EMF_{s, \de, \text{total}, r\geq 11m}[\pr^{\leq \reg}\pmb\phi_{s}] +\sum_{p=0}^2\sup_{\frac{c_*}{2}\leq c\leq c_*}\F_{\Si_{*,c}}[\pr^{\leq \reg}\psis{p}](\tau_1+1, \tau_2-2)\nn\\
\les&\EMFtotalhps{s}{\nab_{\pr_{\tt}}^{\leq \reg}}+\IEde{\pr^{\leq \reg}\pmb\phi_s, \pr^{\leq \reg}\pmb\psi_s} +\sum_{p=0}^2\sup_{\frac{c_*}{2}\leq c\leq c_*}\F_{\Si_{*,c}}[\nab_{\pr_{\tt}}^{\leq \reg}\psis{p}](\tau_1+1, \tau_2-2)\nn\\
&+\sum_{p=0}^{2}\widehat{\mathcal{N}}_{r\geq 10m}''[\pr^{\leq \reg}\psis{p}, \pr^{\leq \reg}\widehat{\F}_{total,s}^{(p)}](\tau_2-3, \tau_2)+\sum_{p=0}^{2}\int_{\RR_*(\tau_1+1, \tau_2-2)}|\pr^{\leq\reg}\widehat{\F}_{total,s}^{(p)}|^2\nn\\
&+\NNttotalph{s}{\pr^{\leq \reg}}+\Big(\EMFtotalhps{s}{\pr^{\leq \reg-1}}\Big)^{\frac{1}{2}}\Big(\EMFtotalhps{s}{\pr^{\leq \reg}}\Big)^{\frac{1}{2}},
\end{align}
with $\EMF_{s,\de, \text{total}}$, $\IEde{\pr^{\leq \reg}\pmb\phi_s, \pr^{\leq \reg}\pmb\psi_s}$, $\widetilde{\NN}_{s, \de, \text{total}}$ and $\widehat{\NN}''[\c,\c](\c,\c)$ given as in \eqref{def:EMFtotalpsnotilde:pm2}, \eqref{eq:definitionofinitialenergyofphiandpsi:IEtotalterm:delta},  \eqref{expression:NNttotalps:pm2} and \eqref{eq:defofNNhat''}, respectively.

Also, under the same assumptions  as in Theorem \ref{th:main:intermediary}, we have,  for  $s=\pm 2$, $1\leq \reg\leq 14$ and $\de\in(0,\frac{1}{3}]$, 
\begin{align}
\lab{eq:imporvedhigherorderenergymorawetzconditionalonlowerrweigth:Teu}
&\EMFtotalhps{s}{\pr^{\leq \reg}}+\sum_{p=0}^2\sup_{\frac{c_*}{2}\leq c\leq c_*}\F_{\Si_{*,c}}[\pr^{\leq \reg}\psis{p}](\tau_1+1, \tau_2-2)\nn\\
\les&\EMFtotalhps{s}{\nab_{(\pr_{\tt},\chi_{0}\pr_{\tilde{\phi}})}^{\leq \reg}}+\IEde{\pr^{\leq \reg}\pmb\phi_s, \pr^{\leq \reg}\pmb\psi_s} +\sum_{p=0}^2\sup_{\frac{c_*}{2}\leq c\leq c_*}\F_{\Si_{*,c}}[\nab_{\pr_{\tt}}^{\leq \reg}\psis{p}](\tau_1+1, \tau_2-2)\nn\\
&+\NNttotalph{s}{\pr^{\leq \reg}}+\sum_{p=0}^{2}\widehat{\mathcal{N}}''[\pr^{\leq \reg}\psis{p}, \pr^{\leq \reg}\widehat{\F}_{total,s}^{(p)}](\tau_2-3, \tau_2)\nn\\
&+\sum_{p=0}^{2}\int_{\RR_*(\tau_1+1, \tau_2-2)}|\pr^{\leq\reg}\widehat{\F}_{total,s}^{(p)}|^2,
\end{align}
with the cut-off function $\chi_0$ given as in \eqref{def:cutoffcuntionchi0}.
\end{lemma}

\begin{proof}
The proof of Lemma 10.2 in \cite{MaSz26}, based on Lemma \ref{lemma:higherorderenergyMorawetzestimates}, immediately applies and yields
\beaa
\EMF_{s, \de, \text{total}, r\geq 11m}[\pr^{\leq \reg}\pmb\phi_{s}] &\les&\EMFtotalhps{s}{\nab_{\pr_{\tt}}^{\leq \reg}}+\IEde{\pr^{\leq \reg}\pmb\phi_s, \pr^{\leq \reg}\pmb\psi_s} \\
&&+\NNttotalph{s}{\pr^{\leq \reg}}+\sum_{p=0}^{2}\widehat{\mathcal{N}}_{r\geq 10m}''[\pr^{\leq \reg}\psis{p}, \pr^{\leq \reg}\widehat{\F}_{total,s}^{(p)}](\tau_2-3, \tau_2)\nn\\
&&
+\Big(\EMFtotalhps{s}{\pr^{\leq \reg-1}}\Big)^{\frac{1}{2}}\Big(\EMFtotalhps{s}{\pr^{\leq \reg}}\Big)^{\frac{1}{2}},
\eeaa
and
\begin{align*}
\EMFtotalhps{s}{\pr^{\leq \reg}} \les&\EMFtotalhps{s}{\nab_{(\pr_{\tt},\chi_{0}\pr_{\tilde{\phi}})}^{\leq \reg}}+\IEde{\pr^{\leq \reg}\pmb\phi_s, \pr^{\leq \reg}\pmb\psi_s}\\
&+\NNttotalph{s}{\pr^{\leq \reg}}+\sum_{p=0}^{2}\widehat{\mathcal{N}}''[\pr^{\leq \reg}\psis{p}, \pr^{\leq \reg}\widehat{\F}_{total,s}^{(p)}](\tau_2-3, \tau_2).
\end{align*}
Then, we conclude the proof of Lemma \ref{lemma:higherorderenergyMorawetzestimates:Teu} by using \eqref{eq:recoveringsupFSi*cforcbetweenc*overc*forallderivativesfromprtaudervativesonly} to control in addition,  for $p=0,1,2$, the terms $\F_{\Si_{*,c}}[\pr^{\leq \reg}\psis{p}](\tau_1+1, \tau_2-2)$.  
\end{proof}

Next, we prove the following analog of Proposition 10.3 in \cite{MaSz26} on the control of EMF estimates for high order unweighted derivatives of $(\phis{p}, \psis{p})$, which is the analog of Theorem \ref{thm:EMF:systemofTeuscalarized:order0:final} in the high order regularity case.

\begin{proposition}[Conditional EMF estimates for unweighted derivatives of $(\phis{p}, \psis{p})$]
\lab{prop:EMFtotalp:finalestimate:pm2:prhighorder}
Under the assumptions of Theorem \ref{th:main:intermediary}, we have, for $s=\pm2$, all $\reg\leq 14$ and $\de\in (0,\frac{1}{3}]$,
\begin{align}
\lab{eq:EMFtotalp:finalestimate:pm2:prhighorder} 
&\EMFtotalh{s}{\nab_{(\pr_{\tt},\chi_{0}\pr_{\tilde{\phi}})}^{\leq \reg}}+\EMFtotalhps{s}{\pr^{\leq \reg}}+\sum_{p=0}^2\sup_{\frac{c_*}{2}\leq c\leq c_*}\F_{\Si_{*,c}}[\pr^{\leq \reg}\psis{p}](\tau_1+1, \tau_2-2) \nn\\
\les& \IEde{\pr^{\leq  \reg}\pmb\phi_s, \pr^{\leq  \reg}\pmb\psi_s} +\sum_{p=0}^2\Errdefect[\nab_{(\pr_{\tt},\chi_{0}\prtphihat)}^{\leq \reg}\psis{p}]+\sum_{p=0}^2\Errdefects[\nab_{\pr_{\tt}}^{\leq \reg}\psis{p}]+\NNttotalph{s}{\pr^{\leq \reg}}\nn\\
& +\sum_{p=0,1,2}\widetilde{\mathcal{N}}_*[\pr^{\leq \reg}\psis{p},\pr^{\leq \reg}\widehat{\F}_{total,s}^{(p)}] +\A[\pmb\psi_{s}](\Iti)+\A[\pmb\phi_{s}](\tau_1,\tau_2) +\A_*[\pmb\psi_s]+\A_*[\pmb\phi_s]\nn\\
&+\sum_{p=0}^{2}\widehat{\mathcal{N}}''[\pr^{\leq \reg}\psis{p}, \pr^{\leq \reg}\widehat{\F}_{total,s}^{(p)}](\tau_2-3, \tau_2)+\sum_{p=0}^{2}\int_{\RR_*(\tau_1+1, \tau_2-2)}|\pr^{\leq\reg}\widehat{\F}_{total,s}^{(p)}|^2,
\end{align}
with $\widetilde{\EMF}_{s,\de, \text{total}}$, $\EMF_{s,\de,\text{total}}$, $\IEde{\pr^{\leq  \reg}\pmb\phi_s, \pr^{\leq  \reg}\pmb\psi_s}$, $\Errdefect[\c]$, $\Errdefects[\c]$, $\widetilde{\NN}_{s, \de, \text{total}}$, $\A$, $\A_*$ and $\widehat{\mathcal{N}}''$ given in \eqref{def:EMFtotalps:pm2}, \eqref{def:EMFtotalpsnotilde:pm2}, \eqref{eq:definitionofinitialenergyofphiandpsi:IEtotalterm:delta}, \eqref{def:Errdefectofpsi}, \eqref{expression:NNttotalps:pm2}, \eqref{def:AandAonorms:tau1tau2:phisandpsis} and \eqref{eq:defofNNhat''}, respectively, and with the cut-off function $\chi_0$ introduced in \eqref{def:cutoffcuntionchi0}.
\end{proposition}

\begin{proof}
Following the procedure in Section 10.2.1 of \cite{MaSz26} which relies on applying Theorem \ref{thm:EMF:systemofTeuscalarized:order0:final} to the equations commuted by $\pr_{\tt}$, we obtain the following analog of (10.68) in \cite{MaSz26}
\begin{align}\lab{eq:EMFtotalp:sumup:finalestimate:pm2:prttandprtphi:sofaronlyprtt}
&\EMFtotalh{s}{\nab_{\pr_{\tt}}^{\leq 1}}+\sum_{p=0,1,2}\sup_{\frac{c_*}{2}\leq c\leq c_*}\F_{\Si_{*,c}}[\nab_{\pr_{\tt}}^{\leq 1}\pmb\psi_p](\tau_1+1, \tau_2-2)\nn\\
\les{}&\IE{\nab_{\pr_{\tt}}^{\leq 1}\pmb\phi_s, \nab_{\pr_{\tt}}^{\leq 1}\pmb\psi_s} +\NNttotalph{s}{\nab_{\pr_{\tt}}^{\leq 1}}+\sum_{p=0,1,2}\widetilde{\mathcal{N}}_*[\nab_{\pr_{\tt}}^{\leq 1}\psis{p},\nab_{\pr_{\tt}}^{\leq 1}\widehat{\F}_{total,s}^{(p)}]\nn\\
&+\ep \EMFtotalhps{s}{\pr^{\leq 1}} +\big(\EMFtotalhps{s}{\pr^{\leq 1}}\big)^{\frac{1}{2}} \big(\EMFtotalps{s}+\Ao[\nab_{\pr_{\tt}}^{\leq 1}\pmb\psi_{s}](\Iti)\big)^{\frac{1}{2}}\nn\\
&
+\sum_{p=0}^2\Errdefect[\nab_{\pr_{\tau}}^{\leq 1}\psis{p}]+\sum_{p=0}^2\Errdefects[\nab_{\pr_{\tau}}^{\leq 1}\psis{p}]+\A[\nab_{\pr_{\tt}}^{\leq 1}\pmb\psi_{s}](\Iti)+\A[\nab_{\pr_{\tt}}^{\leq 1}\pmb\phi_{s}](\tau_1,\tau_2)\nn\\
&+\A_*[\nab_{\pr_{\tt}}^{\leq 1}\pmb\psi_{s}]+\A_*[\nab_{\pr_{\tt}}^{\leq 1}\pmb\phi_{s}],
\end{align}
where we have introduced the following notation, for any $\tau'<\tau''$,
\bea
\Ao[\pmb\psi_s](\tau',\tau''):=\sum_{i,j=1}^3\sum_{p=0}^2\int_{\MM(\tau',\tau'')}\frac{|\psiss{ij}{p}|^2}{r^3},
\eea
and the corresponding quantity $\Ao[\pmb\phi_s](\tau',\tau'')$ with $\psiss{ij}{p}$ replaced by $\phiss{ij}{p}$ in the above formula. The only change w.r.t. to the procedure in Section 10.2.1 of \cite{MaSz26} is for the equation for $\psis{p}$ for which $\widetilde{\widehat{F}_{total,s,ij}^{(p),\pr_{\tt}}}$ in (10.60b) of \cite{MaSz26} has now the term 
$$-\pr_{\tt}(\chi_{\tau_1, \tau_2}\chi_{r_*})f_p\psi^{(p)}_{s,ij}+\pr_{\tt}(\chi_{\tau_1, \tau_2}^{(1)}\chi_{r_*}^{(1)})L^{(p)}_{s,ij}$$
which yields the new term
$$-\chi_{\tau_1, \tau_2}\pr_{\tt}(\chi_{r_*})f_p\psi^{(p)}_{s,ij}+\chi_{\tau_1, \tau_2}\pr_{\tt}(\chi_{r_*}^{(1)})L^{(p)}_{s,ij}.$$
In particular, we need to estimate the following two new terms in the extension of (10.63) in \cite{MaSz26} to our setting 
\beaa
\sum_{p=0}^2\widetilde{\mathcal{N}}\big[\psish{p}{\pr_{\tt}},\chi_{\tau_1,\tau_2}^{(1)}\pr_{\tt}(\chi_{r_*}^{(1)})\L_{s}^{(p)}\big], \quad \sum_{p=0}^2\widetilde{\mathcal{N}}\big[\psish{p}{\pr_{\tt}},-\chi_{\tau_1, \tau_2}\pr_{\tt}(\chi_{r_*})f_p\psis{p}\big].
\eeaa
Now, in view of \eqref{eq:supportpropertiesofcutoffchistar} and \eqref{eq:supportpropertiesofcutoffchistar:(1)}, we have
\beaa
\widetilde{\mathcal{N}}\big[\psi,\chi_{\tau_1,\tau_2}^{(1)}\pr_{\tt}(\chi_{r_*}^{(1)})F\big]+\widetilde{\mathcal{N}}\big[\psi,\chi_{\tau_1, \tau_2}\pr_{\tt}(\chi_{r_*})F\big] &\les& \int_{\RR_*(\tau_1, \tau_2)}|F|(|\pr\psi|+r^{-1}|\psi|)+\int_{\RR_*}|F|^2\\
&\les& \left(\int_{\MM_{r\geq 11m}(\tau_1, \tau_2)}r^{1+\de}|F|^2\right)^{\frac{1}{2}}\left(\M_\de[\psi](\tau_1, \tau_2)\right)^{\frac{1}{2}}
\eeaa
so that the above two new terms are also controlled by the RHS of (10.67) in \cite{MaSz26} which concludes the proof of \eqref{eq:EMFtotalp:sumup:finalestimate:pm2:prttandprtphi:sofaronlyprtt}. 

Then, commuting further with $\pr_{\tt}$, and by an inductive argument, we infer that, for any $\reg\leq 14$,
\begin{align}
\lab{eq:EMFtotalp:sumup:finalestimate:pm2:prtt}
&\EMFtotalh{s}{\nab_{\pr_{\tt}}^{\leq \reg}} +\sum_{p=0,1,2}\sup_{\frac{c_*}{2}\leq c\leq c_*}\F_{\Si_{*,c}}[\nab_{\pr_{\tt}}^{\leq\reg}\pmb\psi_p](\tau_1+1, \tau_2-2)\nn\\
\les{}&\IE{\pr^{\leq \reg}\pmb\phi_s, \pr^{\leq \reg}\pmb\psi_s} +\NNttotalph{s}{\nab_{\pr_{\tt}}^{\leq \reg}}+\sum_{p=0,1,2}\widetilde{\mathcal{N}}_*[\nab_{\pr_{\tt}}^{\leq\reg}\psis{p},\nab_{\pr_{\tt}}^{\leq\reg}\widehat{\F}_{total,s}^{(p)}]\nn\\
&+\ep \EMFtotalhps{s}{\pr^{\leq \reg}}+\big(\EMFtotalhps{s}{\pr^{\leq\reg}}\big)^{\frac{1}{2}} \big(\EMFtotalhps{s}{\pr^{\leq\reg-1}}+\Ao[\nab_{\pr_{\tt}}^{\leq\reg}\pmb\psi_{s}](\Iti)\big)^{\frac{1}{2}}\nn\\
&  +\sum_{p=0}^2\Errdefect[\nab_{\pr_{\tau}}^{\leq \reg}\psis{p}]+\sum_{p=0}^2\Errdefects[\nab_{\pr_{\tau}}^{\leq \reg}\psis{p}] +\A[\nab_{\pr_{\tt}}^{\leq\reg}\pmb\psi_{s}](\Iti)+\A[\nab_{\pr_{\tt}}^{\leq\reg}\pmb\phi_{s}](\tau_1,\tau_2)\nn\\
&+\A_*[\nab_{\pr_{\tt}}^{\leq\reg}\pmb\psi_{s}]+\A_*[\nab_{\pr_{\tt}}^{\leq\reg}\pmb\phi_{s}].
\end{align}

Then, we commute further with $(\pr_{\tt}, \chi_0\prtphihat)$ and follow the procedure in Section 10.2.2 of \cite{MaSz26} which applies immediately since the only novelty compared to \eqref{eq:EMFtotalp:sumup:finalestimate:pm2:prtt} happens in $r\leq 12m$ where our setting is identical to the one in \cite{MaSz26}. This then yields the following EMF estimates for $(\pr_\tau, \chi_{0}\prtphihat)^{\leq \reg}$-derivatives of $(\pmb\phi_s, \pmb\psi_s)$, with $\reg\leq 14$, 
 \begin{align}
\lab{eq:EMFtotalp:sumup:finalestimate:pm2:prttandprtphi:highorder}
&\EMFtotalh{s}{\nab_{(\pr_{\tt},\chi_{0}\prtphihat)}^{\leq \reg}} +\sum_{p=0,1,2}\sup_{\frac{c_*}{2}\leq c\leq c_*}\F_{\Si_{*,c}}[\nab_{\pr_{\tt}}^{\leq \reg}\pmb\psi_p](\tau_1+1, \tau_2-2)\nn\\
\les&\IE{\pr^{\leq  \reg}\pmb\phi_s, \pr^{\leq  \reg}\pmb\psi_s} +\A[\nab_{(\pr_{\tt},\chi_{0}\prtphihat)}^{\leq  \reg}\pmb\psi_{s}](\Iti)+\A[\nab_{(\pr_{\tt},\chi_{0}\prtphihat)}^{\leq  \reg}\pmb\phi_{s}](\tau_1,\tau_2)+\A_*[\nab_{\pr_{\tt}}^{\leq\reg}\pmb\psi_{s}]\nn\\
&+\A_*[\nab_{\pr_{\tt}}^{\leq\reg}\pmb\phi_{s}]+\ep\EMFtotalhps{s}{\pr^{\leq  \reg}}+\sum_{p=0}^2\Errdefect[\nab_{(\pr_{\tt},\chi_{0}\prtphihat)}^{\leq \reg}\psis{p}]+\sum_{p=0}^2\Errdefects[\nab_{\pr_{\tau}}^{\leq \reg}\psis{p}] \nn\\
&+\NNttotalph{s}{\nab_{(\pr_{\tt},\chi_{0}\prtphihat)}^{\leq \reg}} +\sum_{p=0,1,2}\widetilde{\mathcal{N}}_*[\nab_{\pr_{\tt}}^{\leq\reg}\psis{p},\nab_{\pr_{\tt}}^{\leq\reg}\widehat{\F}_{total,s}^{(p)}]+\EMF_{s,\de, \text{total}, 11m\leq r\leq 12m}[\pr^{\leq \reg}\pmb\phi_{s}]\nn\\
&+\big(\EMFtotalhps{s}{\pr^{\leq\reg}}\big)^{\frac{1}{2}} \big(\EMFtotalhps{s}{\pr^{\leq\reg-1}}+\Ao[\nab_{(\pr_{\tt}, \chi_{0}\pr_{\tphi})}^{\leq\reg}\pmb\psi_{s}](\Iti)\big)^{\frac{1}{2}}.
 \end{align}
Finally, the proof of \eqref{eq:EMFtotalp:finalestimate:pm2:prhighorder} follows, as in section 10.2.3 in \cite{MaSz26}, from a simple combination of Lemma \ref{lemma:higherorderenergyMorawetzestimates:Teu} and \eqref{eq:EMFtotalp:sumup:finalestimate:pm2:prttandprtphi:highorder}, noticing in addition the following analog of (10.94) in \cite{MaSz26} for $\A_*$
\beaa
\A_*[\nab_{\pr_{\tt}}^{\leq\reg}\pmb\psi_s] &\les& \sum_{p=0,1,2}\sup_{\frac{c_*}{2}\leq c\leq c_*}\F_{\Si_{*,c}}[\nab_{\pr_{\tt}}^{\leq \reg-1}\pmb\psi_p](\tau_1+1, \tau_2-2), \quad \forall 1\leq\reg\leq 14.
\eeaa
 This concludes the proof of Proposition \ref{prop:EMFtotalp:finalestimate:pm2:prhighorder}.
\end{proof}

We are now in position to prove Theorem \ref{th:main:intermediary}.

\begin{proof}[Proof of Theorem \ref{th:main:intermediary}]
We proceed as in Sections 10.3 to 10.5 of \cite{MaSz26}:
\begin{enumerate}
\item As in section 10.3 of \cite{MaSz26}, we control of the terms ${\bf IE}_{\de}$ and $\Errdefect$ appearing on the RHS of \eqref{eq:EMFtotalp:finalestimate:pm2:prhighorder}. The proof is identical in our setting and we obtain, for $s=\pm2$, all $\reg\leq 14$ and $\de\in (0,\frac{1}{3}]$,
\begin{align}
\lab{eq:localenergyforpsiontau1minus1tau1plus2:initialEMF}
\IEde{\pr^{\leq  \reg}\pmb\phi_{s}, \pr^{\leq  \reg}\pmb\psi_{s}}\les\sum_{p=0}^2\bigg(\E[\pr^{\leq \reg}\phis{p}](\tau_1)+\int_{\MM(\tau_1, \tau_1+2)}r^{1+\de} |\pr^{\leq \reg}\N_{W,s}^{(p)}|^2\bigg),
\end{align}
as well as 
\bea
\lab{esti:Errdefectterm:highorderprderi}
&&\sum_{p=0}^2\Errdefect[\nab_{(\pr_{\tt},\chi_{0}\prtphihat)}^{\leq \reg}\psis{p}]\nn\\
&\les&\sum_{p=0}^2\Big(\E[\pr^{\leq \reg}\phis{p}](\tau_1)+\ep^2\E[\pr^{\leq \reg}\phis{p}](\tau_2-3)\Big)+\sum_{p=0}^2\int_{\MM(\tau_1, \tau_1+2)}r^{1+\de} |\pr^{\leq \reg}\N_{W,s}^{(p)}|^2\nn\\
&&+\ep^2\sum_{p=0}^2\int_{\MM(\tau_2-3, \tau_2-2)}r^{1+\de} |\pr^{\leq \reg}\N_{W,s}^{(p)}|^2.
\eea
In addition, we also need to control the term $\Errdefects$ appearing on the RHS of \eqref{eq:EMFtotalp:finalestimate:pm2:prhighorder}. In view of \eqref{eq:propertyofscalarizationdefect:prop2}, we have, for $\reg\leq 14$,
\begin{align}\lab{esti:Errdefectterm:highorderprderi:Fscase}
& \sum_{p=0}^2\Errdefects[\nab_{\pr_{\tt}}^{\leq \reg}\psis{p}]\nn \\
\les& \sum_{p=0}^2\sup_{c_*/2\leq c\leq c_*}\F_{\Si_{*,c}}[\pr^{\leq \reg}\err_{\textrm{TDefect}}[\psi^{(p)}_{s}]](\tau_1+1, \tau_2-2)\nn\\
\nn\les& \ep^2\sum_{i,j=1}^3\sum_{p=0}^2\bigg(\E[\pr^{\leq \reg}\phi^{(p)}_{s,ij}](\tau_1)+\E[\pr^{\leq \reg}\phi^{(p)}_{s,ij}](\tau_2-3)+\F_{\Si_{*,c_*/2}}[\pr^{\leq\reg}\phi^{(p)}_{s,ij}](\tau_1+1, \tau_2-1)\\
&+\int_{\MM(\tau_1, \tau_1+1)\cup\RR_*(\tau_1+1, \tau_2-3)\cup\MM(\tau_2-3, \tau_2-2)}r^{1+\de}|\pr^{\leq \reg}N_{W,s,ij}^{(p)}|^2\bigg).
\end{align}

\item Next, as in section 10.4 of \cite{MaSz26}, we control of the terms $\NNt_{s, \de, \text{total}}$ appearing on the RHS of \eqref{eq:EMFtotalp:finalestimate:pm2:prhighorder}. In view of \eqref{eq:definitionofwidehatFtotalsij}, the difference with the estimate of $\NNttotalph{s}{\pr^{\leq \reg}}$ in (10.98) of \cite{MaSz26} is that we now need to control $\breve{\F}^{(p)}_{s}$ also in the region $\MM_{r\geq\varsigma_*(\tau-\frac{7c_*}{8})}(\tau_1+1,\tau_2-1)$. To this end, it suffices to control the quantity $\KK_{\reg}$ given for $\reg\leq 14$ by
\beaa
\KK_{\reg} := \sum_{p=0}^2\sum_{i,j=1}^3\int_{\MM_{r\geq\varsigma_*(\tau-\frac{7c_*}{8})}(\tau_1+1,\tau_2-1)}\Big(|\pr^{\leq \reg}\breve{F}^{(p)}_{s,ij}|\big(|\pr\pr^{\leq \reg}\psi^{(p)}_{s,ij}|+r^{-1}|\pr^{\leq \reg}\psi^{(p)}_{s,ij}|\big)+|\pr^{\leq \reg}\breve{F}^{(p)}_{s,ij}|^2\Big).
\eeaa
Together with \eqref{eq:structureofbreveFpsij}, we infer, for $\reg\leq 14$,
\begin{align*}
\KK_{\reg} \les& \sum_{p=0}^2\sum_{i,j=1}^3\int_{\MM_{r\geq\varsigma_*(\tau-\frac{7c_*}{8})}(\tau_1+1,\tau_2-1)}\Big(r^{-1}\big(|\pr\pr^{\leq \reg}(\widetilde{\phi}-\breve{\phi})^{(p)}_{s,kl}|+r^{-1}|\pr^{\leq \reg}(\widetilde{\phi}-\breve{\phi})^{(p)}_{s,kl}|\big)\\
&\times\big(|\pr\pr^{\leq \reg}\psi^{(p)}_{s,ij}|+r^{-1}|\pr^{\leq \reg}\psi^{(p)}_{s,ij}|\big)+r^{-2}\big(|\pr\pr^{\leq \reg}(\widetilde{\phi}-\breve{\phi})^{(p)}_{s,kl}|+r^{-1}|\pr^{\leq \reg}(\widetilde{\phi}-\breve{\phi})^{(p)}_{s,kl}|\big)^2\Big)\\
\les& \M_{0, r\geq\varsigma_*(\tau-7c_*/8)}[\pr^{\leq\reg}(\widetilde{\pmb\phi}-\breve{\pmb\phi})](\tau_1+1, \tau_2-1)\M_{0, r\geq\varsigma_*(\tau-7c_*/8)}[\pr^{\leq\reg}\pmb\psi](\tau_1+1, \tau_2-1)\\
&+\Big(\M_{0, r\geq\varsigma_*(\tau-7c_*/8)}[\pr^{\leq\reg}(\widetilde{\pmb\phi}-\breve{\pmb\phi})](\tau_1+1, \tau_2-1)\Big)^2
\end{align*}
In view of \eqref{eq:localenergyestimateforpsisijponSigmatau2minus3tau2:partialderivatives} and \eqref{eq:controlofwidetildephiminusbreevephiontau2minus1tau2}, we deduce, for $\reg\leq 14$,
\begin{align}\lab{eq:estimateofKKregtowardstheendofgeneralizationproofhigherorderunweigthedenergyMorawetzinsection10dot5}
\KK_{\reg} \les& \ep\sum_{i,j=1}^3\sum_{p=0}^2\bigg(\E[\pr^{\leq \reg}\phi^{(p)}_{s,ij}](\tau_1)+\E[\pr^{\leq \reg}\phi^{(p)}_{s,ij}](\tau_2-3)+\F_{\Si_{*,c_*/2}}[\pr^{\leq\reg}\phi^{(p)}_{s,ij}](\tau_1+1, \tau_2-1)\nn\\
&+\int_{\MM(\tau_1, \tau_1+1)\cup\RR_*(\tau_1+1, \tau_2-3)\cup\MM(\tau_2-3, \tau_2-2)}r^{1+\de}|\pr^{\leq \reg}N_{W,s,ij}^{(p)}|^2\bigg).
\end{align}
The remaining terms in $\NNt_{s, \de, \text{total}}$ are estimated exactly as in section 10.4 of \cite{MaSz26} and we obtain the following analog of (10.98) in \cite{MaSz26}, for any $0<\la \leq 1$, all $\reg\leq 14$ and $\de\in (0,\frac{1}{3}]$, 
\bea
\lab{eq:controlofhighorderNNtspin-2:lastcontrol:pfofthm7.5}
&&\NNttotalph{s}{\pr^{\leq \reg}}\nn\\
&\les&
(\la+\ep)\EMFtotalhps{s}{\pr^{\leq \reg}} +\la^{-3}\sum_{p=0}^2\E[\pr^{\leq \reg}\phis{p}](\tau_1)
\nn\\
&&+\la^{-3}\sum_{p=0}^2{\mathcal{N}}_{\de}[\pr^{\leq \reg}\phis{p}, \pr^{\leq \reg}\N_{W, s}^{(p)}](\tau_1, \tau_2)+\la^{-1}\sum_{p=0}^1\int_{\MM(\tau_1,\tau_2)}r^{-3+\de}\big|\pr^{\leq \reg+1}\N_{T,s}^{(p)}\big|^2
\nn\\
&&+\sum_{p=0}^1\int_{\MM_{r\geq 10m}(\tau_1,\tau_2)} \Big(r^{-1+\de}\big|\pr^{\leq \reg}\N^{(p)}_{W,s}\big|+r^{-2+\de}\big|\pr^{\leq \reg+1}\N^{(p)}_{T,s}\big|\Big)\big|\pr^{\leq \reg}\phis{p}\big|\nn\\
&&+\ep\sum_{i,j=1}^3\sum_{p=0}^2\F_{\Si_{*,c_*/2}}[\pr^{\leq\reg}\phi^{(p)}_{s,ij}](\tau_1+1, \tau_2-1),
\eea
where the last term in \eqref{eq:controlofhighorderNNtspin-2:lastcontrol:pfofthm7.5} is new compared to (10.98) in \cite{MaSz26} and comes from the above estimate for $\KK_{\reg}$. 

\item Finally, we follow Section 10.5 in \cite{MaSz26}. Substituting the estimates \eqref{eq:localenergyforpsiontau1minus1tau1plus2:initialEMF}, \eqref{esti:Errdefectterm:highorderprderi}, \eqref{esti:Errdefectterm:highorderprderi:Fscase} and \eqref{eq:controlofhighorderNNtspin-2:lastcontrol:pfofthm7.5} back into \eqref{eq:EMFtotalp:finalestimate:pm2:prhighorder}, taking both of $\la$ and $\ep$ suitably small such that  the term $(\la+\ep)\EMFtotalhps{s}{\pr^{\leq \reg}}$ in \eqref{eq:controlofhighorderNNtspin-2:lastcontrol:pfofthm7.5} is absorbed by the LHS of \eqref{eq:EMFtotalp:finalestimate:pm2:prhighorder},  and noticing from the definition \eqref{def:AandAonorms:tau1tau2:phisandpsis} of the norms $\A[\cdot](\cdot, \cdot)$ and $\A_*[\cdot]$ that for $s=\pm 2$ we have
\beaa
\A[\pmb\psi_{s}](\Iti)+\A[\pmb\phi_{s}](\tau_1,\tau_2)
\les\sum_{p=0}^2\Big(\EMF[\psis{p}](\Iti)
+\EMF[\phis{p}](\tau_1, \tau_2)\Big),
\eeaa
and 
\beaa
\A_*[\pmb\psi_s]+\A_*[\pmb\phi_s] \les \sum_{p=0}^2\sup_{\frac{c_*}{2}\leq c\leq c_*}\Big(\F_{\Si_{*,c}}[\psis{p}](\tau_1+1, \tau_2-2)+\F_{\Si_{*,c}}[\phis{p}](\tau_1+1, \tau_2-2)\Big),
\eeaa
we deduce the following analog of (10.128) in \cite{MaSz26}
\begin{align}
\lab{eq:EMFtotalpspm2:finalestimatetoprovemaintheom:proof} 
&\EMFtotalh{s}{\nab_{(\pr_{\tt},\chi_{0}\pr_{\tilde{\phi}})}^{\leq \reg}}+\EMFtotalhps{s}{\pr^{\leq \reg}} +\sum_{p=0}^2\sup_{\frac{c_*}{2}\leq c\leq c_*}\F_{\Si_{*,c}}[\pr^{\leq \reg}\psis{p}](\tau_1+1, \tau_2-2)  \nn\\
\les&{\sum_{p=0}^2}\E[\pr^{\leq \reg}\phis{p}](\tau_1)  +\sum_{p=0}^2\Big(\EMF[\psis{p}](\Iti)
+\EMF[\phis{p}](\tau_1, \tau_2)\Big)
\nn\\
&+ \sum_{p=0}^2\sup_{\frac{c_*}{2}\leq c\leq c_*}\Big(\F_{\Si_{*,c}}[\psis{p}](\tau_1+1, \tau_2-2)+\F_{\Si_{*,c}}[\phis{p}](\tau_1+1, \tau_2-2)\Big)\nn\\
&+\sum_{p=0}^2{\mathcal{N}}_{\de}[\pr^{\leq \reg}\phis{p}, \pr^{\leq \reg}\N_{W, s}^{(p)}](\tau_1, \tau_2)
+\sum_{p=0}^1\int_{\MM(\tau_1,\tau_2)}r^{-3+\de}|\pr^{\leq \reg+1}\N_{T,s}^{(p)}|^2
\nn\\
&+\sum_{p=0}^1\int_{\MM_{r\geq 10m}(\tau_1,\tau_2)} \Big(r^{-1+\de}|\pr^{\leq \reg}\N^{(p)}_{W,s}|+r^{-2+\de}|\pr^{\leq \reg+1}\N^{(p)}_{T,s}|\Big)|\pr^{\leq \reg}\phis{p}|\nn\\
&+\sum_{p=0,1,2}\widetilde{\mathcal{N}}_*[\pr^{\leq \reg}\psis{p},\pr^{\leq \reg}\widehat{\F}_{total,s}^{(p)}] +\sum_{p=0}^{2}\widehat{\mathcal{N}}''[\pr^{\leq \reg}\psis{p}, \pr^{\leq \reg}\widehat{\F}_{total,s}^{(p)}](\tau_2-3, \tau_2)\nn\\
&+\sum_{p=0}^{2}\int_{\RR_*(\tau_1+1, \tau_2-2)}|\pr^{\leq\reg}\widehat{\F}_{total,s}^{(p)}|^2 +\ep\sum_{i,j=1}^3\sum_{p=0}^2\F_{\Si_{*,c_*/2}}[\pr^{\leq\reg}\phi^{(p)}_{s,ij}](\tau_1+1, \tau_2-1).
\end{align}
Next:
\begin{enumerate}
\item Proceeding as for $\KK_{\reg}$ in \eqref{eq:estimateofKKregtowardstheendofgeneralizationproofhigherorderunweigthedenergyMorawetzinsection10dot5}, we have 
\begin{align*}
&\sum_{p=0,1,2}\widetilde{\mathcal{N}}_*[\pr^{\leq \reg}\psis{p},\pr^{\leq \reg}\breve{\F}^{(p)}_{s}] +\sum_{p=0}^{2}\int_{\RR_*(\tau_1+1, \tau_2-2)}|\pr^{\leq\reg}\breve{\F}^{(p)}_{s}|^2\\
\les& \ep\sum_{i,j=1}^3\sum_{p=0}^2\bigg(\E[\pr^{\leq \reg}\phi^{(p)}_{s,ij}](\tau_1)+\E[\pr^{\leq \reg}\phi^{(p)}_{s,ij}](\tau_2-3)+\F_{\Si_{*,c_*/2}}[\pr^{\leq\reg}\phi^{(p)}_{s,ij}](\tau_1+1, \tau_2-1)\nn\\
&+\int_{\MM(\tau_1, \tau_1+1)\cup\RR_*(\tau_1+1, \tau_2-3)\cup\MM(\tau_2-3, \tau_2-2)}r^{1+\de}|\pr^{\leq \reg}N_{W,s,ij}^{(p)}|^2\bigg).
\end{align*}
Also, we have
\begin{align*}
&\sum_{p=0,1,2}\widetilde{\mathcal{N}}_*[\pr^{\leq \reg}\psis{p},\pr^{\leq \reg}\chi_{\tau_1, \tau_2}^{(1)}\chi_{r_*}^{(1)}N^{(p)}_{W,s,ij}] +\sum_{p=0}^{2}\int_{\RR_*(\tau_1+1, \tau_2-2)}|\pr^{\leq\reg}\chi_{\tau_1, \tau_2}^{(1)}\chi_{r_*}^{(1)}N^{(p)}_{W,s,ij}|^2\\
\les& \sum_{p=0}^2{\mathcal{N}}_{\de}[\pr^{\leq \reg}\phis{p}, \pr^{\leq \reg}\N_{W, s}^{(p)}](\tau_1, \tau_2)
\end{align*}
and hence
\begin{align*}
&\sum_{p=0,1,2}\widetilde{\mathcal{N}}_*[\pr^{\leq \reg}\psis{p},\pr^{\leq \reg}\widehat{\F}_{total,s}^{(p)}] +\sum_{p=0}^{2}\int_{\RR_*(\tau_1+1, \tau_2-2)}|\pr^{\leq\reg}\widehat{\F}_{total,s}^{(p)}|^2\\
\les&\ep\sum_{i,j=1}^3\sum_{p=0}^2\bigg(\E[\pr^{\leq \reg}\phi^{(p)}_{s,ij}](\tau_1)+\E[\pr^{\leq \reg}\phi^{(p)}_{s,ij}](\tau_2-3)+\F_{\Si_{*,c_*/2}}[\pr^{\leq\reg}\phi^{(p)}_{s,ij}](\tau_1+1, \tau_2-1)\nn\\
&+\sum_{p=0}^2{\mathcal{N}}_{\de}[\pr^{\leq \reg}\phis{p}, \pr^{\leq \reg}\N_{W, s}^{(p)}](\tau_1, \tau_2).
\end{align*}

\item Also, by the local energy estimate \eqref{eq:localenergyestimate:future:unweigthedderivatives:RR*} and \eqref{eq:causlityrelationsforwidetildepsi1}, we have, for $\reg\leq 14$, 
\begin{align*}
\nn&\sum_{p=0,1,2}\sup_{c_*/2\leq c\leq c_*}\F_{\Si_{*,c}}[\pr^{\leq\reg}\pmb\phi_s^{(p)}](\tau_1+1, \tau_2-3)\\
\les &\sum_{p=0}^2\bigg(\E[\pr^{\leq \reg}\pmb\phi^{(p)}_{s}](\tau_1+1) +\F_{\Si_{*,c_*/2}}[\pr^{\leq\reg}\pmb\psi_s^{(p)}](\tau_1+1, \tau_2-3)+ \int_{\RR_*(\tau_1, \tau_2)}r|\pr^{\leq\reg}\N_{W,s}^{(p)}|^2\bigg).
\end{align*}
\end{enumerate}
Combining the above estimates with \eqref{eq:EMFtotalpspm2:finalestimatetoprovemaintheom:proof}, we infer, for $\ep$ small enough, 
\begin{align}
\lab{eq:EMFtotalpspm2:finalestimatetoprovemaintheom:proof:bis} 
&\EMFtotalh{s}{\nab_{(\pr_{\tt},\chi_{0}\pr_{\tilde{\phi}})}^{\leq \reg}}+\EMFtotalhps{s}{\pr^{\leq \reg}}\nn\\
& +\sum_{p=0}^2\sup_{\frac{c_*}{2}\leq c\leq c_*}\F_{\Si_{*,c}}[\pr^{\leq \reg}\psis{p}](\tau_1+1, \tau_2-2)  \nn\\
\les&{\sum_{p=0}^2}\E[\pr^{\leq \reg}\phis{p}](\tau_1)  +\sum_{p=0}^2\Big(\EMF[\psis{p}](\Iti)
+\EMF[\phis{p}](\tau_1, \tau_2)\Big)
\nn\\
&+ \sum_{p=0}^2\sup_{\frac{c_*}{2}\leq c\leq c_*}\Big(\F_{\Si_{*,c}}[\psis{p}](\tau_1+1, \tau_2-2)+\F_{\Si_{*,c}}[\phis{p}](\tau_2-3, \tau_2-2)\Big)\nn\\
&+\sum_{p=0}^2{\mathcal{N}}_{\de}[\pr^{\leq \reg}\phis{p}, \pr^{\leq \reg}\N_{W, s}^{(p)}](\tau_1, \tau_2)
+\sum_{p=0}^1\int_{\MM(\tau_1,\tau_2)}r^{-3+\de}|\pr^{\leq \reg+1}\N_{T,s}^{(p)}|^2
\nn\\
&+\sum_{p=0}^1\int_{\MM_{r\geq 10m}(\tau_1,\tau_2)} \Big(r^{-1+\de}|\pr^{\leq \reg}\N^{(p)}_{W,s}|+r^{-2+\de}|\pr^{\leq \reg+1}\N^{(p)}_{T,s}|\Big)|\pr^{\leq \reg}\phis{p}|\nn\\
& +\sum_{p=0}^{2}\widehat{\mathcal{N}}''[\pr^{\leq \reg}\psis{p}, \pr^{\leq \reg}\widehat{\F}_{total,s}^{(p)}](\tau_2-3, \tau_2) +\ep\sum_{i,j=1}^3\sum_{p=0}^2\F_{\Si_{*,c_*/2}}[\pr^{\leq\reg}\phi^{(p)}_{s,ij}](\tau_2-3, \tau_2-1).
\end{align}
Next, using a slight variation of the local energy estimate \eqref{eq:localenergyestimate:future:unweigthedderivatives:bis}
 on $\MM(\tau_2-3, \tau_2-1)$, we have, for $\reg\leq 14$, 
\beaa
\sum_{i,j=1}^3\sum_{p=0}^2\F_{\Si_{*,c_*/2}}[\pr^{\leq\reg}\phi^{(p)}_{s,ij}](\tau_2-3, \tau_2-1)\les \sum_{p=0}^2\bigg(\E[\pr^{\leq \reg}\phis{p}](\tau_2-3)+\int_{\MM(\tau_2-3, \tau_2-1)}r^{1+\de}|\pr^{\leq\reg}\N_{W,s}^{(p)}|^2\bigg).
\eeaa
Combining with \eqref{eq:EMFtotalpspm2:finalestimatetoprovemaintheom:proof:bis} and taking $\ep>0$ sufficiently small, we infer, for $\reg\leq 14$, 
\begin{align*}
&\EMFtotalh{s}{\nab_{(\pr_{\tt},\chi_{0}\pr_{\tilde{\phi}})}^{\leq \reg}}+\EMFtotalhps{s}{\pr^{\leq \reg}} +\sum_{p=0}^2\sup_{\frac{c_*}{2}\leq c\leq c_*}\F_{\Si_{*,c}}[\pr^{\leq \reg}\psis{p}](\tau_1+1, \tau_2-2)  \nn\\
\les&{\sum_{p=0}^2}\E[\pr^{\leq \reg}\phis{p}](\tau_1)  +\sum_{p=0}^2\Big(\EMF[\psis{p}](\Iti)
+\EMF[\phis{p}](\tau_1, \tau_2)\Big)
\nn\\
&+ \sum_{p=0}^2\sup_{\frac{c_*}{2}\leq c\leq c_*}\F_{\Si_{*,c}}[\psis{p}](\tau_1+1, \tau_2-2)+\sum_{p=0}^{2}\widehat{\mathcal{N}}''[\pr^{\leq \reg}\psis{p}, \pr^{\leq \reg}\widehat{\F}_{total,s}^{(p)}](\tau_2-3, \tau_2)\nn\\
&+\sum_{p=0}^2{\mathcal{N}}_{\de}[\pr^{\leq \reg}\phis{p}, \pr^{\leq \reg}\N_{W, s}^{(p)}](\tau_1, \tau_2)
+\sum_{p=0}^1\int_{\MM(\tau_1,\tau_2)}r^{-3+\de}|\pr^{\leq \reg+1}\N_{T,s}^{(p)}|^2
\nn\\
&+\sum_{p=0}^1\int_{\MM_{r\geq 10m}(\tau_1,\tau_2)} \Big(r^{-1+\de}|\pr^{\leq \reg}\N^{(p)}_{W,s}|+r^{-2+\de}|\pr^{\leq \reg+1}\N^{(p)}_{T,s}|\Big)|\pr^{\leq \reg}\phis{p}|.
\end{align*}
We now control $\widehat{\mathcal{N}}''[\pr^{\leq \reg}\psis{p}, \pr^{\leq \reg}\widehat{\F}_{total,s}^{(p)}](\tau_2-3, \tau_2)$ as in (10.129) in \cite{MaSz26} and conclude as in the end of section 10.5 in \cite{MaSz26}.
\end{enumerate}
This ends the proof of Theorem \ref{th:main:intermediary}.
\end{proof}

%%%%%%%%%%%%%%%%%%%%%%%%%%%%%%%%%%%%%%%%%%%%%%%%%%

\subsection{Energy-Morawetz estimates for the scalar wave equation in perturbations of Kerr}
\lab{sec:energyMorawetzesitmatesforscalarwaveonMM:upto1derivative}

%%%%%%%%%%%%%%%%%%%%%%%%%%%%%%%%%%%%%%%%%%%%%%%%%%

In this section, we derive energy-Morawetz estimates for solutions to the inhomogeneous scalar wave equation on $(\MM, \g)$, i.e.
\begin{equation}
\label{eq:scalarwave}
\Box_{\g} \psi = F, \qquad \MM.
\end{equation}

\begin{theorem}
\label{thm:main:MaSz24}
Let $(\MM, \g)$ as in Section \ref{section:SpacetimeMM-chap6}, with $\g$ satisfying the assumptions \eqref{eq:controloflinearizedinversemetriccoefficients} \eqref{eq:controloflinearizedinversemetriccoefficients:inversegtautau}. Then, for $\ep>0$ small enough, we have for solutions to the inhomogeneous wave equation \eqref{eq:scalarwave} the following energy-Morawetz-flux estimates, for any $1\leq\tau_1<\tau_2\leq\tau_*$ and any $0<\de\leq 1$,
\bea
\label{MainEnerMora:psi}
\EMF^{(1)}_\de[{\psi}](\tau_1,\tau_2) \lesssim \E^{(1)}[\psi](\tau_1)+\mathcal{N}^{(1)}_\de[\psi, F](\tau_1, \tau_2),
\eea
where the norms $\EMF^{(1)}_\de[{\psi}](\tau_1,\tau_2)$, $\E^{(1)}[\psi](\tau_1)$ and $\mathcal{N}^{(1)}_\de[\psi, F](\tau_1, \tau_2)$ for a scalar field $\psi$ are as in Section \ref{subsection:basicnormsforpsi:MMh}, except for the Morawetz norm $\M_\de[\psi](\tau_1, \tau_2)$ which is given as follows
\beaa
\M_\de[\psi](\tau_1, \tau_2) &=& \int_{\MM(\tau_1, \tau_2)}\left(\frac{|\nab_3\psi|^2}{r^{1+\de}} +\frac{|\nabla\psi|^2}{r}\right)+\int_{\MM(\tau_1,\tau_2)}\left(\frac{|\nab_{\Rhat}\psi|^2}{r^{1+\de}} +\frac{|\psi|^2}{r^{3+\de}}\right),
\eeaa
and where the implicit constant in $\lesssim$ only depends on $a$, $m$, $\de$ and $\dec$.
\end{theorem}

Theorem \ref{thm:main:MaSz24} is the analog of Theorem 4.1 in \cite{MaSz24}. Its proof is a slight adaptation of the proof of Theorem 4.1 in \cite{MaSz24} from the case of the spacetime in \cite{MaSz24} which extends to $\II_+$ to the case of $\MM$ which only extends to the spacelike hypersuface $\Si_*$. As in the proof of Theorem \ref{thm:main:MaSz26}, the starting point is an extension to a semi-global problem using the extended metric \eqref{eq:extendedmetricgchitau1tau2}. The adaptation of the proof of Theorem 4.1 in \cite{MaSz24} to the one of Theorem \ref{thm:main:MaSz24} is similar, as well as significantly simpler, to the adaptation of the proof of Theorem 7.1 in \cite{MaSz26} to the one of Theorem \ref{thm:main:MaSz24} in Sections \ref{sect:maintheoremandproof} \ref{sec:proofofth:main:intermediary:00}. Details are left to the reader.


\begin{thebibliography}{99}

\bibitem{AMPW17} L.  Andersson, S.  Ma, C. Paganini and B. Whiting, \textit{Mode stability on the real axis}, J. Math. Phys. \textbf{58} (2017). 
 


\bibitem{Chand2} S. Chandrasekhar, \textit{On the equations governing the perturbations of the Schwarzschild black hole}, P. Roy. Soc. Lond. A Mat. \textbf{343} (1975), 289--298.
            
\bibitem{CB52} Y. Choquet-Bruhat, \textit{Th\'eor\'eme d'existence pour certains syst\'emes d'equations aux d\'eriv\'ees partielles non-lin\'eaires}, Acta Math. \textbf{88} (1952), 141--225.

 \bibitem{CBG69}
 Y. Choquet-Bruhat and R. Geroch, \textit{Global aspects of the Cauchy problem in general relativity}, Communications in Mathematical Physics \textbf{14} (1969),  329--335.

\bibitem{Ch-Kl} D. Christodoulou and S. Klainerman, \textit{The global nonlinear stability of the Minkowski space}, Princeton University Press (1993).
            
\bibitem{D-H-R} M. Dafermos, G. Holzegel and I. Rodnianski, \textit{Linear stability of the Schwarzschild solution to gravitational perturbations}, Acta Math. \textbf{222} (2019), 1--214.


\bibitem{D-H-R-Kerr} M. Dafermos, G. Holzegel and I. Rodnianski, \textit{Boundedness and decay  for the Teukolsky equation on Kerr spacetimes I: The case $|a|\ll M$}, Ann. PDE (2019), 118 pp.


\bibitem{DHRT} M. Dafermos,  G. Holzegel,  I. Rodnianski and M.  Taylor,   \textit{The non-linear stability of the Schwarzschild family of black holes}, arXiv:2104.0822.

\bibitem{GHP} R. Geroch,  A. Held and R.  Penrose, \textit{A space-time calculus based on pairs of null directions}, Journal of Mathematical Physics. \textbf{14} (1973), 874--881.

\bibitem{GKS20} E. Giorgi, S.  Klainerman and J.  Szeftel, \textit{A general formalism for the stability of Kerr}, arXiv:2002.02740.

\bibitem{GKS22} E. Giorgi, S.  Klainerman and J.   Szeftel, \textit{Wave equations estimates and the nonlinear stability of  slowly rotating Kerr black holes}, Pure and Applied Mathematics Quarterly, \textbf{20} (2024), no. 7, 2865--3849. 

\bibitem{Gu} S. Gurriaran, \textit{Non-linear instability of the Kerr Cauchy horizon near $i_+$}, arXiv:2603.17911. 

\bibitem{HeKl2} L. He and S. Klainerman, \textit{Generalized Whiting transform and flux estimates for wave equations in subextremal Kerr}, arXiv:2609.06344.

\bibitem{Hi2} P. Hintz, \textit{Constraint damping on subextremal Kerr spacetimes}, arXiv:2606.27658, 100 pp.

\bibitem{Hi3} P. Hintz, \textit{(Non-)Linear waves on asymptotically flat spacetimes. II: trapping, bound states, nonlinear applications}, arXiv:2606.28008, 285 pp.

\bibitem{Hi1} P. Hintz, \textit{Nonlinear stability of subextremal Kerr black holes}, arXiv:2606.28253, 339 pp.

\bibitem{Kerr63} R.  P. Kerr, \textit{Gravitational field of a spinning mass as an example of algebraically special metrics},  Physical Review Letters \textbf{11} (1963), 237--238.

\bibitem{Kl-review} S. Klainerman, \textit{The black hole stability problem}, Comptes Rendus M\'ecanique \textbf{353} (2025), 551--581. 

\bibitem{KS} S. Klainerman and J.  Szeftel, \textit{Global Non-Linear Stability of Schwarzschild Spacetime under Polarized Perturbations}, Annals of Math Studies, 210. Princeton University Press, Princeton NJ, 2020, xviii+856 pp. 

\bibitem{KS-GCM1}  S. Klainerman and J. Szeftel, \textit{Construction of GCM spheres in perturbations of Kerr}, Ann. PDE, \textbf{8} (2), Art. 17, 153 pp., 2022.

\bibitem{KS-GCM2}  S. Klainerman and J. Szeftel, \textit{Effective  results in uniformization and intrinsic GCM spheres in perturbations of Kerr}, Ann. PDE, \textbf{8} (2), Art. 18, 89 pp., 2022.

\bibitem{KS:Kerr} S. Klainerman and J. Szeftel, \textit{Kerr stability for small angular momentum}, Pure Appl. Math. Q., 19 (2023), 791--1678.

\bibitem{Ma} 
S. Ma, \textit{Uniform energy bound and Morawetz estimate for extreme components of spin fields in the exterior of a slowly rotating Kerr black hole II: linearized gravity}, Comm. Math. Phys. \textbf{377} (2020), 2489--2551.

\bibitem{MaSz24} S. Ma and J.   Szeftel, \textit{Energy-Morawetz estimates for the wave equation in perturbations of Kerr}, arXiv:2410.02341.

\bibitem{MaSz26} 
S. Ma and J. Szeftel, \textit{Energy-Morawetz estimates for Teukolsky equations in perturbations of Kerr}, arXiv:2603.23437.

\bibitem{Millet} P. Millet,  \textit{Optimal decay for solutions of the Teukolsky equation on the Kerr 
 metric for the full subextremal range $|a| < M$}, arXiv:2302.06946.

\bibitem{NP} E.  Newman and R. Penrose, \textit{An approach to gravitational radiation by a method of spin coefficients},  J. Math. Phys. \textbf{3} (1962), 566--578.

\bibitem{RW57} T. Regge and J. A. Wheeler, \textit{Stability of a Schwarzschild singularity}, Phys. Rev. (2), 108:1063--1069, 1957.

\bibitem{Sch16} K. Schwarzschild, \"Uber das Gravitationsfeld eines Massenpunktes nach der Einsteinschen Theorie. Sitzungsber. Preuss. Akad. Wiss. \textbf{1} (1916),  189--196. 

\bibitem{Shen} D. Shen, \textit{Construction of GCM hypersurfaces in perturbations of Kerr}, Annals of PDE, \textbf{9} (1), Art. 11, 2023.

\bibitem{SRTdC20} Y.  Shlapentokh-Rothman and R. Teixeira da Costa, \textit{Boundedness and decay for the Teukolsky equation on Kerr in the full subextremal range $|a|< M$: frequency space analysis}, ARMA \textbf{250} Art. 78, 157 pp, 2026. 

\bibitem{SRTdC23} Y.  Shlapentokh-Rothman and R. Teixeira da Costa, \textit{Boundedness and decay for the Teukolsky equation on Kerr in the full subextremal range $|a|< M$: physical  space analysis}, arXiv:2302.08916, to appear in Analysis \& PDE.

\bibitem{Sze} J. Szeftel, \textit{Kerr stability in the full subextremal range}, arXiv:2609.40068.

\bibitem{TdC20} R. Teixeira da Costa, \textit{Mode stability for the Teukolsky equation on extremal and subextremal Kerr spacetimes}, Comm. Math. Phys. \textbf{378} (2020), 705--781.

\bibitem{Teuk} S. A. Teukolsky, \textit{Perturbations of a rotating black hole. I. Fundamental equations for gravitational, electromagnetic, and neutrino-field perturbations}, Astrophys. J. \textbf{185} (1973), 635--648.
 
\bibitem{Whit}  B. Whiting, \textit{Mode stability of the Kerr black hole},  J. Math. Phys. \textbf{30} (1989), 1301--1305.
\end{thebibliography}
\end{document}